\documentclass{amsart}
\usepackage[T1]{fontenc}
\usepackage{bm}
\usepackage{microtype} %

\usepackage[style=ext-alphabetic,doi=false,isbn=false,maxalphanames=99,maxnames=99,articlein=false,backend=biber]{biblatex}
\DeclareFieldFormat{language}{}

\renewbibmacro*{language}{}

\AtEveryBibitem{%
  \clearfield{language}%
  \ifboolexpr{
    test {\ifentrytype{article}} 
      or test {\ifentrytype{book}}
            or test {\ifentrytype{incollection}}
  }{%
    \clearfield{url}%
  }{}%
  \iffieldundef{eprint}
    {}
    {\clearfield{url}}
}
  \DefineBibliographyStrings{english}{page = {}, pages = {}}
\DefineBibliographyStrings{french}{page = {}, pages = {}}
\DefineBibliographyStrings{german}{page = {}, pages = {}}

\usepackage{tikz}
\usepackage{tikz-cd}
\usepackage{quiver}
\usepackage{adjustbox}
\usepackage{enumitem}
\usepackage[normalem]{ulem}
 \usepackage{stmaryrd}
 \usepackage{mathtools}

\newcommand{\AKModfpA}{A\llbracket K\rrbracket\text{-}\Mod^{\fp}(A)}
\newcommand{\OKtorsModfpA}{\cO\llbracket K\rrbracket_{\mf{a}\text{-}\mathrm{tors}}\text{-}\Mod^{\fp}(A)}
\newcommand{\OKProModfpA}{\cO\llbracket K\rrbracket\text{-}\Pro\Mod^{\fp}(A)}

\newcommand{\HeckeT}{T_p}
\newcommand{\HeckeX}{T_p}

\newcommand{\Cflat}{\C^\flat}
\newcommand{\piss}{\pi_{\semis}}
\newcommand{\catA}{\cA_{U}}
\newcommand{\catB}{\cA_{Z}}
\newcommand{\catD}{\cA}
\newcommand{\uf}{}             %
\newcommand{\mf}{\mathfrak}
\newcommand{\supp}{\operatorname{supp}}

\newcommand{\ProCoh}{\Pro\Coh}
\newcommand{\uHom}{\underline{\Hom}}
\newcommand{\mSpec}{\operatorname{mSpec}}
\newcommand{\Rees}{\operatorname{Rees}}

\newcommand{\Rtilde}{\tld{R}}

\newcommand{\cXhat}{\widehat{\cX}}
\newcommand{\ihat}{\widehat{i}}
\newcommand{\jhat}{\widehat{j}}
\newcommand{\vhat}{\widehat{v}}
\newcommand{\boxdelta}{\boxtimes}

\newcommand{\ResQp}{\Res_{\Qp}}
\renewcommand{\mathbb}{\mathbf}
\newcommand{\chibar}{\overline{\chi}}

\newcommand{\zetabar}{\overline{\zeta}}
\newcommand{\thetabar}{\overline{\theta}}

\newcommand{\sigmacomp}{\sigma^{\operatorname{co}}}
\newcommand{\sigmasigmacomp}{\sigma|\sigmacomp}
\newcommand{\lambdau}{\underline{\lambda}}
\def\crys{\mathrm{crys}}

\def\Dfp{\mathop{D}_{\mathrm{fp}}\nolimits}
\def\coh{\mathrm{coh}}

\def\coh{\mathrm{coh}}
\def\qc{\mathrm{qc}}

\newcommand{\semis}{{\operatorname{ss}}}

\newcommand{\fQ}{\mathfrak{Q}}
\newcommand{\fV}{\mathfrak{V}}
\newcommand{\fW}{\mathfrak{W}}
\newcommand{\fM}{\mathfrak{M}}
\newcommand{\fN}{\mathfrak{N}}
\newcommand{\alphabar}{\overline{\alpha}}
\newcommand{\adm}{\mathrm{adm}}
\newcommand{\coadm}{\mathrm{coadm}}
\newcommand{\tld}[1]{\widetilde{#1}}
\newcommand{\wht}[1]{\widehat{#1}}

\newcommand{\Supp}{\mathop{\mathrm{Supp}}\nolimits}
\newcommand{\Vsigma}{M_{\sigma, \thetabar}}
\newcommand{\Vver}{V^{\ver}}
\newcommand{\Rsigma}{R^{\sigma}}

\newcommand{\Functor}{\mathrm{F}} %
\newcommand{\ProFunctor}{\Pro\mathrm{F}} %
\newcommand {\imm}{k} %

\newcommand{\ps}{\operatorname{ps}}
\newcommand{\Pone}{\mathbb{P}^1}

\newcommand{\Zptimes}{\Zp^\times}

\newcommand{\Qptimes}{\Qp^\times}
\newcommand{\tr}{\operatorname{tr}}

 \newcommand{\Gammat     }{\widetilde{\Gamma}}

\newcommand{\Iw}{\mathrm{Iw}}
\newcommand{\ver}{\mathrm{ver}}
\newcommand{\mss}{\operatorname{ss}}
\newcommand{\cosoc}{\operatorname{cosoc}}

\newcommand{\soc}{\operatorname{soc}}
\newcommand{\Rad}{\operatorname{rad}}
\newcommand{\rad}{\operatorname{rad}}

\newcommand{\Ybad}{Y_{\rm{bad}}}
\newcommand{\cYbad}{\cY_{\rm{bad}}}
\newcommand{\Ugood}{U_{\rm{good}}}
\newcommand{\cUgood}{\cU_{\rm{good}}}
\newcommand{\cUgoodred}{\cU_{\rm{good},\red}}

\usepackage{amssymb,amsmath,amsfonts,amsthm,epsfig,amscd}
\usepackage[all,cmtip,poly]{xy}
\usepackage{color}

\newcommand{\lbar}[1]{\overline{#1}}
\newcommand{\diag}{\mathrm{diag}}
\newcommand{\fourmatrix}[4]{\begin{pmatrix} #1 & #2 \\ #3 & #4 \end{pmatrix}}
\newcommand{\isom}{\xrightarrow{\sim}}
\newcommand{\sm}{{\mathrm{sm}}}
\newcommand{\triv}{\mathbf{1}}
\newcommand{\fm}{\mathfrak{m}}
\newcommand{\fn}{\mathfrak{n}}
\newcommand{\fa}{\mathfrak{a}}

\newcommand{\cbF}{\overline{\bF}}
\newcommand{\ad}{\operatorname{ad}}

\newcommand{\gr}{\mathrm{gr}}

\newcommand{\red}{\operatorname{red}}

\newcommand{\limcommand}{\operatorname{lim}}
\newcommand{\quoteslim}[1]{\operatorname{lim}_{#1}}
\newcommand{\quotescolim}[1]{\operatorname{colim}_{#1}}

\newcommand{\colim}{\operatorname{colim}}

\newcommand{\Ann}{\operatorname{Ann}}

\newcommand{\fB}{\mathfrak{B}}
\newcommand{\fC}{\mathfrak{C}}
\newcommand{\fF}{\mathfrak{F}}
\newcommand{\fX}{\mathfrak{X}}

\newcommand{\univ}{{\operatorname{univ}}}
\newcommand{\To}{\longrightarrow}
\newcommand{\isoto}{\stackrel{\sim}{\To}}
\newcommand{\cotimes}{\, \widehat{\otimes}}
\newcommand{\ccotimes}{\, \widehat{\otimes}}

\newcommand{\cl}{\operatorname{cl}}
\newcommand{\fg}{\operatorname{fg}}
\newcommand{\fp}{\operatorname{fp}}
\newcommand{\ft}{\operatorname{ft}}

\newcommand{\fl}{\operatorname{f.l.}}

\newcommand{\sob}{\operatorname{sob}}
\newcommand{\Rep}{\operatorname{Rep}}

\newcommand{\JH}{\operatorname{JH}}
 \newcommand{\sigmabar   }{\overline{\sigma}}   
\def\iso{\buildrel \sim \over \longrightarrow}

\newcommand{\id}{\operatorname{id}}

\newcommand{\cont}{\operatorname{cont}}

\usepackage{ifpdf}
\ifpdf \usepackage[bookmarksopen,bookmarksdepth=4]{hyperref} \fi

\usepackage{chngcntr}
\counterwithin{equation}{subsection}

\makeatletter
\let\oldsubsubsection\subsubsection

\renewcommand{\subsubsection}{\@ifstar{\subsubsectionstar}{\newsubsubsection}}

\newcommand{\subsubsectionstar}{\oldsubsubsection*}

\newcommand{\newsubsubsection}[1]{%
  \refstepcounter{equation}%
  \@startsection{subsubsection}{3}%
  {\z@}{.5\linespacing\@plus.7\linespacing}{-.5em}%
  {\normalfont\itshape}{#1}%
}
\makeatother

\newtheorem{theorem}[equation]{Theorem}
\newtheorem{thm}[equation]{Theorem}
\newtheorem{lemma}[equation]{Lemma}
\newtheorem{lem}[equation]{Lemma}

\newtheorem{cor}[equation]{Corollary}

\newtheorem{prop}[equation]{Proposition}

\newtheorem{athm}[equation]{Theorem}

\newtheorem{alem}[equation]{Lemma}
\newtheorem{adefn}[equation]{Definition}

\newtheorem{acor}[equation]{Corollary}
\newtheorem{aprop}[equation]{Proposition}

\theoremstyle{definition}
\newtheorem{df}[equation]{Definition}
\newtheorem{defn}[equation]{Definition}

\theoremstyle{remark}
\newtheorem{remark}[equation]{Remark}
\newtheorem{rem}[equation]{Remark}
\newtheorem{arem}[equation]{Remark}
\newtheorem{aremark}[equation]{Remark}

\newtheorem{example}[equation]{Example}

\newtheorem{hypothesis}[equation]{Hypothesis}

\newif\iffinalrun
\iffinalrun
\else
\fi

\iffinalrun
  \newcommand{\need}[1]{}
  \newcommand{\mar}[1]{}
\else
  \newcommand{\need}[1]{{\tiny *** #1}}
  \newcommand{\mar}[1]{\marginpar{\raggedright\tiny  %
		  fixme
      #1}}
\fi

\iffinalrun
  \newcommand{\finalmar}[1]{}
\else
  \newcommand{\finalmar}[1]{\marginpar{\raggedright\tiny   %
      #1}}
\fi

\newcommand{\A}{\AA}
\newcommand{\C}{\CC}
\newcommand{\F}{\FF}

\newcommand{\Q}{\QQ}
\newcommand{\R}{\RR}
\newcommand{\Z}{\ZZ}

\newcommand{\m}{\frakm}

\renewcommand{\AA}{{\mathbb A}}

\newcommand{\CC}{{\mathbb C}}

\newcommand{\FF}{{\mathbb F}}
\newcommand{\GG}{{\mathbb G}}

\newcommand{\QQ}{{\mathbb Q}}
\newcommand{\RR}{{\mathbb R}}

\newcommand{\ZZ}{{\mathbb Z}}

\newcommand{\bA}{\ensuremath{\mathbf{A}}}

\newcommand{\bD}{\ensuremath{\mathbf{D}}}
\newcommand{\bE}{\ensuremath{\mathbf{E}}}
\newcommand{\bF}{\ensuremath{\mathbf{F}}}
\newcommand{\bG}{\ensuremath{\mathbf{G}}}

\newcommand{\bP}{\ensuremath{\mathbf{P}}}
\newcommand{\bQ}{\ensuremath{\mathbf{Q}}}

\newcommand{\bV}{\ensuremath{\mathbf{V}}}

\newcommand{\bZ}{\ensuremath{\mathbf{Z}}}

\newcommand{\cA}{{\mathcal A}}
\newcommand{\cB}{{\mathcal B}}
\newcommand{\cC}{{\mathcal C}}
\newcommand{\cD}{{\mathcal D}}

\newcommand{\cF}{{\mathcal F}}
\newcommand{\cG}{{\mathcal G}}
\newcommand{\cH}{{\mathcal H}}
\newcommand{\cI}{{\mathcal I}}

\newcommand{\cK}{{\mathcal K}}

\newcommand{\cO}{{\mathcal O}}

\newcommand{\cQ}{{\mathcal Q}}

\newcommand{\cT}{{\mathcal T}}
\newcommand{\cU}{{\mathcal U}}
\newcommand{\cV}{{\mathcal V}}
\newcommand{\cW}{{\mathcal W}}
\newcommand{\cX}{{\mathcal X}}
\newcommand{\cY}{{\mathcal Y}}
\newcommand{\cZ}{{\mathcal Z}}
\newcommand{\tE}{{\widetilde E}}

\newcommand{\tP}{{\widetilde P}}

\newcommand{\tR}{{\widetilde R}}
\newcommand{\tX}{{\widetilde X}}

\newcommand{\frakm}{\mathfrak{m}}

\newcommand{\fT}{\mathfrak{T}}

\newcommand{\gM}{\mathfrak{M}}
\newcommand{\gN}{\mathfrak{N}}

\newcommand{\Fbar}{\overline{\F}}
\newcommand{\Qbar}{\overline{\Q}}

\newcommand{\Fp}{\F_p}

\newcommand{\Fpbar}{\Fbar_p}

\newcommand{\Fpbartimes}{\Fpbar^{\times}}

\newcommand{\Zp}{\Z_p}

\newcommand{\Qp}{\Q_p}

\newcommand{\Qpbar}{\Qbar_p}

\DeclareMathOperator{\Aut}{Aut}

\DeclareMathOperator{\CoCh}{CoCh}
\DeclareMathOperator{\Coh}{Coh}
\DeclareMathOperator{\QCoh}{QCoh}

\DeclareMathOperator{\coker}{coker}

\DeclareMathOperator{\End}{End}
\DeclareMathOperator{\REnd}{REnd}
\DeclareMathOperator{\Ext}{Ext}

\DeclareMathOperator{\Fil}{Fil}
\DeclareMathOperator{\Fun}{Fun}
\DeclareMathOperator{\Gal}{Gal}
\DeclareMathOperator{\GL}{GL}

\DeclareMathOperator{\Hom}{Hom}
\DeclareMathOperator{\RHom}{RHom}

\DeclareMathOperator{\Rlim}{Rlim}
\DeclareMathOperator{\Maps}{Maps}

\DeclareMathOperator{\Tor}{Tor}
\DeclareMathOperator{\im}{im}
\DeclareMathOperator{\Ind}{Ind}

\DeclareMathOperator{\Mod}{Mod}

\DeclareMathOperator{\PGL}{PGL}

\DeclareMathOperator{\Pro}{Pro}

\DeclareMathOperator{\SL}{SL}

\DeclareMathOperator{\Sp}{Sp}
\DeclareMathOperator{\Spec}{Spec}

\DeclareMathOperator{\Spf}{Spf}

\DeclareMathOperator{\Sym}{Sym}

\DeclareMathOperator{\Tr}{Tr}

\newcommand\cInd{c\text{-}\!\Ind}
\newcommand{\ab}{\mathrm{ab}}

\newcommand{\Frac}{\mathrm{Frac}}

\newcommand{\Frob}{\mathrm{Frob}}

\newcommand{\HT}{\mathrm{HT}}

\newcommand{\ladm}{\mathrm{l.adm}}
\newcommand{\nr}{\mathrm{nr}}

\newcommand{\St}{\mathrm{St}}

\newcommand{\Zar}{\mathrm{Zar}}

\newcommand{\rhobar}{\overline{\rho}}

\newcommand{\etale}{\'{e}tale~}

\newcommand{\op}{\mathrm{op}}

\newcommand{\into}{\hookrightarrow}
\newcommand{\onto}{\twoheadrightarrow}

\newcommand{\Gm}{\GG_m}

\newcommand{\Ver}{\mathrm{Ver}}

\newcommand{\varepsilonbar}{\overline{\varepsilon}}

\newcommand{\Res}{\operatorname{Res}}

\DeclareMathOperator{\Vcheck}{\check{V}}

\usepackage{crossreftools}

\begin{document}

\title[A categorical $p$-adic Langlands correspondence for~$\mathrm{GL}_2(\mathbf{Q}_p)$]{A categorical $\bm{p}$-adic Langlands correspondence for~$\bm{\mathrm{GL}_2(\mathbf{Q}_p)}$}

\author[A. Dotto]{Andrea Dotto}\email{andrea.dotto@kcl.ac.uk}\address{Department of Mathematics, King's College London,
London~WC2R 2LS,~UK}
\author[M. Emerton]{Matthew Emerton}\email{emerton@math.uchicago.edu}
\address{Department of Mathematics, University of Chicago,
5734 S.\ University Ave., Chicago, IL 60637, USA}
\author[T. Gee]{Toby Gee} \email{toby.gee@imperial.ac.uk} \address{Department of
  Mathematics, Imperial College London,
  London SW7 2AZ,~UK}
\thanks{AD was supported in various stages of this project by the Engineering and Physical Sciences
  Research Council [EP/L015234/1] (The London School of Geometry and Number Theory),
  the James D. Wolfensohn Fund at the Institute for Advanced Study,
  and a Royal Society University Research Fellowship. ME was supported in part by the
  NSF grants %
  DMS-1601871, DMS-1902307, DMS-1952705, and~DMS-2201242. TG was 
  supported in part by an ERC Advanced grant. This project has received funding from the European Research Council
  (ERC) under the European Union’s Horizon 2020 research and
  innovation programme (grant agreement No. 884596). ME and TG were
  both supported in part by the Simons Collaboration on Perfection in Algebra, Geometry and Topology.}

\begin{abstract}%
  Let $p\ge 5$ be a prime. We %
 construct  a fully faithful functor from the
  derived category of all smooth $p$-adic representations of
  $\GL_2(\Qp)$ (with a fixed central character) to a derived category
  of Ind-coherent sheaves on a stack of
  $(\varphi,\Gamma)$-modules. %
\end{abstract}
\maketitle
\setcounter{tocdepth}{2}

\tableofcontents

\section{Introduction}

\subsection{Overview of our results}
We fix a prime $p\ge 5$. Our goal in this paper is to
generalise the $p$-adic local Langlands correspondence for
$\GL_2(\Qp)$.
The existing form of the %
correspondence
(as initiated by Breuil, established in full generality
by Colmez,
 and then promoted to an
equivalence of categories by Pa\v{s}k\={u}nas)
relates locally admissible representations of~$\GL_2(\Qp)$ to finitely
generated modules over various kinds of deformation rings for two-dimensional
mod~$p$ representations of the absolute Galois
group~$\Gal(\Qpbar/\Qp)$. 
Motivated by the ``categorical'' Langlands
program, our generalisation %
relates arbitrary ($p$-power torsion) smooth representations
to complexes of (Ind-) coherent sheaves on a moduli stack of $(\varphi,\Gamma)$-modules.

Our results were announced in the survey~\cite{emerton2023introduction}, and we refer the reader to that paper for
an extensive introduction and motivation; here, we content ourselves with emphasizing
that
the category of smooth representations is significantly %
larger than that of locally admissible representations: for example,
it contains the compactly supported inductions of finite-dimensional
representations of any open and compact (modulo centre) subgroup
of~$\GL_2(\Qp)$, and these are never locally admissible. %

To state our results more precisely, we introduce some notation. 
Let $\cO$ denote the ring of integers in a finite
extension $E$ of $\Q_p$, let~$\F$ be the residue field of~$\cO$, and let $\cA$ denote the
category of smooth representations of $G\coloneqq  \GL_2(\Q_p)$ on $\cO$-modules
on which $p$ acts locally nilpotently, and which have a central
character equal to some fixed character
$\zeta: \Q_p^{\times} \to \cO^{\times}$.
Let $\Dfp^b(\cA)$ denote the full sub-$\infty$-category of the bounded derived
category $D^b(\cA)$ of~$\cA$ consisting of objects with finitely presented cohomology.

Let $\cX$ be the formal algebraic stack over~$\Spf \cO$ which parameterizes rank $2$ projective
\'etale $(\varphi,\Gamma)$-modules %
with fixed determinant~$\zeta\varepsilon^{-1}$, where $\varepsilon$ denotes the
$p$-adic cyclotomic character.
Finally, let $D_{\coh}^b(\cX)$ denote the bounded derived category of Ind-coherent %
sheaves on~$\cX$ with coherent cohomologies. 
(We refer %
to Section~\ref{subsec:formal coherent sheaves} for a precise definition of $D^b_{\coh}(\cX)$.) 

Our main theorem is then the following.

\begin{thm}[Definition~ \ref{defn:definition of F} and Theorem~ \ref{thm:F-is-fully-faithful}]
\label{thm:coherent sheaves intro} There exists an $\cO$-linear fully
faithful exact functor %
$\Functor:\Dfp^b(\cA) \hookrightarrow D_{\coh}^b(\cX).$ 
\end{thm}

Using our assumption that~$p\ge 5$, we show that the unbounded derived category is compactly generated by $\Dfp^b(\cA)$, so for formal reasons, any functor~$\Functor$ 
as in the statement of Theorem~\ref{thm:coherent sheaves intro}
extends to a continuous (i.e.\ colimit-preserving) fully faithful
functor 
\[
\Functor: D(\cA) %
\to \Ind D^b_{\coh}(\cX). 
\]
Thus
Theorem~\ref{thm:coherent sheaves intro} in fact gives a correspondence that
takes account of all smooth representations (and not just the finitely
presented ones).

\begin{rem}\label{rem: more about F}
Although Theorem~\ref{thm:coherent sheaves intro} is only phrased as an existence result, 
the functor~$\Functor$ that we construct in this paper enjoys several additional properties.
For example, Proposition~\ref{finiteness of cohomology} proves that~$\Functor$ has amplitude~$[-1, 0]$, and
Proposition~\ref{concentration and support} describes some properties 
of the image under~$\Functor$ of a natural set of compact generators of~$D(\cA)$
(the compact inductions of Serre weights).
These results are in line with the general expectations described in~\cite[Section~6]{emerton2023introduction}.

The functor~$\Functor$ is %
not an equivalence, and a sketch of the computation of its essential image is given in~\cite[Section~7.5.5]{emerton2023introduction}.
We intend for the full details of this computation to appear elsewhere.
\end{rem}

\subsubsection{Constructing the functor}
We now describe our construction of the functor~$\Functor$. %
We let $\cO\llbracket G\rrbracket _{\zeta}$ denote the quotient of the completed group
ring $\cO\llbracket G\rrbracket $ on which the centre acts via~$\zeta$. 
Then
general principles~\cite[Rem.\ 6.1.23]{emerton2023introduction}
suggest
 that $\Functor$ must be of the form
\begin{equation}
\label{eqn:kernel formula}
\Functor(\pi) = L_{\infty}\otimes^L_{\cO\llbracket G\rrbracket _{\zeta}} \pi,
\end{equation}
for a pro-coherent sheaf $L_{\infty}$ of right $\cO\llbracket G\rrbracket _{\zeta}$-modules
on~$\cX$.
The first step in proving our theorem %
is to directly construct a suitable pro-coherent sheaf~$L_\infty$, 
and then use the formula~\eqref{eqn:kernel formula} to define the functor~$\Functor$. 

Defined in this way, $\Functor$ {\em a priori} takes values in $\Pro D^b_{\coh}(\cX)$,
and so one of our problems is to show that $\Functor$ actually takes values in $D^b_{\coh}(\cX)$.
We then furthermore have to show that $\Functor$ is fully faithful.
Two key techniques used in the solution of these problems are those of {\em localization}
and {\em completion},
both in the geometric sense of localizing or completing at points of~$\cX$, and in a more categorical
sense, in which we pass from $\cA$ to various of its quotient categories,
or to blocks of its full subcategory
$\cA^{\ladm}$ of locally admissible representations.

\subsubsection{Constructing the kernel}
Our construction of $L_{\infty}$ uses
Colmez's constructions from~\cite{MR2642409},
extended to the context of $(\varphi,\Gamma)$-modules
with coefficients in an arbitrary Noetherian $\cO/\varpi^a$-algebra~$A$.
This extension is not formal, and it is the subject of Section~\ref{sec: Colmez}, which includes 
a systematic study of the base change properties of Colmez's functor $D \mapsto D^\natural \boxtimes \Pone$.

We then construct~$L_\infty$ by descent in Section~\ref{subsec: the functor}.
To be precise,
if $\Spf A \to \cX$ 
is a smooth chart, where $A$ is an $I$-adically complete Noetherian $\cO$-algebra,
then we may form the corresponding 
pro-system of \'etale $(\varphi,\Gamma)$-modules~$D_{A/I^n}$.
Via our extension of Colmez's constructions,
we obtain a pro-system
$D_{A/I^n}^{\natural}\boxtimes \Pone$ 
of $\cO\llbracket G\rrbracket $-modules on~$\Spf A$
which admits descent data from $\Spf A$ to~$\cX$.
Descending it (and twisting it appropriately),
we obtain a pro-coherent sheaf of right $\cO\llbracket G\rrbracket _{\zeta}$-modules $L_{\infty}$
on~$\cX$.

\subsubsection{Localization}
The localization theory of~\cite{DEGlocalization} can be briefly summarized as follows. 
In~\cite{MR3150248},
Pa\v{s}k\={u}nas shows that the blocks of $\cA^{\ladm}$
are in bijection with $\Gal(\Fbar/\F)$-conjugacy classes of
$2$-dimensional continuous $\Fbar$-valued pseudorepresentations $\thetabar$
of $G_{\Q_p}$ having determinant~$\zeta\varepsilon^{-1}$,
and we typically denote the block of $\cA^{\ladm}$ corresponding to
such a closed point by~$\cA_{\thetabar}$.
In~\cite{DEGlocalization} we %
explain that these $\thetabar$ are parameterized by (the closed points of)
a scheme $X$ which is a chain of projective lines
over~$\F$, and we show that the category $\cA$ localizes in a natural way
to a sheaf of abelian categories over~$X$.  
Indeed, if $Y$ is a closed subset of~$|X|$, with open complement~$U$,
then in~\cite{DEGlocalization} we define
a certain localizing subcategory $\cA_Y$ of~$\cA$ associated to~$Y$,
and set~$\cA_U \coloneqq  \cA/\cA_Y.$
In the case when $Y$ is finite, $\cA_Y$ consists, by construction, of locally admissible
representations. In particular, when $Y = \{\thetabar\}$ is a single closed point,
then (again by construction) $\cA_Y$ is the block~$\cA_{\thetabar}$.

By design, the process of passing to~$\cA_U$ is related to localization of sheaves on~$\cX$, 
because the closed points of~$|\cX|$ are {\em also} in bijection with the 
(Galois conjugacy classes of) residual
pseudorepresentations~$\thetabar$, and
the resulting bijection
between the closed points of $|\cX|$ and $|X|$ extends
to a continuous morphism $\piss:|\cX| \to |X|$,
characterized by the following condition:
 for each pseudorepresentation~$\thetabar$,  regarded as a closed point of~$|X|$, %
the completion  $\cX_{\thetabar}$ of $\cX$ along its closed
subset~$\piss^{-1}(\thetabar)$
is
the moduli stack of two-dimensional continuous Galois representations  
having determinant~$\zeta \varepsilon^{-1}$ whose underlying pseudorepresentations
deform~$\thetabar$.
(The stack~$\cX_{\thetabar}$ is an example of the moduli stacks of Galois representations constructed and studied by Wang-Erickson in~\cite{MR3831282}.)

\subsubsection{Establishing properties of~\texorpdfstring{$\Functor$}{F}}
As already remarked, having defined~$L_{\infty}$, we are 
able to define a functor
\[
\Functor \coloneqq  L_\infty \otimes^L_{\cO\llbracket G\rrbracket _\zeta}: \Dfp^b(\cA) \to \Pro D^b_{\coh}(\cX).
\]
We then prove the following theorem, which collects a number of results from Section~5 of this paper.

\begin{thm}
  \label{thm:intro-summary-of-Y-U-props}\leavevmode
  \begin{enumerate}%
   \item\label{item:intro finiteness of F} {\em (Prop.\ \ref{finiteness of cohomology}.)} The functor $\Functor$ factors through the full subcategory
$D^b_{\coh}(\cX)$ of~$\Pro D^b_{\coh}(\cX)$, and so may be regarded as a functor
$$\Functor: \Dfp^b(\cA) \to D^b_{\coh}(\cX).$$
\item
{\em (Prop.\ \ref{prop:existence of FY FU overview}.)}\label{item:intro FY FU}
For any closed subset $Y$ of~$X$, with open complement~$U$,
let $\cX_Y$ denote the completion of $\cX$ along~$\pi^{-1}_{\mathrm{ss}}(Y)$, and let $\cU \coloneqq  \pi^{-1}_{\mathrm{ss}}(U)$.
Then the functor $\Functor$ 
induces a commutative diagram defining functors~$\Functor_Y$ and~$\Functor_U$
\begin{equation}\label{eqn:intro recollement diagram}
\begin{tikzcd}
0 \arrow[r] & \Dfp^b(\mathcal{A}_Y) \arrow[r] \arrow[d, "\Functor_Y"] & \Dfp^b(\mathcal{A}) \arrow[r] \arrow[d, "\Functor"] & \Dfp^b(\mathcal{A}_U) \arrow[r] \arrow[d, "\Functor_U"] & 0 \\
0 \arrow[r] & D^b_{\coh}(\cX_Y) \arrow[r] & D^b_{\coh}(\mathcal{X}) \arrow[r] & D^b_{\coh}(\mathcal{U}) \arrow[r] & 0.
\end{tikzcd}
\end{equation}
   
\item\label{item:intro FY fully faithful}
{\em (Cor.\ \ref{cor:overview fully faithful on closed}.)} 
For any finite set of closed points $Y \subset X$,
the functor $\Functor_Y$ is fully faithful. 
\item\label{item:intro FU fully faithful} {\em (Thm.\ \ref{thm:full faithfulness on U}.)} There is a dense open subset $U \subset X$ such that
the functor $\Functor_U$ is fully faithful.  
  \end{enumerate}
\end{thm}

We deduce Theorem~\ref{thm:coherent sheaves intro}
from Theorem~\ref{thm:intro-summary-of-Y-U-props}.
More precisely, part~\eqref{item:intro finiteness of F} of the latter theorem is precisely the statement that
$\Functor$ takes values in~$D^b_{\coh}(\cX)$; while the full faithfulness
of~$\Functor$ follows by combining the full faithfulness
results of parts~\eqref{item:intro FY fully faithful} and~\eqref{item:intro FU fully faithful}
of Theorem~\ref{thm:intro-summary-of-Y-U-props}
with
general category-theoretic arguments.  
These arguments are developed in Appendix~\ref{subsec:semiorthogonal and Ind Pro};
see especially Proposition~\ref{prop:gluingfull-faithfulness-Ind-Pro-version-compact-version}, 
which provides a criterion for full faithfulness
of a functor that preserves a recollement-like decomposition of its source and target categories.
Part~\eqref{item:intro FY FU} of Theorem~\ref{thm:intro-summary-of-Y-U-props} shows that~$\Functor$ 
falls within the scope of Proposition~\ref{prop:gluingfull-faithfulness-Ind-Pro-version-compact-version}, and requires the completion theory
of~\cite{DEGlocalization} in its proof; in fact, we proceed by giving an independent definition of the functor $\Functor_Y$ as a completion of~$\Functor$,
and then proving that~\eqref{eqn:intro recollement diagram} commutes.

The proof of part~\eqref{item:intro FY fully faithful} of Theorem~\ref{thm:intro-summary-of-Y-U-props} immediately reduces to the case in which $Y = \{\thetabar\}$
is a singleton, in which case $\cX_Y$ coincides with~$\cX_{\thetabar}$, and we write~$\Functor_{\thetabar}$ for~$\Functor_Y$. 
As explained in the next paragraph, the full faithfulness of~$\Functor_{\thetabar}$ is a consequence of the results of~\cite{JNWE}.
This %
case of part~\eqref{item:intro FY fully faithful} is then an input for the proof of part~\eqref{item:intro FU fully faithful},
which also requires an explicit quotient presentation of the underlying reduced~$\cU_{\red}$ of the stack~$\cU$, 
which we construct (for this specific choice of~$U$) in Section~\ref{subsec: defining U}.

\subsubsection{Relationship with other works}\label{subsubsec:other works}
For a closed point $\thetabar \in X$,
the functor $\Functor_{\thetabar}$ discussed above turns out to be very closely related to
one of the functors constructed by Johansson--Newton--Wang-Erickson in~\cite{JNWE}. Indeed,
the only (relatively minor) difference
between our functor $\Functor_{\thetabar}$ and the corresponding functor of~\cite{JNWE}
is that while our functor takes values in $D^b_{\coh}(\cX_{\thetabar})$,
where $\cX_{\thetabar}$ is the moduli stack of Galois representations 
introduced above,
theirs takes values in $D^b_{\coh}(\fX_{\thetabar})$, for a certain canonical algebraization
$\fX_{\thetabar}$ of~$\cX_{\thetabar}$;
and our functor is then obtained as a completion of theirs.
Given this fact,
Theorem~\ref{thm:intro-summary-of-Y-U-props}~(4) is essentially a restatement of the full faithfulness results of~\cite{JNWE}.
 
However, establishing this relationship between $\Functor_{\thetabar}$ 
and the corresponding functor of~\cite{JNWE} is non-trivial, since $\Functor_{\thetabar}$ is defined
in terms of Colmez's $D^{\natural}\boxtimes \Pone$ construction,
while the functors of~\cite{JNWE} are defined by certain explicit kernels, using
Morita-theoretic descriptions of the various blocks~$\cA_{\thetabar}$
due to Pa\v{s}k\={u}nas~\cite{MR3150248}.
The relationship between the two 
is given by Proposition~\ref{prop:completed L-infinity},
which in particular provides a conceptual explanation for the particular choice
of kernels made in~\cite{JNWE}. 

More generally, we make full use of previously known results on the $p$-adic local Langlands correspondence for~$\GL_2(\bQ_p)$, including those of~\cite{MR2642409,MR3150248,MR4350140,PaskunasBM,JNWE,KisinFM,HuTan,MR3556448}.
We will often need to restate, or slightly extend, these results to fit our framework, 
and we claim no real novelty in doing so. %

\subsubsection{Sheaf- and category-theoretic techniques}
As this introduction already indicates, we rely on Pro-
and Ind-categorical
techniques to deal with completions,
and to pass from the coherent/finitely presented
setting to the setting of sheaves and modules not satisfying any finiteness conditions.
It seems likely that some, and perhaps all, of the resulting technicalities 
could be ameliorated by using the solid formalism of Clausen and
Scholze;  %
but for various reasons, including a lack of facility in those techniques 
on our part, we have decided to adopt the more traditional Ind/Pro perspective. 

Our arguments also rely heavily on the theory of $\infty$-categories;
in particular, {\em derived categories} are, for us, always understood
in the stable $\infty$-categorical sense.  The most substantial reason
we require these techniques is that we prove our results via gluing
arguments, and the traditional framework of triangulated categories
is inadequate for making such arguments. We note, though,
that we take advantage of this formalism at many other points of
our arguments as well; in particular, we frequently exploit the
flexibility it allows with regard to constructing
and canonically characterizing derived functors.  

While these techniques (Ind/Pro-categories and $\infty$-categories)
are standard tools for experts in the categorical and geometric
aspects of the Langlands program, they may be less familiar to
readers who (like ourselves) are approaching the paper from 
the perspective of the traditional theory 
of the $p$-adic local Langlands correspondence.
To this end, and in order to have a uniform phrasing of several results from the literature, 
we have included two appendices --- Appendix~\ref{app:category-theory}
and Appendix~\ref{sheaves on formal algebraic stacks} --- which
provide background on category-theory (especially
$\infty$-category theory) and on coherent sheaf-theory on formal algebraic stacks.

\subsubsection{A guide to the paper}%
Section~\ref{sec:recollections} collects various background results from Galois representation theory,
smooth representation theory of $\GL_2(\bQ_p)$ in characteristic~$p$, and the existing theory of the local
$p$-adic Langlands correspondence for this group.
We provide several complementary results, usually of a categorical nature, to use later in the paper.
We then describe the localization theory of~\cite{DEGlocalization} and generalize it to derived categories, and,
using results of
Heyer~\cite{Heyerderived, Heyergeometric},
we prove some finiteness results
for $\Ext$-groups between certain non-admissible representations.

Section~\ref{sec: moduli stacks} describes the Galois-theoretic
side of our correspondence, building on the results in Appendix~\ref{sec:stacks-appendix},
and on several previous works in the literature, including especially~\cite{emertongeepicture, MR3831282, JNWE}.
We then specialize the discussion of coherent, Ind-coherent, and Pro-coherent sheaves 
of Appendix~\ref{sheaves on formal algebraic stacks} to the present context of moduli stacks of $(\varphi, \Gamma)$-modules,
and carry out several explicit computations.
We conclude by building a connection between these objects and the localization theory of~\cite{DEGlocalization}.

Section~\ref{sec: Colmez} contains our generalization of a significant part of Colmez's work~\cite{MR2642409}, 
to coefficients in arbitrary Noetherian $\cO/\varpi^a$-algebras.
We begin by proving some technical statements in Drinfeld's theory of Tate modules, which is the natural context for this generalization.
We then work towards one of our main results, which is Theorem~\ref{thm: D natural is G-stable} on the $\GL_2(\Q_p)$-stability of $D^\natural \boxtimes \Pone$,
and develop, along the way, analogues of Colmez's functors $D^+, D^{++}, D^{\nr}, D^\sharp, D^\natural$, and $D \boxtimes -$.
We make significant use of novel base change results, which will also be important in Section~\ref{sec: the functor} to construct~$L_\infty$;
these results can be specialized to the classical context of coefficients in a complete Noetherian local $\cO$-algebra, and yield new results even in that context
(see Theorem~\ref{special fibre of versal deformation}).

Section~\ref{sec: the functor} contains our construction of the functor~$\Functor$, 
and the proof of Theorems~\ref{thm:coherent sheaves intro} and~\ref{thm:intro-summary-of-Y-U-props}.
We begin by constructing a version of $D^\natural \boxtimes \Pone$ for the universal object on the stack~$\cX$, which (as already explained) is done by descent,
as an application of the results of Section~\ref{sec: Colmez};
the pro-coherent sheaf~$L_\infty$ is then obtained as an appropriate twist of this universal
version of~$D^{\natural}\boxtimes \Pone$.
We next make a detailed study of~$L_\infty$, starting from its (completed) pullbacks to the stacks~$\cX_{\thetabar}$ of Galois representations, and of
the functor $\Functor \coloneqq  L_\infty \otimes_{\cO\llbracket G\rrbracket _\zeta} -$.
Sections~\ref{finiteness for Functor} and ~\ref{subsection:localization-to-Ugood} establish the properties of~$\Functor$ that are required to apply
our full faithfulness criterion (Proposition~\ref{prop:gluingfull-faithfulness-Ind-Pro-version-compact-version}), which we do in Section~\ref{subsec:proving-that-F-is-fully-faithful}, where we conclude the proof of our main result.

Appendix~\ref{app:category-theory} assembles some background from $1$-category and $\infty$-category theory, with a significant focus on $\Ind$- and $\Pro$-completions
of categories, which are the context of our main full faithfulness criterion in Section~\ref{subsec:semiorthogonal and Ind Pro}. 
One complication arises from the interaction between $\Pro$-completions and the formation of left derived functors,
which we study systematically in Sections~\ref{subsec:left derived} and~\ref{subsubsec:Pro co-Yoneda}.
As already mentioned, much of this material is standard; but we have found it helpful to collect it here, 
to give a uniform phrasing to several results from the literature,
and to have a uniform foundation for the less standard material in this appendix.

Appendix~\ref{sheaves on formal algebraic stacks} develops the theory of coherent sheaves on formal algebraic stacks, both at the abelian and derived level,
including various functorial operations on these objects, and recasting some of the
coherent completeness results of~\cite{alper2023etale} and~\cite{alper2023coherently}
in our framework. Finally, Appendix~\ref{sec:stacks-appendix} provides some complements to the results of~\cite{emertongeepicture} and~\cite{MR3831282}.

\subsubsection{A comparison with the proof outlined in~{\em \texorpdfstring{\cite{emerton2023introduction}}{the IHES notes}}.}

The results of this paper were announced in~\cite[Section~7]{emerton2023introduction}, which also gave a sketch of the proofs,
emphasizing
the analogy between the results of~\cite[Section~3.8]{DEGlocalization}, and Beauville--Laszlo gluing of sheaves on~$\cX$.
This was presented as an important step in the proof.
Furthermore, we indicated a proof of part~\eqref{item:intro FU fully faithful} of Theorem~\ref{thm:intro-summary-of-Y-U-props}
via an explicit Morita-theoretic argument.  

However, in the course of writing this paper,
we found a more efficient way of formulating these gluing results, avoiding explicit
Beauville--Laszlo-type gluing, and proceeding instead via the
notion of a recollement of $\infty$-categories.
This resulted in the current proof of the full faithfulness of~$\Functor$,
as well as a simplification of the proof of the Beauville--Laszlo-type results of~\cite{DEGlocalization}
(which has been incorporated into that paper).
Furthermore, we found a more efficient proof of 
part~\eqref{item:intro FU fully faithful} of Theorem~\ref{thm:intro-summary-of-Y-U-props},
which avoids the explicit constructions and computations required in the Morita-theoretic approach,
and instead uses the same results related to completions that are used to prove
the full faithfulness of the functors $\Functor_Y$ for finite closed subsets $Y$ of~$X$.

In conclusion,
although the overall structure of the argument sketched in~\cite{emerton2023introduction} has been preserved,
and none of the statements of \emph{loc.~cit.}
need to be corrected, several of the details of our arguments are different from those sketched there.

\subsection{Acknowledgements}\label{subsec: acknowledgements}
We would like to thank Bhargav Bhatt, George Boxer, Ana Caraiani, Pierre Colmez,
Laurent Fargues, Michael Harris, Eugen Hellmann, Christian Johansson, Bao Le Hung, Brandon Levin, Jacob Lurie,
Akhil Mathew,
James Newton, Vytautas Pa\v{s}k\={u}nas, 
David Savitt, Peter Scholze, and Carl Wang-Erickson
for helpful correspondence and conversations.
We are particularly grateful to Christian Johansson, James Newton,
and Carl Wang-Erickson for sharing with us preliminary versions
of the  results of their paper~\cite{JNWE}.
We would also like to thank Vytautas Pa\v{s}k\={u}nas for comments on an
 earlier version of this paper.
\subsection{Notation and conventions}\label{subsec: notation and
  conventions}
We fix throughout the paper a prime~$p\ge 5$. %
Fix an algebraic closure~$\Qpbar$ of~$\Qp$. 
If~$K \subset \Qpbar$ is a finite
extension of~$\bQ_p$, we write~$G_K$ for the absolute Galois
group~$\Gal(\Qpbar/K)$. 
Write~$I_K$ for the inertia subgroup
of~$G_K$. 
We normalise local class field theory so that a uniformizer
corresponds to a geometric Frobenius element. 
Our convention for Hodge--Tate
weights is that the $p$-adic cyclotomic
character~$\varepsilon:G_{\Qp}\to\Q_p^\times$ has Hodge--Tate
weight~$-1$.

Let $\cO$ denote the ring of integers in a fixed finite extension
$E$ of $\Q_p$, let~$\varpi$ be a fixed uniformizer of~$\cO$, and let~$\F=\cO/\varpi$ be the residue field of~$\cO$.
Fix an algebraic closure $\cbF_p$ of~$\bF$.
Unless otherwise stated, all
representations considered in this paper will be on
$\cO$-modules. 
We fix a continuous character
$\zeta:\Qptimes\to\cO^\times$, which we also regard as a character
of~$G_{\Qp}$ via local class field theory. 
Throughout the paper we
will work with representations of~$\GL_2(\Qp)$ of central
character~$\zeta$, and representations of~$G_{\Qp}$ with determinant~$\zeta\varepsilon^{-1}$.
We say that we are in the {\em even} (resp.\ {\em odd}) case if $\zeta$ is
an even (resp.\ odd) character. %

We write~$G=\GL_2(\Q_p)$, $K=\GL_2(\Zp)$, and~$Z=\Qp^\times$ for the centre of~$G$. 
The diagonal torus in~$G$ is denoted by~$T$, and we write~$B$, resp.~$\lbar B$ for the upper-triangular, resp. lower-triangular Borel subgroup containing~$T$, with unipotent radicals~$U$ and~$\lbar U$.
The maximal compact subgroup of~$T$ is denoted by~$T_0$.
The $n$th congruence subgroup of~$K$ is denoted~$K_n \coloneqq  1 + p^nM_2(\bZ_p)$, for $n \geq 1$.
The upper-triangular Iwahori subgroup of~$K$ is denoted~$\Iw$, and we write~$\Iw_1$ for its pro-$p$ Sylow subgroup and $Z_1 = \Iw_1 \cap Z = K_1 \cap Z$ for the maximal pro-$p$ subgroup of~$Z$.
Since~$p \geq 3$, there are isomorphisms
\begin{equation}\label{decomposition of pro-p subgroups}
\Iw_1 \isom Z_1 \times (\Iw_1/Z_1) \text{ and } K_1 \isom Z_1 \times (K_1/Z_1),
\end{equation}
given by the determinant and the projection to the quotient.
Since~$p \geq 5$, the groups~$\Iw_1$ and~$K_1$ are torsion-free, and so the same is true for their direct factors~$\Iw_1/Z_1$ and~$K_1/Z_1$.

If~$Z \subset H \subset G$ is a subgroup, we will write $\cO[H]_\zeta$ for the quotient 
of the group algebra~$\cO[H]$ by the two-sided ideal generated by $\{z-\zeta(z): z \in Z\}$.

\subsubsection*{Stacks}Our conventions on algebraic stacks and formal
algebraic stacks are those of~\cite{stacks-project}
and~\cite{Emertonformalstacks}. The reader who is unfamiliar with this
material may wish to refer to the overview
of~\cite{Emertonformalstacks} in~\cite[App.\
A]{emertongeepicture}. If~$A$ is a topological ring of the kind for which $\Spf A$
is defined (see~\cite[\href{https://stacks.math.columbia.edu/tag/0AIC}{Tag
   0AIC}]{stacks-project} for an extensive discussion of this; we won't
need any examples beyond $I$-adically complete Noetherian rings) 
and~$\cC$ is a stack then we write~$\cC(A)$ for~$\cC(\Spf A)$; if~$A$ has the discrete
topology, then this is equal to~$\cC(\Spec A)$.

\subsubsection*{Module categories}
We will consider several categories of modules over noncommutative rings~$R$.
Unless stated otherwise, we will work with left modules.
We write $\Mod(R)$ for the category of left $R$-modules and $R$-linear maps, $\Mod^{\fl}(R)$ for the full subcategory of modules of finite length, and $\Mod^{\fp}(R)$ for the
full subcategory of finitely presented modules, which coincides with the full subcategory of compact objects of $\Mod(R)$.
We will introduce several variants in the course of the text, such as the category of profinite modules over a profinite topological ring, see Section~\ref{compact modules}.
Recall that the Jacobson radical of~$R$, denoted~$\Rad(R)$, is the intersection of all maximal left ideals of~$R$.
It coincides with the intersection of all maximal right ideals, respectively two-sided ideals, of~$R$.
Unless stated otherwise, we write ``Noetherian'' as a shorthand for ``left and right Noetherian''.

\subsubsection*{Category theory and homological algebra}
We employ the language of $\infty$-categories throughout the paper,
and stable $\infty$-categories are the basic setting in which we
undertake our homological algebra.  We also liberally employ the language
of Ind and Pro completion.  
Appendix~\ref{app:category-theory} %
summarizes much of the material we use.
Of course, most of what we use can be found (often in more general form)
in~\cite{MR2522659}, \cite{LurieHA}, and~\cite{LurieSAG}.
On several occasions we cite~\cite{MR2182076}, which %
provides an extensive survey of many topics in homological algebra and category theory
(in the traditional language of derived categories).  One caution we make is
that its theory of Ind and Pro categories is developed in the context of large categories,
rather than in context of small categories that we work with here. 

We will be applying several of the constructions developed in~\cite{DEGlocalization}.
That paper has a more abelian categorical flavour, 
and we wish to translate its results into the derived context of this paper.
The formation of left derived functors often involves passing to Pro-categories,
but for the applications we have in mind,
we do not wish to work in the derived category of the Pro-category, but
rather in the Pro category of the derived category.
An explanation of how to arrange this is given in \S~\ref{subsec:left derived}. 

\section{Background}\label{sec:recollections}
In this section we gather together a range of results related
to Galois representations as well as to smooth representations
which we will require later on in the paper.

\subsection{Galois representations}
We introduce various notation and concepts related to two-dimensional representations
and pseudorepresentations of~$G_{\Q_p}$.

\subsubsection{Mod~\texorpdfstring{$p$}{p} Galois representations, Serre weights, and the weight part of Serre's
  conjecture}\label{subsubsec: Galois background}%
We write~$\omega=\varepsilonbar:G_{\Q_p} \to \F_p^{\times}$
for the mod~$p$ cyclotomic character.
We sometimes regard~$\omega$ as a character of~$\Q_p^{\times}$
via local class field theory 
(which, as noted above, is
normalised so as to take a uniformizer to a  geometric
Frobenius element). 
Concretely, then, $\omega$ is the reduction mod $p$ of the character $x|x|$
(i.e.\ trivial on $p$, and the reduction map to $\F_p^{\times}$ on $\Z_p^{\times}$).
 
For any~$x\in\cbF_p^\times$, we let
\[
\nr_x:G_{\Qp}\to\cbF_p^\times
\]
be the
unramified character taking a geometric Frobenius element~$\Frob_{p}$ to~$x$. 
Thus every character $G_{\Qp}\to\cbF_p^\times$ is of the
form~$\nr_x\omega^i$ for some uniquely determined~$x\in\cbF_p^\times$
and $0\le i<p-1$.

\begin{defn}
  \label{defn: generic characters}We say that an unordered pair of characters
  $\chi_1,\chi_2:G_{\Qp}\to\Fbar_p^\times$ is \emph{generic} if
  $\chi_1\chi_2^{-1}\ne 1,\omega^{\pm 1}$. Note that by Tate
  local duality, this implies that
  $H^{0}(G_{\Qp},\chi_1\chi_2^{-1})$ and
  $H^{2}(G_{\Qp},\chi_1\chi_2^{-1})$ both vanish, 
  and Tate's local Euler characteristic formula then shows that 
  $H^{1}(G_{\Qp},\chi_1\chi_2^{-1})$ is a one-dimensional
  $\Fbar_p$-vector space.
\end{defn}

We work throughout with the ring of integers $\cO$ in
the finite extension $E$ of~$\Q_p$ as our base ring, and so its residue field
$\F$ plays a special notational role. We often consider representations
defined over~$\F$, or $\Gal(\Fbar_p/\F)$-orbits of representations defined over~$\Fbar_p$. 
However, all the results we discuss are insensitive to base-change from
$\cO$ to $\cO'$ (an extension of rings of integers arising from an inclusion
$E \subseteq E'$ of finite extensions of~$\Q_p$), and so we are always free
to enlarge $\F$ as necessary.  Thus, although we have (of course) tried to
be correct 
in our treatment of questions of fields of definition, we encourage the reader
to not pay too much attention to this purely technical issue.

This being said, 
let $\F'$ now denote the unique quadratic extension of~$\F$, and
write~$\omega_2:I_{\Qp}\to(\F')^\times$ for a choice
of fundamental character of niveau~$2$; writing~$\Q_{p^2}$ for the
quadratic unramified extension of~$\Qp$, we can extend~$\omega_2$ to a
character $G_{\Q_{p^2}}\to(\F')^\times$ in such a way
that~$\omega_2^{p+1}=\omega|_{G_{\Q_{p^2}}}$. Then the 2-dimensional absolutely irreducible
representations $G_{\Qp}\to\GL_2(\F)$ are precisely those of the
form~$\nr_x\otimes\Ind_{G_{\Q_{p^2}}}^{G_{\Qp}}\omega_2^k$ for 
some ~$k\in\Z/(p^2-1)\Z$ with $(p+1)\nmid k$ and some~$x\in(\F')^\times$ with~$x^2\in \F^{\times}$.\footnote{Although
this expression involves objects only defined over~$\F'$, 
the representation itself is defined over~$\F$.  Since, as already noted, we allow
ourselves to enlarge $\F$ as convenient,
questions of precise fields of definition are of little consequence to us in any case.} 
In this case~$x$ is
uniquely determined up to multiplication by~$-1$, and~$k$ is uniquely
determined up to multiplication by~$p$. %

\begin{defn}\label{defn: Serre weight}
A {\em Serre weight} is an irreducible representation
  $\sigma_{a,b} \coloneqq \det{}^{a}\otimes\Sym^b\Fp^2$ of~$\GL_2(\Fp)$, where $0\le
  a<p-1$ and~$0\le b\le p-1$. It is sometimes convenient to view~$a$
  as an element of~$\Z/(p-1)\Z$,
  and we will do so without further
  comment. We say that~$\sigma_{a,b}$ is \emph{Steinberg} if~$b=p-1$, and
  \emph{non-Steinberg} otherwise. %
\end{defn}

We %
now explain what it means to say that ``$\sigma$
is a Serre weight for~$\rhobar$''. %
We do not attempt to motivate this description, which goes
back to~\cite{SerreDuke}, but we remark that one possible motivation
is given by Theorem~\ref{thm: fixed determinant stack for GL2
  Qp}~\eqref{item:points and Serre weights}. %

\begin{defn}\label{defn:Serre weights of rhobar}\leavevmode
\begin{enumerate}
\item If~$\rhobar=\nr_x\otimes\Ind_{G_{\Q_{p^2}}}^{G_{\Qp}}\omega_2^k$ is
irreducible, then~$\sigma_{a,b}$ is a Serre weight for~$\rhobar$ if
and only if $k\equiv (p+1)a+b-p\pmod{p^2-1}$ or $k\equiv
(p+1)a+pb-1\pmod{p^2-1}$. %
\item If~$b\ne 0$, then~$\sigma_{a,b}$ is a Serre weight for a
reducible representation~$\rhobar$ if and only if
\begin{equation}\label{eqn: form of reducible repn Serre weight
  fixed}\rhobar\cong
\begin{pmatrix}
  \nr_x \omega^{a+b} &* \\0&\nr_y\omega^{a-1}
\end{pmatrix}\end{equation}for some~$x,y\in\F^\times$. %
\item
If~$b=0$, then $\sigma_{a,0}$ is a Serre weight for a reducible
representation $\rhobar$ if and only if $\rhobar$ can be written as an extension
as in~\eqref{eqn: form of reducible repn Serre weight fixed}, and if $\rhobar$ is furthermore 
finite flat. If~$x\ne y$ then this further condition actually holds automatically,
while if $x=y$ it is equivalent to requiring %
that the
extension be peu ramifi\'ee (which by definition amounts to asking that %
the corresponding Kummer class is given by an integral unit). %
\end{enumerate}
\end{defn}

We now note a compatibility between central characters and
determinants. %

\begin{defn}\label{defn: compatible Serre weight}%
A Serre weight~$\sigma_{a,b}$ is \emph{compatible with~$\zeta$,} or \emph{$\zeta$-compatible}, if
$\omega^{2a+b} =
\zetabar|_{I_{\Qp}}.$
\end{defn}%
It is immediate from the definitions above that
if~$\det\rhobar=\zetabar\omega^{-1}$, and  ~$\sigma$ is a Serre
weight for~$\rhobar$, then~$\sigma$ is compatible
with~$\zeta$. Accordingly, we assume throughout the paper that all
Serre weights~$\sigma$ that we consider are compatible with~$\zeta$. We will
also sometimes think of a $\zeta$-compatible Serre weight as being an $\cO[KZ]_\zeta$-module, letting
~$K$ act through its quotient $K/K_1 = \GL_2(\F_p)$ and letting~$Z$ act via~$\zeta$.

\begin{remark}
\label{rem:parity}
Via local class field theory, we can view the condition of
Definition~\ref{defn: compatible Serre weight} %
as saying
that~$\zetabar$ (restricted to $\F_p^{\times}$) agrees with the
central character of~$\sigma$. %
In particular, $b$ is even (resp.\ odd) if and only if~$\zeta$ is even (resp.\ odd); furthermore, for each choice of $b$ of the
appropriate parity, there are two corresponding choices of~$a$ for
which $\sigma_{a,b}$ is a compatible Serre weight.
\end{remark}

The following definition will be used throughout the paper.
Again, we do not attempt to motivate this here,
but we remark that it goes back to~\cite{GrossTameness}, and we point to Corollary~\ref{lem:reduced reducible connected components} as one possible motivation.

\begin{defn}\label{defn: companion Serre weights}\leavevmode
  \begin{enumerate}
  \item The \emph{companion} of a Serre weight ~$\sigma=\sigma_{a,b}$
  with~$0\le b\le p-3$ is the Serre weight
  $\sigmacomp\coloneqq \sigma_{a+b+1,p-3-b}$. A Serre weight of the
  form~$\sigma_{a,p-2}$ is defined to be its own companion. (We do not
  define a companion for a Steinberg Serre weight.) %
  We refer to the unordered pair of~$\sigma,\sigmacomp$ as a
  \emph{companion pair}. We say that the unordered pair~$\sigma,\sigmacomp$ is of
  type~\ref{item: non p-distinguished pseudorep} if it is of the form
  $\sigma_{a,p-2}$, of type~\ref{item: Steinberg pseudorep} if it is of the form $\sigma_{a,0},\sigma_{a+1,p-3}$, and otherwise we say that it
  is of type~\ref{item: generic pseudorep}. 
\item Let~$\{\sigma, \sigmacomp\}$ be a companion pair of Serre weights, and write $\sigma= \sigma_{a, b}$.
We write~$\Theta(\sigmasigmacomp)$ for the cuspidal $E[\GL_2(\bF_p)]$-representation denoted~$\Theta(\chi)$ 
in~\cite[\S3]{CDT}, where~$\chi : \bF_{p^2}^\times \to \cO^\times$ lifts $x \mapsto x^{(b+2)+(a-1)(p+1)}$.
Note that by~\cite[Lem.\ 3.1.1(4)]{CDT}, the
  Jordan--H\"older factors of $\Theta(\chi)\otimes_{\cO}\F$
  are~$\{\sigma,\sigmacomp\}$.
  If~$\sigma$ is $\zeta$-compatible, we will often regard~$\Theta(\sigmasigmacomp)$ as an $E[KZ]_\zeta$-module, 
  by letting~$Z$ act via~$\zeta$, and~$K$ via its quotient $K/K_1 =\GL_2(\bF_p)$.
\end{enumerate}
\end{defn}

\subsubsection{Pseudorepresentations}\label{subsec:pseudoreps-recollections}
We will make use of the notion of a %
pseudorepresentation of~$G_{\Qp}$ on an~$\cO$-algebra~$A$. 
We will assume throughout, without further comment, that all
our pseudorepresentations are $2$-dimensional and have
determinant~$\zeta\varepsilon^{-1}$. %
Since~$p>2$, a $2$-dimensional pseudorepresentation is equivalent to the data of a $2$-dimensional determinant in the sense
of~\cite[Defn.\ 1.5]{MR3444227}, so such a pseudorepresentation is
equivalent (see~\cite[Lem.\ 1.9, Prop.\ 1.29]{MR3444227}) to a
continuous function $\theta:G_{\Qp}\to A$ satisfying
\begin{itemize}
\item $\theta(1)=2$.
\item $\theta(gh)=\theta(hg)$.
\item $(\zeta\varepsilon^{-1})(g)\theta(g^{-1}h)-\theta(g)\theta(h)+\theta(gh)=0$.
\end{itemize}

Given a continuous representation $\rho:G_{\Qp}\to\GL_2(A)$, the function
$\tr\rho:G_{\Qp}\to A$ is a pseudorepresentation in this sense. %
If~$A$ is a finite field, or an algebraically closed field,
then any pseudorepresentation comes from a representation, by~\cite[Corollary~2.9]{MR3831282}. 
So we can
and do identify
$\bF$-valued pseudorepresentations with semisimple
representations with coefficients in~$\bF$. 

\begin{defn}\label{defn:pseudorep-types}
We classify the %
$2$-dimensional $\cbF_p$-valued pseudorepresentations $\thetabar: G_{\Qp} \to \cbF_p$
into five ``types'' as follows; %
we explain the labels in Remark~\ref{rem: labelling of blocks}.
  \begin{enumerate}%
  \item[{\crtcrossreflabel{(ssg)}[item: ss pseudorep]}]
    $\thetabar: G_{\Q_p} \to \cbF_p$ is the trace of an absolutely irreducible two-dimensional $\bF$-representation.
  \item[{\crtcrossreflabel{(gen)}[item: generic pseudorep]}]
    $\thetabar = \chi + \zeta \omega^{-1}\chi^{-1}$ for some character $\chi : G_{\bQ_p} \to \bF^\times$ such that the pair $\{\chi, \zeta \omega^{-1}\chi^{-1}\}$
    is generic in the sense of Definition~\ref{defn: generic characters}.
  \item[{\crtcrossreflabel{(scalar)}[item: non p-distinguished pseudorep]}]
    $\thetabar = \chi + \chi$ for some character
    $\chi: G_{\Q_p} \to \F^{\times}$ such that $\chi^2 = \zeta\omega^{-1}$. These only exist if $\zeta$ is odd and $\zeta(p)$ is a square.
  \item[{\crtcrossreflabel{(St)}[item: Steinberg pseudorep]}]
    $\thetabar = \chi + \chi\omega^{-1}$ for some character
    $\chi: G_{\Q_p} \to \F^{\times}$ such that $\chi^2 = \zeta$. These only exist if~$\zeta$ is even and $\zeta(p)$ is a square.
  \item[{\crtcrossreflabel{(gen+)}[item: not abs irred pseudorep]}] $\thetabar = \chi + \zeta\omega^{-1}\chi^{-1}$ %
  for some character $\chi: G_{\bQ_p} \to \cbF_p^\times$ not defined over~$\bF$. 
  \end{enumerate}
\end{defn}

\begin{rem}\label{explain terminology sigmasigmacomp}
The reason for the
  choice of terminology in Definition~\ref{defn: companion Serre weights} is that the semisimple Galois representations of type~\ref{item:
    non p-distinguished pseudorep} (resp.\ of type~\ref{item: Steinberg
    pseudorep}) have Serre weights of the form $\sigma_{a,p-2}$
  (respectively of the form $\sigma_{a,0},\sigma_{a+1,p-3}$, $\sigma_{a,p-1}$).
\end{rem}

Given any $\lbar \bF_p$-valued pseudorepresentation~$\thetabar$ we %
write  $\bF_{\thetabar}$ for the (finite) subfield of~$\Fpbar$ generated by~$\bF$ and the values of~$\thetabar$.
It only depends on the $\Gal(\cbF_p/\bF)$-conjugacy class of~$\thetabar$.
 We write  $(R_{\thetabar}^{\ps}, \theta^\univ)$ for the universal deformation of~$\thetabar$ to complete Noetherian local $\cO \otimes_{W(\bF)} W(\bF_{\thetabar})$-algebras,
 where by a \emph{deformation}
of~$\thetabar$ to a complete local Noetherian $\cO$-algebra~$A$ with residue field~$\F_{\thetabar}$, 
we mean a $2$-dimensional pseudorepresentation $\theta : G_{\bQ_p} \to A$, with determinant~$\zeta\varepsilon^{-1}$, such that the composition $G_{\bQ_p} \xrightarrow{\theta} A \to \F_{\thetabar}$ 
equals~$\thetabar$; the universal deformation exists by~\cite[Prop.\ 3.3]{MR3444227}.
We sometimes refer to~$R_{\thetabar}^{\ps}$ as the (universal) pseudodeformation ring for~$\thetabar$.

\begin{rem}\label{F-rational}
We will sometimes say that a pseudorepresentation~$\lbar \theta$ is \emph{$\bF$-rational} if~$\bF_{\thetabar} = \bF$.
The first four items of Definition~\ref{defn:pseudorep-types} are all $\bF$-rational, and classify pseudorepresentations attached to semisimple representations $G_{\bQ_p} \to \GL_2(\bF)$ 
all of whose irreducible summands are absolutely irreducible. 
The $\bF$-rational pseudorepresentations of type~\ref{item: not abs irred pseudorep} correspond to irreducible representations $\rhobar : G_{\bQ_p} \to \GL_2(\bF)$ that are not absolutely irreducible.
These are precisely the twists of unramified representations such that the characteristic polynomial of Frobenius is irreducible quadratic.
\end{rem}

\begin{defn}\label{defn:CH-alg-Rtilde}
  We write $(\tld R_{\thetabar}, \theta^\univ)$ for the Cayley--Hamilton $R^{\ps}$-algebra associated to~$\theta^\univ$, with anti-involution~$\dagger$. Explicitly, $\tld R_{\thetabar}$ is the quotient $R_{\thetabar}^{\ps}\llbracket G_{\bQ_p}\rrbracket /J$, where $R_{\thetabar}^{\ps}\llbracket G_{\bQ_p}\rrbracket $ denotes the completed group ring, and~$J$ is the closure of the 2-sided ideal in~$R_{\thetabar}^{\ps}\llbracket G_{\bQ_p}\rrbracket $ generated by the elements \[g^2-\theta^{\univ}(g)g+(\zeta\varepsilon^{-1})(g)\]for
  all~$g\in G_{\bQ_p}$.
\end{defn}
Note that the identity
\begin{equation}\label{eqn:Rps in tR}
\theta^{\univ}(g) = g + (\zeta\varepsilon^{-1})(g) g^{-1} \in \widetilde{R}_{\thetabar}
\end{equation}
for $g \in G_{\Q_p}$
shows that the natural map
\begin{equation}\label{CHquotient}
\cO\llbracket G_{\Qp}\rrbracket  \to 
\widetilde{R}_{\thetabar}
\end{equation}
is surjective,
and uniquely determines the embedding
$R_{\thetabar}^{\ps} \hookrightarrow
\widetilde{R}_{\thetabar}$.
The anti-involution $\dagger: g \mapsto
(\zeta\varepsilon^{-1})(g)g^{-1}$ preserves the ideal~$J$ and descends to
an $R_{\thetabar}^{\ps}$-linear anti-involution of~$\tld R_{\thetabar}$. 

Note that for any two pseudorepresentations in the $\Gal(\cbF_p/\bF)$-conjugacy class of~$\thetabar$, the corresponding universal pseudodeformation rings and Cayley--Hamilton algebras 
are canonically identified (using the Galois action on coefficients).

We will need the following structural properties of~$\tld R_{\thetabar}$.

\begin{prop}\label{CHtorsionfree}
Let~$\thetabar$ be a $2$-dimensional $\lbar \bF_p$-valued pseudorepresentation.
The ring $\tld R_{\thetabar}$ is Noetherian, of finite type as an $R_{\thetabar}^{\ps}$-module, and $p$-torsion free.
The natural map $R_{\thetabar}^{\ps} \to Z(\tld R_{\thetabar})$ is an isomorphism.
\end{prop}
\begin{proof}
All statements of the proposition can be checked after making a finite unramified extension $\cO \to \cO'$ (using that $Z(\tld R_{\thetabar}) \otimes_{\cO} \cO' \to Z(\tld R_{\thetabar} \otimes_{\cO} \cO')$
is an isomorphism).
We can therefore assume without loss of generality that~$\thetabar$ does not have type~\ref{item: not abs irred pseudorep}.

Since~$R_{\thetabar}^{\ps}$ is Noetherian (by~\cite[Proposition~3.2]{MR3831282}), the first statement is a consequence of the second.
The second statement is~\cite[Proposition~3.6~(ii)]{MR3831282}.
There remains to prove the claim that~$\tld R_{\thetabar}$ is $p$-torsion free, and the map $R_{\thetabar}^{\ps} \to Z(\tld R_{\thetabar})$ is an isomorphism.

When~$\thetabar$ is residually multiplicity-free (i.e.\ always except in case~\ref{item: non p-distinguished pseudorep}), 
the claim follows from the fact that~$\tld R_{\thetabar}$ is a 
generalized matrix algebra over~$R_{\thetabar}^{\ps}$, by \cite[Proposition~3.6~(5)]{MR3831282},
since~$R^{\ps}_{\thetabar}$ is a $p$-torsion free Noetherian integral domain. 
In more detail, $R^{\ps}_{\thetabar}$ is formally smooth of dimension~3 over~$\cO \otimes_{W(\bF)} W(\bF_{\thetabar})$ in case~\ref{item: ss pseudorep}, by e.g.~\cite[Proposition~6.2, Remark~6.6]{MR3150248}, 
and in case~\ref{item: generic pseudorep}, %
by e.g.~\cite[Proposition~B.17, Remark~B.28]{MR3150248}.
The structure of $R^{\ps}_{\thetabar}$ in case~\ref{item: Steinberg pseudorep} is described
in~\cite[Remark~10.88, Corollary~B.5]{MR3150248}, using
\cite[Remark~10.88]{MR3150248} to identify~$R^{\ps}_{\thetabar}$ with
the ring denoted~$R^\psi$ in \emph{loc.\ cit.}
This makes it clear that $R_{\thetabar}^{\ps}$ is $p$-torsion free in this case as well. 

Finally, when~$\thetabar$ is of type~\ref{item: non p-distinguished pseudorep}, our claim is a consequence of~\cite[Corollary~9.24, Corollary~9.25]{MR3150248}, using \cite[Corollary~9.33]{MR3150248} to identify~$\tld R_{\thetabar}$ with the ring denoted~$R$ in~\emph{loc.\ cit.}
\end{proof}

\begin{lemma}\label{normalization for CH algebra}
Let~$\thetabar$ be a $2$-dimensional $\cbF_p$-valued pseudorepresentation, and assume that~$\thetabar$ does not have type~{\em \ref{item: non p-distinguished pseudorep}}
or~{\em \ref{item: not abs irred pseudorep}}.
Then every $R^{\ps}_{\thetabar}$-algebra automorphism of~$\tR_{\thetabar}$ that fixes all isomorphism classes of simple $\tR_{\thetabar}$-modules is inner.
\end{lemma}
\begin{proof}

The assumptions on~$\thetabar$ mean that its irreducible summands are distinct and absolutely irreducible. 
As recalled during the proof of Proposition~\ref{CHtorsionfree}, $R^{\ps}_{\thetabar}$ is an integral domain.
Writing~$K$ for its field of fractions, it then follows from
\cite[Thm. 1.4.4(ii)]{MR2656025}
that~$\tR_{\thetabar}$ is a generalized matrix $R^{\ps}_{\thetabar}$-subalgebra of~$M_2(K)$,
in such a way that the pseudorepresentation $\theta^{\univ}: \tR_{\thetabar} \to R^{\ps}_{\thetabar}$ coincides %
with the restriction of the trace function.

Each irreducible summand~$\thetabar_i$ of~$\thetabar$ is the pseudorepresentation associated to a unique simple $\tR_{\thetabar}$-module~$M_i$. 
Part of the generalized matrix algebra structure on~$\tR_{\thetabar}$ is
an idempotent~$e_{ii} \in \tR_{\thetabar}$ such that, for all~$x \in \tR_{\thetabar}$, we have 
\[
\lbar \theta(e_{ii}x e_{ii}) = \thetabar_i(x).
\]
Now we claim that the $\varphi(e_{ii})$ satisfy condition~(3) in \cite[Lemma~1.4.3]{MR2656025}, i.e.\ 
\[
\lbar \theta(\varphi(e_{ii})x \varphi(e_{ii})) = \thetabar_i(x) 
\]
for all~$x \in \tR_{\thetabar}$.
The assumption that~$\varphi$ fixes every isomorphism class of simple $\tR_{\thetabar}$-modules implies that $\lbar \theta_i \circ \varphi = \lbar \theta_i$ for all~$i$
(because if~$\thetabar_i$ is the trace of the simple $\tR_{\thetabar}$-module~$M_i$, then $\thetabar_i \circ \varphi$ is the trace of $\varphi^*(M_i) \cong M_i$).
In turn, this implies that $\lbar \theta \circ \varphi = (\sum_i \lbar \theta_i) \circ \varphi = \sum_i (\lbar \theta_i \circ \varphi) = \sum_i \lbar \theta_i = \lbar \theta$, 
and we can then compute that
\[
\lbar \theta(\varphi(e_{ii})x \varphi(e_{ii}))  = (\lbar \theta \circ \varphi)(e_{ii}\varphi^{-1}(x)e_{ii})
= \lbar\theta(e_{ii}\varphi^{-1}(x)e_{ii}) = \lbar \theta_i(\varphi^{-1}(x)) = \lbar \theta_i(x).
\]
Now the final statement of \cite[Lemma~1.4.3]{MR2656025} shows that~$\varphi$ 
is $\tR_{\thetabar}^\times$-conjugate to an automorphism that fixes all the~$e_{ii}$. %
It thus suffices to prove that every $R$-automorphism~$\varphi$ of $\tld R_{\thetabar}$ which fixes all the~$e_{ii}$ is inner. %
This follows by inspection from the presentation of~$\tR_{\thetabar}$ described in Section~\ref{subsubsec:ss CWE stack}---\ref{subsubsec:Steinberg CWE stack} below;
we give details in the case~\ref{item: Steinberg pseudorep}, which is the only case that we will need in the rest of the paper.
Here $R\coloneqq  R^{\ps}_{\thetabar} = \cO\llbracket a_0, a_1, X_0, X_1\rrbracket /(a_0X_1+a_1X_0)$, and
\[
\tld R_{\thetabar} \cong \fourmatrix{R}{RX_0+RX_1}{R}{R} \subset M_2(R), 
\]
compare Remark~\ref{rem:2x2 matrix order}.
The automorphism~$\varphi$ fixes~$e_{11}, e_{22}$, hence it preserves the direct sum decomposition
of~$\tR_{\thetabar}$ according to matrix entries,
and it restricts to the identity on $Re_{11}$ and~$Re_{22}$.
Since $\Aut_R(Re_{21}) = R^\times$, we can compose our automorphism with an inner automorphism and
thereby assume that $\varphi(e_{21}) = e_{21}$.
It must then also restrict to the identity on the remaining summand.
Indeed, we may write $\varphi(X_i e_{12}) = x_ie_{12}$ for some~$x_i \in R$, and we have
\[
  x_{i}e_{11}=x_ie_{12}e_{21}=\varphi(X_i e_{12})e_{21}=\varphi(X_i e_{12}e_{21})=\varphi(X_i)\varphi(e_{11})=X_ie_{11},
\]
so that~$x_i=X_i$, as required.
\end{proof}

\subsection{Representation theory}
We introduce a range of notation and results related to the basics
of the representation theory
of $K \bigl(= \GL_2(\Z_p)\bigr)$
and~$G \bigl(= \GL_2(\Q_p)\bigr)$.

\subsubsection{Abelian categories of representations.}
We say that an $\cO$-module is {\em locally torsion} if each element is annihilated by some power of the uniformizer $\varpi$ of~$\cO$.
As usual,
if~$\Gamma$ is a locally profinite group,
then we say that a representation on a locally torsion $\cO$-module is {\em smooth}
if every element is fixed by an open subgroup of~$\Gamma$,
and that a smooth representation is furthermore {\em admissible} if for any
open subgroup of~$\Gamma$, its submodule of invariants is finite over~$\cO$. 
Since we are working with representations on
$p$-power torsion modules, it suffices to check this
for a single open pro-$p$-subgroup of~$\Gamma$ (assuming that $\Gamma$ admits
such a subgroup).  
Following~\cite[Defn.\ 2.2.15]{MR2667882},
we say that a smooth representation~$\pi$ of $\Gamma$ is %
{\em locally admissible} if every vector~$v\in\pi$ generates an admissible representation,
and is %
{\em locally finite} if every vector~$v\in\pi$ generates a representation of finite length.

We write~$\sm.\, \Gamma$ for the category of smooth $\Gamma$-representations on locally torsion $\cO$-modules.
Then~$\sm.\, \Gamma$ is
a Grothendieck category.
(See e.g.\ the proof of~\cite[Lem.~2.2.3]{DEGlocalization}, 
which extends directly to any~$\Gamma$.%
)
We write~$(\sm.\, \Gamma)^{\adm}$,  resp.\ $(\sm.\, \Gamma)^{\ladm}$, $(\sm.\,
\Gamma)^{\fl}, (\sm.\, \Gamma)^{\fg}, %
$ for the full subcategory
of~$\sm.\, \Gamma$ consisting of admissible objects,  resp.\ locally admissible objects, finite length objects, finitely generated objects.

We will be most interested in the case when $\Gamma = G \coloneqq  \GL_2(\Q_p)$. 
We note that,
since the congruence subgroups $K_n$ of $G$ are open and pro-$p$, an object $\pi$ of $\sm. \, G$
is admissible if and only if $\pi^{K_n}$ is $\cO$-finite for one (and hence every) $n \geq 1$.

Throughout the paper we have fixed a continuous character $\zeta:\Qptimes\to\cO^\times$,
and will focus on representations having fixed central character~$\zeta$.
We therefore introduce notation for the various categories of representations having
this fixed central character, which will be in force throughout the paper.
Namely,
we write~$\cA$ for the full subcategory of~$\sm.\, G$ consisting of objects with
central character~$\zeta$, and similarly~ $\cA^{\adm}$, $\cA^{\ladm}, \cA^{\fl}, \cA^{\fg}$.
Then~$\cA^{\ladm}$ coincides with the category of locally finite objects
of~$\cA$, by~\cite[Thm.\ 2.2.17]{MR2667882}.
We have the following results about its structure; recall (see Section~\ref{subsec:locally finite categories}) that a {\em locally finite category} is a Grothendieck category admitting
a set of generators of finite length. 

\begin{lemma}\label{properties of A^ladm} \leavevmode
\begin{enumerate}
\item $\cA^{\ladm}$ is a locally finite category.
\item
The inclusion $\cA^{\fl} \to \cA^{\ladm}$ induces an equivalence
\begin{equation}
\label{eqn:ladm as Ind}
\Ind \cA^{\fl} \iso \cA^{\ladm}.
\end{equation}
\item Let~$E'/E$ be a finite extension with ring of integers~$\cO'$.
Then extension and restriction of scalars define functors $\cA_\cO^{\ladm} \to \cA_{\cO'}^{\ladm}$ and $\cA_{\cO'}^{\ladm} \to \cA_\cO^{\ladm}$.
If~$\pi_1, \pi_2 \in \cA_{\cO}^{\ladm}$, then 
\[
\Hom_{\cA_\cO^{\ladm}}(\pi_1, \pi_2) \otimes_\cO \cO' \isoto \Hom_{\cA_{\cO'}^{\ladm}}(\pi_1 \otimes_\cO \cO', \pi_2 \otimes_\cO \cO').
\]
\end{enumerate}
\end{lemma}
\begin{proof}
Since every object of~$\cA^{\ladm}$ is locally finite, the set of objects of~$\cA^{\ladm}$ of finite length is a set of generators of~$\cA^{\ladm}$.
Since~$\cA^{\ladm}$ is an exact abelian subcategory of~$\cA$ (in the sense recalled in Appendix~\ref{subsubsec:Serre subcategories}),
 and $\cA^{\ladm}$ is closed under colimits in~$\cA$, and~$\cA$ is a Grothendieck category, we see that~$\cA^{\ladm}$ is a Grothendieck category.
Hence~$\cA^{\ladm}$ is a locally finite category, proving part~(1).
Part~(2) is then a consequence of Lemma~\ref{locally finite categories}(2).
We now prove part~(3).
The same argument as~\cite[Lemma~5.1]{MR3150248} shows that for all~$\pi_1, \pi_2 \in \cA_\cO$ we have
\[
  \Hom_{\cA_\cO}(\pi_1, \pi_2) \otimes_\cO \cO' \isoto \Hom_{\cA_{\cO'}}(\pi_1 \otimes_\cO \cO', \pi_2 \otimes_\cO \cO').
\]
(Note that  since $\cO'$ is finite
free over~$\cO$, we have $\Hom_\cO(V, W) \otimes_\cO \cO' \isoto \Hom_{\cO'}(V \otimes_\cO \cO', W \otimes_\cO \cO')$ for all $\cO$-modules~$V, W$.)
Hence, for all~$\pi \in \cA_\cO$, we have $\pi^K \otimes_\cO \cO' = (\pi \otimes_\cO \cO')^K$.
This immediately implies that $\text{--}\otimes_\cO \cO'$ preserves locally admissible representations.
Finally, if~$\pi \in \cA_{\cO'}^{\ladm}$, then~$\pi$ is also locally admissible when viewed as an object of $\cA_\cO^{\ladm}$ by restriction of scalars: in fact,
we immediately reduce to the case in which~$\pi$ is admissible, and then $\pi^K$ is a finitely generated $\cO'$-module, hence a finitely generated $\cO$-module.
\end{proof}

\subsubsection{Completed group rings.}\label{subsubsec:completed group rings}
If~$\Gamma^\circ$ is a profinite group, we write~$\cO\llbracket \Gamma^\circ\rrbracket $ for the usual Iwasawa algebra (i.e.\
completed group ring) of~$\Gamma^\circ$ with coefficients in~$\cO$.
It is a profinite $\cO$-algebra.
If~$\Gamma^\circ$ is a compact $p$-adic analytic group,
then %
$\cO\llbracket\Gamma^\circ\rrbracket$ is Noetherian~\cite[Prop.\ V.2.2.4]{Lazard}.
We also refer to~\cite[Section~2.1]{MR2667882} for a recollection
of some of the basic theory of these rings.

In fact, we will frequently consider Iwasawa coefficients with coefficients,
especially in the case when $\Gamma^{\circ} = K$. We recall the definition.

\begin{defn}\label{defn:Iwasawa-algebra-coefficients}
If $A$ is an~$\cO$-algebra, then
$A\llbracket \Gamma^{\circ}\rrbracket $ denotes the topological $A$-algebra 
  \[A\llbracket \Gamma^{\circ} \rrbracket \coloneqq \varprojlim_{H}A[\Gamma^{\circ}/H],\]
where $H$ runs over a neighbourhood basis of normal open subgroups
of~$\Gamma^{\circ}$, each of the group rings $A[\Gamma^{\circ}/H]$ is endowed
with the discrete topology, and~$A\llbracket \Gamma^{\circ} \rrbracket $ is endowed
with the projective limit topology.
\end{defn}

We next record some properties %
of~$A\llbracket K\rrbracket $ (under the additional hypothesis that $A$ is a finite
type $\cO/\varpi^a$-algebra for some $a \geq 1$).
Before doing so,
we
recall from e.g.\ \cite[Definition~II.2.1.1]{HvO} 
that a filtered $\cO$-algebra $(R, \Fil^i R)_{i \in \bZ}$ (possibly non-commutative, decreasingly filtered, and with $\Fil^0 = R$) 
is called \emph{Zariskian} if $\Fil^{1}R \subset \Rad(R)$ and the Rees ring
\[
\Rees(R) \coloneqq  \bigoplus_{i \in \bZ} \Fil^i R
\]
is Noetherian.
We will use~\cite{HvO} as our main reference for this material, with the caveat that we will work with decreasing filtrations,
whereas \emph{loc.\ cit.} uses increasing filtrations.
A filtered $R$-module~$M$ is \emph{good} if
$\Rees(M)$ is finitely generated over $\Rees(R)$,
which,
by~\cite[Lem.\ I.5.4]{HvO}, is equivalent to
the existence of
generators $m_1, \ldots, m_n$ of~$M$ as an $R$-module, and integers~$k_1, \ldots, k_n$, such that
\[
\Fil^i M = (\Fil^{i-k_1} R)m_1 + \cdots + (\Fil^{i-k_n}R)m_n.
\]

\begin{lemma}\label{lem: properties of completed group ring O K}
Let~$\fa$ be the kernel of $\cO\llbracket K_1\rrbracket  \to \F$, i.e.\ the two-sided ideal generated by~$\varpi$ and the augmentation ideal of~$\cO\llbracket K_1\rrbracket $.
Then the profinite topology on $\cO\llbracket K_1\rrbracket $ and~$\cO\llbracket K\rrbracket $ is induced by the $\fa$-adic filtration, which is Zariskian.
\end{lemma}
\begin{proof}
First, recall that since $K_1$ is a uniform pro-$p$ group, the rings $\cO\llbracket K_1 \rrbracket$ and $\gr_\fa(\cO\llbracket K_1 \rrbracket)$ are Noetherian: 
since the $\fa$-adic filtration on $\cO\llbracket K_1 \rrbracket$ is induced by a $p$-valuation on~$G$, by~\cite[Lem.\ 3.24, 3.25]{MR1924402},
this is a consequence of~\cite[Thm.\ III.2.3.3]{Lazard}.
Since $\fa = \Rad(\cO\llbracket K_1\rrbracket )$, we now see that the topology on $\cO\llbracket K_1\rrbracket $ is $\fa$-adic, by Lemma~\ref{properties of compact modules}~(8).
Finally, since $\cO\llbracket K_1\rrbracket $ is $\fa$-adically complete, \cite[Prop.\ II.2.2.2]{HvO} shows that $\cO\llbracket K_1\rrbracket $ is Zariskian, because the associated graded
ring $\gr_\fa \cO\llbracket K_1\rrbracket $ is Noetherian.
This concludes the proof for~$K_1$.

Since $\cO\llbracket K\rrbracket $ is a finite $\cO\llbracket K_1\rrbracket $-module, its $\fa$-adic filtration induces the profinite topology, again by Lemma~\ref{properties of compact modules}~(8), 
and is a good filtration (by definition). 
By~\cite[Thm.\ I.5.7]{HvO}, we conclude that $\gr_\fa(\cO\llbracket K\rrbracket )$ is a finite $\gr_\fa(\cO\llbracket K_1\rrbracket )$-module.
Hence $\gr_\fa(\cO\llbracket K\rrbracket )$ is a Noetherian ring, which by another application of \cite[Prop.\ II.2.2.2]{HvO} shows 
that~$\cO\llbracket K\rrbracket $ is Zariskian with respect to the
$\fa$-adic filtration.
\end{proof}

\begin{lem}
  \label{lem:topology-Iwasawa-algebra}
  Let~$A$ be a finite type $\cO/\varpi^a$-algebra for some~$a \geq 1$.
  \begin{enumerate}
    \item\label{item: Iwasawa 1} $A\llbracket K\rrbracket $ is topologically isomorphic to the $\fa$-adic completion of $A \otimes_\cO \cO\llbracket K\rrbracket $.
    \item\label{item: Iwasawa 2} The rings $A \otimes_\cO \cO\llbracket K\rrbracket $ and $A\llbracket K\rrbracket $ are Noetherian.
    \item\label{item: Iwasawa 3} $A\llbracket K\rrbracket $ is flat over~$A \otimes_\cO \cO\llbracket K\rrbracket $.
    \item\label{item: Iwasawa 4} The projective limit topology on~$A\llbracket K\rrbracket $ coincides with its $\mf{a}$-adic topology. 
    \item\label{item: Iwasawa 6} Every morphism in $\Mod^{\fp}(A\llbracket K\rrbracket )$ is strict with respect to the $\fa$-adic topology 
    (in the sense that the quotient topology on the coimage coincides with the subspace topology on the image). 
    The functor 
    \begin{equation}\label{eqn: tilde functor on Iwasawa modules}
    \Mod^{\fp}(A\llbracket K\rrbracket ) \to \Pro \Mod^{\fp}(A),\ M \mapsto \tld M \coloneqq \quoteslim{i} M/\fa^i M 
    \end{equation}
    is exact.
    \item\label{item: Iwasawa 7} 
Let $\AKModfpA$ denote the abelian category of~$A\llbracket K \rrbracket$-module objects in~$\Mod^{\fp}(A)$,
so that the functor~\eqref{eqn: tilde functor on Iwasawa modules}
factors through the forgetful functor $\Pro \AKModfpA \to \Pro \Mod^{\fp}(A)$.
Then
if \[X_1 \to X_2 \to X_3 \to X_4 \to X_5\] is an exact sequence
in~$\Pro \AKModfpA$
    and~$X_1, X_2, X_4, X_5$ are in the essential image of~\eqref{eqn: tilde functor on Iwasawa modules}, the same is true of~$X_3$.
  \end{enumerate}
\end{lem}
\begin{proof}

Proof of~\eqref{item: Iwasawa 1}: %
By Lemma~\ref{lem: properties of completed group ring O K},
the sequences of ideals~$J_n$ and %
$\fa^n(\cO/\varpi^a)\llbracket K\rrbracket $ are cofinal in $(\cO/\varpi^a)\llbracket K\rrbracket $.
So the pro-objects
\[
\quoteslim{n} (\cO/\varpi^a)[K/K_n] \quad \text{ and } \quad
\quoteslim{n} (\cO/\varpi^a)\llbracket K\rrbracket /\fa^n(\cO/\varpi^a)\llbracket K\rrbracket 
\]
are canonically isomorphic.
Tensoring with~$A$, we see that the same is true of
\[
\quoteslim{n} A[K/K_n] \quad \text{ and } \quad \quoteslim{n} A \otimes_\cO \cO\llbracket K\rrbracket /\fa^n\cO\llbracket K\rrbracket .
\]
Passing to the inverse limit, we see that~$A\llbracket K\rrbracket $ is isomorphic to the $\fa$-adic completion of $A \otimes_{\cO} \cO\llbracket K\rrbracket $, as desired.
This concludes the proof of~(1).

Proof of~\eqref{item: Iwasawa 2}: Observe first that $\Rees_\fa(A \otimes_\cO \cO\llbracket K\rrbracket )$ is Noetherian, since it is a quotient of $A \otimes_\cO \Rees_\fa(\cO\llbracket K\rrbracket )$,
which is Noetherian 
(because~$\Rees_\fa(\cO\llbracket K\rrbracket )$ is Noetherian, by Lemma~\ref{lem: properties of completed group ring O K}, and~$A$ has finite type over~$\cO$).
Thus $A \otimes_\cO \cO\llbracket K\rrbracket $ and $\gr_\fa(A \otimes_\cO \cO\llbracket K\rrbracket )$ are also Noetherian, because they are quotients of $\Rees_\fa(A \otimes_\cO \cO\llbracket K\rrbracket )$.
Finally, \cite[Prop.\ I.7.1.2]{HvO} implies that~$A\llbracket K\rrbracket $ is also Noetherian,
since it is complete with respect to a filtration whose associated graded is 
$\gr_\fa(A \otimes_\cO \cO\llbracket K\rrbracket )$.

Proof of~\eqref{item: Iwasawa 3}: By~\cite[Cor.\ I.5.5, Rem.\ II.1.1.2(1)]{HvO}
and the already observed fact that
$\Rees_\fa(A \otimes_\cO \cO\llbracket K\rrbracket )$ is Noetherian,
we deduce that the $\fa$-adic filtration on $A \otimes_\cO \cO\llbracket K\rrbracket $ has the Artin--Rees property.
The claim is then a consequence of~\cite[Thm.\ II.1.2.4]{HvO}, given that 
$\gr_\fa(A \otimes_\cO \cO\llbracket K\rrbracket )$ is (as we have noted) also Noetherian.

Proof of~\eqref{item: Iwasawa 4}: By part~\eqref{item: Iwasawa 1}, it suffices to prove that
\[
\fa^n A\llbracket K\rrbracket  = \ker(A\llbracket K\rrbracket  \to A \otimes_\cO \cO\llbracket K\rrbracket /\fa^n\cO\llbracket K\rrbracket ).
\]
Since~$\fa^n$ is a finitely generated ideal of~$\cO\llbracket K\rrbracket $, a choice of (right) generators gives an exact sequence
\[
(A \otimes_\cO \cO\llbracket K\rrbracket )^{\oplus i} \to A \otimes_\cO \cO\llbracket K\rrbracket  \to A \otimes_\cO \cO\llbracket K\rrbracket /\fa^n\cO\llbracket K\rrbracket  \to 0
\]
for some~$i \geq 0$.
This sequence stays exact after passing to the $\fa$-adic completion, because the $\fa$-adic filtration on~$A \otimes_\cO \cO\llbracket K\rrbracket $ has the Artin--Rees property
(compare~\cite[Lem.\ II.1.2.5]{HvO}).
This concludes the proof of~\eqref{item: Iwasawa 4}.

Proof of~\eqref{item: Iwasawa 6}: By~\cite[Prop.\ II.2.2.1]{HvO}, the $\fa$-adic filtration on $A\llbracket K\rrbracket $ is Zariskian.
By~\cite[Thm.\ 2.1.2(5)]{HvO}, every inclusion $N \to M$ of finitely generated $A\llbracket K\rrbracket $-modules is a topological embedding with respect to the
$\fa$-adic topology on source and target.
This immediately implies both statements of~\eqref{item: Iwasawa 6}.

Proof of~\eqref{item: Iwasawa 7}: By part~\eqref{item: Iwasawa 6}, it suffices to prove this when~$X_1 = X_5 = 0$, and so, writing~$X \coloneqq X_3$, we have finitely presented 
$A\llbracket K \rrbracket$-modules $M \coloneqq X_2, M' \coloneqq X_4$, and an exact sequence %
in~$\Pro \AKModfpA$
 \begin{equation}\label{eqn: to prove this comes from tilde}
0 \to \tld M \to X \to \tld M' \to 0.
 \end{equation}
Now both $\tld M$ and $\tld M'$ are countably indexed; indeed,
$M = \quoteslim{n} M/\fa^n$, and similarly for~$M'$.
Lemma~\ref{lem:countably indexed}
shows that~\eqref{eqn: to prove this comes from tilde}
can be obtained as a countably indexed limit of short exact sequences in~$\AKModfpA$
$$0 \to M/\fa^n \to X_n \to M'_n \to 0.$$
(Here the projective system $(M'_n)$ need not coincide with 
$(M'/\fa^n)$, but is isomorphic to it in 
$\Pro \AKModfpA$.)
We now pass to the actual limit over~$n$ of these short exact sequences,
to obtain a short exact sequence
$$0 \to M \to Y \to M' \to 0$$
of~$A\llbracket K\rrbracket$-modules.
(The exactness on the right follows from the fact that the transition morphisms
in $(M/\fa^n)$ are surjective.)
We may equip $Y$ with its inverse limit topology (each $X_n$ being endowed 
with its discrete topology), and then (from its construction as a limit)
we see that we have a strict short exact sequence of topological $A\llbracket K\rrbracket$-modules,
in which each of $M$ and $M'$ are endowed with their $\fa$-adic topologies.
Consequently, $Y$ is an object of $\Mod^{\fp}(A\llbracket K\rrbracket)$ (being 
an extension of two such objects).  Furthermore, since 
$Y$ is constructed as a countable limit of discrete spaces,
the topology on $Y$ is completely metrizable, and so~\cite[Prop.~C.6]{emertongeepicture}
shows that it must coincide with the $\fa$-adic topology on~$Y$
(cf.\ the discussion of Remark~\ref{rem:Polish Iwasawa} below).
Thus we obtain a morphism of short exact sequences
$$\xymatrix{
0 \ar[r] & \tld M \ar[r] \ar@{=}[d] & \tld Y \ar[r] \ar[d] & \tld M' \ar[r]\ar@{=}[d] & 0 \\
0 \ar[r] & \tld M  \ar[r]  & X \ar[r]  & \tld M' \ar[r] & 0 }
$$
in $\Pro \AKModfpA$. The five lemma implies that
$\tld Y \iso X$, concluding the proof.    
\end{proof}

\begin{remark}
\label{rem:Polish Iwasawa}
As in the statement 
of Lemma~\ref{lem:topology-Iwasawa-algebra}, assume that $A$
is of finite type over $\cO/\varpi^a$ for some~$a\geq 1$. 
Then since the ring $A$ is then countable,
the completed group ring $A\llbracket K\rrbracket $, endowed with its $\fa$-adic topology, is Polish.
The general theory of finitely generated modules over Noetherian Polish topological
rings (recapitulated, for example, in~\cite[Prop.~C.6]{emertongeepicture})
shows that any finitely generated $A\llbracket K\rrbracket $-module
is endowed with a canonical completely metrizable topology,
and that morphisms between such modules are automatically continuous and strict
with respect to this topology.

Since the topology on $A\llbracket K\rrbracket $ is its $\fa$-adic topology,
this canonical topology on any finitely generated $A\llbracket K\rrbracket $-module
is again the $\fa$-adic topology.
One may then interpret 
Lemma~\ref{lem:topology-Iwasawa-algebra}~\eqref{item: Iwasawa 6}
as providing another proof of the strictness of morphisms with respect to canonical topologies.
\end{remark}

\begin{remark}\label{rem:Polish Iwasawa II}
Suppose that~$A$ is a finite type $\cO/\varpi^a$-algebra for some~$a \geq 1$, that~$M$ is an $A\llbracket K \rrbracket$-module, and that~$M$ is finitely presented over~$A$.
Then~$M$ is also finitely presented over $A\llbracket K \rrbracket$, and its canonical topology is discrete. 
 Indeed, we can write the complete metric space~$M$ as a countable union of closed subsets $M = \bigcup_{x \in M} \{x\}$, and if~$M$ were not discrete, each  subset~$\{x\}$ would have empty interior, 
contradicting the Baire category theorem.
\end{remark}

\begin{remark}\label{rem:Polish Iwasawa III}
Note that the forgetful functor induces an equivalence between $\AKModfpA$ %
and the category $\OKtorsModfpA$ of $\cO\llbracket K \rrbracket$-module objects~$M$ in~$\Mod^{\fp}(A)$ which are $\fa$-power torsion (in the sense that the $\cO\llbracket K \rrbracket$-action
factors through $\cO\llbracket K \rrbracket/\fa^i$ for some~$i$, depending on the module). 
Indeed, the $\fa$-adic topology on any $M \in \AKModfpA$ is necessarily discrete, by Remark~\ref{rem:Polish Iwasawa II}, hence~$M$ is $\fa$-power torsion.
The inverse functor endows each %
$N \in \OKtorsModfpA$ 
with its $A \otimes_\cO \cO\llbracket K \rrbracket$-module structure, which uniquely extends to an $A\llbracket K \rrbracket$-module
structure, since~$N$ is $\fa$-power torsion.
Hence, in the statement of Lemma~\ref{lem:topology-Iwasawa-algebra}~\eqref{item: Iwasawa 7}, we may replace~$\AKModfpA$ with the category of 
$\fa$-power torsion $\cO\llbracket K \rrbracket$-modules in~$\Mod^{\fp}(A)$.
\end{remark}

\begin{remark}\label{rem:Polish Iwasawa IV}
Write  $\OKProModfpA$ for the category of $\cO\llbracket K \rrbracket$-modules in~$\Pro \Mod^{\fp}(A)$.
The forgetful functor $\AKModfpA \to \Mod^{\fp}(A)$ is exact and faithful, so its Pro-extension factors through an exact and faithful functor
\[
\Pro \AKModfpA \to \OKProModfpA.
\]
We claim that this functor is also full.
Note that~$\OKProModfpA$ has cofiltered limits, and that the forgetful functor $\OKProModfpA \to \Pro \Mod^{\fp}(A)$ preserves them;
so it suffices to prove that if $M = \quoteslim{i} M_i \in \Pro \AKModfpA$, and $N \in \AKModfpA$, then every $\cO\llbracket K \rrbracket$-linear morphism
$\varphi: M \to N$ in $\OKProModfpA$ factors through a morphism $M_i \to N$ in~$\AKModfpA$ for some~$i$.

By definition, $\varphi$ factors through an $A$-linear morphism $\varphi_j: M_j \to N$ for some~$j$, which may not be $\cO\llbracket K \rrbracket$-linear (i.e.\ it may not be a morphism in~$\AKModfpA$).
However, the $\cO\llbracket K \rrbracket$-linearity of~$\varphi$ 
shows that for all~$x \in K$ there exists~$i_x \geq j$ such that
$x\varphi_j - \varphi_j x$ vanishes on~$M_{i_x}$.
Since~$M_j$ and~$N$ are $\fa$-power torsion, there exists~$n$ such that~$K_{n}$ acts trivially on both, and so~$\varphi_j$ is automatically $K_n$-linear.
Choose now representatives~$x_1, \ldots, x_t$ of~$K/K_n$ in~$K$, as well as some $i \geq \sup\{i_{x_1}, \ldots, i_{x_t}\}$.
Then $M_i \to M_j \xrightarrow{\varphi_j} N$ is~$\cO\llbracket K \rrbracket$-linear, as desired.
\end{remark}

\subsubsection{Smooth representations of compact \texorpdfstring{$p$}{p}-adic analytic groups}\label{subsubsec: contragredient}
If~$\Gamma^\circ$ is a compact $p$-adic analytic group,
then Pontrjagin duality induces an equivalence of categories
\begin{equation}\label{eqn: Pontrjagin duality}
(\text{--})^\vee: (\sm.\, \Gamma^\circ)^{\op} \iso \Mod_c(\cO\llbracket\Gamma^\circ\rrbracket^\op) \cong \Mod_c(\cO\llbracket\Gamma^\circ\rrbracket).
\end{equation}
More precisely, when we apply Pontrjagin duality to
a 
smooth $\Gamma^\circ$-representation,
we typically %
convert the resulting right $\cO\llbracket\Gamma^\circ\rrbracket$-module structure
to a left $\cO\llbracket\Gamma^\circ\rrbracket$-module structure by applying the continuous anti-involution
of $\cO\llbracket\Gamma^\circ\rrbracket$ induced by $g \mapsto g^{-1}$ 
(and this anti-involution is what provides the isomorphism 
between the second and third of the categories in~\eqref{eqn: Pontrjagin duality}).
However, in some contexts (particularly when forming tensor products,
as we will do extensively below), it will be convenient to work directly with
right $\cO\llbracket\Gamma^\circ\rrbracket$-modules.

The equivalence~\eqref{eqn: Pontrjagin duality} restricts to an equivalence
\begin{equation}\label{eqn: admissible Pontrjagin duality}
\bigl((\sm.\, \Gamma^\circ)^{\adm}\bigr)^{\op} \iso \Mod^{\fg}(\cO\llbracket\Gamma^\circ\rrbracket),
\end{equation}
where we regard the target as a full subcategory
of $\Mod_c(\cO\llbracket\Gamma^\circ\rrbracket)$ via~\eqref{eqn:mod-fp-to-mod-c}. %

\subsubsection{Smooth \texorpdfstring{$\Gamma$}{\Gamma}-representations as
\texorpdfstring{$\cO\llbracket \Gamma\rrbracket $}{O[\Gamma]}-modules
in the non-compact case}
If~$\Gamma$ is a $p$-adic analytic group (compact or not), if $\Gamma^\circ \subset \Gamma$ is a compact open subgroup, 
and if $\pi$ is an object of $\sm.\, \Gamma$, then the $\Gamma^\circ$-action on~$\pi$
induces a unique continuous $\cO\llbracket\Gamma^\circ\rrbracket$-action on~$\pi$ (where we endow $\pi$
with its discrete topology).  
This action can be extended to a suitably completed group ring
$\cO\llbracket\Gamma \rrbracket$ of~$\Gamma$.
We recall its definition, in the more general context where we allow coefficients.  

\begin{defn}\label{defn: Iwasawa algebra of G}
If~$A$ is an~$\cO$-algebra
then we let 
\[
A\llbracket \Gamma\rrbracket \coloneqq \cO[\Gamma] \otimes_{\cO[\Gamma^\circ]} A\llbracket \Gamma^\circ \rrbracket
\]
with its natural ring structure, as explained in 
~\cite[\S3]{MR4106887}.
\end{defn}

\begin{lemma}
\label{lem:A[[Gamma]] actions}
If $M$ is a Hausdorff topological $A$-module equipped with  
an action of $\Gamma$ by continuous automorphisms, 
and with a continuous $A \llbracket \Gamma^{\circ}\rrbracket$-module structure,
such that the two resulting $A[\Gamma^{\circ}]$-module structures on $M$ coincide,
then the induced $A[\Gamma]$-module structure on $M$ extends uniquely
to an $A\llbracket \Gamma\rrbracket$-module structure.
\end{lemma}
\begin{proof}
This is easily verified; see e.g.~\cite[Prop.~2.23]{MR4551876}, whose 
proof carries over verbatim.
\end{proof}

As a special case of~Lemma~\ref{lem:A[[Gamma]] actions},
we deduce the claim made above, namely that
if~$\pi \in \sm.\, \Gamma$, then the action of $\cO\llbracket \Gamma^\circ \rrbracket$ on~$\pi$ induces an $\cO\llbracket \Gamma\rrbracket $-action 
on~$\pi$. Thus we obtain a fully faithful embedding $\sm.\, \Gamma
\hookrightarrow \Mod(\cO\llbracket \Gamma\rrbracket )$.
The following lemma shows that this realizes $\sm.\, \Gamma$ as a localizing subcategory of ~$\Mod(\cO\llbracket \Gamma\rrbracket )$ (see Appendix~\ref{subsec:1-categories} for a brief recollection of this notion).

\begin{lemma}\label{embedding of smooth in modules}
The fully faithful embedding $\sm.\, \Gamma \hookrightarrow \Mod(\cO\llbracket \Gamma\rrbracket )$
realizes $\sm.\, \Gamma$ as a localizing subcategory of $\Mod(\cO\llbracket \Gamma\rrbracket )$.
\end{lemma}

\begin{rem}
The main result of~\cite{Heyersmooth} shows that $\sm.\, \Gamma$ is even a localizing subcategory of $\Mod(\cO[\Gamma])$, i.e.\ that extensions of smooth representations are automatically smooth.
However, we will not need this more general fact, but will satisfy ourselves
with recalling the (standard) proof of the easier result that we will be using.  
\end{rem}

\begin{proof}[Proof of Lemma~{\em \ref{embedding of smooth in modules}}]
It is immediate that the essential image of~$\sm. \, \Gamma$ is closed under
subobjects, quotients, and colimits, and in particular under arbitrary direct sums.
So it suffices to prove that if~$V$ is an extension in~$\Mod(\cO\llbracket \Gamma\rrbracket )$ of two objects~$\pi_1, \pi_2 \in \sm.\, \Gamma$, then~$V$ is smooth.
Choose~$v \in V$, let~$\Gamma^\circ \subset \Gamma$ be a pro-$p$ open subgroup that fixes the image of~$v$ in~$\pi_2$, %
and let~$n \geq 0$ be such that $\varpi^n v = 0$.
Replacing~$\pi_2$ by its cyclic $\cO\llbracket\Gamma^\circ\rrbracket$-submodule generated by the image of~$v$, which is isomorphic to $\cO/\varpi^n$ with the trivial $\Gamma^\circ$-action, 
it suffices to prove that 
\[
\Ext^1_{\sm.\, \Gamma^\circ}(\cO/\varpi^n, \pi_1) \to \Ext^1_{\Mod\cO\llbracket\Gamma^\circ\rrbracket}(\cO/\varpi^n, \pi_1)
\]
is an isomorphism (equivalently, a surjection); indeed, it follows that~$v$ is contained in a smooth $\Gamma^\circ$-submodule of~$V$.

We claim that it suffices to consider the case that~$\pi_1$ is finitely generated in $\sm.\, \Gamma^\circ$ (and so is a finitely generated torsion $\cO$-module);
indeed, $\pi_1$ is the filtered colimit of its finitely generated subobjects in~$\sm.\, \Gamma^\circ$, and $\Ext^1_{\sm.\, \Gamma^\circ}(\cO/\varpi^n,\text{--})$, $ \Ext^1_{\Mod\cO\llbracket\Gamma^\circ\rrbracket}(\cO/\varpi^n, \text{--})$ commute with filtered colimits,
because the categories $\sm.\,\Gamma^\circ$ and $\Mod \cO\llbracket\Gamma^\circ\rrbracket$ are locally Noetherian (see e.g.~\cite[Proposition~A.1.1]{DEGlocalization}).
If~$\pi_1$ is finitely generated, then every extension of $\cO/\varpi^n$ by~$\pi_1$ in $\Mod(\cO\llbracket\Gamma^\circ\rrbracket)$ is a finitely generated torsion $\cO$-module, and so its canonical topology as an $\cO\llbracket\Gamma^\circ\rrbracket$-module 
is discrete, which implies that the action of~$\Gamma^\circ$ is smooth, as required. 
\end{proof}

\subsubsection{Smooth \texorpdfstring{$G$}{G}-representations as \texorpdfstring{$\cO\llbracket G\rrbracket $}{O[G]}-modules}
\label{subsubsec:smooth to O-G }
When~$\Gamma = G = \GL_2(\bQ_p)$
and $A$ is an~$\cO$-algebra,
we write~$A\llbracket G\rrbracket _\zeta$ for the quotient of~$A\llbracket G\rrbracket $ by the ideal generated by $z-\zeta(z)$ for~$z \in Z(G) \cong \bQ_p^\times$.
The category $\cA$ may thus be regarded as a full subcategory of
$\Mod(\cO\llbracket G\rrbracket _\zeta)$, and the inclusion $\cA\subseteq \Mod(\cO\llbracket G\rrbracket _\zeta)$ is
exact and colimit-preserving
(see e.g.\ \cite[\S E.2]{emerton2023introduction}).
We recall the following result about the structure of $\cO\llbracket G \rrbracket$.

\begin{lem}
  \label{lem:group-ring-coherent} The ring~$\cO\llbracket G\rrbracket _\zeta$ is coherent.
\end{lem}
\begin{proof}
  This is a consequence of~\cite[Theorem~9.7]{MR4551876} (see
  also~\cite[Theorem~1.2]{MR4106887}), which states that~$\cO\llbracket G\rrbracket $ is coherent,
  together with the fact that the quotient of a coherent ring by an ideal that
  is generated by finitely many central elements is
  coherent~\cite[Corollary~3.6]{MR3981694}.
\end{proof}

\begin{remark}\label{embedding of A in modules}
In the context of~$\Gamma = G = \GL_2(\bQ_p)$, an evident variant
of Lemma~\ref{embedding of smooth in modules} shows that the essential image of the fully faithful embedding 
$\cA \hookrightarrow \Mod(\cO\llbracket G\rrbracket _{\zeta})$
is again a localizing subcategory of its target.
\end{remark}

\begin{rem}\label{fg equals fp for cA}
An object~$\pi$ of~$\cA$ is \emph{finitely presented} if it is the cokernel of a map $\cInd_{KZ}^G(\sigma_1) \to \cInd_{KZ}^G(\sigma_2)$ where~$\sigma_1, \sigma_2$ are smooth $KZ$-representations
on $\cO$-finite modules.
By~\cite[Proposition~3.4]{MR4106887}, $\pi$ is finitely presented if and only if the image of~$\pi$ under the embedding in Remark~\ref{embedding of A in modules} is finitely presented as an $\cO\llbracket G\rrbracket _\zeta$-module.
Since $G = \GL_2(\bQ_p)$, $\pi$ is finitely presented if and only if it is finitely generated (see e.g.\ \cite[Theorem~2.2.1]{DEGlocalization} for a proof of this fact).
For this reason, and for consistency with our notation for modules, we will often write~$\cA^{\fp}$ for the category we have denoted~$\cA^{\fg}$.
We can thus summarize our discussion by saying that the inclusion of $\cA$ in $\Mod(\cO\llbracket G \rrbracket_\zeta)$ induces an inclusion 
\begin{equation}
\label{eq:Afp-in-Modfp-OG}
\cA^{\fp}\subseteq  \Mod^{\fp}(\cO\llbracket G\rrbracket _\zeta).
\end{equation}  
\end{rem}

\subsubsection{Smooth \texorpdfstring{$G$}{G}-representations and duality}
Pontrjagin duality induces an equivalence of
$(\sm.\, G)^{\op}$
with the category $\Mod_G^{\mathrm{pro\,aug}}(\cO)$ introduced in~\cite[Section~2.1]{MR2667882}.
The objects of~$\Mod_G^{\mathrm{pro\,aug}}(\cO)$ are %
(left or right, compare Section~\ref{subsubsec: contragredient}) $\cO\llbracket G\rrbracket $-modules 
endowed 
with a profinite topology that admits a neighborhood basis of the identity consisting of open
topological $\cO\llbracket K\rrbracket $-modules, 
while the
morphisms in this category
are the continuous $\cO\llbracket G\rrbracket $-linear maps.
The forgetful functor $\Mod_G^{\mathrm{pro\,aug}}(\cO) \to \Mod(\cO\llbracket G\rrbracket )$
is exact.
There are also forgetful functors $\Mod_G^{\mathrm{pro\,aug}}(\cO) \to \Mod_c(\cO\llbracket K\rrbracket ) \to \Mod_c(\cO)$,
which factor through the category of $\cO\llbracket G\rrbracket $-modules in $\Mod_c(\cO\llbracket K\rrbracket )$, resp.\ $\Mod_c(\cO)$.
When working with the category $\Mod_G^{\mathrm{pro\,aug}}(\cO)$,
or one of its full subcategories,
we will sometimes write $\Hom_{\cO\llbracket G\rrbracket }^{\cont}$
rather than~$\Hom_{\Mod_G^{\mathrm{pro\,aug}}(\cO)}$.

We denote the image of $\cA^{\op}$ in~$\Mod_G^{\mathrm{pro\,aug}}(\cO)$ by $\Mod_{G, \zeta}^{\mathrm{pro\,aug}}(\cO)$,
and the image of
$\cA^{\ladm, \op}$ by~$\fC$. 
This category is denoted~$\fC(\cO)$ in~\cite{MR3150248}, but as we
will not work with more general coefficients, we shorten the notation accordingly. %
It is a Serre subcategory of~$\Mod_{G, \zeta}^{\mathrm{pro\,aug}}(\cO)$. %
\begin{remark}
\label{rem:contragredient again}
If we regard the objects of~$\fC$ as left~$\cO\llbracket G\rrbracket $-modules,
following the discussion in Section~\ref{subsubsec: contragredient}, then they have central character~$\zeta^{-1}$.
If we regard them as right~$\cO\llbracket G\rrbracket $-modules (again as in that discussion),
then they have central character~$\zeta$, i.e.\ they are right $\cO\llbracket G\rrbracket _{\zeta}$-modules.
\end{remark}

Since locally admissible representations of~$G$ are locally finite, every object~$P \in \fC$ is the cofiltered limit of its finite length quotients in~$\fC$.
We refer to the image of $\cA^{\adm,\op}$ in~$\fC$ as the (full) subcategory
of {\em coadmissible}~modules,
and denote it by~$\fC^{\coadm}.$
Equivalently, this is the full subcategory of~$\fC$ consisting of objects
that are finitely generated over~$\cO\llbracket K\rrbracket $,
or (again equivalently) over $\cO\llbracket K_n\rrbracket $ for one, or any,~$n\geq 1$.
Note that the topology on an object of $\fC^{\coadm}$ coincides with the canonical 
compact topology it inherits as a finitely generated module over the compact Noetherian
ring~$\cO\llbracket K\rrbracket $, by the uniqueness part of Lemma~\ref{properties of compact modules}~\eqref{item: compact 5}.
By~\eqref{eqn:ladm as Ind},
we see that
formation of projective limits (equipped with their projective limit topologies) induces
an equivalence
\begin{equation}
\label{eqn:fC as a Pro category}
\Pro \fC^{\fl} \iso \fC,
\end{equation}
where, as usual, $\fC^{\fl}$ denotes the full subcategory of finite length objects of~$\fC$.
The following ``automatic continuity'' result is often useful.

\begin{lemma}\label{continuousadmissible}
If $M$ and $N$ are objects of~$\fC$,
with $M$ furthermore lying in~$\fC^{\coadm},$
then the evident inclusion
\[
\Hom_{\fC}(M,N) \, \bigl(=\Hom_{\cO\llbracket G\rrbracket }^{\cont}(M, N)\bigr)
\subseteq \Hom_{\cO\llbracket G\rrbracket }(M, N)
\]
is an equality.
\end{lemma}
\begin{proof}
It suffices to prove that $\cO\llbracket K\rrbracket $-linear morphisms between~$M$ and~$N$ are automatically continuous.
Since~$M$ and~$N$ are objects of $\Mod_c(\cO\llbracket K\rrbracket )$, and~$M$ is finitely generated over~$\cO\llbracket K\rrbracket $, this is a consequence of~\cite[Proposition~3.5]{MR1474172}.
\end{proof}

We also have the following related result.

\begin{lemma}
\label{lem:fC to O-G- }
The forgetful functor $\fC^{\coadm} \to \Mod(\cO\llbracket G\rrbracket _{\zeta})$
is fully faithful, and its essential image is a Serre subcategory of 
its target.
\end{lemma}
\begin{proof}
The claimed full faithfulness is an immediate consequence
of Lemma~\ref{continuousadmissible}.
If $M$ is an object of~$\fC^{\coadm}$,
then (by definition) $M$ is finitely generated over the Noetherian profinite $\cO$-algebra~$\cO\llbracket K\rrbracket $,
and (as noted above) is equipped with the canonical compact topology 
it inherits by virtue of this.
Thus,
any $\cO\llbracket G\rrbracket _{\zeta}$-submodule of $M$ is again finitely generated over~$\cO\llbracket K\rrbracket $,
and so is an object of~$\fC^{\coadm}$. Similar arguments in the case of quotients and extensions
show that $\fC^{\coadm}$ is indeed a Serre subcategory of~$\Mod(\cO\llbracket G\rrbracket _{\zeta}).$
\end{proof}

\subsubsection{Induced representations}\label{induced representations}
We write $\Ind_B^G: \sm.\, T \to \sm.\, G$ for the functor of parabolic induction with respect to~$B$, and $\cInd_{KZ}^G: \sm.\, KZ \to \sm.\, G$ for the functor of compact induction.
We adopt similar notation for other groups, and for the induction functors on the categories of smooth representations with fixed central character 
(note that~$\Ind_B^G$ and~$\cInd_{KZ}^G$ preserve central characters when these are defined).
Then~$\Ind_B^G$ is an exact functor, and~$\cInd_{KZ}^G$ is an exact left adjoint to~$\Res^G_{KZ}$.
The parabolic induction~$\Ind_B^G(\triv_T)$ of the trivial character of~$T$ is a reducible representation, and we write~$\St$ for its unique irreducible quotient, the Steinberg representation of~$G$.
There is an exact sequence
\[
0 \to \triv_G \to \Ind_B^G(\triv_T) \to \St \to 0.
\]
If~$\sigma$ is a $\zeta$-compatible Serre weight, and~$\chi = \sigma^{\Iw_1}$ is the highest weight of~$\sigma$ viewed as a character of~$T_0 Z$, then we write
\begin{gather*}
\cH_G(\sigma) \coloneqq  \End_G(\cInd_{KZ}^G\sigma)\\
\cH_T(\chi) \coloneqq  \End_T(\cInd_{T_0 Z}^T \chi).
\end{gather*}
There are canonical isomorphisms
\begin{equation}\label{eqn: Satake isomorphism for the torus}
\cH_T(\chi) \iso \cH_T(\triv_{T_0Z}) \iso \F[\HeckeX^{\pm 1}]
\end{equation}
(so in particular $\cH_T(\chi)$ is 
canonically independent of $\chi$).
We normalize our choice of generator~$\HeckeX$ in such a way that the Satake map
$\cH_G(\sigma) \to \cH_T(\chi)$
of \cite{HerzigSatake}
is an injection which sends the usual spherical Hecke operator to~$\HeckeX$, and
determines an isomorphism \begin{equation}\label{eqn:Satake-iso-with-coordinate-X} \cH_G(\sigma)\isoto\F[\HeckeX].\end{equation}
If $\sigma$ is not a twist of $\Sym^0$, %
then~\cite[Theorem~25]{BarthelLivneDuke} gives an isomorphism
\begin{equation}\label{universal parabolic induction}
\left ( \cInd_{KZ}^G\sigma \right ) [1/\HeckeX] \iso \Ind_{\lbar B}^G\left ( \cInd_{T_0 Z}^T\chi \right ),
\end{equation}
where $\lbar B$ is the lower-triangular Borel subgroup.
This isomorphism is compatible with the Satake isomorphism,
in the sense that it is an isomorphism of $\F[\HeckeX^{\pm 1}]$-modules (acting on the right-hand side by $\Ind_{\lbar B}^G$-functoriality),
and so can be regarded as characterizing the Satake morphism (for $\sigma$ not a twist 
of $\Sym^0$).

If $\sigma\coloneqq \sigma_{a,0}$ is a twist of $\Sym^0$,
so that $\sigma_{a,p-1}$ is the corresponding twist
of $\Sym^{p-1}$,
then there is a short exact sequence
\begin{equation}\label{Steinberg change of weight}
0 \to \cInd_{KZ}^G \sigma_{a,0} \to \cInd_{KZ}^G \sigma_{a,p-1} \to
\bigl(\omega^a\nr_{\HeckeX^2-\zeta(p)}\circ \det\bigr) \otimes \St
\to 0,
\end{equation}
with the first embedding being~$\F[\HeckeX]$-linear.
Similarly, there is a short exact sequence
\begin{equation}\label{trivial change of weight}
0 \to \cInd_{KZ}^G \sigma_{a,p-1} \to \cInd_{KZ}^G \sigma_{a,0} \to
\bigl(\nr_{\HeckeX^2-\zeta(p)}\circ \det\bigr) \to
\to 0,
\end{equation}
where again the first embedding is $\F[\HeckeX]$-linear.
The composite of the injections $\cInd_{KZ}^G \sigma_{a,0} \hookrightarrow
\cInd_{KZ}^G \sigma_{a, p-1}$ and $\cInd_{KZ}^G \sigma_{a,p-1}\hookrightarrow
\cInd_{KZ}^G \sigma_{a,0}$, in either order, is a non-zero scalar multiple of
multiplication by $\HeckeX^2 -\zeta(p)$.

\subsubsection{Irreducible objects of~\texorpdfstring{$\cA$}{A}.}\label{irreducibles}
Recall the following classification of the irreducible objects
of~$\cA$. 
\begin{thm}\label{thm:classification-of-irreps}\leavevmode
\begin{enumerate}
\item\label{item:35} Every irreducible object of~$\cA$ is isomorphic to one of the following
representations: 
\begin{enumerate}
\item\label{item:32} 
   $\cInd_{KZ}^G\sigma_{a, b}/f(\HeckeX)$ for some
    $\zeta$-compatible~$\sigma_{a, b}$ and monic irreducible
    $f \in \cH_G(\sigma_{a, b})$, such that
    $f$ is coprime to $\HeckeX^2-\zeta(p)$ if $b = 0, p-1$;
  \item\label{item:33} an irreducible subquotient of $\cInd_{KZ}^G\sigma_{a, 0}/f(\HeckeX)$, or equivalently
  of $\cInd_{KZ}^G\sigma_{a, p-1}/f(\HeckeX)$, for some monic irreducible factor of $\HeckeX^2-\zeta(p)$. 
  \end{enumerate}
  \item\label{item:36} The representations of part~\eqref{item:32} are pairwise
non-isomorphic, with the  exception of the isomorphisms
\[\cInd_{KZ}^G \sigma_{a, 0}/f(\HeckeX) \cong \cInd_{KZ}^G\sigma_{a, p-1}/f(\HeckeX)\]
when $f(\HeckeX)$ is coprime to $\HeckeX^2 -\zeta(p)$, %
and
\[\cInd_{KZ}^G\sigma_{a, b}/\HeckeX \cong \cInd_{KZ}^G\sigma_{a+b, p-1-b}/\HeckeX.\]
\item\label{item:37} A representation $\cInd_{KZ}^G\sigma_{a, b}/f(\HeckeX)$ as in~\eqref{item:32} is
  absolutely irreducible if and only if~$f$ has degree one.
\end{enumerate} 
\end{thm}
\begin{proof}
  See~\cite[Theorem~2.12, Corollary~2.14]{DEGlocalization}, which builds upon the results of~\cite{BarthelLivneDuke, BreuilGL2I} (see also~\cite[Section~5.3]{MR3150248}).
\end{proof}

\begin{rem}
  \label{rem:irreps-as-parabolic-induction} By~\eqref{universal parabolic induction}, 
  the representations $\cInd_{KZ}^G\sigma_{a, b}/f(\HeckeX)$ for~$f(\HeckeX) \ne \HeckeX$ 
  occurring in Theorem~\ref{thm:classification-of-irreps} can be written as parabolic inductions from 
the lower-triangular Borel subgroup~$\lbar B$, or equivalently (after a twist) from~$B$.
For example, %
\[\cInd_{KZ}^G\sigma_{a, b}/(\HeckeX-\lambda) \cong \Ind_B^G(\nr_{\lambda^{-1}\zeta(p)}\omega^a \otimes \nr_\lambda\omega^{a+b})\]
for all~$\lambda \in \bF^\times$.
\end{rem}

\subsubsection{Tensoring between \texorpdfstring{$\fC$}{C} and \texorpdfstring{$\cA$}{A}}
If $M$ and $N$ are respectively right and left $\cO\llbracket G\rrbracket _{\zeta}$-modules, 
then we can form their tensor product $M\otimes_{\cO\llbracket G\rrbracket _{\zeta}}N \in \Mod(\cO)$.
There are the usual trifunctorial bijections
\begin{equation}\label{balancedmaps}
\Hom_{\cO}(M \otimes_{\cO\llbracket G\rrbracket _\zeta} N, W) \to \Hom_{\cO\llbracket G\rrbracket _\zeta}\bigl(N, \Hom_\cO(M, W)\bigr)
\end{equation}
and
\begin{equation}\label{balancedmaps bis}
\Hom_{\cO}(M \otimes_{\cO\llbracket G\rrbracket _\zeta} N, W) \to \Hom_{\cO\llbracket G\rrbracket _{\zeta}}\bigl(M, \Hom_\cO(N, W)\bigr )
\end{equation}
given by evaluation maps.
On the right hand side of~\eqref{balancedmaps bis}, contrary to our usual conventions, we are
forming $\Hom$ of right $\cO\llbracket G\rrbracket _{\zeta}$-modules.
Of course, we can always use the anti-involution $g \mapsto g^{-1}$ on $\cO\llbracket G\rrbracket $
to convert $M$ and $\Hom_{\cO}(N,W)$ to left $\cO\llbracket G\rrbracket $-modules,
and then write the target of~\eqref{balancedmaps bis} as $\Hom_{\cO\llbracket G\rrbracket }\bigl(M,\Hom_{\cO}(N,W)\bigr),$
now formed in the category of left $\cO\llbracket G\rrbracket $-modules. 

A key example of the tensor products over~$\cO\llbracket G\rrbracket _{\zeta}$ that
we will consider are those of the form
$M \otimes_{\cO\llbracket G\rrbracket _\zeta} \pi$, 
where $M$ is an object of~$\fC$ (thought of as a right $\cO\llbracket G\rrbracket _{\zeta}$-module;
see Remark~\ref{rem:contragredient again}) and~$\pi$ is an object of~$\cA$.
Since the formation of tensor products is compatible  with filtered colimits,
we will primarily focus on the case when $\pi$ is an object of $\cA^{\fp}$,
and then it will be convenient to also consider the more general
case when $\pi$ is an object of $\Mod^{\fp}(\cO\llbracket G\rrbracket _{\zeta}).$
We are going to describe some general properties of this construction; see~\cite[Section~6.1]{JNWE} for related results.

The functor $M \otimes_{\cO\llbracket G\rrbracket _\zeta} - : \Mod^{\fp}(\cO\llbracket G\rrbracket _\zeta) \to \Mod(\cO)$ is of course a special case of the general construction in Section~\ref{subsubsec:tensor products in abelian categories}, 
viewing~$M$ as a right $\cO\llbracket G\rrbracket _\zeta$-module in the abelian category $\Mod(\cO)$. 
Since $M \in \fC$, we can also regard it as a right $\cO\llbracket G\rrbracket _\zeta$-module in $\Mod_c(\cO)$,
and then Lemma~\ref{lem:abelian-composition-right-exact-tensor} constructs a lift of $M \otimes_{\cO\llbracket G\rrbracket _\zeta} -$ through $\Mod_c(\cO) \to \Mod(\cO)$.
As explained in Remark~\ref{rem:algebra-version-EW}, for all $\pi \in \Mod^{\fp}(\cO\llbracket G\rrbracket _\zeta)$ and~$W \in \Mod_c(\cO)$ there is an isomorphism
\begin{equation}\label{eqn:continuous adjunction}
\Hom_{\Mod_c(\cO)}(M\otimes_{\cO\llbracket G\rrbracket _{\zeta}} \pi, W)
\iso \Hom_{\cO\llbracket G\rrbracket _{\zeta}}\bigl(\pi,\Hom_{\Mod_c(\cO)}(M,W)\bigr).
\end{equation}
Furthermore, by Lemma~\ref{lem:tensor and projective limits}, 
the formation of $M\otimes_{\cO\llbracket G\rrbracket _{\zeta}} \pi$,
for $M$ an object of~$\fC$ and $\pi$ an object of~$\Mod^{\fp}(\cO\llbracket G\rrbracket _{\zeta})$,
is compatible with the formation of cofiltered limits in~$\fC$.
We next describe the interaction of tensor products with the Pontrjagin duality functor~$(\text{--})^\vee$.

\begin{lemma}
\label{lem:tensor facts}
Let $M$ be an object of~$\fC$ and $\pi$ be an object of $\Mod^{\fp}(\cO\llbracket G\rrbracket _{\zeta})$.

\begin{enumerate}
\item\label{item:38} 
There is a natural isomorphism of discrete $\cO$-modules
$$(M\otimes_{\cO\llbracket G\rrbracket _{\zeta}}\pi)^{\vee} \iso \Hom_{\cO\llbracket G\rrbracket _{\zeta}}(\pi,M^{\vee}),$$
or, equivalently, a natural isomorphism of compact $\cO$-modules 
$$M\otimes_{\cO\llbracket G\rrbracket _{\zeta}}\pi \iso \Hom_{\cO\llbracket G\rrbracket _{\zeta}}(\pi,M^{\vee})^{\vee}.$$
\item
If $M$ is coadmissible, 
then $M\otimes_{\cO\llbracket G\rrbracket _{\zeta}} \pi$ is a discrete~$\cO$-module of finite cardinality.
\item
If $\pi$ is of finite length, then there is  
a natural isomorphism of compact $\cO$-modules
$$M\otimes_{\cO\llbracket G\rrbracket _{\zeta}}\pi \iso \Hom_{\fC}(M,\pi^{\vee})^{\vee}$$
{\em (}where $\Hom_{\fC}(M,\pi^{\vee})$ is endowed with its discrete topology,
so its Pontrjagin dual is a compact~$\cO$-module{\em )}.
\end{enumerate}
Similarly, if $\pi \in \Mod^{\fl}(\cO\llbracket KZ\rrbracket _\zeta)$, and~$M \in \Mod_c(\cO\llbracket KZ\rrbracket _\zeta)$, then there is a natural isomorphism of compact $\cO$-modules
\[
M \otimes_{\cO\llbracket KZ\rrbracket _\zeta} \pi \isoto \Hom_{\Mod_c(\cO\llbracket KZ\rrbracket _\zeta)}(M, \pi^\vee)^\vee.
\] 
\end{lemma}
\begin{proof}
Part~(1) follows from applying~\eqref{eqn:continuous adjunction} with $W = \varpi^{-n}\cO/\cO$ and taking the colimit over~$n$ (using the assumption that~$\pi$ is a
compact object of~$\Mod \cO\llbracket G \rrbracket_\zeta$).

For part~(2), note that if $M$ is coadmissible, then its Pontrjagin dual $M^{\vee}$
is an object of $\cA^{\adm}$.  Since furthermore
$\pi$ is finitely generated, we see that
$\Hom_{\cO\llbracket G\rrbracket _{\zeta}}(\pi,M^{\vee})$ is a 
(literally) finite~$\cO$-module.
Thus $M\otimes_{\cO\llbracket G\rrbracket _{\zeta}} \pi$ has a finite Pontrjagin dual, by (1),
and so is itself finite. This proves~(2).

Finally we turn to~(3), and so assume that $\pi$ is finitely presented and of finite length. 
Write $M \iso \varprojlim_i M_i$ as a cofiltered limit of
objects $M_i \in \fC$ of finite length.
Then the~$M_i$ have finite $\cO\llbracket G\rrbracket _\zeta$-length, and so
\begin{multline*}
\Hom_{\cO}^{\cont}(M\otimes_{\cO\llbracket G\rrbracket _{\zeta}} \pi, E/\cO)
 \iso 
\Hom_{\cO}^{\cont}(\varprojlim_i M_i \otimes_{\cO\llbracket G\rrbracket _{\zeta}} \pi, E/\cO)
\\
\iso
\varinjlim_i\Hom_{\cO}(M_i \otimes_{\cO\llbracket G\rrbracket _\zeta} \pi, E/\cO)
\isom
\varinjlim_i\Hom_{\cO\llbracket G\rrbracket }(M_i, \pi^{\vee})
\\
\iso \varinjlim_i \Hom_{\fC}(M_i,\pi^{\vee})
\iso 
\Hom_{\fC}(M,\pi^{\vee}).
\end{multline*}
Here the first isomorphism follows from the compatibility 
with projective limits provided by Lemma~\ref{lem:tensor and projective limits};
the second isomorphism follows from the fact that the Pontrjagin dual of a
projective limit of finite length $\cO$-modules is the direct limit of their Pontrjagin
duals (note that by part~(2), the tensor products $M_i\otimes_{\cO\llbracket G\rrbracket _{\zeta}} \pi$
are finite length $\cO$-modules, since the $M_i$ are coadmissible, being of finite length);
the third isomorphism is provided by~\eqref{balancedmaps bis};
the fourth isomorphism
is given by Lemma~\ref{continuousadmissible};
and the fifth isomorphism follows from the fact that $\pi^{\vee}$ 
is of finite length (so that any homomorphism $M \to \pi^{\vee}$
in $\fC$ factors through some~$M_i$).
This proves~(3), and the same proof works for~$\cO\llbracket KZ\rrbracket _\zeta$.%
\end{proof}

\begin{lemma}
\label{lem:projective Tor vanishing}
If $P$ is a projective object of~$\fC,$
then:
\begin{enumerate}
\item $P$ is a projective object of $\Mod_c(\cO\llbracket KZ\rrbracket _{\zeta}^\op)$ {\em (}and so,
in particular, it is topologically $\cO\llbracket KZ\rrbracket _{\zeta}$-flat, and $\cO\llbracket KZ\rrbracket _{\zeta}$-flat{\em )}.
\item
$\Tor_i^{\cO\llbracket G\rrbracket _{\zeta}}(P,\pi) = 0$ for any object $\pi$ of~$\cA$
and any~$i\geq 1$.
\end{enumerate}
\end{lemma}
\begin{proof}
Since~$P^\vee$ is injective in the category~$\cA^{\ladm}$, it follows from~\cite[Corollary~3.10]{MR2667892} 
that~$P^\vee$ is also injective in the category $(\sm.\,KZ)_\zeta$ of smooth $KZ$-representations with central character~$\zeta$.
Hence~$P$ is projective in the dual category to $(\sm.\,KZ)_\zeta$, which is $\Mod_c(\cO\llbracket KZ\rrbracket _\zeta^\op)$, compare Remark~\ref{rem:contragredient again}.
Lemma~\ref{topologically flat implies flat over Noetherian} now implies that~$P$ is topologically $\cO\llbracket KZ\rrbracket _{\zeta}$-flat and $\cO\llbracket KZ\rrbracket _{\zeta}$-flat.
This establishes the first part of the lemma.

We now prove the second part.
Since~$\Tor_i$ commutes with filtered colimits, we can assume without loss of generality that~$\pi \in \cA^{\fp}$.
Since~$\cA^{\fp}$ is an abelian subcategory of~$\cA$, there exists a resolution 
\[\cInd_{KZ}^G(V_\bullet) \to \pi,\]
where each~$V_i$ is a smooth
$\cO\llbracket KZ\rrbracket _\zeta$-module of finite $\cO$-length. %
By the first part of the lemma, $\cInd_{KZ}^GV_i$ is acyclic for the functor $P \otimes_{\cO\llbracket G\rrbracket _{\zeta}} (\text{--})$, since
\[
P \otimes^L_{\cO\llbracket G\rrbracket _\zeta} \cInd_{KZ}^GV_i
\cong
P \otimes^L_{\cO\llbracket G\rrbracket _\zeta}
(\cO\llbracket G\rrbracket _{\zeta}\otimes_{\cO\llbracket KZ\rrbracket _{\zeta}}V_i)
\cong P \otimes^L_{\cO\llbracket KZ\rrbracket _\zeta} V_i.
\]
Hence $\Tor_i^{\cO\llbracket G\rrbracket _\zeta}(P, \pi)$ is the homology of 
\[P \otimes_{\cO\llbracket G\rrbracket _\zeta} \cInd_{KZ}^GV_\bullet,\]
which by Lemma~\ref{lem:tensor facts}~(1) is isomorphic to the homology of
\[\Hom_{\cO\llbracket G\rrbracket _\zeta}(\cInd_{KZ}^GV_\bullet, P^\vee)^\vee.\]
This is concentrated in degree zero as a consequence of~\cite[Corollary~5.18]{MR3150248}, which implies that~$P^\vee$ is an injective object of~$\cA$.
\end{proof}

\subsection{Morita theory for blocks \texorpdfstring{of $\cA^{\ladm}$ and $\fC$}{}}\label{VP Morita theory}
Unlike $\cA$, the category $\cA^{\ladm}$ 
is locally finite, by Lemma~\ref{properties of A^ladm}.
By the structure theory of locally finite categories, recalled in Section~\ref{subsec:locally finite categories}, $\cA^{\ladm}$ admits a decomposition into blocks.
An extensive analysis of the various blocks is made in~\cite{MR3150248}, 
and in what follows we recall, and slightly expand, some of these results.
The arguments in~\cite{MR3150248} rely crucially on Colmez's functor ``$V$'',
but we postpone our discussion of this functor until the following section;
in this section we focus on those results whose statements can be made without
reference to Colmez's functor. 

\subsubsection{Classification of blocks}
Recall that,
by definition, a block of~$\cA^{ \ladm}$
is an equivalence class of (isomorphism classes of) irreducible objects\footnote{We note that  it
follows from the results of \cite{BarthelLivneDuke, BreuilGL2I}
that the irreducible objects of $\cA$ are automatically
admissible, and hence lie in $\cA^{\ladm}.$  Thus
we can equally well regard this as an equivalence relation
on the irreducible objects of~$\cA$.}
 under the equivalence relation generated by
\[
\pi_1 \sim \pi_2 \text{ if } \Ext^1_{\cA}(\pi_1, \pi_2) \ne 0 \text{ or } \Ext^1_{\cA}(\pi_2, \pi_1) \ne 0.
\]
The blocks~$\mathfrak{B}$ containing absolutely irreducible objects were classified in~\cite{MR3150248}.
The remaining blocks were classified in~\cite[Proposition~2.4.8, Remark~2.4.9]{DEGlocalization}.
In summary, blocks come in five ``types'' as follows (we
explain the labels in Remark~\ref{rem: labelling of blocks} below):
\begin{enumerate}
\item[{\crtcrossreflabel{(ssg)}[item: ss block]}] $\fB =\{\pi\}$ where $\pi = \cInd_{KZ}^G \sigma/\HeckeT$, for some $\zeta$-compatible Serre weight~$\sigma$,
is an irreducible supersingular representation.
\item[{\crtcrossreflabel{(gen)}[item: generic block]}] $\fB = \{\Ind_B^G(\chi_1 \otimes \omega^{-1}\chi_2), \Ind_B^G(\chi_2 \otimes \chi_1\omega^{-1})\}$ 
for characters~$\chi_1, \chi_2 : \bQ_p^\times \to \F^\times$ such that $\chi_1\chi_2 = \lbar \zeta\omega$ and $\chi_1\chi_2^{-1} \ne 1, \omega^{\pm 1}$.
\item[{\crtcrossreflabel{(scalar)}[item: non p-distinguished block]}] $\fB= \{\Ind_B^G(\chi \otimes \omega^{-1} \chi)\}$ for a character $\chi: \bQ_p^\times \to \F^\times$
such that $\chi^2 = \lbar \zeta \omega$.
\item[{\crtcrossreflabel{(St)}[item: Steinberg block]}] $\fB=\{\chi, \chi \otimes \St, \Ind_B^G(\omega\chi \otimes \omega^{-1} \chi)\}$ for a character $\chi : \bQ_p^\times \to \F^\times$
such that $\chi^2 = \zetabar$.
\item[{\crtcrossreflabel{(gen+)}[item: not abs irred block]}] 
$\fB$ is the set of irreducible subquotients of 
\[
\cInd_{KZ}^G \sigma/f_\fB(\HeckeT) \oplus \cInd_{KZ}^G\sigmacomp/f_\fB^*(\HeckeT)
\]
for some companion pair~$\sigmasigmacomp$ of $\zeta$-compatible Serre weights,
and some irreducible monic polynomial
$f_\fB \in \F[\HeckeT]$ of degree~$>1$.
Here we have written
$f_\fB^*(\HeckeT) \coloneqq  f_\fB(0)^{-1}\HeckeT^{\deg f_\fB}f_\fB(\zeta(p)/\HeckeT)$
for the irreducible monic polynomials whose roots are $\lambda^{-1}\zeta(p)$, where~$\lambda$ runs through the roots of~$f_\fB$.
\end{enumerate}

\begin{rem}
\label{rem:comparing gen and gen+}
Note that blocks of type~\ref{item: generic block}
can also be written in the form
$\{\cInd_{KZ}^G \sigma / (\HeckeT-\lambda), \cInd_{KZ}^G \sigmacomp/(\HeckeT-\lambda^{-1}\zeta(p))\}$
for an appropriately chosen companion pair of weights~$\sigmasigmacomp$ %
and $\lambda \in \F^{\times}$.
Thus the key point in the description of blocks of type~\ref{item: not abs irred block}
is that the polynomial $f_\fB$ be of degree strictly greater than~$1$.
\end{rem}

\subsubsection{The full subcategory associated to a block}
If~$\fB$ is a block of irreducible objects of~$\cA$, the discussion of Section~\ref{subsec:locally finite categories} yields a direct factor~$\cA_\fB$ of~$\cA$.
By Lemma~\ref{lem:fl equals Noetherian}, the category~$\cA_\fB$ is locally finite, and its subcategory $\cA_{\fB}^{\fl}$ of objects of finite length coincides with its
subcategory of compact objects, resp.\ Noetherian objects.
These coincide furthermore with the finitely generated objects, by the following lemma.

\begin{lem}\label{fg equals fl for Athetabar}
The category~$\cA_{\fB}^{\fl}$ coincides with the full subcategory of finitely presented objects of~$\cA_{\fB}$.
\end{lem}
\begin{proof}
Lemma~\ref{lem:fl equals Noetherian} implies that~$\cA_\fB^{\fl}$ is the subcategory of Noetherian objects of~$\cA_{\fB}$.
On the other hand, an object of~$\cA_{\fB}$ is Noetherian in~$\cA_{\fB}$ if and only if it is Noetherian in~$\cA$, and the Noetherian objects of~$\cA$ are precisely the finitely
generated objects, by~\cite[Corollary~2.2.4]{DEGlocalization}. We conclude because finitely generated objects are finitely presented (see Remark~\ref{fg equals fp for cA}). %
\end{proof}

\subsubsection{Blocks of type~{\em\ref{item: not abs irred block}}}
If~$\fB$ is a block of type~\ref{item: not abs irred block}, by~\cite[Prop.\ 2.4.8, Rem.\ 2.4.9]{DEGlocalization} %
we have the following (mutually exclusive and exhaustive) possibilities:
\begin{itemize}
  \item $\sigma = \sigmacomp = \sigma_{a, p-2}$ and $f_\fB = f_\fB^*$. Then~$\fB$ is a singleton.
  \item $\{\sigma, \sigmacomp\} = \{\sigma_{a, 0}, \sigma_{a+1, p-3}\}$, and $f_\fB = \HeckeT^2-\zeta(p)$ (hence $\zeta(p)$ is not a square in $\cO^\times$). Then~$\fB$ has three elements.
  \item $\fB$ has two elements.
\end{itemize} 

If~$\F'$ is a splitting field of~$f_\fB$ over~$\F$,
then the set
\[
\fB \otimes_\F \F'\coloneqq \{\JH(\pi \otimes_\F \F'): \pi \in \fB\}
\]
is a union of blocks of $\cA_{\F'}^{\ladm}$.
The group $\Gal(\F'/\F)$ acts transitively on~$\fB \otimes_\F \F'$ and preserves its partition as a union of blocks; 
the stabilizer of each part is trivial, except in the following case:
\begin{itemize}
  \item $\sigma = \sigmacomp = \sigma_{a, p-2}$, $f_\fB = f_\fB^*$, and $f_\fB \ne \HeckeT^2-\zeta(p)$.
\end{itemize}
In fact, in this case, the element of order two of $\Gal(\F'/\F)$ sends each root~$\lambda$ of~$f_\fB$ to $\lambda^{-1}\zeta(p) \ne \lambda$, and so it fixes each part.
(This exception can be explained in Galois-theoretic terms: for example,
when $\deg f_\fB = 2$, $f_\fB = f_\fB^*$, and $f_\fB \ne \HeckeT^2-\zeta(p)$, the block $\{\cInd_{KZ}^G(\sigma_{a, p-2})/f_\fB\}$ will correspond under Definition~\ref{defn: pseudorep labelled by block}
to an irreducible but not absolutely irreducible Galois pseudorepresentation. 
Note that these are all twists of unramified representations, hence have Serre weight~$\sigma_{a, p-2}$, as explained in Remark~\ref{F-rational}.)

In keeping with Remark~\ref{F-rational}, for any block of~$\cA$ (of type~\ref{item: not abs irred block} or not) 
we write~$\bF_\fB$ for the splitting field of~$f_\fB$, except in the case just described, in which we write $\bF_\fB$ for 
the index-two subextension of the splitting field (so that, in this case, $f_\fB$ splits in~$\bF_\fB[\HeckeX]$ as a product of irreducible quadratic polynomials,
each fixed by the involution $f \mapsto f^*$).

\begin{defn}\label{defn:F-rational block}
We say that a block~$\fB$ of~$\cA^{\ladm}$ is $\F$-rational if $\F_\fB = \F$.
\end{defn}

The next result shows that every block~$\fB$ of~$\cA$ is equivalent to an $\bF_\fB$-rational block of $\cA_{\bF_\fB}$.

\begin{lemma}\label{reducing gen+ to gen}
Let~$\fB$ be a block of type~{\em\ref{item: not abs irred block}},
and let $\fB'$ be a block of $\cA_{\F_\fB}$ contained in~$\fB \otimes_\F \F_\fB$.
Write~$\cO'/\cO$ for the unramified extension with residue field~$\F_\fB$.
Then $\fB'$ is $\bF_\fB$-rational, and restriction of scalars from~$\cO'$ to~$\cO$ is an equivalence of categories $\cA_{\cO', \fB'} \to \cA_{\cO, \fB}$.
\end{lemma}
\begin{proof}

By definition, $\F_\fB$ is a splitting field of $f_\fB$, except when $\sigma = \sigma_{a, p-2}, f_\fB = f_\fB^*$, and $f_\fB \ne \HeckeT^2-\zeta(p)$.
In this case, $f_\fB$ splits as a product of irreducible quadratic polynomial, which are fixed by the involution $f \mapsto f^*$.
Thus we have the following mutually exclusive possibilities:
\begin{itemize}
  \item $\sigma = \sigma_{a, p-2}$, and $f_\fB = \HeckeT^2-\zeta(p)$.
  Then $\fB'$ has type~\ref{item: non p-distinguished block}, and it has the form $\{\pi_0\}$ for some absolutely irreducible object~$\pi_0$ of $\cA_{\cO'}$.
  \item $\sigma = \sigma_{a, p-2}$, $f_\fB = f_\fB^*$, and $f_\fB \ne \HeckeT^2-\zeta(p)$.
  Then $\fB'$ is an $\F_\fB$-rational block of type~\ref{item: not abs irred block}, and it has the form $\{\pi_0\}$ for some irreducible object~$\pi_0$ of $\cA_{\cO'}$.
  \item $\sigma = \sigma_{a, 0}$, and $f_\fB = \HeckeT^2-\zeta(p)$.
  Then $\fB'$ has type~\ref{item: Steinberg block}, and it has the form $\{\pi_0, \pi_1, \pi_2\}$ for some absolutely irreducible objects~$\pi_0, \pi_1, \pi_2$ of $\cA_{\cO'}$.
  \item $\fB'$ %
has type~\ref{item: generic block},
and it has the form $\{\pi_0, \pi_1\}$ for some absolutely irreducible objects~$\pi_0, \pi_1$ of $\cA_{\cO'}$.
\end{itemize} 

In all of these cases, $\fB'$ is $\F_\fB$-rational, and 
for all $\alpha\ne \beta \in \Gal(\cO'/\cO)$, we have $\alpha^*\fB' \ne \beta^* \fB'$.
Furthermore, restriction of scalars preserves 
local admissibility, by Lemma~\ref{properties of A^ladm}~(3),
and so we may also regard each $\pi_i$ as a locally admissible object of~$\cA_{\cO}$.

We now prove that~$\pi_i$ is an irreducible object of~$\cA_{\cO, \fB}$.
Choose an irreducible subobject $\pi \to \pi_i$
of~$\pi_i$ in~$\cA_\cO^{\ladm}$.
Then~$\pi$ is an irreducible subquotient of $\cInd_{KZ}^G\sigma'/f'(\HeckeT)$ for some irreducible monic~$f' \in \F[\HeckeT]$.
Since $\pi \otimes_\cO \cO'$ is a subobject of $\pi_i \otimes_\cO \cO' = \oplus_{\gamma \in \Gal(\cO'/\cO)} \gamma^* \pi_i$, we see that
$\cInd_{KZ}^G\sigma'/f'(\HeckeT)$ and 
\[
\cInd_{KZ}^G \sigma/f_\fB(\HeckeT) \oplus \cInd_{KZ}^G\sigmacomp/f_\fB^*(\HeckeT)
\]
have an irreducible subquotient in common after extending scalars to~$\cO'$.
By~\cite[Corollary~2.1.14]{DEGlocalization} we deduce that $(\sigma', f') \in \{(\sigma, f_\fB), (\sigmacomp, f_\fB^*)\}$, and so~$\pi$ is an irreducible object 
of~$\cA_{\cO, \fB}$.
Now, by construction, we know that $\pi \otimes_\cO \cO'$ has $\cA_{\cO'}$-length equal to $[\cO':\cO]$.
This is the same as the length of $\pi_i \otimes_\cO \cO'$, and so the inclusion $\pi \to \pi_i$ is an isomorphism, as desired.

It follows from the previous paragraph that restriction of scalars defines a functor $\cA_{\cO', \fB'} \to \cA_{\cO, \fB}$, and there remains to prove that it is an equivalence.
We will do this by an application of Morita theory.
Let~$\cI$ be an injective object of~$\cA_{\cO', \fB'}$ with 
socle equal to the direct sum of the objects of~$\fB'$.
Then $\cI$ is an injective object of~$\cA_{\cO, \fB}$:
in fact, $\Hom_{\cA_{\cO}^{\ladm}}(\text{--}, \cI) \otimes_\cO \cO'$ is exact, because it 
is naturally isomorphic to $\Hom_{\cA_{\cO'}^{\ladm}}(\text{--} \otimes_\cO \cO', \cI \otimes_\cO \cO')$ 
(by Lemma~\ref{properties of A^ladm}~(3))
and $\cI \otimes_\cO \cO' = \oplus_{\gamma \in \Gal(\cO'/\cO)} \gamma^* \cI$
is injective.
Furthermore, $\cI$ is an injective cogenerator of~$\cA_{\cO, \fB}$: this is a consequence of the fact that~$\fB$ and~$\fB'$ have the same number of elements, and
non-isomorphic irreducible objects~$\pi_i, \pi_j$ of~$\fB'$ are not isomorphic in $\cA_{\cO}$
(as can be seen by computing $\Hom_{\cA_\cO}(\pi_i, \pi_j) \otimes_\cO \cO' = \Hom_{\cA_{\cO'}}(\pi_i \otimes_\cO \cO', \pi_j \otimes_\cO \cO') = 0$, using that 
$\pi_i \ne \gamma^* \pi_j$ for any~$\gamma \in \Gal(\cO'/\cO)$).
By Morita theory (see e.g.\ Section~\ref{subsubsec:theta Morita})
there remains to prove that the restriction map
\begin{equation}\label{to prove isomorphism for gen+}
\End_{\cA_{\cO'}^{\ladm}}(\cI) \to \End_{\cA_\cO^{\ladm}}(\cI)
\end{equation}
is an isomorphism, or equivalently, is surjective.
We can verify this after applying~$\text{--}\otimes_\cO \cO'$.
Now
\[
\End_{\cA_\cO^{\ladm}}(\cI) \otimes_\cO \cO' = \bigoplus_{\alpha, \beta \in \Gal(\cO'/\cO)}\Hom_{\cA^{\ladm}_{\cO'}}(\alpha^*\cI, \beta^* \cI) = 
\bigoplus_{\gamma \in \Gal(\cO'/\cO)} \End_{\cA_{\cO'}^{\ladm}}(\gamma^*\cI)
\]
since~$\alpha^*\cI$ and~$\beta^* \cI$ are in different blocks of~$\cA_{\cO'}^{\ladm}$ if~$\alpha \ne \beta$,
and $\eqref{to prove isomorphism for gen+} \otimes_{\cO} \cO'$ is the diagonal inclusion.
Since the image of $\eqref{to prove isomorphism for gen+} \otimes_{\cO} \cO'$ also contains the idempotents in the target (because it contains $\cO' \otimes_{\cO} \cO'$), 
this concludes the proof.
\end{proof}

Because of Lemma~\ref{reducing gen+ to gen}, one can often assume without loss of generality that~$\fB$ is $\bF$-rational.
As already mentioned, these blocks have been thoroughly studied in~\cite{MR3150248}, with the exception of $\bF$-rational blocks~$\fB$ of type~\ref{item: not abs irred block};
note that these are not equivalent to blocks of type~\ref{item: generic block},
e.g.\ because they contain a unique irreducible object.
However, many of the invariants of~$\fB$ that we define in the following sections become isomorphic to invariants of blocks of type~\ref{item: generic block} after a base extension, which will be enough
for our intended applications.

\subsubsection{The dual category associated to a block}
We write~$\fC_\fB$ for the image of~$\cA_\fB$ under Pontrjagin duality.
By Proposition~\ref{Gabriel decomposition}, there exists a pseudocompact $\cO$-algebra~$E_\fB$ such that $\fC_\fB$ is equivalent to the category of pseudocompact right $E_\fB$-modules.
The ring~$E_{\fB}$ is only well-defined up to its category of pseudocompact modules,
and one of the main results of Pa\v{s}k\={u}nas in~\cite{MR3150248} is the computation
of an explicit choice of $E_{\fB}$, and of its centre.
We will recall these results in detail in the next section. For now,
we content ourselves by noting
the following particular lemma, which follows directly from those more
precise results,
and then deriving some consequences of it.
We do this because~\cite{MR3150248}
only considers blocks that contain absolutely irreducible representations (i.e.\ not of type~\ref{item: not abs irred block})
whereas it will be important for us to have a uniform statement for all blocks.
This will recur throughout our discussion
of ``classical'' $p$-adic local Langlands for~$\GL_2(\Q_p)$
(for example, 
the absolutely irreducible case of Lemma~\ref{lem:projective generators are flat over their endos}~\eqref{item:projective 4} 
can be found in~\cite[Corollary~6.7]{MR4350140}), 
and as we already intimated in Section~\ref{subsubsec:other works},
these amplifications of the existing literature are straightforward consequences
of those existing results. 

\begin{lemma}\label{finite cosocle implies Noetherian endomorphisms}
Let~$M$ be an object of~$\fC_{\fB}$ with cosocle of finite length.
Then $\End_{\fC_\fB}(M)$ is finitely generated over the Bernstein centre~$\cZ_\fB$ of~$\fC_\fB$, 
which is a Noetherian, local, profinite $\cO$-algebra.
Hence~$\End_{\fC_\fB}(M)$ is a Noetherian profinite $\cO$-algebra.
\end{lemma}
\begin{proof}
It suffices to construct a projective generator~$P_\fB$ 
of~$\fC_\fB$, with cosocle of finite length, such that 
$E_\fB\coloneqq  \End_{\fC_\fB}(P_\fB)$ is finitely generated over its centre~$Z(E_\fB)$, and~$Z(E_\fB)$ is Noetherian and local.
In fact, the equivalence~\eqref{equivalence from Gabriel's thesis} implies that $\cZ_\fB \cong Z(E_\fB)$, and so~$\cZ_\fB$ is Noetherian, local and profinite.
Finally, if~$M$ has cosocle of finite length, then $\End_{\fC_\fB}(M)$ is an $E_\fB$-submodule of $\Hom_{\fC_\fB}(P_\fB^{\oplus n}, M)$ for some~$n \geq 0$, which is a quotient of
$\Hom_{\fC_\fB}(P_\fB^{\oplus n}, P_\fB^{\oplus n}) \cong E_\fB^{\oplus n^2}$.
Hence $\End_{\fC_\fB}(M)$ is a subquotient of a finitely generated $\cZ_\fB$-module, and so it is finitely generated over~$\cZ_\fB$.

If~$\fB$ contains absolutely irreducible representations, then the existence of~$P_\fB$ with these properties is established in~\cite{MR3150248}.
Hence the lemma is true when~$\fB$ contains absolutely irreducible representations.
By Lemma~\ref{reducing gen+ to gen}, there remains to prove the lemma in the case that~$\fB = \{\pi_0\}$ is $\bF$-rational of type~\ref{item: not abs irred block}.
If~$P_\fB \to \pi_0^\vee$ is a projective envelope, and $\cO \to \cO'$ is unramified quadratic, then $P_\fB \otimes_\cO \cO'$ is a projective generator of a block of $\fC_{\cO'}$ 
of type~\ref{item: generic block}.
Hence $\End(P_\fB \otimes_\cO \cO')$ is finitely generated over its centre, which is Noetherian and local. 
Since $\End(P_\fB \otimes_\cO \cO') = \End(P_\fB) \otimes_\cO \cO'$ and $Z(\End(P_\fB) \otimes_\cO \cO') = Z(\End(P_\fB)) \otimes_\cO \cO'$, we deduce that 
$Z(\End(P_\fB))$ is Noetherian and local, and that $\End(P_\fB)$ is finitely generated over $Z(\End(P_\fB))$, as desired.
\end{proof}

\begin{lemma}
\label{lem:projective generators are flat over their endos}
Let~$P$ be a projective object of~$\fC_{\fB}$ with cosocle of finite length.
\begin{enumerate}
\item\label{item:projective 1} The endomorphism ring $E\coloneqq  \End_{\fC_{\fB}}(P)$ is a Noetherian profinite $\cO$-algebra, and $P$ is a Noetherian object of $\fC$.
\item\label{item:projective 2} The module~$P$, with its natural topology as an object of~$\fC_{\fB}$, is an object of~$\Mod_c(E)$.
\item\label{item:projective 4} For all~$i > 0$, the quotient $P/\rad(E)^i P$ is a finite length, hence coadmissible, object of~$\fC_\fB$.
\item\label{item:projective 3} If~$P$ is furthermore a projective generator of~$\fC_{\fB}$, then 
$P$ is topologically flat over~$E$ {\em (}hence projective in~$\Mod_c(E)$, by Lemma~{\em\ref{lem: topologically flat implies projective}}{\em )}
and~$E$ has finite global dimension. %
\end{enumerate}
\end{lemma}
\begin{proof}
Proof of~\eqref{item:projective 1}: An application of Lemma~\ref{finite cosocle implies Noetherian endomorphisms} shows that~$E$ is a (left and right) Noetherian profinite $\cO$-algebra.
To see that~$P$ is Noetherian in $\fC$, or equivalently in~$\fC_\fB$, note that there exists a projective generator $P_\fB$ of~$\fC_\fB$, with finite cosocle and Noetherian endomorphism ring~$E_\fB$, 
such that~$P$ is a quotient of $P_\fB$.
Under the equivalence $\Hom_{\fC_\fB}(P_\fB, \text{--}) : \fC_\fB \to \Mod_c(E_\fB^\op)$, $P$ goes to a quotient of~$E_\fB$, which is a Noetherian $E_\fB^\op$-module.
Hence~$P$ is Noetherian in~$\fC_\fB$.

Proof of~\eqref{item:projective 2}:
Since the natural topology on~$P$ is profinite,
it suffices to prove that every open $\cO$-submodule of~$P$ of finite index
in~$P$ contains an open left $E$-submodule, which is exactly what is proved in~\cite[Lemma~2.7]{MR3150248}.
(Strictly speaking, this reference works under the additional assumption that $\cosoc(P)$ is multiplicity free, but this assumption is only used to ensure that $E/\Rad(E)$
has finite $\bF$-dimension, which follows from the fact that $E$ is compact and has the $\rad(E)$-adic topology, by Lemma~\ref{compact modules}~\eqref{item: compact 7}).

Proof of~\eqref{item:projective 4}: By Lemma~\ref{finite cosocle implies Noetherian endomorphisms}, $E$ is finitely generated over the Bernstein centre $\cZ_\fB$, 
which is a Noetherian local ring.
Writing $\fm \coloneqq  \rad(\cZ_\fB)$, 
we deduce that~$E$ has the $\fm$-adic topology, and so
the sequences $\rad(E)^i$ and~$\fm^i E$ are cofinal in~$E$. 
It thus suffices to prove that $P/\fm^i P$ has finite length.
Choose a projective generator~$P_\fB$ of~$\fC_\fB$ with finite cosocle, and write~$E_\fB$ for its endomorphism algebra.
Since~$P$ has finite cosocle, there is a surjection $P_\fB^{\oplus n} \to P/\fm^i P$ for some~$n$, and so
$\Hom_{\fC_\fB}(P_\fB, P/\fm^i P)$ is an $E^{\op}_\fB$-quotient of~$E_\fB^{\oplus n}$.
Hence $\Hom_{\fC_\fB}(P_\fB, P/\fm^i P)$ is a finite~$E_\fB/\fm^i E_\fB$-module.
Since $\fm^i E_\fB$ is open in~$E_\fB$, 
and~$E_\fB$ is profinite, the quotient $E_\fB/\fm^i E_\fB$ is a finite set, and so we conclude that
$\Hom_{\fC_\fB}(P_\fB, P/\fm^i P)$ has finite $E_\fB^{\op}$-length. 
Since $\Hom_{\fC_\fB}(P_\fB, \text{--}): \fC_\fB \to \Mod_c(E_\fB^\op)$ is an equivalence, we conclude that~$P/\fm^i P$ has finite $\fC_\fB$-length, as desired.

Proof of~\eqref{item:projective 3}: Assume that~$P$ is a projective generator of~$\fC_{\fB}$.
Then, by Proposition~\ref{compact Morita theory}, $P$ is a complete left $E$-module in~$\fC_\fB$, and the corresponding functor
\[
\text{--}\cotimes_E P : \Mod_c(E^\op) \to \fC_\fB
\]
is an equivalence.
Since the forgetful functor $i: \fC_\fB \to \Mod_c(\cO)$ is exact and cofiltered limit-preserving, 
the composition $i \circ (\text{--} \cotimes_E P)$ is an exact and cofiltered limit-preserving functor $\Mod_c(E^\op) \to \Mod_c(\cO)$.
By Lemma~\ref{lem:abelian-composition-right-exact-tensor}, this composition is naturally isomorphic to the completed tensor product
associated to~$i(P)$.
By Lemma~\ref{compact modules are compact}, this is the usual completed tensor product~$\text{--}\cotimes_E P$, which is therefore exact, as desired.

Finally, the claim that~$E$ has finite global dimension can be verified as follows.
By~\cite[Proposition~4.2.2]{JNWE} and Lemma~\ref{finite cosocle implies Noetherian endomorphisms}, it suffices to prove that
every simple left $E^{\op}$-module~$M$ has a finite projective resolution (since there are finitely many isomorphism classes of simple modules).
For this, it suffices to prove that~$M$ has a finite projective resolution in~$\Mod_c(E^{\op})$, all of whose terms have finite length $E^{\op}$-cosocle, or equivalently,
that $M$ has a finite injective resolution in $\Mod_c(E^{\op})^{\op}$, all of whose terms have socle of finite length.
Since~$\Mod_c(E^{\op})^{\op}$ is a locally finite category, general theory shows that~$M$ has an injective resolution~$M \to J^\bullet$ such that, for all simple~$N \in \Mod_c(E^{\op})$, we have
\[
\operatorname{length}_\cO \Hom_{\Mod_c(E^{\op})^{\op}}(N, J^t) = \operatorname{length}_\cO \Ext^t_{\Mod_c(E^{\op})^{\op}}(N, M)
\]
(see the explanation in~\cite[Remark~10.11]{MR3150248}).
It thus suffices to prove that for any two simple modules~$M, N$, the Ext-group 
$\Ext^t_{\Mod_c(E^{\op})^{\op}}(N, M)$ is $\cO$-finite, and vanishes for large enough~$t$.
Since $\Mod_c(E^{\op})^{\op} \isoto \cA_{\fB}$, this is a consequence of known properties of $\Ext$-groups between irreducible objects of~$\cA_{\fB}$.
\end{proof}

In Section~\ref{socle and radical} we have defined socle and radical filtrations for objects of a complete and cocomplete abelian category,
such as $\cA$, $\fC$, or~$\Mod(\cO\llbracket G\rrbracket _\zeta)$.
The next lemma shows that, for finite length objects, the natural inclusions $\cA \to \Mod(\cO\llbracket G\rrbracket _\zeta)$ and $\fC \to \Mod(\cO\llbracket G\rrbracket _\zeta)$ preserve these filtrations.

\begin{lemma}
If~$M \in \fC$ has finite length, then the $\fC$-socle filtration and the $\Mod(\cO\llbracket G\rrbracket _\zeta)$-socle filtration of~$M$ coincide.
The same statement is true with~$\fC$ replaced by~$\cA$, or ``socle'' replaced by ``radical''. 
\end{lemma}
\begin{proof}
This is an immediate consequence of Lemma~\ref{lem:fC to O-G- },
resp.\ Lemma~\ref{embedding of smooth in modules}.
\end{proof}

We now prove two finiteness results for the radical filtration %
of objects of~$\fC$, or equivalently, the socle filtration of objects of~$\cA$.

\begin{lemma}
\label{lem:Loewy facts}
If $M$ is an object of $\fC$ with cosocle of finite length, 
then $M$ is Noetherian in~$\fC$, and $M/\operatorname{rad}_n M$ is of finite length for each~$n\geq 0$.
Furthermore, $M \iso \varprojlim_n (M/\operatorname{rad}_n M).$
\end{lemma}
\begin{proof}
If~$\cosoc(M)$ has finite length, then~$M$ is a quotient of a projective envelope~$P$ of~$\cosoc(M)$, and~$P$ is Noetherian by 
Lemma~\ref{lem:projective generators are flat over their endos}~\eqref{item:projective 1}.
Hence~$M$ is Noetherian.
If~$\rad^i M / \rad^{i+1} M$ has infinite length for some~$i$, then (by the dual to Corollary~\ref{associated graded of socle filtration is semisimple}~(1)) 
it is an infinite direct product of simple objects, and so it contains an infinite ascending chain
of subobjects.
This contradicts the fact that~$M$ is Noetherian and concludes the proof of the first statement of the lemma.

We now prove that $M \iso \varprojlim_n M/\rad^n M$.
Since~$\fC$ is dual to~$\cA^{\ladm}$, we know that 
\[
M \iso \varprojlim_{M' \subseteq M} M/ M'
\]
where the limit is over the set of closed submodules~$M' \subseteq M$ of finite $\fC$-colength.
Since we have just proved that~$\rad^n M$ has finite colength for all~$n$, it suffices to prove that if
$q: M \to N$ is a surjection, and~$N$ has finite length, then there exists~$n \geq 0$ such $\rad^n M \subseteq \ker(q)$. 
However, $\rad^n N = 0$ for~$n$ large enough, and then Lemma~\ref{general Loewy facts}(2) implies
that $\rad^n M \subseteq \ker(q)$, as desired.
This concludes the proof.
\end{proof}

\begin{lemma}\label{bounded quotients of finite length}%
Let~$\pi \in \cA^{\fp}$, let $\fB$ be a block of~$\cA^{\ladm}$, and let~$n \geq 0$.
Then there exists a unique subobject~$\pi^{(n)}_{\fB} \subseteq \pi$ such that
\begin{enumerate}
\item $\pi/\pi^{(n)}_{\fB} \in \cA^{\fp}_{\fB}$ and has Loewy length~$\leq n$.
\item if $\pi' \subseteq \pi$, and $\pi/\pi' \in \cA^{\fp}_{\fB}$ and has Loewy length~$\leq n$, then $\pi^{(n)}_{\fB} \subseteq \pi'$. 
\end{enumerate}
The set $\{\pi^{(n)}_{\fB}: n \geq 0\}$ is
cofinal in the set of quotients of~$\pi$ contained in~$\cA^{\fp}_{\fB}$.
\end{lemma}
\begin{proof}
The uniqueness part is immediate.
For the existence part, let~$\cI$ be a direct sum of injective envelopes of the simple objects of~$\cA_{\fB}$.
The dual of Lemma~\ref{lem:Loewy facts} implies that~$\soc_{\cA, n} \cI$ has finite $\cA$-length.
Hence $\Hom_\cA(\pi, \soc_{\cA,n}\cI)$ is finitely generated over~$\cO$. 
Choose generators~$\varphi_1, \ldots, \varphi_r$ of this $\cO$-module, and let
\[
\pi^{(n)}_{\fB} = \bigcap_{i} \ker(\varphi_i).
\]
Then
\[\pi/\pi^{(n)}_{\fB} \subseteq \bigoplus_i \soc_{\cA, n} \cI,\]
and the right-hand side is an object of~$\cA_{\fB}^{\fp}$ and has Loewy length $\leq n$.
Hence~(1) holds, by Lemma~\ref{general Loewy facts}.
To see that~(2) holds, we need to prove that $\pi^{(n)}_{\fB}$ maps to zero in~$\pi/\pi'$. Since
$\cI$ is an injective cogenerator of~$\cA_{\fB}$,
 it suffices to prove that for all
$\varphi \in \Hom_\cA(\pi/\pi', \cI)$, the composite $\pi^{(n)}_{\fB}\to\pi \to \pi/\pi' \xrightarrow{\varphi} \cI$ is zero. %
Now Lemma~\ref{general Loewy facts} implies that~$\varphi$ factors through $\soc_{\cA, n} \cI$, and all elements of $\Hom_\cA(\pi, \soc_{\cA, n} \cI)$ are zero on
$\pi^{(n)}_{\fB}$, by construction.
This concludes the proof of the existence of~$\pi^{(n)}_{\fB}$.

Finally, the cofinality claim follows from the second property of~$\pi^{(n)}_{\fB}$, since any quotient of~$\pi$ contained in $\cA^{\fp}_{\fB}$ has finite $\cA$-length, 
hence finite Loewy length.
\end{proof}

We now consider the case of a projective generator with finite cosocle.

\begin{lemma}
\label{lem:recapitulated Morita}
Let~$\fB$ be a block of~$\cA^{\ladm}$, and let~$P_{\fB}$ be a projective generator of~$\fC_{\fB}$ with cosocle of finite length.
Let~$E_\fB \coloneqq \End_{\fC_\fB}(P_\fB)$.
Then:
\begin{enumerate}
\item
the functor $M \mapsto \Hom_{\fC_\fB}(P_{\fB}, M) = \Hom_{G}^{\cont}(P_{\fB},M)$
gives an equivalence
$$\fC_{\fB} \iso \Mod_c(E_{\fB}^{\op}),$$
with quasi-inverse given by
$N \mapsto N \cotimes_{E_{\fB}} P_{\fB}.$
\item\label{item:63} 
$\pi \mapsto P_{\fB} \otimes_{\cO\llbracket G\rrbracket _{\zeta}} \pi$ gives an equivalence
$$\cA_{\fB}^{\fp} \iso \Mod^{\fl}(E_{\fB}),$$
with quasi-inverse given by
$M \mapsto \Hom_{E_{\fB}}^{\cont}(P_{\fB}, M)$.
\item If~$M$ is an object of~$\fC_{\fB}$ with cosocle of finite length, then the module 
$\Hom_{\fC_{\fB}}(P_{\fB}, M)$ is finitely presented over $E_{\fB}^\op$ and over~$\End_{\fC_\fB}(M)$.
\item If~$M$ is an object of $\fC_\fB$ with cosocle of finite length, then there is a natural isomorphism
\begin{equation}
  \label{eq:Pthetabar-tensor-pi-formula}
  M%
  \iso
  \Hom_{\fC_{\fB}}(P_{\fB},M)
  \otimes_{E_{\fB}} P_{\fB}.
\end{equation}
\end{enumerate}
\end{lemma}
\begin{proof}
Part~(1) is Proposition~\ref{compact Morita theory}, whose assumptions are met because of Lemma~\ref{lem:projective generators are flat over their endos}~\eqref{item:projective 1}.
We now prove part~(2).
Restricting the equivalences in part~(1) to the subcategories of finite length objects, and recalling that 
\[
\Mod_c(E_{\fB}^\op)^{\fl} \isoto \Mod^{\fl}(E_{\fB}^\op)
\]
by Lemma~\ref{properties of compact modules}~\eqref{item: compact 7}
(which applies since~$E_{\fB}$ is Noetherian, by Lemma~\ref{lem:projective generators are flat over their endos}~\eqref{item:projective 1}),
we obtain inverse anti-equivalences
\[\cA_{\fB}^{\fl} \isoto \Mod^{\fl}(E_{\fB}^\op), \; \pi \mapsto \Hom_{\fC_{\fB}}(P_{\fB}, \pi^\vee)\]
and
\[\Mod^{\fl}(E_{\fB}^\op) \to \cA_{\fB}^{\fl}, N \mapsto (N \cotimes_{E_{\fB}} P_{\fB})^\vee.\]
Now Pontrjagin duality gives an anti-equivalence $\Mod^{\fl}(E_{\fB}^\op) \to \Mod^{\fl}(E_{\fB})$, and
we have
\[\Hom_{\fC_{\fB}}(P_{\fB}, \pi^\vee)^\vee \cong P_{\fB} \otimes_{\cO\llbracket G\rrbracket _\zeta} \pi\]
by Lemma~\ref{lem:tensor facts}.
On the other hand, if~$M \in \Mod^{\fl}(E_{\fB})$ then
$N \coloneqq M^\vee$ is finitely presented over~$E_{\fB}^\op$, and so $N \cotimes_{E_{\fB}} P_{\fB} = N \otimes_{E_{\fB}} P_{\fB}$, 
and 
\[(N \otimes_{E_{\fB}} P_{\fB})^\vee \cong \Hom_{E_{\fB}^\op}(N, P_{\fB}^\vee) = \Hom_{E_{\fB}}^{\cont}(P_{\fB}, M),\]
where the isomorphism is because both sides are right exact functors of~$N$, and send~$N = E_{\fB}^\op$ to~$P_{\fB}^\vee$.
This concludes the proof of part~(2).

We now prove part~(3).
By assumption, $M$ is a quotient of~$P_{\fB}^{\oplus n}$ for some~$n$, and so
the module $\Hom_{\fC_{\fB}}(P_{\fB}, M)$ is a quotient of~$(E^{\op}_{\fB})^{\oplus n}$, i.e.\ it is a finitely generated $E^{\op}_{\fB}$-module.
By Lemma~\ref{finite cosocle implies Noetherian endomorphisms}, it follows that $\Hom_{\fC_{\fB}}(P_{\fB}, M)$ is finitely generated over the Bernstein centre of~$\fC_\fB$,
and so also over~$\End_{\fC_\fB}(M)$.
This concludes the proof because~$E_{\fB}$ and~$\End_{\fC_\fB}(M)$ are Noetherian, again by Lemma~\ref{finite cosocle implies Noetherian endomorphisms}.

Finally, part~(4) follows because part~(1) implies that $$M \iso
\Hom_{\fC_{\fB}}(P_{\fB},M)
\cotimes_{E_{\fB}} P_{\fB},$$
and we can use part~(3) to identify the completed tensor product with the uncompleted tensor product.
\end{proof}

\subsection{The \texorpdfstring{$p$}{p}-adic local Langlands correspondence for locally admissible \texorpdfstring{$\GL_2(\bQ_p)$-}{}representations.}\label{subsec:classical p-adic Langlands correspondence}
We now recall the results of Pa\v{s}k\={u}nas from~\cite{MR3150248} more precisely.
A key ingredient in both the statement and proof of these results is Colmez's functor
from $\GL_2(\Q_p)$-representations to $G_{\Q_p}$-representations, and so we begin
by briefly recalling some of the basic properties of this functor.

\subsubsection{Colmez's functor~\texorpdfstring{$V$}{V}}\label{ColmezV}
We refer to the fundamental work~\cite{MR2642409} for the definition of a covariant exact functor
\begin{equation}\label{eqn:definition of V}
(\sm.\, G)^{\fl}\to (\sm.\,G_{\bQ_p})^{\fl}
\end{equation}
whose values on absolutely irreducible objects can be found in \cite[Section~VII.4]{MR2642409}, see also \cite[Section~5.7]{MR3150248}.
This functor is usually denoted~$V$.
Since we wish to work with Galois representations of determinant~$\zeta \varepsilon^{-1}$, 
we will use the notation~$V$ to denote the twist of this functor by the Galois character~$\varepsilon^{-1}$,
and so its values on absolutely irreducible objects
(as classified in
Theorem~\ref{thm:classification-of-irreps} and
Remark~\ref{rem:irreps-as-parabolic-induction} above)
are as follows.
Note that if~$E'/E$ is a finite extension with ring of integers~$\cO'$, then~$V$ commutes with $\text{--}\otimes_\cO \cO'$, and with restriction of scalars from~$\cO'$ to~$\cO$:
this allows us to compute~$V$ on all irreducible objects.

\begin{lemma}\label{properties of V}\leavevmode
  \begin{enumerate}
    \item $V(\Ind_B^G(\chi_1 \otimes \chi_2\omega^{-1})) \cong \chi_2\omega^{-1}$ for any smooth characters $\chi_i: \bQ_p^\times \to \bF^\times$. %
    \item $V(\cInd_{KZ}^G\sigma_{a, b} /\HeckeT) \cong \Ind_{G_{\bQ_p^2}}^{G_{\bQ_p}}(\nr_{-\zeta(p)}\omega_2^{b+1}) \otimes \omega^{a-1}$.
    \item $V(\chi \circ \det) = 0$ and $V((\chi \circ \det) \otimes \St) = \chi$ for any smooth character $\chi : \bQ_p^\times \to \bF^\times$.
    \end{enumerate}    
\end{lemma}
\begin{proof}
See~\cite[Section~5.7]{MR3150248} and~\cite[Section~VII.4]{MR2642409}. 
\end{proof}

Recall that $\cO\llbracket G_{\bQ_p}\rrbracket $ has an anti-involution~$\dagger$ induced by $g \mapsto (\zeta\varepsilon^{-1})(g)g^{-1}$.
We write
\[
V^\dagger: \cA^{\fl} \to \Mod^{\fl}(\cO\llbracket G_{\bQ_p}\rrbracket ^{\op})
\]
for the functor obtained by precomposing the $\cO\llbracket G_{\bQ_p}\rrbracket $-action on $V$ with~$\dagger$.
Following~\cite[Section~5.7]{MR3150248}, we also introduce a covariant exact functor
\[
\Vcheck : \fC \to \Mod_c(\cO\llbracket G_{\bQ_p}\rrbracket )
\]
as follows.
We first introduce an autoequivalence $W \mapsto W^*$ of~$\Mod^{\fl}(\cO\llbracket G_{\bQ_p}\rrbracket )$
by letting~$W^*$ be the Pontrjagin dual of~$W$, with the left $\cO\llbracket G_{\bQ_p}\rrbracket $-module structure obtained by composing the natural right $\cO\llbracket G_{\bQ_p}\rrbracket $-module structure with
$\dagger$.
In other words, $W^*$ is the contragredient of~$W$, twisted by~$\zeta\varepsilon^{-1}$.
If~$M$ is an object of~$\fC$ of finite length, we then define
\[
\Vcheck(M) = V(M^\vee)^{*} \in \Mod^{\fl}(\cO\llbracket G_{\bQ_p}\rrbracket ),
\]
Then we extend~$\Vcheck$ to~$\fC$ by imposing compatibility with cofiltered limits, i.e.\ we use the identification $\fC = \Pro(\fC^{\fl})$ of Lemma~\ref{locally finite categories}
and the fact that $\Mod_c(\cO\llbracket G_{\bQ_p}\rrbracket )$ is complete.

\begin{rem}\label{Vcheck and Vdagger}
If~$M \in \fC$ has finite length, and we regard~$\Vcheck(M)^\vee$ as a right $\cO\llbracket G_{\bQ_p}\rrbracket $-module via contravariant Pontrjagin functoriality, 
then it is naturally isomorphic to $V^\dagger(M^\vee)$.
\end{rem}

\begin{rem}\label{comparison with Montreal paper}
Similarly to~$V$, our functor~$\Vcheck$ also differs from the~$\Vcheck$ in~\cite[Section~5.7]{MR3150248}
by a twist by~$\varepsilon^{-1}$.
We also caution the reader that on \cite[p.~320]{PaskunasBM},
Pa\v{s}k\={u}nas also writes
$\Vcheck$ to denote the contravariant functor on admissible Banach space representations obtained
by first dualizing to land in the isogeny category associated to~$\fC$,
and then applying the functor~$\Vcheck$ as we've recalled it above. 
We will not do this.
\end{rem}

\subsubsection{Blocks and pseudorepresentations}

The functor~$V$ induces a bijection between 
the blocks of~$\cA^{\ladm}$, and the $\Gal(\cbF_p/\bF)$-conjugacy classes of two-dimensional $\cbF_p$-valued pseudorepresentations
$\thetabar$ of~$G_{\Q_p}$ with determinant $\zetabar\omega^{-1}$. 
This bijection is given as follows.

\begin{defn}\label{defn: pseudorep labelled by block}\leavevmode
\begin{enumerate}
\item Assume that~$\fB$ is an $\bF$-rational block.
If~$\fB = \{\pi\}$ has type~\ref{item: ss block} or type~\ref{item: not abs irred block}, 
we set $\thetabar(\fB) \coloneqq  V(\pi)$.  
In all other cases, we set 
\[
\thetabar(\fB) \coloneqq 
V(\pi) + V(\pi)^{*}
\] 
for any infinite-dimensional irreducible~$\pi \in \fB$.
\item Assume that~$\fB$ is not an $\bF$-rational block. By Lemma~\ref{reducing gen+ to gen}, there exists an $\bF_{\fB}$-rational block~$\fB'$ 
of~$\cA_{\bF_{\fB}}$ equivalent to~$\fB$.
We define $\thetabar(\fB)$ to be the $\Gal(\cbF_p/\bF)$-conjugacy class of $\thetabar(\fB')$. 
\end{enumerate}
\end{defn}

Note that when~$\fB = \{\pi\}$ is an $\bF$-rational block of type~\ref{item: not abs irred pseudorep},
the compatibility of~$V$ with extension of scalars, together with Lemma~\ref{properties of V}~(1), 
shows that $V(\pi)$ is an irreducible, not absolutely irreducible representation $G_{\bQ_p} \to \GL_2(\bF)$. 
This implies that $\bF_{\thetabar(\fB)} = \bF_\fB$ for all blocks~$\fB$.
Lemma~\ref{properties of V} also implies that~$\thetabar(\text{--})$ is a bijection, with the following explicit description on $\bF$-rational objects.

\begin{lemma}\label{thetabar on blocks}\leavevmode
\begin{enumerate}
\item $\thetabar\{\cInd_{KZ}^G \sigma_{a, b}/\HeckeT\} = \omega^{a-1}\otimes\Ind_{\bQ_{p^2}}^{\bQ_p}\nr_{-\zeta(p)}\omega_2^{b+1}$.
\item $\thetabar\{\Ind_B^G(\chi_1 \otimes \chi_2\omega^{-1}), \Ind_B^G(\chi_2 \otimes \chi_1\omega^{-1})\} = \omega^{-1}(\chi_1+\chi_2)$.
\item $\thetabar\{\Ind_B^G(\chi \otimes \omega^{-1}\chi)\} = 2\omega^{-1}\chi$.
\item $\thetabar\{\chi \circ \det, (\chi \circ \det) \otimes \St, \Ind_B^G(\omega\chi \otimes \omega^{-1}\chi)\} = (1 + \omega^{-1})\chi$.
\item $\thetabar\{\cInd_{KZ}^G(\sigma_{a, b})/(\HeckeT^2-t\HeckeT+\lbar \zeta(p))\} = \omega^{a-1} \otimes \nr_{\HeckeT^2-t\HeckeT+\lbar \zeta(p)}$.
\end{enumerate}
\end{lemma}
\begin{proof}
  This is a direct computation. In case~(5), we have written $\nr_{\HeckeT^2-t\HeckeT+\lbar \zeta(p)}$ for the two-dimensional unramified representation 
with characteristic polynomial of Frobenius given by $\HeckeT^2-t\HeckeT+\lbar \zeta(p)$.
\end{proof}

From now on we write~$\mf{B}_{\thetabar}$ for the block corresponding to~$\thetabar$.
Slightly more informally, we will often label data associated to blocks in terms
of pseudorepresentations~$\thetabar$, rather than using the notation
$\fB$ or $\fB_{\thetabar}$.
For example, from now on we'll typically 
write $\fC_{\thetabar}$ in place of $\fC_{\fB_{\thetabar}}$.
Note that since $\bF_{\fB_{\thetabar}} = \bF_{\thetabar}$, this will not lead to ambiguity when discussing fields of definition.

 \begin{rem}\label{rem: labelling of blocks}
   Our labelling of the blocks and pseudorepresentations is now seen to be justified as
   follows: %
 \begin{itemize}
 \item[\ref{item: ss block}] is short for ``supersingular'', because these blocks
   contain the supersingular representations of~$\GL_2(\Qp)$. 
 \item[\ref{item: generic block}] is short for
   ``generic'', because of Definition~\ref{defn: generic characters}. Note that these pseudorepresentations are the traces of the generic representations on the Emerton--Gee stack
   for~$G_{\Qp}$.
 \item[\ref{item: non p-distinguished block}] is for the scalar pseudorepresentations.
 \item[\ref{item: Steinberg block}] is short for ``Steinberg'', because these blocks contain
a twist of the Steinberg representation of~$G$.
\item[\ref{item: not abs irred block}] becomes ``gen'' after extension of scalars.
 \end{itemize} 
\end{rem}

\subsubsection{Construction of projective objects}
We now exhibit for every block~$\fB_{\thetabar}$ a
particular projective object~$P_{\thetabar}$ of~$\fC_{\thetabar}$ with finite cosocle. 
If~$\thetabar$ is not of type \ref{item: Steinberg block}, then~$P_{\thetabar}$ will furthermore be a generator of~$\fC_{\thetabar}$.

\begin{defn}
\label{def:P theta}
We define a projective object~$P_{\thetabar}$ of~$\fC_{\thetabar}$ in the following way:
\begin{enumerate}
\item If~$\thetabar$ has type~\ref{item: ss pseudorep}, let~$\pi$ be the unique irreducible object of~$\fB_{\thetabar}$,
and let~$P_{\pi^\vee}$ be a projective envelope of~$\pi^\vee$ in~$\fC_{\thetabar}$.
Then $P_{\thetabar} \coloneqq  P_{\pi^\vee}^{\oplus 2}$.
\item If~$\thetabar$ has type~\ref{item: Steinberg pseudorep}, let
$\pi_1 = (\chi \circ \det) \otimes \pi_\alpha$ and~$\pi_2 = (\chi \circ \det) \otimes \St$ be the irreducible, infinite-dimensional objects of~$\fB_{\thetabar}$.
Let~$P_{\pi_1^\vee}$ and~$P_{\pi_2^\vee}$ be their projective envelopes in~$\fC_{\thetabar}$.
Then~$P_{\thetabar} \coloneqq  P_{\pi_1^\vee} \oplus P_{\pi_2^\vee}$.
\item Otherwise, $P_{\thetabar}\coloneqq  \oplus_{\pi} P_{\pi^{\vee}}$ is the direct sum of projective envelopes of the duals of the irreducible objects~$\pi$ of~$\fB_{\thetabar}$. 
\end{enumerate}
We also write
\begin{gather*}
E_{\thetabar} \coloneqq  \End_{\cO\llbracket G\rrbracket _\zeta}^{\cont}(P_{\thetabar}) = \End_{\fC_{\thetabar}}(P_{\thetabar})\\
\fT_{\thetabar} \coloneqq  \text{ right orthogonal to~$P_{\thetabar}$ in~$\fC_{\thetabar}$}\\
\fQ_{\thetabar} \coloneqq  \fC_{\thetabar}/\fT_{\thetabar}.
\end{gather*}

\end{defn}

\subsubsection{The functor \texorpdfstring{$\Vcheck$}{V} and Cayley--Hamilton modules.}\label{subsubsec:Vcheck and CH}
In this subsection we study the restriction of~$\Vcheck$ to~$\fC_{\thetabar}$ for all pseudorepresentations~$\thetabar$.
When~$\thetabar$ has type~\ref{item: not abs irred pseudorep}, and $\cO'/\cO$ is the finite unramified extension with residue field~$\bF_{\thetabar}$, 
Lemma~\ref{reducing gen+ to gen} produces an $\bF_{\thetabar}$-rational block~$\fC_{\thetabar'}$ of~$\fC_{\cO'}$,
and shows that restriction of scalars $\fC_{\thetabar'} \to \fC_{\thetabar}$ is an equivalence.
Since~$V$ commutes with restriction of scalars, the diagram
\[
\begin{tikzcd}
\fC_{\thetabar'} & \Mod_c(\cO'\llbracket G_{\bQ_p}\rrbracket )\\
\fC_{\thetabar} & \Mod_c(\cO\llbracket G_{\bQ_p}\rrbracket )
\arrow["\Vcheck", from=1-1, to=1-2]
\arrow["\Vcheck", from=2-1, to=2-2]
\arrow["\sim", from=1-1, to=2-1]
\arrow[from=1-2, to=2-2]
\end{tikzcd}
\]
commutes.
In the rest of this subsection, we will therefore be able to assume without loss of generality that~$\thetabar$ is $\bF$-rational.

We begin by forming the profinite module $\Vcheck(P_{\thetabar}),$  
noting that it 
has commuting $\cO$-linear left actions of~$E_{\thetabar}$ (by functoriality) and~$\cO\llbracket G_{\bQ_p}\rrbracket $.
By Lemma~\ref{lem:projective generators are flat over their endos}(2), $P_{\thetabar}$ is a pseudocompact left $E_{\thetabar}$-module, and so $\Vcheck(P_{\thetabar})$
can be seen as an object of $\Mod_c(E_{\thetabar})$. 
The following result is essentially a restatement of some of the results 
of~\cite{MR4350140}, which in turn, since we are assuming $p\geq 5$, are
essentially a restatement of results from~\cite{MR3150248}. 

\begin{prop}\label{VcheckP}
\leavevmode
Let~$\thetabar$ be a two-dimensional $\cbF_p$-valued pseudorepresentation of~$G_{\bQ_p}$.
\begin{enumerate}
\item
The $\cO\llbracket G_{\bQ_p}\rrbracket $-action on~$\Vcheck(P_{\thetabar})$ factors through $\tld R_{\thetabar}$, and induces an isomorphism
\[
\tld R_{\thetabar} \xrightarrow{\sim} \End_{\Mod_c(E_{\thetabar})}\bigl(\Vcheck(P_{\thetabar})\bigr).
\]
\item $\Vcheck(P_{\thetabar})$ 
is free of rank one over each of $\tld R_{\thetabar}$ and $E_{\thetabar}$; in particular there is an isomorphism $\End_{\Mod_c(E_{\thetabar})}\bigl(\Vcheck(P_{\thetabar})\bigr) \cong E^{\op}_{\thetabar}$, 
well-defined up to inner automorphism. 
\item
The composition of the isomorphisms in parts~(1) and~(2) is an isomorphism $\tld R_{\thetabar} \iso E_{\thetabar}^{\op}$,
well-defined up to inner automorphism.
\end{enumerate}
\end{prop}
\begin{proof}
Without loss of generality, $\thetabar$ is $\bF$-rational.
Furthermore, all statements in the proposition may be checked after a finite unramified base extension $\text{--} \otimes_\cO \cO'$. 
We can therefore assume without loss of generality that~$\thetabar$ does not have type~\ref{item: not abs irred pseudorep}.
The statement of~(1) is then a consequence of~\cite[Theorem~6.3, Theorem~6.13]{MR4350140}, together with the assertion in Proposition~\ref{CHtorsionfree}
that $\tld R_{\thetabar}$ is $\cO$-torsion free.
The statement in~(2) concerning the $E_{\thetabar}$-structure of~$\Vcheck(P_{\thetabar})$ is contained in the last three paragraphs of the proof of~\cite[Proposition~4.18]{MR4350140}.

If we choose an isomorphism $\imath:\Vcheck(P_{\thetabar}) \iso E_{\thetabar}$
in $\Mod_c(E_{\thetabar})$, 
then we find that
$$\End_{\Mod_c(E_{\thetabar})}\bigl(\Vcheck(P_{\thetabar})\bigr) \iso E_{\thetabar}^{\op}$$
(with the isomorphism depending 
upon the choice of $\imath$, and thus being well-defined up to an inner automorphism)
and that $\Vcheck(P_{\thetabar})$ is also free of rank one over its endomorphism algebra
$\End_{\Mod_c(E_{\thetabar})}\bigl(\Vcheck(P_{\thetabar})\bigr)$. 
Combining this with the statement of~(1) proves the remainder of the proposition.
\end{proof}

Part~(2) of Proposition~\ref{VcheckP} allows us to regard $\Vcheck(P_{\thetabar})$ as an $(\Rtilde_{\thetabar},E_{\thetabar}^{\op})$-bimodule, free of rank one with respect to the action of either ring.
We now use Morita theory to give an alternative description of the functor~$\Vcheck.$
We recall that the functor $\Vcheck$ is trivial on 
$\fT_{\thetabar}$: this is a vacuous statement except in case~\ref{item: Steinberg block},
in which case it amounts to the statement that $\Vcheck\bigl((\chi\circ \det)^{\vee}\bigr) = 0$.
Thus~$\Vcheck$ induces a functor $\fQ_{\thetabar} \to \Mod_c(\cO\llbracket G_{\bQ_p}\rrbracket )$
(again exact and cofiltered limit-preserving, because so is $\Vcheck$)
which we continue to denote by~$\Vcheck$. 
Because our projective object $P_{\thetabar}$ has been chosen precisely
so that it detects the quotient $\fQ_{\thetabar}$, 
we can then use Morita theory with respect to $P_{\thetabar}$ to 
describe~$\Vcheck$.  This is the subject of the following theorem,
which is essentially due to~Pa\v{s}k\={u}nas.
In particular, the statement of~(1), to the effect that $\Vcheck$ induces 
the indicated equivalence, is one of the main results of~\cite{MR3150248}.
\begin{theorem}
\label{thm:VP equiv}
Let $\fB_{\thetabar}$ be a block of~$\cA^{\ladm}$,
and let
$\thetabar$ be the associated
$\Gal(\cbF_p/\bF)$-conjugacy class of pseudorepresentations.
Then the following are true:
\begin{enumerate}
\item The functor $\Vcheck: \fC_{\thetabar} \to \Mod_c(\cO\llbracket G_{\bQ_p}\rrbracket )$ factors through~$\Mod_c(\tld R_{\thetabar})$, is naturally isomorphic to 
\begin{equation}\label{formula for Vcheck to prove}
M \mapsto
\Hom_{\fC_{\thetabar}}(P_{\thetabar}, M)
\otimes_{E_{\thetabar}}
\Vcheck(P_{\thetabar}),
\end{equation}
and induces an equivalence
$\fQ_{\thetabar} \iso \Mod_c(\tld R_{\thetabar}).$
\item The functor $V^\dagger: \cA_{\thetabar}^{\fp} \to
\Mod^{\fl}(\cO\llbracket G_{\bQ_p}\rrbracket )$ factors through $\Mod^{\fl}(\widetilde{R}_{\thetabar}^{\op})$, and 
is naturally isomorphic to
\[
\pi \mapsto \Hom_{E_{\thetabar}}(\Vcheck(P_{\thetabar}), E_{\thetabar}) \otimes_{E_{\thetabar}} (P_{\thetabar} \otimes_{\cO\llbracket G\rrbracket _{\zeta}} \pi),
\]
where~$E_{\thetabar}$ acts by right multiplication 
on the $\Hom$-space, 
and~$\tld R_{\thetabar}$ acts through~$\Vcheck(P_{\thetabar})$.
\end{enumerate}
\end{theorem}
\begin{proof}
As in Section~\ref{subsubsec:theta Morita}, we write $\overline{M}$ for the
image of an object $M$ of $\fC_{\thetabar}$ in~$\fQ_{\thetabar}$.
Then Lemma~\ref{properties of projective objects}(2) shows that $\Hom_{\fC_{\thetabar}}(P_{\thetabar}, M) = \Hom_{\fQ_{\thetabar}}(
\overline{P}_{\thetabar},\overline{M})$ for all~$M \in \fC_{\thetabar}$.
Hence $\End_{\fQ_{\thetabar}}(\lbar P_{\thetabar}) \cong E_{\thetabar}$ is a Noetherian profinite $\cO$-algebra.
By Lemma~\ref{properties of projective objects}~(3), $\lbar P_{\thetabar}$ is a projective generator of~$\fQ_{\thetabar}$, and it has finite cosocle in~$\fQ_{\thetabar}$.
Hence Proposition~\ref{compact Morita theory} applies, and shows that $\Hom_{\fQ_{\thetabar}}(\overline{P}_{\thetabar},\text{--})$
an equivalence $\fQ_{\thetabar} \iso \Mod_c(E_{\thetabar}^{\op}),$ with quasi-inverse %
given by $-\cotimes_{E_{\thetabar}} \lbar P_{\thetabar}$.
For any object $M$ of~$\fC_{\thetabar},$
we now see that 
$$\Vcheck(M) = \Vcheck( \overline{M})
\iso \Vcheck\bigl( \Hom_{\fQ_{\thetabar}}(\overline{P}_{\thetabar}, 
\overline{M})\cotimes_{E_{\thetabar}} \overline{P}_{\thetabar} \bigr) 
\iso \Hom_{\fC_{\thetabar}}(P_{\thetabar},M)\otimes_{E_{\thetabar}} \Vcheck(P_{\thetabar}).
$$
For the second isomorphism, we use the fact that $\Vcheck$ is exact and cofiltered limit-preserving,
and also the fact that $\Vcheck(P_{\thetabar})$ is finitely generated (in fact,
free of rank one) over $E_{\thetabar}$, so that $\text{--}\cotimes_{E_{\thetabar}} 
\Vcheck(P_{\thetabar})$ can be computed as the usual tensor product.
This gives the factorization and the natural isomorphism of~(1). 
Since, as we've already noted, 
$\Hom_{\fQ_{\thetabar}}(\overline P_{\thetabar}, \text{--})$ is an equivalence,
it then follows from Proposition~\ref{VcheckP}~(2)
that $\Vcheck$ induces an equivalence $\fQ_{\thetabar} \to \Mod_c(\tld R_{\thetabar}).$
This completes the proof of~(1). 

We now prove part~(2). 
By Remark~\ref{Vcheck and Vdagger}, the functor~$V^\dagger$ is the Pontrjagin dual of~$\Vcheck$.
So it follows from part~(1) that
  \[
    V^\dagger(\pi) \cong
\left ( 
\Hom_{\fC_{\thetabar}}(P_{\thetabar}, \pi^\vee)
      \otimes_{E_{\thetabar}}
 \Vcheck(P_{\thetabar})
\right )^{\vee}
  \]
  functorially in~$\pi$.  Since $\Vcheck(P_{\thetabar})$ is a free
  left $E_{\thetabar}$-module of finite rank, we deduce from this that
  \[
    V^\dagger(\pi) \cong \Hom_{E_{\thetabar}}\bigl(\Vcheck(P_{\thetabar}),
    E_{\thetabar}\bigr) \otimes_{E_{\thetabar}}
    \Hom_{\fC_{\thetabar}}(P_{\thetabar}, \pi^\vee)^\vee.
  \]
The result now follows from Lemma~\ref{lem:tensor facts}~(1).
\end{proof}

\begin{rem}
\label{description of Bernstein centres}
Write~$\cZ_{\thetabar}$ for the centre of the category~$\fC_{\thetabar}$. 
Then there exists a unique isomorphism $Z(\tld R_{\thetabar}) \to \cZ_{\thetabar}$ 
such that the functor~$\Vcheck$ is $Z(\tld R_{\thetabar})$-linear.
This is part of the results of~\cite{MR3150248}; from our point of view, it is a direct consequence of Theorem~\ref{thm:VP equiv}, except in case~\ref{item: Steinberg pseudorep}, 
in which case we need the additional fact that $\fC_{\thetabar} \to \fQ_{\thetabar}$
induces an isomorphism of centres, which follows from~\cite[Theorem~10.80, Corollary~10.77]{MR3150248}.

Note also that by Proposition~\ref{CHtorsionfree}, the natural map $R_{\thetabar}^{\ps} \to Z(\tld R_{\thetabar})$ is an isomorphism.
Hence we can reformulate the above as the existence of a unique isomorphism $R_{\thetabar}^{\ps} \to \cZ_{\thetabar}$ such that the functor~$\Vcheck$ is $R_{\thetabar}^{\ps}$-linear.
It then follows that the functors~$V$ and~$V^\dagger$ are also $R_{\thetabar}^{\ps}$-linear.
In fact, since $V(\pi) = \Vcheck(\pi^\vee)^*$, to prove the $R_{\thetabar}^{\ps}$-linearity of~$V$ it suffices to check that for all $z \in R_{\thetabar}^{\ps} \subset \tld R_{\thetabar}$ 
and all $M \in \Mod^{\fl}(\tld R_{\thetabar})$,
we have the equality $z_{M^*} = z_{M}^*$ as central endomorphisms of~$M^*$.
The equality $z_{M^\vee} = z_M^\vee$ holds by definition (here~$\vee$ denotes the Pontrjagin dual).
Now, again by definition, we have $z_{M^*} = (z^\dagger)_{M^\vee}$, and $z_M^* = z_M^\vee$ (since the functors $M \mapsto M^*$ and~$M \mapsto M^\vee$ act in the same way on morphisms). 
Since the anti-involution~$\dagger$ is trivial on the centre of~$\tld R_{\thetabar}$, we conclude that~$z_{M^*} = z_M^*$, as desired.
\end{rem}

Since $\Vcheck:\fQ_{\thetabar} \to \Mod_c(\tld R_{\thetabar})$ is an equivalence
of categories, it also admits a Morita-theoretic description.   It is easy to make
this explicit,
since the natural isomorphism~\eqref{formula for Vcheck to prove}
shows that the equivalences $\Hom_{\fQ_{\thetabar}}(\overline{P}_{\thetabar},\text{--})$
and $\Vcheck(\text{--})$ differ by tensoring with the
free rank one bimodule~$\Vcheck(P_{\thetabar}).$

\begin{defn}\label{defn: tld P}
We define
$\tld P_{\thetabar} \coloneqq 
\Hom_{E_{\thetabar}}\bigl(\Vcheck(P_{\thetabar}),E_{\thetabar}
\bigr)\otimes_{E_{\thetabar}} P_{\thetabar}$
and $\tld E_{\thetabar} \coloneqq  \End_{\fC_{\thetabar}}(\tld P_{\thetabar})$.
\end{defn}

If~$\thetabar$ is not of type \ref{item: Steinberg block}, then~$\tld P_{\thetabar}$ is a generator of~$\fC_{\thetabar}$ (since the same is true of~$P_{\thetabar}$).
Since $\Vcheck(P_{\thetabar})$ is an $(\tld R_{\thetabar}, E_{\thetabar}^{\op})$-bimodule 
which is free of rank one over both rings,
we see that 
$\Hom_{E_{\thetabar}}\bigl(\Vcheck(P_{\thetabar}),E_{\thetabar}\bigr)$
is an $(\tld R_{\thetabar}^{\op}, E_{\thetabar})$-bimodule which is again
free of rank one over both rings. 
Thus $\tld P_{\thetabar}$ is a projective object of $\fC_{\thetabar}$,
with a natural isomorphism
\begin{equation}
\label{eqn:endo identification}
\tld R_{\thetabar}^{\op} \iso \tld E_{\thetabar} = \End_{\fC_{\thetabar}}(\tld P_{\thetabar}).
\end{equation}
Hence we can regard the functor~$V^\dagger$ 
as a functor
\begin{equation}\label{refactored Vdagger Ethetabar}
V^\dagger: \cA_{\thetabar}^{\fp} \to \Mod^{\fl}(\tld E_{\thetabar}),
\end{equation}%
and the functor~$V$ as a functor
\begin{equation}\label{refactored V Rthetabar}
V: \cA_{\thetabar}^{\fp} \to \Mod^{\fl}(\tR_{\thetabar}).
\end{equation}

\begin{rem}
  The main reason to introduce~$\tld P_{\thetabar}$ is to obtain the canonical isomorphism~\eqref{eqn:endo identification}, and therefore 
  avoid making a choice of isomorphism as in Proposition~\ref{VcheckP}~(3).
  Throughout the paper we will mostly work with~$\tld P_{\thetabar}$, rather than~$P_{\thetabar}$. 
  As an exception, in Section~\ref{subsec:steinberg computations} we will need an explicit presentation of $\tld E_{\thetabar}$ when~$\thetabar = 1 + \omega$, 
  and so it will be convenient to choose an isomorphism $P_{\thetabar} \isoto \tld P_{\thetabar}$.
  By construction, this choice is equivalent to the choice of an $E_{\thetabar}$-basis of the free rank-one $E_{\thetabar}$-module~$\Vcheck(P_{\thetabar})$.
  A convenient choice of basis is described in Lemma~\ref{lem:P_thetabar normalization}.
\end{rem}

\begin{lemma}\label{needed for universal deformation}
\leavevmode
\begin{enumerate}
\item
There is a natural isomorphism
$$\Vcheck(\text{--}\, ) \iso \Hom_{\fC_{\thetabar}}({\tld P_{\thetabar}}, \text{--})$$
of functors $\fC_{\thetabar} \to \Mod_c(\tld R_{\thetabar})$. 
In particular,
the functor
$$M \mapsto M\cotimes_{\tld R_{\thetabar}} \overline{\tld P_{\thetabar}}$$
provides a quasi-inverse to the equivalence $\Vcheck:\fQ_{\thetabar} \to \Mod_c(\tld R_{\thetabar}).$
\item
There is a natural isomorphism
$$V^\dagger(\text{--}\,) \iso \tld P_{\thetabar} \otimes_{\cO\llbracket G\rrbracket _{\zeta}}\text{--}
$$
of functors
$\cA_{\thetabar}^{\fp} \to \Mod^{\fl}(\tld E_{\thetabar}).$
\item\label{item:64}  If $\thetabar$ is not of type~{\em\ref{item: Steinberg block}}, then $V^\dagger:\cA_{\thetabar}^{\fp} \to \Mod^{\fl}(\tld E_{\thetabar})$ is an equivalence.
\end{enumerate}
\end{lemma}
\begin{proof}
The claimed natural isomorphism in~(1) follows from
the definition of~$\tld P_{\thetabar}$ together with~\eqref{formula for Vcheck to prove}.
The rest of part~(1) is then 
a special case of 
Proposition~\ref{compact Morita theory},

The statement of~(2) is dual to the statement of~(1) in
the same way that the statement of 
part (2) of Theorem~\ref{thm:VP equiv} is dual to the statement of part~(1)
of that theorem, and is proved in an identical manner.

The final part is immediate from Lemma~\ref{lem:recapitulated Morita}~\eqref{item:63}, since~$\tld P_{\thetabar}$ is a projective generator of~$\fC_{\thetabar}$.
\end{proof}

\subsection{Localization of smooth representations \texorpdfstring{of~$\GL_2(\bQ_p)$}{}}\label{subsec:recollections-of-localization-paper}
In this and the next section we recall various results from our paper~\cite{DEGlocalization},
and explain how to translate them from the abelian categorical context (which
is the focus of that paper) to the stable $\infty$-categorical context in which
we are working.  
Along the way, we establish various relationships between the stable
$\infty$-categories of representations that we will be considering.

\subsubsection{A chain of projective lines.}\label{chain of projective lines} 
One of the main results of~\cite{DEGlocalization} is a localization theory for~$\cA$ over the Zariski site of a scheme~$X$ over~$\bF$, which is a chain of projective lines of length~$(p \pm 1)/2$, with ordinary double points (where the sign is equal to~$-\zeta(-1)$).
There is a bijection between 
companion pairs of $\zeta$-compatible weights and the irreducible
components\footnote{In~\cite[Prop.\ 2.5.4]{DEGlocalization}, we
label the components as $X(\tau)$, where $\tau$ runs over the set of isomorphisms
classes of irreducible cuspidal representations of
$\GL_2(\bF_p)$ over~$E$ with central character~$\zeta$. %
The map $\tau \mapsto \JH(\overline{\tau})$ (i.e.\ passing to the Jordan--H\"older factors
of the semisimplified mod $p$ reduction) induces a bijection between the set of such~$\tau$ 
and the set of companion pairs.  This bijection is described in more detail
in the proof of Proposition~\ref{prop: explicit description of underlying reduced of U} below. Compare also~\cite[Rem.\ 2.5.7]{DEGlocalization}.}
of~$X$,
and we write $X(\sigmasigmacomp)$ for the component corresponding to the pair~$\{\sigma,
\sigmacomp\}$.

If~$\sigma = \sigma_{a, b}$ is $\zeta$-compatible, and~$b \ne p-1$, there is a
morphism (\cite[Definition~2.6.1]{DEGlocalization})
\begin{equation}\label{eqn:defn-of-fsigma}
f_\sigma : \Spec \cH_G(\sigma) = \bA^1  \to X(\sigmasigmacomp)
\end{equation}
which is an open immersion when~$b \ne p-2$ (given by $x\mapsto x^{\pm 1}$, 
the sign depending on~$\sigma$), and is a degree two morphism when~$b = p-2$ (given by~$x\mapsto (x+x^{-1})^{\pm 1}$).
The identification of $\Spec \cH_G(\sigma)$ with~$\bA^1$ uses the isomorphism~\eqref{eqn:Satake-iso-with-coordinate-X}.
When~$b = p-1$, we furthermore define~$f_{\sigma_{a, p-1}} = f_{\sigma_{a, 0}}$.

There is a bijection $x \mapsto \fB_x$ between the set of closed points of~$X$
and the set of blocks of~$\cA^{\ladm}$, given by
\begin{equation}\label{closed points to blocks}
\fB_x = \bigcup_{\sigma,y} %
\JH\left ( \cInd_{KZ}^G\sigma \otimes_{\cH_G(\sigma)} y \right ),
\end{equation}
i.e.\ the union of the Jordan--H\"older factors of $\cInd_{KZ}^G\sigma
\otimes_{\cH_G(\sigma)} y$ for those~$\sigma,y$ with $f_{\sigma}(y)=x$. 
Composing with the bijection in Definition~\ref{defn: pseudorep labelled by block},
we see that there is also a bijection between the set of closed points of~$X$,
and the set of $\Gal(\cbF_p/\bF)$-conjugacy classes of two-dimensional $\cbF_p$-pseudorepresentations 
$\thetabar$
of~$G_{\bQ_p}$, having determinant~$\lbar \zeta \omega^{-1}$. Accordingly, we
will typically denote a closed point $x \in X$ by~$\thetabar$. 
Note that in the notation of Section~\ref{subsec:pseudoreps-recollections}, the residue field of~$\thetabar$ is~$\F_{\thetabar}$.

\begin{rem}\label{rem:explicit formula for f_sigma}
Using the identification of closed points of~$X$ with $\cbF_p$-pseudorepresentations that we have just described, 
we have the following explicit formula for the morphism~$f_\sigma$ from~\eqref{eqn:defn-of-fsigma}:
if~$\sigma=\sigma_{a,b}$, then for ~$t\in\Fpbartimes$ (regarded as an $\Fpbar$-point of $\Spec \cH_G(\sigma)$ via the isomorphism~\eqref{eqn:Satake-iso-with-coordinate-X}) we have
\begin{equation}
  \label{eq:formula-for-f-sigma-t}
 f_{\sigma}(t)= \nr_{t^{-1}}\omega^{a-1}+\nr_{t}\zetabar\omega^{-a},
\end{equation}
while
\begin{equation}
  \label{eq:formula-for-f-sigma-t-at-zero}
 f_{\sigma}(0)= \Ind_{G_{\bQ_{p^2}}}^{G_{\bQ_p}}(\nr_{-\zeta(p)}\omega_2^{b+1}) \otimes \omega^{a-1}.
\end{equation}
Hence the
points $0,\infty$ of $X(\sigmasigmacomp)$ correspond to irreducible
pseudorepresentations.
\end{rem} 

\subsubsection{Localization of abelian categories.}
We associate to every closed subset $Y \subseteq |X|$ the
full Serre subcategory~$\cA_Y \subset \cA$ of objects all of whose irreducible
subquotients are contained in $\bigcup_{\thetabar \in Y} \fB_{\thetabar}$. We
write $i_{Y,*}:\cA_Y\to\cA$ for the fully faithful inclusion functor. Write~$U \coloneqq  X \setminus Y$,  let~$\cA_U\coloneqq \cA/\cA_Y$ denote the Serre
quotient, and write $j_U^{*}:\cA \to \cA_U$ for the quotient functor. We write
$\cA_Y^{\fp}$ for the full subcategory of finitely generated objects of~$\cA_{Y}$, and~$\cA_U^{\fp}$ for
the essential image of~$\cA^{\fp}$ in~$\cA_U$.
When~$Y = \{\thetabar\}$ is a singleton, the category~$\cA_{\{\thetabar\}}$
agrees with the category~$\cA_{\fB_{\thetabar}}$
of~\eqref{eqn:product-decomposition-blocks}, and we will denote it
by~$\cA_{\thetabar}$ from now on. %
On the other hand, if~$Y=X$ then~$\cA_Y=\cA$, so all the results we recall
below about the categories~$\cA_Y$ apply in particular to~$\cA$ itself.

The following proposition summarizes some of the fundamental properties of the
categories~$\cA_Y$ established in~\cite{DEGlocalization} (see
Appendix~\ref{subsec:1-categories} for a brief recollection of some material on
locally Noetherian categories). %
\begin{prop}
  \label{prop:recollections-of-localization-paper}\leavevmode
  \begin{enumerate}
  \item The categories $\cA_Y$ and~$\cA_U$ are locally
Noetherian. %
\item The compact objects of~$\cA_Y$ {\em (}resp.\ $\cA_U${\em )} are the finitely generated
  objects~$\cA_Y^{\fp}$ {\em (}resp.\ $\cA_U^{\fp}${\em )}.
\item The inclusion $i_{Y,*}:\cA_{Y}\to\cA$ is exact, preserves colimits,
  preserves injectives, and preserves
  compact objects.
\item  The quotient functor  $j_U^{*}:\cA \to \cA_U$ admits an exact and fully faithful
  right adjoint $j_{U,*}:\cA_U\to\cA$. %
  In particular, $\cA_Y$ is a localizing subcategory of~$\cA$.
\item The assignment $U \mapsto \cA_U$ defines a stack~$\cA_{\bullet}$ of abelian
categories on~$X_{\Zar}$.
  \end{enumerate}
\end{prop}
\begin{proof}The first four points are immediate from~\cite[Lem.\ 3.1.2, 3.1.7,
  Prop.\ 3.1.13~(2), Cor.\ 3.5.8]{DEGlocalization}, together with the observation that by~\cite[Proposition~A.1.1]{DEGlocalization}, 
  the Noetherian objects of~$\cA$, resp.\ $\cA_Y$, are precisely the
  compact objects, and since~$\cA_Y \subset \cA$ is a Serre
  subcategory, the inclusion preserves Noetherian objects. 
The final point is~\cite[Theorem~3.3.1]{DEGlocalization}.
\end{proof}

The exact functor $i_{Y,*}$ induces an exact functor
  $\ihat_{Y,*}:\Pro(\cA_Y^{\fp})\to\Pro(\cA^{\fp})$, which by the adjoint functor theorem
  has a left adjoint
  $\ihat_Y^*: \Pro(\cA^{\fp})\to\Pro(\cA_Y^{\fp})$. 
  This functor was denoted~$\widehat{(\text{--})}_Y$ in~\cite{DEGlocalization}.
  \begin{lem}
    \label{i-hat-is-exact}%
    The functor
    $\ihat_Y^{*}:\Pro(\cA^{\fp})\to\Pro(\cA_Y^{\fp})$ is exact.
  \end{lem}
  \begin{proof}%
    By~\cite[Cor.\ 3.5.5]{DEGlocalization}, the
    functor 
    \[\ihat_Y^*|_{\cA^{\fp}}: \cA^{\fp} \to \Pro(\cA_Y^{\fp})\] is
    exact. The lemma then follows immediately from
    Lemma~\ref{lem:left-adjoint-to-Pro-functor-abelian}~(2).\qedhere
  \end{proof}
  
\begin{remark}
If $Y \subseteq X$ consists of a single closed point~$\thetabar$, then we
will typically write $i_{\thetabar,*}$ and $\ihat_{\thetabar}^*$
rather than $i_{Y,*}$ and $\ihat_Y^*$.
\end{remark}

\begin{rem}%
  \label{rem:completion-to-subset-of-open}
  Passing to $\Ind$-extensions, we obtain
  an adjoint pair of functors, still denoted~$(\ihat^*_Y, \ihat_{Y, *})$, between $\Ind \Pro \cA^{\fp}$ and~$\Ind \Pro \cA^{\fp}_Y$.
  The
  functor~$\ihat_Y^{*}$ 
  was also denoted $\wht{(\text{--})}_Y$ in~\cite{DEGlocalization}.

  If~$\thetabar$ is a closed point of~$U$, then the functor $\ihat^*_{\thetabar} : \cA^{\fp} \to \Pro \cA_{\thetabar}^{\fp}$
  factors through~$\cA^{\fp}_{U}$, %
  and we continue to denote the resulting functor $\cA_U^{\fp} \to \Pro \cA_{\thetabar}^{\fp}$ by~$\ihat^{*}_{\thetabar}$.
Bearing in mind this slight abuse of notation, we have a natural isomorphism %
\begin{equation}\label{eqn:unit localization and completion on U}
\ihat^{*}_{\thetabar}\isoto \ihat^{*}_{\thetabar}j_{U}^{*}:\cA^{\fp} \to \Pro \cA_{\thetabar}^{\fp}.
\end{equation}
Similarly, the composite $\ihat^*_{\thetabar}j_{U, *}$, which is \emph{a priori} a functor $\cA_U^{\fp} \to \Ind \Pro \cA_{\thetabar}^{\fp}$, 
is actually valued in~$\Pro \cA_{\thetabar}^{\fp}$,
and there is a natural isomorphism
\begin{equation}\label{eqn:unit localization and completion on U II}
\ihat^*_{\thetabar}\isoto \ihat_{\thetabar}^*j_{U,*} : \cA^{\fp}_U \to \Pro \cA_{\thetabar}^{\fp};
\end{equation}
to see this, it suffices to precompose the $\Ind$-extension of~\eqref{eqn:unit localization and completion on U} with $j_{U, *}$, and use the fact that the counit
$j_U^*j_{U, *} \to 1$ is an isomorphism.
See also~\cite[Section~3.5.21]{DEGlocalization} for related material.
\end{rem}

The functors~$\ihat_{Y}^{*}$ and~$j_{U,*}$ can be explicitly computed in
some cases. By~\cite[Lem.\ 3.5.2]{DEGlocalization}, if~$\pi$ is an object
of~$\cA^{\fp}$, there is a natural isomorphism
\begin{equation}\label{explicit description of completion}
\ihat_Y^{*}\pi \iso \quoteslim{} \pi',
\end{equation}
where $\pi'$ runs over the cofiltered %
set of quotients of~$\pi$ 
lying in~$\cA_Y^{\fp}$. 
Furthermore, the unit of adjunction $\pi \to \ihat_{Y, *}\ihat_Y^* \pi$ corresponds under~\eqref{explicit description of completion} to the inverse limit of the corresponding quotient maps.
When~$\sigma$ is a Serre weight, and $\pi = \cInd_{KZ}^G \sigma$, we can produce
an explicit cofinal set in~\eqref{explicit description of completion} in the following way.

\begin{defn}\label{f_Y}
Let~$\sigma$ be a Serre weight, and let $Y \subset X$ be a closed subset. %
\begin{enumerate}
\item If~$f_\sigma^{-1}(Y)$ is finite, we define~$f_Y \in \cH(\sigma)= \F[\HeckeT]$ to be the unique monic squarefree generator of the ideal of $f_{\sigma}^{-1}(Y)$.
If~$Y = \{\thetabar\}$ is a singleton, we write~$f_{\thetabar}$ for~$f_Y$.
\item If~$f_\sigma^{-1}(Y)$ is infinite, we define~$f_Y \coloneqq  0$.
\end{enumerate}
\end{defn}

When~$Y = \{\thetabar\}$ is a single $\bF$-rational point, 
the polynomial~$f_{\thetabar}$ has degree one, unless $\sigma = \sigma_{a, p-2}$,
in which case~$f_{\thetabar}$ can have degree two (and can be reducible).
Furthermore, by~\cite[Lemma~3.1.8(2)]{DEGlocalization}, $(\cInd_{KZ}^G \sigma)/f_{\thetabar}(\cInd_{KZ}^G \sigma)$ is the maximal multiplicity-free quotient of $\cInd_{KZ}^G \sigma$
which is an object of~$\cA_{\thetabar}$.

\begin{lem}%
  \label{lem:localization-of-cInd-sigma}Let~$\sigma \coloneqq  \sigma_{a, b}$ be a Serre
  weight, and let~$Y \subset X$ be a closed subset.  
  \begin{enumerate}
  \item 
 \begin{enumerate}
 \item\label{item:30} The set 
 $\{\cInd_{KZ}^G\sigma/f_Y^n\cInd_{KZ}^G\sigma: n \geq 1\}$ is cofinal in~\eqref{explicit description of completion},
and so we have
\begin{equation}\label{eqn:recollection-of-explicit-completion-of-cInd-sigma}\ihat_Y^{*}\cInd_{KZ}^G \sigma\iso \quoteslim{n}(\cInd_{KZ}^G \sigma )/f_Y^{n}\cInd_{KZ}^G\sigma.\end{equation}
\item\label{item:31} The  unit morphism \[\cInd_{KZ}^G \sigma \to
j_{U,*}j_U^* \cInd_{KZ}^G \sigma\] is given by the natural map %
\begin{equation}\label{eqn:recollection-of-description-of-j*-cInd-sigma}\cInd_{KZ}^G \sigma \rightarrow
(\cInd_{KZ}^G \sigma)[1/f_Y].\end{equation}
\end{enumerate}
\item\label{item:support of cInd sigma} 
Let~$\thetabar \in X$ be a closed point, and let~$\rhobar : G_{\bQ_p} \to \GL_2(\cbF_p)$ be a semisimple Galois representation with trace~$\thetabar$.
Then~$\ihat_{\thetabar}^* \cInd_{KZ}^G \sigma \ne 0$ if and only if~$\sigma$ is a Serre weight of~$\rhobar$.
\end{enumerate}
\end{lem}
\begin{proof}
  Part~(1a) is~
  \cite[Cor.\ 3.5.3]{DEGlocalization}, while Part~(1b)  is~ \cite[Prop.~3.1.13(1)]{DEGlocalization}.  
We now prove part~\eqref{item:support of cInd sigma}.
By part~\eqref{item:30},
$\ihat^*_{\thetabar}(\cInd_{KZ}^G \sigma) \ne 0$ if and only if $f_{\{\thetabar\}}$ is not a unit in $\cH(\sigma)$.
By Definition~\ref{f_Y}, this occurs if and only if 
there exists $t \in \Spec \cH(\sigma)$ such that $f_\sigma(t) = \thetabar$.
By~\eqref{eq:formula-for-f-sigma-t} and~\eqref{eq:formula-for-f-sigma-t-at-zero}, this occurs if and only if $t \in \cbF_p^\times$ and $\thetabar = \nr_{t^{-1}}\omega^{a-1} + \nr_t\lbar \zeta \omega^{-a}$,
or $t = 0$ and $\thetabar = \Ind_{G_{\bQ_{p^2}}}^{G_{\bQ_p}}(\nr_{-\zeta(p)}\omega_2^{b+1}) \otimes \omega^{a-1}$.
Bearing in mind Definition~\ref{defn:Serre weights of rhobar}, this is equivalent to~$\sigma$ being a Serre weight of~$\rhobar$.
\end{proof}

\begin{rem}\label{rem: what happens when f_Y = 0}
  In the setting of Lemma~\ref{lem:localization-of-cInd-sigma}~\eqref{item:30}, if~$f_Y = 0$ then $\cInd_{KZ}^G \sigma$ is an object of~$\cA_Y^{\fg}$, and so~\eqref{eqn:recollection-of-explicit-completion-of-cInd-sigma}
  exhibits its completion as a constant pro-object.
  Similarly, the localization $j_U^* \cInd_{KZ}^G \sigma$ is equal to zero, as is the right-hand side of~\eqref{eqn:recollection-of-description-of-j*-cInd-sigma}, which is by definition the 
  $\cH_G(\sigma)$-module localization 
  \[
  \cInd_{KZ}^G \sigma \otimes_{\cH_G(\sigma)} \cH_G(\sigma)[1/f_Y].
  \]
\end{rem}

\subsubsection{Tensor product and completion}
\label{subsubsec:tensor and completion}
If~$P$ is a projective object of~$\fC_{\thetabar}$ with finite length cosocle,
then we have seen in Lemma~\ref{lem:projective generators are flat over their endos}
that~$E \coloneqq  \End_{\fC}(P)$ is a compact Noetherian $\cO$-algebra
and that~$P \in \Mod_c(E)$.
Hence~$P$ is a right $\cO\llbracket G\rrbracket _\zeta$-module in~$\Mod_c(E)$, and so Lemma~\ref{lem:abelian-composition-right-exact-tensor}
shows that the functor
\[
P\otimes_{\cO\llbracket G\rrbracket _{\zeta}} \text{--}: \Mod^{\fp}(\cO\llbracket G\rrbracket _\zeta) \to \Mod_c(\cO)
\]
can be lifted through the forgetful functor $\Mod_c(E) \to \Mod_c(\cO)$ to a functor 
\begin{equation}
\label{what we are factoring}
P \otimes_{\cO\llbracket G\rrbracket _\zeta} -: \Mod^{\fp}(\cO\llbracket G\rrbracket _\zeta) \to \Mod_c(E).
\end{equation}
Since $\Mod_c(E) = \Pro \Mod^{\fl}(E)$, this is an instance of Lemma~\ref{lem:derived-tensor-product-EW}~\eqref{item:86a}.
Since $\cA^{\fp}$ is a full subcategory of $\Mod^{\fp}(\cO\llbracket G\rrbracket _\zeta)$,
the restriction of \eqref{what we are factoring} to~$\cA^{\fp}$ extends  uniquely to a cofiltered limit-preserving functor
\begin{equation}\label{what we are factoring II}
P \ccotimes_{\cO\llbracket G\rrbracket _\zeta} - : \Pro(\cA^{\fp}) \to \Mod_c(E).
\end{equation}
This is an instance of Lemma~\ref{lem:derived-tensor-product-EW}~\eqref{item:86}, which strictly speaking
considers the Pro-extension to all of $\Pro\bigl(\Mod^{\fp}(\cO\llbracket G\rrbracket _{\zeta})\bigr)$,
but in what follows we will be focussed on its full subcategory $\Pro(\cA^{\fp}).$
In particular, in that context we have the following exactness result.

\begin{lemma}
\label{lem:tensor exactness}
The functor $P \otimes_{\cO\llbracket G\rrbracket _\zeta} (\text{--})$ is exact on $\cA^{\fp}$,
and hence the Pro-extended functor
$P\ccotimes_{\cO\llbracket G\rrbracket _{\zeta}} (\text{--}): \Pro(\cA^{\fp}) \to \Mod_c(E)$ 
is also exact.
\end{lemma}
\begin{proof}
The forgetful functors $\Mod_c(E) \to \Mod(E)$ and $\Mod(E) \to \Mod(\cO)$ are exact,
and so it suffices to verify the claimed exactness after applying their composite,
i.e.\ after forgetting the topology and $E$-module structure on
$P\otimes_{\cO\llbracket G\rrbracket _{\zeta}} \pi.$
Doing this, we are reduced to considering the tensor product with $P$ and $\pi$
thought of as abstract $\cO\llbracket G\rrbracket _{\zeta}$-modules, and the claimed
exactness follows from the Tor-vanishing
of Lemma~\ref{lem:projective Tor vanishing}~(2).
\end{proof}

Regarding $\thetabar$ as a closed point of~$X$,
we may consider the completion functor~$\ihat_{\thetabar}^*$; %
we will use it together with~\eqref{what we are factoring II} to relate~\eqref{what we are factoring}
to its restriction to~$\cA_{\thetabar}^{\fp}$.
We first note the following.

\begin{lem}
\label{lem:finite length tensor}
If $\pi$ is an object of~$\cA_{\thetabar}^{\fp}$,
then $P\otimes_{\cO\llbracket G\rrbracket _{\zeta}} \pi$ is a finite length $E$-module.
\end{lem}
\begin{proof}
It suffices to show that $P\otimes_{\cO\llbracket G\rrbracket _{\zeta}} \pi$ has finite $\cO$-length, or equivalently that $(P\otimes_{\cO\llbracket G\rrbracket _{\zeta}}\pi)^{\vee}$
has finite $\cO$-length.  
By Lemma~\ref{lem:tensor
  facts}~\eqref{item:38}, $(P\otimes_{\cO\llbracket G\rrbracket _{\zeta}}\pi)^{\vee}$
 is isomorphic to $\Hom_{\cO\llbracket G\rrbracket _{\zeta}}(\pi,P^{\vee}) = \Hom_\cA(\pi, P^\vee)$. 
 By Lemma~\ref{fg equals fl for Athetabar}, $\pi$ has finite $\cA$-length, and so
 $\Hom_\cA(\pi, P^\vee) = \Hom_\cA(\pi, \soc_{\cA, n} P^\vee)$ for some~$n \geq 0$, by Lemma~\ref{general Loewy facts}(1).
 This concludes the proof since~$\soc_{\cA, n} P^\vee$ has finite $\cA$-length, by the dual to Lemma~\ref{lem:Loewy facts}.
\end{proof}

Lemma~\ref{lem:finite length tensor} shows that~\eqref{what we are factoring}
restricts to a functor
\begin{equation}\label{what we are extending from}
P \otimes_{\cO\llbracket G\rrbracket _\zeta} (\text{--}): \cA_{\thetabar}^{\fp} \to \Mod^{\fl}(E).
\end{equation}
Because of the equivalence~\eqref{improved compact modules and Pro},
the restriction of~\eqref{what we are factoring II} to~$\Pro(\cA_{\thetabar}^{\fp})$ therefore coincides with the $\Pro$-extension 
\[
\Pro(\cA_{\thetabar}^{\fp}) \to \Pro\bigl(\Mod^{\fl}(E)\bigr) \iso \Mod_c(E)
\]
of~\eqref{what we are extending from},
which we will also denote $P \ccotimes_{\cO\llbracket G\rrbracket _\zeta} (\text{--})$.

\begin{lemma}\label{tensor and completion}%
If~$P$ is a projective object of~$\fC_{\thetabar}$ with finite length cosocle 
and~$E = \End_{\fC_{\thetabar}}(P)$,
then the natural transformation 
\[
P \otimes_{\cO\llbracket G\rrbracket _\zeta} (\text{--}) \to P \ccotimes_{\cO\llbracket G\rrbracket _\zeta} \ihat_{\thetabar,*}\ihat_{\thetabar}^*(\text{--})
\]
 of functors $\cA^{\fp} \to \Mod_c(E)$ 
{\em (}induced by the unit of adjunction $\id_{\Pro \cA^{\fp}} \to \ihat_{\thetabar,*}\ihat_{\thetabar}^*${\em )}
is an isomorphism. %
\end{lemma}
\begin{proof}
Choose~$\pi \in \cA^{\fp}$, and let $P_n = P/\operatorname{rad}_n P$. 
By Lemma~\ref{lem:Loewy facts}, $P_n$ has finite $\fC$-length, and
\[
P \iso \varprojlim_n P_n.
\]
Similarly,
let $\pi_n \coloneqq  \pi^{(n)}_{\thetabar}$ be
the quotient of~$\pi$ defined in Lemma~\ref{bounded quotients of finite length}.
By the cofinality statement in Lemma~\ref{bounded quotients of finite length}, the isomorphism~\eqref{explicit description of completion} can be rewritten as
\[\ihat_{\thetabar, *}\ihat_{\thetabar}^*\pi \iso \quoteslim{n} \pi_n.\] 
By Lemma~\ref{bounded quotients of finite length} and Lemma~\ref{general Loewy facts} respectively, 
we have isomorphisms of finite torsion $\cO$-modules
$$\Hom_{\cA}(\pi_n, \operatorname{soc}_\cA^n P^\vee) 
\iso
\Hom_{\cA}(\pi, \operatorname{soc}_\cA^n P^\vee)$$ 
and
$$\Hom_{\cA}(\pi_n, \operatorname{soc}_\cA^n P^\vee) 
\iso
\Hom_{\cA}(\pi_n, P^\vee)$$ 
(induced by pullback and pushforward, respectively).
Applying $\Hom_{\cO}(\text{--}, E/\cO)$, we obtain by Lemma~\ref{lem:tensor facts}~(1)
a pair of isomorphisms
\[
P_n \otimes_{\cO\llbracket G\rrbracket _\zeta} \pi \iso P_n \otimes_{\cO\llbracket G\rrbracket _\zeta} \pi_n,
\]
and
\[
P \otimes_{\cO\llbracket G\rrbracket _\zeta} \pi_n
\iso
P_n \otimes_{\cO\llbracket G\rrbracket _\zeta} \pi_n.
\]
We thus find that 
\begin{multline*}
P\otimes_{\cO\llbracket G\rrbracket _{\zeta}} \pi \iso
\varprojlim_n (P_n\otimes_{\cO\llbracket G\rrbracket _{\zeta}} \pi)
\iso \varprojlim_n (P_n\otimes_{\cO\llbracket G\rrbracket _{\zeta}} \pi_n) 
\\
\iso \varprojlim_n (P\otimes_{\cO\llbracket G\rrbracket _{\zeta}} \pi_n)
\iso P \ccotimes_{\cO\llbracket G\rrbracket _{\zeta}} \ihat_{\thetabar, *}\ihat_{\thetabar}^*\pi
\end{multline*}
(the first isomorphism being an application 
of Lemma~\ref{lem:tensor and projective limits} in the category~$\Mod_c(E)$, 
and the final isomorphism holding because
$P\ccotimes_{\cO\llbracket G\rrbracket _{\zeta}} (\text{--})$ is cofiltered limit-preserving by construction),
as required. %
\end{proof}

\subsection{Derived categories of smooth representations}\label{subsec:derived-categories-smooth-representations}
We now write $D(\cA_Y)$ for the (unbounded) derived stable $\infty$-category of
the Grothendieck category~$\cA_Y$. (See~\cite[App.\
A.1]{emerton2023introduction} and Appendix~\ref{subsec: category theory
  background} for a brief recollection of the material on stable
$\infty$-categories that we use below, and Appendix~\ref{subsec:t-structures}
for a recollection of $t$-structures.) 

Since the inclusion $i_{Y,*}: \cA_Y \to \cA$ is exact and preserves colimits, it
induces (by Lemmas~\ref{lem:derived-adjunction-unbounded-Grothendieck} and~\ref{lem:exact-functor-induces-t-exact}) a continuous $t$-exact
functor %
\begin{equation}\label{definition of i_Y*}
i_{Y,*} : D(\cA_Y) \to D(\cA).
\end{equation}
(Here and below, we will use the same notation for exact functors between
abelian categories, and their $t$-exact extensions to derived categories.)
Recall that in the notation of Definition~\ref{defn:D-sub-A-B} we write
$D_{\cA_Y}(\cA)$ for the full subcategory of $D(\cA)$ consisting of objects all
of whose cohomologies lie in~$\cA_Y$, and similarly
$D_{\cA_Y}(\cO\llbracket G\rrbracket _{\zeta})$ (where~$D(\cO\llbracket G\rrbracket _{\zeta})$ is the derived
category of~$\cO\llbracket G\rrbracket _{\zeta}$-modules).

\begin{lemma}\label{unbounded i_* fully faithful}%
The $t$-exact continuous functor $i_{Y,*} : D(\cA_Y) \to D(\cA)$ %
is fully faithful,
and induces an equivalence
$D(\cA_Y) \iso D_{\cA_Y}(\cA)$. Furthermore we have a commutative diagram of
 t-exact 
fully faithful continuous functors
\[\begin{tikzcd}
 D(\cA_Y)	 & D_{\cA_Y}(\cA) & D(\cA) \\
D_{\cA_{Y}}(\cO\llbracket G\rrbracket _{\zeta})	 & D_{\cA_{Y}}(\cO\llbracket G\rrbracket _{\zeta})  & D_{\cA}(\cO\llbracket G\rrbracket _{\zeta})
\arrow["\simeq", from=1-1, to=1-2]
	\arrow[ hook, from=1-2, to=1-3]
	\arrow[ "\simeq", from=1-1, to=2-1]
	\arrow[ "\simeq", from=1-2, to=2-2]
        	\arrow[ "\simeq", from=1-3, to=2-3]
	\arrow["=", from=2-1, to=2-2]
        	\arrow[ hook, from=2-2, to=2-3]
\end{tikzcd}\]
\end{lemma}
\begin{proof}
 We begin by noting that since we have a colimit-preserving, exact and fully faithful functor $\cA\to
\Mod(\cO\llbracket G\rrbracket _\zeta)$, %
we
have a continuous $t$-exact functor $D(\cA)\to D(\cO\llbracket G\rrbracket _\zeta)$. %
By an identical argument to the proof of ~\cite[Prop.\
E.2.2]{emerton2023introduction}, 
taking into account Remark~\ref{embedding of A in modules},
this functor satisfies the hypotheses of
Theorem~\ref{thm:unbounded-adjoints-localizing}, so it induces
the right-most vertical equivalence in the commutative diagram above.

In turn, the right-most equivalence restricts to the middle equivalence.
It therefore suffices to show that the composite
\[D(\cA_{Y})\to D(\cA)\to D(\cO\llbracket G\rrbracket _\zeta)\]is fully faithful and induces an
equivalence $D(\cA_{Y})\to D_{\cA_Y}(\cO\llbracket G\rrbracket _\zeta)$. Since $i_{Y,*}: \cA_Y
\to \cA$ preserves injectives by \cite[Cor.\ 3.5.8]{DEGlocalization}, 
and its right
adjoint has finite cohomological
amplitude by~\cite[Cor.\ 3.6.3]{DEGlocalization}, this follows from an
application of Theorem~\ref{thm:unbounded-adjoints-localizing} to the fully faithful
functor $\cA_{Y}\to \Mod(\cO\llbracket G\rrbracket _\zeta)$. %
\end{proof}

We recall the following standard lemma.
\begin{lemma}
\label{lem:modules of finite projective dimension}
If $A$ is a coherent ring and $M$ is a finitely generated
$A$-module for which $\RHom_A(M,\text{--})$ is of bounded amplitude,
then $M$ is compact as an object of~$D(A)$.
\end{lemma}
\begin{proof}By for example~\cite[Lem.\ 4.1.6]{MR1269324}, $M$ 
 has a finite length resolution
by finitely generated projective $A$-modules, and is therefore compact in~$D(A)$
(see for example~\cite[Prop.\ 4.6]{MR1214458}). %
\end{proof}

\begin{prop}\label{prop: t structure on sm G is regular}%
  The objects of~$\cA^{\fp}$ are compact in~$D(\cA)$.
\end{prop}
\begin{proof} 
By Lemma~\ref{unbounded i_* fully faithful}, there is a fully faithful
continuous functor
 $D(\cA)\into
D(\cO\llbracket G\rrbracket _{\zeta})$. By Lemma~\ref{lem:group-ring-coherent}, the
ring~$\cO\llbracket G\rrbracket _\zeta$ is coherent, so it follows from Lemma~\ref{lem:modules of
  finite projective dimension} 
that it
suffices to show that for any finitely generated object $\pi$ of~$\cA$,
the functor $\RHom_{\cO\llbracket G\rrbracket _\zeta}(\pi,\text{--})$ has bounded
amplitude.

To this end, recall that by for
example~\cite[Section~2.4.14]{inertialp-adic} we have a short exact
sequence 
$$0 \to \cInd_{N}^G \delta\pi \to \cInd_{KZ}^G \pi \to \pi \to 0$$
where~$N$ is the normalizer of the Iwahori subgroup~$\Iw$,
and~$\delta$ is the nontrivial quadratic character of~$N/\Iw Z$. We
can therefore reduce to showing that for any finitely generated~$\pi$, and any
subgroup~$H$ of~$G$ which contains~$\Iw Z$ with finite coprime-to-$p$
index, 
the functor $\RHom_{\cO\llbracket G\rrbracket _\zeta}(\cInd_H^G\pi,\text{--})$ has bounded
amplitude. %

Since~$\pi$ is finitely generated, and~$G/KZ$ is countable, $\pi$ is countably generated over~$\cO$.
Hence we can write~$\pi$ as a countably-indexed filtered colimit $\colim_iV_i$ of finitely
generated~$H$-subrepresentations. 
Then
\[
\RHom_{\cO\llbracket G\rrbracket _\zeta}(\cInd_H^G\pi,\text{--})=\Rlim_i\RHom_{\cO\llbracket H\rrbracket _\zeta}(V_i,\text{--}).
\]
Using the fact that~$R^j\lim_i = 0$ for~$j \geq 2$ (since the limit is countably indexed) we are reduced to checking that
$\RHom_{\cO\llbracket H\rrbracket _\zeta}(V_i,\text{--})$ has bounded amplitude. 
Since~$V_i$ is an $H$-linear direct summand of $\cInd_{\Iw_1 Z}^H(V_i)$ (since 
the index of $\Iw_1 Z$ in $H$ is prime-to-$p$), it
suffices to check that $\RHom_{\cO\llbracket \Iw_1 Z\rrbracket _\zeta}(V_i, \text{--})$ has
bounded amplitude.
This is a consequence of the fact that the ring~$\cO\llbracket \Iw_1Z\rrbracket _{\zeta}$ has finite global
dimension, which follows from~\cite[Thm.\ 3.26]{MR1924402}. 
More precisely, the decomposition~\eqref{decomposition of pro-p subgroups} shows that
we have an isomorphism of rings
$$\cO\llbracket \Iw_1 Z\rrbracket _{\zeta} \iso \cO\llbracket \Iw_1/Z_1\rrbracket $$ 
and since $p \geq 5,$ we see that $\Iw_1$ is torsion-free, 
and thus so is $\Iw_1/Z_1$ (being a direct factor of the former group).
Thus~\cite[Thm.\ 3.26]{MR1924402} applies to the target of this isomorphism.
\end{proof}
\begin{rem}
  \label{rem:regular-used-p-5}If~$p=2$ or~$3$ then $\Iw_1$ has a nontrivial $p$-torsion element, namely $-1 \in \bQ_p$ when~$p = 2$, 
  and a conjugate of $\zeta_3 \in \bQ_p(\zeta_3) \subset M_2(\bQ_p)$ when~$p = 3$.
  So we anticipate that the analogue of Proposition~\ref{prop: t structure
    on sm G is regular} would fail in these cases.
\end{rem}

\begin{defn}\label{defn:Dfp}
  We write~$D^b_{\fp}(\cA_Y)$ for~$D^b_{\cA_Y^{\fp}}(\cA_{Y})$, and similarly ~$D^b_{\fp}(\cA_U)$ for~$D^b_{\cA_U^{\fp}}(\cA_{U})$.
\end{defn}

\begin{cor}
  \label{cor:all-the-nice-properties-of-D(A)}Let  $Y \subseteq X$ be a closed
  subset with open complement $U \subseteq X$. Then
  \begin{enumerate}
  \item The canonical functors
$D^b(\cA_Y^{\fp}) \to D^b_{\fp}(\cA_Y)$
and $D^b(\cA_U^{\fp}) \to D^b_{\fp}(\cA_U)$ %
are equivalences.
  \item  The categories~$D(\cA_{Y})$ and~$D(\cA_U)$ are compactly
  generated, and we have $D(\cA_Y)^c=D^b_{\fp}(\cA_Y)$ and
  $D(\cA_U)^c=D^b_{\fp}(\cA_U)$.
  Hence the natural maps $\Ind D^b_{\fp}(\cA_Y) \isoto D(\cA_Y)$ and $\Ind D^b_{\fp}(\cA_U) \to D(\cA_U)$ are equivalences.
\item $j_U^*:\cA \to \cA_U$
induces  equivalences
$$D(\cA)/D(\cA_Y) \iso D(\cA_U)$$ and $$D^b_{\fp}(\cA)/D^b_{\fp}(\cA_Y) \iso D^b_{\fp}(\cA_U).$$
\item The functor~ $j_{U,*} : D(\cA_U) \to D(\cA)$ is fully faithful, continuous, and
is right adjoint to $j_U^*: D(\cA) \to D(\cA_U)$.
  \end{enumerate}
\end{cor}
\begin{proof}
  This is immediate from
  Proposition~\ref{prop:localizing-derived-compact-objects} (whose hypotheses
  hold by   %
  Lemma~\ref{unbounded i_* fully
    faithful}, Proposition~\ref{prop:recollections-of-localization-paper}, and Proposition~\ref{prop: t structure on sm G is regular}).
\end{proof}
We now give an explicit set of compact generators  (in the sense recalled in
Remark~\ref{rem:compact generation}) for the compactly generated
categories~$D(\cA)$ and~$D(\cA_{U})$.
\begin{cor}
  \label{cor:explicit-generators-D(A)}
  Letting $\sigma$ range over the $\zeta$-compatible
  Serre weights, the collection $\{\cInd_{KZ}^G \sigma\}$ is a set of compact generators
  of~$D(\cA)$, and the collection $\{j_U^*\cInd_{KZ}^G \sigma\}$ is a set of compact
  generators of~$D(\cA_{U})$.
\end{cor}
\begin{proof}%
 The second statement follows from the first. For the first, by Proposition~\ref{prop:abelian-compact-implies-derived-compactly-generated}
  it suffices to show that  $\{\cInd_{KZ}^G \sigma\}$ is a set of weak generators
  of~$\cA$, i.e.\ that for each $\pi\in\cA$, there exists a non-zero map
  $\cInd_{KZ}^G \sigma\to\pi$ for some~$\sigma$. This is standard; see e.g.\ the
  proof of~\cite[Lem.\ 2.2.3]{DEGlocalization}.
\end{proof}

Finally we turn to the exact functor $\ihat_Y^{*}:\Pro(\cA^{\fp})\to\Pro(\cA_Y^{\fp})$
from Lemma~\ref{i-hat-is-exact}.
Note that since~$i_{Y,*}$ is exact,
it induces a $t$-exact functor 
\begin{equation}\label{definition of i_Y*-hat}
\ihat_{Y,*}:\Pro D^b_{\fp}(\cA_Y) \to \Pro D^b_{\fp}(\cA). 
\end{equation}

\begin{lem}
  \label{lem:derived-t-exact-ihat-Y}
  Let
  \begin{equation}
  \ihat_Y^{*}:\Pro D^b_{\fp}(\cA) \to \Pro D^b_{\fp}(\cA_Y)
  \end{equation}
  be the $\Pro$-extension of the $t$-exact functor $D^b_{\fp}(\cA) \to \Pro D^b_{\fp}(\cA_Y)$ whose restriction to hearts is $\ihat^*_Y$.
  Then~$\ihat^*_Y$ is left adjoint to~$\ihat_{Y,*}$. %
\end{lem}
\begin{proof}
  This is immediate from  Lemma~\ref{i-hat-is-exact} and Proposition~\ref{prop:pro-adjoints} (applied with $\cC = \cA^{\fp}, \cC' = \cA_Y^{\fp}$, and~$f = \ihat^*_Y$).
\end{proof}

\subsection{Finiteness of some \texorpdfstring{$\Ext$}{}
  groups}\label{subsec:Ext-group-computations}%
In this subsection we prove some results on extension groups between parabolic inductions which we will use in the proof of our main theorem (see Proposition~\ref{properties of Hecke modules sigma
  version statements we use}). %

\subsubsection{Removing the central character.}
We will make significant use of the results of
Heyer~\cite{Heyerderived, Heyergeometric}
on parabolic induction and its adjoints. Since Heyer's papers
work with the category~$\sm.\,G$ of all smooth $\cO$-representations of~$G$,
possibly without a central character, we begin by explaining how to relate 
$\Ext$ groups in~$\cA$ and~$\sm.\,G$.

By Lemma~\ref{unbounded i_* fully faithful}, the inclusion $\cA \to \Mod(\cO\llbracket G\rrbracket _\zeta)$ induces a $t$-exact equivalence
\[
D(\cA) \to D_{\sm}(\cO\llbracket G\rrbracket _\zeta),
\]
where the right-hand side denotes the full subcategory of complexes with smooth cohomology.
By ~\cite[Proposition~E.2.2]{emerton2023introduction}, we have a commutative
diagram \begin{equation}\label{comparing with group algebra}
\begin{tikzcd} 
D(\sm.\, G) \arrow[r, "\sim"] & D_{\sm}(\cO\llbracket G\rrbracket )\\
D(\cA) \arrow[u] \arrow[r, "\sim"] & D_{\sm}(\cO\llbracket G\rrbracket _\zeta) \arrow[u].
\end{tikzcd}
\end{equation}
Writing~$A \coloneqq \cO\llbracket \bQ_p^\times\rrbracket $ and~$B \coloneqq \cO$, we can thus study the relationship between~$\Ext^i_{\cA}$ and~$\Ext^i_{\sm.\, G}$ as a particular case of the following problem:
Suppose given a morphism $\zeta: A\to B$ of commutative rings
and a (not necessarily commutative) $A$-algebra $R$,
and write $S \coloneqq  B\otimes_A R.$
Then, if $M$ and $N$ are two $S$-modules, which can then also be regarded as $R$-modules
via the canonical morphism $R \to S$,
what is the relationship between
$\Ext^{\bullet}_R(M,N)$ 
and
$\Ext^{\bullet}_S(M,N)$? 

To begin with, if $M$ is an $R$-module and $N$ is an $S$-module, then deriving the tensor--$\Hom$
adjunction gives
$$\RHom_R(M,N) \iso \RHom_S( B\otimes^L_A M, N).$$
Now, if $M$ is itself an $S$-module, then
\[B\otimes^L_A M \iso (B\otimes^L_A B)\otimes_B^L M.\]
Thus we obtain, for $S$-modules $M$ and~$N$, 
an isomorphism
\begin{equation}\label{eqn:computing-derived-tensor}\RHom_R(M,N) \iso \RHom_S\bigl( (B\otimes^L_A B)\otimes^L_B M,N).\end{equation}

\begin{lemma}\label{Koszul central character} %
There exists an isomorphism of functors $D(\cA)^\op \times D(\cA) \to D(\cO)$:
\[
\RHom_{\sm.\, G}(\text{--}, \text{--}) \cong \RHom_{\cA}(\text{--}, \text{--}) \oplus \RHom_{\cA}(\text{--}, \text{--})^{\oplus 2}[-1] \oplus \RHom_{\cA}(\text{--}, \text{--})[-2]. 
\]
\end{lemma}
\begin{proof}
We apply the above discussion with~$A=\cO\llbracket \Q_p^{\times}\rrbracket $, $B=\cO$,  and the morphism $A\to B$ being $\zeta:\cO\llbracket \Q_p^{\times}\rrbracket  \to \cO$;
and we take $R=\cO\llbracket G\rrbracket $ and $S=\cO\llbracket G\rrbracket _{\zeta}$. 
Choosing a uniformizer of $\bQ_p$ and an isomorphism $1+p\bZ_p \isoto \bZ_p$ then induces an isomorphism %
\[A \cong \cO\llbracket U\rrbracket [V^{\pm 1}][\bF_p^\times],\] %
and we note that the group ring (denoted $[\bF_p^\times]$) of the cyclic group~$\bF_p^\times$ is semisimple in residue characteristic~$p$.
Computing with a Koszul complex in the variables~$U$ and~$V$, we see that
\[
B\otimes^L_A B \iso B[2] \oplus B[1]^{\oplus 2} \oplus B.
\]
Using~\eqref{eqn:computing-derived-tensor} we have
\begin{multline*}%
\RHom_{\cO\llbracket G\rrbracket }(M,N)
\iso
\\
\RHom_{\cO\llbracket G\rrbracket _{\zeta}}(M,N)
\oplus
\RHom_{\cO\llbracket G\rrbracket _{\zeta}}(M,N)^{\oplus 2}[-1] 
\oplus
\RHom_{\cO\llbracket G\rrbracket _{\zeta}}(M,N)[-2].
\end{multline*}
The result now follows from~\eqref{comparing with group algebra}.
\end{proof}

\subsubsection{Extensions of universal parabolic inductions.}
We begin by introducing some notation involving the connected components of~$X$.
We refer to Section~\ref{subsec: defining U} for motivation and some related material.

\begin{defn}\label{defn:f(t)}%
  Let~$\sigmasigmacomp$ be a companion pair of Serre weights.
  \begin{enumerate}
  \item Define a polynomial $f \in \F[t]$ by
  \[
    f(t) \coloneqq
    \begin{cases}
      t & \text{if } \sigmasigmacomp \text{ is of type~\ref{item:
          generic pseudorep}}  \\
      t(t^2-\zeta(p)) & \text{if } \sigmasigmacomp \text{ is of type~\ref{item:
                 non p-distinguished pseudorep} or~\ref{item: Steinberg pseudorep}}.
    \end{cases}
  \]
  \item Let~$\Ybad\subset X$ be the closed
subset determined by the finite set of
points %
which become of type~\ref{item: ss
  pseudorep}, \ref{item: non p-distinguished pseudorep}
or~\ref{item: Steinberg pseudorep} after a finite extension of~$\cO$. 
Let~$\Ugood$ be the complement of~$\Ybad$.  
\item   If $\{\sigma,\sigmacomp\}$ is a pair of companion weights, let
  \[U(\sigmasigmacomp) \coloneqq  \Ugood\cap X(\sigmasigmacomp),\]
  \end{enumerate}
\end{defn}

We fix for the rest of
this section a companion pair of Serre weights~$\{\sigma,\sigmacomp\}$, and a choice of $\sigma_1, \sigma_2 \in \{\sigma, \sigmacomp\}$ (possibly with $\sigma_1 = \sigma_2$).
Recall from Section~\ref{induced representations} that we write
\begin{gather*}
\cH_G(\sigma_i) = \End_G(\cInd_{KZ}^G\sigma_i),\\
\cH_T = \End_{T}(\cInd_{T_0 Z}^T 1)\iso \F[X^{\pm 1}],
\end{gather*}%
where $T_0 = T(\bZ_p)$ is the maximal compact subgroup of the diagonal torus~$T$. 
We let $\chi_i = (\sigma_i)^{\Iw_1}$, viewed as a character of~$T_0 Z$. In Section~\ref{induced representations}, we have
canonically identified the endomorphism algebras $\End_T(\cInd_{T_0 Z}^T
\chi_i)$ with~$\cH_T$, and via the Satake morphisms, we have identified each of
the~$\cH_G(\sigma_i)$ with~$\F[\HeckeX]$ (with~$\HeckeX$ being identified with the spherical
Hecke operator~$\HeckeT \in \cH_G(\sigma_i)$).

In the rest of this section we will denote the element~$f(\HeckeX)$ of~$\cH_G(\sigma_i)$, where the polynomial~$f$ is defined in Definition~\ref{defn:f(t)}, by~$f$.
In the notation of Definition~\ref{defn:f(t)},
the vanishing set of~$f$ in $\Spec \cH_G(\sigma_i)$ coincides with $f_{\sigma_i}^{-1}(\Ybad \cap X(\sigmasigmacomp))$.

We now 
write~$j^*, j_*$ for the functors~$j_{U(\sigmasigmacomp)}^*, j_{U(\sigmasigmacomp), *}$ defined in Proposition~\ref{prop:recollections-of-localization-paper}.
Then~\eqref{eqn:recollection-of-description-of-j*-cInd-sigma}
implies that
\begin{equation}\label{identifying localization}
j_*j^* \cInd_{KZ}^G\sigma_i = \cInd_{KZ}^G\sigma_i[1/f].
\end{equation}
Since
\[
\left (\cInd_{KZ}^G\Sym^{0} \otimes \operatorname{det}^a \right )[1/(\HeckeX^2-\zeta(p))] \cong \left ( \cInd_{KZ}^G\Sym^{p-1} \otimes \operatorname{det}^a \right )[1/(\HeckeX^2-\zeta(p))],
\]
in the case that~$\sigma_1$ or~$\sigma_2$ is a twist of~$\Sym^0$, we can and do
replace it by the corresponding twist of~$\Sym^{p-1}$. 
(We do this in order that the isomorphism~\eqref{universal parabolic induction} is valid.)

Since~\eqref{universal parabolic induction} is equivariant for the Satake map, and $\Ind_{\lbar B}^G$ preserves colimits,
we may combine~\eqref{universal parabolic induction} with~\eqref{identifying localization} to obtain an isomorphism
\begin{equation}\label{universal parabolic induction II}
j_*j^*\cInd_{KZ}^G\sigma_i \cong \Ind_{\lbar B}^G\left ( (\cInd_{T_0 Z}^T\chi_i)[1/f] \right ).
\end{equation}
In particular we see that $\cH_T[1/f]$ acts on~$j_*j^*
\cInd_{KZ}^G\sigma_2$ and thus on 
the $\Ext$ groups \[\Ext^i_{\cA}(j_*j^* \cInd_{KZ}^G\sigma_1, j_*j^* \cInd_{KZ}^G\sigma_2).\] %

\begin{prop}\label{properties of Hecke modules}
The $\cH_T[1/f]$-modules~ $\Ext^i_{\cA}(j_*j^* \cInd_{KZ}^G\sigma_1, j_*j^* \cInd_{KZ}^G\sigma_2)$ are finitely generated and torsion-free. 
\end{prop}
\begin{proof}
  We begin by noting that~\eqref{identifying localization} 
  induces an isomorphism of~$\cH_T[1/f]$-modules
\begin{equation}\label{reformulating E^i}
\Ext^i_{\cA}(j_*j^* \cInd_{KZ}^G\sigma_1, j_*j^* \cInd_{KZ}^G\sigma_2) \cong \Ext^i_{\cA}(\cInd_{KZ}^G\sigma_1[1/\HeckeX], \cInd_{KZ}^G\sigma_2[1/f]).
\end{equation}
Indeed, it suffices to show that for all ~$i$ we
have \[\Ext^i_{\cA}(\cInd_{KZ}^G\sigma_1[1/f]/\cInd_{KZ}^G\sigma_1[1/\HeckeX],
  \cInd_{KZ}^G\sigma_2[1/f])=0,\]and this is immediate from
~\cite[Corollary~3.1.15]{DEGlocalization} (with~$Y$ taken to be the vanishing
locus of~$f$ in~$X$, as before).
Furthermore, by \eqref{universal parabolic induction} we have 
\begin{equation}\label{eqn:instance-of-universal-parabolic-induction}\cInd_{KZ}^G\sigma_1[1/\HeckeX]\cong \Ind_{\lbar B}^G\left ( \cInd_{T_0 Z}^T\chi_{1}
  \right ).\end{equation}
We now use some results of Heyer~\cite{Heyerderived, Heyergeometric}. 
By~\cite[Theorem~4.1.1]{Heyerderived}, the functor $\Ind_{\lbar B}^G : D(\sm.\, T) \to D(\sm.\, G)$ of parabolic
induction has a left adjoint $L(\lbar U,
\text{--})$. By~\cite[Corollary~4.1.3]{Heyerderived} there is a spectral sequence  of
$\cH_T[1/f]$-modules %
\begin{multline}\label{left adjoint spectral sequence}
E_2^{i, j} = \Ext^i_{\sm.\, T}(L^{-j}(\lbar U, \Ind_{\lbar B}^G(\cInd_{T_0
  Z}^T\chi_1)), (\cInd_{T_0 Z}^T\chi_2)[1/f]) \\ \Rightarrow \Ext^{i+j}_{\sm.\, G}\left (\Ind_{\lbar B}^G(\cInd_{T_0
  Z}^T\chi_1), \Ind_{\lbar B}^G\left ((\cInd_{T_0 Z}^T\chi_2)[1/f] \right ) \right ).
\end{multline}
As previously recalled, the target of this spectral sequence is isomorphic to
\[\Ext^{i+j}_{\sm.\, G}(\cInd_{KZ}^G(\sigma_1)[1/\HeckeX], \cInd_{KZ}^G(\sigma_2)[1/f]).\]
(In fact, this is~\eqref{eqn:instance-of-universal-parabolic-induction} for the first term, and~\eqref{identifying localization} and~\eqref{universal parabolic induction II} for the second term.)

We now compute the~$E_2$ page of~\eqref{left adjoint spectral sequence}.
For all~$w \in W(G, T)$, the intersection $T \cap w\lbar U w^{-1}$ is trivial.
Using this, the statement and proof of %
~\cite[Example~4.2.1]{Heyergeometric} go through unchanged to
show that $L(\lbar U, \Ind_{\lbar B}^G ( \cInd_{T_0 Z}^T\chi_1))$ is
concentrated in degrees~$[-1, 0]$, and that there are natural
isomorphisms %
\begin{equation}\label{eqn:explicit-description-of-Li-functors}
\begin{gathered}
L^0(\lbar U, \Ind_{\lbar B}^G(\cInd_{T_0 Z}^T\chi_1)) \cong \cInd_{T_0 Z}^T\chi_1,\\
L^{-1}(\lbar U, \Ind_{\lbar B}^G (\cInd_{T_0 Z}^T\chi_1) ) \cong \delta_s
\otimes (\cInd_{T_0 Z}^T \operatorname{ad}(s)^* \chi_1),
\end{gathered}
\end{equation}
where $s \in W(G, T)$ is the nontrivial element, and $\delta_s = \omega^{-1} \otimes \omega$ is the $\lbar B$-positive root. Here we have used~\cite[Lemma~4.1.6]{Heyergeometric}
to compute~$\delta_s$, and we have also used the natural isomorphism of objects of~$\sm.\, T$
\[
\ad(s)^*(\cInd_{T_0 Z}^T)(\chi_1) \cong \cInd_{T_0 Z}^T(\ad(s)^*(\chi_1)).
\]

We now analyse the spectral sequence~\eqref{left adjoint spectral sequence}. By~\eqref{eqn:explicit-description-of-Li-functors} and
Lemma~\ref{lem:Ext-groups-over-T} below, the $E_2^{i,j}$-terms are finite free
$\cH_T[1/f]$-modules, so the $\Ext$-groups \begin{equation}\label{eqn:Ext-groups-without-fixed-central-character}\Ext^{i}_{\sm.\,
    G}(\cInd_{KZ}^G(\sigma_1)[1/\HeckeX], \cInd_{KZ}^G(\sigma_2)[1/f])\end{equation} are
finitely generated $\cH_T[1/f]$-modules.
By~\eqref{reformulating E^i} and
Lemma~\ref{Koszul central character}, the same is true of the
$\Ext^i_{\cA}(j_*j^* \cInd_{KZ}^G\sigma_1, j_*j^*
\cInd_{KZ}^G\sigma_2)$.
Furthermore, we see from Lemma~\ref{Koszul central
  character} that if some $\Ext^j_{\cA}(j_*j^* \cInd_{KZ}^G\sigma_1, j_*j^*
\cInd_{KZ}^G\sigma_2)$ has $\cH_T[1/f]$-torsion, then the same is true of the
groups~\eqref{eqn:Ext-groups-without-fixed-central-character} for~$j=i,i+1,i+2$.
However, the vanishing result of Lemma~\ref{lem:Ext-groups-over-T} shows that  ~\eqref{left adjoint spectral sequence} degenerates at~$E_3$, %
and since all the objects on the~$E_2$ page are free $\cH_T[1/f]$-modules,  %
it  furthermore shows that the groups~\eqref{eqn:Ext-groups-without-fixed-central-character} could
only have non-zero torsion in degrees~$2$ or~$3$. In particular, it is impossible
for them to have non-zero torsion for a range of degrees ~$j=i,i+1,i+2$, so they
are torsion-free, as required.
\end{proof}

\begin{lem}\label{lem:Ext-groups-over-T}For any $\zeta$-compatible smooth character $\psi: T_0 \to \bF^\times$, the
  $\Ext$-groups \[\Ext^i_{\sm.\,T}(\cInd_{T_0 Z}^T\psi,\cInd_{T_0 Z}^T
  \chi_2 [1/f])\]  are finite free ~$\cH_T[1/f]$-modules, and are zero
for~$i\ge 4$.
\end{lem}%
\begin{proof}
  By Frobenius reciprocity, we have a natural isomorphism
\[\Ext^i_{\sm.\,T}(\cInd_{T_0 Z}^T\psi ,\cInd_{T_0 Z}^T
  \chi_2 [1/f])\cong \Ext^i_{\sm.\ T_0 Z}(\psi, \Res^T_{T_0 Z}\cInd_{T_0 Z}^T
  \chi_2 [1/f]).\]
 Since the ring $\cO\llbracket T_0Z\rrbracket $ is Noetherian (e.g.\ by
~\cite[Section~2.1]{MR4551876}), the category $\sm.\, T_0 Z$ is locally Noetherian
(e.g.\ because it is a localizing subcategory of
$\Mod \cO\llbracket T_0 Z\rrbracket $, by Lemma~\ref{embedding of smooth in modules}). Thus $\Ext^i_{\sm.\, T_0 Z}(\psi, \text{--})$
commutes with filtered colimits (e.g.\ by~
\cite[Proposition~A.1.1(3)]{DEGlocalization}). 
On the other hand, we have an isomorphism of $\cH_T[1/f]$-modules
\[
\Res^T_{T_0 Z} \cInd_{T_0 Z}^T \chi_2 [1/f] \cong \chi_2  \otimes_{\bF} \cH_T[1/f],
\]
and thus an identification
\[
\Ext^i_{\sm.\,T}(\cInd_{T_0 Z}^T\psi ,\cInd_{T_0 Z}^T
  \chi_2 [1/f])\cong \Ext^i_{\sm.\, T_0 Z}(\psi, \chi_2 ) \otimes_{\bF} \cH_T[1/f], 
\]
which shows that the $\Ext$-groups are finite free over~$\cH_T[1/f]$. 
Since~$T_0Z \cong \bQ_p^\times
\times \bZ_p^\times$, these~$\Ext$ groups are furthermore
concentrated in degrees~$i \in [0, 3]$, as required.
\end{proof}

\begin{rem}\label{two Hecke actions}
Assume that~$\chi_1 \ne \delta_s \otimes \ad(s)^*\chi_1$.
There are two actions of~$\cH_T$ on $\Ext^i_{\cA}(j_*j^* \cInd_{KZ}^G\sigma_1, j_*j^* \cInd_{KZ}^G\sigma_2)$, arising from the action of~$\cH_G(\sigma_i)$ on the two factors.
The proof of Proposition~\ref{properties of Hecke modules} shows that they coincide when~$\chi_2 = \chi_1$, and are related by the automorphism $\HeckeX \mapsto \HeckeX^{- 1}$ of~$\cH_T$ when~$\chi_2 = \delta_s\otimes\operatorname{ad}(s)^*\chi_1 $.
\end{rem}

\subsubsection{\texorpdfstring{$\Ext$}{Ext} and completion.}%
We now show the compatibility with completion of the formation of
$\Ext$-modules.
Recall that we have fixed a companion pair of Serre weights~$\{\sigma,\sigmacomp\}$ and a choice
 $\sigma_1=\sigma$, $ \sigma_2 =\sigmacomp$. %
We continue to write~$j_*, j^*$ for~$j_{U(\sigmasigmacomp),*}, j_{U(\sigmasigmacomp)}^*$.
Fix a closed point $\thetabar$ of~$U(\sigmasigmacomp)$, and let~$f_{\thetabar} \in \cH_G(\sigma_2)$ be the
squarefree monic polynomial
associated to~$\{\thetabar\} \subset X$ in Definition~\ref{f_Y}. %

\begin{lemma}\label{A.3.14 revisited}
Let $\pi_i \coloneqq  \cInd_{KZ}^G\sigma_i$.
Then the natural map
\begin{equation}\label{eqn:natural map to limit of Ext}
\Ext^i_{\Pro D^b_{\fp}(\cA_{\thetabar})}(\ihat^*_{\thetabar}(\pi_1), \ihat^*_{\thetabar}(\pi_2)) \to 
\varprojlim_n \Ext^i_{\Pro D^b_{\fp}(\cA_{\thetabar})}(\ihat^*_{\thetabar}(\pi_1), \pi_2/f_{\thetabar}^n\pi_2)
\end{equation}%
is an isomorphism.
\end{lemma}
\begin{proof}%
By~\eqref{eqn:recollection-of-explicit-completion-of-cInd-sigma}, we have
\[
\ihat^*_{\thetabar}(\pi_2) = \quoteslim{n} \pi_2/f_{\thetabar}^n\pi_2 \in \Pro(\cA_{\thetabar}^{\fp}).
\] 
By Lemma~\ref{lem:formal properties of RHom}~(1), applied to the compactly generated
$\infty$-category $\Ind \Pro D^b_{\fp}(\cA_{\thetabar})$, 
we have natural isomorphisms %
\begin{multline*}
\RHom_{\Pro D^b_{\fp}(\cA_{\thetabar})}(\ihat^*_{\thetabar}(\pi_1), \ihat^*_{\thetabar}(\pi_2)) = \RHom_{\Ind\Pro D^b_{\fp}(\cA_{\thetabar})}(\ihat^*_{\thetabar}(\pi_1), \ihat^*_{\thetabar}(\pi_2)) \\ 
= \lim_n\RHom_{\Ind\Pro D^b_{\fp}(\cA_{\thetabar})}(\ihat^*_{\thetabar}(\pi_1), \pi_2/f^n_{\thetabar}\pi_2) = \lim_n\RHom_{\Pro D^b_{\fp}(\cA_{\thetabar})}(\ihat^*_{\thetabar}(\pi_1), \pi_2/f^n_{\thetabar}\pi_2)
\end{multline*}
(where~$\lim_n$ is formed in the derived $\infty$-category $D(\F)$ of $\F$-vector spaces, and we have used that $\Pro D^b_{\fp}(\cA_{\thetabar}) \to \Ind \Pro D^b_{\fp}(\cA_{\thetabar})$
is fully faithful and limit-preserving). %
The map~\eqref{eqn:natural map to limit of Ext}
is thus an edge map in a %
spectral sequence
\begin{equation}\label{eqn:Rlim spectral sequence}
E_2^{p, q} = R^p\lim_n\Ext^q_{\Pro D^b_{\fp}(\cA_{\thetabar})}(\ihat^*_{\thetabar}(\pi_1), \pi_2/f^n_{\thetabar}\pi_2) \Rightarrow 
\Ext^{p+q}_{\Pro D^b_{\fp}(\cA_{\thetabar})}(\ihat^*_{\thetabar}(\pi_1), \ihat_{\thetabar}^*(\pi_2)).
\end{equation}
We now claim that, for all $\tau \in \cA_{\thetabar}^{\fp}$ and~$i \in \bZ$, the $\F$-vector space
\begin{equation}\label{eqn:need-to-prove-fd}
\Ext^i_{\Pro D^b_{\fp}(\cA_{\thetabar})}(\ihat^*_{\thetabar}(\pi_1), \tau)
\end{equation}
is finite-dimensional.
Assuming the claim, the Mittag-Leffler criterion shows that~$E_2^{p, q} = 0$ whenever~$p \geq 1$, and so~\eqref{eqn:Rlim spectral sequence} degenerates at~$E_2$, 
and~\eqref{eqn:natural map to limit of Ext} is an isomorphism, as desired.

We now prove the claim.
By Lemma~\ref{lem:derived-t-exact-ihat-Y}, the functor %
\[
\ihat^*_{\thetabar} : \Pro D^b_{\fp}(\cA) \to \Pro D^b_{\fp}(\cA_{\thetabar})
\]
is left adjoint to the inclusion, and so~\eqref{eqn:need-to-prove-fd}
is isomorphic to
\[
\Ext^i_{\Pro D^b_{\fp}(\cA)}(\pi_1, \tau) = \Ext^i_{D^b_{\fp}(\cA)}(\pi_1, \tau).
\]
In turn, since~$\cA = \Ind(\cA^{\fp})$, the equivalence~\eqref{eqn:Yoneda equivalence} induces an isomorphism
\[
\Ext^i_{D^b(\cA)}(\pi_1, \tau) \isoto \Ext^i_{\Ind D^b_{\fp}(\cA)}(\pi_1, \tau) = \Ext^i_{D^b_{\fp}(\cA)}(\pi_1, \tau),
\] 
and so we have reduced to proving that $\Ext^i_{D^b(\cA)}(\pi_1, \tau)$ has finite $\F$-dimension.
By d\'evissage, we can assume that~$\tau$ is irreducible. %
Then the claim is a consequence
of~\cite[Proposition~4.2.4]{DEGlocalization}. 
\end{proof}

Continue to write $\pi_i \coloneqq  \cInd_{KZ}^G\sigma_i$. By Remark~\ref{rem:completion-to-subset-of-open}, the  exact functor $\ihat_{\thetabar}^*$ induces a map
\begin{equation*}
\Ext^i_{D^b_{\fp}(\cA_U)}(j^*\pi_1, j^*\pi_2) \to \Ext^i_{D^b(\Pro\cA^{\fp}_{\thetabar})}(\ihat_{\thetabar}^*\pi_1, \ihat_{\thetabar}^*\pi_2)
\end{equation*}
which can be composed with the $t$-exact functor $p : D^b(\Pro \cA^{\fp}_{\thetabar}) \to \Pro D^b(\cA^{\fp}_{\thetabar})$ of~\eqref{eqn:bounded pro functor}
to produce a map
\begin{equation}\label{eqn: Ext and completion}
  \ihat^*_{\thetabar}: \Ext^i_{D^b_{\fp}(\cA_U)}(j^*\pi_1, j^*\pi_2) \to \Ext^i_{\Pro D^b(\cA^{\fp}_{\thetabar})}(\ihat_{\thetabar}^*\pi_1, \ihat_{\thetabar}^*\pi_2)
\end{equation}

\begin{prop}\label{Ext and completion}Let $\pi_i \coloneqq  \cInd_{KZ}^G\sigma_i$. Then
  ~\eqref{eqn: Ext and completion} 
induces an isomorphism
\begin{equation}\label{eqn: Ext and completion II}
\Ext^i_{D^b_{\fp}(\cA_U)}(j^*\pi_1, j^*\pi_2)^{\wedge}_{f_{\thetabar}}
 \isom \Ext^i_{\Pro D^b(\cA_{\thetabar}^{\fp})}(\ihat_{\thetabar}^*\pi_1, \ihat_{\thetabar}^*\pi_2).
\end{equation}
\end{prop}
\begin{proof}

By Remark~\ref{rem:completion-to-subset-of-open}, as well as the fact that $j_*$ is fully faithful and $t$-exact, 
it suffices to prove that~$\ihat^*_{\thetabar}$ induces an isomorphism
\begin{equation}\label{eqn: Ext and completion III}
\Ext^i_{D^b(\cA)}(j_*j^*\pi_1, j_*j^*\pi_2)^{\wedge}_{f_{\thetabar}}
 \isom \Ext^i_{\Pro D^b(\cA_{\thetabar}^{\fp})}(\ihat_{\thetabar}^*\pi_1, \ihat_{\thetabar}^*\pi_2).
\end{equation}
By Proposition~\ref{properties of Hecke modules}, the group
$\Ext^i_{D^b(\cA)}(j_*j^*\pi_1, j_*j^*\pi_2)$ is $\cH_T[1/f]$-torsion free and in
particular $f_{\thetabar}$-torsion free, so that
for any~$n\ge 1$,
the exact sequence
\begin{equation}\label{completing pi_2}
0 \to j_*j^*\pi_2 \xrightarrow{f_{\thetabar}^n} j_*j^*\pi_2 \to \pi_2/f_{\thetabar}^n\pi_2 \to 0
\end{equation}
induces an isomorphism
\begin{equation}\label{eqn:another-iso-of-Exts-to-complete}
\Ext^i_{D^b(\cA)}(j_*j^*\pi_1, j_*j^*\pi_2)/f^n_{\thetabar} \cong \Ext^i_{D^b(\cA)}(j_*j^*\pi_1, \pi_2/f^n_{\thetabar}\pi_2).
\end{equation}
By Lemma~\ref{lem:derived-t-exact-ihat-Y},
we see that $\ihat_{\thetabar}^*$ induces an isomorphism
\begin{equation*}
\Ext^i_{D^b_{\fp}(\cA)}(\pi_1, \pi_2/f^n_{\thetabar}\pi_2) \isoto 
\Ext^i_{\Pro D^b(\cA^{\fp}_{\thetabar})}(\ihat_{\thetabar}^*\pi_1, \pi_2/f_{\thetabar}^n\pi_2),
\end{equation*}
hence (bearing in mind Remark~\ref{rem:completion-to-subset-of-open}) an isomorphism
\begin{equation}\label{completing pi_2 II}
  \ihat^*_{\thetabar}: \Ext^i_{D^b(\cA)}(j_*j^*\pi_1, \pi_2/f^n_{\thetabar}\pi_2) \isoto 
\Ext^i_{\Pro D^b(\cA^{\fp}_{\thetabar})}(\ihat_{\thetabar}^*\pi_1, \pi_2/f_{\thetabar}^n\pi_2).
\end{equation}

On the other hand, applying~$\Ext^i_{\Pro D^b(\cA^{\fp}_{\thetabar})}(\ihat^*_{\thetabar}\pi_1, \text{--})$ to the exact sequence
\[
0 \to \ihat^*_{\thetabar}\pi_2 \xrightarrow{f_{\thetabar}^n} \ihat^*_{\thetabar}\pi_2 \to \pi_2/f_{\thetabar}^n\pi_2 \to 0
\]
obtained by applying~$\ihat^*_{\thetabar}$ to~\eqref{completing pi_2}, %
we obtain an
injection 
\begin{equation}\label{eqn:
      Ext and completion-again}
      \Ext^i_{\Pro D^b(\cA^{\fp}_{\thetabar})}(\ihat^*_{\thetabar}\pi_1,
  \ihat^*_{\thetabar}\pi_2)/f^n_{\thetabar}\into \Ext^i_{\Pro D^b(\cA^{\fp}_{\thetabar})}(\ihat^*_{\thetabar}\pi_1, \pi_2/f^n_{\thetabar}\pi_2).\end{equation}
We conclude from the discussion above that the composition of~\eqref{eqn:another-iso-of-Exts-to-complete} and~\eqref{completing pi_2 II} is an isomorphism
\[
\Ext^i_{D^b(\cA)}(j_*j^*\pi_1, j_*j^*\pi_2)/f^n_{\thetabar} \isom \Ext^i_{\Pro D^b(\cA^{\fp}_{\thetabar})}(\ihat_{\thetabar}^*\pi_1, \pi_2/f_{\thetabar}^n\pi_2),
\]
which can also be factored as a composition
\begin{multline*}
  \Ext^i_{D^b(\cA)}(j_*j^*\pi_1, j_*j^*\pi_2)/f^n_{\thetabar} \stackrel{\ihat^*_{\thetabar}}{\to}
  \Ext^i_{\Pro D^b(\cA^{\fp}_{\thetabar})}(\ihat^*_{\thetabar}\pi_1,
  \ihat^*_{\thetabar}\pi_2)/f^n_{\thetabar}\\
  \stackrel{\eqref{eqn:
      Ext and completion-again}}\into \Ext^i_{\Pro D^b(\cA^{\fp}_{\thetabar})}(\ihat^*_{\thetabar}\pi_1, \pi_2/f^n_{\thetabar}\pi_2),
\end{multline*}
It follows that both arrows in this composition are isomorphisms, and so 
\[
\ihat^*_{\thetabar} : \Ext^i_{D^b(\cA)}(j_*j^*\pi_1, j_*j^*\pi_2) \to
  \Ext^i_{\Pro D^b(\cA^{\fp}_{\thetabar})}(\ihat^*_{\thetabar}\pi_1, \ihat^*_{\thetabar}\pi_2)
\]
becomes an isomorphism after quotienting out by~$f_{\thetabar}^n$.
Hence, to conclude the proof that~\eqref{eqn: Ext and completion III} is an isomorphism, it suffices to prove that the natural map
\[
\Ext^i_{\Pro D^b(\cA^{\fp}_{\thetabar})}(\ihat_{\thetabar}^*\pi_1, \ihat_{\thetabar}^*\pi_2) \to 
\varprojlim_n \Ext^i_{\Pro D^b(\cA^{\fp}_{\thetabar})}(\ihat^*_{\thetabar}\pi_1, \ihat^*_{\thetabar}\pi_2)/f^n_{\thetabar}
\]
is an isomorphism.
Since~\eqref{eqn:
      Ext and completion-again} is an isomorphism, the right-hand side is isomorphic to
\[\varprojlim_n \Ext^i_{\Pro D^b(\cA^{\fp}_{\thetabar})}(\ihat^*_{\thetabar}\pi_1, \pi_2/f^n_{\thetabar}\pi_2).\]
So it suffices to prove that
\[
\Ext^i_{\Pro D^b(\cA^{\fp}_{\thetabar})}(\ihat_{\thetabar}^*\pi_1, \ihat_{\thetabar}^*\pi_2) \to \varprojlim_n \Ext^i_{\Pro D^b(\cA^{\fp}_{\thetabar})}(\ihat^*_{\thetabar}\pi_1, \pi_2/f^n_{\thetabar}\pi_2)
\]
is an isomorphism. This is 
Lemma~\ref{A.3.14 revisited}.
\end{proof}

\begin{rem}\label{properties of Hecke modules sigma version 2}%
Proposition~\ref{properties of Hecke modules} (and the fact that $\cH_T[1/f]$
is a principal ideal domain) shows that the  %
$\F[\HeckeX^{\pm
  1}, f(\HeckeX)^{-1}]$-module \[\Ext^i_{\cA}\bigl( j_*j^*\cInd_{KZ}^G \sigma_1,
j_*j^*\cInd_{KZ}^G \sigma_2\bigr)\] is finite free. 
A further analysis of the spectral sequence~\eqref{left adjoint spectral
  sequence}  then allows one to determine its rank:
\begin{itemize}
\item If $i=0$ then the rank is~$1$ if ~$\sigma_1 =\sigma_2 $, and zero
  otherwise.
\item If ~$i=1$ then the rank is~$1$ unless~$\{\sigma,\sigmacomp\}$ is of 
  type~\ref{item: non p-distinguished pseudorep} in which case the rank is~$2$.
\item If $i=2$ then the rank is~$0$ if  $\sigma_1 =\sigma_2$ and
  $\{\sigma,\sigmacomp\}$ is not of 
  type~\ref{item: non p-distinguished pseudorep}, and otherwise the rank is~$1$.
\item If $i=3$ then the rank is~$0$.
\end{itemize}
\end{rem}

\section{The moduli stack of rank~2 \'etale \texorpdfstring{$(\varphi,\Gamma)$}{(ϕ,Γ)}-modules}\label{sec: moduli stacks}%
In this section we will recall some of the main results
of~\cite{emertongeepicture,EGaddenda}, specialised to the case of~$\GL_2(\bQ_p)$,
and prove some additional results in this setting. %

\subsection{\'Etale \texorpdfstring{$(\varphi,\Gamma)$}{(ϕ,Γ)}-modules}\label{subsec: the coefficient
  rings}We begin by recalling some definitions and notation
from~\cite{emertongeepicture}.
       We write
       $\Gamma\coloneqq  \Gal(\Qp(\zeta_{p^{\infty}})/\Qp)$, so that the
       cyclotomic character induces an isomorphism
       $\chi:\Gamma \isoto \Z_p^{\times}$. For each $a\in\Zptimes$, we
       write~$\sigma_a\in\Gamma$ for the element with $\chi(\sigma_a)=a$. %
       If we choose a compatible system of $p^n$-th roots of $1$,
       then these give rise in the usual way to an element
       $\varepsilon \in (\widehat{\Qp(\zeta_{p^{\infty}})})^{\flat}.$
       As usual, let $[\varepsilon]$ denote the Teichm\"uller 
       lift of $\varepsilon$ to an element of
       $W\bigl( \cO_{\widehat{\Qp(\zeta_{p^{\infty}})}}^\flat).$
       There is then a continuous embedding
       $$\Z_p\llbracket T\rrbracket  \hookrightarrow 
       W( \cO_{\widehat{\Qp(\zeta_{p^{\infty}})}}^\flat)$$
       (the source being endowed with its $(p,T)$-adic topology,
       and the target with its weak topology),
       defined via $T \mapsto [\varepsilon]- 1$.
       We denote the image of this embedding by $\A^+$.
       This embedding extends to an embedding
       $$\widehat{\Z_p((T))} \hookrightarrow W\bigl( (\widehat{\Qp(\zeta_{p^{\infty}})})^\flat
       \bigr)$$
       (here the source is the $p$-adic completion of the Laurent
       series ring $\Z_p((T))$),
       whose image we denote by $\A$. %

       We have commuting actions of~$\varphi$ and~$\Gamma$ on~$\A$
which are given by the explicit formulae
\[\varphi(T)=(1+T)^p-1,\]
\[\sigma_a(T)=(1+T)^{a}-1.\]
Note that $\A^+$ is
visibly~$(\varphi,\Gamma)$-stable.

There is a left inverse~$\psi$ to~$\varphi$ defined as follows. We
have a decomposition~$\A=\oplus_{i=0}^{p-1}(1+T)^i\varphi(\A)$,
and~$\psi:\A\to\A$ is defined to be the projection onto the $i=0$
factor; so in particular~$\psi(\varphi(x))=x$. This restricts to a
surjection~$\psi:\A^+\to\A^+$.

\begin{rem}
  \label{rem: change of notation and meaning versus EG paper}Our
  notation and conventions differ from that
  of~\cite{emertongeepicture}, for two reasons. Firstly, since we are
  only working over~$\Qp$, and not over an extension $K/\Qp$, we have
  dropped~$K$ from the notation. Secondly, and more significantly, we
  use $(\varphi,\Gamma)$-modules for the full cyclotomic extension
  $\Qp(\zeta_{p^{\infty}})/\Qp$, and we make no use of
  $(\varphi,\Gamma)$-modules for the $\Zp$-subextension. %

  The theories of $(\varphi,\Gamma)$-modules for the two extensions
  are equivalent (via taking invariants for the group of roots of
  unity in~$\Zp^\times$), but in most of~\cite{emertongeepicture} the
  $\Zp$-theory was used, because it was convenient when proving the
  existence of various moduli stacks to work with a pro-cyclic pro-$p$ group.
However, for applications to the $p$-adic Langlands program
  (or more generally in arguments involving explicit formulas for the
  actions of $\varphi$ and~$\Gamma$), it is much more convenient to
  work with the full cyclotomic extension. Since we make no use of the $\Zp$-subextension, and we want to avoid
  notational clutter where possible, %
  we have written~$\Gamma$
  for the group denoted~$\Gammat$ in~\cite{emertongeepicture}, and  $\A,\A^+$ for the rings denoted~$\A'_{\Qp}$ and $(\A')^+_{\Qp}$
  respectively in~\cite{emertongeepicture}.
\end{rem}

We now introduce coefficients. Recall that~$\cO$ is the ring of integers in our fixed finite
extension $E/\Qp$. Let~$A$ be a $p$-adically complete $\cO$-algebra.
Frequently, we will work modulo a fixed power $\varpi^a$ of $\varpi$,
and thus assume that $A$ is actually an $\cO/\varpi^a$-algebra
(and sometimes we impose further conditions on $A$,
such as that of being Noetherian, or even of finite type over
$\cO/\varpi^a$). %
We write
\[\A^+_{A}=\A^+\cotimes_{\Zp}A \coloneqq  \varprojlim_{n} \A^+/(p,T)^n
	\otimes_{\Z_p} A = \varprojlim_{m} (\varprojlim_n \A^+/(p^m,T^n)
	\otimes_{\Z_p} A), \]
and
\[\A_{A} =\A\cotimes_{\Zp}A \coloneqq 
	\varprojlim_{m} ((\varprojlim_n \A^+/(p^m,T^n)
	\otimes_{\Z_p} A) [1/T]). \]
      If~$A$ is a Noetherian $\cO/\varpi^{a}$-algebra for some $a\ge 1$,
      then $\A^+_A=A\llbracket T\rrbracket $ and $\A_A=A((T))$ are flat $A$-algebras, and
      the morphism $\A^+_A\to\A_A$ is flat and injective. If~$A\to B$
      is a (faithfully) flat homomorphism of Noetherian
      $\cO/\varpi^a$-algebras, then the morphisms $\A^+_A\to\A^+_B$ and
      $\A_A\to\A_B$ are (faithfully) flat. (All of these statements
      follow easily from the facts that completions and localizations
      of Noetherian rings are flat, see~\cite[Lem.\
      5.1.7]{EGstacktheoreticimages}.)

      The rings $\A^+_A$, $\A_A$ have natural topologies, which are
      discussed further in Section~\ref{subsec: Tate module}. By~\cite[Lem.\ 2.2.17]{emertongeepicture}, the actions
      of~$\varphi,\Gamma$ on~$\A^+$ and~$\A$ extend to continuous
      $A$-linear actions on~$\A^+_A$ and~$\A_A$. Similarly, the action
      of~$\psi$ extends to continuous $A$-linear surjections
      $\psi:\A_A\to\A_A$, $\psi:\A^+_A\to\A^+_A$, satisfying $\psi(\varphi(x))=x$. %

  By
  definition, a projective \'etale $(\varphi,\Gamma)$-module with $A$-coefficients
  is a finitely generated projective $\A_A$-module $D$, equipped with
\begin{itemize}
\item a $\varphi$-linear morphism $\varphi:D\to D$ with the property
  that the corresponding morphism $\Phi_D:\varphi^*D\to D$ is an
  isomorphism (i.e.\ $D$ is given the structure of an \'etale
  $\varphi$-module over $\A_A$), and 
\item a continuous semi-linear action
  of~$\Gamma$ that commutes with~$\varphi$.
\end{itemize}Here the notion of continuity is with respect to the
canonical topology on~$D$ inherited from that on~$\A_A$ 
(see~\cite[Rem.\ D.2]{emertongeepicture}). Our \'etale
$(\varphi,\Gamma)$-modules will usually be assumed to be projective, and we
will often write ``\'etale $(\varphi,\Gamma)$-module'' for
``projective \'etale
$(\varphi,\Gamma)$-module'' when no confusion should arise.
We will also sometimes write ``\'etale $(\varphi,\Gamma)$-module over~$A$'' for ``\'etale $(\varphi,\Gamma)$-module with $A$-coefficients''.

\begin{rem}\label{rem:Fontaine equivalence}
In the case that~$A$ is Artinian, there is the usual equivalence of
categories (see \cite{MR1106901}) between the category of (projective) \'etale
$(\varphi,\Gamma)$-module with $A$-coefficients and the category of
finite projective $A$-modules with an action of~$G_{\Qp}$, given by
the functor $D\mapsto
(W(\C^\flat)\otimes_{\A}D)^{\varphi=1}$.
\end{rem}

\subsection{Moduli stacks of \'etale \texorpdfstring{$(\varphi,\Gamma)$}{(ϕ,Γ)}-modules}\label{subsec: moduli stacks}

\begin{defn}\label{defn:stack-cX-fixed-det-2D}
  We %
  write $\cX$ for the stack of projective \'etale %
  $(\varphi,\Gamma)$-modules of rank $2$ with
  determinant~$\zeta\varepsilon^{-1}$, as defined
  in %
  Definition~\ref{defn: Xdchi}.  Explicitly, if~$A$ is a $p$-adically
  complete $\cO$-algebra, then $\cX(\Spf A)$
  is the groupoid of pairs $(D,\theta)$ where~$D$ is a rank~$2$
  projective \'etale $(\varphi,\Gamma)$-module with $A$-coefficients,
  and $\theta$ is an identification of~$\wedge^2D$ with (the
  $(\varphi,\Gamma)$-module corresponding to)~$\zeta\varepsilon^{-1}$.
\end{defn}

Theorem~\ref{thm: fixed determinant stack for GL2 Qp} below summarises the basic
properties of~$\cX$. In order to state it, we make the following definitions. Our
notation for Serre weights is as in Section~\ref{subsec: notation and
  conventions}; recall in particular that we assume throughout this paper that our Serre
weights~$\sigma$ are
compatible with~$\zeta$ %
in the sense of Definition~\ref{defn: compatible Serre weight}.

\begin{defn}
  \label{defn:potentially crystalline EG stacks}
  Let $\lambda \in \{(a, b) \in \bZ^2: a > b\}$ be a regular Hodge type, and let
  $\tau: I_{\bQ_p} \to \GL_2(E)$ be an inertial type (i.e.\ a
representation with open kernel that extends to the Weil group of~$\bQ_p$) such that the pair~$(\lambda, \tau)$ is compatible with~$\zeta \varepsilon^{-1}$, in the sense of Definition~\ref{adefn: Hodge type compatible with chi }.
 Then we write~$\cX^{\lambda, \tau}$ for the $\varpi$-adic formal algebraic substack of~$\cX$ denoted $\cX_2^{\crys, \uline \lambda, \tau, \zeta\varepsilon^{-1}}$ 
  in~\eqref{eqn:defn-of-crys-lambda-tau-chi}.
 If~$\tau$ is the trivial representation, we will often write $\cX^{\lambda, \crys}$ for~$\cX^{\lambda, \tau}$.
\end{defn}

\begin{defn}
  \label{defn:special-fibre-crystalline-weight-sigma}If~$\sigma=\sigma_{a,b}$ is a Serre
  weight, then we 
  define 
  \[
  \cZ(\sigma)
  \coloneqq (\cX^{(1-a, -a-b), \crys}\times_{\cO}\F)_{\red}.
  \] 
This is a closed algebraic substack of~$\cX$ (since
$\cX^{(1-a,-a-b),\crys}$ is $\varpi$-adically formal).
\end{defn}

The next theorem collects some of the basic properties of the moduli stacks we have just introduced.
Recall that by~\cite[Theorem~6.6.3]{emertongeepicture}, the finite type points of~$\cX$ are in bijection
with $\Gal(\cbF_p/\F)$-conjugacy classes of continuous representations
$\rhobar: G_{\bQ_p} \to \GL_2(\cbF_p)$ with determinant~$\zeta\varepsilon^{-1}$. 
This bijection sends~$\rhobar$ to the image of the classifying map $\Spec \cbF_p \to \cX$ of the $(\varphi, \Gamma)$-module corresponding to~$\rhobar$ under the equivalence of Remark~\ref{rem:Fontaine equivalence}.

\begin{thm}%
  \label{thm: fixed determinant stack for GL2 Qp}\leavevmode
  \begin{enumerate}
  \item $\cX$ is a Noetherian
  formal algebraic stack. 
\item The underlying reduced substack $\cX_{\red}$
  is of finite type over~$\F$, and is equidimensional of
  dimension~$1$.
\item   The irreducible components of~$\cX_{\red}$ admit a
  natural surjection to the set of $\zeta$-compatible Serre weights. 
  The fibres of this surjection are single irreducible
  components~$\cX(\sigma)$ unless $\sigma=\sigma_{a,p-1}$ for some~$a$, in
  which case we have two irreducible components~$\cX(\sigma)^{\pm}$.
\item\label{item:cycle-of-Z(sigma)} 
  The stack~$\cZ(\sigma)$ is algebraic, and coincides with~$\cX(\sigma)$ unless $\sigma=\sigma_{a,p-1}$ for some~$a$, in
  which case it is the scheme-theoretic union of~$\cX(\sigma_{a,0})$ and the two irreducible components~$\cX(\sigma)^{\pm}$. %
  \item\label{item:points and Serre weights}  Let $x \in |\cX|$ be a finite type point, with corresponding representation $\rhobar : G_{\bQ_p} \to \GL_2(\cbF_p)$.
    Then $x \in |\cZ(\sigma)|$ if and only if~$\sigma$ is a Serre weight of~$\rhobar$. %
\end{enumerate}%
  \end{thm}
\begin{proof}
  The first three parts are a special case of
  Theorem~\ref{athm: fixed determinant stack in general}. 
  Part~\eqref{item:cycle-of-Z(sigma)} is immediate from~\cite[Thm.\ 8.6.2]{emertongeepicture}.
  Finally, Definition~\ref{defn:Serre weights of rhobar} was formulated in such a way that part~\eqref{item:points and Serre weights} is true by definition;
  see e.g.\ \cite[\S 8.5]{emertongeepicture} %
  for an elaboration of this point. 
\end{proof}

We now describe a choice of~$(\lambda, \tau)$ so that $\cX^{\lambda, \tau} \otimes_{\cO} \bF = \cZ(\sigma)$. %

\begin{lemma}\label{Z(sigma) as special fibre}
Let~$\sigma$ be a Serre weight.
\begin{enumerate}%
\item If $\sigma = \Sym^b \otimes \det^a$, and $b \ne p-2$, we let $\lambda = (1-a, -a-b)$, and we let~$\tau$ be trivial.
\item If $\sigma = \Sym^{p-2} \otimes \det^a$, 
we let
$\lambda = (1, 0)$, and we let~$\tau = \omega_2^{p + (a-1)(p+1)} \oplus \omega_2^{p^2 + (a-1)(p+1)}$. 
\end{enumerate}
Then the special fibre $\cX^{\lambda, \tau} \otimes_{\cO} \bF$ is reduced, and $\cX^{\lambda, \tau} \otimes_{\cO} \bF = \cZ(\sigma)$.
\end{lemma}%
\begin{proof}
The top-dimensional cycle of $\cX^{\lambda, \tau}\otimes_{\cO}\F$ is always~$[\cZ(\sigma)]$; this is immediate from the definitions in case~(1), and  from  the geometric 
Breuil--M\'ezard conjecture~\cite[Theorem~8.6.2]{emertongeepicture} in case~(2).
(In case~(2), note that~$\tau$ is the inertial type corresponding to the representation~$\Theta(\sigmasigmacomp)$ of Definition~\ref{defn: companion Serre weights}~(2), whose
semisimplified mod $\varpi$ reduction is~$\sigma$.)
It thus suffices to show that $\cX^{\lambda, \tau}\otimes_{\cO}\F$ is reduced. 

In case~(1), if furthermore $b \ne p-1$, the reducedness follows from Fontaine--Laffaille theory.
If $b = p-1$, it suffices to prove that for any semisimple $\rhobar: G_{\bQ_p} \to \GL_2(\bF)$,
there exists a versal ring to~$\cX^{\lambda, \crys}\otimes_{\cO}\F$ at~$\rhobar$ which is reduced.
(This follows for example from an application of Lemma~\ref{lem:pro-coh-vanishes-if-pullbacks-do} to the nilradical of the structure sheaf of $\cX^{\lambda, \tau}\otimes_{\cO}\F$.) %
When $\rhobar$ is not a twist of~$1 \oplus \omega^{-1}$, 
we can take the versal ring to be $R^{\square, \lambda, \crys}_{\rhobar}/\varpi$, i.e.\ the special fibre of the crystalline lifting ring of~$\rhobar$ with weight~$\lambda$, 
which is reduced by~\cite[Corollary~1.7.14]{KisinFM}.
When~$\rhobar$ is a twist of~$1 \oplus \omega^{-1}$, we can apply e.g.\ Corollary~\ref{cor:explicit presentation of crystalline deformation ring II} below.
We thus conclude that $\cX^{\lambda, \crys}\otimes_{\cO}\F = \cZ(\sigma)$, as desired.

In case~(2), 
we know %
by~\cite[Theorem~1.3]{CEGSBM} that $\cX^{\lambda, \tau}\otimes_{\cO}\F$ is generically reduced. 
Furthermore, by~\cite[Theorem~1.1.2]{hung2023local}, $\cX^{\lambda, \tau}$ is normal, so that 
$\cX^{\lambda, \tau}\otimes_{\cO}\F$ is ~$S_1$, and thus reduced, as claimed. %
\end{proof}

\subsection{A continuous map from \texorpdfstring{$|\cX|$}{X} to
  \texorpdfstring{$|X|$}{X}}\label{subsec-cX-to-X}
We now construct a continuous map $\piss:|\cX|\to|X|$ from the underlying topological space of the stack~$\cX$ to the underlying topological space of the chain~$X$ of projective lines 
constructed in Section~\ref{chain of projective lines}.
\subsubsection{Soberization}
Recall that a topological space is \emph{sober} if every irreducible closed subset has a unique generic point. The inclusion of the full subcategory of sober
spaces into the category of all topological spaces (and continuous maps) admits
a left adjoint, the {\em soberization} of a space, which we will denote $\sob(\text{--}).$
Concretely (see e.g.\ \cite[\href{https://stacks.math.columbia.edu/tag/0A2N}{Tag 0A2N}]{stacks-project}),
for any topological space~$S$, we define $\sob(S)$ to be the set of irreducible closed 
subspaces of~$S$.  The topology on $\sob(S)$ is defined by letting its closed subsets
be the collection of subsets~$\sob(T)$, where $T$ ranges over the closed subsets of~$S$.
Of course, the natural morphism (unit of adjunction) $S \to \sob(S)$ is a homeomorphism
if and only if~$S$ itself is sober.

\subsubsection{Underlying topological spaces
of stacks}\label{subsubsec: underlying}
If $\cZ$ is an algebraic stack of finite presentation over a field~$k$ (and so is in particular quasi-compact and quasi-separated),
then we write (as usual) $|\cZ|$ for the  
underlying topological space of~$\cZ$,
$|\cZ|_{\ft}$ for its subset of finite type points
and $|\cZ|_{\cl}$ for its subset of closed points.
We endow each of these subsets with the subspace topology.
We always have that $|\cZ|_{\cl} \subseteq |\cZ|_{\ft}$,
with equality for schemes (or more generally Deligne--Mumford stacks),
but not in general.
We have that $|\cZ|_{\ft}$ is dense in $|\cZ|$,
and that $|\cZ|$ is
sober~\cite[\href{https://stacks.math.columbia.edu/tag/0DQP}{Tag 0DQP}]{stacks-project}.
This implies that 
the inclusion $|\cZ|_{\ft} \hookrightarrow |\cZ|$ induces a homeomorphism
\begin{equation}\label{eqn:soberization isomorphism}
\sob(|\cZ|_{\ft}) \iso |\cZ|.
\end{equation}

\subsubsection{Finite type points of~\texorpdfstring{$\cX$}{X}}\label{subsubsec:points of X}
As already recalled, by \cite[Thm.\ 6.6.3]{emertongeepicture}, the finite type points of~$\cX$ are in bijection with $\Gal(\cbF_p/\F)$-conjugacy classes of continuous representations
$\rhobar: G_{\bQ_p} \to \GL_2(\cbF_p)$ with determinant~$\zeta\varepsilon^{-1}$. 
Furthermore, the finite type point associated to~$\rhobar$ is closed if and only if~$\rhobar$ is semisimple; 
and the closure of~$\{\rhobar\}$ contains a unique closed point, which is associated to the semisimplification $\rhobar^{\mathrm{ss}}$.
Thus we obtain a map (of sets!)
\begin{equation}
\label{eqn:specialization to closed points}
|\cX|_{\ft} \to |\cX|_{\cl}, \rhobar \mapsto \rhobar^{\mathrm{ss}}.
\end{equation}
We caution the reader that this map 
is {\em not} continuous.

\subsubsection{Constructing the morphism}
We let $X$ denote the chain of projective lines over $\F$ recalled in Section~\ref{chain of projective lines}, %
so that $|X|_{\cl}$ is in bijection with the set of blocks of~$\cA^{\ladm}$,
and hence is also in bijection with the set of residual pseudorepresentations.
Thus we obtain a bijection $|\cX|_{\cl} \iso |X|_{\cl}$,
which, when composed with~\eqref{eqn:specialization to closed points},
yields a morphism
\begin{equation}
\label{eqn:cX to X on ft points}
|\cX|_{\ft} \to |X|_{\cl}.
\end{equation}

\begin{lemma}
The map~{\em \eqref{eqn:cX to X on ft points}} is continuous.
\end{lemma}
\begin{proof}
The closed subsets of $|X|_{\cl}$ consist 
of finite unions of points and irreducible components, and so it suffices
to show that the preimage of 
each of these in $|\cX|_{\ft}$ is closed. %

To see that the preimage of a point is closed, it suffices to prove that it is a finite union of closed sets, 
and so it suffices to prove that if~$\rhobar,\rhobar' \in |\cX|_{\ft}$, and $\rhobar' \in \lbar{\{\rhobar\}}$ (the closure in~$|\cX|$ of~$\{\rhobar\}$),
then~$(\rhobar')^{\mathrm{ss}}=\rhobar^{\mathrm{ss}}$. 
This is a consequence of~\cite[Thm.\ 6.6.3~(3)]{emertongeepicture}.
In more detail, we find that $(\rhobar')^{\mathrm{ss}} \in \lbar{\{\rhobar'\}} \subset \lbar{\{\rhobar\}}$, and so 
$(\rhobar')^{\mathrm{ss}}$ must coincide with the unique closed point of
$\lbar{\{\rhobar\}}$, which is $\rhobar^{\mathrm{ss}}$, as desired.

On the other hand, the preimage of an irreducible component of $X$, corresponding to the companion pair
$\{\sigma,\sigmacomp\}$, is equal either to $|\cX(\sigma)|_{\ft} \cup |\cX(\sigmacomp)|_{\ft}$,
if the pair $\{\sigma,\sigmacomp\}$ is not of the form
$\{\sigma_{a,0},\sigma_{a+1,p-3}\}$, 
or to $|\cX(\sigma_{a,0})|_{\ft} \cup |\cX(\sigma_{a,p-1})^+|_{\ft} \cup 
|\cX(\sigma_{a,p-1})^-|_{\ft} \cup 
|\cX(\sigma_{a+1,p-3})|_{\ft}$
otherwise.  
In either case, this is a union of irreducible components of $|\cX|_{\ft}$,
and thus closed.
\end{proof}

Note that $|\cX| = |\cX_{\red}|$ (essentially by definition), and so,
since $\cX_{\red}$ is an algebraic stack of finite presentation over~$\F$, the discussion in Section~\ref{subsubsec: underlying}
applies to $|\cX|$.
In particular, we have a canonical homeomorphism $\sob(|\cX|_{\ft}) \iso |\cX|$. 
Since $X$ is a finite type scheme over~$\F$,
we also have 
a canonical homeomorphism $\sob(|X|_{\cl}) \iso |X|$.
Soberizing the map~\eqref{eqn:cX to X on ft points},
and  post- and pre-composing with the second of these homeomorphisms and 
the inverse of the first,
we obtain the continuous morphism
\begin{equation}\label{eqn:defn-of-cts-map-f}\piss:|\cX| \iso \sob(|\cX|_{\ft}) \buildrel
\sob\text{\eqref{eqn:cX to X on ft points}} \over \longrightarrow
\sob(|X|_{\cl}) \iso |X|\end{equation}
that we wished to construct.%

\subsection{Stacks of Galois representations}\label{subsec:CWE stacks}
As we recall in greater generality in
Appendix~\ref{subsubsec:stacks with fixed pseudorepresentation}, 
Wang--Erickson has associated 
to each pseudorepresentation
$\thetabar$  of $G_{\Q_p}$ %
over~$\Fpbar$, with determinant~$\zetabar\omega^{-1}$, 
a Noetherian formal algebraic stack~$\cX_{\thetabar}$.
If~$A$ is an $\cO$-algebra in which~$p$ is nilpotent, then $\cX_{\thetabar}(\Spf A)$ is the groupoid of
continuous representations
of $G_{\Q_p}$ on rank two projective $A$-modules of determinant~$\zeta\varepsilon^{-1}$,
whose associated residual pseudorepresentation
is $\Gal(\cbF_p/\F)$-conjugate to~$\thetabar$.
Passage to the associated pseudorepresentation yields then a morphism
\begin{equation}
\label{eqn:passage to pseudoreps}
\cX_{\thetabar} \to \Spf R_{\thetabar}^{\ps}.
\end{equation}
On the other hand, passage to the associated \'etale $(\varphi,\Gamma)$-module
yields a morphism
\begin{equation}
\label{eqn:CWE to EG}
i'_{\thetabar}: \cX_{\thetabar} \to \cX.
\end{equation}%

\begin{thm}\label{thm: thetabar substacks}
\leavevmode
\begin{enumerate}
\item
The morphism~\eqref{eqn:passage to pseudoreps}
is representable by algebraic stacks.
\item The morphism~\eqref{eqn:CWE to EG} is a completion of~$\cX$ at the closed subset $|\cX_{\thetabar}| \subset |\cX|$, which coincides with $\pi_{\mathrm{ss}}^{-1}\{\thetabar\}$.
\end{enumerate}
\end{thm}
\begin{proof}
Part~(1) is a special case of Theorem~\ref{athm:CWE-main-result}.
The first claim of part (2) is a special case of Theorem~\ref{thm: thetabar substacks app version} (bearing in mind Remark~\ref{arem:fixed-det-versions-of-Galois-stacks}). %
Finally, since $|\cX_{\thetabar}|$ and $\pi_{\mathrm{ss}}^{-1}\{\thetabar\}$ are closed subsets of~$|\cX|$ with the same set of finite type points, 
the equality $|\cX_{\thetabar}| = \pi_{\mathrm{ss}}^{-1}\{\thetabar\}$ follows from~\eqref{eqn:soberization isomorphism}.
\end{proof}

Furthermore, as explained in Appendix~\ref{subsubsec:stacks with fixed pseudorepresentation},
Wang-Erickson shows in ~\cite[Section~3.2]{MR3831282} that there is
an algebraic stack~$\fX_{\thetabar}$ of finite type over
$\Spec R_{\thetabar}^{\ps}$, such that~$\cX_{\thetabar}$ can be identified
with the completion of~$\fX_{\thetabar}$ at the maximal ideal
of~$R_{\thetabar}^{\ps}$. By definition, $\fX_{\thetabar}$ represents the moduli
problem of $2$-dimensional ``compatible representations''
of the Cayley--Hamilton algebra~$\Rtilde_{\thetabar}$ (see Definition~\ref{defn:CH-alg-Rtilde}) with determinant~$\zeta\varepsilon^{-1}$.
In particular, there is a canonical rank two vector bundle $\fV_{\thetabar}$ lying
over~$\fX_{\thetabar}$, endowed with an action of~$\tR_{\thetabar}$. 

Following~\cite{JNWE}, we now give descriptions 
of $\fX_{\thetabar}$,
together with the bundle $\fV_{\thetabar}$ lying over it,
and the $\tR_{\thetabar}$-action on this bundle.  
We will also define a versal ring to~$\cX_{\thetabar}$ at its unique closed point, and make the following definitions. 

\begin{defn}\label{defn:various versal rings}
Let~$\thetabar$ be a $2$-dimensional $\cbF_p$-valued pseudorepresentation.
\begin{enumerate}
\item We write~$R^{\ver}_{\thetabar}$ for the versal ring to $\cX_{\thetabar}$ defined in Sections~\ref{subsubsec:ss CWE stack}---\ref{subsubsec:not abs irred CWE stack}, and $\Vver_{\thetabar}$ for
the versal object on~$R^{\ver}_{\thetabar}$.
\item For all pairs~$(\lambda, \tau)$
as in Definition~\ref{defn:potentially crystalline EG stacks}, we will write $R^{\lambda, \tau}_{\thetabar}$ for the quotient of~$R^{\ver}_{\thetabar}$ 
defined by $\Spf R_{\thetabar}^{\lambda, \tau} \coloneqq  \Spf R^{\ver}_{\thetabar} \times_{\cX} \cX^{\lambda, \tau}$.
\item If~$\sigma$ is a Serre weight, 
and~$(\lambda, \tau)$ is the pair associated to~$\sigma$ in Proposition~\ref{Z(sigma) as special fibre}, 
we will write $\Rsigma_{\thetabar}$ for the ring $R^{\lambda, \tau}_{\thetabar}$. By Lemma~\ref{Z(sigma) as special fibre}, we thus have
$\Spf \Rsigma_{\thetabar}/\varpi = \Spf R^{\ver}_{\thetabar}/\varpi \times_{\cX/\varpi} \cZ(\sigma)$.
\end{enumerate}
When~$\thetabar$ is clear from the context, we will sometimes omit it from this notation.
\end{defn}

\subsubsection{The case~\emph{\ref{item: ss pseudorep}}}\label{subsubsec:ss CWE stack}
In this case, $\thetabar$ is the trace of an absolutely irreducible representation $\rhobar: G_{\bQ_p} \to \GL_2(\bF)$.
We let~$R^{\ver}_{\thetabar}$ be the universal deformation ring of~$\rhobar$, to complete Noetherian local $\cO$-algebras, with determinant~$\zeta \varepsilon^{-1}$.
Then
$\fX_{\thetabar} \cong [\Spec R^{\ver}_{\thetabar}/\bm{\mu}_2]$
for the trivial action of~$\bm{\mu}_2$,
$\fV_{\thetabar} = (R^{\ver}_{\thetabar})^{\oplus 2}$ placed in degree~$-1$, 
and
$\tR_{\thetabar}$ is isomorphic to $M_2(R^{\ver}_{\thetabar})$ with its natural action on~$\fV_{\thetabar}$.
See for example~\cite[Section~3.1]{JNWE}.

\subsubsection{The case~\emph{\ref{item: generic pseudorep}}}\label{subsubsec:gen CWE stack}
After a twist (which possibly changes~$\zeta$) ,
we can assume without loss of generality that $\thetabar = \nr_\lambda\omega^r + \nr_{\lambda^{-1}\zeta(p)}\omega^{-1}$ for some $\lambda \in \bF^\times$ and some $0 \leq r \leq p-2$, 
such that $\lambda \ne \pm \zeta(p)^{1/2}$ if~$r \in \{0, p-2\}$.
As explained in~\cite[Section~3.2]{JNWE}, the stack $\fX_{\thetabar}$ admits the following
description as a torus quotient:
Let
\[
S = \cO\llbracket a_0, a_1, X\rrbracket [b, c]/(bc-X),
\]
with grading
\[
\deg(b) = 2, \deg(c) = -2, \deg(a_0) = \deg(a_1) = 0.
\] 
Then $\fX_{\thetabar} \cong [\Spec S/\Gm]$ (the $\Gm$-action being determined by the grading
on~$S$).
The map $R^{\ps}_{\thetabar} \to S$,
corresponding to the morphism $\fX_{\thetabar} \to \Spec R^{\ps}_{\thetabar}$,
induces an isomorphism
$$R^{\ps}_{\thetabar} \iso \cO\llbracket a_0,a_1,X\rrbracket ,$$
identifying $R^{\ps}_{\thetabar}$ with the degree zero component of~$S$.

This description of $\fX_{\thetabar}$ allows us to identify
the category
$\Coh(\fX_{\thetabar})$ with the category of finitely generated graded $S$-modules.
As usual, if~$M$ is a graded $S$-module, we write $M(i)$ for the graded shift given by
\[
M(i)_n \coloneqq  M(n+i).
\]
The universal rank two vector bundle $\fV_{\thetabar}$ on $\fX_{\thetabar}$
corresponds to the graded module $S(1) \oplus S(-1)$.
The action of $\tld R_{\thetabar}$ on $\fV_{\thetabar}$ 
is described by identifying $\tR_{\thetabar}$ (as an $R^{\ps}_{\thetabar}$-algebra) with the matrix order
\begin{equation}\label{2x2 matrix order (gen)}
\End_{\text{S-gr}}(\fV_{\thetabar}) =
\fourmatrix{R^{\ps}_{\thetabar}}{R^{\ps}_{\thetabar}b}
{R^{\ps}_{\thetabar}c}{R^{\ps}_{\thetabar}}
\subseteq M_2(S)
\end{equation}
acting by left multiplication on~$\fV_{\thetabar}$ (thought of as column
vectors) (see \cite[Thm.~3.2.1]{JNWE}). %
We write~$e_{ij}$ for the matrix in~$M_2(S)$ with~$1$ in the $ij$-th entry, and zero elsewhere, and we normalize our presentation of~$S$ in such a way that
$\tR_{\thetabar}e_{11}$ (i.e.\ the first column of~$\tR_{\thetabar}$) is the $\tR_{\thetabar}$-projective envelope of the simple module~$\nr_\lambda\omega^r$.
Then the completion~$R^{\ver}_{\thetabar}$ of~$S$ at the maximal homogeneous ideal is a versal ring to~$\cX$ at 
$\rhobar \coloneqq  \nr_\lambda \omega^r \oplus \nr_{\lambda^{-1}\zeta(p)}\omega^{-1}$.

\subsubsection{The case~\emph{\ref{item: non p-distinguished pseudorep}}}\label{subsubsec:scalar CWE stack}
In this case, $\thetabar$ is not multiplicity-free, and so $\fX_{\thetabar}$
does not have a natural presentation as a $\Gm$-quotient of an affine scheme.
Rather,
$\fX_{\thetabar}$ is the stack denoted~$\Rep(E)$
in~\cite[Section~3.3]{JNWE},
where~$E = \tld R_{\thetabar}$, and so it is a global quotient $\Rep(E) = [\Rep^\square(E)/\SL_2]$
for some affine scheme $\Rep^\square(E)~=~\Spec A$.
A presentation of~$A$ as an $R_{\thetabar}^{\ps}$-algebra is given after the proof of~\cite[Proposition~3.3.9]{JNWE};
we remark that
\[
A \otimes_{R^{\ps}_{\thetabar}} R^{\ps}_{\thetabar}/\fm \cong \bF[a, b, c_1, c_2, d_1, d_2]/(a^2+c_1c_2, b^2+d_1d_2, 2ab + (c_1d_2+c_2d_1)),
\] 
which is not a reduced ring (since for example $c_1d_2-c_2d_1$ is a nonzero element which squares to zero).
There is a unique maximal $\SL_2$-invariant ideal~$\fn$ of~$A$, and it is a maximal ideal;
the corresponding morphism
$$
\Spec A/\fn \hookrightarrow \Spf A \to \cX_{\thetabar} \hookrightarrow \cX
$$
classifies the unique scalar Galois representation
$\rhobar: G_{\Q_p} \to \GL_2(A/\fn)$ having trace~$\thetabar$.
We write~$R^{\ver}_{\thetabar}$ for the $\fn$-adic completion of~$A$, which is a versal ring to~$\cX$ at the closed point corresponding to this~$\rhobar$.

\subsubsection{The case~{\em \ref{item: Steinberg pseudorep}}}
\label{subsubsec:Steinberg CWE stack}
After a twist (which possibly changes~$\zeta$),
 we can assume without loss of generality that $\thetabar = 1 + \omega^{-1}$.
As explained in~\cite[Section~3.4]{JNWE}, the stack $\fX_{\thetabar}$ admits the following
description as a torus quotient:
Let
\[
S = \cO\llbracket a_0, a_1, X_0, X_1\rrbracket [b_0, b_1, c]/(X_0-b_0c, X_1-b_1c, a_1b_0 + a_0b_1),
\]
with grading
\[
\deg(b_0) = \deg(b_1) = 2, \deg(c) = -2, \deg(a_0) = \deg(a_1) = 0.
\] 
Then $\fX_{\thetabar} \cong [\Spec S/\Gm]$ (the $\Gm$-action being determined by the grading
on~$S$).
The map $R^{\ps}_{\thetabar} \to S$,
corresponding to the morphism $\fX_{\thetabar} \to \Spec R^{\ps}_{\thetabar}$,
induces an isomorphism
$$R^{\ps}_{\thetabar} \iso \cO\llbracket a_0,a_1,X_0,X_1\rrbracket /(a_0 X_1+ a_1 X_0),$$
identifying $R^{\ps}_{\thetabar}$ with the degree zero component of~$S$.
Note that our~$a_1$ is denoted $a'_1=a_1 +p$ in~\cite{JNWE}, and our $X_i$ is denoted~$b_ic$.

The universal rank two vector bundle $\fV_{\thetabar}$ on $\fX_{\thetabar}$
corresponds to the graded module $S(1) \oplus S(-1)$.
The action of $\tld R_{\thetabar}$ on $\fV_{\thetabar}$ 
is described by identifying $\tR_{\thetabar}$ (as an $R^{\ps}_{\thetabar}$-algebra) with the matrix order
\begin{equation}\label{2x2 matrix order}
\End_{\text{S-gr}}(\fV_{\thetabar}) =
\fourmatrix{R^{\ps}_{\thetabar}}{R^{\ps}_{\thetabar}b_0+R^{\ps}_{\thetabar}b_1}
{R^{\ps}_{\thetabar}c}{R^{\ps}_{\thetabar}}
\subseteq M_2(S)
\end{equation}
acting by left multiplication on~$\fV_{\thetabar}$ (thought of as column
vectors) (see \cite[Thm.~3.4.1]{JNWE}). %
Again we write~$e_{ij}$ for the matrix in~$M_2(S)$ with~$1$ in the $ij$-th entry, and zero elsewhere.
Then the modules $\tR_{\thetabar}e_{11}, \tR_{\thetabar}e_{22}$ are the $\tR_{\thetabar}$-projective envelopes of the simple
modules~$1$,  $\omega^{-1}$ respectively.
We let~$R^{\ver}$ denote
the completion %
of~$S$ at its maximal homogeneous ideal; then $R^{\ver}$ is a versal ring to~$\cX$ at 
the closed point corresponding to $\rhobar \coloneqq  1 \oplus \omega^{-1}$.

It will also be useful to have some notation for the scheme-theoretic union of the (finitely many) substacks~$\cX_{\thetabar} \subset \cX$
where $\thetabar$ is of type~\ref{item: Steinberg pseudorep}, and so we make the following definition.

\begin{df}
\label{def:Steinberg locus}
We write $\cX(\St) \coloneqq  \bigcup_{\thetabar \text{ of type \ref{item: Steinberg pseudorep}}} 
\cX_{\thetabar}$. 
\end{df}

\begin{rem}
  \label{rem:Steinberg-components-are-CWE-stacks}For~$\thetabar$ of type~\ref{item: Steinberg
  pseudorep}, the stack $\cX_{\thetabar,\red}$ coincides with one of the irreducible
components~$\cX(\sigma_{a,p-1})^{\pm}$
(with~$a$ and the sign~$\pm$
depending on~$\thetabar$).
Thus $\cX(\St)_{\red}$ is the disjoint union of four irreducible
components~$\cX(\sigma_{a,p-1})^{\pm}$ (assuming that  it is non-empty, i.e.\ that
$\zeta$ is even).
\end{rem}

\subsubsection{The case~{\em\ref{item: not abs irred pseudorep}}}\label{subsubsec:not abs irred CWE stack}
If~$\thetabar$ is a $\cbF_p$-pseudorepresentation of type~\ref{item: not abs irred pseudorep}, then there exists a smooth character $\chi : G_{\bQ_p} \to \cbF_p^\times$
such that $\thetabar = \chi+\zeta\omega^{-1}\chi^{-1}$, and~$\chi$ does not factor through~$\F^\times$.
Writing $\F_{\thetabar}/\F$ for the extension generated by the values of~$\thetabar$, we can identify $\thetabar$ with an $\F_{\thetabar}$-valued pseudorepresentation~$\thetabar'$,
and by definition, we have equalities
\[
R_{\thetabar}^{\ps} = R_{\thetabar'}^{\ps}, \tR_{\thetabar} = \tR_{\thetabar'},
\fX_{\thetabar} = \fX_{\thetabar'}, \cX_{\thetabar} = \cX_{\thetabar'}.
\]
Up to natural isomorphism, these objects only depend on the $\Gal(\cbF_p/\F)$-conjugacy class of~$\thetabar$.

Note that~$\F_{\thetabar}$ almost always coincides with the extension~$\F_{\thetabar}'/\F$ generated by the values of~$\chi$, and then $\thetabar'$ is the trace of a 
semisimple reducible Galois representation, classified by a morphism
$\Spec \F_{\thetabar} \to \cX_{\thetabar}$.
We then let $R^{\ver}_{\thetabar}$ be the versal ring at $\Spec \F_{\thetabar} \to \cX_{\thetabar}$ defined in Section~\ref{subsubsec:gen CWE stack}.

The only exception is when $\chi$ and~$\zeta \omega^{-1}\chi^{-1}$ are Galois conjugates over~$\F$, 
in which case $\thetabar'$ is the trace of an %
irreducible Galois representation~$\rhobar:G_{\Q_p} \to \GL_2(\F_{\thetabar})$,
for which $\rhobar \otimes_{\F_{\thetabar}}\F_{\thetabar}'$
is reducible.
We have a morphism $\Spec \F_{\thetabar}' \to \cX_{\thetabar}$ classifying
$\rhobar \otimes_{\bF_{\thetabar}} \bF'_{\thetabar}$, and
we let
\[
\Spf R_{\thetabar}^{\ver} \to \cX_{\thetabar} \times_{W(\F_{\thetabar})} W(\F'_{\thetabar})
\]
be the versal morphism to $\Spec \F'_{\thetabar} \to \cX_{\thetabar} \times_{W(\F_{\thetabar})} W(\F'_{\thetabar})$ defined in
Section~\ref{subsubsec:gen CWE stack}.
Note that the composite $\Spf R^{\ver}_{\thetabar} \to \cX_{\thetabar}$ is versal to $\Spec \F_{\thetabar}' \to \cX_{\thetabar}$.

\subsection{Coherent sheaves on stacks of Galois representations}\label{subsec:sheaves on CWE stacks}
We fix throughout this section a 2-dimensional $\cbF_p$-valued pseudorepresentation~$\thetabar$ of~$G_{\bQ_p}$.
We will discuss some coherent sheaf theory on the stacks~$\fX_{\thetabar}$ and~$\cX_{\thetabar}$ defined in Section~\ref{subsec:CWE stacks}.
To ease notation, we write $R\coloneqq  R^{\ps}_{\thetabar}$ and $\fm\coloneqq  \fm_R$.
We write $\fX_0$ to denote the underlying reduced substack of~$\cX_{\thetabar}$;
this coincides with the underlying reduced substack of the
base-change~$\fX_{\thetabar}\times_{\Spec R}\Spec R/\fm.$

\begin{remark}\label{rem:nonreduced case}
We saw above (see in particular Section~\ref{subsubsec:scalar CWE stack} that in case~\ref{item: non p-distinguished pseudorep}) that in fact 
$\fX_{\thetabar} \times_{\Spec R} \Spec R/\fm$ is reduced, except in case~\ref{item: non p-distinguished pseudorep}.
\end{remark}

Adopting notation from Appendix~\ref{subsubsec:stacks with fixed pseudorepresentation}, 
we write $k_{\thetabar} : \cX_{\thetabar} \to \fX_{\thetabar}$ for the completion map, and
\[
k_{\thetabar, *}: D^b_{\coh}(\cX_{\thetabar}) \to D^b_{\coh}(\fX_{\thetabar}),
\]
\[
\widehat k^*_{\thetabar} : D^b_{\coh}(\fX_{\thetabar}) \to \Pro D^b_{\coh}(\cX_{\thetabar})
\]
for the pushforward and completed pullback functors.
Their basic properties are summarized in Theorem~\ref{thm:coherent-completeness-CWE-stacks} (see also Remark~\ref{arem:fixed-det-versions-of-Galois-stacks}).

\subsubsection{Tensoring pro-coherent sheaves with compact modules}
\label{subsubsec:tensoring-pro-coherent-compact}
We now fix a (not necessarily commutative) $R$-algebra
$E$ which is finite and flat as an $R$-module.
We regard~$E$ as a topological $R$-algebra
by endowing it with its $\fm$-adic topology. 
With this topology, $E$ is a Noetherian profinite $\cO$-algebra, and so we can apply the material of Section~\ref{compact modules}.
We also suppose given a coherent sheaf $\fF$
on $\fX_{\thetabar}$ endowed with a right action %
of~$E$ extending
the natural $R$-action on~$\fF$.
We write $\cF \coloneqq  \widehat{\imm}_{\thetabar}^* \fF$,
so $\cF$ is the $\fm$-adic completion of~$\fF$;
then $\cF$ is an object of $\Pro \Coh(\cX_{\thetabar})$ 
endowed with a right $E$-action
extending the $R$-action on~$\fF$.

Our goal in this subsection is to discuss the various functors  
we can obtain by tensoring $\fF$ and $\cF$ with different flavours of $E$-modules.
We first discuss things on the abelian level, before turning to the analogous
constructions on the level of stable $\infty$-categories.

By Proposition~\ref{prop:Eilenberg-Watts} we have a right exact functor %
\begin{equation}
\label{eqn:tensoring with fF}
\fF
\otimes_E\text{--} : \Mod^{\fp}(E) \to \Coh(\fX_{\thetabar}),
\end{equation} 
characterized up to unique isomorphism by sending~$E$ to~$\fF$.
Similarly,
we have a unique right exact functor
\begin{equation}
\label{eqn:tensoring with cF}
\cF\otimes_E\text{--} : \Mod^{\fp}(E) \to \Pro \Coh(\cX_{\thetabar})
\end{equation}
sending~$E$ to~$\cF$.
Since $\cF = \widehat{\imm}_{\thetabar}^* \fF$, we thus have a natural isomorphism
\begin{equation}
\label{eqn:tensor completion compatibility}
\cF\otimes_E\text{--}  \iso \widehat{\imm}_{\thetabar}^* \circ (\fF \otimes_E\text{--} ) . 
\end{equation}

On the other hand, by definition, $\cF$ is complete
as an $R$-module object of $\Pro\Coh(\cX_{\thetabar})$.
Remark~\ref{complete if and only if derived complete} then shows that~$\cF$ is also derived complete,
and in turn, by Lemma \ref{lem:equivalence-of-R-and-E-complete}, we see that
$\cF$ is complete and derived complete as an~$E$-module. 
We thus obtain, by Lemma~\ref{lem:EW-gives-completed-tensor-product},
a unique right exact and cofiltered limit-preserving functor
\begin{equation}\label{eqn:completed tensoring with cF}
\cF\cotimes_E\text{--}  : \Mod_c(E) \to \ProCoh(\cX_{\thetabar})
\end{equation}
sending~$E$ to~$\cF$.

\begin{lem}
  \label{lem:comparing-frak-and-cal-tensors --- abelian case}
With the notation of the previous paragraph, there exists a unique (up to isomorphism) right exact functor
\begin{equation}
\label{eqn:fl tensoring with cF}
\cF\otimes_E\text{--} : \Mod^{\fl}(E) \to \Coh(\cX_{\thetabar})
\end{equation}
such that the following diagram commutes: 

\[\begin{tikzcd}
	 \Mod^{\fl}(E) &
        \Mod^{\fp}(E) & \Mod_{c}(E) \\
	 \Coh(\cX_{\thetabar}) &
        \Coh(\mf{X}_{\thetabar}) & \ProCoh(\cX_{\thetabar}).
	\arrow[hook, from=1-1, to=1-2]
	\arrow["\cF\cotimes_E\text{--} ", "\eqref{eqn:completed tensoring with cF}"',from=1-3, to=2-3]
	\arrow[hook, from=1-2, to=1-3]
	\arrow["\mf{F}\otimes_E\text{--} ", "\eqref{eqn:tensoring with fF}"', from=1-2, to=2-2]
	\arrow["\cF\otimes_E\text{--} ", "\eqref{eqn:fl tensoring with cF}"', from=1-1, to=2-1]
  \arrow[hook, "k_{\thetabar, *}", from=2-1, to=2-2]
	\arrow[hook, "\widehat k^*_{\thetabar}", from=2-2, to=2-3]
\end{tikzcd}\]
Furthermore, $\cF \cotimes_E\text{--} $ is the $\Pro$-extension of $\cF \otimes_E\text{--} $.
\end{lem}
\begin{proof}
By Corollary~\ref{cor:restricting-completed-tensor}, the restriction of~\eqref{eqn:completed tensoring with cF}
through $\Mod^{\fp}(E) \to \Mod_c(E)$ is~\eqref{eqn:tensoring with cF}.
Hence~\eqref{eqn:tensor completion compatibility} shows that the rightmost square commutes.

We now construct the left-hand square.
Each of the quotients %
\[
\fF/\fm^n \fF \coloneqq  \fF \otimes_E (E/\fm^n E) = \cF \otimes_E (E/\fm^n E) 
\]
is an object of~$\Coh_{\fX_0}(\fX_{\thetabar})$, since $\fX_0$ was defined to be
the underlying reduced substack of the vanishing locus of~$\fm$.
If $N$ is a finite length $E$-module, then it is annihilated by some power of~$\m$,
so we find more generally that 
$\fF\otimes_E N$ lies in $\Coh_{\fX_0}(\fX_{\thetabar}),$
which is to say that~\eqref{eqn:tensoring with fF} restricts to a right exact functor
\begin{equation}
\label{eqn:fl tensoring with fF}
\fF\otimes_E\text{--} : \Mod^{\fl}(E) \to \Coh_{\fX_0}(\fX_{\thetabar}).
\end{equation}
Theorem~\ref{thm:coherent-completeness-CWE-stacks} shows that
$\imm_{\thetabar, *}$ is an equivalence
$\Coh(\cX_{\thetabar}) \iso \Coh_{\fX_0}(\fX_{\thetabar})$
We deduce the existence and uniqueness of~\eqref{eqn:fl tensoring with cF}, and the commutativity of the left-hand square.
Finally, bearing in mind that $\Mod_c(E) = \Pro \Mod^{\fl}(E)$ and $\cF \cotimes_E\text{--} $ preserves cofiltered limits, we see that the commutativity of the outer rectangle implies that $\cF \cotimes_E\text{--}$ is the
Pro-extension of $\cF \otimes_E\text{--} $, as desired.
\end{proof}

We now apply
Lemma~\ref{lem:derived-tensor-product-EW}
to construct derived versions of the various tensor products just introduced.
Namely, 
Lemma~\ref{lem:derived-tensor-product-EW}~\eqref{item:83} shows that~\eqref{eqn:tensoring
with fF} induces a right $t$-exact functor
\begin{equation}
\label{eqn:derived tensoring with fF}
\fF\otimes^L_E\text{--}  : D^b_{\fp}(E) \to D^-\bigl( \Coh(\fX_{\thetabar})\bigr),
\end{equation}
while
Lemma~\ref{lem:derived-tensor-product-EW}~\eqref{item:84} shows that~\eqref{eqn:tensoring
with cF} induces a right $t$-exact functor
\begin{equation}
\label{eqn:derived tensoring with cF}
\cF\otimes^L_E\text{--}  : D^b_{\fp}(E) \to \Pro D^b_{\coh}(\cX_{\thetabar}),
\end{equation}
where we implicitly use the equivalence
$D^b\bigl(\Coh(\cX_{\thetabar})\bigr) \iso D^b_{\coh}(\cX_{\thetabar})$.
Finally, 
Lemma~\ref{lem:derived-tensor-product-EW}~\eqref{item:85} shows that~\eqref{eqn:fl tensoring
with cF} induces a right $t$-exact limit preserving functor
\begin{equation}
\label{eqn:completed derived tensoring with cF}
\cF\cotimes^L_E\text{--}   : \Pro D^b_{\fl}(E) \to \Pro D^b_{\coh}(\cX_{\thetabar}).
\end{equation}
For simplicity, from now on we will make the additional assumption that~$E$ has finite global dimension; in our applications, this will be true 
by Lemma~\ref{lem:projective generators are flat over their endos}~\eqref{item:projective 3}. 
The functor~\eqref{eqn:derived tensoring with fF}
then factors through $D^b(\Coh(\fX_{\thetabar}))$, by Lemma~\ref{derived tensor and finite global dimension}.

\begin{lem}
  \label{lem:comparing-frak-and-cal-tensors --- derived case}
Assume that~$E$ has finite global dimension.
Then there exists a unique (up to isomorphism) right $t$-exact functor
\begin{equation}
\label{eqn:fl derived tensoring with cF}
\cF\otimes^L_E\text{--}  : D^b_{\fl}(E) \to D^b\bigl(\Coh(\cX_{\thetabar})\bigr)
\end{equation}
such that the following diagram commutes:
\[\begin{tikzcd}
	 D^b_{\fl}(E) &
        D^b_{\fp}(E) & \Pro D^b_{\fl}(E) \\
         D^b\bigl(\Coh(\cX_{\thetabar})\bigr) & D^b\bigl(\Coh(\fX_{\thetabar})\bigr)
      & \Pro D^b_{\coh}(\cX_{\thetabar}).
	\arrow["\eqref{eqn:fl to fp}", hook, from=1-1, to=1-2]
	\arrow["\cF\cotimes^L_E\text{--} ",from=1-3, to=2-3]
	\arrow[hook, from=1-2, to=1-3, "\eqref{eqn:derived fp to pro fl}"]
	\arrow["\mf{F}\otimes^L_E\text{--} ",from=1-2, to=2-2]
	\arrow["\cF\otimes^L_E\text{--} ",from=1-1, to=2-1]
        	\arrow[hook, from=2-1, to=2-2, "k_{\thetabar, *}"]
        \arrow["\cF\otimes^L_E\text{--} ", from=1-2, to=2-3]
	\arrow[hook, from=2-2, to=2-3, "\widehat k_{\thetabar}^*"]
\end{tikzcd}\]
\end{lem}
\begin{proof}
  The upper triangle commutes by Lemma~\ref{upper triangle in diagram}.
  The two directions in the lower triangle are right $t$-exact functors $D^b_{\fp}(E) \to \Pro D^b_{\coh}(\cX_{\thetabar})$, sending~$E$ to~$\cF$, resp.\ 
  $\widehat k_{\thetabar}^* \fF = \cF$.
  Hence they are naturally isomorphic, by Lemma~\ref{lem:derived-tensor-product-EW}~(2).

  We now construct~\eqref{eqn:fl derived tensoring with cF}, and prove the commutativity of the left-hand square.
  Since
  $\widehat k^*_{\thetabar} : D^b(\Coh(\cX_{\thetabar})) \to D^b_{\coh, \fX_0}(\fX_{\thetabar})$
  is an equivalence,
  by Theorem~\ref{thm:coherent-completeness-CWE-stacks},
  it suffices to
  prove that the composite 
  \[
  D^b_{\fl}(E) \to D^b_{\fp}(E) \xrightarrow{\fF \otimes_E^L\text{--}} D^b(\Coh(\fX_{\thetabar}))
  \]
  factors through $D^b_{\coh, \fX_0}(\fX_{\thetabar})$, i.e.\ the full subcategory of $D^b_{\coh}(\fX_{\thetabar})$ 
  whose objects are the bounded complexes of coherent sheaves whose cohomology is set-theoretically supported on
  the vanishing locus~$\fX_0$ of~$\fm \coloneqq  \rad(R)$ in~$\fX$.
  By induction on the amplitude of objects in $D^b_{\fl}(E)$, it suffices to prove that $\fF \otimes^L_E N \in D^b_{\coh, \fX_0}(\Coh(\fX_{\thetabar}))$ whenever~$N$ is an object of~$\Mod^{\fl}(E)$,
  which is a consequence of the fact that~$N$ is annihilated by a power of~$\fm$.
\end{proof}

\subsubsection{Tensoring pro-coherent sheaves with \texorpdfstring{$G$}{G}-representations}
\label{subsubsec:pro-coherent tensor G reps}
We now fix a projective generator
$P$ of $\fC_{\thetabar}$, chosen to have finite length cosocle,
and write 
\[
E \coloneqq  \End_{\fC_{\thetabar}}(P).
\]  
Recall from Lemma~\ref{finite cosocle implies Noetherian endomorphisms} and Lemma~\ref{lem:projective generators are flat over their endos}~\eqref{item:projective 3}
that $E$ has finite type as a module over its centre (which is the Bernstein centre of~$\fC_{\thetabar}$), has finite global dimension,
and is Noetherian.
Hence, by Lemma~\ref{properties of compact modules}~\eqref{item: compact 7}, the natural topology on $E$ (that it inherits as 
the endomorphism ring of the compact $\cO\llbracket G\rrbracket _{\zeta}$-module~$P$)
coincides with its $\fm$-adic topology, and makes
$E$ a Noetherian profinite $\cO$-algebra.
Recall also from Remark~\ref{description of Bernstein centres} that there is a unique isomorphism $R \isoto Z(E) = \cZ(\fC_\fB)$ such that the functor $\Vcheck$ is $R$-linear.
Furthermore, we assume given a coherent sheaf $\fF$ on $\fX_{\thetabar}$
with a left $E^{\op}$-action, %
so that we are in the situation of Section~\ref{subsubsec:tensoring-pro-coherent-compact}.

Applying the functor~\eqref{eqn:completed tensoring with cF} to~$P$, we obtain
an object $\cF \cotimes_E P$ of $\ProCoh(\cX_{\thetabar})$.
Since $P$ is projective in $\Mod_c(E)$
(by Lemma~\ref{lem:projective generators are flat over their endos}~\eqref{item:projective 3}),
this coincides
with the object $\cF\cotimes^L_E P$ obtained by
applying~\eqref{eqn:completed derived tensoring with cF} to~$P$.
The right $\cO\llbracket G\rrbracket _{\zeta}$-action on $P$ induces a right $\cO\llbracket G\rrbracket _{\zeta}$-action
on $\cF\cotimes_E P$.  
An application of Lemma~\ref{lem:derived-tensor-product-EW}~\eqref{item:84}
then produces a right $t$-exact functor
$$(\cF\cotimes_E P)\otimes^L_{\cO\llbracket G\rrbracket _{\zeta}}\text{--} : D^b_{\fp}(\cO\llbracket G\rrbracket _{\zeta}) \to
\Pro D^b_{\coh}(\cX_{\thetabar}).$$ 

Another application of Lemma~\ref{lem:derived-tensor-product-EW}~\eqref{item:84},
this time to $P$ itself, yields a right $t$-exact functor
$$P\otimes^L_{\cO\llbracket G\rrbracket _{\zeta}}\text{--}: D^b_{\fp}(\cO\llbracket G\rrbracket _{\zeta}) \to \Pro D^b_{\fl}(E).$$
In fact this functor may be factored as a composite
\begin{equation}
\label{eqn:composite derived tensor}
P\otimes^L_{\cO\llbracket G\rrbracket _{\zeta}}\text{--}: D^b_{\fp}(\cO\llbracket G\rrbracket _{\zeta}) \to
D^-\bigl(\Mod_c(E)\bigr) \to \Pro D^b_{\fl}(E),
\end{equation}
where the first arrow is the restriction to 
\[
D^b_{\fp}(\cO\llbracket G\rrbracket _{\zeta})
\iso D^b\bigl(\Mod^{\fp}(\cO\llbracket G\rrbracket _{\zeta})\bigr) \subset
D^-\bigl(\Mod^{\fp}(\cO\llbracket G\rrbracket _{\zeta})\bigr) 
\]
of the (right $t$-exact) derived
functor of the functor
$$P\otimes_{\cO\llbracket G\rrbracket _{\zeta}}\text{--}:
\Mod^{\fp}(\cO\llbracket G\rrbracket _{\zeta}) \to \Mod_c(E)$$
arising from an application of Theorem~\ref{thm:left derived functors};
and the second arrow is the $t$-exact functor
$$D^-\bigl(\Mod_c(E)\bigr) \iso D^-\bigl(\Pro \Mod^{\fl}(E)\bigr) \xrightarrow{p} \Pro D^b_{\fl}(E)$$
obtained as a special case of the diagram~\eqref{eqn:Pro comparison},
taking $\cC$ there to be $\Mod^{\fl}(E)$, and recalling the equivalence
$\Mod_c(E) \iso \Pro\Mod^{\fl}(E)$. 
That this composite coincides, up to natural isomorphism, with the functor
$P\otimes^L_{\cO\llbracket G\rrbracket _{\zeta}}\text{--}$ follows from
Lemma~\ref{lem:derived-tensor-product-EW}~\eqref{item:84},
since both functors are right $t$-exact,
and both take $\cO\llbracket G\rrbracket _{\zeta}$ to the $E$-module object $P$ in~$\Mod_c(E)
\iso \Pro\Mod^{\fl}(E)$.

We will often make implicit use of the following lemmas involving these two functors.

\begin{lemma}\label{lem:t-exactness on cA}
The composite
\[
D^b_{\fp}(\cA) \xrightarrow{\eqref{eq:Afp-in-Modfp-OG}} D^b_{\fp}(\cO\llbracket G\rrbracket _\zeta) \xrightarrow{P \otimes^L_{\cO\llbracket G\rrbracket _\zeta}\text{--}} \Pro D^b_{\fl}(E)
\]
is $t$-exact.
\end{lemma}
\begin{proof}
Let $F$ denote the composite functor appearing in the statement of the lemma.
By Corollary~\ref{cor:derived functors of exact functors}, it suffices to prove that the restriction of~$F$ to~$\cA^{\fp}$
(which is the heart
of $D^b(\cA^{\fp})\iso D^b_{\fp}(\cA)$) takes values in
the heart $\Pro \Mod^{\fl}(E)$ of $\Pro D^b_{\fl}(E)$.

By the alternative description of $P\otimes^L_{\cO\llbracket G\rrbracket _{\zeta}}\text{--}$ 
provided by~\eqref{eqn:composite derived tensor},
it suffices to prove that the composite
\[
\cA^{\fp} \subset D^b_{\fp}(\cA) \xrightarrow{\eqref{eq:Afp-in-Modfp-OG}} D^b_{\fp}(\cO\llbracket G\rrbracket _\zeta) \to D^-(\Mod_c(E))
\]
has essential image in $\Mod_c(E)$.
To do so, it suffices to prove that the composite 
\[
\cA^{\fp} \subset D^b_{\fp}(\cA) \xrightarrow{\eqref{eq:Afp-in-Modfp-OG}} D^b_{\fp}(\cO\llbracket G\rrbracket _\zeta) \to D^-(\Mod_c(E)) \to D^-(\Mod(E))
\]
has essential image in~$\Mod(E)$, where the final arrow is the $t$-exact extension of the forgetful functor $\Mod_c(E) \to \Mod(E)$ (which is exact and conservative, by Lemma~\ref{properties of compact modules}~\eqref{item: compact 1}).
Now the composite $D^b_{\fp}(\cO \llbracket G \rrbracket_\zeta) \to D^-(\Mod(E))$
is right $t$-exact and takes $\cO\llbracket G \rrbracket_\zeta$ to~$P$, 
hence coincides with the left derived functor of the usual tensor product $P \otimes_{\cO\llbracket G \rrbracket_\zeta} -$, by Theorem~\ref{thm:left derived functors}.
Hence our claim follows from Lemma~\ref{lem:projective Tor vanishing}~(2).
\end{proof}

\begin{lemma}
\label{lem:tensor product comparison}
There is a natural isomorphism of right $t$-exact functors $D^b_{\fp}(\cO\llbracket G\rrbracket _{\zeta}) \to
\Pro D^b_{\coh}(\cX_{\thetabar})$
\begin{equation}\label{eqn:derived tensor associativity I}
(\cF\cotimes_E P) \otimes^L_{\cO\llbracket G\rrbracket _{\zeta}}\text{--} = (\cF\cotimes^L_E P) \otimes^L_{\cO\llbracket G\rrbracket _{\zeta}}\text{--}
\iso \cF \cotimes^L_E (P \otimes^L_{\cO\llbracket G\rrbracket _{\zeta}}\text{--}),
\end{equation}
and there is a natural isomorphism of right $t$-exact functors %
$\Pro D^b_{\fp}(\cO\llbracket G\rrbracket _\zeta) \to \Pro D^b_{\coh}(\cX_{\thetabar})$
\begin{equation}\label{eqn:derived tensor associativity II}
(\cF\cotimes_E P) \ccotimes^L_{\cO\llbracket G\rrbracket _{\zeta}}\text{--} = (\cF \cotimes^L_E P) \ccotimes_{\cO\llbracket G\rrbracket _\zeta}^L\text{--} \isoto \cF \cotimes^L_E (P \ccotimes_{\cO\llbracket G\rrbracket _\zeta}^L\text{--}).
\end{equation}
\end{lemma}
\begin{proof}Recalling again that $P$ is projective in $\Mod_c(E)$,
Lemma~\ref{lem:derived-tensor-product-EW}~\eqref{item:84}
shows that to prove~\eqref{eqn:derived tensor associativity I}, it suffices to note that each functor takes~$\cO\llbracket G\rrbracket _{\zeta}$ to the object $\cF \cotimes^L_E P=\cF \cotimes_E P$ of $\ProCoh(\cX_{\thetabar})$. %
On the other hand, since $P \ccotimes_{\cO\llbracket G\rrbracket _\zeta}^L\text{--}$, resp.\ $(\cF \cotimes^L_E P) \ccotimes_{\cO\llbracket G\rrbracket _\zeta}^L\text{--}$, is the Pro-extension of
$P \otimes_{\cO\llbracket G\rrbracket _\zeta}^L\text{--}$, resp.\ $(\cF \cotimes^L_E P) \otimes_{\cO\llbracket G\rrbracket _\zeta}^L\text{--}$,
and $\cF \cotimes^L_E -$ is cofiltered limit-preserving,
the isomorphism~\eqref{eqn:derived tensor associativity II} can be taken to be the Pro-extension of~\eqref{eqn:derived tensor associativity I}.
\end{proof}

\subsubsection{The universal vector bundle}\label{subsubsec:the universal vector bundle}
Recall from Section~\ref{subsec:CWE stacks} that the universal object on~$\fX_{\thetabar}$ is denoted~$\fV_{\thetabar}$.
It is a free $\cO_{\fX_{\thetabar}}$-module of rank two with a left action of~$\tld R_{\thetabar}$.
Its pullback via $\widehat{\imm}_{\thetabar}^*$
to~$\ProCoh(\cX_{\thetabar})$ %
is denoted 
\begin{equation}\label{eqn:definition of V_thetabar}
\cV_{\thetabar} \coloneqq  \widehat{\imm}_{\thetabar}^*\fV_{\thetabar}. 
\end{equation}
We may apply the discussion of Section~\ref{subsubsec:tensoring-pro-coherent-compact}
to~$\fV_{\thetabar}$ and~$\cV_{\thetabar}$, viewing them as right $\tR^{\op}_{\thetabar}$-modules.
In particular, from the commutative diagram 
of Lemma~\ref{lem:comparing-frak-and-cal-tensors --- abelian case},
we obtain the
commutative diagram 
\begin{equation}
\label{eq:comparing-tensor-cV-fV}
\begin{tikzcd}
	\Mod^{\fl}(\tR_{\thetabar}^{\op}) & \Coh(\cX_{\thetabar})  \\
	& \Coh(\fX_{\thetabar})
	\arrow["\cV_{\thetabar}\otimes_{\tR_{\thetabar}^{\op}}\text{--}", from=1-1, to=1-2]
	\arrow["\fV_{\thetabar}\otimes_{\tR_{\thetabar}^{\op}}\text{--}"', from=1-1, to=2-2]
	\arrow[hook, from=1-2, to=2-2]
\end{tikzcd}
\end{equation}
We now describe the main properties of~$\fV_{\thetabar}$ and of~$\cV_{\thetabar}$ that we will use in Section~\ref{sec: the functor}.
These results are essentially due to Johansson--Newton--Wang-Erickson~\cite{JNWE}.

\begin{prop}\label{topological flatness of V}
Assume that~$\thetabar$ does not have type~{\em \ref{item: Steinberg pseudorep}}.
Then~$\cV_{\thetabar}$, which is a complete right $\tR^{\op}_{\thetabar}$-module object in $\Pro \Coh(\cX_{\thetabar})$, is topologically flat 
{\em (}in the sense of Definition~{\em \ref{topologically flat object}}{\em )}.
\end{prop}
\begin{proof}
By the discussion in Section~\ref{subsubsec:not abs irred CWE stack}, we can assume without loss of generality that $\F = \F_{\thetabar}$.
We need to prove that
the functor 
\[
\cV_{\thetabar} \cotimes_{\tR_{\thetabar}^{\op}}\text{--} : \Mod_c(\tR_{\thetabar}^{\op}) \to \Pro\Coh(\cX_{\thetabar})
\]
is exact.
Since $\Mod_c(\tR_{\thetabar}^{\op}) = \Pro \Mod^{\fl}(\tR_{\thetabar}^{\op})$, 
it suffices to prove that its restriction to $\Mod^{\fl}(\tR_{\thetabar}^{\op})$ is exact.
By Lemma~\ref{lem:comparing-frak-and-cal-tensors --- abelian case}, %
this restricted functor is the same as the composite
\[
\Mod^{\fl}(\tR_{\thetabar}^{\op}) \to \Mod^{\fp}(\tR_{\thetabar}^{\op}) \xrightarrow{\mathfrak{V}_{\thetabar} \otimes_{\tR_{\thetabar}^{\op}}\text{--}} 
\Coh(\fX_{\thetabar}) \to \Pro \Coh (\cX_{\thetabar}).
\]
The last arrow is exact, so it suffices to prove that the composite of the first two arrows is exact.

Assume first that~$\thetabar$ does not have type~\ref{item: non p-distinguished pseudorep}.
We will prove that $\fV_{\thetabar}$ is flat as a right $\tR_{\thetabar}^{\op}$-module, or equivalently, flat as a left $\tR_{\thetabar}$-module.
This can be verified after base change to a quadratic unramified extension $\cO'/\cO$, hence we can assume without loss of generality that~$\thetabar$ has type~\ref{item: ss pseudorep}
or type~\ref{item: generic pseudorep}.
For type~\ref{item: ss pseudorep}, the flatness of~$\fV_{\thetabar}$ over~$\tR_{\thetabar}$
   is clear from Section~\ref{subsubsec:ss CWE stack}.
  Indeed, we
  have $\fX_{\thetabar} \cong [\Spec R_{\thetabar}^{\ps}/\bm{\mu}_2]$, while $\fV_{\thetabar} \cong R_{\thetabar}^{\ps}(1) \oplus R_{\thetabar}^{\ps}(1)$, and 
  $\tld R_{\thetabar} \cong M_2(R_{\thetabar}^{\ps})$ with its standard left action on~$\fV_{\thetabar}$ (thought of as column vectors).
Thus~$\fV_{\thetabar}$ is a direct summand of $\tld R_{\thetabar}$ and is in
particular $\tR_{\thetabar}$-flat, hence $\tR_{\thetabar}$-flat. 
For type~\ref{item: generic pseudorep}, the flatness of~$\fV_{\thetabar}$ over~$\tR_{\thetabar}$ 
is~\cite[Prop.\ 5.3.1]{JNWE}.

Finally, assume that~$\thetabar$ has type~\ref{item: non p-distinguished pseudorep}. 
Let $p: \Spec A \to \fX_{\thetabar}$ be the presentation described in~\cite[Section~3.3]{JNWE} (and recalled in Section~\ref{subsubsec:scalar CWE stack}). %
Write $V_{\thetabar} \coloneqq  p^*\fV_{\thetabar}$.
Since $p^*$ is exact and faithful, it suffices to prove that the composite 
\[
V_{\thetabar} \otimes_{\tR_{\thetabar}^{\op}}\text{--}: \Mod^{\fl}(\tR^{\op}_{\thetabar}) \to \Mod^{\fp}(\tR^{\op}_{\thetabar}) 
\xrightarrow{\mathfrak{V}_{\thetabar} \otimes_{\tR^{\op}_{\thetabar}}\text{--}} \Coh(\fX_{\thetabar}) \xrightarrow{p^*} \Mod^{\fp}(A)
\]
is exact.
This follows from
the proof of 
\cite[Prop.~5.4.2]{JNWE},
which shows that if~$N$ is the unique simple left $\tR^{\op}_{\thetabar}$-module, then 
$\Tor_1^{\tR^{\op}_{\thetabar}}(V_{\thetabar}, N) = 0.$
(More precisely, \emph{loc.\ cit.} writes~$E$ for~$\tR_{\thetabar}$, and proves that $\Tor_1^E(p^*\cV^*, M) = 0$, where~$M$ is the unique simple $E$-module,
$\cV$ is $\fV_{\thetabar}$, and~$\cV^*$ is the coherent dual of~$\fV_{\thetabar}$, viewed as a right $\tR_{\thetabar}$-module via the dual action.
Now~\cite[Prop.\ 2.2.4]{JNWE} shows that $\cV^*$ with the dual right $\tR_{\thetabar}$-action is isomorphic to the pullback of~$\cV$ under 
$\dagger : \tR_{\thetabar}^{\op} \isoto \tR_{\thetabar}$, and so we deduce that 
$\Tor_1^{\tR_{\thetabar}^{\op}}(p^*\cV, M^\dagger) = 0$.
Since~$M^\dagger$ is the unique simple left $\tR^{\op}_{\thetabar}$-module, this is what we wanted.)
\end{proof}

\begin{thm}\label{thm:JNWE-properties-of-V}%
Let~$\thetabar$ be a 2-dimensional $\cbF_p$-valued pseudorepresentation of~$G_{\bQ_p}$.
\begin{enumerate}
\item\label{item:20} The natural map
\[  
\Rtilde_{\thetabar}\to \REnd_{\Pro D^b_{\coh}(\cX_{\thetabar})}(\cV_{\thetabar}) 
\]
is an isomorphism, and
the functor
\[
\cV_{\thetabar} \otimes^L_{\tR_{\thetabar}^{\op}}\text{--}: D^b_{\fp}(\tR_{\thetabar}^{\op}) \to \Pro D^b_{\coh}(\cX_{\thetabar})
\]
defined in {\em Lemma~\ref{lem:derived-tensor-product-EW}}~\eqref{item:87}, is fully faithful.
\item\label{item:22} If~$\thetabar$ does {\em not} have type~{\em \ref{item: Steinberg pseudorep}},
the composite
\[
D^b_{\fl}(\tR_{\thetabar}^{\op}) \xrightarrow{\eqref{eqn:fl to fp}} D^b_{\fp}(\tR_{\thetabar}^{\op}) \xrightarrow{\cV_{\thetabar} \otimes^L_{\tR_{\thetabar}^{\op}}\text{--}} \Pro D^b_{\coh}(\cX_{\thetabar})
\]
is $t$-exact.
\end{enumerate}
\end{thm}
\begin{proof}%
  By Lemma~\ref{lem:fp-tensor-full-faithfulness-criterion}, the second statement of part~(1) is a consequence of the first statement.
  By the discussion in Section~\ref{subsubsec:not abs irred CWE stack}, 
  we can assume without loss of generality that~$\F = \F_{\thetabar}$.
  
  We begin by proving the first part.
  By Theorem~\ref{thm:coherent-completeness-CWE-stacks},
  the completion map 
  \[
  \REnd_{D^b_{\coh}(\fX_{\thetabar})}(\fV_{\thetabar}) \to \REnd_{\Pro D^b_{\coh}(\cX_{\thetabar})}(\cV_{\thetabar}) 
  \]
  is an isomorphism.
  So it suffices to prove that
  the natural map 
  \begin{equation}\label{eqn:to prove isomorphism for REnd}
  \Rtilde_{\thetabar}\to\REnd_{D^b_{\coh}(\fX_{\thetabar})}(\fV_{\thetabar})
  \end{equation}
  is an isomorphism. %
  If~$\cO'/\cO$ is a quadratic unramified extension, then the formation of
$\Ext^i_{D^b_{\coh}(\fX_{\thetabar})}$ commutes with base-change to $\cO'$ %
(as can be seen e.g.\ by computing with free resolutions after pullback
  by a smooth surjection $p: \Spec S \to \fX_{\thetabar}$).
  Thus we can assume without loss of generality that
  $\thetabar$ does not have type~\ref{item: not abs irred pseudorep}.
  It then follows from~\cite[Thms.\ 3.1.1, 3.2.3,
  3.3.1, 3.5.1]{JNWE} that~\eqref{eqn:to prove isomorphism for REnd} is an isomorphism. 
  This concludes the proof of the first part.

  We now prove the second part. 
  By Lemma~\ref{upper triangle in diagram}, it suffices to prove that the functor
  \[
  \cV_{\thetabar} \cotimes_{\tR_{\thetabar}^{\op}}\text{--} : \Pro D^b_{\fl}(\tR_{\thetabar^{\op}}) \to \Pro D^b_{\coh}(\cX_{\thetabar}),
  \]
  defined in Lemma~\ref{lem:derived-tensor-product-EW}~\eqref{item:94}, is $t$-exact.
  By Lemma~\ref{topological flatness and t-exactness}, it suffices to prove that $\cV_{\thetabar}$ is a topologically flat $\tR_{\thetabar}^{\op}$-module.
  This is Proposition~\ref{topological flatness of V}.\qedhere
\end{proof}%

\subsubsection{Completion at closed points}
\label{subsubsec:completion at closed points}
In Section~\ref{finiteness for Functor} we will need 
the following results about coherent sheaves on versal rings to~$\cX$ at closed points.

Recall from Section~\ref{subsubsec:points of X}
that the closed points of~$\cX$ are in bijection with $\Gal(\cbF_p/\F)$-conjugacy classes of semisimple representations $\rhobar: G_{\bQ_p} \to \GL_2(\cbF_p)$ with determinant~$\zeta\varepsilon^{-1}$,
or equivalently, $\Gal(\cbF_p/\F)$-conjugacy classes of 2-dimensional $\cbF_p$-valued pseudorepresentations~$\thetabar$.
Given a 2-dimensional $\cbF_p$-valued pseudorepresentation~$\thetabar$, we will adopt without further comments the notation in Definition~\ref{defn:various versal rings} for versal rings.

\begin{defn}\label{V_sigma}
For any Serre weight~$\sigma$ and $2$-dimensional $\cbF_p$-valued 
pseudorepresentation~$\thetabar$, we introduce the notation %
\[
\Vsigma \coloneqq  (\Pro V)(\ihat_{\thetabar}^* \cInd_{KZ}^G \sigma) \in \Pro \Mod^{\fl}(\tR_{\thetabar}) \isoto \Mod_c(\tR_{\thetabar}),
\]
where we write $\Pro V : \Pro \cA_{\thetabar}^{\fp} \to \Pro \Mod^{\fl}(\tR_{\thetabar})$ for the $\Pro$-extension of~\eqref{refactored V Rthetabar}. 
Similarly, we write 
\[
\Vsigma^\dagger \coloneqq  (\Pro V^\dagger)(\ihat_{\thetabar}^* \cInd_{KZ}^G \sigma) \in \Pro \Mod^{\fl}(\tR_{\thetabar}^{\op}) \isoto \Mod_c(\tR_{\thetabar}^{\op}).
\]
\end{defn}

The $\tR_{\thetabar}$-modules~$\Vsigma$ have been computed in~\cite{KisinFM} following work of Berger--Breuil~\cite{MR2642406}; we recall these results in the next lemma.

\begin{lemma}\label{lem: computing Vsigma}
Let~$\sigma = \sigma_{a, b}$ be a Serre weight, and let~$\thetabar$ be a 2-dimensional $\cbF_p$-valued pseudorepresentation.
Assume that~$\Vsigma \ne 0$, and that~$\thetabar$ does not have type~{\em \ref{item: not abs irred pseudorep}} {\em (}i.e.\ it is $\F$-valued with absolutely irreducible summands{\em )}.
Then:
\begin{enumerate}
\item If $\thetabar$ has type~{\em\ref{item: ss pseudorep}}, then $\Vsigma$ is $\tR_{\thetabar}$-linearly isomorphic to the
Galois representation $\Vver_{\thetabar} \otimes_{R_{\thetabar}^{\ver}} \Rsigma_{\thetabar}/\varpi$.

\item If~$\thetabar = \nr_\lambda \omega^{a+b} + \nr_{\lambda^{-1}\zeta(p)} \omega^{a-1}$,
then~$M_{\sigma, \thetabar}$ is $\tR_{\thetabar}$-linearly isomorphic to a free $\bF\llbracket S\rrbracket $-module of rank one on which~$\tR_{\thetabar}$
acts via the character %
\begin{equation}\label{Galois action on V_sigma}
\nr_{S+\lambda}\omega^{a+b} : \tR_{\thetabar} \to \bF\llbracket S\rrbracket .
\end{equation}
The restriction of~\eqref{Galois action on V_sigma} to $R^{\ps}_{\thetabar} = Z(\tR_{\thetabar})$ is a ring homomorphism
\begin{equation}
\label{eqn:beta map}
\beta : R_{\thetabar}^{\ps} \to \bF\llbracket S\rrbracket ,
\end{equation}
which is surjective whenever $\thetabar$ does not have type~{\em \ref{item: non p-distinguished pseudorep}}.
\end{enumerate}
\end{lemma}
\begin{proof}
Part~(1) follows from \cite[Lemma~1.5.3]{KisinFM} (and its proof).
Part~(2) follows from \cite[Lemma~1.5.9, Lemma~1.5.11]{KisinFM}, where the morphism~$\beta$ is denoted~$\theta$.
\end{proof}

\begin{rem}\label{rem:what beta classifies}
Using the description~\eqref{eqn:Rps in tR} of the central 
embedding $R^{\ps}_{\thetabar} \to \tR_{\thetabar}$, one sees that~$\beta$ classifies the deformation
\[
\nr_{S+\lambda}\omega^{a+b} + \nr_{(S+\lambda)^{-1}\zeta(p)}\omega^{a-1}
\]
of~$\thetabar$ to~$\bF\llbracket S\rrbracket $.
\end{rem}

We now deduce an alternative description of $\Vsigma$.
It will be useful to employ the following notation: if~$J \subset \tR_{\thetabar}$ is a subset, we will write $\langle J \rangle$ for the two-sided ideal
generated by~$J$, and $J\tR_{\thetabar}$ for the right ideal generated by~$J$.
We also note that in part~(2) of the following lemma, which encompasses the cases 
when~$\thetabar$ has type~\ref{item: generic pseudorep} or~\ref{item: Steinberg pseudorep},
we employ the explicit presentation of~$\tR_{\thetabar}$ given above
in Section~\ref{subsubsec:gen CWE stack}, resp.\ Section~\ref{subsubsec:Steinberg CWE stack}.

\begin{lemma}\label{V of completed compact induction}Let~$\sigma = \sigma_{a, b}$ be a Serre weight, and let~$\thetabar$ be a 2-dimensional $\cbF_p$-valued pseudorepresentation.
\begin{enumerate}
\item If~$\thetabar$ has type~{\em \ref{item: ss pseudorep}}, and $\Vsigma \ne 0$, then $\Vsigma$ is $R_{\thetabar}^{\ver}$-linearly isomorphic to $(\Rsigma_{\thetabar}/\varpi)^{\oplus 2}$.
\item Suppose that~$0 \leq b \leq p-2$, let $\lambda \in \F^\times$, and assume that $(b, \lambda) \ne (p-2, \pm \zeta(p)^{1/2})$.
Let $\sigma = \sigma_{a, b}$ and $\thetabar = \nr_\lambda \omega^{a+b} + \nr_{\lambda^{-1}\zeta(p)} \omega^{a-1}$.
Then $\Vsigma \cong \tR_{\thetabar}/\langle \ker \beta, e_{22} \rangle$,
and $M_{\sigmacomp, \thetabar} \cong \tR_{\thetabar}/\langle \ker \beta, e_{11} \rangle$.
\end{enumerate}
\end{lemma}
\begin{proof}
Recalling that when $\thetabar$ has type~\ref{item: ss pseudorep} we have
$R^{\ps}_{\thetabar} = R^{\ver}_{\thetabar}$,
so that $\tR_{\thetabar}$ is in particular an~$R^{\ver}_{\thetabar}$-algebra,
we see that 
part~(1) follows from 
Lemma~\ref{lem: computing Vsigma}~(1).

We now prove part~(2).
By Lemma~\ref{lem: computing Vsigma}, we have $\Vsigma = \tR_{\thetabar}/\ker(\nr_{S+\lambda}\omega^{a+b})$.
Since~$\bF\llbracket S\rrbracket $ is a local ring, at least one of the orthogonal idempotents 
$e_{11}$ and~$e_{22}$ goes to zero under $\nr_{S+\lambda}\omega^{a+b}$, 
and since the kernel is a two-sided ideal, it contains at most one of~$e_{11}, e_{22}$.
Hence there exists~$i = 1$ or~$2$ such that $\langle \ker \beta, e_{ii} \rangle \subset \ker(\nr_{S+\lambda}\omega^{a+b})$.
The inclusion must then be an equality, since
$\tR_{\thetabar}/\langle \ker \beta, e_{ii} \rangle \cong R^{\ps}_{\thetabar}/\ker(\beta) = \bF\llbracket S\rrbracket $.

It remains to prove that~$i = 2$.
If $e_{22} : \Vsigma \to \Vsigma$ is not zero, then there exists a nonzero $\tR_{\thetabar}$-linear map
$\tR_{\thetabar}e_{22} \to \Vsigma$. 
Bearing in mind that our presentation of $\tR_{\thetabar}$ is such that
$\tR_{\thetabar}e_{22}$ is a projective envelope of~$\nr_{\lambda^{-1}\zeta(p)}\omega^{a-1}$,
this contradicts our assumption that either $\omega^{a+b} \ne \omega^{a-1}$ (i.e.\ $b \ne p-2$) or
$\lambda \ne \lambda^{-1}\zeta(p)$ (when $b = p-2$).
\end{proof}

\begin{rem}\label{explicit kernel of theta}\leavevmode
If~$\thetabar = 1+ \omega^{-1}$,
then the kernel of~$\beta$ is~$(\varpi, a_1, X_0, X_1)$.
This is proved
in~\cite[Lemma~3.9]{HuTan}, bearing in mind that $(a_0, a_1, X_0, X_1)$ 
are denoted~$(d_0, -d_1, c_0, c_1)$ in \emph{loc.\ cit.}, as explained in~\cite[Section~5.5]{JNWE}.
\end{rem}

\begin{rem}\label{a_0 and T}
If $\thetabar = 1 + \omega^{-1}$ we will also need the following connection with the Hecke operator~$\HeckeT$.
If $\sigma \in \{\Sym^0, \Sym^{p-1}, \Sym^{p-3} \otimes \det\}$, then~$\HeckeT$ acts on~$\Vsigma$ by $V$-functoriality.
Lemma~\ref{lem: computing Vsigma} provides an identification
$\End_{\Mod_c (\tR_{\thetabar})}(\Vsigma) = \bF\llbracket S\rrbracket$,
and~$\HeckeT-1$ is a uniformizer of this ring:
in fact, the exactness of $\ihat^*_{\thetabar}$ and~$\Pro V$ implies that
\[
\Vsigma/(\HeckeT-1)\Vsigma = V\bigl(\ihat^*_{\thetabar}\cInd_{KZ}^G \sigma/(\HeckeT-1)\cInd_{KZ}^G \sigma\bigr),
\]
and hence the left-hand side is a one-dimensional $\F$-vector space;
since $M_{\sigma,\thetabar}$ is free of rank one over~$\F\llbracket S \rrbracket$
(by Lemma~\ref{lem: computing Vsigma}), we conclude that indeed $\HeckeT-1$
is identified with a uniformizer of~$\F\llbracket S \rrbracket$.
On the other hand, the map
\[
\beta : R_{\thetabar}^{\ps} \to \End_{\Mod_c (\tR_{\thetabar})}(\Vsigma) = \F\llbracket S \rrbracket,
\]
describing the restriction of the $\tR_{\thetabar}$-action on~$\Vsigma$ to $Z(\tR_{\thetabar}) = R_{\thetabar}^{\ps}$
is surjective, by Lemma~\ref{lem: computing Vsigma}, and so it sends~$a_0$ to a uniformizer of~$\F\llbracket S \rrbracket$, by Remark~\ref{explicit kernel of theta}.
Putting the preceding observations together, we conclude that
$\beta(a_0)$ is an $\F\llbracket S \rrbracket^\times$-multiple of~$\HeckeT-1$.
\end{rem}

Finally, we have the right $\tR_{\thetabar}$-modules $\Vsigma^\dagger \coloneqq  \Pro V^\dagger(\ihat_{\thetabar}^* \cInd_{KZ}^G \sigma)$, which are obtained from~$\Vsigma$ by
composing with the anti-involution~$\dagger$.

\begin{lemma}\label{V of completed compact induction II}Let~$\sigma = \sigma_{a, b}$ be a Serre weight, and let~$\thetabar$ be a 2-dimensional $\cbF_p$-valued pseudorepresentation.
\begin{enumerate}
\item If~$\thetabar$ has type~{\em \ref{item: ss pseudorep}}, then $\Vsigma^\dagger$ is $R_{\thetabar}^{\ver}$-linearly isomorphic to $(\Rsigma_{\thetabar}/\varpi)^{\oplus 2}$.
\item Assume that~$0 \leq b \leq p-2$, let $\lambda \in \F^\times$, and assume that $(b, \lambda) \ne (p-2, \pm \zeta(p)^{1/2})$.
Let $\thetabar = \nr_\lambda \omega^{a+b} + \nr_{\lambda^{-1}\zeta(p)} \omega^{a-1}$.
Then $\Vsigma^\dagger \cong \tR_{\thetabar}/\langle \ker \beta, e_{11}\rangle$, and
$M_{\sigmacomp, \thetabar}^\dagger \cong \tR_{\thetabar}/\langle \ker \beta, e_{22}\rangle$
(as right $\tR_{\thetabar}$-modules).
\end{enumerate}
\end{lemma}
\begin{proof}
By definition, if~$g \in \tR_{\thetabar}^\times$ is the image of an element of $G_{\bQ_p}$, then $g^\dagger = (\zeta\varepsilon^{-1})(g)g^{-1}$.
Our assumptions imply that~$\thetabar$ does not have type~\ref{item: non p-distinguished pseudorep} or~\ref{item: not abs irred pseudorep}, and so $\tR_{\thetabar}$ is an order in $M_2(R^{\ps}_{\thetabar})$,
as described in Section~\ref{subsec:CWE stacks}.
It follows from this description of~$\tR_{\thetabar}$ that~$\dagger$ coincides, at elements of $G_{\bQ_p}$, with the restriction to~$\tR_{\thetabar}$ of the anti-involution
\begin{equation}\label{eqn:canonical involution}
\fourmatrix a b c d \mapsto \fourmatrix{d}{-b}{-c}{a}
\end{equation}of  $M_2(R^{\ps}_{\thetabar})$.
Since the images of elements of~$G_{\bQ_p}$ are dense in~$\tR_{\thetabar}$, we conclude that~$\dagger$ coincides with~\eqref{eqn:canonical involution}.
Thus the lemma follows from Lemma~\ref{V of completed compact induction}. %
\end{proof}

\begin{rem}\label{rem:tensoring with V_sigma}
In case~(2) of Lemma~\ref{V of completed compact induction II}, it follows from~\eqref{Galois action on V_sigma}
that the inertia group~$I_{\bQ_p}$ acts on $\Vsigma^\dagger$ via the character $(\omega^{a+b})^\dagger|_{I_{\bQ_p}} = \omega^{a-1}|_{I_{\bQ_p}}$.
Hence, if~$W$ is a left $\tR_{\thetabar}$-module, then $W \cotimes_{\tR^{\op}_{\thetabar}} \Vsigma^\dagger$ is an $\tR_{\thetabar}$-quotient of~$W$ on which~$I_{\bQ_p}$ acts by~$\omega^{a-1}$.
\end{rem}

\begin{prop}\label{prop:support of M(sigma)}
Let~$\sigma = \sigma_{a, b}$ be a Serre weight, and let $\thetabar$ be a $2$-dimensional $\cbF_p$-valued pseudorepresentation. 
Let~$R^{\ver}_{\thetabar}, \Rsigma_{\thetabar}$ be the versal rings defined in
Definition~{\em \ref{defn:various versal rings}}.
Assume that $b \ne 0$, and $\Vsigma^\dagger \ne 0$.
Then 
\begin{enumerate}
\item $\Vver_{\thetabar} \cotimes_{\tR_{\thetabar}^{\op}} \Vsigma^\dagger$ has scheme-theoretic support $\Spec \Rsigma_{\thetabar}/\varpi$, and
\item if~$\thetabar$ does not have type~{\em \ref{item: non p-distinguished pseudorep}}, then
$\Vver_{\thetabar} \cotimes_{\tR_{\thetabar}^{\op}} \Vsigma^\dagger \cong \Rsigma_{\thetabar}/\varpi$.
\end{enumerate}
\end{prop}

\begin{rem}\label{rem:Vsigma change of weight}
The assumption that~$b \ne 0$ in Proposition~\ref{prop:support of M(sigma)} is due to the fact that, for any~$\thetabar$, we have an exact sequence
\[
0 \to \ihat_{\thetabar}^* \cInd_{KZ}^G \sigma_{a, p-1} \to \ihat_{\thetabar}^* \cInd_{KZ}^G \sigma_{a, 0} \to \ihat^*_{\thetabar}\bigl(\omega^a\nr_{X^2-\zeta(p)}\circ \det\bigr) \to 0 
\]
obtained by applying~$\ihat_{\thetabar}^*$ to~\eqref{trivial change of weight}.
Since $V^\dagger$ is zero on $\SL_2(\bQ_p)$-invariant representations of~$G$, and~$V^\dagger$ is exact, we deduce that $M_{\sigma_{a, p-1}, \thetabar}^\dagger = M_{\sigma_{a, 0}, \thetabar}^\dagger$.
Hence $\Vver_{\thetabar} \cotimes_{\tR_{\thetabar}^{\op}} V_{\sigma_{a, 0}}^\dagger \cong R^{\sigma_{a, p-1}}_{\thetabar}/\varpi$, 
  which can be larger than $R^{\sigma_{a, 0}}_{\thetabar}/\varpi$.
\end{rem}

\begin{proof}[Proof of Proposition~{\em \ref{prop:support of M(sigma)}}]
By definition, the Galois representation $\Vver/\fm \Vver$ is semisimple with absolutely irreducible summands.
Hence, after possibly replacing~$\F$ with a finite extension, we can assume that~$\thetabar$ does not have type~\ref{item: not abs irred pseudorep}.
Furthermore, by Lemma~\ref{lem:localization-of-cInd-sigma}~\eqref{item:support of cInd sigma}, $\Vsigma^\dagger \ne 0$ if and only if~$\sigma$ is a Serre weight of $\Vver/\fm \Vver$.
We now proceed by cases, according to the type of~$\thetabar$.

Case~\ref{item: ss pseudorep}: Recall that $\tR_{\thetabar} = M_2(R^{\ver}_{\thetabar})$ acting on $\Vver_{\thetabar} = (R^{\ver}_{\thetabar})^{\oplus 2}$ in the standard representation. 
Hence $\Vsigma^\dagger$ is $R^{\ver}_{\thetabar}$-linearly isomorphic to the direct sum of 
$\Vver_{\thetabar} \cotimes_{\tR_{\thetabar}^{\op}} \Vsigma^\dagger$ with itself.
On the other hand, by Lemma~\ref{V of completed compact induction II}~(1), $\Vsigma^\dagger$ is $R^{\ver}_{\thetabar}$-linearly isomorphic to 
$(\Rsigma_{\thetabar}/\varpi)^{\oplus 2}$, and so $\Vver_{\thetabar} \cotimes_{\tR_{\thetabar}^{\op}} \Vsigma^\dagger = \Rsigma_{\thetabar}/\varpi$,
as desired.

Case~\ref{item: generic pseudorep}: After a twist, we may assume that $\sigma = \Sym^b$ and 
$\thetabar = \nr_\lambda\omega^b + \nr_{\lambda^{-1}\zeta(p)}\omega^{-1}$, for some~$\lambda \in \bF^\times$
and~$0 \leq b \leq p-2$.
Then Lemma~\ref{V of completed compact induction II}~(2) shows that
\[
\Vver_{\thetabar} \cotimes_{\tR_{\thetabar}^{\op}} \Vsigma^\dagger = \Vver_{\thetabar}/(\ker \beta, e_{11}, ce_{21})\Vver_{\thetabar},
\]
where we have used that $\langle \ker \beta, e_{11}\rangle = (\ker \beta, e_{11}, ce_{21})\tR_{\thetabar}$.
Since $\Vver_{\thetabar}/e_{11}\Vver_{\thetabar}$
is free of rank one as an $R^{\ver}_{\thetabar}$-module,
we see that
$\Vver_{\thetabar} \cotimes_{\tR_{\thetabar}^{\op}} \Vsigma^\dagger$ is a cyclic $R^{\ver}_{\thetabar}$-module with annihilator $(\ker\beta, c)R^{\ver}_{\thetabar}$.

The proposition will thus follow in this case once we show that
$R^{\ver}_{\thetabar}/(\ker \beta, c)R^{\ver}_{\thetabar}$ coincides with
$\Rsigma_{\thetabar}/\varpi$.
Since these rings are abstractly isomorphic
(being power series rings over~$\bF$ in the same number of variables), 
it suffices to prove that 
the quotient map $R^{\ver}_{\thetabar} \to R^{\ver}_{\thetabar}/(\ker \beta, c)R^{\ver}_{\thetabar}$ factors through~$\Rsigma_{\thetabar}/\varpi$.
To this end, note that
we have an exact sequence of $R^{\ver}_{\thetabar}$-modules
\[
0 \to \Vver_{\thetabar}/(\ker \beta, c, e_{22}) \to \Vver_{\thetabar} \otimes_{R^{\ver}_{\thetabar}} R^{\ver}_{\thetabar}/(\ker\beta, c)R^{\ver}_{\thetabar} \to \Vver_{\thetabar}/(\ker \beta, c, e_{11}) \to 0
\]
which exhibits the free rank two $R^{\ver}_{\thetabar}/(\ker \beta, c)$-module
$\Vver_{\thetabar}/(\ker \beta, c)$ as an extension of free rank one modules.
The explicit form of~$\tR_{\thetabar}$ given in~\eqref{2x2 matrix order (gen)} shows that
this is actually a sequence of $\tR_{\thetabar}$-modules.
Furthermore, since the quotient 
$\Vver_{\thetabar}/(\ker \beta, c, e_{11})$
has been proven in the previous paragraph to be isomorphic to
$\Vver_{\thetabar} \cotimes_{\tR_{\thetabar}^{\op}} \Vsigma^\dagger$,
Remark~\ref{rem:tensoring with V_sigma} shows that $I_{\bQ_p}$ acts by~$\omega^{-1}$ on $\Vver_{\thetabar}/(\ker \beta, c, e_{11})$.
Since $\det \Vver_{\thetabar}$ is an unramified twist of $\omega^{b-1}$,
this shows that~$\Vver_{\thetabar} \otimes_{R^{\ver}_{\thetabar}} R^{\ver}_{\thetabar}/(\ker \beta, c)R^{\ver}_{\thetabar}$ is an ordinary deformation of 
$\nr_\lambda \omega^b \oplus \nr_{\lambda^{-1}\zeta(p)}\omega^{-1}$, 
and its classifying map $R^{\ver}_{\thetabar} \to R^{\ver}_{\thetabar}/(\ker \beta, c)R^{\ver}_{\thetabar}$ 
factors through $\Rsigma_{\thetabar}/\varpi$,
as desired.

Case~\ref{item: Steinberg pseudorep}: 
This is similar to case~\ref{item: generic pseudorep}.
After a twist, we may assume that $\thetabar = 1 + \omega^{-1}$ and $\sigma \in \{\Sym^{p-1}, \Sym^{p-3} \otimes \det\}$.
Bearing in mind the isomorphism $M_{\Sym^0, \thetabar}^\dagger \cong M_{\Sym^{p-1}, \thetabar}^\dagger$ 
from Remark~\ref{rem:Vsigma change of weight},
the module $\Vsigma^\dagger$ has thus been computed in Lemma~\ref{V of completed compact induction II}.
Taking into account the explicit form of~$\ker \beta$ given in Remark~\ref{explicit kernel of theta}, we find that
\begin{gather}
\Vver_{\thetabar} \cotimes_{\tld R_{\thetabar}^{\op}} M_{\Sym^{p-1}, \thetabar}^\dagger = %
R^{\ver}_{\thetabar}/(\varpi, a_1, X_0, X_1, c)R^{\ver}_{\thetabar} = 
\bF\llbracket a_0, b_0, b_1\rrbracket /(a_0b_1)\label{eqn:computing Vsigma Steinberg}\\
\Vver_{\thetabar} \cotimes_{\tld R_{\thetabar}^{\op}} M_{\Sym^{p-3} \otimes \det, \thetabar}^\dagger = %
R^{\ver}_{\thetabar}/(\varpi, a_1, X_0, X_1, b_0, b_1)R^{\ver}_{\thetabar} = \bF\llbracket a_0, c\rrbracket \label{eqn:computing Vsigma Sym^p-3}
\end{gather}
as $R^{\ver}_{\thetabar}$-modules.
Arguing as in case~\ref{item: generic pseudorep}, and using the fact 
that $R^{(p, 0), \crys}/\varpi$ is a power series ring over $\bF\llbracket x, y\rrbracket /(xy)$, %
we thus see that $V \cotimes_{\tR_{\thetabar}^{\op}} \Vsigma^\dagger \cong \Rsigma_{\thetabar}/\varpi$, as desired. %
(This presentation of $R^{(p, 0), \crys}/\varpi$ can be deduced from a more careful analysis of the presentation of~$\fX_{\thetabar}$ 
described in Section~\ref{subsubsec:Steinberg CWE stack}, or alternatively,
it can be deduced from~\cite[Lemma~7.3.7]{BCGP}.)

Case~\ref{item: non p-distinguished pseudorep}: 
After a twist, we may assume that $\thetabar = \omega^{-1} + \omega^{-1}$, and $\sigma = \Sym^{p-2}$. 
This case was excluded in Lemma~\ref{V of completed compact induction II}, so we will argue in a different way, by applying the results of~\cite{MR3556448}.
In Definition~\ref{defn: companion Serre weights} we have defined an irreducible $E[KZ]_\zeta$-module~$\Theta(\sigmasigmacomp)$ 
whose mod $\varpi$ reduction is isomorphic to~$\sigma$,
and we now fix an $\cO[KZ]_\zeta$-lattice~$\Theta$ in~$\Theta(\sigmasigmacomp)$. 
We have isomorphisms
\begin{multline*}
\Vver_{\thetabar} \cotimes_{\tR_{\thetabar}^{\op}} \Vsigma^\dagger \isoto 
\Vver_{\thetabar} \cotimes_{\tE_{\thetabar}} (\tP_{\thetabar} \ccotimes_{\cO\llbracket G\rrbracket _\zeta} \ihat_{\thetabar, *}\ihat^*_{\thetabar}\cInd_{KZ}^G \sigma)\\ 
\isoto \Vver_{\thetabar} \cotimes_{\tE_{\thetabar}} (\tP_{\thetabar} \otimes_{\cO\llbracket G\rrbracket _\zeta} \cInd_{KZ}^G \sigma) 
\isoto \Hom_{KZ}^{\cont}(\Vver_{\thetabar} \cotimes_{\tE_{\thetabar}} \tP_{\thetabar}, \sigma^\vee)^\vee,
\end{multline*}
where the first arrow is Lemma~\ref{needed for universal deformation}~(2), the second arrow is Lemma~\ref{tensor and completion}, and the third arrow is Lemma~\ref{lem:tensor facts}.
Now~$\Vver_{\thetabar} \cotimes_{\tE_{\thetabar}} \tP_{\thetabar}$ is the object denoted~$N$ 
in~\cite{MR3556448} (and we employ this same notation in what follows),
and so we see that $\Vver_{\thetabar} \cotimes_{\tR_{\thetabar}^{\op}} \Vsigma^\dagger$ 
is the module denoted~$M^\square(\sigma)$ in~\cite[(2.12)]{MR3556448}.
Following~\emph{loc.\ cit.} and~\cite[Section~2]{PaskunasBM}, we also consider the module
\[
M(\Theta) \coloneqq  \Hom_{\cO}(\Hom_{KZ}^{\cont}(N, \Theta^d), \cO).
\]
Since $\Vver_{\thetabar} \cotimes_{\tE_{\thetabar}} -$ is exact, by Proposition~\ref{topological flatness of V}, 
and $\tP_{\thetabar} \cotimes_{\cO\llbracket KZ\rrbracket _\zeta} -$ is exact, by Lemma~\ref{lem:projective generators are flat over their endos}~\eqref{item:projective 3}, 
we see that $N$ is topologically flat in~$\Mod_c(\cO\llbracket KZ\rrbracket _\zeta)$,
hence projective in $\Mod_c(\cO\llbracket KZ\rrbracket _\zeta)$, by~\ref{lem: topologically flat implies projective}.
By~\cite[Lemma~2.14]{PaskunasBM},
we conclude that
$M^\square(\sigma)$ is isomorphic to $M(\Theta)/\varpi M(\Theta)$.
We claim that $M(\Theta)$ is a faithful finitely presented $\Rsigma_{\thetabar}$-module.

Assuming the claim, we conclude the proof of the proposition as follows.
By the claim, the natural map 
\[
\Rsigma_{\thetabar} \to \End_{R^{\ver}_{\thetabar}}(M(\Theta)) 
\]
is injective,
and so the natural map 
\[
\Rsigma_{\thetabar}/\varpi \to \End_{R^{\ver}_{\thetabar}}(M(\Theta)/\varpi M(\Theta)) 
\]
has nilpotent kernel.
Since $\Rsigma_{\thetabar}/\varpi$ is reduced, by Lemma~\ref{Z(sigma) as special fibre}, 
we conclude that the scheme-theoretic support of $M^\square(\sigma) = M(\Theta)/\varpi M(\Theta)$ is $\Rsigma_{\thetabar}/\varpi$, as desired.

To prove the claim, observe first that~$M(\Theta)$ is finitely presented over~$R^{\ver}_{\thetabar}$, by~\cite[Proposition~2.15]{PaskunasBM}.
Then let $\fa \coloneqq  \Ann_{R^{\ver}_{\thetabar}} M(\Theta)$. 
As shown during the proof of~\cite[Theorem~4.10]{MR3556448}, the radical of~$\fa$ is the kernel of $R^{\ver}_{\thetabar} \to \Rsigma_{\thetabar}$.
However, it is also shown in~\emph{loc.\ cit.} that the assumptions of~\cite[Theorem~2.1]{MR3556448}, or equivalently~\cite[Theorem~2.42]{PaskunasBM}, are met.
Hence $\fa$ is a radical ideal, which concludes the proof.
\end{proof}

\begin{cor}\label{cor:explicit presentation of crystalline deformation ring II}
  Let~$\thetabar = 1 + \omega^{-1}$.
  Then $R_{\thetabar}^{(p, 0), \crys}/\varpi = \F\llbracket a_0, b_0, b_1\rrbracket /(a_0b_1)$, and $R_{\thetabar}^{(1, 0), \crys}/\varpi = R_{\thetabar}^{(p, 0), \crys}/(\varpi, b_1)$.
\end{cor}
\begin{proof}
The claim about $R_{\thetabar}^{(p, 0), \crys}/\varpi$ is an immediate consequence of~\eqref{eqn:computing Vsigma Steinberg}.
We now prove the claim about $R_{\thetabar}^{(1, 0), \crys}/\varpi$.
It follows from e.g.\ \cite[Theorem~5.11]{HuTan}
that the quotient map $R_{\thetabar}^{\ver}/\varpi \to R_{\thetabar}^{(1, 0), \crys}/\varpi$ factors through a surjection 
$R_{\thetabar}^{(p, 0), \crys}/\varpi \to R_{\thetabar}^{(1, 0), \crys}/\varpi$ 
whose kernel is a minimal prime, and so it is generated by either~$a_0$ or~$b_1$.
But it cannot be generated by~$a_0$: 
if this were the case, 
then the map $R_{\thetabar}^{\ps} \to R_{\thetabar}^{(p, 0), \crys}/\varpi$ would factor through~$\bF$, 
since the maximal ideal of~$R^{\ps}_{\thetabar}$ is generated by~$(a_0, a_1, X_0, X_1)$.
This would contradict the existence of lifts of $1 \oplus \omega^{-1}$ to~$\bF[\epsilon]/\epsilon^2$
that factor through $R_{\thetabar}^{(p, 0), \crys}/\varpi$ and have nontrivial pseudocharacter 
(for example, the lift $\nr_{1-\epsilon} \oplus \nr_{1+\epsilon}\omega^{-1}$).
\end{proof}

\subsection{Complements in the Steinberg case}
\label{subsec:steinberg computations}
We assume throughout this section that $\thetabar$ is of type~\ref{item: Steinberg block}.
Making a twist if necessary, we assume that in fact $\thetabar = 1 + \omega^{-1}$.
We will shorten notation to $R \coloneqq  R^{\ps}_{\thetabar}$.
Recall that in Section~\ref{subsubsec:Steinberg CWE stack}
we have given an explicit description 
of~$R, \tld R_{\thetabar}, \fX_{\thetabar}$, and of the bundle $\fV_{\thetabar}$ lying over~$\fX_{\thetabar}$, together with its
left $\tR_{\thetabar}$-action.  
Using the notation of that section, %
we write $$I \coloneqq  (X_0,X_1)R \subseteq R;$$
this is the reducibility ideal in~$R$.

\subsubsection{The projective generator~\texorpdfstring{$\bP_{\thetabar}$}{P theta} of \texorpdfstring{$\fC_{\thetabar}$}{C theta}}
\label{subsubsec:Steinberg complements}
Recall from Definition~\ref{def:P theta} and Definition~\ref{defn: tld P} the projective objects $P_{\thetabar}, \tP_{\thetabar} \in \fC_{\thetabar}$.
Neither of these %
is a projective generator,
which leads us to introduce yet another projective object.
 
\begin{defn}\label{defn: bP in Steinberg case}%
Let~$\thetabar = 1 + \omega^{-1}$.
We choose a $\fC_{\thetabar}$-projective envelope~$P_{\triv_G^\vee}$ of the trivial character of~$G$,
we define the projective generator (with finite cosocle)~$\bP_{\thetabar}$ of~$\fC_{\thetabar}$ as
$$\bP_{\thetabar} \coloneqq  \tld P_{\thetabar} \oplus P_{\triv_G^\vee},
$$
and we write $\bE_{\thetabar} \coloneqq  \End_{\fC_{\thetabar}}(\bP_{\thetabar})$.
\end{defn}
We can now apply the Morita theory of Lemma~\ref{lem:recapitulated Morita}, 
taking the projective generator in the statement of that lemma to be~$\bP_{\thetabar}$.
Then Lemma~\ref{lem:recapitulated Morita}~(3)
shows that $\Hom_{\fC_{\thetabar}}(\bP_{\thetabar}, \tP_{\thetabar})$
is of finite type over $\tE_{\thetabar}$, and~\eqref{eq:Pthetabar-tensor-pi-formula} gives an isomorphism
\begin{equation}
\label{eqn:Steinberg case rewrite}
\tP_{\thetabar} %
\iso
\Hom_{\fC_{\thetabar}}(\bP_{\thetabar},\tP_{\thetabar})
\otimes_{\bE_{\thetabar}} \bP_{\thetabar}.
\end{equation}

\subsubsection{The ideal~\texorpdfstring{$J$}{J} of \texorpdfstring{$\bE_{\thetabar}$}{E theta}}\label{subsubsec:the ideal J}

Let $\cO_{\triv_G^\vee}$ 
denote $\cO$ regarded as an object of~$\fC_{\thetabar}$ by 
letting~$G$ act trivially.
Since the $\fC_{\thetabar}$-cosocle of~$\cO_{\triv_G^{\vee}}$ is equal to~$\triv_G^\vee$, with
multiplicity one,
we find that $\cO_{\triv_G^{\vee}}$ is a quotient of $P_{\triv_G^\vee}$,
and that 
\[
\Hom_G^{\cont}(\bP_{\thetabar},\cO_{\triv_G^{\vee}}) =
\Hom_G^{\cont}(P_{\triv_G^{\vee}},\cO_{\triv_G^{\vee}})
\]
is free of rank one over~$\cO$.
Since $\cO_{\triv_G^{\vee}}$ is a quotient
of~$P_{\triv_G^{\vee}}$, which is in turn a quotient of~$\bP_{\thetabar}$,
we see that $\Hom_G^{\cont}(\bP_{\thetabar},\cO_{\triv_G^{\vee}})$
is also a quotient
of $\bE_{\thetabar}$.
Thus we have an isomorphism of $\bE_{\thetabar}^{\op}$-modules
$$\Hom_G^{\cont}(\bP_{\thetabar},\cO_{\triv_G^{\vee}}) \iso \bE_{\thetabar}^{\op}/J$$
for some closed right ideal $J$ of~$\bE_{\thetabar}$ with the property
that $\cO \iso \bE_{\thetabar}^{\op}/J$ (as~$\cO$-modules). 
More precisely, $J$ can be taken to be the annihilator of any surjection $\bP_{\thetabar} \to \cO_{\triv_G^\vee}$ in~$\fC_{\thetabar}$.

We claim that~$J$ is in fact a two-sided ideal.
In fact, since $\Hom_G^{\cont}(\tP_{\thetabar},\cO_{\triv_G^{\vee}}) = 0$, $J$ contains the two-sided ideal~$J'$ generated by the 
idempotent of~$\bE_{\thetabar}$ associated to~$\tP_{\thetabar}$.
Since the quotient $\bE_{\thetabar}/J'$ is commutative,
we see that
$J$ is a two-sided ideal, and 
$\cO \iso \bE_{\thetabar}^{\op}/J$ (as $\cO$-algebras).
It also follows from 
Lemma~\ref{lem:recapitulated Morita}~(1)
that we have the following isomorphism in~$\fC_{\thetabar}$:
$$\cO_{\triv_G^{\vee}} \iso (\bE^{\op}_{\thetabar}/J) \cotimes_{\bE_{\thetabar}} \bP_{\thetabar} 
\iso \bP_{\thetabar}/J\bP_{\thetabar}.$$ %

\begin{lemma}
\label{lem:SL_2 invariants via J torsion}
If $M$ is an object in $\Mod_c(\bE_{\thetabar}^{\op})$, then
the natural morphism
$$M[J] \cotimes_{\bE_{\thetabar}} \bP_{\thetabar}
\to
M\cotimes_{\bE_{\thetabar}} \bP_{\thetabar}$$
is injective, and induces an isomorphism
$$
M[J] \cotimes_{\bE_{\thetabar}} \bP_{\thetabar}
\iso
(M\cotimes_{\bE_{\thetabar}} \bP_{\thetabar})^{\SL_2(\Q_p)}.
$$
\end{lemma}
\begin{proof}
To ease notation, write $N \coloneqq  M\cotimes_{\bE_{\thetabar}} \bP_{\thetabar}$;
note that this is the object of $\fC_{\thetabar}$ that corresponds to
$M$ under the Morita equivalence of Lemma~\ref{lem:recapitulated Morita}~(1).
The symbol~$M[J]$ denotes the $J$-torsion in~$M$; since~$J$ is a left ideal of $\bE_{\thetabar}^{\op}$, this is also an instance of the construction in~ Remark~\ref{rem:Hom-version-of-EW}, i.e.
$M[J] = \Hom_{\bE_{\thetabar}^{\op}}^{\cont}(\bE_{\thetabar}^{\op}/J, M)$.

By Lemmas~\ref{technical remark for Steinberg block} and ~\ref{trivial action} below, 
$G$ acts trivially on $N^{\SL_2(\Q_p)}$, since the trivial character is the only square root of $\zeta^{-1}: \bQ_p^\times \to \cO^\times$ contained in~$\fC_{\thetabar}$.  
Hence 
$N^{\SL_2(\bQ_p)}$
is the image of
the evaluation map
\begin{equation}
\label{eqn:SL_2 inv embedding}
\Hom_{\fC_{\thetabar}}(\cO_{\triv_G^{\vee}}, N) \otimes_{\cO} \cO_{\triv_G^{\vee}} \hookrightarrow N.
\end{equation}
The Morita equivalence gives an isomorphism
$$\Hom_{\fC_{\thetabar}}(\cO_{\triv_G^{\vee}}, N) \iso \Hom^{\cont}_{\bE_{\thetabar}^{\op}}(\bE_{\thetabar}^{\op}/J,
M) = M[J],$$
and so we may rewrite~\eqref{eqn:SL_2 inv embedding} in the form
\begin{equation}
\label{eqn:SL_2 inv embedding bis}
M[J] \otimes_{\cO} \cO_{\triv_G^{\vee}} \isoto
(M\cotimes_{\bE_{\thetabar}}\bP_{\thetabar})^{\SL_2(\Q_p)}.
\end{equation}
Recalling that $\bE_{\thetabar}/J = \cO,$
and that $\bP_{\thetabar}/J\bP_{\thetabar} = \cO_{\triv_G^{\vee}}$,
we may rewrite the source
as $M[J]\otimes_{\bE_{\thetabar}} \bP_{\thetabar},$
and hence rewrite~\eqref{eqn:SL_2 inv embedding bis} as
$$M[J]\otimes_{\bE_{\thetabar}} \bP_{\thetabar} \iso 
(M\cotimes_{\bE_{\thetabar}}\bP_{\thetabar})^{\SL_2(\Q_p)}.
$$
Tracing the maps through, one sees that the composition of this isomorphism with
the inclusion
$$(M\cotimes_{\bE_{\thetabar}}\bP_{\thetabar})^{\SL_2(\Q_p)} \hookrightarrow
M\otimes_{\bE_{\thetabar}} \bP_{\thetabar}$$
coincides with the natural morphism in the statement of the lemma.\qedhere

\end{proof} 

\begin{lemma}\label{technical remark for Steinberg block}
Let~$\fB$ be a block of type~\emph{\ref{item: Steinberg block}}.
If $M$ is an object of~$\fC_{\fB}$, then $M^{\SL_2(\Q_p)}$ is a subobject
of~$M$ in~$\fC_{\fB}$.
Furthermore, we have the equality
\[M^{\SL_2(\bQ_p)} = \bigcup_{i \in \cI} N_i,\]
where~$\{N_i\}_{i\in \cI}$ is the set of $\fC_{\fB}$-subobjects of~$M$ of finite $\cO$-type.
\end{lemma}
\begin{proof}
The first statement is true because the condition of being $\SL_2(\bQ_p)$-invariant is a closed condition.
For the second statement, observe first that 
if~$v \in M^{\SL_2(\bQ_p)}$ then $\langle G \cdot v \rangle$ is a
cyclic module over $\cO[G/\SL_2(\bQ_p)]_\zeta$, which is finitely generated over~$\cO$. %
There remains to prove that if~$N \subset M$ is a $\fC_\fB$-subobject of finite $\cO$-type, then $N \subset M^{\SL_2(\bQ_p)}$.
This is a consequence of Lemma~\ref{trivial action}.
\end{proof}

\begin{lemma}\label{trivial action}Let $R$ be a complete Noetherian local
$\cO$-algebra with finite residue field.
Let~$W$ be a
finite type $R$-module equipped with a 
continuous {\em (}for the $\m_R$-adic topology on $W${\em )} $G$-action,
admitting a central character, say ~$\alpha: \Q_p^\times \to R^\times$.
Then~$\alpha$ is a square and~$W$ decomposes as a direct sum
of eigenspaces $W_{\beta}$, where $\beta$ runs over the various square roots of~$\alpha$
in the group of continuous $R^{\times}$-valued characters of~$\Q_p^{\times}$,
and $G$ acts on $W_{\beta}$ via $\beta \circ \det$.
\end{lemma}
\begin{proof}
Replacing~$R$ by its Artinian quotients, we can assume that~$R$ is an Artinian ring.
Then~$W$ has the discrete topology, and so the $G$-action is smooth.
Since~$W$ has finite type, the upper unipotent subgroup~$U$ has an open subgroup that acts trivially on~$W$, and conjugating by~$\diag(p^{-1}, 1)$ we see that~$U$ itself also acts trivially.
It follows that~$wUw^{-1}$ acts trivially, hence so does~$\SL_2(\bQ_p)$, which is generated by~$U$ and~$wUw^{-1}$ (see~\cite[Proposition~III.1.1]{MR2642409} for a proof).
Hence~$G$ acts on~$W$ via the determinant, and so~$\alpha$ is a square.
Twisting by~$\alpha^{-1/2}$ (for some choice of this square root),
we can assume that~$\alpha$ is trivial.
But then the $G$-action factors through $\bQ_p^\times/(\bQ_p^\times)^2$, and the lemma follows.
\end{proof}

\subsubsection{The sheaves~\texorpdfstring{$\fW_{\thetabar}$}{W theta} and~\texorpdfstring{$\cW_{\thetabar}$}{W theta}}\label{subsubsec:the sheaves fW and cW}
In this subsection we give an alternative description of the pro-coherent sheaf
\[
\cV_{\thetabar} \cotimes_{\tE_{\thetabar}} \tP_{\thetabar} \in \Pro \Coh(\cX_{\thetabar}), 
\]
which will be useful in Section~\ref{subsec:completing L-infinity}.
As usual, we have formed the completed tensor product using the canonical isomorphism $\tE_{\thetabar} \isoto \tR_{\thetabar}^\op$.

Recall %
from~\eqref{eqn:Steinberg case rewrite}, and the surrounding discussion, that 
$\Hom_{\fC_{\thetabar}}(\bP_{\thetabar},\tP_{\thetabar})$
is an $\bE_{\thetabar}^{\op}$-module and an $\tE_{\thetabar}$-module, finitely presented over both.
It is thus a complete right $\bE_{\thetabar}$-module in $\Mod_c(\tE_{\thetabar})$,
and we obtain a right exact, cofiltered limit-preserving functor
\[
\Hom_{\fC_{\thetabar}}(\bP_{\thetabar}, \tP_{\thetabar}) \cotimes_{\bE_{\thetabar}} \text{--} : \Mod_c(\bE_{\thetabar}) \to \Mod_c(\tE_{\thetabar})
\]
whose composition with $\Mod_c(\tE_{\thetabar}) \to \Mod_c(\cO)$ is the usual completed tensor product.
The evaluation map
\begin{equation}
\label{eqn:Steinberg case rewrite bis}
\Hom_{\fC_{\thetabar}}(\bP_{\thetabar},\tP_{\thetabar}) \cotimes_{\bE_{\thetabar}} \bP_{\thetabar}
\iso \tP_{\thetabar}
\end{equation}
is then an $\tE_{\thetabar}$-linear isomorphism.
Formula~\eqref{eqn:Steinberg case rewrite bis} motivates the following definition.
\begin{df}\label{defn:W-thetabar}
If $\thetabar$ is of type~\ref{item: Steinberg pseudorep},
we write 
\[
\cW_{\thetabar} \coloneqq  
\cV_{\thetabar}\cotimes_{\tE_{\thetabar}} 
\Hom_{\fC_{\thetabar}}(\bP_{\thetabar},\tP_{\thetabar}) %
\]
and $$\fW_{\thetabar} \coloneqq  
\fV_{\thetabar}\otimes_{\tE_{\thetabar}} 
\Hom_{\fC_{\thetabar}}(\bP_{\thetabar},\tP_{\thetabar}).$$
\end{df}
Then $\cW_{\thetabar}$ is a complete right $\bE_{\thetabar}$-module in $\Pro \Coh(\cX_{\thetabar})$, $\fW_{\thetabar}$ is a right $\bE_{\thetabar}$-module in
$\Coh(\fX_{\thetabar})$, and the restriction of the $\bE_{\thetabar}$-action to~$R$ coincides with the action of~$R$ through the structure map $\cX_{\thetabar} \to \Spf R$, 
resp.\ $\fX_{\thetabar} \to \Spec R$.
By~\eqref{eqn:Steinberg case rewrite bis}, we have an isomorphism 
\begin{equation}\label{comparing cV and cW}
\cW_{\thetabar} \cotimes_{\bE_{\thetabar}} \bP_{\thetabar}\iso
\cV_{\thetabar}\cotimes_{\tE_{\thetabar}} \tP_{\thetabar}.
\end{equation}
Since $\cV_{\thetabar}$ algebraizes to the universal
rank two bundle $\fV_{\thetabar}$ on~$\fX_{\thetabar}$,
we see by Lemma~\ref{lem:comparing-frak-and-cal-tensors --- abelian case} that $\cW_{\thetabar}$ algebraizes to 
$\fW_{\thetabar}$,
in the sense that $\cW_{\thetabar}$ is recovered as~$\widehat{\imm}_{\thetabar}^* \fW_{\thetabar}$.

 \begin{rem}%
The sheaf~$\fW_{\thetabar}$ will often intervene in our discussion through its quotient $\fW_{\thetabar}/\fW_{\thetabar}[J]$, where $J \subset \bE_{\thetabar}$ is
the ideal defined in Section~\ref{subsubsec:the ideal J}.
This is because of Theorem~\ref{thm:alg-W-REnd-JWNE} below (which %
is essentially due to Johansson--Newton--Wang-Erickson) and Proposition~\ref{prop:completed L-infinity}~(2).

\end{rem}

\subsubsection{Describing \texorpdfstring{$\bE_{\thetabar}$}{E theta}}
We next recall the structure of $\bE_{\thetabar}$.
This is the subject of~\cite[\S 10.5]{MR3150248},  
and we refer to the discussion there for details of the computations, and for detailed
explanations of the notation that we use. (See also~\cite[\S 5.5]{JNWE}.)%

We first recall from Remark~\ref{description of Bernstein centres} that there is an identification of
the centre of $\bE_{\thetabar}$
(equivalently, the centre of~$\fC_{\thetabar}$;
or again equivalently, 
the centre of the category~$\cA_{\thetabar}$)
with the pseudodeformation ring~$R=R^{\ps}_{\thetabar}$,
such that the functors~$V, V^\dagger$ and~$\Vcheck$ become $R$-linear.

Next, to be precise,
we note that the computations of~\cite[\S 10.5]{MR3150248} concern the endomorphism ring of
$P_{\pi_{\alpha}^{\vee}} \oplus P_{\St^{\vee}} \oplus P_{\triv^{\vee}}$, 
whereas we wish to describe the endomorphism ring~$\bE_{\thetabar}$
of $\bP_{\thetabar} \coloneqq  \tP_{\thetabar} \oplus P_{\triv^{\vee}}$, 
where $\tP_{\thetabar}$ is the twist of 
$P_{\thetabar}=P_{\pi_{\alpha}^{\vee}} \oplus P_{\St^{\vee}}$ 
described in Definition~\ref{defn: tld P}.
(We remind the reader that where we write $P_{\pi^{\vee}_{\alpha}}$, etc.,
Pa\v{s}k\={u}nas in~\cite{MR3150248} writes $\tP_{\pi^{\vee}_{\alpha}}$, etc; while the
notation $P_{\thetabar}$ and $\tP_{\thetabar}$ doesn't appear in {\em loc.\ cit.})
There are therefore three endomorphism rings to be considered in this subsection, namely $E_{\thetabar}  \coloneqq  \End_{\fC_{\thetabar}}(P_{\thetabar}), 
\tE_{\thetabar} \coloneqq  \End_{\fC_{\thetabar}}(\tP_{\thetabar})$ and 
$\bE_{\thetabar} \coloneqq  \End_{\fC_{\thetabar}}(\bP_{\thetabar})$.
One of the ingredients in our description of~$\bE_{\thetabar}$  
will be a choice of basis of the free rank-one $E_{\thetabar}$-module~$\Vcheck(P_{\thetabar})$, which provides an identification
$\tR_{\thetabar} \iso E_{\thetabar}^{\op}$ as in Proposition~\ref{VcheckP}~(3),
and an identification $P_{\thetabar} \iso \tP_{\thetabar}$.
We will describe a convenient choice of basis in Remark~\ref{eqn:identifying E and tE}.

In preparation for this, we recall that \cite[(252)]{MR3150248}, and the displayed equation immediately above
it on~\cite[p.~160]{MR3150248}, 
yield identifications
\begin{equation}
\label{eqn:E_thetabar}
E_{\thetabar} \coloneqq  
\End_{\fC_{\thetabar}}(P_{\thetabar}) 
=
\begin{pmatrix} Re_1 & R\varphi_{12} \\R\varphi_{21}^0 + R\varphi_{21}^1&Re_2
\end{pmatrix}
\end{equation}
and
\begin{equation}\label{first presentation of bE}
\End_{\fC_{\thetabar}}(P_{\pi_{\alpha}^{\vee}} \oplus P_{\St^{\vee}} \oplus P_{\triv^{\vee}})
=
\begin{pmatrix} Re_1 & R\varphi_{12} & R\varphi_{13}^0 + R\varphi_{13}^1\\R\varphi_{21}^0 + R\varphi_{21}^1&Re_2&R\varphi_{23}^0+R\varphi_{23}^1\\R\varphi_{31}&R\beta + R\varphi_{32}&Re_3 \end{pmatrix},
\end{equation}
where the various $\varphi$, as well as $\beta$,
denote certain homomorphisms which are constructed in~\emph{loc.\ cit.},
and~$e_i$ denotes the idempotent of the endomorphism ring corresponding to the $i$-th
direct summand of $P_{\pi_{\alpha}^{\vee}} \oplus P_{\St^{\vee}} \oplus P_{\triv^{\vee}}$.
Of course, the isomorphism~\eqref{eqn:E_thetabar} is just the restriction
to the upper left $2\times 2$ block of the isomorphism~\eqref{first presentation of bE}.

\begin{lemma}\label{3x3 matrix order}
There is a unique $R$-algebra isomorphism
\begin{multline}\label{second presentation of bE}
\End_{\fC_{\thetabar}}(P_{\pi_{\alpha}^{\vee}} \oplus P_{\St^{\vee}} \oplus P_{\triv^{\vee}})
\\
\isom \begin{pmatrix} R & Rc & RX_0+RX_1\\Rb_0+Rb_1&R&Rb_0+Rb_1\\R&Rb_0^{-1}a_0 + Rb_0^{-1}X_0&R \end{pmatrix} \subset M_{3 \times 3}(\Frac(S))
\end{multline}
that preserves the displayed generators of~\eqref{first presentation of bE} and~\eqref{second presentation of bE}.
\end{lemma}
\begin{proof}See~\cite[\S 5.5]{JNWE}.
\end{proof}

\begin{remark}
\label{rem:3x3 matrix order}
We can rewrite the isomorphism~\eqref{second presentation of bE} 
more succinctly (but in a manner that suppresses the explicit identification of generators)
as
$$
\End_{\fC_{\thetabar}}(P_{\pi_{\alpha}^{\vee}} \oplus P_{\St^{\vee}} \oplus P_{\triv^{\vee}})
\cong
\begin{pmatrix} R & Rc & I\\I c^{-1}&R&I c^{-1}\\R& I^{\vee} c &R \end{pmatrix}, 
$$
where as usual we write $I = (X_0,X_1) \subset R$ to denote the reducibility ideal,
and we write $I^{\vee} \coloneqq  \{ a \in \Frac{R}: a I \subseteq R\}$ for the fractional ideal-theoretic
inverse of~$I$.
\end{remark}

\begin{remark}%
\label{rem:2x2 matrix order}
Taken together, \eqref{2x2 matrix order} and \eqref{second presentation of bE} yield isomorphisms
\begin{equation}
\label{eqn:collected isos}
\tld R_{\thetabar} \cong \fourmatrix{R}{Rb_0+Rb_1}{Rc}{R} \text{ and } E_{\thetabar} \cong \fourmatrix{R}{Rc}{Rb_0+Rb_1}{R}.
\end{equation}
Analogously to the description of $\End_{\fC_{\thetabar}}(P_{\pi_{\alpha}^{\vee}} \oplus P_{\St^{\vee}} \oplus P_{\triv^{\vee}})$ provided by Remark~\ref{rem:3x3 matrix order},
we can also write
\begin{equation}
\label{eqn:2x2 matrix order}
E_{\thetabar} 
\cong \fourmatrix{R}{Rc}{Ic^{-1}}{R}.
\end{equation}
\end{remark}

\begin{lemma}
\label{lem:P_thetabar normalization}
We can choose a basis of $\Vcheck(P_{\thetabar})$ as an
$E_{\thetabar}$-module
so that the isomorphism $R_{\thetabar}
\iso E_{\thetabar}^{\op}$ 
of Proposition~{\em \ref{VcheckP}~(3)} 
is given, in terms of the descriptions of its source and target
provided by~\eqref{eqn:collected isos}, by matrix transpose.
\end{lemma}
\begin{proof}
As noted in 
Proposition~\ref{VcheckP}~(3), 
choosing a basis of $\Vcheck(P_{\thetabar})$ gives rise to an $R_{\thetabar}^{\ps}$-algebra isomorphism
\begin{equation}
\label{eqn:iso from basis choice}
\tR_{\thetabar} \iso E_{\thetabar}^{\op},
\end{equation}
and changing this choice
of basis composes~\eqref{eqn:iso from basis choice} with an inner automorphism (of either
its source or target; either one amounts to the same thing).

A consideration of~\eqref{eqn:collected isos} shows that 
matrix transposition provides an $R_{\thetabar}^{\ps}$-algebra isomorphism
$\tR_{\thetabar} \iso E_{\thetabar}^{\op}$, 
which is necessarily the composition of~\eqref{eqn:iso from basis choice}
with an automorphism (again, of either its source or its target). 
It thus suffices to show that  any $R_{\thetabar}^{\ps}$-algebra automorphism of $\tR_{\thetabar}$ is inner.
By Lemma~\ref{normalization for CH algebra}, it in turn suffices to note that any $R_{\thetabar}^{\ps}$-algebra automorphism of $\tR_{\thetabar}$ necessarily fixes the isomorphism classes of the two simple 
modules~$1, \omega^{-1}$ of~$\tR_{\thetabar}$,  because $\Ext^1_{\tR_{\thetabar}}(1, \omega^{-1})$ and~$\Ext^1_{\tR_{\thetabar}}(\omega^{-1}, 1)$ are not isomorphic.
\end{proof}

\begin{remark}
\label{eqn:identifying E and tE}
Fixing an $E_{\thetabar}$-basis for $\Vcheck(P_{\thetabar})$ as in
Lemma~\ref{lem:P_thetabar normalization}
yields isomorphisms $P_{\thetabar} \iso \tP_{\thetabar}$ and $E_{\thetabar} \iso \tE_{\thetabar}$,
such that the isomorphism~\eqref{eqn:endo identification}
is given by matrix transpose (when its source and target
are described $E_{\thetabar} \iso \tE_{\thetabar}$ and~\eqref{eqn:collected isos}).
We make such a choice of basis from now on.
Thus we also obtain
an isomorphism
$\bP_{\thetabar} \iso P_{\pi_{\alpha}^{\vee}} \oplus P_{\St^{\vee}} \oplus P_{\triv^{\vee}},$
and so~\eqref{first presentation of bE} and~\eqref{second presentation of bE}
can be interpreted as an isomorphism
\begin{equation}\label{presentation of bE we use}
\bE_{\thetabar}
\isom \begin{pmatrix} R & Rc & RX_0+RX_1\\Rb_0+Rb_1&R&Rb_0+Rb_1\\R&Rb_0^{-1}a_0 + Rb_0^{-1}X_0&R \end{pmatrix} = \begin{pmatrix} R & Rc & I\\I c^{-1}&R&I c^{-1}\\R& I^{\vee} c &R \end{pmatrix},
\end{equation}
and~\eqref{eqn:collected isos} and~\eqref{eqn:2x2 matrix order} can be interpreted as an isomorphism
\begin{equation}\label{presentation of tE we use}
\tE_{\thetabar} \isom \fourmatrix{R}{Rc}{Rb_0+Rb_1}{R} = \fourmatrix{R}{Rc}{Ic^{-1}}{R}.
\end{equation}
Accordingly, from now on, we drop~$E_{\thetabar}$ from the notation, and work exclusively with~$\tE_{\thetabar}$ and~$\bE_{\thetabar}$.
\end{remark}

Recall from Section~\ref{subsubsec:the ideal J} that $J\subset \bE_{\thetabar}$ denotes the two-sided ideal for which
the quotient $\bE_{\thetabar}/J$ corresponds to $\cO_{\triv_G^\vee}$, i.e.\ the module~$\cO$ with its trivial $G$-action,
under the Morita equivalence of Lemma~\ref{lem:recapitulated Morita}.

\begin{lemma}
\label{lem:J description}
If we write $\mathfrak a \coloneqq  (a_0,a_1,X_0,X_1) \subset R,$ and describe $\bE_{\thetabar}$ via~\eqref{presentation of bE we use},
then
we have
$$
J =
\begin{pmatrix} R & Rc & I\\I c^{-1}&R&I c^{-1}\\R& I^{\vee} c &\mathfrak a \end{pmatrix} 
\subset
\begin{pmatrix} R & Rc & I\\I c^{-1}&R&I c^{-1}\\R& I^{\vee} c &R \end{pmatrix}. 
$$
\end{lemma}
\begin{proof}
The surjection
$\bE_{\thetabar} \to \bE_{\thetabar}/J$ is obtained
by applying
$\Hom_{\fC_{\thetabar}}(\bP_{\thetabar}, \text{--})$
to the composite
$\bP_{\thetabar} \to P_{\triv_G^{\vee}} \to \cO_{\triv_G^{\vee}}.$
Applying 
$\Hom_{\fC_{\thetabar}}(\bP_{\thetabar}, \text{--})$ to the first
arrow,
we obtain the projection of
$\begin{pmatrix} R & Rc & I\\I c^{-1}&R&I c^{-1}\\R& I^{\vee} c &R \end{pmatrix}$ 
onto its third column.
The top two entries of this third column are given by
$\Hom_{\fC_{\thetabar}}(P_{\pi_{\alpha}^{\vee}}, P_{\triv_G^{\vee}})$
and
$\Hom_{\fC_{\thetabar}}(P_{\St^{\vee}}, P_{\triv_G^{\vee}})$
respectively,
each of which map to zero when we compose with the surjection
$P_{\triv_G^{\vee}} \to \cO^{\triv_G^\vee}.$
Thus in fact 
$\bE_{\thetabar} \to \bE_{\thetabar}/J$ is the composition
\begin{multline*}
\begin{pmatrix} R & Rc & I\\I c^{-1}&R&I c^{-1}\\R& I^{\vee} c &R \end{pmatrix} 
\to R =
\Hom_{\fC}(P_{\triv_G^{\vee}},P_{\triv_G^{\vee}})
\\
\to 
\Hom_{\fC}(P_{\triv_G^{\vee}},\cO_{\triv_G^{\vee}}) \iso
\Hom_{\fC}(\cO_{\triv_G^{\vee}}, \cO_{\triv_G^{\vee}}) = \cO
\end{multline*}
(the first arrow being projection onto the bottom right entry; hence this factorization only occurs in the category of $R$-modules, and not of $R$-algebras).
Thus the lemma will follow once we show that
the kernel of the given morphism $R \to \cO$ is indeed equal to~$\mathfrak a$.
Now~\cite[Lem.~10.75]{MR3150248}
shows that this kernel is equal
to the image of the natural morphism
$$
\Hom_{\fC_{\thetabar}}(P_{\pi_{\alpha}^{\vee}} \oplus P_{\St^{\vee}}, P_{\triv_G^{\vee}})
\otimes 
\Hom_{\fC_{\thetabar}}(P_{\triv_G^{\vee}},P_{\pi_{\alpha}^{\vee}}\oplus P_{\St^{\vee}})
\to \Hom_{\fC_{\thetabar}}(P_{\triv_G^{\vee}},P_{\triv_G^{\vee}}),$$
which we can recast more concretely (using our explicit description of $\bE_{\thetabar}$
as a matrix order) as 
the image of the morphism
$$(R \oplus I^{\vee} c) \otimes (I \oplus I c^{-1}) \to R$$
given by performing summandwise multiplication, then adding the results.
This image is equal to the ideal $I + I^{\vee} I$, which one computes to be~$\mathfrak a$. 
\end{proof}

\subsubsection{Describing \texorpdfstring{$\fW_{\thetabar}/\fW_{\thetabar}[J]$}{W mod W[J]}}
Recall that $\tR_{\thetabar}$ acts on $\fV_{\thetabar} = S(1) \oplus S(-1)$
through its identification with a matrix order given by~\eqref{2x2 matrix order},
acting via left multiplication on the elements of $\fV_{\thetabar}$ regarded as column
vectors.
Since we have chosen our identification of $P_{\thetabar}$ and $\tP_{\thetabar}$
to satisfy the conclusion of Lemma~\ref{lem:P_thetabar normalization},
the right action of $\tE_{\thetabar}$ on $\fV_{\thetabar}$ is via the
identification of $\tE_{\thetabar}$ with a matrix order given by~\eqref{presentation of tE we use},
acting by right multiplication on elements of $\fV_{\thetabar}$ (thought of as row vectors).

Recall also that in Section~\ref{subsubsec:the sheaves fW and cW} we defined
$$
\fW_{\thetabar}
\coloneqq  \fV_{\thetabar}
\otimes_{\tld E_{\thetabar}}
\Hom_{\fC_{\thetabar}}(\bP_{\thetabar}, \tld P_{\thetabar})
\iso
\fV_{\thetabar}
\otimes_{\tE_{\thetabar}}
\Hom_{\fC_{\thetabar}}(\bP_{\thetabar},  P_{\thetabar})
$$
(the isomorphism being induced by our identification of $P_{\thetabar}$ with $\tP_{\thetabar}$),
which we can think as a coherent sheaf on $\fX_{\thetabar}$ endowed with
an~$\bE^{\op}_{\thetabar}$-action,
or equivalently as a graded
$S \otimes_R \bE^\op_{\thetabar}$-module
(regarded as a graded ring via the grading on $S$ and the trivial grading
on $\bE^{\op}_{\thetabar}$).

Note that 
$\Hom_{\fC_{\thetabar}}(\bP_{\thetabar},  P_{\thetabar})$
coincides, as an $(\tE_{\thetabar},\bE_{\thetabar})$-bimodule,
with the top two rows of~$\bE_{\thetabar}$ when written as a matrix order as
in~\eqref{presentation of bE we use};
here~$\bE_{\thetabar}$ acts via right multiplication, while $\tE_{\thetabar}$
acts via left multiplication, using the description of~\eqref{presentation of tE we use}.
In particular, reading off the first two rows of~$\bE_{\thetabar}$ in this description,
we find that
\begin{equation}
\label{eqn:explicit iso}
\Hom_{\fC_{\thetabar}}(\bP_{\thetabar},  P_{\thetabar})
\iso
\begin{pmatrix} R & Rc & I\\I c^{-1}&R&I c^{-1} \end{pmatrix}. 
\end{equation}
If we write (as usual) $e_2 = \left( \begin{smallmatrix}0 & 0 \\ 0 & 1 \end{smallmatrix}\right)
\in \tE_{\thetabar},$
then we see that 
$$\Hom_{\fC_{\thetabar}}(\bP_{\thetabar},  P_{\thetabar})
\iso \tE_{\thetabar} \oplus \bigl( (\tE_{\thetabar}e_2) \otimes_R I c^{-1} \bigr)$$
as an $\tE_{\thetabar}$-module.
In fact, $\tE_{\thetabar} e_2$ is simply the second column of $\tE_{\thetabar}$,
regarded as a left $\tE_{\thetabar}$-module. 
So this isomorphism
simply expresses the fact that the third column of the right hand side 
of~\eqref{eqn:explicit iso} is obtained from the second by tensoring over $R$ with~$Ic^{-1}$.

From this discussion, we compute that
\begin{multline}
\label{eqn:fW as a graded module}
\fW_{\thetabar} \coloneqq 
\bigl(S(1)\oplus S(-1)\bigr) \otimes_{\tE_{\thetabar}} \bigl(\tE_{\thetabar} \oplus (\tE_{\thetabar} e_2)
\otimes_R Ic^{-1}\bigr)
\\
\iso
S(1) \oplus S(-1) \oplus \bigl(S(-1) \otimes_R Ic^{-1}\bigr)
\end{multline}
as graded $S$-modules,
where the grading on~$S(-1)\otimes_R Ic^{-1}$ comes from the first factor; i.e.\ we take the tensor product grading with~$R$ and~$Ic^{-1}$ concentrated in degree zero.
In fact, we can make $\bE_{\thetabar}$ (as described in~\eqref{presentation of bE we use})
act
by right multiplication on the target of the isomorphism~\eqref{eqn:fW as a
graded module}: to see that this makes sense,
recall that $b_i = X_ic^{-1}$, so that $Ic^{-1}$ naturally embeds in~$S$,
and there is a product morphism of graded $S$-modules
\begin{equation}
\label{eqn:tensor to product}
S(-1) \otimes_R I c^{-1} \to S(1),
\end{equation}
with image equal to $(b_0,b_1)(1)$ 
(where $(b_0,b_1)$ denotes the indicated homogeneous ideal of~$S$ with its natural grading).
Then ~\eqref{eqn:fW as a graded module} is an isomorphism 
of graded $S\otimes_R \bE_{\thetabar}^{\op}$-modules.
Finally, taking the direct sum of the identity morphisms on the first two summands of~\eqref{eqn:fW as a graded module}, and 
of~\eqref{eqn:tensor to product} on the third summand, we obtain a surjection
\begin{equation}
\label{eqn:fW surjection}
\fW_{\thetabar} 
\buildrel \text{\eqref{eqn:fW as a graded module}} \over = 
S(1) \oplus S(-1) \oplus \bigl( S(-1)\otimes_R I \bigr) \to S(1) \oplus S(-1)
\oplus (b_0,b_1) (1).
\end{equation}

\begin{lemma} 
\label{lem:tensor to product}
The kernel of~{\em\eqref{eqn:fW surjection}} is equal to~$\fW_{\thetabar}[J]$.
Consequently, we have an isomorphism 
of graded ~$S\otimes_R \bE_{\thetabar}^{\op}$-modules
\begin{equation}
\label{eqn:fW/fWJ description}
\fW_{\thetabar}/\fW_{\thetabar}[J] \iso S(1)\oplus S(-1) \oplus (b_0,b_1)(1).
\end{equation}
\end{lemma}
\begin{proof}
Since~\eqref{eqn:fW surjection} is given by the identity on its first two summands,
its kernel coincides with the kernel of~\eqref{eqn:tensor to product}.  As $S$ 
is an integral domain, multiplication by $c^{-1}$ induces an isomorphism
$S(-1) \iso c^{-1} S(1).$
We may then factor~\eqref{eqn:tensor to product} as the composite of this isomorphism
with the morphism
\begin{equation}
\label{eqn:tensor to product bis}
S(-1) \otimes_R I \to S(-1)
\end{equation}
induced by multiplication,
and so the kernel of~\eqref{eqn:fW surjection} coincides with the kernel 
of~\eqref{eqn:tensor to product bis}. 

If we consider the result of tensoring the short exact sequence of $R$-modules
$$ 0 \to I \to R \to R/I \to 0$$
with $S$ over~$R$, then we see that the kernel of~\eqref{eqn:fW surjection}
is isomorphic to $\Tor_1^R(S,R/I).$
Now $R/I$ admits a free resolution with initial terms
$$ R^{\oplus 2} \buildrel \left ( \begin{smallmatrix}  a_1 & -X_1 \\ a_0 & X_0 
\end{smallmatrix} \right )
\over \longrightarrow 
R^{\oplus 2} \buildrel \left ( \begin{smallmatrix} X_0 & X_1 \end{smallmatrix} \right ) \over \longrightarrow R$$
(we regard the elements of $R^{\oplus 2}$ as column vectors,
and perform matrix multiplication on the left),
and so $\Tor_1^R(S,R/I)$ can be computed as the cohomology of the three term
complex
$$ S^{\oplus 2} \buildrel \left ( \begin{smallmatrix}  a_1 & -X_1 \\ a_0 & X_0 
\end{smallmatrix} \right )
\over \longrightarrow 
S^{\oplus 2} \buildrel \left ( \begin{smallmatrix} X_0 & X_1 \end{smallmatrix} \right ) \over \longrightarrow S.$$
Now since $X_i = b_i c,$ and since multiplication by~$c$
is injective, this coincides with the cohomology of the complex
$$ S^{\oplus 2} \buildrel \left ( \begin{smallmatrix}  a_1 & -b_1 c \\ a_0 & b_0 c 
\end{smallmatrix} \right )
\over \longrightarrow 
S^{\oplus 2} \buildrel \left ( \begin{smallmatrix} b_0 & b_1 \end{smallmatrix} \right ) \over \longrightarrow S.$$
One easily confirms that
$$ S^{\oplus 2} \buildrel \left ( \begin{smallmatrix}  a_1 & -b_1  \\ a_0 & b_0  
\end{smallmatrix} \right )
\over \longrightarrow 
S^{\oplus 2} \buildrel \left ( \begin{smallmatrix} b_0 & b_1 \end{smallmatrix} \right ) \over \longrightarrow S$$
is exact,
and so
\begin{equation}
\label{eqn:Tor description}
\ker \text{\eqref{eqn:fW surjection}}
\iso
\Tor_1^R(S,R/I) \iso
\bigl(
S\left( \begin{smallmatrix} a_1 \\ a_0 \end{smallmatrix} \right) 
+
S\left( \begin{smallmatrix} -b_1 \\ b_0 \end{smallmatrix} \right) 
\bigr)
/
\bigl(
S\left( \begin{smallmatrix} a_1 \\ a_0 \end{smallmatrix} \right) 
+
Sc \left( \begin{smallmatrix} -b_1 \\ b_0 \end{smallmatrix} \right) 
\bigr),
\end{equation}
the sums being taken in $S^{\oplus 2}$.

With this explicit description of~$\ker \text{\eqref{eqn:fW surjection}}$
in hand, it is easy to prove the lemma.  Firstly, one directly
checks from the description~\eqref{eqn:Tor description}
that $\fa \coloneqq  (a_0,a_1,X_0,X_1) \subset R$ annihilates 
$ \ker \text{\eqref{eqn:fW surjection}}$,
and then, using the description of $J$ given by Lemma~\ref{lem:J description},
one sees that  
$ \ker \text{\eqref{eqn:fW surjection}}$ is an $\bE_{\thetabar}^{\op}$-submodule
of $\fW_{\thetabar}$ which is annihilated by~$J$.
To complete the proof of the lemma, it remains to show that
\[\fW_{\thetabar}/ \ker \text{\eqref{eqn:fW surjection}}\isoto
S(1) \oplus S(-1) \oplus (b_0,b_1) (1)\]
is $J$-torsion free.
Again,
the $\bE_{\thetabar}$-action on this module is given by right multiplication
by the matrix order of~\eqref{presentation of bE we use},
and so one easily checks this; indeed, extending scalars to~$\Frac(S)$,
it amounts to the fact that $M_3\bigl(\Frac(S)\bigr)$ acts
faithfully on~$\Frac(S)^{\oplus 3}$. 
\end{proof}

\subsubsection{Some constructions from~\texorpdfstring{\cite{JNWE}}{JNWE}}
In~\cite[Prop.~3.5.6]{JNWE} there is introduced a certain finitely generated
graded $S$-module $Q^*$. %
Here $(\text{--})^*$ denotes the internal $\Hom$ with~$S$ on the category of graded $S$-modules, 
i.e.\ the graded duality $\underline{\Hom}_{S\text{-}\gr}(\text{--},S)$; in particular,
the module $Q^*$ is
constructed as
 $\underline{\Hom}_{S\text{-}\gr}(Q,S)$ for some other graded $S$-module~$Q$.
The module $Q$ is actually reflexive 
(in fact it is maximal Cohen--Macaulay over the Gorenstein ring~$S$),
and so $Q$ can be recovered as~$(Q^*)^*$.
For our purposes, however, the most convenient description of $Q^*$ is the one coming from the
displayed formula above~\cite[Prop.~3.5.6]{JNWE},
which exhibits $Q^*$ as the image of the morphism
\begin{equation}
\label{eqn:Q* description}
S(-1) \oplus S(-1) \xrightarrow{\fourmatrix{a_1}{a_0}{-b_1}{b_0}} S(-1) \oplus S(1).
\end{equation}

\begin{lemma}
\label{lem:Q* description}
If we regard $Q^*$ as a submodule of~$S(-1) \oplus S(1)$ {\em (}by regarding
it as the image of  the morphism~{\em \eqref{eqn:Q* description})},
then projection onto the second factor induces an isomorphism
$Q^* \iso (b_0,b_1) (1).$
\end{lemma}
\begin{proof}
If we regard elements of the direct sum as column vectors,
then $Q^*$ is the span of the column
vectors~$\left(\begin{smallmatrix} a_1 \\ - b_1 \end{smallmatrix} \right)$
and~$\left(\begin{smallmatrix} a_0 \\ b_0 \end{smallmatrix} \right)$.
Thus, the image of $Q^*$ under the projection is equal to~$(b_0,b_1).$

Now if $s_0$ and $s_1$ are elements of~$S$,
then
$$ (s_1 a_1 + s_0 a_0)b_0 = a_0(s_1 (-b_1) + s_0 b_0).$$ 
Thus, since $S$ is an integral domain,
we see that kernel of the projection is trivial.  This completes the proof of the lemma. 
\end{proof}

We now recall one of the main results of~\cite[\S 5.5]{JNWE}. 

\begin{defn}\label{defn:X*} Write $X^* \coloneqq  S(1)\oplus S(-1)\oplus Q^*$. \end{defn}\begin{prop}
\label{prop:JNWE End computation}
The graded $S$-module $X^*$ is maximal Cohen--Macaulay,
and there are isomorphisms
$$\bE_{\thetabar}^{\op} \iso \End_{S\text{-}\gr}(X^*) \iso \REnd_{S\text{-}\gr}(X^*).$$ 
The action of $\bE_{\thetabar}^{\op}$ is described
by using the identification~{\em \eqref{presentation of bE we use}} of $\bE_{\thetabar}$ with an 
order in $M_3(S)$, and having it act by matrix multiplication on the right
of $X^* \coloneqq  S(1) \oplus S(-1) \oplus Q^*$ {\em (}regarded as a module of row vectors{\em )}.
\end{prop}
\begin{proof}As the notation suggests, $X^{*}$ is the dual of $X\coloneqq S(-1)\oplus S(1)\oplus Q$.
It therefore suffices to show 
that 
$\Ext^i_{S\text{-}\gr}(X, X) = 0$ for~$i \ne 0$, and that
$\bE_{\thetabar} \iso  \End_{S\text{-}\gr}(X)$.
The first statement is proved in~\cite[Prop.\ 3.5.2]{JNWE},
and the second statement is proved immediately before~\cite[Rem.\ 5.5.2]{JNWE}.
\end{proof}
\begin{thm}\label{thm:alg-W-REnd-JWNE}
The $S \otimes_R \bE_{\thetabar}^{\op}$-module 
$\fW_{\thetabar}/\fW_{\thetabar}[J]$ is isomorphic to~$X^*$,
and the natural map\[\bE_{\thetabar}^{\op}\to
  \REnd_{\Coh(\fX_{\thetabar})}(\fW_{\thetabar}/\fW_{\thetabar}[J])\] is an isomorphism.
Consequently,
the functor defined via {\em Lemma~\ref{lem:derived-tensor-product-EW}}~\eqref{item:87}:
\[
\cW_{\thetabar}/\cW_{\thetabar}[J] \otimes_{\bE_{\thetabar}}^L -: D^b_{\fp}(\bE_{\thetabar}) \to \Pro D^b_{\coh}(\cX_{\thetabar}),
\]
 is fully faithful.

\end{thm}
\begin{proof}
Taken together,
Lemmas~\ref{lem:tensor to product} 
and~\ref{lem:Q* description} 
show that each of
$\fW_{\thetabar}/\fW_{\thetabar}[J]$ and $X^*$
is isomorphic to $S(1) \oplus S(-1) \oplus (b_0,b_1)(1) $ as graded ~$S\otimes_R \bE_{\thetabar}^{\op}$-modules.
This proves the first claim.
The second claim follows from the first claim and Proposition~\ref{prop:JNWE End computation}.
  By Theorem~\ref{thm:coherent-completeness-CWE-stacks},
  the completion map 
  \[
  \REnd_{D^b_{\coh}(\fX_{\thetabar})}(\fW_{\thetabar}/\fW_{\thetabar}[J]) \to \REnd_{\Pro D^b_{\coh}(\cX_{\thetabar})}(\cW_{\thetabar}/\cW_{\thetabar}[J]) 
  \]
  is an isomorphism,
  the final claim is a consequence of the second one and Lemma~\ref{lem:fp-tensor-full-faithfulness-criterion}.
\end{proof}

\subsubsection{\texorpdfstring{$\Tor$}{Tor}-dimension of \texorpdfstring{$X^*$}{X*} over \texorpdfstring{$\bE_{\thetabar}$}{E theta}}

It follows from \cite[Prop.~5.5.3]{JNWE} that $X^*\otimes^L_{\bE_{\thetabar}}\text{--}$
has amplitude $[-1,0]$ when evaluated on objects of $\Mod^{\fl}(\bE_{\thetabar}).$
Presumably there is a ``$\Tor$-dimension $1$''
variant of the local criterion for flatness %
which allows us to deduce that $\fW_{\thetabar}/\fW_{\thetabar}[J]\cong X^{*}$ is of $\Tor$-dimension~$1$, but we give a direct proof as follows.
We use the notion of $\Tor$-dimension from Definition~\ref{topologically flat object}.

\begin{lemma}
\label{lem:X* tor-dim'n}
The $\bE_{\thetabar}^{\op}$-module  $\fW_{\thetabar}/\fW_{\thetabar}[J]$ is of $\Tor$-dimension~$1$.
\end{lemma}

\begin{proof}%
  By Theorem~\ref{thm:alg-W-REnd-JWNE}, the $\bE_{\thetabar}^{\op}$-modules   $\fW_{\thetabar}/\fW_{\thetabar}[J]$ and~$X^{*}$ are isomorphic, so it is equivalent to show that~$X^{*}$ has $\Tor$-dimension~$1$.
The explicit computations of Example~\ref{example:functor on irreducible objects in the Steinberg case}
below imply that the $\Tor$-dimension of $X^*$ is at least~$1$.  Thus we need only
show that it is also at most~$1$, which we now do.

To begin,
recall the isomorphism~\eqref{presentation of bE we use}:
$$
\bE_{\thetabar} \iso
\begin{pmatrix} R & Rc & RX_0+RX_1\\Rb_0+Rb_1&R&Rb_0+Rb_1\\R&Rb_0^{-1}a_0 + Rb_0^{-1}X_0&R \end{pmatrix} 
=
\begin{pmatrix} R & Rc & I\\I c^{-1}&R&I c^{-1}\\R& I^{\vee} c &R \end{pmatrix} 
$$
acting by right multiplication on $X^* = S(1) \oplus S(-1) \oplus (b_0,b_1)(1).$
Note also that
$$
S = \cdots \oplus R c^{-n} \oplus Rc^{1-n}
\oplus \cdots \oplus R c \oplus R \oplus I c^{-1} \oplus I^2 c^{-2} \oplus \cdots \oplus
I^n c^{-n} \oplus \cdots,$$
while
$$(b_0,b_1)
= \cdots \oplus I c^{-n} \oplus Ic^{1-n}
\oplus \cdots \oplus I c \oplus I \oplus I c^{-1} \oplus I^2 c^{-2} \oplus \cdots \oplus
I^n c^{-n} \oplus \cdots.$$
Thus (using subscripts to denote graded pieces) we have

\[
(X^*)_{2n-1} =
\begin{cases}
  Rc^{-n} \oplus Rc^{-n+1} \oplus Ic^{-n} 
    & \text{if } n \leq 0, \\[0.5ex]
  I^n c^{-n} \oplus I^{n-1} c^{-n+1} \oplus I^n c^{-n}
    & \text{if } n \geq 1.
\end{cases}
\]

In order to show that~$X^{*}$ has $\Tor$-dimension~$\leq 1$, it suffices to show that it has a projective resolution of length two.
In fact, the composite
\[
\Mod^{\fp}(\bE_{\thetabar}) \xrightarrow{X^* \otimes_{\bE_{\thetabar}}} \Coh(\fX_{\thetabar}) = \operatorname{GrMod}(S) \to \Mod(S)
\]
is the usual tensor product, and so its derived functor can be computed by projective resolutions on either variable; and the forgetful functor
$\operatorname{GrMod}(S) \to \Mod(S)$ is faithful and exact. %
Thus it suffices in turn to show that each $(X^*)_{2n-1}$ has such a resolution.

If $n \leq 1,$ then in fact
$(X^*)_{2n-1}$ is isomorphic (as a right $\bE_{\thetabar}$-module)
to one of the rows of $\bE_{\thetabar}$, and is thus projective over~$\bE_{\thetabar}$.

If $n \geq 2,$ then
$$(X^*)_{2n-1} = I^n c^{-n} \oplus I^{n-1}  c^{-n + 1} \oplus I^n c^{-n}.$$
This has the following two step projective resolution:%
\[0\to  P_1 ^{\oplus n-1}c^{-n}\to  P_2 ^{\oplus n}c^{1-n}\to
  (X^*)_{2n-1}\to 0 \] where $P_1 =(R,I^{\vee}c,R)$, $P_2 =(Ic^{-1},R,Ic^{-1})$,
the first map is given by \[\bigl((-X_1 , X_0 ,0,\dots,0),(0,-X_1 ,X_0
  ,\dots,0),\dots,(0,\dots,0,-X_1 ,X_0 )\bigr),\] and the second map by $(X_0 ^{n-1},X_0 ^{n-2}X_1 ,\dots,X_1^{n-1})$.
\end{proof}

\begin{remark}
\label{rem:lifting S structure}
Note that the projective resolution of $X^*$ as an $\bE_{\thetabar}^{\op}$-module
constructed in the preceding lemma is in fact a graded resolution.
Furthermore,
it is not difficult to describe
explicit lifts of the action of each of $b_0,$ $b_1$, and $c$
to this resolution. 
Indeed, for $n \leq 1$, the ``resolution'' of $X^*_{2n-1}$ is of length one,
i.e.\ $X^*_{2n-1}$ is already its own projective resolution,
and so nothing need be said about the action of $c$ on these graded 
pieces, or the action of the $b_i$ on $X^*_{2n-1}$ for $n \leq 0$.

If $n \geq 1$, then we have the following commutative diagram (where we regard the terms as column vectors):
$$\xymatrix{
0 \ar[r] & P_1^{\oplus n-1}c^{-n} \ar^-{
\left(\begin{smallmatrix} -X_1 & 0 & \cdots & 0\\
X_0 & - X_1 & \cdots & 0 \\
\vdots & \vdots & \cdots & \vdots \\
0 & 0 & \cdots & -X_1 \\
0 & 0 & \cdots & X_0 \end{smallmatrix}\right)}[rrr] 
\ar^-
{\left(\begin{smallmatrix} c^{-1} & 0 & \cdots & 0\\
0 & c^{-1} & \cdots & 0 \\
\vdots & \vdots & \cdots & \vdots \\
0 & 0 & \cdots & c^{-1} \\ 
0 & 0 & \cdots & 0 
\end{smallmatrix}\right)}[dd] 
&&&
P_2^{\oplus n} c^{1-n} 
\ar^-{\left(\begin{smallmatrix} X_0^{n-1} & X_0^{n-2}X_1 & \cdots & X_1^{n-1}\end{smallmatrix}\right)}[rrr]
\ar^-
{\left(\begin{smallmatrix} c^{-1} & 0 & \cdots & 0\\
0 & c^{-1} & \cdots & 0 \\
\vdots & \vdots & \cdots & \vdots\\
0 & 0 & \cdots & c^{-1} \\ 
0 & 0 & \cdots & 0 
\end{smallmatrix}\right)}[dd] 
&&&
X_{2n-1}^* \ar^-{b_0}[dd]
\ar[r]
& 0
\\
&&&&&&&&
\\
0 \ar[r] & P_1^{\oplus n}c^{-n-1} \ar_-{
\left(\begin{smallmatrix} -X_1 & 0 & \cdots & 0\\
X_0 & - X_1 & \cdots & 0 \\
\vdots & \vdots & \cdots & \vdots \\
0 & 0 & \cdots & -X_1 \\
0 & 0 & \cdots & X_0 \end{smallmatrix}\right)}[rrr] 
&&&
P_2^{\oplus n+1} c^{-n} 
\ar_-{\left(\begin{smallmatrix} X_0^{n} & X_0^{n-1}X_1 & \cdots & X_1^{n}\end{smallmatrix}\right)}[rrr]
&&&
X_{2n+1}^* \ar[r] 
& 0
}
$$
which gives a lift of multiplication by $b_0$. 
One easily constructs 
analogous diagrams which give lifts of multiplication by $b_1$, and of~$c$.
These then allow one to compute $\Tor_i^{\bE_{\thetabar}}(X^*,M)$
as a graded $S$-module for any $\bE_{\thetabar}$-module~$M$.
\end{remark}

\begin{example}\label{example:functor on irreducible objects in the Steinberg case}
If $\pi$ is an irreducible object of $\cA_{\thetabar}$,
then there is an $\bE_{\thetabar}$-linear isomorphism
$$\bP_{\thetabar}\otimes_{\cO\llbracket G\rrbracket _{\zeta}} \pi
\iso \bE_{\thetabar}/\mathbf{m}_{\pi},$$
where  
$\mathbf{m}_{\pi}$ is a maximal left ideal
in~$\bE_{\thetabar}$.
Concretely, writing~$\fm$ for the maximal ideal of~$R$, we have
$$\mathbf{m}_{\pi_{\alpha}}
= 
\begin{pmatrix} \fm & Rc & I\\I c^{-1}&R&I c^{-1}\\R& I^{\vee} c &R \end{pmatrix},
$$
$$\mathbf{m}_{\St}
= 
\begin{pmatrix} R & Rc & I\\I c^{-1}&\fm&I c^{-1}\\R& I^{\vee} c &R \end{pmatrix},
$$
and
$$\mathbf{m}_{\triv}
= 
\begin{pmatrix} R & Rc & I\\I c^{-1}&R&I c^{-1}\\R& I^{\vee} c &\fm \end{pmatrix}.
$$
Following the prescription of Remark~\ref{rem:lifting S structure},
one can then compute
$$\Tor_i^{\bE_{\thetabar}}(X^*, \bP_{\thetabar}\otimes_{\cO\llbracket G\rrbracket _{\zeta}}\pi)
= \Tor_i^{\bE_{\thetabar}}(X^*, \bE_{\thetabar}/\mathbf{m}_{\pi_\alpha})
$$
as a graded $S$-module
for each irreducible~$\pi$.  
Since these $\Tor$-values are already computed
in~\cite[Prop.~5.5.3]{JNWE}
by a different method
(namely, by computing resolutions in $\fC_{\thetabar}$,
rather than by resolving~$X^*$ as we have done) 
we omit the full details
of the calculation,
and simply record the result for future reference.
As it turns out, these~$\Tor$-values vanish for all but one value of~$i$, and are then as follows:
\[
\Tor_0^{\bE_{\thetabar}}(X^*, \bP_{\thetabar}\otimes_{\cO\llbracket G\rrbracket _{\zeta}}\pi_\alpha) \cong S(1)/(a_0, a_1, b_0, b_1,
  \varpi),\] 
\[\Tor_0^{\bE_{\thetabar}}(X^*, \bP_{\thetabar}\otimes_{\cO\llbracket G\rrbracket _{\zeta}}\St)  \cong S(-1)/(a_0, a_1, c, \varpi),\] 
\[\Tor_{1}^{\bE_{\thetabar}}(X^*, \bP_{\thetabar}\otimes_{\cO\llbracket G\rrbracket _{\zeta}}\triv_G)\cong S(-3)/(a_0, a_1, c, \varpi).\]
\end{example}

   \subsection{Defining an open substack of~\texorpdfstring{$\cX$}{X}}\label{subsec: defining U}
In this section, we shift our attention away from the various stacks~$\cX_{\thetabar}$,
and back to the stack~$\cX$.  Our main goal is to describe a certain ``good'' open
substack of~$\cX$ --- to be denoted $\cUgood$ --- and to describe its underlying
reduced substack explicitly.

   We begin by noting that since~$\thetabar$ has fixed determinant $\zetabar\omega^{-1}$,
   there are only finitely many~$\thetabar$ of types~\ref{item: ss
     pseudorep}, \ref{item: non p-distinguished pseudorep}
   or~\ref{item: Steinberg pseudorep}, even if we allow ourselves to replace~$\F$ by~$\Fpbar$.
   Note furthermore that~$\thetabar$ of
   type~\ref{item: non p-distinguished pseudorep} (resp.\ of
   type~\ref{item: Steinberg pseudorep}) exist if and only
   if~$\zeta$ is odd (resp.\ even).

\begin{defn}%
\label{defn:closed substack of bad points}
\leavevmode
\begin{enumerate}
\item We let~$\cYbad$ denote the
reduced closed substack of~$\cX$ with ~$|\cYbad|=\piss^{-1}(|\Ybad|)$, where~$\Ybad$ is defined in Definition~\ref{defn:f(t)}. 
Equivalently,
$\cYbad$ is the (finite) union of
the~$\cX_{\thetabar,\red}$ for~$\thetabar$ of types~\ref{item: ss
  pseudorep}, \ref{item: non p-distinguished pseudorep}
or~\ref{item: Steinberg pseudorep}.  
\item We let~$\cUgoodred$ be the open
complement of~$\cYbad$ in~$\cX_{\red}$, and we let~$\cUgood$ be the open
substack of~$\cX$ with underlying reduced substack~$\cUgoodred$.
Note that $|\cUgood| = \pi_{\mathrm{ss}}^{-1}(|\Ugood|)$, where~$\Ugood$ is defined in Definition~\ref{defn:f(t)}.
\end{enumerate}
\end{defn}

\begin{rem}
  \label{rem: what are the points of U}
By the definition of~$\cUgood$, the finite type points of~$\cUgood$
correspond to those
  $\rhobar:G_{\Qp}\to\GL_2(\cbF_p)$ which are reducible, and furthermore have
  corresponding pseudorepresentation  of type~\ref{item: generic
    pseudorep}. 
\end{rem}

\begin{defn}\label{defn:U sigmasigmaco}
  If $\{\sigma,\sigmacomp\}$ is a pair of companion weights, we
  set   
\[\cU(\sigmasigmacomp)_{\red} \coloneqq  
\cUgoodred\cap \bigl(\cX(\sigma)\cup\cX(\sigmacomp)\bigr),\] \[\cU(\sigma)_{\red} \coloneqq  
\cUgoodred\cap \cX(\sigma).\]
Note that~$\cU(\sigmasigmacomp)_{\red}$ is the open
substack of~$\cX_{\red}$ whose underlying topological space is $\piss^{-1}(U(\sigmasigmacomp))$, where~$U(\sigmasigmacomp)$ is defined in Definition~\ref{defn:f(t)}.
\end{defn}

\begin{lemma}
\label{lem:reduced reducible connected components}
The stack $\cUgoodred$ decomposes as a disjoint union of connected components
$$\cUgoodred= \coprod_{\sigmasigmacomp} \cU(\sigmasigmacomp)_{\red},$$
where $\{\sigma,\sigmacomp\}$ runs over all companion pairs of Serre weights.
\end{lemma}
\begin{proof}Since $\cYbad$ contains all the~$\cX_{\red,\thetabar}$
  for~$\thetabar$ of type~\ref{item: Steinberg pseudorep}, it in
  particular contains all irreducible components of~$\cX_{\red}$ of
  the form~$\cX(\sigma_{a,p-1})^{\pm}$, so $\cUgoodred$ is
  the union of the $\cU(\sigmasigmacomp)_{\red}$. %
  Since~$\cYbad$ also
  contains the $\cX_{\red,\thetabar}$ for~$\thetabar$ of
  type~\ref{item: ss pseudorep}, it follows from the description of
  the Serre weights associated to a representation~$\rhobar$ (as
  recalled in Section~\ref{subsec: notation and
  conventions}) that the
  $\cU(\sigmasigmacomp)_{\red}$ are disjoint, and that each $\cU(\sigmasigmacomp)_{\red}$ is
  connected, because the components~$\cX(\sigma)$ and~$\cX(\sigmacomp)$
  intersect along the split locus, i.e.\ at those~$\rhobar$ which are
  direct sums of characters. %
\end{proof}

\begin{defn}We let %
$\cU(\sigmasigmacomp)$
be the open substack of~$\cX$ with underlying reduced substack
$\cU(\sigmasigmacomp)_{\red}$. 
\end{defn}

\begin{cor}\label{cor: U itself is a disjoint union of companion pairs} The stack $\cUgood$ decomposes as a disjoint union of connected components
$$\cUgood= \coprod_{\sigmasigmacomp} \cU(\sigmasigmacomp),$$
where $\{\sigma,\sigmacomp\}$ runs over all companion pairs of Serre weights.
\end{cor}
\begin{proof}
  This is immediate from Lemma~\ref{lem:reduced reducible connected
    components} and the definitions of the stacks~$\cUgood$, $\cU(\sigmasigmacomp)$.
\end{proof}

\begin{defn}
  \label{defn: explicit U sigmasigmacomp notation}Let
  $\{\sigma,\sigmacomp\}$ be a companion pair of Serre weights. As in
  Definition~\ref{defn:f(t)}, we write
\[
f(t) =
\begin{cases}
  t & \text{if } \sigmasigmacomp \text{ is of type~\ref{item:
            generic pseudorep}}  \\
  t(1-t^2) & \text{if } \sigmasigmacomp \text{ is of type~\ref{item:
            non p-distinguished pseudorep} or~\ref{item: Steinberg pseudorep}}.
\end{cases}
\]
We define a group scheme~$\bG$ over~$\Spec \cO$ by
\[
\bG =
\begin{cases}
  \Gm & \text{if } \sigmasigmacomp \text{ is of type~\ref{item:
            generic pseudorep} or ~\ref{item: Steinberg pseudorep}}  \\
 C_2\ltimes \Gm  & \text{if } \sigmasigmacomp \text{ is of type~\ref{item:
            non p-distinguished pseudorep}},
\end{cases}
\]where $C_2$ denotes the cyclic group $\{1,\tau\}$ of order $2$,
and the semi-direct product
$C_2\ltimes \Gm$ is defined by
$\tau\cdot u = u^{-1}$.

\end{defn}

  \begin{prop}  \label{prop: explicit description of underlying
      reduced of U}%

  We have an isomorphism

\begin{equation}\label{eqn:description of underlying reduced of U sigma sigmaco}\cU(\sigmasigmacomp)_{\red} \cong [ 
  \Spec \F [t,f(t)^{-1},x,y]/(xy) /\mathbf{G}],\end{equation}
where 
\begin{itemize}  \renewcommand{\labelitemi}{--} 
\item the action of $\Gm\le\mathbf{G}$ is given by
\[  u\cdot (t,x,y)   =   (t, u^2 x, u^{-2} y),\]
\item and if $\{\sigma,\sigmacomp\}$ is of type ~\emph{\ref{item:
            non p-distinguished pseudorep}}, then the action of
        $C_2$ is given by \[  \tau\cdot (t,x,y)   =
          (t^{-1}, y, x).\]
\end{itemize}

If~$\sigma=\sigma_{a,b}$ then the representation~$\rhobar:G_{\Qp}\to\GL_2(\Fpbar)$ corresponding to
a finite type point of $\cU(\sigmasigmacomp)_{\red}$ is
of the form \begin{equation}\label{eqn:points-of-explicit-stack-U-fixed-det}\rhobar\cong
  \begin{pmatrix}
    \nr_{t^{-1}}\zetabar\omega^{-a} & *_x\\ *_y&\nr_t\omega^{a-1}
  \end{pmatrix}.
\end{equation}More precisely,
by~\eqref{eqn:points-of-explicit-stack-U-fixed-det} we mean the following:
if~$y=0$, then~$\rhobar$ is an extension of $\nr_{t }\omega^{a-1}$ by
$\nr_{t^{-1} }\zetabar\omega^{-a}$, and is split if and only if~$x=0$ \emph{(}note that
since~$t^2\ne 1$ if $\{\sigma,\sigmacomp\}$ is of type~\emph{\ref{item: non p-distinguished
  pseudorep}} or~\emph{\ref{item: Steinberg pseudorep}}, there is a unique non-split
extension up to isomorphism\emph{)}; and similarly if~$x=0$, it is an extension of
$\nr_{t^{-1}}\zetabar\omega^{-a} $ by $\nr_{t}\omega^{a-1}$, and is split if
and only if~$y=0$.

\end{prop}
\begin{proof}
This is a straightforward consequence of the relationship
between~$\cX(\sigmasigmacomp)$ and the special fibre of a tamely potentially
crystalline stack of cuspidal type, and in particular of the explicit
description of these special fibres in ~\cite{hung2023local}. Since
~\cite{hung2023local} does not fix the determinant of the stacks that
it considers, we begin by explaining how to reduce to that context. To
that end, we write $\cX_2$ for the stack of projective \'etale 
  $(\varphi,\Gamma)$-modules of rank $2$ as
  in~\cite{emertongeepicture}, and we write~$\cX_2(\sigmasigmacomp)$ for the
  union of the irreducible components~$\cX_2(\sigma)$ and~$\cX_2
  (\sigmacomp)$ of~$\cX_{2,\red}$.

  We will construct a certain
  open substack~$\cU$ of~$\cX_2(\sigmasigmacomp)$.
  If~$\sigma$ is not of type~\ref{item:
            non p-distinguished pseudorep}, then $\cU$ will admit the description \[\cU\cong [ 
  \Spec \F [t_1 ^{\pm 1},t_2 ^{\pm 1},x,y]/(xy)
  /(\Gm\times\Gm)], \]where the action of $\Gm\times\Gm$ is given by
\[  (u_1 ,u_2 )\cdot (t_1 ,t_2 ,x,y)   =   (t_1 ,t_2 , u_1 u_2 ^{-1}
  x, u_1 ^{-1}u_2  y).\] If $\sigma$ is of type~\ref{item:
            non p-distinguished pseudorep}, then $\cU$ will admit the description  \[\cU\cong [ 
  \Spec \F [t_1 ^{\pm 1},t_2 ^{\pm 1},(t_1 -t_2 )^{-1},x,y]/(xy)
  /C_2\ltimes(\Gm\times\Gm)], \]where~$C_2 $ exchanges the two copies
of~$\Gm$, the action of $\Gm\times\Gm$ is as above, and  the action of
        $C_2$ is via \[  \tau\cdot (t_1 ,t_2 ,x,y)   =
          (t_2 ,t_1 , y, x).\]
Furthermore, in either case, the representation~$\rhobar:G_{\Qp}\to\GL_2(\Fpbar)$ corresponding to
a finite type point of $\cU(\sigmasigmacomp)_{\red}$ is
of the form  
\[\rhobar\cong
  \begin{pmatrix}
    \nr_{t_1 }\omega^{a+b} & *_{x} \\ *_{y}&\nr_{t_2 }\omega^{a-1}
  \end{pmatrix}.\]

Admitting the existence of such a~$\cU$, %
we deduce the proposition as follows: by the
definition of~$\cX$ (i.e.\ Definition~\ref{defn: Xdchi}), there is an
open substack~$\cU'$ of~$\cX(\sigmasigmacomp)$ determined by the pullback diagram   \[ \begin{tikzcd}%
      \cU'\ar{d}\ar{rr}\arrow[dr, phantom, "\ulcorner", very near start]&{}&\cU\ar{d}{t_1 t_2 }\\
      \Spec\F\ar{rr}{(\zetabar\omega^{-2a-b})(\Frob_p)}&{}&\Gm
    \end{tikzcd}\] It follows from the definitions that~$\cU'$ is
  given by the right hand side of~\eqref{eqn:description of underlying
    reduced of U sigma sigmaco}, and to see
  that~$\cU'=\cU(\sigmasigmacomp)_{\red}$, it suffices to note that
  they have the same $\Fpbar$-points. %

  We now turn to the construction of~$\cU$. %
Replacing~$\sigma$ by
  ~$\sigmacomp$ if necessary, we may assume that $0\le b\le
  (p-3)/2$ %
  unless~$\sigma$ is of type~\ref{item: non p-distinguished pseudorep}. After twisting, we can and do furthermore assume
  that~$a=1$. %
  We write~$\Theta(\sigmasigmacomp)$ for 
  the cuspidal representation of~$\GL_2 (\Fp)$ 
  in Definition~\ref{defn: companion Serre weights}~(2).
  By ~\cite[Lem.\ 3.1.1(4)]{CDT}, the
  Jordan--H\"older factors of $\Theta(\chi)\otimes_{\cO}\F$
  are~$\{\sigma,\sigmacomp\}$. 
  It then follows (for example\footnote{This 
is a particular case (possibly a motivating one) of the Breuil--M\'ezard
conjecture~\cite{BreuilMezard} 
(in a weak form in which we pay no attention to the size of the multiplicities, but only
to whether or not they are non-zero), which was proved in general by Kisin~\cite{KisinFM}.
For a reformulation in the framework of moduli stacks of \'etale $(\varphi,\Gamma)$-modules,
see~\cite[\S 8.5]{emertongeepicture}.})  from~\cite[Lem.\ 4.2.4]{CDT} 
 that $\rhobar:G_{\Qp}\to\GL_2(\Fpbar)$ has a potentially
  crystalline lift of type ~$\tau\coloneqq \omega_2
  ^{b+2}\oplus\omega_2 ^{p(b+2)}$ %
  and Hodge--Tate weights $0,1$ if
  and only if $\rhobar$ corresponds to an $\Fpbar$-point of~$\cX_2
  (\sigmasigmacomp)$.

  Consequently, we see that in the notation of~\cite[Defn.\
4.8.8]{emertongeepicture}, %
any open substack of $(\cX_2 ^{\crys,(1,0),\tau})_{\F}$ which is also reduced is
in fact an open substack of $\cX_2 (\sigmasigmacomp)$. We will
now %
exhibit~$\cU$ as such an open substack. Firstly, note that by \cite[Thm.\
4.5]{bellovin2023irregular} it is equivalent to produce an open substack of the
special fibre of the stack denoted $\cZ^{\tau}$ in~\cite{hung2023local}, where
in the notation of that paper, we take~$f=1$, so that~$K=\Qp$, and we
have~$\tau=\tau((12),(b+2,0))$. %

Suppose firstly that we are not in the ~\ref{item: non p-distinguished
  pseudorep} case. Write $G=\GL_2/\F$ with its usual upper-triangular Borel
subgroup~$B=TU$, where~$T$ is the diagonal maximal torus, and~$U$ the subgroup of upper-triangular unipotent matrices. Write ~$LG$ for its loop group in
the variable~$v$, $L^+G$ for the positive loop group, and respectively write
$L^+I$ %
and $L^+_1G$ for the Iwahori subgroup of matrices which
are upper-triangular modulo~$v$ and for the principal congruence subgroup of
matrices congruent to the identity modulo~$v$.

As in~\cite[\S 3.1]{hung2023local}, write~$\mathcal{A}(\eta)\subset
LG$ for the subfunctor with $\cA(\eta)(R)$ given by those $A\in
\GL_2(R((v)))$ satisfying
\begin{enumerate}
\item $\det A\in v\cdot R\llbracket v\rrbracket ^{\times}$.
\item $A\in M_2(R\llbracket v\rrbracket )$, and $A$ is upper-triangular modulo~$v$.
\item $vA^{-1}$ is upper-triangular modulo~$v$.
\end{enumerate}
Write $LG^{\tau}\subset LG$ for the subfunctor \[LG^{\tau}(R)=\mathcal{A}(\eta)(R)\begin{pmatrix}
    & 1\\1&
  \end{pmatrix}
  \begin{pmatrix}
    v^{b+2}&\\&1
  \end{pmatrix}.\]Then %
by \cite[Lem.\
3.2.1, 3.3.7]{hung2023local}  (and in the case $p=5$ and~$b=1$, \cite[Cor.\
3.3.15]{hung2023local} and its proof) %
we can identify the special fibre of~$\cZ^{\tau}$ %
with the quotient stack %
\[ [(L^+_1G\backslash LG^{\tau})/B-\textrm{conj}]\] %

By~\cite[Lem.\ 3.1.2(3)]{hung2023local}, this has an open substack~$\cU$
given by %

\[\left[\left (L^+_1G\backslash L^+I \begin{pmatrix}
    X & 1\\v&Y
  \end{pmatrix}\begin{pmatrix}
    & 1\\1&
  \end{pmatrix}
  \begin{pmatrix}
    v^{b+2}&\\&1
  \end{pmatrix}\right)/B-\textrm{conj}\right]\]with $XY=0$. 

Note that \[\begin{pmatrix}
    X & 1\\v&Y
  \end{pmatrix}\begin{pmatrix}
    & 1\\1&
  \end{pmatrix}
  \begin{pmatrix}
    v^{b+2}&\\&1
  \end{pmatrix}=\begin{pmatrix}
    v^{b+2}&X \\v^{b+2}Y&v
                \end{pmatrix}.\]  Conjugation by an element of~$U$ acts as follows:

              {\footnotesize
              \begin{align*}
                \begin{pmatrix}
                  1 & x\\ 0&1
                \end{pmatrix}
                \begin{pmatrix}
                  v^{b+2}&X \\v^{b+2}Y&v
                \end{pmatrix} \begin{pmatrix}
                                1 & -x\\ 0&1
                              \end{pmatrix}&=
                \begin{pmatrix}
                  1+xYv^{b+1}+x^2Y^2v^{b+1}  & x-v^{b+1}x-v^{b+1}x^2Y \\ xY^2v^{b+1}& 1-xYv^{b+1}
                \end{pmatrix}
                \begin{pmatrix}
                  v^{b+2}&X \\v^{b+2}Y&v
                \end{pmatrix}\\ &\in L^{+}I\cdot \begin{pmatrix}
                  v^{b+2}&X \\v^{b+2}Y&v
                \end{pmatrix}
              \end{align*}\par}%

Since we have $B=T\ltimes U$ and $L^+G=\GL_2\ltimes L_1^+G$,
$L^+I=B\ltimes L_1^+G$, %
it follows that ~$\cU$ is
equal to the quotient stack%

\[\left[\left (\begin{pmatrix}
                  a&0 \\0&d
                \end{pmatrix} \begin{pmatrix}
                  v^{b+2}&X \\v^{b+2}Y&v
                \end{pmatrix}\right)/T-\textrm{conj}\right]\]with $XY=0$. By construction, the
          \'etale $\varphi$-module corresponding to an $\Fpbar$-point of this
          stack is given by taking~$\varphi$ to be the matrix $\big(\begin{smallmatrix}
                  a&0 \\0&d
               \end{smallmatrix}\big) \big(\begin{smallmatrix}
                  v^{b+2}&X \\v^{b+2}Y&v
                \end{smallmatrix}\big)$, from which we can immediately read off~\eqref{eqn:points-of-explicit-stack-U-fixed-det}.

Suppose now that we are in the ~\ref{item: non p-distinguished
  pseudorep} case. %
We follow the analysis
of this case in~\cite[\S 5.5.1]{hung2023local}. Accordingly, we set
$\widetilde{w}=w_0 \eta$, $\widetilde{z}=t_{(1,1)}$, and it follows from the discussion
in \emph{loc.\ cit.} that in this case we have an open substack~$\cU$ given by the quotient by
simultaneous conjugation of the pairs of matrices $(\varphi,N)$ satisfying the equation \[N\cdot\varphi N\varphi^{-1}=0,\] where~$\varphi$ is regular semisimple and $N$ is nilpotent (more precisely, $N$ is in the nilpotent cone of~$\GL_2 $); and the corresponding \'etale $\varphi$-module is given by the matrix $\varphi\left(v 1_2+N\right)$.
 Indeed in the notation of ~\cite[\S
5.5.1]{hung2023local}, we take \[\varphi=
  \begin{pmatrix}
    \alpha&\beta\\\gamma&\delta
  \end{pmatrix}^{-1},\
N=
  \begin{pmatrix}
    -D&B\\ C&D
  \end{pmatrix}.
\]%

Since the regular semisimple matrix~$\varphi$ is conjugate to an
element of~$T$, this stack is isomorphic to the quotient of the pairs
of matrices $(\varphi,N)$ as above where now~$\varphi\in T$ is regular
semisimple, by conjugation by the normalizer~$N(T)=T\rtimes C_2$. This
is the required description. %
\end{proof}

Write  $j':\cUgood\to\cX$ for the open immersion.
\begin{lem}\label{lem:vanishing-higher-derived-functors-j}%
  The functor
$j'_*: D^b_{\coh}(\cUgood)
\to \Ind D^b_{\coh}(\cX)$ 
is $t$-exact.
\end{lem}
\begin{proof}%
By Corollary~\ref{cor: U itself is a disjoint union of companion pairs}, we can replace~$\cUgood$ by $\cU\coloneq \cU(\sigmasigmacomp)$ for a companion pair~$\{\sigma,\sigmacomp\}$. Write $\cX \iso \colim \cX_n$ as a colimit of algebraic stacks with the 
transition maps being thickenings, and define $j'_n : \cU_n \hookrightarrow \cX_n$
to be the base-change
of the open immersion $j':\cU \hookrightarrow \cX.$
Then (as noted in Remark~\ref{notational remark about j star}),
the functor $j'_*:D^b_{\coh}(\cU) \to\Ind  D^b_{\coh}(\cX)$
is equivalent to the colimit of functors $(j'_n)_*: D^b_{\coh}(\cU_n) \to
\Ind D^b_{\coh}(\cX_n),$
these functors (or, more precisely, their Ind-extensions) being instances of 
the general construction~\eqref{eqn:Ind derived pushforward}
associated to an open immersion of Noetherian
algebraic stacks admitting affine diagonals.

To prove the lemma
it thus suffices to show that each $(j_n')_*$ is $t$-exact.
By Remark~\ref{rem:cohomologically affine case} (and the discussion
that precedes it), it suffices to show that each of the open immersions 
$j_n'$ is cohomologically affine; this we now do. (Recall that a morphism $f:\cX\to\cY$ of algebraic stacks is
\emph{cohomologically affine} (\cite[Defn.\ 3.1]{MR3237451}) if it is quasi-compact, and if furthermore
$f_{*}:\QCoh(\cX)\to\QCoh(\cY)$ is exact.)

By Proposition~\ref{prop: explicit description of underlying reduced of U}, $\cU_{\red}$ is  the quotient of an affine scheme by a linearly reductive group (note that since we are assuming
  that~$p>2$, the group scheme~$\bG/\Spec \F$ is linearly reductive  by~\cite[Ex.\ 12.4(2)]{MR3237451};  recall that an \emph{fppf} group scheme~$G/S$ 
is \emph{linearly reductive} if
the corresponding morphism $BG\to S$ is cohomologically
affine).
Thus $\cU_{\red}$ is cohomologically
affine over~$\Spec \F$ (see~\cite[Ex.\ 12.9]{MR3237451}), and consequently each~$\cU_n$ is 
cohomologically affine over~$\Spec \cO$, by~\cite[Prop.\ 3.9]{MR3237451}. Furthermore each $\cX_n$ has affine (and hence cohomologically affine) diagonal.
We now factor $j'_n$ as
$$\cU_n \buildrel \Gamma_{j_n'} \over \longrightarrow\cU_n \times_{\Spec \cO} \cX_n \to \cX_n,$$
the first arrow being the graph of $j_n'$, and the second arrow being the second projection.
Since
$\Gamma_{j_n'}$ is a base-change of $\Delta_{\cX_n}$, it is affine (and in particular
cohomologically affine).  
Since $\cU_n \to \Spec \cO$ is cohomologically affine, so is its pullback
over~$\cX_n$.
Thus $j'_n$ is the composite of cohomologically affine morphisms, and hence is itself
cohomologically affine, as required. %
\end{proof}

Although we don't need it in what follows, we end this section by %
showing that the morphism~$\piss:|\cX| \to |X|$ defined in~\eqref{eqn:defn-of-cts-map-f} is a quotient map.

\begin{lemma}\label{preparation for quotient map}
Let~$\{\sigma,\sigmacomp\}$ be a companion pair of $\zeta$-compatible Serre weights.
Then the restriction
\[
\piss: |\cU(\sigmasigmacomp)| = \piss^{-1}(|U(\sigmasigmacomp)|) \to |U(\sigmasigmacomp)|
\]
is a quotient map.
\end{lemma}
\begin{proof}
Since (as noted in the proof of Lemma~\eqref{lem:vanishing-higher-derived-functors-j}) the group~$\bG$ in Definition~\ref{defn: explicit U sigmasigmacomp notation} is linearly reductive over~$\bF$, it follows from~\eqref{eqn:description of underlying reduced of U sigma sigmaco} and~\cite[Theorem~13.2]{MR3237451} that the algebraic stack~$\cU(\sigmasigmacomp)_{\red}$
has a good moduli space $p: \cU(\sigmasigmacomp)_{\red} \to U$, given by the spectrum~$U$ of the ring of $\bG$-invariants of $\bF[t, f(t)^{-1}, x, y]/(xy)$.
This ring of invariants is equal to $\bF[t,f(t)^{-1}],$ resp.\ $\bF[t+t^{-1}, 1/(t^2+t^{-2}-2)]$,
if $\{\sigma,\sigmacomp\}$ is of type~\ref{item: generic pseudorep}, resp.\ of
type~\ref{item: non p-distinguished pseudorep}, 
and the induced morphism on $\Fbar_p$-valued points
sends an $\Fbar_p$-valued point of $\cU(\sigmasigmacomp)$ of 
the form~\eqref{eqn:points-of-explicit-stack-U-fixed-det}
to the point $\Fbar_p$-valued point $t$, resp.\ $t+t^{-1}$, of~$U$.

Now we may (uniquely) identify $\Spec U$ with $U(\sigmasigmacomp)$ in such a way that 
an $\Fbar_p$-valued point~$t$, resp.\ $t+t^{-1}$, of $U$ is identified with the
$\Fbar_p$-valued point of~$U(\sigmasigmacomp)$ corresponding to the pseudorepresentation
$\nr_{t^{-1}}\zetabar\omega^{-a} + \nr_t\omega^{a-1}$.
With this identification, we see (using the above description
of its effect on $\Fbar_p$-valued points) that the morphism
\[
|\cU(\sigmasigmacomp)|_{\ft} = |\cU(\sigmasigmacomp)_{\red}|_{\ft} \to |U|= |U(\sigmasigmacomp)|_{\ft}
\]
induced by the good moduli space morphism~$p$ 
is precisely (the restriction to~$\cU(\sigmasigmacomp)$ of)
the morphism~\eqref{eqn:cX to X on ft points}
which maps each finite type point to its associated pseudorepresentation.
Soberizing, we find that the morphism induced by $p$ on
underlying topological spaces in fact coincides with~$\piss$.
By~ \cite[Theorem~4.16(v)]{MR3237451}, $p$ induces a quotient map, so we are done.
\end{proof}

\begin{lemma}
The morphism $\piss:|\cX| \to |X|$ is a quotient map of topological spaces.
\end{lemma}
\begin{proof}
Since $|\cX|$ is Noetherian, and since $|\cX|_{\ft}\to |X|_{\ft}$ is surjective,
the morphism $|\cX| \to |X|$ is also surjective,
by Lemma~\ref{lem:Noetherian soberization surjectivity} below. 
Hence it suffices to show that $\piss : \piss^{-1}(|X(\sigmasigmacomp)|) \to
|X(\sigmasigmacomp)|$ is a quotient map for each companion pair~$\{\sigma,\sigmacomp\}$.
Let $S \subset |X(\sigmasigmacomp)|$, and assume that $\piss^{-1}(S)$ is closed in $\piss^{-1}(|X(\sigmasigmacomp)|)$, or equivalently in~$|\cX|$.
We need to prove that~$S$ is closed.
By Lemma~\ref{preparation for quotient map}, we see that $S \cap |U(\sigmasigmacomp)|$ is closed in~$|U(\sigmasigmacomp)|$.
Hence either $S \cap |U(\sigmasigmacomp)|$ is a finite set of closed points, or $|U(\sigmasigmacomp)| \subset S$.
Since $X(\sigmasigmacomp) \setminus U(\sigmasigmacomp)$ is a finite set of closed points, the first possibility implies that~$S$ is a finite set of closed points,
which concludes the proof in this case.
Similarly, the second possibility implies that $X(\sigmasigmacomp) \setminus S$ is a finite set of closed points, and so~$S$ is open in~$X(\sigmasigmacomp)$.
Hence $\piss^{-1}(S)$ is a closed and open subset of $\piss^{-1}(|X(\sigmasigmacomp)|)$, which is a connected topological space.
This implies that $\piss^{-1}(S) = \piss^{-1}(|X(\sigmasigmacomp)|)$ and so $S = X(\sigmasigmacomp)$, since~$\piss$ is surjective. 
This concludes the proof that~$S$ is closed.
\end{proof}

\begin{lem}
\label{lem:Noetherian soberization surjectivity}
If $f:T \to S$ is a continuous surjective map of topological spaces, and $T$ is Noetherian,
then the induced map $\sob(T)\to \sob(S)$ is surjective. 
\end{lem}
\begin{proof}
Let $A$ be an irreducible closed subspace of $S$.  Let $\{B_i\}$ 
denote the (finite, because $T$ is Noetherian) set of irreducible
components of~$f^{-1}(A)$.  Then, since $f$ is surjective,
we have $A = \bigcup \overline{f(B_i)}$, and so (since $A$ is irreducible),
we have $A = \overline{f(B_{i_0})}$ for one of the components $B_{i_0}$.
Thus the induced map $\sob(T) \to \sob(S)$ is indeed surjective.
\end{proof}

\begin{remark}
In fact the morphism $f$ has a lift to a morphism $\cX_{\red} \to X$,
which in turn admits a thickening to a morphism $\cX \to \tX$,
where $\tX$ is a certain formal scheme over $\Spf \cO$ with underlying reduced
equal to~$X$.
The formal scheme~$\tX$, and the morphism from $\cX$ to it,
are characterized as being initial in the category of morphisms from
$\cX$ to formal algebraic spaces over~$\Spf \cO$.  In other words,
$\tX$ is the formal moduli space associated to~$\cX$.
Since we do not need this result in this paper, we don't give a proof here
(although we intend to give a construction elsewhere). 
\end{remark}%

\section{\texorpdfstring{$D^\natural\boxtimes\Pone$}{D natural box P1}}%
\label{sec: Colmez}
In this section we
study analogues for general coefficients %
of some of the constructions of~\cite{MR2642409}.
Throughout the discussion, all \'etale $(\varphi,\Gamma)$-modules are
automatically understood to be projective, with the exception that from
Section~\ref{subsec:formal phi Gamma modules} onward
we exploit the relationship between 
\'etale $(\varphi,\Gamma)$-modules
and Galois representations,
and so allow ourselves at times to consider \'etale $(\varphi,\Gamma)$-modules
(and especially so-called {\em formal} \'etale $(\varphi,\Gamma)$-modules) that are 
not necessarily projective.

Our ring of coefficients is typically taken to be a Noetherian~$\cO/\varpi^a$-algebra,
for some $a \geq 1$; the main exception to this occurs in our use 
of formal \'etale $(\varphi,\Gamma)$-modules, %
where the coefficients are 
instead taken to be a complete Noetherian local $\cO$-algebra with finite residue field.
At certain points we also restrict our coefficients to be a finite type~$\cO/\varpi^a$-algebra;
we do this sometimes for topological reasons (the various topological modules
under consideration will then be Polish topological groups, which lets us apply
the open mapping theorem), and sometimes because we wish to complete at maximal ideals
and pass to the formal context. 

\subsection{Tate modules and topologies}\label{subsec: Tate module}The
coefficient rings~$\A_A$ and~$\A_A^+$ introduced in Section~\ref{subsec: the coefficient
  rings}, and many of the various
modules over %
these rings that we study below, are examples of \emph{Tate $A$-modules} in
the sense of~\cite[\S3]{MR2181808}. We now recall this notion and some
of its basic properties; see also~\cite[\S5]{EGstacktheoreticimages}
and~\cite[App.\ D]{emertongeepicture} for some related results and
discussion.

\begin{df}
\label{def:Tate modules}
Let $A$ be a commutative ring. An \emph{elementary Tate
  $A$-module} is a topological $A$-module which is isomorphic to
$P\oplus Q^*$, where $P,Q$ are discrete projective $A$-modules,
and $Q^* \coloneqq  \Hom_A(Q,A)$ equipped with its natural projective limit
topology (where we write~$Q^*$ as the projective limit of the
$(Q')^*$, where~$Q'$ is a finite direct summand of~$Q$, and give each
$(Q')^*$ the discrete topology). A \emph{Tate $A$-module} is a direct summand of an elementary Tate
$A$-module.

A morphism of Tate modules is a continuous morphism of the underlying $A$-modules.
\end{df}

\begin{df}
\label{def:lattice}
A submodule $L$ of a Tate $A$-module $M$ is a \emph{lattice} if it is
open, and if furthermore for every open submodule $U\subseteq L$, the
quotient $L/U$ is a finitely generated $A$-module. %
A subset of a Tate
$A$-module is \emph{bounded} if it is contained in some lattice.
\end{df}

\begin{rem}
  \label{rem: open and bounded equals lattice}If $A$ is Noetherian,
  then a submodule~$L$ of a Tate $A$-module~$M$ is a lattice if and
  only if it is open and bounded. Indeed if~$L$ is open and bounded,
  then by definition there is a lattice~$L'\supseteq L$, and for any
  open submodule $U\subseteq L$, the quotient~$L'/U$ is a finitely
  generated $A$-module; so the submodule $L/U\subseteq L'/U$ is also a
  finitely generated $A$-module.
\end{rem}

If~$A$ is an $\cO/\varpi^a$-algebra for some $a\ge 1$, then $\A_A^+$
and~$\A_A$ are elementary Tate modules, and~$\A^+_A$ is a lattice
in~$\A_A$; furthermore every finitely generated
projective~$\A_A$-module has a natural topology, making it a Tate
$A$-module. (See for example~\cite[Ex.\ 5.1.13]{EGstacktheoreticimages}.)

 In fact, by~\cite[Thm.\ 3.10]{MR2181808}, there is a natural bijection between finitely generated projective $\A_A$-modules, and pairs $(M,T)$
consisting of a Tate $A$-module $M$ and a topologically nilpotent
automorphism $T:M\to M$, by giving $M$ the $\A_A$-module structure
determined by
$Tm\coloneqq T(m)$. (Here $T$ is topologically nilpotent if and only if for
each pair of lattices $L,L'\subseteq M$, we have $T^nL\subseteq L'$ for
all sufficiently large~$n$.)

\begin{defn}\label{defn: A+ lattice}
  If~$M$ is a finitely generated projective $\A_A$-module, then we
  refer to an~$\A^+_A$-submodule $\gM\subset M$ which is also a
  lattice as an \emph{$\A^+_A$-lattice}  (or \emph{$\A^+$-lattice}
  if~$A$ is clear from the context).
\end{defn}

\begin{rem}
  \label{rem: all lattices are Aplus lattices after rescaling}Let $A$
  be an $\cO/\varpi^a$-algebra. If ~$M$
  is a finitely generated projective $\A_A$-module, and $\gM\subset M$
  is a lattice, then~$\gM$ need not be an $\A_A^+$-lattice. However,
  since multiplication by~$T$ is topologically nilpotent, for all
  sufficiently large~$n$ we have $T^n\gM\subset\gM$, so  that $\gM$ is
  an $A\llbracket T^n\rrbracket $-module. Since $M$ is a finite projective
  $A((T^n))$-module, and the rings $A\llbracket T^n\rrbracket $, $A((T^n))$ are
  respectively isomorphic to~$\A_A^+$ and~$\A_A$, this means that we
  can often reduce questions about lattices to the case of
  $\A^+$-lattices. 
\end{rem}

We recall the following lemma due to Drinfeld~\cite[Lem.~5.1.20]{EGstacktheoreticimages}.

\begin{lem}
\label{lem:coprojective}
If $M$ is a finitely generated projective~$\A_A$-module,
and if $\gM$ is a lattice in~$M$, then $M/\gM$ is $A$-flat if and only
if it is $A$-projective.   If $\gM$ is furthermore an~$\A_A^+$-lattice,
then these conditions are in turn equivalent to $\gM$ being~$\A_A^+$-projective.
\end{lem}

\subsubsection{The weak topology on lattices}\label{subsubsec:weak topology}

If~$A$ is a complete Noetherian local ring, 
and $L$ is a lattice in an $A$-Tate module,
then in addition to its Tate-module topology,
it has a {\em weak topology}, as we now explain.
Lemma~\ref{lattices in Tate modules are complete} below shows that we have an isomorphism (of $A$-modules, i.e.\ disregarding
topologies) $L \iso \varprojlim_k L/\m_A^k L$.
If we endow $L/\m_A^k L$ with its quotient topology induced by the Tate module
topology on $L$, 
then the right-hand side inherits an inverse limit topology,
which we then can transport to $L$ via this isomorphism;
this is (by definition) the weak topology on~$L$.

\begin{lem}\label{lattices in Tate modules are complete}
If~$A$ is a complete Noetherian local ring, and $L \subset M$ is a lattice in a Tate $A$-module, then~$L$ is $\fm_A$-adically complete.
The $A$-module~$L$ with the weak topology is an object of $\Mod_c(A)$, hence it is profinite if $A$ has finite residue field.
\end{lem}
\begin{proof}
Since $M$ is a direct summand of an elementary Tate module, and the direct sum of lattices is a lattice, it suffices to prove the lemma under the assumption
that $M = P \oplus Q^*$ is elementary.

We begin by proving the lemma when~$L = Q^*$. %
Since~$A$ is a local ring, $Q$ is free, and so
$Q^*$ is isomorphic to a product $\prod_{i \in I} A$ for some index set~$I$.
Then $Q^*$ is $\fm_A$-adically complete because $\otimes_A A/\fm_A^n$ commutes with arbitrary products, since $A/\fm_A^n$ is finitely presented.
For the second claim of the lemma, note that the Tate module topology %
on~$\prod_{i \in I} A$ is the product topology with respect to the discrete topology on~$A$, hence the quotient topology on
$\prod_{i \in I}A/\fm_A^n$ is the product topology with respect to the discrete topology on $A/\fm_A^n$.
It follows that the weak topology on~$\prod_{i \in I}A$ %
is the product topology with respect to the $\fm_A$-adic topology on~$A$, 
which makes $\prod_{i \in I} A$ an object of~$\Mod_c(A)$.

Now, if $L \subset Q^*$ is any lattice, then by definition there exists an $A$-finite direct summand $Q' \subset Q$ such that $(Q')^\perp \subset L$.
Writing $Q''$ for the complement of~$Q'$, so that $Q = Q' \oplus Q''$,
we see that $Q^* = Q'^* \oplus Q''^*$ and $(Q')^\perp = Q''^*$.
Hence there exists an $A$-submodule $V \subset Q'^*$ such that $L = V \oplus Q''^*$.
Since~$Q'^*$ is finite, we see that~$V$ is also finite.
By the results of the previous paragraph, this immediately implies that~$L$ is $\fm_A$-adically complete, and that the weak topology makes~$L$ an object of $\Mod_c(A)$.

Finally, if $L \subset M$ is any lattice, then it is contained in $P' \oplus Q^*$ for some finite $A$-summand $P' \subset P$.
Since $M = (P/P') \oplus (P' \oplus Q^*)$, and $(P' \oplus Q^*) \cong (P'^* \oplus Q)^*$, we can conclude the proof by an application
of the results in the previous paragraph (applied with $Q^*$ replaced by~$(P'^* \oplus Q)^*$).
\end{proof}

\begin{example}
In the case when $L$ is an $\A_A^+$-lattice in a finitely generated projective $\A_A$-module,
the weak topology on $L$ coincides with its $(\m_A,T)$-adic topology,
whereas the Tate module topology coincides with its $T$-adic topology.
\end{example}

\subsubsection{Completions and completed tensor products}
\label{subsubsec:completions}
Suppose that $A$ is a ring,
and that $M$ is a topological $A$-module endowed with an $A$-linear topology,
i.e.\ $M$ is an $A$-module with a topological abelian group structure,   
such that $0 \in M$ admits a neighbourhood basis consisting of $A$-submodules of~$M$.
Note that the product morphism $A \times M \to M$ is then jointly continuous, when
$A$ is equipped with its discrete topology.
We may then consider the completion of~$M$ (either as a topological group,
or as a topological
$A$-module; the result is the same).  In symbols, if we let $U$ run over
the partially ordered set of open $A$-submodules of~$M$,
then
$$\widehat{M} \coloneqq  \varprojlim_U M/U.$$
Suppose now that $A \to B$ is a ring homomorphism.
If $M$ is a topological $A$-module endowed with an $A$-linear topology,
then we equip $M\otimes_A B$ with the $B$-linear topology defined by the images 
of the various $U\otimes_A B$, as $U$ runs over the open $A$-submodules of~$M$.
We then define $M\cotimes_A B$ to be the completion of $M\otimes_A B$ with respect
to this topology. 
Thus, if we let $\overline{U\otimes_A B}$ denote the image of $U\otimes_A B$
in~$M\otimes_A B$,
then 
\begin{equation}\label{eqn: recognizing completed tensor}
M\cotimes_A B \coloneqq  \varprojlim_U (M\otimes_A B)/(\overline{U\otimes_A B})
= \varprojlim_U (M/U)\otimes_A B
\end{equation}
(each of the terms in each of the inverse limits being endowed with its discrete topology).
The symbol~$\cotimes_A$ was also used in Section~\ref{subsubsec:theta Morita} to denote the completed tensor product on categories of compact modules,
but this should not lead to ambiguities, since in the current context, the ring~$A$ is not a profinite topological ring.

We will most often apply this notion in the case when~$A\to B$ is a morphism of Noetherian
$\cO/\varpi^a$-algebras, in situations related to the above discussion of finite projective
$\A_A$-modules and their lattices.
If $M_A$ is a finite projective $\A_A$-module,
then its Tate module topology is $A$-linear,
and $M_A$ is complete, as is any lattice in~$M_A$.
The tensor product $M_B\coloneqq  M_A\otimes_{\A_A}\A_B$ is a finite
projective~$\A_B$-module, and hence is again complete.
There is a canonical continuous morphism $M_A \otimes_A B \to M_A \otimes_{\A_A} \A_B$.
Since the target is complete, 
it induces a morphism from the completion of the source to the target, which is in fact
an isomorphism
\begin{equation}\label{recognizing completion of module}
M_A\cotimes_A B \iso M_A \otimes_{\bA_A} \bA_B.
\end{equation}

If~$\gM_A\subseteq M_A$ is a lattice, then we may also consider the completed
tensor product~$\gM_A\cotimes_AB$.
If~$\gM_A$ is an $\A_A^+$-lattice, then this completion naturally
identifies with $\gM_A\otimes_{\A^+_A}\A^+_B$,
in the sense that the natural map
\begin{equation}\label{recognizing completion of lattice}
\gM_A \cotimes_A B \to \gM_A \otimes_{\A^+_A} \A^+_B
\end{equation} 
is a topological isomorphism. %
Indeed, the natural $A$-linear topology on~$\gM_A$ is the $T\bA_A^+$-adic topology, hence the $B$-linear topology on 
\[
\gM_A \otimes_A B = \gM_A \otimes_{\bA_A^+} (\bA_A^+ \otimes_A B)
\]
is the
$T(\bA_A^+ \otimes_A B)$-adic topology.
Since $\A^+_A \otimes_A B$ is a Noetherian ring with $T$-adic completion
$\A_B^+$, and~$\gM_A$ is a finite $\A_A^+$-module, we see that
$\gM_A\otimes_{\A^+_A}\A^+_B$ is the $T$-adic completion of
$\gM_A\otimes_AB$.
The inclusion $\gM_A \subseteq M_A$ induces a continuous morphism
$\gM_A\cotimes_A B\to M_B$, which is open (by Lemma~\ref{lem: collection of facts about lattices}~\eqref{item: lattice 4} below), but which need not be an inclusion in general.

The following lemma collects together various basic statements about
lattices in finite projective $\A_A$-modules, and their behaviour 
under completed tensor product.

\begin{lem}
  \label{lem: collection of facts about lattices}
Let~$A\to B$ be a
  morphism of Noetherian $\cO/\varpi^a$-algebras. %
  Let $M_A$ be a finite
  projective $\A_A$-module, and as above write $M_B = M_A\cotimes_AB = M_A\otimes_{\A_A}\A_B$.
  \begin{enumerate}
  \item\label{item: lattice 1} An $\A^+_A$-submodule $\gM_A\subset M_A$ is an
    $\A^+_A$-lattice if and only if it is finitely generated and
    the $\A_A$-span of~$\gM_A$ is~$M_A$.
      \item\label{item: lattice 2} $M_A$ contains an $\A^+_A$-lattice $\gM_A$.
  \item\label{item: lattice 3} If $L_1$ is a  lattice in~$M_A$,
    and~$L_2$ is an $A$-submodule of~$M_A$, then~$L_2$ is a lattice if
    and only if it is commensurate to~$L_1$, in the sense that for~$n$ sufficiently
    large we have $T^nL_1\subseteq L_2\subseteq
    T^{-n}L_1$. 
  \item \label{item: lattice 4}
If $\gM_A\subseteq M_A$ is a lattice then the evident sequence \begin{equation}
\label{eqn:desired s.e.s.}
  \gM_A\cotimes_A B \to M_B \to (M_A/\gM_A)\otimes_A B \to 0,
\end{equation}
  is exact, 
and
the image of the
    morphism $\gM_{A}\cotimes_AB\to M_{A}\cotimes_AB=M_B$ is a lattice in~$M_B$,
while its kernel is finitely generated over~$B$.
In addition, this kernel is discrete {\em(}with respect to the topology induced on it
by~$\gM_A\cotimes_A B$), and the morphism $\gM_A \cotimes_A B \to M_B$ is open. 

If furthermore either $A \to B$ is flat or $M_A/\gM_A$ is flat as an $A$-module
then \eqref{eqn:desired s.e.s.} is exact on the left as well, i.e.\ we have a short
exact sequence
  $$0 \to \gM_A\cotimes_A B \to M_B \to (M_A/\gM_A)\otimes_A B \to 0.$$
In particular,  $\gM\cotimes_AB$ is then a lattice in $M_B$.
  \item\label{item: lattice 5} If $A\into B$ is injective, and~$\gM$ is a lattice in~$M_B$,
    then $\gM\cap M_A$ is a lattice in~$M_A$. %
  \end{enumerate}

\end{lem}
\begin{proof}
  Part~\eqref{item: lattice 1} is~\cite[Lem.\
  D.9]{emertongeepicture}. Part~\eqref{item: lattice 2} follows from part~\eqref{item: lattice 1},
  because we can choose a surjection $f:(\A_A)^r\to M_A$ for some
  sufficiently large~$r$, and set $\gM_A=f((\A_A^+)^r)$. The ``only
  if'' direction of~\eqref{item: lattice 3} is immediate from the
  topological nilpotence of~$T$. For the converse, since~$L_2$
  contains $T^nL_1$, it is open, while if $U\subseteq L_2$ is open,
  then $L_2/U$ is a submodule of the finitely generated $A$-module
  $(T^{-n}L_1)/U$, and is therefore finitely generated, so~$L_2$ is a
  lattice by definition.  %

  We turn to proving part~\eqref{item: lattice 4}.  
By Remark~\ref{rem: all lattices are Aplus lattices after rescaling}
we can, without loss of generality, assume that~$\gM_A$ is an $\A^+$-lattice. 
As explained above, in this case we have that
\begin{equation*}
\gM_A\cotimes_AB=\gM_A\otimes_{\A^+_A}\A^+_B
\end{equation*}
and that
\begin{equation*}
M_A\cotimes_AB = M_A\otimes_{\A_A}\A_B
= M_A\otimes_{\A_A^+[1/T]}\A^+_B[1/T] = M_A\otimes_{\A_A^+}\A^+_B.
\end{equation*}
Now $M_A/\gM_A = \bigcup_{n\geq 0} (T^{-n}\gM_A)/\gM_A,$
so that 
\begin{multline*}
(M_A/\gM_A) \otimes_{\A_A^+}\A_B^+
= \bigcup_{n\geq 0} \bigl((T^{-n}\gM_A)/\gM_A\bigr) \otimes_{\A_A^+}\A_B^+
\\
= \bigcup_{n\geq 0} \bigl((T^{-n}\gM_A)/\gM_A \bigr)\otimes_A B
= (M_A/\gM_A) \otimes_AB.
\end{multline*}
Thus tensoring the short exact sequence
$$0 \to \gM_A \to M_A \to M_A/\gM_A \to 0$$
with $\A_B^+$ over $\A_A^+$, and reinterpreting the various terms according
to the preceding isomorphisms, we obtain the exact sequence~\eqref{eqn:desired
  s.e.s.}, establishing the first claim of~\eqref{item: lattice 4}.
Furthermore, since the image of $\gM_A\cotimes_AB = \gM_A\cotimes_{\A_A^+}\A_B^+$
in $M_B$ is finitely generated over $\A_B^+$ and generates $M_B$ over $\A_B$,
it is a lattice in~$M_B$, by~\eqref{item: lattice 1}. 

If $A \to B$ is flat, then so is $\A_A^+ \to \A_B^+$, and thus
the preceding right exact sequence becomes short exact.
On the other hand, if $M_A/\gM_A$ is $A$-flat,
then Lemma~\ref{lem:coprojective} shows that $\gM_A$ is projective
over~$\A_A^+$. 
Thus $\gM_A\otimes_{\A_A^+}\A_B^+$ is projective over~$\A_B^+$,
and in particular is $T$-torsion free, and so again we
find that it embeds into its localization
$(\gM_A\otimes_{\A_A^+}\A_B^+)[T^{-1}] = M_B.$

To complete the proof of~\eqref{item: lattice 4},
we must prove our claims about the kernel of the left-most arrow
in~\eqref{eqn:desired s.e.s.},
and also prove that $\gM\otimes_A B \to M_B $ is an open mapping.
Since $M_A$ is a finite projective $\A_A$-module by assumption,
we may find a second finite projective $\A_A$-module $M'_A$ so that
$M_A \oplus M'_A$ is free of finite rank over~$\A_A$.  If we choose
an $\A^+_{A}$-lattice $\gM'_A$ in $M'_A$, then $\gM_A\oplus \gM'_A$ is a lattice
in $M_A \oplus M'_A$.
The formation of the various completed tensor products, kernels, and images 
in play is compatible with taking direct sums, and so replacing $M_A$
by $M_A \oplus M_A'$,
and $\gM_A$ by $\gM_A\oplus \gM_A'$,
we may assume without loss of generality that~$M_A$ 
is free.
In this case, we may choose a free $\A_A^+$-lattice
$\gM_A\subseteq \gM''_A\subset M_A$.
Furthermore, part~(3) shows that there exists $n_0 \geq 0$ 
so that $T^n \gM''_A \subseteq \gM_A$ for
$n \geq n_0$, and the sublattices $T^n \gM''_A$ ($n \geq n_0$)
then form a basis of neighbourhoods of $0$
in~$\gM_A$. 

Since $M_A$ is flat over $\A_A^+$, we see that the kernel of the morphism
$$\gM_A\cotimes_A B = \gM_A \otimes_{\A_A^+}\A_B^+ \to M_B$$
is identified with $\Tor_1^{\A_A^+}(M_A/\gM_A, \A_B^+),$
and we begin by explaining why this is a finite~$B$-module.
From what we have already proved, 
we know that
$$\gM''_A\otimes_{\A_A^+}\A_B^+ \to M_B \to (M_A/\gM_A'')\otimes_{\A_A^+}\A_B^+ \to 0$$
is exact on the left,
so that $\Tor_1^{\A_A^+}(M_A/\gM''_A, \A_B^+) = 0.$
A consideration of the short exact sequence
$$0 \to \gM''_A/\gM_A \to M_A/\gM_A \to M_A/\gM''_A \to 0$$
then yields a surjection
$$\Tor_1^{\A_A^+}(\gM''_A/\gM_A , \A_B^+) \to \Tor_1^{\A_A^+}(M_A/\gM_A,\A_B^+).$$
Now the source of this surjection is finite over $\A_B^+$ (since $\gM''_A/\gM_A$
is finite over~$\A_A^+$) and killed by some power of~$T$ (since $\gM''_A/\gM_A$~is).
Thus it is finite over~$B$, and thus so is $\Tor_1^{\A_A^+}(M_A/\gM_A,\A_B^+)$,
which is what we wanted to prove.

Next, consider the neighbourhood basis $T^n \gM''_A \subseteq \gM_A$ of $0$ in~$\gM_A$.
By definition, the images
$$\overline{(T^n \gM''_A)\cotimes_A B}
\coloneqq  \im\bigl((T^n \gM''_A)\cotimes_A B \to \gM''_A\cotimes_A B\bigr)
$$
form a neighbourhood basis of the topology on~$\gM''_A\cotimes_A B.$
Again, by what we have already proved, for each $n$ the morphism
$(T^n\gM''_A )\cotimes_A B \to M_B$ is injective, and so 
$\overline{(T^n \gM''_A)\cotimes_A B)}$
has trivial intersection with  $\ker(\gM_A\cotimes_A B \to M_B)$.
This shows that this kernel is discrete in $\gM_A\cotimes_A B$.
Also, the image of $(T^n\gM''_A)\cotimes_A B$
is a lattice in $M_B$, and so $\gM_A\cotimes_A B$ has a neighbourhood basis of zero
whose members have open image in~$M_B$.  Thus $\gM_A\cotimes_A B \to M_B$ is open,
as claimed. This completes the proof of~(4).

Finally for~\eqref{item: lattice 5}, by
  parts~\eqref{item: lattice 2} and~\eqref{item: lattice 3} there is
  an $\A_B^+$-lattice $\gM'\subseteq M_B$ and an integer $n\ge 0$ such
  that $T^n\gM'\subseteq\gM\subseteq T^{-n}\gM'$. By~\cite[Lem.\
  D.11]{emertongeepicture}, $\gM'\cap M_A$ is a lattice in~$M_A$, and
  $(T^i\gM')\cap M_A=T^i(\gM'\cap M_A)$ for all integers~$i$. In
  particular we have
  \[T^n(\gM'\cap M_A)\subseteq\gM\cap M_A\subseteq T^{-n}(\gM'\cap
    M_A),\]so~$\gM\cap M_A$ is a lattice in~$M_A$ by~\eqref{item:
    lattice 3}.
\end{proof}

We also note the following result, showing that certain tensor products
are automatically complete.

\begin{lem}
\label{lem:tensor with finite algebra is automatically complete}
If~$A\to B$ is a
  finite morphism of Noetherian $\cO/\varpi^a$-algebras
  and $M_A$ is a finite
  projective $\A_A$-module, 
then the natural morphism $M_A \otimes_A B \to M_A \cotimes_A B$
is a topological isomorphism.  Furthermore, if $\gM_A$ is a lattice in~$M_A$, 
then the natural morphism $\gM_A \otimes_A B \to \gM_A \cotimes_A B$
is a topological isomorphism.
\end{lem}
\begin{proof}
We begin by proving the second statement.
By Remark~\ref{rem: all lattices are Aplus lattices after rescaling}, without loss of generality $\fM_A$ is an $\bA_A^+$-lattice.
By~\eqref{recognizing completion of lattice}, the completed tensor product $\fM_A \cotimes_A B$ identifies with $\fM_A \otimes_{\bA_A^+} \bA_B^+$.
Since~$\bA_A^+ = A\llbracket T\rrbracket $ is Noetherian and $\fM_A$ is a lattice in~$M_A$, $\fM_A$ is finitely presented over~$\bA_A^+$.
The map
\[
\fM_A \otimes_A B \to \fM_A \otimes_{\bA_A^+} \bA_B^+
\]
is part of a natural transformation of right-exact functors on finitely presented $\bA_A^+$-modules.
Hence, to prove that it is an isomorphism, it suffices to prove that it is an isomorphism when~$\fM_A = \bA_A^+$, i.e.\ that $\A^+_A\otimes_AB=\A^+_B$.
Since~$B[T]$ is finitely presented as an $A[T]$-module, this follows because $\bA_A^+ \otimes_A B = A\llbracket T\rrbracket  \otimes_{A[T]} B[T]$ is the $T$-adic completion of~$B[T]$.

We now prove the first statement.
By~\eqref{recognizing completion of module} we have $M_A \cotimes_A B = M_A \otimes_{\bA_A} \bA_B$, so we need to prove that the natural map
\begin{equation}\label{to prove bijective for automatically complete}
M_A \otimes_A B \to M_A \otimes_{\bA_A} \bA_B
\end{equation} 
is a topological isomorphism.
Since it is continuous and open (by Lemma~\ref{lem: collection of facts about lattices}~\ref{item: lattice 4}), it suffices to prove that it is bijective.
Since~$M_A$ is a direct summand of a finite free $\bA_A$-module, it suffice to prove that~\eqref{to prove bijective for automatically complete} is bijective when~$M_A = \bA_A$.
This holds because $\bA_A = \bA_A^+[1/T]$ and  $ \bA_B^+=\bA^+_A \otimes_A B $, as proved in the previous paragraph.
\end{proof}

We close this section by establishing some additional properties related
to the topology of finite projective $\A_A$-modules.

\begin{lemma}
\label{lem:f.g. discrete}
Let $A$ be a Noetherian ring,
let $M$ be a topological $A$-module endowed with
an $A$-linear topology, and suppose furthermore that
there is a sequence $(\mu_n)_{n\geq 0}$ of continuous
morphisms $\mu_n: M \to A$ for which the induced morphism
$M \to \prod_{n=0}^{\infty} A$ is injective.
Then the induced topology on any finitely generated 
$A$-submodule of $M$ is discrete.
\end{lemma}
\begin{proof}
Consider a morphism $f:A^m \to M,$
and write $\nu_n \coloneqq  \mu_n \circ f.$
Let $e_i$ ($i = 1,\ldots,m$) denote the standard basis vectors of~$A^m$,
and write
$$S_r \coloneqq  \{ \bigl(\nu_0(e_1), \ldots, \nu_0(e_m)\bigr), \ldots,
\bigl(\nu_r(e_1), \ldots, \nu_r(e_m)\bigr) \} \subseteq A^m.$$
Then, if $K_r \coloneqq  \ker(\nu_0\times \cdots \times \nu_r)$,
we find that $K_r = S_r^{\perp}$ (the ``orthogonal complement''
of $S_r$ with respect to the usual pairing $A^m \times A^m \to A$).

The $S_r$ form an increasing sequence of subsets of $A^m$, %
and hence the $S_r^{\perp}$
form a decreasing sequence of orthogonal complements.  The Artinian property of 
orthogonal complements\footnote{If $S$ is a subset of $A^m$, then the formation
of $S^{\perp}$ is order reversing, and $S \subseteq S^{\perp\perp}$.  These properties then
imply that $S^{\perp} = S^{\perp\perp\perp}$, and that $S \mapsto S^{\perp\perp}$ is
order preserving. 
Thus $S_r^{\perp\perp}$ is an increasing sequence of submodules of
the Noetherian module~$A^m$, and so stabilizes.
Thus the decreasing sequence $S_r^{\perp} = S_r^{\perp\perp\perp}$ also stabilizes.}
implies that this decreasing sequence stabilizes.  
If we recall that $M \to \prod_{n=0}^{\infty} A$ is injective by assumption, we then see that
$K \coloneqq  \ker(f) = \bigcap_{r \geq 0} K_r = K_{r_0}$ for some sufficiently large~$r_0$.
Hence, if
$N \coloneqq  \im(f)$ then 
\[
N \cap \ker(\prod_{n=0}^\infty A \to \prod_{n =0}^{r_0} A) = 0.
\]
Hence $N$ inherits the discrete topology, as claimed.
\end{proof}

\begin{example}%
\label{ex:f.g. discrete}\leavevmode
\begin{enumerate}
\item
If $A$ is a Noetherian $\cO/\varpi^a$-algebra,
then Lemma~\ref{lem:f.g. discrete} applies to finite projective $\A_A$-modules,
and hence also to lattices in such modules.
\item
If $(M_n)_{n \geq 0}$ is an inverse system of topological modules satisfying the hypotheses
of Lemma~\ref{lem:f.g. discrete}, 
then the same is true of the inverse limit $M \coloneqq  \varprojlim_n M_n.$  
\end{enumerate}
\end{example}

\begin{lemma}
\label{lem:lattice completeness}
Let $A$ be a Noetherian $\cO/\varpi^a$-algebra, $I$-adically complete
with respect to some ideal $I\subseteq A$.  If $M_A$ is a finite projective $\A_A$-module,
then any lattice in $M_A$ is closed in the $I$-adic topology of~$M_A$.
\end{lemma}
\begin{proof}
Let $\gM$ be any lattice in~$M_A$. 
Applying Remark~\ref{rem: all lattices are Aplus lattices after rescaling},
we see that it is no loss of generality to assume that $\gM$ is an~$\A_A^+$-lattice.
Choose another finite projective $\A_A$-module~$M_A'$ such that $M_A \oplus M_A'$
is free, and choose an $\A_A^+$-lattice $\gM'$ in $M_A'$.  Then it suffices
to prove the result for the $\A_A^+$-lattice $\gM\oplus \gM'$ in $M_A\oplus M_A'$,
and so it is no loss of generality to further assume that $M_A$ is free.
Then we may choose a free $\A_A^+$-lattice $L$ in~$M_A$, and by part~(3) 
of Lemma~\ref{lem: collection of facts about lattices}, we may assume (scaling $L$
by a power of $T$ if necessary) that $L$ contains~$\gM$. 

Now the $I$-adic topology on $\A_A \coloneqq  A((T))$
evidently induces the $I$-adic topology on $\A_A^+ \coloneqq  A\llbracket T\rrbracket $,
and so the $I$-adic topology on $M_A$ induces the $I$-adic topology on~$L$.
Thus it suffices to prove that the $I$-adic topology on $L$ induces the $I$-adic topology
on~$\gM$.  This follows directly from Artin--Rees (applied to the ideal $I\A_A^+$
in the Noetherian ring~$\A_A^+$, using the fact that $\bA_A^+$-lattices are finitely generated over~$\bA_A^+$, by Lemma~\ref{lem: collection of facts about lattices}~\eqref{item: lattice 1}).
\end{proof}

\begin{remark}
\label{rem:fake A((T))-modules}
Note that the statements of
parts (4) and (5) of Lemma~\ref{lem: collection of facts about lattices},
and the statements of Lemmas~\ref{lem:f.g. discrete} and~\ref{lem:lattice completeness},
depend only on the structure of $M_A$ as a Tate $A$-module.
Thus if $M_A$ is any Tate $A$-module which is topologically isomorphic
(as a topological $A$-module) to a finite projective $\A_A$-module,
these statements apply to~$M_A$.
\end{remark}

\subsection{\texorpdfstring{$D^\natural$ and $D^\sharp$}{D natural and D sharp}}\label{subsec: D natural D
  sharp}
Unless stated otherwise, in this section we assume that~$A$ is a Noetherian
$\cO/\varpi^a$-algebra for some~$a\ge 1$. %
We prove some analogues with coefficients of some results
from~\cite{MR2642409,ColmezMirabolique}. Our overall strategy, and many of our arguments, is
based on those of Colmez, but the presence of the coefficient ring~$A$ means that we sometimes
have to find new ways to argue.

\begin{rem}
  \label{rem: not using Gamma}Throughout this section we do not make
  any use of the action of~$\Gamma$, and our results are equally valid
  for \'etale $\varphi$-modules over~$\A_A$. We will occasionally
  exploit this in our proofs.
\end{rem}

We begin by defining the~$\psi$-operator on $(\varphi,\Gamma)$-modules.

\begin{prop}%
  \label{prop: existence of psi on phi Gamma modules}%
  Let $D$ be an
  \'etale $(\varphi,\Gamma)$-module with $A$-coefficients. Then there
  is a unique $A$-linear morphism~$\psi:D\to D$ 
  such that
  \[\psi(\varphi(a)m)=a\psi(m), \] \[\psi(a\varphi(m))=\psi(a)m\]
  for any $a\in \A_{A}$, $m\in D$.
  The morphism $\psi$ is furthermore continuous, open, and surjective, and commutes with~$\Gamma$. 
\end{prop}
\begin{proof}
Since $D = \A_A\varphi(D),$ we see that the required properties of $\psi$ uniquely
determine it if it exists.  In fact, the same observation, expressed in the more precise
form of the isomorphism $\Phi_D: \A_A\otimes_{\A_A,\varphi} D \iso D,$
allows us to construct~$\psi$, namely as the
  composite
  \[D\stackrel{\Phi_D^{-1}}{\to} \A_{A}\otimes_{\A_{A},\varphi}D
    \stackrel{\psi\otimes 1}{\to}\A_{A}\otimes_{\A_{A}}D=D.  \]
  The relations
  $\psi(\varphi(a)m)=a\psi(m)$, $\psi(a\varphi(m))=\psi(a)m$ 
  then follow immediately. %
That ~$\psi$ is continuous and open follows
  from the continuity of $\Phi_D^{-1}$ and the continuity and openness\footnote{The
  openness of $\psi$  on $\A_A$
  is clear from the fact that $\psi(\varphi(T^n)\A_A^+) = T^n\psi(\A_A^+) = T^n \A_A^+$
  for any $n\geq 0$.} of~ $\psi$ on~$\A_{A}$, 
  the surjectivity of~$\psi$ follows similarly (and
  is also immediate from the relation $\psi(\varphi(m))=m$),
  and the fact that $\psi$ commutes with the $\Gamma$-action again follows by the same argument
  (or by the uniqueness).
 \end{proof}

\begin{lem}\label{lem: basic psi inclusions}%
  Let~$\gM$ be an $\A_A^+$-lattice in an \'etale $(\varphi,\Gamma)$-module~$D$
  with $A$-coefficients. Then:
  \begin{enumerate}
  \item $\psi(\gM)$ is an $\A^+_A$-module.
  \item If $\varphi(\gM)\subseteq \gM$ then $\gM\subseteq \psi(\gM)$.
  \item If $\gM\subseteq \A^+_A\cdot\varphi(\gM)$, then
    $\psi(\gM)\subseteq \gM$.
  \item If $\psi(\gM)\subseteq \gM$, then $\psi(T^{-1}\gM)\subseteq
    T^{-1}\gM$, and for any $x\in D$ there exists $N\ge 0$ such that
    for all~$n\ge N$, we have $\psi^n(x)\in T^{-1}\gM$.
  \end{enumerate}
\end{lem}
\begin{proof}
In the case $A=\cO/\varpi^a$ this is~\cite[Lem.\
II.4.1]{ColmezMirabolique}, and the same proof works more generally,
as we now briefly sketch. The first two parts follow from the
identities $a\psi(x)=\psi(\varphi(a)x)$
and~$\psi(a\varphi(x))=\psi(a)x$ for $a\in\A_A$, $x\in D$. For the
third part, note that
$\psi(\sum_ia_i\varphi(x_i))=\sum_i\psi(a_i)x_i$.

Finally for the fourth part, note firstly that
\[\psi(T^{-1}\gM)\subseteq \psi(\varphi(T)^{-1}\gM)\subseteq
  T^{-1}\psi(\gM)\subseteq T^{-1}\gM.\]Furthermore, for each~$k\ge 1$
we have
\[\psi(\varphi^k(T)^{-1}\gM)\subseteq
\varphi^{k-1}(T)^{-1}\psi(\gM)\subseteq \varphi^{k-1}(T)^{-1}\gM.\] For
any $x\in D$ we have $x\in \varphi^{k}(T)^{-1}\gM$ for all~$k$
sufficiently large, and the result follows.
\end{proof}

\begin{lem}
  \label{lem: existence of psi stable lattices and bounded gaps}Let~$D$ be an \'etale
  $(\varphi,\Gamma)$-module with ~$A$-coefficients. %
  Then there exists an
  $\A^+$-lattice~$\gM\subset D$
  with $\psi(\gM)=\gM$. Furthermore if~$\gM'\subseteq D$ is another
  $\A^+$-lattice with~$\psi(\gM')=\gM'$ then $T\gM\subseteq\gM'\subseteq T^{-1}\gM$.
\end{lem}
\begin{proof}
  By, for example, \cite[Lem.\ 5.2.7]{EGstacktheoreticimages}, we can
  choose $\A^+$-lattices $\gN_0,\gN_1\subset D$ with
  $\varphi(\gN_0)\subseteq\gN_0\subseteq\gN_1\subseteq\A^+_A\varphi(\gN_1)$. For
  each~$n\ge 0$, set~$\gM_n=\psi^n(\gN_0)$; then by parts (1)-(3) of
  Lemma~\ref{lem: basic psi inclusions}, $(\gM_n)_{n \geq 0}$ is an increasing
  sequence of $\A^+$-lattices, each of which contains~$\gN_0$ and is contained
  in~$\gN_1$. By Lemma~\ref{lem: collection of facts about lattices}~\eqref{item: lattice 3}
  $\gM\coloneqq \cup_{n\ge 0}\gM_n$ is an $\A^+$-lattice in~$D$, and by construction
  we have $\psi(\gM)=\gM$.

  We now show that $\gM'\subseteq T^{-1}\gM$, by arguing as in the
  proof of Lemma~\ref{lem: basic psi inclusions}~(4). Indeed since any
  two $\A^+$-lattices are commensurate (by Lemma~\ref{lem: collection of facts about lattices}~\eqref{item: lattice 3}), there
  exists some integer $n\ge 0$ such that
  $\gM'\subseteq \varphi^n(T)^{-1}\gM$, so that
  $\gM'=\psi^n(\gM')\subseteq T^{-1}\gM$. %

  Finally, the inclusion $T\gM\subseteq\gM'$ is equivalent to
  $\gM\subseteq T^{-1}\gM'$, and thus follows from the previous
  paragraph upon reversing the roles of~$\gM,\gM'$.
\end{proof}

\begin{prop}\label{prop: D sharp natural exist} %
  Let~$D$ be an \'etale $(\varphi, \Gamma)$-module with $A$-coefficients.
  \begin{enumerate}
  \item\label{item: D sharp exists} $D$ contains a maximal $\A^+$-lattice~$D^\sharp\subset D$ with the
    property that %
    $\psi(D^\sharp)=D^\sharp$.
Furthermore, if $\gN$ is any bounded $\A^+$-submodule 
for which $\psi(\gN) = \gN$, then $\gN \subseteq D^\sharp$.
  \item\label{item: D natural exists} $D$ contains a minimal $\psi$-stable $\A^+$-lattice $D^\natural$;
    furthermore $\psi(D^\natural)=D^\natural$.
  \item\label{item: bounding D natural sharp} We have $D^\natural\subseteq D^{\sharp}\subseteq
    T^{-1}D^\natural$.
  \end{enumerate}
\end{prop}
\begin{proof}By Lemma~\ref{lem: existence of psi stable lattices and
    bounded gaps}, there exists an $\A^+$-lattice $\gM\subset D$ with
  $\psi(\gM)=\gM$. Let~$D^\sharp$ denote the sum of all the
  $\A^+$-lattices~$\gM'\subset D$ with $\psi(\gM')=\gM'$; then $\gM\subseteq
  D^\sharp\subseteq T^{-1}\gM$ by Lemma~\ref{lem: existence of psi
    stable lattices and bounded gaps}, so $D^\sharp$ is an $\A^+$-lattice
  by Lemma~\ref{lem: collection of facts about lattices}~\eqref{item: lattice 3}. By construction
  $\psi(D^\sharp)=D^\sharp$, and (again by construction) $D^\sharp$ is
  maximal with this property.
  If $\gN$ is any bounded $\A^+$-submodule for which $\psi(\gN) = \gN $,
  then $D^\sharp + \gN$ is an $\bA^+$-lattice 
  satisfying $\psi(D^\sharp + \gN) = D^\sharp + \gN$.  Thus $D^\sharp + \gN \subseteq D^\sharp$
  by the property we've just proved, and so $\gN \subseteq D^\sharp$.
  This concludes the proof of part~\eqref{item: D sharp exists}.

  We claim that if $\gM\subseteq D^\sharp$ is an $\A^+$-lattice with
  $\psi(\gM)\subseteq\gM$, then in fact $\psi(\gM)=\gM$. Granting
  this, parts~\eqref{item: D natural exists} and~\eqref{item: bounding D natural sharp} follow easily. Indeed, we
  let $D^\natural$ be the intersection of all $\psi$-stable $\A^+$-lattices
  in~$D$; then $D^\natural\subseteq D^\sharp$ (because $D^\sharp$ is a
  $\psi$-stable $\A^+$-lattice), so that (by the claim) $D^\natural$ is equal to the
  intersection of all the $\A^+$-lattices $\gM\subseteq D^\sharp$ with
  $\psi(\gM)=\gM$. It follows from Lemma~\ref{lem: existence of psi
    stable lattices and bounded gaps} that for any such $\A^+$-lattice~$\gM$
  we have $T D^\sharp\subseteq \gM$, so we have $TD^{\sharp}\subseteq
  D^{\natural}\subseteq D^\sharp$ (which establishes part~\eqref{item:
    bounding D natural sharp}), so $D^\natural$ is an $\A^+$-lattice by
  another application of Lemma~\ref{lem: collection of facts about
    lattices}~\eqref{item: lattice 3}. Part~\eqref{item: D natural
    exists} follows from another application of the claim, since~$D^\natural$ is $\psi$-stable (by construction).
  
  It remains to prove the claim. To this end, suppose that
  $\gM\subseteq D^\sharp$ and $\psi(\gM)\subseteq\gM$. Then we argue
  as in the proof of~\cite[Prop.\ II.5.11]{ColmezMirabolique}:
  since~$\psi(D^\sharp)=D^\sharp$, the composite
  \[D^\sharp/\gM\stackrel{\psi}{\to}D^\sharp/\psi(\gM)\to
    D^\sharp/\gM\] is a surjective endomorphism of the finitely
  generated $A$-module~$D^\sharp/\gM$, and is thus an
  isomorphism. %
  In particular $D^\sharp/\psi(\gM)\to D^\sharp/\gM$ is injective, so
  $\psi(\gM)=\gM$, as required.
\end{proof}

\begin{lem}
  \label{lem: psi pushes into T inverse D natural}Let $D$ be an \'etale $(\varphi, \Gamma)$-module with $A$-coefficients, and let $\gM\subset D$ be any
  $\A^+$-lattice. Then there exists an~$N\ge 0$ such that for all~$n\ge N$,
  we have $\psi^n(\gM)\subseteq T^{-1}D^\natural$. 
\end{lem}
\begin{proof}Since any
  two $\A^+$-lattices are commensurate (by Lemma~\ref{lem: collection of facts about lattices}~\eqref{item: lattice 3}), there
  exists some integer $N\ge 0$ such that
  $\gM\subset \varphi^N(T)^{-1}D^{\natural}$. Then for any $n\ge N$
  we have $\psi^n(\gM)\subseteq \psi^{n-N}(T^{-1}D^\natural)\subseteq
  T^{-1}D^\natural$ (recalling that $T^{-1}D^\natural$
  is~$\psi$-stable, by Lemma~\ref{lem: basic psi inclusions}). 
\end{proof}

We close this subsection by describing the behaviour of $D^{\sharp}$ and $D^{\natural}$
under morphisms of $(\varphi,\Gamma)$-modules, and under change of coefficients.

\begin{lem}
  \label{lem: D natural induced maps}If $f:D_1\to D_2$ is a
  morphism of \'etale $(\varphi,\Gamma)$-modules with $A$-coefficients, then~$f$ restricts to
  morphisms $f^\sharp:D_1^\sharp\to D_2^\sharp$
  and $f^\natural:D_1^\natural\to D_2^\natural$.
  If~$f$ is furthermore injective \emph{(}resp.\ surjective\emph{)}
  then so are~$f^\natural$ and $f^{\sharp}$ {\em (}resp.\ $f^{\natural}${\em )}.
\end{lem}
\begin{proof}
Once we know that $f^\sharp$ and~$f^\natural$ exist, the injectivity statement is
obvious, while the surjectivity statement for $f^\natural$
follows from the defining property
of~$D_2^\natural$, since $f(D_1^\natural)$ is a $\psi$-stable $\A^+$-lattice in~$D_2$ if~$f$ is surjective.

We next note that~$f(D_1^\sharp)$ is a bounded
$\A^+$-submodule of $D_2$ on which $\psi$ acts surjectively
(since $D_1^\sharp$ is an $\bA^+$-lattice  in $D_1$ on which $\psi$ acts  surjectively),
and so $f(D_1^\sharp) \subset D_2^{\sharp}$ by Proposition~\ref{prop: D sharp natural exist}~\eqref{item: D sharp exists}, as claimed.
Since $D_1^\natural \subseteq D_1^\sharp,$
we then find that $f(D_1^\natural) \subseteq D_2^\sharp.$
Let~$\gM\coloneqq f^{-1}(D_2^\natural)\cap D_1^\sharp$. Then~$\gM$ is an
$\A^+$-lattice, because it contains~$TD_1^\sharp$ (note that
$f(TD_1^\sharp)=Tf(D_1^\sharp)\subseteq TD_2^\sharp\subseteq
D_2^\natural$ by Proposition~\ref{prop: D sharp natural exist}~\eqref{item: bounding D natural sharp}).
Furthermore~$\gM$
is~$\psi$-stable by definition, so it contains~$D_1^\natural$, as
required.
\end{proof}

\begin{lem}%
  \label{lem: D natural base change}Suppose that $A\to B$ is a morphism
  of Noetherian $\cO/\varpi^a$-algebras, and that~$D_A$ is an
  \'etale $(\varphi,\Gamma)$-module with
  $A$-coefficients.
  Let~$D_B := D_A \cotimes_A B$.
  The following then hold:
  \begin{enumerate}
  \item
  $D_B^\sharp$ contains the image of the 
  natural map $D_A^\sharp\cotimes_AB\to D_B$.
  \item
  $D_B^\natural$ is equal to the image of the
  natural map $D_A^\natural\cotimes_AB\to D_B$.
  \item  If $A \to B$ is surjective, then there is a canonical surjection
  $ D_A^\natural \to D_B^\natural$.
  \item %
  If~$A\to B$ is flat, then there is a canonical isomorphism
  $D_A^\natural\cotimes_AB \iso D_B^\natural$.
  \end{enumerate}
\end{lem}
\begin{proof}
  Since $D_A^\sharp$ is a $\psi$-stable $\A^+_A$-lattice satisfying
  $\psi(D_A^\sharp) = D_A^{\sharp}$,
  we see that the image~$\gN$ of the map
  $D_A^\sharp\cotimes_AB\to D_B$ is a $\psi$-stable 
  $\A_B^+$-lattice for which $\psi(\gN)$ is dense in~$\gN$.
  Since $\psi$ is open, we see that in fact $\psi(\gN) = \gN$,
  and so $\gN \subseteq D_B^\sharp$, proving~(1).
  
  Now write~$\gN'$ for the image of the map
  $D_A^\natural\cotimes_AB\to D_B$; to prove~(2), we must show that $\gN'$
  is equal to $D_B^{\natural}$. By Lemma~\ref{lem: collection
    of facts about lattices}~\eqref{item: lattice 4}, $\gN'$ is a $\psi$-stable
  $\A^+_B$-lattice in~$D_B$, and therefore contains $D_B^\natural$. To see that
  $\gN'$ is contained in $D_B^\natural$, we argue as in the proof of
  Lemma~\ref{lem: D natural induced maps}. Namely, let~$\gM$ be the preimage
  of~$D_B^\natural$ in~$D_A^\natural$ under the natural morphism
  $D_A^\natural\subseteq D_A\to D_B$; we must show that $\gM=D_A^\natural$.

  To see this, note that~$\gM$ is an $\A_A^+$-lattice, because it contains~$T D_A^\sharp$, by
  Proposition~\ref{prop: D sharp natural exist}~\eqref{item: bounding D natural sharp}
  together with part~(1) of the present lemma.
  Since~$\gM$ is $\psi$-stable, we must then have
  $\gM=D_A^\natural$, as required.
  This concludes the proof of part~(2).
  
  As discussed prior to the statement of Lemma~\ref{lem: collection of facts about lattices},
  the completed tensor product
  of~(2) can be reinterpreted as the usual tensor product $\text{--}\otimes_{\A_A^+}{\A_B^+}.$
  Since $\A_A^+ \to \A_B^+$ is surjective when $A\to B$ is, part~(3) then follows.
  
  Finally, part~(4) follows from part~(2) and Lemma~\ref{lem: collection of facts about lattices}~(4).
\end{proof}

\begin{remark}
Part~(4) of Lemma~\ref{lem: D natural base change} shows that the formation of~$D^\natural$ satisfies flat base-change in~$A$.
Unlike~$D^\natural$,
$D^\sharp$ need not satisfy flat base-change
in general.
\end{remark}

\subsection{\texorpdfstring{$D^+$, $D^{++}$, and $D^{\nr}$}{D plus, D plus plus, and D unramified}}
We continue to assume that~$A$ is a Noetherian
$\cO/\varpi^a$-algebra for some~$a\ge 1$,
and write $D$ to denote an \'etale $(\varphi,\Gamma)$-module with coefficients in~$\A_A$.
As in the preceding subsection, 
we will prove various  analogues with coefficients of some results
from~\cite{MR2642409,ColmezMirabolique}. 
Just as was noted in Remark~\ref{rem: not using Gamma}
regarding the results of the preceding subsection, the $\Gamma$-action
continues to play no role in the discussion of this subsection.

\begin{defn}
  \label{defn: Dplus}Let $D^+$ be the $A$-submodule of $D$ consisting
  of those $z\in D$ for which the sequence
$\{\varphi^n(z)\}_{n\ge 0}$ is bounded.
  Let $D^{++}$ be the $A$-submodule of $D$ consisting of those $z \in D$ for which
  the sequence
  $\{\varphi^n(z)\}_{n\ge 0}$ tends to zero (i.e.\ the submodule consisting of those
  $z$ on which $\varphi$ acts topologically nilpotently).
\end{defn}

\begin{lem}
\label{lem:Dplus properties}
$T D^+ \subseteq D^{++} \subseteq D^+$.
\end{lem}
\begin{proof}  The second inclusion is obvious from the definitions,
while the first follows from the fact that if $(a_n)$ is a bounded sequence,
then $(\varphi^n(T) a_n)$ is a sequence converging to~$0$.
\end{proof}

\begin{lem}
  \label{lem: Dplus contained in Dnatural}$D^+\subseteq D^\natural$.
\end{lem}
\begin{proof}
  We follow the proof of~\cite[Prop.\ II.5.14]{ColmezMirabolique}. Fix
  some~$z\in D^+$.
  Then also $\varphi(z) \in D^+$, and so 
  Lemma~\ref{lem:Dplus properties} shows that $T\varphi(z) \in D^{++}$; 
  i.e.\ that the sequence
  $\{\varphi^n(T\varphi(z))\}_{n\ge 0}$ tends to zero.  Since $D^\natural$ is a lattice,
  we thus find that
  $\varphi^n(T\varphi(z))\in D^\natural$ for some sufficiently large value of~$n$.
  Since $D^\natural$ is
  $\psi$-stable, and since
  $\psi^{n+1}(\varphi^n(T\varphi(z)))=\psi(T\varphi(z))=\psi(T)z=-z$,
  we thus have $z\in D^\natural$, as required.
\end{proof}

The following corollary provides a useful characterization of~$D^+$.

\begin{cor}
  \label{cor: Dplus is a lattice}$D^+$ is an $\A^+_A$-lattice in~$D$,
and is the {\em (}unique{\em )} maximal $\varphi$-invariant bounded subset of~$D$. 
\end{cor}
\begin{proof}
  We saw in Lemma~\ref{lem:Dplus properties} that $TD^+ \subseteq D^+.$
  Thus, if we show that $D^+$ is a lattice in~$D$, it is in fact an~$\A^+_A$-lattice. 

  By definition any $\varphi$-bounded subset of~$D$ is contained in~$D^+$.
  In particular, if~$\gM$ is any
  $\varphi$-stable lattice in~$D$, then $\gM\subseteq D^+$. Such a
  lattice exists, by~\cite[Lem.\ 5.2.7]{EGstacktheoreticimages},
  and thus $D^+$ contains a lattice, and so is open.
  By Remark~\ref{rem: open and bounded equals lattice}, this proves the lemma, provided that we show that~$D^+$ is also bounded;
  but this follows from
  Lemma~\ref{lem: Dplus contained in Dnatural}.
\end{proof}

Our goal in what follows is to show that the formation
of each of $D^{++}$
and $D^+$ is compatible with flat base change.  In the case
of $D^{++}$ this is relatively straightforward, but the case of $D^+$
is more involved.

A general observation is that,
when analyzing constructs such as $D^+$
and~$D^{++}$ (and the related construct $D^{\nr}$ of
Definition~\ref{def:Dnr} below) whose definitions involve dynamical aspects of the action
of~$\varphi,$ complications can sometimes arise from the existence of
a descending chain of ideals
in $A$ that do not stabilize. 
The following lemma (which is presumably standard)
sometimes allows us to get around this problem, by reducing our
analysis to the case of an Artinian coefficient ring, for which this
complication cannot occur.

\begin{lemma}
\label{lem:Artinian replacement}
If $A$ is a Noetherian ring, then there
exists an embedding $A \hookrightarrow B$ with $B$ Artinian.
\end{lemma}
\begin{proof}
First recall that if $f$ is any element of~$A$, then 
for $n$ large enough, 
the natural map $A \to A_f \times A/f^n$ is an embedding.
Now, if $A$ is non-zero, 
and if $\mathfrak p$ is a minimal prime of~$A$,
then we may choose an element $f \in A$ such that the distinguished open $D(f)$ is
a neighbourhood of $\mathfrak p$ in $\Spec A$ which contains no other associated
primes of~$A$.  Thus $\mathfrak p$ is the unique associated prime of $A_f$,
and so the natural map $A_f \to A_{\mathfrak p}$ is an embedding.
Thus, for $n$ large enough,
the natural map $A \to A_{\mathfrak p} \times A/f^n$ is an embedding. 
Since $\mathfrak p$ was chosen to be minimal, the localization $A_{\mathfrak p}$ is Artinian.
 
We continue by applying the same argument to $A/f^n$. 
By Noetherian induction (in $\Spec A$) the process eventually terminates,
and we obtain an embedding of $A$ into a finite product of Artinian local rings,
as required.
\end{proof}

We now return to the setting of \'etale $(\varphi,\Gamma)$-modules
over Noetherian $\cO/\varpi^a$-algebras.

\begin{lemma}
\label{lem:phi is injective on quotients}
The induced action of $\varphi$ on each of $D/D^+$ and $D/D^{++}$ is injective.
\end{lemma}
\begin{proof}
This follows from the fact that each of $D^+$ and $D^{++}$ is $\varphi$-saturated in~$D$
(i.e.\ if $\varphi(x)$ is in $D^+$, resp.\ $D^{++}$, then so is $x$ itself),
as is evident from their definitions.
\end{proof}

The preceding lemma admits a kind of converse, characterizing~$D^{++}$.

\begin{lemma}
\label{lem:characterizing Dplusplus}
 $D^{++}$ is an $\A_A^+$-lattice in~$D$,
and is furthermore
the {\em (}unique{\em )} minimal $\varphi$-invariant
lattice $\gM$ of~$D$ 
for which $\varphi$ acts injectively on~$D/\gM.$
\end{lemma}
\begin{proof}
It follows from
Lemma~\ref{lem:Dplus properties}
and Corollary~\ref{cor: Dplus is a lattice}
that $D^{++}$ is a $T$-invariant lattice,
and thus is an $\A^+_A$-lattice in~$D$.
The injectivity of~$\varphi$ on $D/D^{++}$  was proved in Lemma~\ref{lem:phi is injective
on quotients}.
  If $\gM$ is any $\varphi$-invariant lattice in~$D$,
and if $x \in D^{++}$, then by definition we see that
$\varphi^n(x) \in \gM$ for some value of~$n$.  Thus,
if $\varphi$ is injective on~$D/\gM$,
then we see that
$x \in \gM$, and hence that $D^{++} \subseteq \gM$.
\end{proof}

We also note the following technical property 
of~$D^{++}$.

\begin{lemma}
\label{lem:Dplusplus uniformity}
The action of $\varphi$ on~$D^{++}$ is {\em uniformly} topologically nilpotent,
i.e.\ given any lattice $\gM \subset D$, there is some $n$ such that $\varphi^n(D^{++}) \subset \gM$.
\end{lemma}
\begin{proof}
Any lattice $\gM$
contains an $\A_A^+$-invariant sublattice, and replacing $\gM$ by such a sublattice
if necessary, we may assume that~$\gM$ is itself an $\A_A^+$-lattice.
Since $D^{++}$ is a lattice, 
by Lemma~\ref{lem:characterizing Dplusplus},
it is finitely generated over~$\A_A^+$, say by elements
$x_1,\ldots,x_m$. %
By definition, we may find some $n$ such that $\varphi^n(x_i) \in \gM$ for each~$i$.
Then
$$\varphi^n(D^{++}) = \varphi^n( \A_A^+\langle x_1,\ldots,x_m\rangle)
\subseteq \A_A^+\langle \varphi^n(x_1),\ldots,\varphi^n(x_m)\rangle \subseteq \gM,$$
as required.
\end{proof}

It is now easy to establish flat base-change for~$D^{++}$.

\begin{lemma}
\label{lem:flat base-change for Dplusplus}
If $A \to B$ is a flat morphism of Noetherian $\cO/\varpi^a$-algebras, if $D_A$ is an \'etale
$(\varphi,\Gamma)$-module over $\A_A$, and if we write
$D_B \coloneqq  D_A\cotimes_A B,$ then the resulting morphism
$D_A^{++} \cotimes_A B \to D_B$ induces an isomorphism
$D_A^{++} \cotimes_A B \iso D_B^{++}$.
\end{lemma}
\begin{proof}
Since $A \to B$ is flat, 
and since $D_A^{++}$ is a lattice in~$D_A$ by Lemma~\ref{lem:characterizing Dplusplus},
it follows from
Lemma~\ref{lem: collection of facts about lattices}~\eqref{item: lattice 4}
that $D_A^{++}\cotimes_A B
\to D_B$ identifies its source with a lattice in~$D_B$.
Using the uniform topological nilpotency of Lemma~\ref{lem:Dplusplus uniformity},
we see that $\varphi$ acts topologically nilpotently on the completed tensor
product~$D_A^{++}\cotimes_AB$,
and so this lattice is contained in~$D_B^{++}$.

Furthermore, since $B$ is flat over $A$, and since 
$\varphi$ is injective on $D_A/D_A^{++}$ by Lemma~\ref{lem:phi is injective on quotients},
we see that $\varphi$ is injective on $(D_A/D_A^{++})\otimes_A B
= D_B/(D_A^{++}\cotimes_A B).$
(The indicated identification is provided by
Lemma~\ref{lem: collection of facts about lattices}~\eqref{item: lattice 4}.)
Lemma~\ref{lem:characterizing Dplusplus}
now shows the reverse inclusion,
namely that $D_B^{++} \subseteq D_A^{++}\cotimes_A B.$
This proves the lemma.
\end{proof}

Our next goal is to prove Lemma~\ref{lem:flat base-change for Dplus},  which shows that~$D^+$ is compatible with flat base-change for~$D^+$.

\begin{lemma}
\label{lem:Cartesian diagram for Dplus}
Let $A\hookrightarrow B$ is an embedding of Noetherian $\cO/\varpi^a$-algebras,
and  $D_A$ is an \'etale $(\varphi, \Gamma)$-module with~$A$ coefficients, 
so that the embedding of $A$ into  $B$  induces an embedding
$D_A \hookrightarrow D_B \coloneqq  D_A \cotimes_A B,$
then we
have $D_A^+ = D_A \cap D_B^+.$
Equivalently, the natural map $D_A/D_A^+ \to D_B/D_B^+$ is injective.
\end{lemma}
\begin{proof}
Since $D_A^+$ is a bounded $\varphi$-invariant subset of $D_B$, it is contained
in $D_B^+$, and thus in the indicated intersection. 
On the other hand, this intersection is a bounded $\varphi$-invariant
subset of $D_A$ (using Lemma~\ref{lem: collection of facts about lattices}~\eqref{item: lattice 5} to see the boundedness), and thus is contained in~$D_A^+$, by Corollary~\ref{cor: Dplus is a lattice}.
\end{proof}

\begin{lemma}
\label{lem:phi expands}
Let $D$ be an \'etale $(\varphi,\Gamma)$-module over the Noetherian %
$\cO/\varpi^a$-algebra~$A$.
Then for all sufficiently large $n$,
we have $\varphi^{-n}\bigl( (D/D^+)[T] \bigr) \cap \bigl((D/D^+)[T]\bigr) = 0.$
\end{lemma}
\begin{proof}
Let $M_n$ denote the kernel of the endomorphism $T \varphi^n$ of $D/D^+$;
we must show that $M_n[T] = 0$ if $n$ is sufficiently large.
Applying Lemmas~\ref{lem:Artinian replacement} 
and~\ref{lem:Cartesian diagram for Dplus},
we see that we may replace $A$ by an Artinian overring,
and thus we assume that $A$ is Artinian for the remainder of the proof.

Since $\varphi$ is injective on $D/D^+$, by Lemma~\ref{lem:phi is injective on quotients},
we find that
$$M_n \coloneqq  \ker T\varphi^n = \ker \varphi(T \varphi^n) =  
\ker \varphi(T) \varphi^{n+1} \supseteq \ker T \varphi^{n+1} \eqcolon M_{n+1}.$$
Thus $M_n[T]$ is a descending sequence of $A$-submodules of the finite type
$A$-module $(D/D^+)[T]$, which therefore stabilizes.
Suppose that the sequence stabilizes at~$n_0$,
and write $N \coloneqq  M_{n_0}[T] = M_n[T]$ for $n \geq n_0.$
Then we see (by definition of the $M_n$) that $\varphi^n(N) \subset (D/D^+)[T]$
for $n\geq n_0$.  Thus $\bigcup_n \varphi^n(N)$ is a bounded subset of $D/D^+$ (i.e.\ it is contained in a finite type $A$-module)
and so $N = 0$ (by the definition of $D^+$).
This is what we had to prove.  
\end{proof}

\begin{lemma}
\label{lem:flat base-change for Dplus}
If $A \to B$ is a flat morphism of Noetherian $\cO/\varpi^a$-algebras, if $D_A$ is an \'etale
$(\varphi,\Gamma)$-module over $\A_A$, and if we write
$D_B \coloneqq  D_A\cotimes_A B,$ then the resulting morphism
$D_A^+ \cotimes_A B \to D_B$ induces an isomorphism
$D_A^+ \cotimes_A B \iso D_B^+$.
\end{lemma}
\begin{proof}
Since $A \to B$ is flat, 
and since $D_A^{+}$ is a lattice in~$D_A$ by Corollary~\ref{cor: Dplus is a lattice},
it follows from
Lemma~\ref{lem: collection of facts about lattices}~\eqref{item: lattice 4}
that $D_A^{+}\cotimes_A B
\to D_B$ identifies its source with a lattice in~$D_B$.
This lattice is furthermore~$\varphi$-stable, 
and thus is contained in~$D_B^+$ (by another application of Corollary~\ref{cor: Dplus is a lattice}).
Our goal is to show that in fact $D_A^+\cotimes_A B = D_B^+$,
or equivalently that 
$D_B^+/(D_A^+\cotimes_A B) = 0.$
Since every element of $D_B/(D_A^+\cotimes_A B)$ is $T$-power torsion,
it in fact suffices to show that
$\bigl(D_B^+/(D_A^+\cotimes_A B)\bigr)[T] = 0,$
and this is what we will do.

Lemma~\ref{lem:phi expands}
shows that the morphism
$(D_A/D_A^+)[T] \to D_A/D_A^+$ induced by $T \varphi^n$
is injective if~$n$ is sufficiently large.
Base-changing over the flat $A$-algebra~$B$,
we find that the morphism
$\bigl(D_B/(D_A^+\cotimes_A B)\bigr)[T] \to D_B/(D_A^+\cotimes_A B)$ induced by $T \varphi^n$
is again injective if~$n$ is sufficiently large.
(Here and below we use 
Lemma~\ref{lem: collection of facts about lattices}~\eqref{item: lattice 4}
to make the evident manipulations of various tensor products and completed
tensor products with~$B$.)
Composing with $\varphi^r$ for any~$r$, and recalling Lemma~\ref{lem:phi is injective on
quotients},
we that the same is true of the morphism
$\bigl(D_B/(D_A^+\cotimes_A B)\bigr)[T] \to D_B/(D_A^+\cotimes_A B)$ induced by $\varphi^r(T) \varphi^n$,
if $n$ is sufficiently large (depending on~$r$).
Thus no non-zero element of $\bigl(D_B/(D_A^+\cotimes_A B)\bigr)[T]$ has a bounded $\varphi$-orbit;
in fact, since $D_A^+ \cotimes_A B$ is a lattice in $D_B$, the bounded subsets of $D_B^+/D_A^+ \cotimes_A B$
are precisely those contained in some $B$-module of finite type.
So we see that indeed
$\bigl(D_B^+/(D_A^+\cotimes_A B)\bigr)[T] = 0$,  
as required.
\end{proof}

Finally, we introduce and study some basic properties of~$D^{\nr}$ in our context of $(\varphi, \Gamma)$-modules with coefficients.

  \begin{defn}
\label{def:Dnr}
If $D$ is an \'etale $(\varphi,\Gamma)$-module with coefficients in~$\A_A$,
for some Noetherian $\cO/\varpi^a$-algebra~$A$,
then we write~$\displaystyle D^{\nr}\coloneqq \bigcap_n\varphi^n(D)$, equipped with its
actions of~$\varphi$ and~$\Gamma$.
  \end{defn}

\begin{rem}
In the case that~$A=\cO/\varpi^a$, the $A$-module $D^{\nr}$ is related to the
maximal abelian subrepresentation %
of the corresponding Galois representation 
(see~\cite[Rem.\ II.1.2]{ColmezMirabolique}, and Lemma~ \ref{Dnr and V^ab} %
below).
Since this notion does not behave well in the context of general families
of \'etale $(\varphi,\Gamma)$-modules (e.g.\ because of the fundamental
phenomenon of generically reducible families specializing to
irreducible objects), it should not come as a surprise that %
$D^{\nr}$ does
not behave well in families; 
for example, it is not compatible with
flat base-change in general. %
\end{rem}

\begin{lemma}
\label{lem: Dnr is finite}
$D^{\nr}$ is a finite type $A$-submodule of~$D^+$.
\end{lemma}
\begin{proof}
Clearly $D^{\nr}$ is $\varphi$-stable, while
since $\varphi$ is injective on $D$ we see that $D^{\nr}$ is also $\varphi$-saturated.
Since $D^{\nr} \subseteq \varphi(D)$ by definition,
we then see that $\varphi(D^{\nr}) = D^{\nr}$.

Since $D^{++}$  is also both $\varphi$-stable and $\varphi$-saturated,
we see that $D^{\nr}  \cap D^{++}$  is again $\varphi$-stable and $\varphi$-saturated.
It then follows that also $\varphi(D^{\nr} \cap D^{++}) = D^{\nr}\cap D^{++}.$
The uniform topological nilpotency of Lemma~\ref{lem:Dplusplus uniformity}
then shows that $D^{\nr}\cap D^{++} = 0.$

Now  consider the image $M$ of $D^{\nr}$ in~$D/D^+$;  
evidently $\varphi(M) = M$.
Suppose that $m \in M[T],$ and let $m' \in M$ be such
that $\varphi(m') = m.$  Then $\varphi(Tm') = \varphi(T) \varphi(m') = \varphi(T) m = 0,$
and thus $Tm' = 0$ (by Lemma~\ref{lem:phi is injective on quotients}). 
Consequently we see that $\varphi(M[T]) \supseteq M[T]$, %
and hence that $\varphi^n(M[T]) \supseteq M[T]$ for all $n \geq 0.$
Lemma~\ref{lem:phi expands} then implies that 
$M[T] = 0,$  and so we see that in fact $M = 0$.
Thus $D^{\nr} \subseteq D^+$.
(In more detail: if $x \in M[T]$, then the inclusion $\varphi^n(M[T]) \supseteq M[T]$ shows that there exists $y \in M[T]$ such that $\varphi^n(y) = x$.
Hence $y \in \varphi^{-n}\left ( (D/D^+)[T] \right ) \cap (D/D^+)[T]$, which by Lemma~\ref{lem:phi expands} is zero for~$n$ large enough, 
and so $x = \varphi^n(y)$ is also equal to~$0$.)

Putting what we've shown together, we see that $D^{\nr}$ embeds into
$D^+/D^{++}$, which is a finite type $A$-module by Lemma~\ref{lem:Dplus properties}
together with Corollary~\ref{cor: Dplus is a lattice}, 
and hence $D^{\nr}$ is itself a finite type $A$-module.
\end{proof}

\begin{cor}
  \label{cor: Dnr contained in Dnatural}We have $D^{\nr}\subset D^\natural$.
\end{cor}
\begin{proof} %
This is immediate from Lemmas~\ref{lem: Dnr is finite} and  %
~\ref{lem:
    Dplus contained in Dnatural}.
  \end{proof}

\subsection{From \texorpdfstring{$(\varphi,\Gamma)$}{(ϕ,Γ)}-modules to equivariant sheaves}
\label{subsec:sheaves}
We now explain how some of the constructions of~\cite[\S III.1, \S
V]{ColmezMirabolique} extend to our setting.
We apply some results from~\cite[\S 3]{MR3444236},  %
in which some of the more formal aspects of Colmez's constructions
are presented in a natural level of generality.
With these general results in hand,
the extension of Colmez's constructions to the case of \'etale
$(\varphi,\Gamma)$-modules with coefficients is for the
most part immediate, and we frequently refer
to~\cite{ColmezMirabolique} for the proofs of results when they are
literally identical in our setting.

We write $P^+$ for the monoid $
\big(\begin{smallmatrix}
  \Zp\setminus\{0\}&\Zp\\0&1
\end{smallmatrix}\big)
$, and~$P$ for the group~$
\big(\begin{smallmatrix}
  \Qptimes &\Qp\\0&1
\end{smallmatrix}\big)
$. Then there is a natural action of~$P^+$ on~$\Zp$ (respectively
of~$P$ on~$\Qp$) via $
\big(\begin{smallmatrix}
  a &b\\0&1
\end{smallmatrix}\big)x=ax+b
$.
The basic motivation for Colmez's constructions is as follows:  we can regard
an \'etale $(\varphi,\Gamma)$-module $D$ as being the global sections
of a $P^+$-equivariant sheaf of~$A$-modules on~$\Z_p$.  
We are then %
able to extend its $P^+$-equivariance to an equivariance
under the category of all non-degenerate\footnote{In the sense of having non-zero
derivative at every point.} piecewise affine-linear maps between open subsets 
of~$\Z_p$, %
and then (by an appropriate limiting process) 
to the category of all {\em local diffeomorphisms} (in the
sense of Definition~\ref{def:local diffeo} below) 
between open subsets of~$\Z_p$. 

With these constructions in place, it becomes possible to localize~$D$ 
in an equivariant  manner over other locally analytic~$\Z_p$-manifolds, such as $P$-equivariantly
over $\Q_p$, and (our ultimate goal) $\GL_2(\Q_p)$-equivariantly over~$\Pone(\Q_p)$.
In fact, the exposition becomes simpler if we construct the localization over $\Q_p$ 
en route to the construction of the equivariant structure under local diffeomorphisms, 
and so we do this.
We discuss localization over $\Pone(\Q_p)$ in the next subsection. %

Suppose now that $M$ is any $P^+$-module with coefficients in a Noetherian $\cO/\varpi^a$-algebra~$A$.
Write $U^+ \coloneqq  
\big(\begin{smallmatrix}
  1 &\Zp\\0&1
\end{smallmatrix}\big).
$
Let $\varphi$ denote the endomorphism 
$$\varphi:\left(\begin{smallmatrix} 1 & x \\ 0 & 1 \end{smallmatrix}\right)
\mapsto 
\left(\begin{smallmatrix} 1 & px \\ 0 & 1 \end{smallmatrix}\right)$$
of~$U^+$,
as well as the endomorphism of $M$ induced by the action of
$\left(\begin{smallmatrix} p & 0 \\ 0 & 1 \end{smallmatrix}\right).$
Then if $u \in U^+$ and $m\in M$, we find that $\varphi(um)  = \varphi(u)\varphi(m).$
Consequently for each~ $n\ge 1$ we see that $\varphi^n(M)$ is $\varphi^n(U^+)$-stable,
and so we obtain an induced morphism
\begin{equation}
\label{eqn:Phi maps}
A[U^+]\otimes_{A[\varphi^n(U^+)]} \varphi^n(M) \to M.
\end{equation}

\begin{defn}
  \label{defn:etale-Pplus-module}We say that a  $P^+$-module~$M$ with
  coefficients in~$A$ is \emph{\'etale} if the morphisms~\eqref{eqn:Phi maps}
  are isomorphisms for all~$n\ge 1$.
\end{defn}

Our interest in \'etale $P^+$-modules comes from the following lemma.
\begin{lem}
\label{lemma:Pplus action on D}
If~$D$ is an \'etale $(\varphi,\Gamma)$-module over a Noetherian $\cO/\varpi^a$-algebra~$A$,
then $D$ is an \'etale $P^+$-module with $A$-coefficients, where the
action of~$P^+$ on~$D$ is via the continuous maps\[
  \begin{pmatrix}
    p^ka&b\\0&1
  \end{pmatrix}z=(1+T)^b\varphi^k\circ\sigma_a(z)
\] for $a\in\Zptimes$, $b\in\Zp$,
$k\in\Z_{\ge 0}$, and $z \in D$. \emph{(}Recall that~$\sigma_a\in\Gamma$ is the element with
$\chi(\sigma_a)=a$, where~$\chi$ is the cyclotomic character.\emph{)}
\end{lem}
\begin{proof}
  This is immediate from the definitions. Indeed by definition, the condition
  that  $D$ is \'etale as a $(\varphi,\Gamma)$-module is equivalent to the
   morphisms~\eqref{eqn:Phi maps} being isomorphisms.
 \end{proof}
 By~\cite[Thm.\ 3.32]{MR3444236}, there is an equivalence of abelian categories
 between the category of \'etale $P^+$-modules with $A$-coefficients, and the
 category of $P$-equivariant sheaves of $A$-modules on~$\Qp$; the inverse
 functor is given by passage to the module of sections over~$\Zp$. Following
 Colmez, for any open set~$U\subset\Qp$, we write $D\boxtimes U$ for the
 sections over~$U$ of the sheaf corresponding to~$D$, so that in particular we
 have  $D\boxtimes\Zp=D$. As is the case for any $P$-equivariant
 sheaf of~$A$-modules on $\bQ_p$, there is a $P$-equivariant action on~$D\boxtimes\Qp$ of
 $\cC(\Qp,A)$, the ring of $A$-valued locally constant functions
 on~$\Qp$, and thus there is in particular a $P$-equivariant action of the ring
 $\cC_c(\Qp,A)$  of compactly supported $A$-valued locally constant
 functions. Similarly for each open ~ $U\subseteq \Qp$ there is a natural action
 of $\cC_c(U,A)$ on~$D\boxtimes U$.

 If~$U\subseteq \Qp$ is open then we write $\Res_U:D\boxtimes\Qp\to D\boxtimes
 U$ for the restriction map. If~$U\subseteq V\subseteq \Qp$ are open subsets
 then we have the restriction map $\Res^V_U:D\boxtimes V\to D\boxtimes U$,
 which we will usually %
 denote simply by~$\Res_U$.  If~$U$ is compact then the
 restriction map~$\Res_U:D\boxtimes\Qp\to D\boxtimes U$ admits a section
 given by extension by zero, and we accordingly identify the action of the locally constant function
 ~$1_U\in \cC_c(\Qp,A)$ with~$\Res_U$.

 By  ~\cite[Prop.\
  3.26]{MR3444236} (and its proof), the global sections~ $D\boxtimes\Qp$ are given explicitly by \begin{equation}\label{eqn:D-box-Qp-as-inverse-limit}D\boxtimes\Q_p = \varprojlim_{\psi} D,\end{equation}
where the right hand side denotes the inverse limit of the inverse system
$$\cdots \buildrel \psi \over \longrightarrow D
\buildrel \psi \over \longrightarrow D
 \buildrel \psi \over \longrightarrow \cdots;$$
equivalently,
$$D\boxtimes\Q_p \coloneqq  \{(z_n)_{n \geq N} \, | \, \psi(z_{n+1}) = z_n \text{ for all } n \}.$$(Here $N$ can be any element of $\Z\cup\{-\infty\};$
we take advantage of this flexibility in some of the constructions introduced
below.) We can therefore endow~$D\boxtimes\Qp$ with the projective limit
topology (where~$D$ has its canonical Tate module topology). Note that if~$U$ is compact
and open, then ~$D\boxtimes U$ also has a canonical Tate module topology; indeed
$D\boxtimes U$ is a direct summand of some~$D\boxtimes\frac{1}{p^n}\Zp$, which
is homeomorphic to~$D\boxtimes\Zp=D$.

Furthermore, the action of~$P$ on  $D\boxtimes\Qp$  is determined by
  the following properties. %
  \begin{enumerate}
  
  \item The restriction map $\Res_{\Zp}:D\boxtimes\Qp\to D = D\boxtimes \Z_p $
    is given by \[(z_n)_{n\in\Z}\mapsto
     z_0\in D=D\boxtimes\Zp.\]  \item The extension by zero map
    $D = D\boxtimes \Z_p \to D\boxtimes \Q_p$ is given by \[z \mapsto
      \left\{\bigl(\varphi^n(z)\bigr)_{n \geq 0}\right\};\] this  is continuous and $P^+$-equivariant.
  \item  $\left(\begin{smallmatrix} p & 0 \\ 0 & 1\end{smallmatrix}\right)$ acts
    via $(z_n)_{n \in\Z} \mapsto (z_{n+1})_{n \in \Z}.$ 
      \item the elements
    $\left(\begin{smallmatrix} 1 & b \\ 0 & 1\end{smallmatrix}\right)$, for
    $b \in p^{-N}\Z_p$, act via the transformation \[(z_n)_{n \geq N} \mapsto \bigl((1+T)^{p^n b}z_n\bigr)_{n \geq N}.\]
  \end{enumerate}
  It follows from the definitions that this action of~$P$ is continuous. 

\begin{defn}
  \label{defn:restriction-maps-and-m-alpha}If $\alpha\in \cC_c(U,A)$,  we
  write~$m_{\alpha}$ %
  for the corresponding continuous endomorphism of~$D\boxtimes U$. %
\end{defn}

\begin{lem}
  \label{lem:properties-of-Res-map-construction}\leavevmode
  \begin{enumerate}
  \item For any~$n\ge 0$ the restriction map $\Res_{p^n\Zp}:D\boxtimes\Zp\to
    D\boxtimes p^n\Zp$ is equal to  $\varphi^n\circ\psi^n$.
  \item Let  $\alpha:U\to A$ be locally constant, and let~$n\ge 0$ be large
    enough that for
  all~$a\in U$, we have $a+p^{n}\Zp\subseteq U$, and
  $\alpha|_{a+p^n\Zp}$ is constant. Let~$I_n(U)$
for a system of coset representatives for~$U$ modulo~$p^n\Zp$. Then  \[m_\alpha=\sum_{i\in
      I_n(U)}\alpha(i)\Res_{i+p^n\Zp}.\]%
  \end{enumerate}
\end{lem}
\begin{proof}
  The first part follows from the identity
  $\Res_{p^n\Zp}=\left(\begin{smallmatrix} p^n & 0 \\ 0 &
                                                          1\end{smallmatrix}\right)\circ
                                                      1_{\Zp}\circ
                                                      \left(\begin{smallmatrix}
                                                        p^{-n} & 0 \\ 0 &
                                                                          1\end{smallmatrix}\right)$,
                                                                      and the
                                                                      second is
                                                                      immediate
                                                                      from the
                                                                      identity $\alpha=\sum_{i\in
      I_n(U)}\alpha(i)1_{i+p^n\Zp}$
\end{proof}

\begin{example}
For future reference,
note that~$D\boxtimes\Zptimes=D^{\psi=0}$, with the identification being
induced via the embedding $D\boxtimes\Zptimes \hookrightarrow D\boxtimes \Z_p \coloneqq  D$
given by extension by zero. 
Indeed, the sections
of~$D=D\boxtimes\Zp$ which are supported on~$\Zptimes$ are precisely
those in the kernel of $\Res_{p\Zp}=\varphi\circ\psi$, and~$\varphi$ is injective.
\end{example}%

Our next goal is to extend the $P$-equivariant localization of $D$ over $\Q_p$ 
to an equivariant structure under the category of local diffeomorphisms.
As already noted, this will then allow us to localize $D$ over more general
one-dimensional $p$-adic manifolds, such as $\Pone(\Q_p)$.

\begin{df}
\label{def:local diffeo}
If~$U,V$ are compact open subsets of~$\Qp$, then we say that a map
$f:U\to V$ is a local diffeomorphism if it is~$\cC^1$ (in the sense of~\cite[Section~I.5.1]{Colmezfonctions}),
and if its derivative is non-vanishing on~$U$.
\end{df}

\begin{df}
  \label{def:uniform convergence}
  For all~$n \geq 0$, let $X_n  \subset D$ be a subset.
  We say that the sequence $X_n$ tends uniformly to 0 as $n \to \infty$
  if for any $A$-lattice $\fM \subset D$ there exists~$N>0$ such that $X_n \subset \fM$ for all~$n > N$.
  (This is a slight rephrasing of the definition in~\cite[Section~V.1.2]{ColmezMirabolique}.)

  Fix~$m > 0$, and for all~$n \geq 0$, let $X_n \subset D \boxtimes \frac 1 {p^m} \bZ_p$ be a subset.
  We say that the sequence $X_n$ tends uniformly to~$0$ as $n \to \infty$ if its image under $\varphi^m: D \boxtimes \frac 1 {p^m} \bZ_p \isoto D$ tends uniformly to~$0$. 

\end{df}

\begin{prop}
  \label{prop: direct image of local diffeomorphism}
  Suppose that
  $f:U\to V$ is a local diffeomorphism. Then there is a continuous
  $A$-linear morphism $f_*:D\boxtimes U\to D\boxtimes V$, defined as follows:
For each sufficiently large~$n$, write~$I_n(U)$
for a system of coset representatives for~$U$ modulo~$p^n\Zp$,
and then for each $z\in D\boxtimes U$, define
  \[f_*(z) \coloneqq \lim_{n\to\infty}\sum_{i\in I_n(U)}
    \begin{pmatrix}
      f'(i)&f(i)\\0&1
    \end{pmatrix}\Res_{p^n\Zp}\left (
    \begin{pmatrix}
      1&-i\\0&1
    \end{pmatrix}z
\right );\] this limit exists, and is independent of the choices
of~$I_n(U)$.

Furthermore, %
\[
\left \{ \Res_{j+p^n\Zp}\left(f_*(z)-\sum_{i\in I_n(U)}
    \begin{pmatrix}
      f'(i)&f(i)\\0&1
    \end{pmatrix}\Res_{p^n\Zp}\left(
    \begin{pmatrix}
      1&-i\\0&1
    \end{pmatrix}z
\right)\right): j \in V \right \}\]
tends uniformly to zero as~$n\to\infty$.
\end{prop}
\begin{example}
\label{ex:undoing restriction to pZp} 
Let $f : \bZ_p \to p \bZ_p$ be multiplication by~$p$.
Then
\[
f_*(z) = \lim_{n \to \infty} \sum_{i\in I_n(U)} \fourmatrix p 0 0 1 \Res_{i+p^n\bZ_p}(z) = \fourmatrix p 0 0 1 z = \varphi(z)
\]
and so $f_* = \varphi : D \to D \boxtimes p\bZ_p = \varphi(D)$.

Similarly, let $f : p\Z_p \buildrel \over \longrightarrow \Z_p$ be multiplication by~$p^{-1}$.
Then 
\[
f_*: \varphi(D) \eqcolon  D\boxtimes p\Z_p \to D\boxtimes \Z_p \coloneq D
\]
is given by $\varphi^{-1}: \varphi(D) \to D.$ 
Thus the composite
$$D \eqcolon  D\boxtimes \Z_p \buildrel \Res_{p\Z_p} \over
\longrightarrow D\boxtimes p\Z_p \buildrel f_* \over \longrightarrow D\boxtimes\Z_p \coloneq D$$
is given by $\varphi^{-1}\circ \varphi \psi,$ i.e.\ by the operator~$\psi$.
\end{example}

\begin{proof}[Proof of Proposition~\ref{prop: direct image of local diffeomorphism}]Note firstly that if the limit exists, it is
  necessarily independent of the choices of~$I_n(U)$. Indeed
  if~$I_n(U)$, $I'_n(U)$ are two such choices, this follows by
  considering the limit for both choices together with the third
  choice~$I''_n(U)$ defined by ~$I''_n(U)=I_n(U)$ if~$n$ is even and
  $I'_n(U)$ if~$n$ is odd.
  
  Exactly as in the proof of~\cite[Prop.\ V.1.3]{ColmezMirabolique},
  we can reduce the statement of the proposition to the case that
  $U=V=\Zp$, and that~$f$ is \emph{regular} in the sense that
  $v_p(f'(x))=0$ for all~$x\in\Zp$. If~$A$ is a finite $\cO/\varpi^a$-algebra, this
  case is~\cite[Lem.\ V.1.2]{ColmezMirabolique}, and essentially the
  same proof works in our setting, as we now explain. Write \[u_n=\sum_{i\in I_n(U)}
    \begin{pmatrix}
      f'(i)&f(i)\\0&1
    \end{pmatrix}\Res_{p^n\Zp}\left (
    \begin{pmatrix}
      1&-i\\0&1
    \end{pmatrix}z
\right).\] We will show next that for any $\A^+$-lattice~$\gM$ containing~$z$,
$u_n-u_{n-1}$ tends to~$0$ uniformly in~$z \in \gM$; it then follows that
$f_*(z)$ exists, and that~$f_*$ is $A$-linear and continuous.

For~$i\in\Zp$, write $r_{n,i}(z)=\psi^n(
\big(\begin{smallmatrix}
  1&-i\\0&1
\end{smallmatrix}\big)
z)$. The argument of the second paragraph of the proof of ~\cite[Lem.\
V.1.2]{ColmezMirabolique} goes through unchanged, and shows that we
can write \begin{equation}\label{un minus u n minus 1} u_n-u_{n-1}=\sum_{j\in I_n(\Zp)}
  \begin{pmatrix}
    1&f(j)\\0&1
  \end{pmatrix}\varphi^n((g_j-h_{n,j})\cdot r_{n,j}(z)),\end{equation}
where $g_j,h_{n,j}\in P^+$ (and the difference
$g_j-h_{n,j}$ is as elements of the group ring, rather than as
elements of~$P^+$). Furthermore, $g_j,h_{n,j}$ enjoy the following
properties: we have $g_j\in
\big(\begin{smallmatrix}
  \Zptimes &0\\0&1
\end{smallmatrix}\big)
$, and $g_j^{-1}h_{n,j}\in
\big(\begin{smallmatrix}
  1+p^{a(n)}\Zp & p^{a(n)}\Zp\\ 0&1
\end{smallmatrix}\big)
$ for some integer~$a(n)$, and $a(n)\to \infty$ as $n\to\infty$.

Let~$\gM$ be an $\A^+$-lattice containing~$z$. Then there is a $\Gamma$-stable
$\A^+$-lattice~$\gM'\supseteq\gM$ such that $r_{n,j}(\gM)\subseteq\gM'$ for
all~$j\in\Zp$, $n\ge 0$. 
Indeed $\big(\begin{smallmatrix}
  1&-j\\0&1
\end{smallmatrix}\big)\gM=(1+T)^j\gM=\gM$ for all~$j\in\Zp$, and by
Lemma~\ref{lem: psi pushes into T inverse D natural}, there is an
$\A^+$-lattice~$\gM''$ containing all the~$\psi^n(\gM)$.
Then we can take $\fM'$ to be the $\bA^+$-submodule~$\Gamma \fM''$ of~$D$ generated by $\gamma m$ for~$\gamma \in \Gamma, m \in \fM''$, which is a lattice by compactness of~$\Gamma$
(see~\cite[Lemma~5.1.5]{emertongeepicture}).

Writing \[(g_j-h_{n,j})\cdot
  r_{n,j}(z)= g_j\cdot (1-g_j^{-1}h_{n,j})\cdot r_{n,j}(z), \] we see
that it suffices to show that there is a nested sequence of $\A^+$-lattices~$\gM_n$
with the properties that for all~$j$ we have~$(1-g_j^{-1}h_{n,j})(\gM')\subseteq \gM_n$,
that the intersection of the~$\gM_n$ is~$0$, and that~$\gM_n$ is
$\varphi$-stable and $\Gamma$-stable for~$n$ sufficiently large. Indeed it then follows
from~\eqref{un minus u n minus 1} that $u_n-u_{n-1}\in\gM_n$ for ~$n$
sufficiently large, uniformly in~$z\in\gM$.

It is therefore enough to show that
there exists a sequence $m(n)$, tending to~$\infty$ with~$n$, 
such that if we put $\fM_n \coloneqq   T^{m(n)} \fM'$ then $(1-g_j^{-1}h_{n,j})(\gM')\subseteq \gM_n$.
In fact, $\fM_n$ is $\Gamma$-stable for all~$n$, since $\fM'$ is $\Gamma$-stable; 
and $\fM_n$ is $\varphi$-stable for all~$n$ large enough, by~\cite[Lemma~5.2.7]{EGstacktheoreticimages}.
Recalling that $g_j^{-1}h_{n,j}\in
\big(\begin{smallmatrix}
  1+p^{a(n)}\Zp & p^{a(n)}\Zp\\ 0&1
\end{smallmatrix}\big)
$ and that $a(n)\to\infty$, we see in turn that it is enough to show
that for each~$m\ge 0$, there exists an~$M\ge 0$ such that if $h\in \big(\begin{smallmatrix}
  1+p^{M}\Zp & p^{M}\Zp\\ 0&1
\end{smallmatrix}\big)$, then $(1-h)(\gM')\subseteq T^m(\gM')$. Writing $h= \big(\begin{smallmatrix}
  1+p^{M}a & p^{M}b\\ 0&1
\end{smallmatrix}\big)$, we have
\[(1-h)(z)=(1-(1+T)^{p^Mb})z+(1+T)^{p^Mb}(z-\sigma_{1+p^Ma}(z)).\]Certainly
$(1-(1+T)^{p^Mb})\in T^m\A^+_A$ for~$M$ sufficiently large uniformly
in~$b$ (because
$A$ is an $\cO/\varpi^a$-algebra), so it suffices to show that
$z-\sigma_{1+p^Ma}(z)\in T^m\gM'$ for~$M$ sufficiently large, for
all~$z\in \gM'$ and~$a\in\Zp$. This follows from~\cite[Lem.\
D.28~(3)]{emertongeepicture}, since~$A$ is an $\cO/\varpi^a$-algebra.

Finally, we note that by~\eqref{un minus u n minus 1} we have for
any~$k\ge 1$ and~$z\in\gM$  \[
  \Res_{i+p^n\Zp}(u_{n+k}-u_{n+k-1})=\sum_{j\in I_{n+k}(\Zp),f(j)\in i+p^n\Zp}
  \begin{pmatrix}
    1&f(j)\\0&1
  \end{pmatrix}\varphi^{n+k}((g_j-h_{n+k,j})\cdot
  r_{n+k,j}(z)), \]which is contained in~$\gM_{n+k}$ for~$n$
sufficiently large (uniformly in $z\in\gM$ and~$i\in\Zp$). Thus we
have in particular that $\Res_{i+p^n\Zp}(f_*(z)-u_n)\in\gM_n$, and the
result follows.
\end{proof}

The operators~$m_\alpha$ and~$f_*$ enjoy the following properties.

\begin{prop}
  \label{prop: properties of analytic operations}Let $U,V$ be compact
  open subsets of~$\Zp$. %
  \begin{enumerate}
  \item For all locally constant maps $\alpha_1,\alpha_2:U\to A$, we have
    $m_{\alpha_1}\circ m_{\alpha_2}=m_{\alpha_1\alpha_2}$.
  \item If $f:U\to V$ is a local diffeomorphism and $\alpha:V\to A$ is
    locally constant, then \[f_*\circ m_{\alpha\circ f}=m_\alpha\circ f_*.\]
  \item If $f:U\to V$ and $g:V\to W$ are local diffeomorphisms, then
    $(g\circ f)_*=g_*\circ f_*$.
  \item If $\alpha:U\to A$ is locally constant and $V\subseteq U$,
    then~$m_\alpha$ commutes with~$\Res_V$.
  \item If $\alpha:\Zptimes\to A$ is constant on $a+p^n\Zp$ for all
    $a\in\Zptimes$, then \[m_\alpha=\sum_{i\in(\Z/p^n\Z)^\times}\alpha(i)\Res_{i+p^n\Zp}.\]
  \end{enumerate}

\end{prop}
\begin{proof} All but the second and third properties were already established
  above.  In the case that~$A$ is a finite $\cO/\varpi^a$-algebra, the 
  second and third results are
  respectively ~\cite[Prop.\  V.2.4, V.1.6]{ColmezMirabolique}, %
  the proofs %
  of which go over
  unchanged in our setting. %
\end{proof}%

In summary, the sheaf $U \mapsto D\boxtimes U$ is equivariant for
the action of local diffeomorphisms.

\subsubsection{Defining \texorpdfstring{$\text{--}\boxtimes \Q_p$}{box Qp} more generally }
For any~$\psi$-stable
$A$-submodule ~$M\subset D$, we let $M\boxtimes\Zp=M$, and we let $M\boxtimes\Qp$ denote the set of
sequences~$(x_n)_{n\ge 0}$ with $x_n\in M$ and $\psi(x_{n+1})=x_n$ for
all ~$n\ge 0$. Equivalently, we have
$M\boxtimes\Qp\coloneqq \varprojlim_\psi M$. %
We endow $M\boxtimes\Qp$ with its usual projective limit topology,
where $M\subseteq D$ has the subspace topology, and
then~$M\boxtimes\Qp$ is closed in the infinite product $\prod_{n\ge
  0}M$, and has the subspace topology (with the product having the
product topology). 
To see that it is closed, one can for example note that it is
the intersection of the sets of sequences~$(x_n)$
with~$\psi(x_{n+1})=x_n$ for~$n\le N$, which are closed for any~$N$
because~$\psi$ is continuous.
Note that $M \boxtimes \bQ_p$ need not be an $A[P]$-submodule of~$D \boxtimes \bQ_p$, e.g.\ if~$M$ is not an $\bA_A^+$-submodule of~$D$.
However, this will be the case when $M = D^\natural, D^\sharp$: see Lemma~\ref{lem: D sharp mod D natural box Qp}.

\begin{remark}
As Colmez explains in~\cite[Section~V.1.2]{ColmezMirabolique},
Proposition~\ref{prop: direct image of local diffeomorphism} is motivated by considering
the case of $\A^+_A \coloneqq  A\llbracket T\rrbracket $ with its $(\varphi,\Gamma)$-structure.
Indeed, we may interpret $\A^+_A$ as the $A$-module of $A$-valued measures
on~$\Z_p$, and 
Proposition~\ref{prop: direct image of local diffeomorphism} in the case of a local diffeomorphism $f: \bZ_p \to \bZ_p$
computes the pushforward of measures under~$f$.
Here we should note that $\A^+_A$ is also $\psi$-invariant, because $\A_A^+ = \A_A^{\natural}$,  
and so
the preceding discussion can be applied to it: indeed, by~\cite[Corollary~3.30]{MR3444236},
an $A[P^+]$-submodule of an \'etale $(\varphi, \Gamma)$-module~$D$ with $A$-coefficients is \'etale in the sense of Definition~\ref{defn:etale-Pplus-module} 
if and only if it is $\psi$-stable.

A general \'etale $(\varphi,\Gamma)$-module
$D$ need not 
contain an $\A_A^+$-lattice that is simultaneously $\varphi$ and $\psi$-invariant ---
it contains such a lattice if and only if $D^+ = D^{\natural}$ --- and so in general it is only 
$D$ that localizes, rather than any of its lattices.
Nevertheless, as
we will see
in what follows, the lattices $D^+$, $D^{\natural}$, and $D^{\sharp}$
do interact in important ways with the localization of~$D$.   
\end{remark}%

\subsubsection{Key examples: \texorpdfstring{$D^{\natural}\boxtimes \Q_p$}{D natural box Qp} and \texorpdfstring{$D^{\sharp}\boxtimes \Q_p$}{D sharp box Qp}}
Since each of $D^{\natural}$ and $D^{\sharp}$ is $\psi$-stable, both
$D^{\natural}\boxtimes \Q_p$ and $D^{\sharp}\boxtimes \Q_p$ are defined.
Furthermore, the inclusion $D^{\natural} \subseteq D^{\sharp}$ evidently induces an embedding
\begin{equation}
\label{eqn:embedding D natural box Qp into D sharp box Qp}
D^{\natural}\boxtimes \Q_p \hookrightarrow D^{\sharp}\boxtimes \Q_p.
\end{equation}

\begin{lem}\label{lem: D sharp mod D natural box Qp}
\leavevmode
\begin{enumerate}
\item The embedding~{\em \eqref{eqn:embedding D natural box Qp into D sharp box Qp}}
is open.
\item 
  We have natural identifications of finite $A$-modules
  \[(D^\sharp\boxtimes\Qp)/(D^\natural\boxtimes\Qp)=(D^\sharp/D^\natural)\boxtimes\Qp=D^\sharp/D^\natural.\]
\item $D^\natural \boxtimes \bQ_p$ and~$D^\sharp \boxtimes \bQ_p$ are $B(\bQ_p)$-stable $A$-submodules of~$D \boxtimes \bQ_p$, and the embedding~\eqref{eqn:embedding D natural box Qp into D sharp box Qp}
is $A[B(\bQ_p)]$-equivariant.
\end{enumerate}
\end{lem}
\begin{proof}We have a short exact sequence
  \[0\to D^\natural\to D^\sharp\to D^\sharp/D^\natural\to 0.\]
  Since~$\psi$ acts surjectively on each term, %
  the passage
  to~$\varprojlim_\psi$ is exact, and we have a short exact
  sequence
  \[0\to D^\natural\boxtimes\Qp\to D^\sharp\boxtimes\Qp\to
    (D^\sharp/D^\natural)\boxtimes\Qp\to 0.\] Since~$D^\sharp/D^\natural$ %
  is a finite $A$-module, its surjective endomorphism~$\psi$ is bijective, so we may identify
  $(D^\sharp/D^\natural)\boxtimes\Qp$ with $D^\sharp/D^\natural$.  This proves~(2),
  and also shows that~\eqref{eqn:embedding D natural box Qp into D sharp box Qp}
identifies $D^{\natural}\boxtimes \Q_p$ with the kernel of the composite
$$\xymatrix{D^{\sharp}\boxtimes \Q_p \ar^-{\text{proj.\ to 1st factor}}[rr] &&
  D^{\sharp} \ar[r]&
D^{\sharp}/D^{\natural}},$$
proving~(1).
Finally, part~(3) is immediate from the definitions. %
\end{proof}%

\subsection{\texorpdfstring{The $\GL_2(\Q_p)$-representations $D\boxdelta\mathbb{P}^1$}{D box P1}}
\label{subsec: D box P1}
Continue to let~$A$ denote a Noetherian $\cO/\varpi^a$-algebra, and let ~$D$ be
an \'etale $(\varphi,\Gamma)$-module with
$A$-coefficients.
In the preceding subsection we have explained how to localize $D$ to a sheaf
over~$\Q_p$, equivariantly with respect to the action of local diffeomorphisms.
By general principles, this structure then allows 
us to localize $D$ over other $1$-dimensional locally analytic $p$-adic
manifolds, such as $\Pone(\Q_p)$.  Since $\PGL_2(\Q_p)$ acts
on $\Pone(\Q_p)$ via diffeomorphisms, the resulting sheaf on $\Pone(\Q_p)$
will be $\PGL_2(\Q_p)$-equivariant,
and in particular its global sections will afford a $\PGL_2(\Q_p)$-representation.

In fact, we follow Colmez in making a slightly more involved construction
(via explicit formulas)
which depends on our choice of character~$\zeta$, 
and leads to a $G$-representation with central character~$\zeta \varepsilon^{-2}$. 
In the case when $\zeta = \varepsilon^2$
we obtain the construction of the preceding paragraph.

Let~$w \coloneqq
\big(\begin{smallmatrix}
  0&1\\1&0
\end{smallmatrix}\big)
$. Then~$w$ acts on~$\Zptimes$ via the local
diffeomorphism~$z\mapsto 1/z$, and we let
$H_w:D\boxtimes\Zptimes\to D\boxtimes\Zptimes$ be the composite
$w_*\circ m_{\zeta\varepsilon^{-2}}$ (where we think of~
$\zeta\varepsilon^{-2}$ as the function $\zeta\varepsilon^{-2}:\Zp^{\times}\to\cO^{\times}\to A^{\times}\subset A$). 
We set \[D\boxdelta\Pone \coloneqq \{z=(z_1,z_2)\in
  D\times D,\
  \Res_{\Zptimes}(z_2)=H_w(\Res_{\Zptimes}(z_1))\}.\] This is a closed
subspace of~$D\times D$. %

If $g=
\big(\begin{smallmatrix}
  a&b\\c&d
\end{smallmatrix}\big)
\in G$, and~$U$ is a compact open subset of~$\Qp$ not
containing $-d/c$, then we write~$\alpha_g:U\to A$ for the function
$\alpha_g(x)=(\zeta\varepsilon^{-2})(cx+d)$. The action of~$g$ on~$\Pone(\Qp)$ via M\"obius
transformations induces a diffeomorphism $g:U\to gU$, and we define
$H_g:D\boxtimes U\to D\boxtimes gU$ to be $g_*\circ
m_{\alpha_g}$. 
Note that this specialises to the definition of $H_w$
above.

We have a restriction map $\Res_U: D\boxdelta\Pone\to D\boxtimes U$
defined by \[\Res_U(z)=\Res_{U\cap\Zp}(z_1)+H_w(\Res_{wU\cap
    p\Zp}(z_2))=\Res_{U\cap p\Zp}(z_1)+H_w(\Res_{wU\cap \Zp}(z_2)).\]

\begin{thm}
  \label{thm: existence of G action on D box P1}There is a unique
  action $(g,z)\mapsto g\cdot z$ of $\GL_2(\Qp)$ on $D\boxdelta\Pone$
  such that for all open compact subsets~$U$ of~$\Qp$, and all
  elements $z=(z_1,z_2)$ of~$D\boxdelta\Pone$, we have
  \begin{equation}\label{formula for G action}
  \Res_U(g\cdot z)=H_g(\Res_{g^{-1}U\cap\Zp}(z_1))+H_{gw}(\Res_{(gw)^{-1}U\cap
    p\Zp}(z_2)).
  \end{equation}
\end{thm}%
\begin{proof}
  In the case that~$A$ is a finite $\cO/\varpi^a$-algebra this
  is~\cite[Thm.\ II.1.4]{MR2642409}. The proof only makes use of the
  properties of~$f_*$ and~$m_\alpha$ in Proposition~\ref{prop:
    properties of analytic operations}, and therefore goes over
  immediately in our setting.
\end{proof}

\begin{rem}%
  \label{rem: central char}It is immediate from the definition that %
~$D\boxdelta\Pone$ has central character~$\zeta\varepsilon^{-2}$. %
  Furthermore, comparing our definitions with~\cite[Construction~II.1]{MR2642409}, we see that our $D \boxtimes \Pone$ coincides (in the case that $A$ is a finite $\cO/\varpi^{a}$-algebra) with the representation denoted 
  $D \boxtimes_{\zeta\varepsilon^{-2}} \Pone$ in the notational scheme of \emph{loc.\ cit.}
  \end{rem}

As explained in~\cite[Rem.\ II.1.7]{MR2642409}, it follows from the
proof of Theorem~\ref{thm: existence of G action on D box P1} that if
$z=(z_1,z_2)\in D\boxdelta\Pone$, then $z_1=\Res_{\Zp}z$,
and~$z_2=\Res_{\Zp}(w\cdot z)$.
These relations imply that $w(z_1, z_2) = (z_2, z_1)$ for all~$(z_1, z_2) \in D \boxtimes \Pone$.
Furthermore, \cite[Rem.\ II.1.7~(iii)]{MR2642409} shows that if $U \subset \bQ_p$ is compact open, 
then~$\Res_U$ admits a continuous $A$-linear section~$\iota_U$, defined by
\begin{equation}\label{section to restriction}
\begin{aligned}
\iota_U: D \boxtimes U &\to D \boxtimes \mathbb{P}^1 \\
z &\mapsto \left( \Res_{U \cap \mathbb{Z}_p}(z), H_w(\Res_{U \cap w\mathbb{Z}_p}(z)) \right).
\end{aligned}
\end{equation}
The $A$-module~$D \boxtimes U$ is therefore a direct summand of $D \boxtimes \Pone$, and we will sometimes identify it with this direct summand, and
write $\Res_U$ for~$\iota_U \circ \Res_U$.
Similarly, if $\fM \subseteq D$ is an $A$-submodule, we will write $\fM \boxtimes \bZ_p$ for the image of~$\fM$ under~$\iota_{\bZ_p}$.
We single out a special case of this discussion in the next lemma.

\begin{lem}\label{lem: seeing D box P1 as a Tate module}
 The map
  \begin{equation}\label{eqn: explicit D box Pone}D\boxdelta\Pone\to (D\boxtimes\Zp)\oplus(D\boxtimes p\Zp)\end{equation} given
  by $z\mapsto (\Res_{\Zp}z,\Res_{p\Zp}(w\cdot z))$ \emph{(}or equivalently
  by $(z_1,z_2)\mapsto (z_1,\Res_{p\Zp}(z_2))$\emph{)} is a homeomorphism,
  the inverse morphism being given
  by 
  \[\iota_{\bZ_p} + w \circ \iota_{p\bZ_p}: (y_1,y_2)\mapsto(y_1,y_2+H_w(\Res_{\Zptimes}(y_1))).\]
  The induced action of~$w$ on $(D \boxtimes \bZ_p) \oplus (D \boxtimes p\bZ_p)$ is
  \[
  w(y_1, y_2) = (y_2+H_w(\Res_{\bZ_p^\times}(y_1)), \Res_{p\bZ_p}(y_1)).
  \]
\end{lem}
\begin{proof}
The verification that the maps are mutually inverse is formal, and
their continuity follows from the continuity of~$H_w$,
$\Res_{\Zptimes}$, and $\Res_{p\Zp}$. 
The formula for~$w$ is an immediate consequence of the fact that $w(z_1, z_2) = (z_2, z_1)$ for all $(z_1, z_2) \in D \boxtimes \Pone$.
\end{proof}

\begin{cor}\label{pointwise continuity of G-action}
For all $g \in G$, the map $D \boxtimes \Pone \to D \boxtimes \Pone, z \mapsto g \cdot z$ is continuous.
\end{cor}
\begin{proof}
By Lemma~\ref{lem: seeing D box P1 as a Tate module}, it suffices to prove that $z \mapsto \Res_U(g \cdot z)$ is continuous for all
open compact subsets $U \subseteq \bQ_p$.
This is true by the formula for~$\Res_U(g \cdot z)$ given in Theorem~\ref{thm: existence of G action on D box P1}, together with the fact that the operators
$f_*$ and~$m_\alpha$ are continuous (using Proposition~\ref{prop: direct image of local diffeomorphism} for~$f_*$).
\end{proof}

\begin{rem}
  Since $D\boxtimes p\Zp=(\varphi\circ\psi)(D) = \varphi(D)$ is a direct summand of $D = \varphi(D) \oplus D^{\psi = 0}$, the subspace topology on $D \boxtimes p\bZ_p$ makes it a Tate $A$-module.  The map $\varphi: D \to D \boxtimes p\bZ_p$ is a topological isomorphism with inverse~$\psi$.  Similarly, $D \boxtimes p^{n}\bZ_p$ is a Tate $A$-module for all $n \in \bZ$, and it is topologically isomorphic to~$D$ via the action of $\diag(p^\bZ, 1) \subseteq P$.  The homeomorphism \eqref{eqn: explicit
    D box Pone} then shows that $D\boxdelta\Pone$ is a Tate $A$-module. In particular, it makes sense to talk about $A$-submodules of $D\boxdelta\Pone$ being lattices, and we will do so without further comment.
\end{rem}

\begin{lemma}\label{restriction preserves lattices}
Let~$n \in \bZ$.
Then $\Res_{p^{n}\bZ_p} : D \boxtimes \Pone \to D \boxtimes p^{n} \bZ_p$ is open, and sends lattices to lattices.
The same is true for $\Res_{p^n \bZ_p} : D \to D \boxtimes p^n \bZ_p$ if~$n \geq 0$.
\end{lemma}
\begin{proof}
By construction, the map $\Res_{p^{n}\bZ_p}$ is continuous, surjective and $A$-linear.
We now show that it is open. 
The equivariance properties of~$\Res_{p^n\bZ_p}$ under the action of~$G$ show that it suffices to prove this when~$n = 0$, in which case $\Res_{\bZ_p}$ is the projection onto
a direct summand, by Lemma~\ref{lem: seeing D box P1 as a Tate module}.
Hence it is open.

There remains to prove that $\Res_{p^n \bZ_p}$ sends lattices to lattices.
It suffices to show that if $f: M \to N$ is a surjective, continuous and open morphism between $A$-Tate modules, and $\fM \subset M$ is a
lattice, then $f(\fM)$ is a lattice in~$N$.
Let~$U \subset f(\fM)$ be an open $A$-submodule; we need to prove that $f(\fM)/U$ is a finite $A$-module.
Since~$f$ induces an isomorphism %
\[
\fM/(\fM \cap f^{-1}(U))\cong f(\fM)/U,
\]
this is a consequence of the fact that~$\fM$ is a lattice and $\fM \cap f^{-1}(U)$ is open.
\end{proof}

\begin{lemma}\label{lattices in D box P^1}
Let~$\fM \subset D$ be a lattice. Then: 
\begin{enumerate}
\item $(\fM \boxtimes \bZ_p) + w(\fM \boxtimes \bZ_p)$ is open in $D \boxtimes \Pone$, and 
\item $\{z \in D \boxtimes \Pone: \Res_{\bZ_p}z \in \fM, \Res_{\bZ_p}(wz) \in \fM\}$ is a lattice in $D \boxtimes \Pone$.
\end{enumerate}
\end{lemma}
\begin{proof}
For any $A$-submodule $L \subset D \boxtimes \Pone$, it follows from Lemma~\ref{lem: seeing D box P1 as a Tate module} that $L$ is open 
if and only if the image of~$L$ under
\begin{align*}
    \rho : D \boxtimes \Pone &\to D \oplus \varphi(D),\\ z &\mapsto (z, \Res_{p\Zp}(w\cdot z)),\\ (z_1 ,z_2 ) &\mapsto (z_1, \Res_{p\Zp}(z_2))
\end{align*}
is open.

Let $L \coloneqq (\fM \boxtimes \bZ_p) + w(\fM \boxtimes \bZ_p)$.
Recall from~\eqref{section to restriction} that 
\[
\fM \boxtimes \bZ_p \coloneqq  \iota_{\bZ_p}(\fM) = \{(m, H_w \Res_{\bZ_p^\times} m): m \in \fM\} \subset D \boxtimes \Pone.
\]
Then 
\[
\rho(\fM \boxtimes \bZ_p) = \{(m, 0): m \in \fM\} \subset D \oplus \varphi(D)
\]
and 
\[
\rho(w(\fM \boxtimes \bZ_p)) = \{(H_w\Res_{\bZ_p^\times}m, \Res_{p \bZ_p} m): m \in \fM\} \subset D \oplus \varphi(D).
\]
To see that $\rho(L)$ is open, consider the $A$-linear map
\[
\lambda : D \oplus D \to D \oplus \varphi(D), (x, y) \mapsto (x + H_w \Res_{\bZ_p^\times} y, \Res_{p \bZ_p} y).
\]
Then $\rho(L) = \lambda (\fM \oplus \fM)$, and since $\fM \oplus \fM$ is open in $D \oplus D$, it suffices to prove that~$\lambda$ is open.
To do so, it suffices to find neighborhood bases~$(M_n),(L_n)$ of~$0\in D$ such that $\lambda(M_n \oplus L_n)$ is open for all~$n$.
Choose~$M_n$ arbitrarily, and choose~$L_n$ such that $H_w \Res_{\bZ_p^\times}(L_n) \subseteq M_n$, using the continuity of $H_w \circ \Res_{\bZ_p^\times}$.
Then
\[
\lambda(M_n \oplus L_n) \supseteq M_n \oplus \Res_{p\bZ_p}(L_n).
\]
Since~$\Res_{p\bZ_p}: D \to \varphi(D)$ is open (because it coincides with~$\varphi \psi$), %
this concludes the proof of part~(1).

Now let $L \coloneqq \{z \in D \boxtimes \Pone: \Res_{\bZ_p}z \in \fM, \Res_{\bZ_p}(wz) \in \fM\}$.
Since $\Res_{\bZ_p}$ and~$w$ are continuous (using Corollary~\ref{pointwise continuity of G-action} for~$w$), $L$ is open.
Hence $\rho(L) \subset D \oplus \varphi(D)$ is also open.
On the other hand, $\rho(L)$ is contained in $\fM \oplus \Res_{p \bZ_p}\fM$, which is a lattice in $D \oplus \varphi(D)$, by Lemma~\ref{restriction preserves lattices}, so that $\rho(L)$ is also bounded, and so it is a lattice, by Remark~\ref{rem: open and bounded equals lattice}.
\end{proof}

We now study the continuity of the $G$-action on $D \boxtimes \Pone$.

\begin{lem}%
  \label{lem: G action on G box P1 is continuous}The action
  of~$G$ on~$D\boxdelta\Pone$ that was constructed in
Theorem~{\em \ref{thm: existence of G action on D box P1}} is continuous.
\end{lem}
\begin{proof}

  We claim %
  that it is enough to show that the action map \begin{equation}\label{eqn: H D box P1 continuity}H \times  ( D \boxdelta \Pone)\to
    D\boxdelta\Pone\end{equation} is jointly continuous for some open subgroup~$H$
  of~$G$.

  Indeed, assume that this is the case, and let
  $(g,x)\in G\times (D\boxdelta\Pone)$ be arbitrary. 
  Any
  open neighbourhood of $gx\in D\boxdelta\Pone$ contains a
  neighbourhood of the form $gx+\gM$, where
  $\gM\subset D\boxdelta\Pone$ is a lattice. 
  Since $g\in G$ acts continuously
  on~$D\boxdelta\Pone$,
  $g^{-1}\gM\subset D\boxdelta\Pone$ is open. By the continuity
  of~\eqref{eqn: H D box P1 continuity}, there is an open subgroup~$H'$
  of~$H$ and a lattice $L\subset D\boxdelta\Pone$ such that
  $H'(x+L)\subseteq x+g^{-1}\gM$; so $(gH')(x+L)\subseteq gx+\gM$, and
  $gH'\times (x+L)\subset G\times D\boxdelta\Pone$ is an open
  neighbourhood of~$(g,x)$.

 We now turn to establishing the continuity of~\eqref{eqn: H D box P1
   continuity}. If $H=\cup_iU_i$ and $D\boxdelta\Pone=\cup_j V_j$ are open
  covers of $H$ and $D\boxdelta \Pone$ respectively, then it suffices to check
  that each $U_i \times V_j \to D \boxdelta \Pone$ is jointly continuous.
  We claim that it suffices to show that
  $H \times D \to D \boxdelta \Pone$ is continuous. Indeed, since~$w$ is a continuous automorphism of
  $D \boxdelta \Pone$, if we write 
  \[D \boxdelta \Pone = D \boxtimes \Zp \cup w (D \boxtimes \Zp) = D
    \cup w D,\]  then we deduce that $wHw\times wD\to D\boxdelta\Pone$ is continuous, 
  and so
  the action map 
  \[(H\cap wHw) \times D\boxdelta\Pone\to D\boxdelta\Pone\]
  is
  continuous, and since $H\cap wHw$ is also an open subgroup
  of~$G$, we will be done.

  We now take~$H$ to be the usual Iwahori subgroup of matrices which
  are upper triangular modulo~$p$, so that by the Iwahori
  decomposition we may write \[H = w U (p\Zp) w B(\Zp)\]where $B$ is
  the (upper triangular) Borel subgroup, and $U(p\Zp)\coloneqq 
\big(\begin{smallmatrix}
  1&p\Zp\\0&1
\end{smallmatrix}\big)
$.
In summary, we have reduced to showing that 
    \begin{equation}\label{eqn: first joint continuity equation}w U( p \Zp) w  B(\Zp) \times D \to D \boxdelta \Pone \end{equation}
is jointly continuous.    

Now, the continuity of
\begin{equation}\label{eqn: joint Borel continuity}B(\Zp) \times D \to D = D \boxtimes \Zp \subset D \boxdelta \Pone\end{equation}
is more or less immediate from the definition of the $B(\Zp)$-action
on~$D$ and the continuity assumptions in the definition of a
$(\varphi,\Gamma)$-module. 
More precisely, it follows from the
existence of a neighborhood basis of zero in~$D$ consisting of $(\varphi,\Gamma)$-stable lattices, which in turn
follows from the compactness of~$\Gamma$ (see~\cite[Lemma~5.1.5]{emertongeepicture}
for a proof in the case that~$A$ is a Noetherian $\Fp$-algebra, which goes over
unchanged to the case of Noetherian $\cO/\varpi^a$-algebras).

Accordingly,  we can factor~\eqref{eqn: first joint continuity equation} as 
  \begin{equation}\label{eqn: second joint continuity equation} w U( p
  \Zp) w   B(\Zp) \times D  \to    w U(p\Zp) w \times D  \to
  D\boxdelta \Pone\end{equation} with the first arrow continuous; so we just have to show that the second arrow is continuous.    Since~$w$ is a continuous automorphism, this amounts to showing that
  \[U (p\Zp)   \times   w (D \boxtimes \Zp)  \to  D \boxdelta \Pone\]
is jointly continuous. Writing $D\boxtimes\Zp=D\boxtimes\Zptimes\cup
D\boxtimes p\Zp$, and taking into account the continuity of
~\eqref{eqn: joint Borel continuity}, we see that it is enough to show
that \[U (p\Zp)   \times   w (D \boxtimes p\Zp)  \to  D \boxdelta
  \Pone\]is jointly continuous.

It follows from the definitions (see the proof
of~\cite[Prop.\ II.1.8]{MR2642409}) that~$U(p\Zp)$ preserves
$w(D\boxtimes p\Zp)$, and more precisely if~$b\in p\Zp$, then  $
\big(\begin{smallmatrix}
  1&b\\0&1
\end{smallmatrix}\big)
$ takes~$(0,z)\in w(D\boxtimes p\Zp)\subset D\boxdelta\Pone$ to
$(0,z')$ with %
\[z' = (\zeta\varepsilon^{-2})^{-1}(1+b)
  \begin{pmatrix}
    1&-1\\0&1
  \end{pmatrix} w
  \begin{pmatrix}
    (1+b)^2 & b(1+b)\\0 & 1
  \end{pmatrix}w
  \begin{pmatrix}
    1 & 1/(1+b)\\0& 1
  \end{pmatrix}z.
\]Noting that $\big(\begin{smallmatrix}
    1 & 1/(1+b)\\0& 1
  \end{smallmatrix}\big)$ takes $p\Zp$ to~$1+p\Zp$, and that~$w$ is a continuous
    automorphism of $D\boxtimes(1+p\Zp)$, the result follows from the
    continuity of~\eqref{eqn: joint Borel continuity}.
  \end{proof}

In our later applications, rather than appealing to this continuity directly,
we will use the following consequence of it.

\begin{cor}%
  \label{cor: O-G  action on G box P1}The $G$-action
  on~$D\boxdelta\Pone$ extends to an
  $A\llbracket G\rrbracket $-module structure on $D\boxtimes\Pone$,
  uniquely characterized by the requirement that 
  the induced $A\llbracket K\rrbracket $-action on $D\boxtimes \Pone$ is continuous.
  Furthermore, the endomorphisms of $D\boxdelta\Pone$ induced
  by the elements of~$A\llbracket G\rrbracket $ are all continuous.
\end{cor}
\begin{proof}%
  
  We will show that the action of~$K$ on
  ~$D\boxdelta\Pone$ extends to a continuous action of~$A\llbracket K\rrbracket $.
  Such an extension is unique if it exists (by continuity, together with
  the density of~$A[K]$ in~$A\llbracket K\rrbracket $),
  and by the definition of~$A\llbracket G\rrbracket $, this will yield the required~$A\llbracket G\rrbracket $-action.
  The final claim of the lemma will also then follow, by combining the continuity
  of the $G$-action and the $\cO\llbracket K\rrbracket $-action on~$D\boxdelta \Pone$. 

In order to construct the required~$A\llbracket K\rrbracket $-action, choose a lattice~$L\subset D\boxtimes\Pone$.
  Lemma~\ref{lem: G action on G box P1 is continuous}
  shows that the action map
  $K\times (D\boxtimes \Pone) \to D\boxtimes \Pone$ 
  is jointly continuous, and so we see that for each~$k\in K$, there
  is an open subgroup~$U_k$ of~$K$ and a lattice $L_k\subseteq L$ such that
  $(kU_k)\cdot L_k\subseteq L$. Since~$K$ is compact, it is covered by finitely many
  of the $kU_k$, and if we let~$L'$ denote the intersection of the
  corresponding lattices~$L_{k}$, we see that $K\cdot L'\subseteq L$. 
  The $A$-submodule generated by $K \cdot L'$ is therefore open and bounded, and so it is a lattice.
  Replacing~$K \cdot L'$
  by this lattice, %
  and recalling that~$L$ was
  arbitrary, we see that $D\boxtimes\Pone$ contains a basis of $K$-stable
  lattices~$L_n$. Then each %
  $(D\boxtimes\Pone)/L_n$ is a discrete
  $A$-module with a continuous action of~$K$, hence is a smooth $K$-representation,
  and thus is canonically an $A\llbracket K\rrbracket $-module (with jointly continuous action map). Consequently
   \[D\boxtimes\Pone=\varprojlim_n(D\boxtimes\Pone)/L_n\] inherits 
  a continuous $A\llbracket K\rrbracket $-action, as required.
\end{proof}

We now need to consider the restriction map
$\Res_{\Qp}:D\boxdelta\Pone\to D\boxtimes\Qp$, which by definition is
given by \[z\mapsto\left(\Res_{\Zp}\left(
  \left(\begin{smallmatrix}
    p^n&0\\0&1
  \end{smallmatrix}\right)
z\right)\right)_{n\ge 0}. \]The verification that this defines a morphism
to~$D\boxtimes\Qp$ is formal; see the proof of~\cite[Prop.\
II.1.14]{MR2642409}. %
The action of~$P$ on~$D\boxtimes\Qp$ extends to an action of
the Borel subgroup~$B(\Qp)$ of~$G$ by letting the diagonal matrix~$
\big(\begin{smallmatrix}
  z&0\\0&z
\end{smallmatrix}\big)
$ act via~$(\zeta\varepsilon^{-2})(z)$. %
\begin{lem}
  \label{Res Qp is continuous and P equivariant}$\Res_{\Qp}:D\boxdelta\Pone\to
  D\boxtimes\Qp$ is continuous and $B(\Qp)$-equivariant. %
\end{lem}
\begin{proof}The $B(\Qp)$-equivariance is formal; see the proof
  of~\cite[Prop.\ II.1.14]{MR2642409}. Since $D\boxtimes\Qp$ has the
  inverse limit topology, the continuity of~$\Res_{\Qp}$ follows from
  the continuity of~$\Res_{\Zp}$ and of the action of~$\big(\begin{smallmatrix}
    p^n&0\\0&1
  \end{smallmatrix}\big)$.  
\end{proof}

We then define \[D^\natural\boxdelta\Pone \coloneqq \{z\in D\boxdelta\Pone,\
  \Res_{\Qp}(z)\in D^\natural\boxtimes\Qp\},\]
\[D^\sharp\boxdelta\Pone \coloneqq \{z\in D\boxdelta\Pone,\ \Res_{\Qp}(z)\in
  D^\sharp\boxtimes\Qp\}.\] 
By Lemma~\ref{Res Qp is continuous and P equivariant}, $D^\natural\boxdelta\Pone$ and
$D^\sharp\boxdelta\Pone$ are closed $A$-submodules of
$D\boxdelta\Pone$. In the case that~$D$ has rank two and determinant~$\zeta \varepsilon^{-1}$,
we will show that they are in fact lattices, but
this is not obvious; indeed our proof will be closely intertwined with
the proof that $D^\natural\boxdelta\Pone$ is $G$-stable, and
as in the original work of Colmez, this is ultimately established via
$p$-adic interpolation from the crystalline %
case.
We end this section by establishing two properties of these modules.
\begin{lem}
  \label{D sharp D natural box Pone agree up to
    finite}$(D^\sharp\boxdelta\Pone)/(D^\natural\boxdelta\Pone)$ is a
  finite $A$-module, and $D^\natural \boxtimes \Pone$ is an open submodule of $D^\sharp \boxtimes \Pone$.
\end{lem}
\begin{proof}
  By definition, $\Res_{\Qp}$ induces an
  injection
  \[(D^\sharp\boxdelta\Pone)/(D^\natural\boxdelta\Pone)\into
    (D^\sharp\boxtimes\Qp)/(D^\natural\boxtimes\Qp),\] and the target
  is a finite $A$-module by Lemma~\ref{lem: D sharp mod D natural box
    Qp}.
    This proves the first statement, and since~$\Res_{\bQ_p}$ is continuous, 
    the second statement also follows from Lemma~\ref{lem: D sharp mod D natural box Qp}, which asserts that $D^\natural \boxtimes \bQ_p$ is open
    in $D^\sharp \boxtimes \bQ_p$.
\end{proof}

\begin{lem}\label{lem: Dplus box Zp contained in Dnatural box P1}
$D^\natural \boxtimes \Pone$ contains $D^+ \boxtimes \bZ_p$.
\end{lem}
\begin{proof}
Let~$t \coloneqq \big(\begin{smallmatrix} p & 0\\0 & 1\end{smallmatrix}\big)$.
By definition, we need to prove that if~$z \in D^+ \boxtimes \bZ_p$ and~$n \geq 0$ then $\Res_{\bZ_p}(t^n z) \in D^\natural$.
Recall that $D^+ \boxtimes \bZ_p \coloneqq  \iota_{\bZ_p}(D^+)$. 
Applying~\eqref{section to restriction}, we see that there exists~$z_0 \in D^+$ such that 
$z = (z_0, H_w \Res_{\bZ_p^\times} z_0)$ as elements of $D \boxtimes \Pone$.
Since~$D^+$ is $\varphi$-stable, and $D^+ \subset D^\natural$ by Lemma~\ref{lem: Dplus contained in Dnatural}, 
it thus suffices to prove that $\Res_{\bZ_p}(t^n z)= \varphi^n(z_0)$ for all~$n \geq 0$.

Applying~\eqref{formula for G action} we see that
\[
\Res_{\bZ_p}(t^n z) = H_{t^n}(z_0):
\]
indeed, the second summand in~\eqref{formula for G action} vanishes, because $\Res_{p\bZ_p}(H_w \Res_{\bZ_p^\times} z_0) = 0$.
We now observe that $H_{t^n} = t^n_*$, where $t^n : \bZ_p \to p^n \bZ_p$ is the diffeomorphism of multiplication by~$p^n$:
indeed, by definition 
\[
H_{t^n} = t^n_* \circ m_{\alpha_{t^n}},
\] 
and the function $\alpha_{t^n}$
is identically~$1$ on~$\bZ_p$.
Now Example~\ref{ex:undoing restriction to pZp} shows that
\[
t^n_* = \varphi^n : D \to D \boxtimes p^n \bZ_p = \varphi^n(D),
\]
which concludes the proof.
\end{proof}

\subsection{The action of \texorpdfstring{$\GL_2(\Qp)$}{GL2} on lattices in \texorpdfstring{$D\boxdelta\Pone$}{D box P1}}\label{subsec: D natural box P1 GL2
  stable}%
We continue to assume that $A$ is a Noetherian $\cO/\varpi^a$-algebra.
Our goal in this section is to prove Proposition~\ref{prop: G-stable
  lattice implies D natural box P1 stable}, which gives a criterion
for $D^\natural\boxdelta\Pone$ to be
$G$-stable, or (equivalently) $\cO\llbracket G\rrbracket $-stable (with respect to
the $\cO\llbracket G\rrbracket $-action on $D\boxdelta \Pone$ given by Corollary~\ref{cor: O-G  action on G box P1}).
The study of this question is intertwined with the study of certain
lattices in~$D\boxdelta\Pone$.

We embed $\Q_p^{\times}$ into $G$ via
\begin{equation}
\label{eqn:Qp-times embedding into G}
 \Q_p^{\times} \iso \begin{pmatrix} \Q_p^{\times} & 0 \\ 0 & 1 \end{pmatrix};
\end{equation}
in particular, the composite of this embedding with $\det:G \to \Q_p^{\times}$
is the identity on~$\Q_p^{\times}$. 

\begin{rem}\label{rem:Dnr-is-Qptimes-stable}
  Note that since~ $D^{\nr}\subseteq D$ is $\varphi$- and ~$\Gamma$-stable, $D^{\nr}\subseteq D=D\boxtimes\Zp\subset D\boxtimes\Qp\subset D\boxtimes\Pone$ is $\Qptimes$-stable.
\end{rem}
\begin{lem}
  \label{lem: kernel of Res Qp}
  The kernel of $\ResQp:D\boxtimes\Pone\to
  D\boxtimes\Qp$ is $w\cdot D^{\nr}$. In
  particular, $\ker\ResQp$ is a finite $A$-module. %
\end{lem}
\begin{proof} Since $D^{\nr}$ is a finite $A$-module by
  Lemma~\ref{lem: Dnr is finite}, it is enough to prove that
  $\ker\ResQp=w\cdot D^{\nr}$.  This follows formally from the
  definitions, exactly as in the proof of~\cite[Prop.\
  II.1.14]{MR2642409}. %
\end{proof}

We have the following variant on ~\cite[Lem.\
  III.3.6]{ColmezMirabolique}, which is proved in the same way.
\begin{lem}
  \label{lem: minimal stable lattice in D box Qp-prelim-version}If
  $M\subseteq D\boxdelta\Qp$ is a closed $B(\Qp)$-stable $A$-submodule, then
  \[M=M_0 \boxtimes\Qp,\] where $M_0\coloneqq \Res_{\Zp}(M)\subseteq D$ is
  a $\Gamma$- and $\psi$-stable $\A_A^+$-submodule on which $\psi$ acts surjectively. 
\end{lem}
\begin{proof}%
  Since~$M$ is $B(\Qp)$-stable, it follows from the definition of the
  action of~$P$ on $D\boxtimes\Qp$ 
  that~$M_0 $ is $\A_A^+$-, $\Gamma$-, and $\psi$-stable,
  and that $\psi$ acts surjectively on~$M_0$.  %
It also follows that, if we write $D\boxtimes\Q_p \coloneqq  \varprojlim_{\psi} D$,
then the image of $M$ under the $k$th projection (for any $k\ge 0$) is equal to~$M_0$.
We therefore have
  $M\subseteq M_0\boxtimes\Qp$, and it remains to show
  that this inclusion is an equality. 
  Thus, let~$z=(z_n)_{n\ge 0}\in M_0\boxtimes\Qp$.
  The just-established equality of images shows that for 
  each~$k\ge 0$ we can find~$u_k=(u_{k,n})_{n\ge 0}\in M$
  with $u_{k,k}=z_k$. By the definition of the 
  projective limit
  topology on $D\boxtimes\Qp$, we see that $u_k\to z$ as $k\to\infty$,
  and we are done, since $M$ is closed
  in~$D\boxtimes\Qp$ by hypothesis.   
\end{proof}%

\begin{lem}
  \label{lem: minimal stable lattice in D box Qp}If
  $\gM\subset D\boxdelta\Pone$ is a $B(\Qp)$-stable lattice, then
  \[D^\natural\boxtimes\Qp\subseteq\ResQp(\gM)\subseteq D^\sharp\boxtimes\Qp.\]%
\end{lem}
\begin{proof}
Since each
  $\Res_{p^{-n}\Zp} : D \boxtimes \Pone \to D \boxtimes p^{-n}\bZ_p$ %
  takes lattices to lattices, by Lemma~\ref{restriction preserves lattices},
we see that
$$\ResQp(\gM) \coloneqq  \varprojlim_n \Res_{p^{-n}\Z_p} \gM \hookrightarrow D\boxtimes\Q_p \coloneqq 
\varprojlim_n D\boxtimes p^{-n}\Z_p,$$
is an inverse limit of open, and hence closed,
$A$-submodules of the various $D\boxtimes p^{-n}\Z_p$, 
and so is a closed
  $A$-submodule of~$D\boxtimes\Qp$.
Now Lemma~\ref{lem: minimal stable lattice in D box Qp-prelim-version}
shows that $\ResQp(\gM) \iso \gM_0 \boxtimes \Q_p,$
where $\gM_0 \coloneqq  \Res_{\Z_p} \gM$ is an $\A_A^+$-submodule of $D\boxtimes \Z_p = D$
which is $\psi$-stable
and on which $\psi$ acts surjectively. 
As already noted,
$\gM_0$ is also a lattice in~$D$. 
Thus $D^\natural\subseteq \gM_0\subseteq D^\sharp$,
by Proposition~\ref{prop: D sharp natural exist}, and the lemma follows. 
\end{proof}

\begin{lem}\label{lem: B stable implies G stable}If $\gM\subset D\boxdelta\Pone$ is a $B(\Qp)$-stable lattice, then
  it is
  $G$-stable, and hence in fact $\cO\llbracket G\rrbracket $-stable.  %
\end{lem}
\begin{proof}Since the group~$G$ has a neighbourhood basis of
  the identity given by open subgroups (e.g.\ the congruence subgroups
  of~$\GL_2(\Zp)$), it follows from Lemma~\ref{lem: G action on G box P1 is continuous} and \cite[Lem.\
  D.13]{emertongeepicture} that there is an open subgroup
  $H\subseteq G$ such that $H\gM=\gM$. Since~$B(\Qp)\gM=\gM$ we
  have $B(\Qp)HB(\Qp)\gM=\gM$, and since we have $B(\Qp)HB(\Qp)=G$ by the Bruhat decomposition, we conclude that~$\fM$ is $G$-stable.

  To see the final assertion, it suffices to note that since $\cO[K]$ is dense in
  $\cO\llbracket K\rrbracket $, and since $\gM$ is closed in $D\boxdelta\Pone$, the $K$-stability of
  $\gM$ implies the $\cO\llbracket K\rrbracket $-stability of~$\gM$, and $G$-stability and $\cO\llbracket K\rrbracket $-stability
  together imply~$\cO\llbracket G\rrbracket $-stability. 
\end{proof}

\begin{prop}
  \label{prop: G-stable lattice implies D natural box P1
    stable}Suppose that there is a $B(\Qp)$-stable lattice
  $\gM\subset D\boxdelta\Pone$. Then $D^\natural\boxdelta\Pone$ is an 
  $\cO\llbracket G\rrbracket $-stable lattice in $D\boxdelta\Pone$.
\end{prop}
\begin{proof} %
  By Lemma~\ref{lem: minimal stable lattice in D box Qp}, we
  have
  \[D^\natural\boxdelta\Pone\subseteq\gM+\ker\ResQp\subseteq
    D^\sharp\boxdelta\Pone.\]
    Since~$\ker\ResQp$ is a finite $A$-module, by Lemma~\ref{lem: kernel of Res Qp}, it follows that $\fM + \ker \Res_{\bQ_p}$ is a lattice.
    By Lemma~\ref{D sharp D natural box Pone agree up to finite},
  $(D^\sharp\boxdelta\Pone)/(D^\natural\boxdelta\Pone)$ is a
  finite $A$-module, hence $D^\sharp \boxtimes \Pone$ is a lattice,
  and $D^\natural \boxtimes \bP^1$ is open in $D^\sharp \boxtimes \Pone$, hence
  $D^\natural\boxdelta\Pone$ is a lattice. 
  Since
  $D^\natural\boxdelta\Pone$ is $B(\Qp)$-stable, the proposition follows
  from Lemma~\ref{lem: B stable implies G stable}.
\end{proof}

We also have the following ``converse'' to Proposition~\ref{prop: G-stable lattice implies D natural box P1
  stable}.

\begin{lem}
  \label{lem: G stable implies lattice}If $D^\natural\boxdelta\Pone$
  is $G$-stable, then it is a lattice in $D\boxdelta\Pone$.
\end{lem}
\begin{proof}
  Under the assumption that $D^\natural\boxdelta\Pone$ is
  $G$-stable, we have
  \[ D^\natural\boxdelta\Pone
    \subseteq \{x\in D\boxdelta\Pone\mid \Res_{\Zp}x\in D^\natural
    \textrm{ and } \Res_{\Zp}wx\in D^\natural\}.\] 
 Since~$D^\natural$ is a lattice in~$D$, it follows
    from
    Lemma~\ref{lattices in D box P^1}~(2) that 
    $D^\natural\boxdelta\Pone$ is a closed
    $A$-submodule of a lattice.
By Remark~\ref{rem: open and bounded equals lattice},
it thus suffices to show that
  $D^\natural\boxdelta\Pone$ is open. 
  To this end, note that it follows from 
  Lemma~\ref{lem: Dplus box Zp contained in Dnatural box P1}  
  that $D^\natural\boxdelta\Pone$
  contains 
  $D^+\boxtimes\Zp$, and since
  $D^\natural\boxdelta\Pone$ is assumed to be
  $G$-stable, it also contains $w\cdot (D^+\boxtimes\Zp)$. It
  therefore contains $D^+\boxtimes\Zp+w\cdot
  (D^+\boxtimes\Zp)$, which by Lemma~\ref{lattices in D box P^1}~(1) is 
  an open neighborhood of the identity, as required. %
\end{proof}

\subsection{Formal \texorpdfstring{$(\varphi,\Gamma)$}{(ϕ,Γ)}-modules}
\label{subsec:formal phi Gamma modules}
In this subsection, which is an interlude in our main discussion, 
we briefly extend some of the previous theory to the context of what we
call {\em formal} \'etale $(\varphi,\Gamma)$-modules. %

Let~$R$ be a complete local Noetherian $\cO$-algebra with finite residue
field. Following Dee~\cite{MR1805474} we let~$\widehat{\A}_{R}$ denote the
$\m_R$-adic completion of~${\A}_{R}$, and we define a \emph{formal \'etale
  $(\varphi,\Gamma)$-module} with $R$-coefficients to be an \'etale
$(\varphi,\Gamma)$-module over $\widehat{\A}_{R}$; that is, a finitely
generated $\widehat{\A}_{R}$-module equipped with commuting semilinear actions
of~$\varphi$ and~$\Gamma$, with the underlying $\varphi$-module being \'etale in
the usual sense. 
We do not demand that a formal \'etale
$(\varphi,\Gamma)$-module is projective; if it is, we explicitly refer to it
as a {\em projective} formal \etale $(\varphi,\Gamma)$-module; note that
it is then in fact {\em free}, since $\widehat{\A}_R$ is a local ring.
Note that if~$R$ is furthermore a
finite $\cO$-module, then $\widehat{\A}_{R}=\A_R$, so that the categories of (projective) formal
\'etale $(\varphi,\Gamma)$-modules with $R$-coefficients and \'etale
$(\varphi,\Gamma)$-modules with $R$-coefficients are equivalent. %

Of course, in the case that~$R$ is Artinian this agrees with the usual
definition of a (not necessarily projective) \'etale $(\varphi,\Gamma)$-module,
and as we recalled in Section~\ref{subsec: the coefficient rings}, in this case
there is an equivalence of categories between the category of \'etale
$(\varphi,\Gamma)$-modules with $R$-coefficients and the category of finitely
generated $R$-modules equipped with a continuous action of~$G_{\Qp}$. We write $V\mapsto\bD(V)$ for
the functor taking a representation of~$G_{\Qp}$ to the corresponding
$(\varphi,\Gamma)$-module, and $D\mapsto \mathbf{V}(D)$ for the quasi-inverse functor.

Returning to the general case that $R$ is a complete
local Noetherian $\cO$-algebra with finite residue field, it is shown in
~\cite{MR1805474} that this equivalence extends to an exact equivalence of categories
between the category of formal \'etale $(\varphi,\Gamma)$-modules with
$R$-coefficients and the category of finitely generated $R$-modules with a
continuous action of~$G_{\Qp}$, where
now \[\bD(V)\coloneqq \varprojlim_k\bD(V/\m_R^{k}V)\]and similarly for~$\mathbf{V}(D)$. (It is also shown in
~\cite{MR1805474} that this construction is compatible with extension of
scalars, so we do not record~$R$ in the notation~$\bD$.)

If~$D$ is a formal \'etale $(\varphi,\Gamma)$-module with $R$-coefficients
then we write~$D_{k}\coloneqq D/\m_R^{k} D$. 
Then $D_k$ is a not-necessarily-projective \'etale $(\varphi,\Gamma)$-module
over $R/\m_R^k$, which however is projective if $D$ is,
and $D=\varprojlim_kD_k$. We extend many
of the definitions that we have made in earlier parts of Section~\ref{sec: Colmez} to the case of formal
\'etale $(\varphi,\Gamma)$-modules by passage to projective limits.
Strictly speaking,
in case $D$ is not projective, so that the $D_k$ need not be projective,
the various definitions and constructions that
 we have made will not apply to the $D_k$;
however, since $R/\m_R^k$ is a finite $\cO$-algebra,
$D_k$ is of finite length over $\A$, and we may ignore the coefficients
and directly apply the definitions of~\cite{ColmezMirabolique}
and~\cite{MR2642409}.
To be precise,
we make the following definition.

\begin{defn}\label{def: various functors on formal phi Gamma modules}
Let~$R$ be a complete Noetherian local $\cO$-algebra with finite residue field, and let~$D$ be a formal \'etale $(\varphi, \Gamma)$-module with $R$-coefficients.
We set:
\begin{itemize}%
\item $D^\sharp \coloneqq  \varprojlim_k D_k^\sharp$ and
  $D^\natural \coloneqq  \varprojlim_k D_k^\natural$. These are $\psi$-stable submodules
  of~$D$.
\item %
  $D \boxtimes \bP^1 \coloneqq  \varprojlim_k D_k \boxtimes \bP^1$, with its natural action of~$\cO\llbracket G\rrbracket $. %

\item $D^\sharp \boxtimes \bQ_p \coloneqq  \varprojlim_\psi D^\sharp$ with its natural
  action of~$B$. %
Hence
  $D^\sharp \boxtimes \bQ_p \cong \varprojlim_k D_k^\sharp \boxtimes \bQ_p$.  We
  define $D^\natural \boxtimes \bQ_p$ and $D \boxtimes \bQ_p$ in a similar way as inverse limits over~$k$,
  so there is a $B$-equivariant map %
  \[
    \Res_{\bQ_p} : D \boxtimes \bP^1 \to D \boxtimes \bQ_p
  \]
  defined as the inverse limit of the restriction maps for~$D_k$.
\item %
  $D^\sharp \boxtimes \bP^1 \coloneqq \varprojlim_k D_k^\sharp \boxtimes \bP^1$, and $D^\natural \boxtimes \bP^1 \coloneqq \varprojlim_k D_k^\natural \boxtimes \bP^1$.
\item $D^{\nr}\coloneqq \varprojlim_kD_k^{\nr}$.
\end{itemize}
\end{defn}

Each of the objects appearing in these various inverse limits has a natural topology,
and so each of the objects 
being defined has a natural inverse limit topology.
In fact, except for $D \boxtimes \Pone$, each of the objects appearing in the various inverse systems in these definitions
is profinite, and so these inverse limit topologies are also profinite.
Recall also that the formation of inverse limits of profinite $R$-modules 
is exact; we will use this fact frequently below.

\begin{lem}
  \label{lem:formal-etale-phi-gamma-behave-well}Let~$R$ be a complete local
  Noetherian $\cO$-algebra with finite residue field, and let~$D$ be a 
  formal \'etale $(\varphi,\Gamma)$-module with $R$-coefficients.
  Then:
  \begin{enumerate}
  \item
    $(D^{\sharp}\boxtimes\Qp)/(D^{\natural}\boxtimes\Qp)=(D^{\sharp}/D^{\natural})\boxtimes\Qp
    =D^{\sharp}/D^{\natural}=\varprojlim_kD_k^\sharp
    / D_k^\natural$.
  \item $D^{\nr}=\bigcap_n\varphi^n(D)$.
  \end{enumerate}
\end{lem}
\begin{proof}
  The first claim follows from Lemma~\ref{lem: D sharp mod D natural box Qp} and the
exactness of the formation of inverse limits of profinite modules.
Noting that since~$\varphi$ is injective we may
  identify $\bigcap_n\varphi^n(D)$ with $\varprojlim_{\varphi}D$, the second claim
  is clear, as limits commute with limits.
\end{proof}

We next recall the description of 
$D^{\nr}$ and $D^{\sharp}/D^{\natural}$ in the context of the 
comparison between formal \'etale $(\varphi,\Gamma)$-modules and 
Galois representations.
Using this, we will see in Lemma~\ref{finiteness in formal case} below
that $D^{\nr}$ and $D^{\sharp}/D^{\natural}$ are both finite $R$-modules.  

\begin{defn}\label{defn:category-of-abelian-representations}
  We equip every finitely generated $R$-module with its canonical (i.e.\ $\fm_R$-adic) topology.
  We write~$\Rep_R(G_{\Qp})$ for the category of continuous representations of~$G_{\Qp}$ on finitely
generated $R$-modules, and ~$\Rep_R^{\ab}(G_{\Qp})$ for the subcategory of abelian representations, i.e.\ representations for which the action of~$G_{\Qp}$ factors through~$G_{\Qp}^{\ab}$.
We also write $\Rep_R(\Qptimes)$ for the category of continuous representations of~$\Qptimes$ on finitely generated $R$-modules.
\end{defn}

We recall the following result, which is simply a reformulation of local class field 
theory for~$\Q_p$ in our context.

\begin{lemma}
\label{lem:local Artin}
Pullback along the local Artin map $\Q_p^{\times} \to G_{\Q_p}^{\ab}$ induces 
an equivalence of categories
$\Rep_R^{\ab}(G_{\Q_p}) \iso \Rep_R(\Qptimes)$.
\end{lemma}
\begin{proof}
The local Artin map identifies $G_{\Q_p}^{\ab}$
with the profinite completion of $\Q_p^{\times}$.
The lemma follows immediately from this, together with the fact that
any continuous action of $\Q_p^{\times}$ on a finitely generated $R$-module
extends to its profinite completion (since $R$ is profinite).
\end{proof}

\begin{defn}
  \label{defn:upper-and-lower-ab}If~$V$ is an object of~$\Rep_R(G_{\Qp})$, then we write~$V^\ab$, respectively~$V_\ab$, for the maximal abelian $R[G_{\Q_p}]$-submodule of~$V$, respectively abelian $R[G_{\Q_p}]$-quotient module of~$V$
(where, as above, ``abelian'' 
has the meaning that the action of $G_{\Q_p}$ factors through~$G_{\Q_p}^{\ab}$).
\end{defn}

\begin{defn}%
  \label{defn:Dnr-functor}
We let $\bD^\nr:\Rep_R(G_{\Qp})\to\Rep_R(\Qptimes)$ be the functor given by $\bD^\nr(V)\coloneqq \bD(V)^{\nr}$ (which has an action of~$\Qptimes$ by Remark~\ref{rem:Dnr-is-Qptimes-stable}).%
\end{defn}

\begin{lemma}\label{Dnr and V^ab}\leavevmode
\begin{enumerate}
\item If~ $V\in \Rep_R(G_{\Qp})$ then $\bD^\nr(V) = \bD^\nr(V^\ab)$.
\item The restriction of the functor $\bD^\nr$
to ~$\Rep_R^{\ab}(G_{\Qp})$ %
is naturally isomorphic to the equivalence of categories $\Rep_R^{\ab}(G_{\Qp})\iso\Rep_R(\Qptimes)$ 
of Lemma~\em{\ref{lem:local Artin}}.
\end{enumerate}

\end{lemma}
\begin{proof}
For the first part, since~$\varprojlim_k$ commutes with limits, we have $\varprojlim_k (V_k)^\ab \cong V^\ab$.
Indeed, $V^\ab$ is naturally isomorphic to $\ker(V \xrightarrow{x_j-1} \prod_{j \in J}V)$, where~$\{x_j: j \in J\}$ is a set of generators of the commutator subgroup of~$G_{\bQ_p}$.
Hence it suffices to prove that the natural map $\bD^\nr((V_k)^\ab) \to \bD^\nr(V_k)$ is an isomorphism.
This is proved in~\cite[Rem.\ II.1.2]{ColmezMirabolique}.

Turning to the second part, if we regard $\Z_{p^n}$ as a $\Z_p[\Z/n\Z]$-module (via the action of geometric Frobenius), 
then there is an isomorphism of $\Z_p[\Z/n\Z]$-modules 
$\Z_p[\Z/n\Z] \iso \Z_{p^n}.$
We can choose these isomorphisms compatibly in $n$, in the sense that, for $m \mid  n$,  the natural 
surjection $\Z_p[\Z/n\Z] \to \Z_p[\Z/m\Z]$ becomes identified with the trace map
$\Tr^n_m: \Z_{p^n} \to \Z_{p^m}$. Equivalently, we may find elements $\delta_n$ of $\Z_{p^n}$
such that $\delta_n$ freely generates $\Z_{p^n}$ as a $\Z_p[\Z/n\Z]$-module, and such that
$\Tr^n_m(\delta_n) = \delta_m$ if~$m \mid n$. %

Next, if $V$ is a finite-cardinality $\Z_p$-module with a continuous action of $\widehat{\Z}$,
choose $n$ such that this action factors through~$\Z/n\Z$,
and define an isomorphism 
\begin{equation}
\label{eqn:Dnr iso}
V \iso (V \otimes_{\Z_p} \Z_{p^n})^{\widehat{\Z}}
= (V\otimes_{\Z_p} \breve \Z_p)^{\widehat{\Z}}
\end{equation}
(invariants being taken with respect to the diagonal action of~$\widehat{\Z}$)
via
$$v \mapsto \sum_{i = 0}^{n-1} \langle i \rangle ( v \otimes \delta_n),$$ 
where $\langle i\rangle$ denotes the action of the coset $i \bmod n \in \Z/n\Z$.
The trace relations between the various $\delta_n$ ensure that this isomorphism is
independent of the choice of~$n$ (provided that the action of $\widehat{\Z}$
on $V$ factors through $\Z/n\Z$).
This isomorphism intertwines the action of $\widehat{\Z}$ on~$V$ 
with the inverse of the action of $\widehat{\Z}$ on the second factor
of~$(V\otimes_{\Z_p} \breve \Z_p)^{\widehat{\Z}}.$

Now, if $V$ is endowed with a continuous action of 
$G_{\Q_p}^{\ab} \iso \widehat{\Q_p^{\times}} \iso p^{\widehat{\Z}} \times \Gamma$,
where $p$ corresponds to geometric Frobenius, and the projection to $\Gamma = \Z_p^{\times}$
is the cyclotomic character,
then
$\bD^{\nr}(V) \iso (V \otimes_{\Z_p} \breve \Z_p)^{\Frob},$
with the $\varphi$-action being induced by the action of Frobenius on~$\breve \Z_p$
and the $\Gamma$-action being induced by the action of $\Gamma$ on~$V$.
The isomorphism~\eqref{eqn:Dnr iso} then induces a $\Q_p^{\times}$-equivariant isomorphism
$$V \iso \bD^{\nr}(V).$$
By construction it is natural in~$V$, and so by passing to inverse limits,
it induces a corresponding isomorphism for any finite type $R$-module
endowed with a continuous $G_{\Q_p}^{\ab}$-representation.
\end{proof}

\begin{lemma}\label{Tate duality and Pontrjagin duality}%
There is a natural isomorphism of functors $\Rep_R(G_{\Qp})\to\Rep_R(\Qptimes)$ 
\[
\bD(V)^\sharp/\bD(V)^\natural \cong \bD^\nr(V_\ab) \otimes \varepsilon^{-1}.
\]
\end{lemma}
\begin{proof}
Since $D \mapsto D^\sharp/D^\natural$ and $D \mapsto D^\nr$ are both compatible with passage to inverse limits, it suffices to prove this when $R$ is an Artin ring.
Write~$(\text{--})^\vee = \Hom_{\bZ_p}(\text{--}, \bQ_p/\bZ_p)$.
It is shown in~\cite[Prop.\ I.2.3]{ColmezMirabolique} that if we put $D \coloneqq \bD(V)$ and $\check D \coloneqq \bD(V^\vee \otimes \varepsilon)$ then there is a $\bZ_p$-linear map
\[
D \otimes_{R} \check D \to \bQ_p/\bZ_p
\]
which identifies~$\check D$ with the Pontrjagin dual of~$D$ (i.e.\ the space of
continuous $\bZ_p$-linear maps $D \to \bQ_p/\bZ_p$, equipped with the weak
topology), and vice versa.
This map is invariant under~$\varphi$ and~$\Gamma$, and by~\cite[Prop.\ II.5.19]{ColmezMirabolique} it descends to a map
\[
D^\sharp/D^\natural \otimes_{R} \check D^\nr \to \bQ_p/\bZ_p
\]
which once again is $(\varphi, \Gamma)$-invariant, and puts the two factors in Pontrjagin duality (which now coincides with $\bZ_p$-duality, since the two sides have finite $\bZ_p$-length).
By Lemma~\ref{Dnr and V^ab} we have
\[
\check D^\nr = \bD^\nr((V^\vee \otimes \varepsilon)^\ab) \cong \bD^\nr((V_\ab)^\vee \otimes \varepsilon),
\]
so it suffices now to prove that~$\bD^\nr$ is compatible with duals and tensor products.
This is standard: for example, in order to prove that the natural map
\[
\bD^\nr(V_1) \otimes_{\bZ_p} \bD^\nr(V_2) \to \bD^{\nr}(V_1 \otimes_{\bZ_p} V_2)
\]
is an isomorphism, it suffices to extend scalars to~$\breve \bZ_p$, where it yields the isomorphism
\[
(\breve \bZ_p \otimes_{\bZ_p} V_1) \otimes_{\breve \bZ_p} (\breve
\bZ_p \otimes_{\bZ_p} V_2) \isom \breve\bZ_p \otimes_{\bZ_p} (V_1
\otimes_{\bZ_p} V_2). \qedhere
\]
\end{proof}

\begin{rem}\label{example trivial character}
To illustrate the twist in Lemma~\ref{Tate duality and Pontrjagin duality} with an example, recall from~\cite[Ex.\ II.4.5, II.5.16]{ColmezMirabolique} %
that when~$V$ is the trivial character of $\bF[G_{\bQ_p}]$ we have
\[
\bD(V)^\sharp = T^{-1}\bF\llbracket T\rrbracket  \text{ and } \bD(V)^\natural = \bF\llbracket T\rrbracket .
\] 
Then the action of $\Gamma$ on~$T^{-1}$ is via
\[
\gamma(T^{-1}) = ((1+T)^{\varepsilon(\gamma)}-1)^{-1} = (\lbar \varepsilon(\gamma)T + O(T^2))^{-1} = \lbar \varepsilon^{-1}(\gamma)T^{-1}u(T)
\]
for some $u(T) \in 1 + T\bF\llbracket T\rrbracket $.
\end{rem}%

\begin{lemma}\label{finiteness in formal case}
Let~$R$ be a complete Noetherian local $\cO$-algebra with finite residue field, and let~$D$ be a formal \'etale $(\varphi, \Gamma)$-module with $R$-coefficients.
Then $D^\sharp/D^\natural$ and $D^\nr$ are finite $R$-modules.
\end{lemma}
\begin{proof}
By Lemma~\ref{Tate duality and Pontrjagin duality} it suffices to prove this
for~$D^\nr$. %
Write~$V\coloneqq\mathbf{V}(D)$, %
so that by Lemma~\ref{Dnr and V^ab} the lemma is equivalent to the claim that $\bD^\nr(V^\ab)$ is $R$-finite.
By definition $\bD^\nr(V^\ab)$ is $\fm_R$-adically complete, and since~$D^\nr$ is exact on abelian $R[G_{\bQ_p}]$-modules (by Lemma~\ref{Dnr and V^ab}), %
 we know that $\bD^\nr(V^\ab) \otimes_R \bF \cong \bD^\nr(V^\ab \otimes_R \bF)$.
By Lemma~\ref{lem: Dnr is finite}, the module $\bD^\nr(V^\ab \otimes_R \bF)$ is $\bF$-finite.
This implies the required finiteness by~\cite[\href{https://stacks.math.columbia.edu/tag/031D}{Tag~031D}]{stacks-project}.%
\end{proof}

We also note the following application of the preceding results, which we will need below.
If~$x \in \Spec(R)$, we use the notation~$k(x)$ to denote the residue field of the localization of~$R$ at the prime ideal corresponding to~$x$.

\begin{lemma}\label{no D^ab}
Let $R$ be a reduced $\cO$-flat complete Noetherian local $\cO$-algebra with finite residue field.
Let~$V$ be a continuous $R[G_{\bQ_p}]$-module which is finite free over~$R$ of rank at least~$2$, 
and such that there exists a Zariski dense set of closed points~$x \in \Spec(R[1/p])$ 
such that $V \otimes_R k(x)$ 
is an absolutely irreducible $k(x)[G_{\bQ_p}]$-module.
Then 
$\bD^\nr(V)= 0$.
\end{lemma}
\begin{proof}%
By Lemma~\ref{Dnr and V^ab} we need to prove that $V^\ab = 0$.
Since $V^\ab$ is $\cO$-torsion free, and $V^\ab[1/p] \subset V[1/p]^\ab$ (where $V[1/p]^\ab$  is the maximal abelian submodule of $V[1/p]$), %
it suffices to prove that $V[1/p]^\ab = 0$.
However, $V[1/p]^\ab$ is contained in a finite free $R[1/p]$-module, hence its associated primes are minimal in $R[1/p]$, since $R[1/p]$ is reduced. 
To show that $V[1/p]^{\ab} = 0$ it thus suffices to show that it vanishes after localization at any minimal prime~$\eta$ of~$R[1/p]$.
To do so, it suffices to prove that $V[1/p] \otimes_R k(\eta)$ is an absolutely irreducible $k(\eta)[G_{\bQ_p}]$-module.
However, if this is not the case, then the properness of the Grassmannian of~$V[1/p]$ over~$\Spec R[1/p]$ shows that $V[1/p] \otimes_R k(x)$ is an absolutely 
reducible $k(x)[G_{\bQ_p}]$-module
for all~$x \in \lbar{\{ \eta \}}$: in fact, the $G_{\bQ_p}$-fixed locus is closed in the Grassmannian.
This contradicts our assumption on~$V$.
\end{proof}

We next establish some properties of the functor $D^{\sharp}\boxtimes \Q_p$.  It turns
out that this functor is very well behaved in the formal context,
thanks to the following result of Colmez, which shows that the formation of $D^{\sharp}\boxtimes \Q_p$ is exact
on \'etale $(\varphi,\Gamma)$-modules of finite $\A$-length.

\begin{thm}
  \label{thm:Colmez-exactness-Dsharp-box-Qp}
Let~$R$ be a complete Noetherian local $\cO$-algebra with finite residue field.
For any short exact sequence \[0\to D_1 \to D_2 \to D_3 \to 0\]  of (not necessarily projective, formal) \'etale $(\varphi,\Gamma)$-modules with $R$-coefficients, each of which is of finite length as an $\A$-module,  the sequence  \[0\to D_1^{\sharp}\boxtimes\Qp \to D_2^{\sharp}\boxtimes\Qp \to D_3^{\sharp}\boxtimes\Qp \to 0\] is also short exact.
\end{thm}
\begin{proof}
  This is immediate from~\cite[Thm.\ III.3.5]{ColmezMirabolique}.
\end{proof}%

\begin{cor}\label{cor:Colmez-exactness-Dsharp-box-Qp}
Let~$R$ be a complete Noetherian local $\cO$-algebra with finite residue field.
Then $D \mapsto D^\sharp \boxtimes \bQ_p$ is a right exact functor on the category of formal \'etale $(\varphi, \Gamma)$-modules with $R$-coefficients.
\end{cor}
\begin{proof}
Given an exact sequence $0 \to D_1 \to D \to D_2 \to 0$ of formal \'etale $(\varphi, \Gamma)$-modules with $R$-coefficients,
the sequence $(D_1)_k \to D_k \to (D_2)_k \to 0$ is exact for all~$k \geq 1$, and its terms have finite $\bA$-length.
By Theorem~\ref{thm:Colmez-exactness-Dsharp-box-Qp}, it stays exact after applying $(\text{--})^\sharp \boxtimes \bQ_p$, and then the corollary follows from the exactness of $\varprojlim_k$
on profinite $R$-modules. 
\end{proof}

Before formulating our next result,
we note that
if $D$ is a formal \'etale $(\varphi,\Gamma)$-module over~$R$,
and if $M$ is a finitely generated $R$-module,
then $D\otimes_R M$ 
(which agrees with $D\cotimes_{\widehat \bA_R} (\widehat \bA_R \otimes_R M)$;
see Lemma~\ref{lem:automatic tensor completeness}
and Remark~\ref{rem:automatic tensor completeness}) %
is again a formal \'etale $(\varphi,\Gamma)$-module,
and
$$(D\otimes_R M)_k = D_k\otimes_R M = D_k \otimes_{R/\m_R^k} (M/\m_{R}^k M)$$
is a not necessarily projective
\'etale $(\varphi,\Gamma)$-module over $R/\m_R^k$, for each $k\geq 1$.

\begin{lemma}
\label{Dsharp box Qp flatness}
Let~$R$ be a complete Noetherian local $\cO$-algebra with finite residue field,
and let $D$ be a formal \'etale $(\varphi,\Gamma)$-module over~$R$.
\begin{enumerate}
\item\label{item:65} If $M$ is a finitely generated $R$-module,
then there is a natural isomorphism
$$(D^{\sharp}\boxtimes \Qp)\otimes_R M \iso (D\otimes_R M)^{\sharp}\boxtimes \Qp.$$
\item\label{item:66} %
If $D$ is furthermore projective,
then $D^{\sharp}\boxtimes \Q_p$ is topologically flat over~$R$ {\em (}and
so in particular flat over~$R${\em )}.%
\end{enumerate}
\end{lemma}
\begin{proof}
For the first part, by Proposition~\ref{prop:Eilenberg-Watts}~\eqref{item:56} it suffices to show that the functor  $M \mapsto (D\otimes_R M)^{\sharp} \boxtimes \Q_p$ is  right exact
(on the category of finitely generated $R$-modules).
This follows from the right exactness of $M \mapsto D\otimes_R M$ and the right exactness of $D\mapsto D^{\sharp}\boxtimes\Qp$ 
(Corollary~\ref{cor:Colmez-exactness-Dsharp-box-Qp}). %

Now suppose 
that $0 \to M_1 \to M_2 \to M_3 \to 0$ is a short exact sequence of profinite 
$R$-modules.  Write this as the inverse limit of a short exact sequences
of finite-cardinality $R$-modules
$0 \to M_{1,i} \to M_{2,i} \to M_{3,i} \to 0.$
Since $D$ is projective, it is $R$-flat, and so
each of the sequences
$$0 \to D\otimes_{R} M_{1,i} \to D\otimes_{R} M_{2,i} \to D\otimes_{R} M_{3,i} \to 0$$
is again short exact.
By Theorem~\ref{thm:Colmez-exactness-Dsharp-box-Qp}, %
we find that 
$$0 \to (D\otimes_{R} M_{1,i})^{\sharp}\boxtimes\Qp \to (D\otimes_{R} M_{2,i})^{\sharp}
\boxtimes\Qp \to (D\otimes_{R} M_{3,i})^{\sharp}\boxtimes\Qp \to 0$$
is short exact.
Applying the isomorphism of~(1) (noting that it does indeed apply,
since the finite $R$-modules $M_{j,i}$
are in particular finitely generated), we find that the sequence
$$0
\to (D^{\sharp}\boxtimes \Qp)\otimes_R M_{1,i} 
\to (D^{\sharp}\boxtimes \Qp)\otimes_R M_{2,i} 
\to (D^{\sharp}\boxtimes \Qp)\otimes_R M_{3,i} 
\to 0
$$
is short exact.  Passing to the inverse limit over~$i$,
and recalling that the formation of inverse limits 
of profinite modules is exact,
we find that the sequence
$$0 
\to (D^{\sharp}\boxtimes \Qp)\cotimes_R M_{1} 
\to (D^{\sharp}\boxtimes \Qp)\cotimes_R M_{2} 
\to (D^{\sharp}\boxtimes \Qp)\cotimes_R M_{3} 
\to 0$$
is exact. Bearing in mind Lemma~\ref{topologically flat implies flat over Noetherian}, this proves~(2).
\end{proof}

\begin{lemma}
  \label{lem:formal D natural base change}
Let $R\to S$ be a morphism
of complete local Noetherian $\cO$-algebras with finite residue fields,
let~$D$ be a formal \'etale $(\varphi, \Gamma)$-module over~$R$,
and write $D_S \coloneqq  D\cotimes_R S.$
Then the natural map
$$(D^{\natural}\boxtimes \Q_p)\cotimes_R S \to D_S^{\natural}\boxtimes \Q_p$$
is surjective.
\end{lemma}
\begin{proof}
Suppose to begin with that $R$ and $S$ are both Artinian.
It follows from Lemma~\ref{lem: D natural base change}
that the morphism $D^{\natural}\otimes_R S \to D_S^{\natural}$ is surjective.
(Since $S$, being Artinian, is finite over~$R$, the tensor product is automatically
complete. Also, although the cited Lemma assumes projectivity of~$D$, the
claimed surjectivity follows in general from~\cite[Prop.~II.5.17]{ColmezMirabolique}.)
Passing to the inverse limit over~$\psi,$ and taking into account (as always)
that these inverse limits of profinite $\cO$-modules are exact,
we find that
$$\varprojlim_{\psi} (D^{\natural}\otimes_R S) \to D_S^{\natural}\boxtimes \Q_p$$
is surjective.
Now the natural morphism
$$(D^{\natural}\boxtimes {\Q_p})\cotimes_R S \to
\varprojlim_{\psi} (D^{\natural}\otimes_R S)$$
evidently has dense image, and hence is surjective (being a continuous
map of profinite $\cO$-modules).  
This completes the proof of the lemma in the Artinian case.

The general case follows by writing the morphism $R \to S$ as
the inverse limit of the morphisms $R/\fm_R^k \to S/\fm_S^k$,
and then passing to the inverse limit from the Artinian case (yet again
using the fact the formation of inverse limits of profinite $\cO$-modules
is exact). 
\end{proof}

\subsection{\texorpdfstring{Base-change and the action of $\GL_2(\Qp)$ on
  $D^\natural\boxdelta\Pone$}{Base change}}\label{subsec: base change}We now
return to the setting of a morphism $A\to B$ of Noetherian
$\cO/\varpi^a$-algebras, and a projective \'etale
$(\varphi,\Gamma)$-module with $A$-coefficients $D_A$. As in
Section~\ref{subsec: Tate module}, we write $D_B \coloneqq  D_A\cotimes_AB$;
then~$D_B$ is naturally an \'etale
$(\varphi,\Gamma)$-module with $B$-coefficients.

\begin{remark}
\label{rem:P box P1 lattices}
  Lemma~\ref{lem: seeing D box P1 as a Tate module} 
gives a topological isomorphism of Tate $A$-modules
$$D_A\boxdelta \Pone \iso D_A \oplus (D_A \boxtimes p\Z_p).$$
Now $D_A$ is a finite projective $\A_A$-module,
while $D_A\boxtimes p\Z_p = \varphi  (D_A)$
is a finite projective $\varphi(\A_A)$-module.   Since $\varphi$ induces
an isomorphism of $\A_A$ onto its image, we see that $D_A\boxtimes p\Z_p$ 
is isomorphic (as a Tate $A$-module) to a finite projective $\A_A$-module.
Hence $D_A\boxdelta \Pone$ itself is isomorphic to a finite projective $\A_A$-module,
and hence Example~\ref{ex:f.g. discrete} and Remark~\ref{rem:fake A((T))-modules} apply to~$D_A\boxdelta \Pone.$
\end{remark}

\begin{lem}
  \label{lem: base change for D box P1}Let $A\to B$ be a morphism of
  Noetherian $\cO/\varpi^a$-algebras for some $a\ge 1$. Let~$D_A$ be an
   \'etale $(\varphi,\Gamma)$-module with $A$-coefficients. Then
there is a natural isomorphism
  \[(D_A\boxdelta\Pone)\cotimes_AB \iso
  D_B\boxdelta\Pone,\]
which induces a morphism
$$(D_A^{\natural}\boxdelta \Pone )\cotimes_A B \to D_B^{\natural}\boxdelta \Pone.$$
\end{lem}
\begin{proof}%
  The isomorphism $D_A^{\oplus 2}\cotimes_AB\isoto D_B^{\oplus 2}$ induces a
  morphism  $(D_A\boxdelta\Pone)\cotimes_AB \to
  D_B\boxdelta\Pone$.  To show that this is an isomorphism,
  consider the following commutative diagram, where the vertical maps are the
  homeomorphisms  $z\mapsto (\Res_{\Zp}z,\Res_{p\Zp}(w\cdot z))$ of Lemma~\ref{lem: seeing D box P1 as a Tate module}:
  \[
  \begin{tikzcd}
     (D_A\boxdelta\Pone)\cotimes_AB\arrow[r]\arrow[d]  & D_B\boxdelta\Pone \arrow[d] \\
   (D_A\boxtimes\Zp)\cotimes_AB\oplus(D_A\boxtimes p\Zp)\cotimes_AB\arrow[r] & (D_B\boxtimes\Zp)\oplus(D_B\boxtimes p\Zp)    
  \end{tikzcd}
\]
Now $(D_A\boxtimes \Z_p)\cotimes_A B \coloneqq  D_A \cotimes_A B$ maps isomorphically
to $D_B \eqcolon D_B\boxtimes \Z_p$ by the very definition of~$D_B$, and so
  it suffices
  to show that also \[D_B\boxtimes p\Zp=(D_A\boxtimes
    p\Zp)\cotimes_AB,\]i.e.\
  that \[\varphi\psi(D_B)=(\varphi\psi(D_A))\cotimes_AB.\]Since
  $\varphi\psi$ is the projection
  onto a direct summand as a $\varphi(\A_A)$-module, we are done. %

To see that this induces a morphism
$(D_A^{\natural}\boxdelta \Pone )\cotimes_A B \to D_B^{\natural}\boxdelta \Pone$
as claimed,
it suffices (by the definition of $D^{\natural}$) to consider the commutativity of
the diagram
$$\xymatrix{
(D_A\boxdelta \Pone)\cotimes_A B \ar^-{\sim}[r] \ar[d] &  D_B\boxdelta \Pone
\ar[d] \\
(D_A\boxtimes\Q_p)\cotimes_A B \ar[r] &
D_B\boxtimes\Q_p}
$$
induced by $\Res_{\Q_p}$,
together with 
the fact that, by Lemma~\ref{lem: D natural base change}~(2),
the bottom horizontal arrow induces a morphism
$(D_A^{\natural}\boxtimes\Q_p)\cotimes_A B \to D_B^{\natural}\boxtimes\Q_p.$
\end{proof}

\begin{lem}\label{lem: D boxtimes P1 as A((T)) module}%
  \leavevmode
  Let~$A\to B$ be a morphism of Noetherian $\cO/\varpi^a$-algebras for
  some $a\ge 1$, and let~$D_A$ be an
   \'etale $(\varphi,\Gamma)$-module with $A$-coefficients.
  \begin{enumerate}
  \item\label{item: base change lattice box Pone} If $\gM_A\subseteq D_A\boxdelta\Pone$ is a lattice then the image of the
    natural morphism $\gM_{A}\cotimes_AB\to (D_A\boxdelta\Pone)\cotimes_AB=D_B\boxdelta\Pone$ is a lattice
    in~$D_B\boxdelta\Pone$. In particular if $A\to B$ is flat then $\gM_{A}\cotimes_AB$
    is a lattice in $D_B\boxdelta\Pone$.
  \item\label{item: descend lattice box times Pone} If $A\into B$ is injective, and~$\gM$ is a lattice in~$D_B\boxdelta\Pone$,
    then $\gM\cap (D_A\boxdelta\Pone)$ is a lattice in~$D_A\boxdelta\Pone$. 
  \end{enumerate}
\end{lem}
\begin{proof}
  This follows from %
  Lemma~\ref{lem: base change for D box P1} together with 
  parts~\eqref{item: lattice 4} and~\eqref{item: lattice 5} of
  Lemma~\ref{lem: collection of facts about lattices} (taking into account
  Remark~\ref{rem:P box P1 lattices}).
\end{proof}

\begin{lem}
  \label{lem: base change results we can use to prove G-stability}Let 
  $A\to B$ be a morphism of
  Noetherian $\cO/\varpi^a$-algebras for some $a\ge 1$. Let~$D_A$ be an \'etale
   $(\varphi,\Gamma)$-module with $A$-coefficients, and
  let~$D_B=D_A\cotimes_AB$.
  \begin{enumerate}%
     \item\label{item: going down for Dnatural P1 lattice} If $A\into B$ is injective,
and if $D_B^\natural\boxdelta\Pone$
    is an $\cO\llbracket G\rrbracket $-stable lattice in $D_B\boxdelta\Pone$, then
    $D_A^\natural\boxdelta\Pone$ is an $\cO\llbracket G\rrbracket $-stable lattice in
    $D_A\boxdelta\Pone$.
  \item\label{item: going up for Dnatural P1 lattice} If $D_A^\natural\boxdelta\Pone$ is an $\cO\llbracket G\rrbracket $-stable lattice in
    $D_A\boxdelta\Pone$, then $D_B^\natural\boxdelta\Pone$ is a
    $\cO\llbracket G\rrbracket $-stable lattice in
    $D_B\boxdelta\Pone$.
\item\label{item: flat case}
In the situation of~{\em \eqref{item: going up for Dnatural P1 lattice}},
 if~$A\to B$ is furthermore flat,
then $(D_A^\natural\boxdelta\Pone)\cotimes_AB$
maps isomorphically to $D_B^{\natural}\boxdelta\Pone.$
   \end{enumerate}
\end{lem}
\begin{proof}We begin with~\eqref{item: going down for Dnatural P1
    lattice}.  By Proposition~\ref{prop: G-stable lattice implies D
    natural box P1 stable}, $D_A^\natural\boxdelta\Pone$ is a
  $\cO\llbracket G\rrbracket $-stable lattice if and only if
  $D_A^\natural\boxdelta\Pone$ is a lattice if and only if there
  exists some $B(\Qp)$-stable lattice in $D_A\boxdelta\Pone$. Now,
  $(D_B^\natural\boxdelta\Pone)\cap (D_A\boxdelta\Pone)$ is a lattice
  in $D_A\boxdelta\Pone$ by Lemma~\ref{lem: D boxtimes P1 as A((T)) module}~\eqref{item: descend lattice box times Pone}, and it is $B(\Qp)$-stable
  because both $D_B^\natural\boxdelta\Pone$ and $D_A\boxdelta\Pone$
  are $\cO\llbracket G\rrbracket $-stable, so we are done. %

  We now turn to~\eqref{item: going up for Dnatural P1 lattice}. By
Lemma~\ref{lem: D boxtimes P1 as A((T)) module}~\eqref{item: base change lattice box Pone},
the image~$\gM_B$ of the
  morphism  \[(D_A^\natural\boxdelta\Pone)\cotimes_AB\to
    D_B\boxdelta\Pone\]is a lattice in $D_B\boxdelta\Pone$. It is
  $\cO\llbracket G\rrbracket $-stable by construction, so it follows from Proposition~\ref{prop: G-stable lattice implies D
    natural box P1 stable} that $D_B^\natural\boxdelta\Pone$ is a
  $\cO\llbracket G\rrbracket $-stable lattice.

Finally, let us put ourselves in the situation of~\eqref{item: flat case},
so that $A\to B$ is furthermore assumed to be flat.  
Then
Lemma~\ref{lem: D boxtimes P1 as A((T)) module}~\eqref{item: base change lattice box Pone}
shows that $(D_A^\natural\boxdelta\Pone)\cotimes_AB$ maps isomorphically onto its image in~$D_B\boxdelta \Pone$,
i.e.\ onto~$\fM_B$,
while Lemma~\ref{lem: base change for D box P1} shows that
$\gM_B \subseteq D_B^{\natural} \boxdelta \Pone.$
Combining this last fact with the conclusion of
Lemma~\ref{lem: minimal stable lattice in D box Qp},
we find that $\Res_{\Q_p}(\gM_B) = D_B^{\natural}\boxtimes \Q_p.$
Thus to show that $\gM_B = D_B^{\natural}\boxdelta \Pone,$
 it suffices by Lemma~\ref{lem: kernel of Res Qp} to show that $\ker \Res_{\Q_p} = w D_B^{\nr}$
is contained in $\gM_B$.  Since $\gM_B$ is $G$-invariant,
it is equivalent to show that $D_B^{\nr} \subseteq \gM_B$.
Since $D_B^{\nr} \subseteq D_B^+$, it suffices in turn to show
that $D_B^+$ 
(regarded as a submodule of~ $D_B=D_B\boxdelta\Zp$ and thus of~$D_B\boxdelta\Pone$)
is contained in $\gM_B.$

Now, 
by Lemma~\ref{lem: Dplus box Zp contained in Dnatural box P1}, we have an inclusion
$D_A^+ = D_A^+\boxdelta \Z_p \subseteq D_A^{\natural}\boxdelta\Pone,$
which induces a morphism
$(D_A^+)\cotimes_A B \to \gM_B$.  
By
Lemma~\ref{lem:flat base-change for Dplus}, this morphism
is an injection whose image is identified with $D_B^+$
when we identify $\gM_B$ with its image in $D_B^{\natural}\boxdelta \Pone$.
This completes the proof.
\end{proof}

\begin{lem}
  \label{lem: D natural in CNL reduces to checking it for
    quotients}Let~$(A,\m)$ be a complete local Noetherian
  $\cO/\varpi^a$-algebra, %
  and suppose that ~$D_A$ is an \'etale
  $(\varphi,\Gamma)$-module %
  with $A$-coefficients. Suppose
  that for each~$n\ge 1$, $(D_{A/\m^n})^\natural\boxdelta\Pone$ is an
  $\cO\llbracket G\rrbracket $-stable lattice in $D_{A/\m^n}\boxdelta\Pone$. Then
  $D_A^\natural\boxdelta\Pone$ is an $\cO\llbracket G\rrbracket $-stable lattice in
  $D_A\boxdelta\Pone$.
\end{lem}
\begin{proof}By Lemma~\ref{lem: G stable implies lattice}, it is
  enough to prove that $D_A^\natural\boxdelta\Pone$ is
  $\cO\llbracket G\rrbracket $-stable; so by definition we need to prove that if
  $f \in \cO\llbracket G\rrbracket $
  and if $x\in D_A\boxdelta\Pone$ satisfies
  $\Res_{\Zp}\left(
  \big(\begin{smallmatrix}
    p^n&0\\0&1
  \end{smallmatrix}\big)
x\right)\in D_A^\natural$ for all~$n\ge 0$, then
  $\Res_{\Zp}\left( \big(\begin{smallmatrix}
    p^n&0\\0&1
  \end{smallmatrix}\big)fx\right)\in D_A^\natural$ for all~$n\ge 0$.

Now, for each~$m\ge 1$ the assumption that $(D_{A/\m^m})^\natural\boxdelta\Pone$ is
$\cO\llbracket G\rrbracket $-stable means that the image of  $\Res_{\Zp}\left( \big(\begin{smallmatrix}
    p^n&0\\0&1
  \end{smallmatrix}\big)fx\right)$ in~$D_{A/\m^m}$ is contained in
$D^\natural_{A/\m^m}$. By Lemma~\ref{lem: D natural base change}~(3), this
means that for each~$n\ge 0$, $\Res_{\Zp}\left( \big(\begin{smallmatrix}
    p^n&0\\0&1
  \end{smallmatrix}\big)fx\right)$ is contained in
$D_A^\natural+\m^mD_A$. Since this holds for all~$m$,
the lemma follows by an application of Lemma~\ref{lem:lattice completeness}.
\end{proof}

We now specialize to the case that~$D$ has rank 2. %

\begin{thm}
  \label{thm: D natural is G stable in the CNL case}
  Let~$A$ be a
  complete local Noetherian $\cO/\varpi^a$-algebra with finite residue
  field, and suppose that ~$D_A$ is an
   \'etale $(\varphi,\Gamma)$-module of rank~$2$ with
  $A$-coefficients and determinant~$\zeta\varepsilon^{-1}$. Then $D_A^\natural\boxtimes\Pone$ is a
  $\cO\llbracket G\rrbracket $-stable lattice in $D_A\boxtimes\Pone$.
\end{thm}%
\begin{proof}By Lemma~\ref{lem: D natural in CNL reduces to checking it for
    quotients}, we may assume (by replacing $A$ by $A/\m^n$ where $\m$
  is the maximal ideal of~$A$) that~$A$ is Artinian, in which
  case~$D_A$ corresponds to a 2-dimensional representation
  $\rho_A:G_{\Qp}\to\GL_2(A)$. By Lemma~\ref{lem: base change results
    we can use to prove G-stability}~\eqref{item: going up for Dnatural
    P1 lattice}, we can assume that~$A$ is a quotient of the universal (fixed determinant)
  framed deformation ring~$S$ for~$\rho_A\pmod \m$.

  As explained in~\cite[\S II.3.4]{MR2642409}, there is a universal
  formal \'etale $(\varphi,\Gamma)$ module~$D$ over~$S$, %
  and it is proved there
  that $D^\natural\boxtimes\Pone$ is
  $G$-stable; indeed, by~\cite[Lem.~II.3.6]{MR2642409}, this follows from the
  results of Berger--Breuil~\cite{MR2642406} together with %
  the Zariski density of the so-called crystalline benign points in
  the generic fibre $\Spec S[1/p]$, which %
  footnote 7 of~\cite{MR2642409} explains is known for $p\geq 3$.
  (In fact, 
  Zariski density statements of this kind are now known in much greater generality; see~\cite[Thm.\ 6.1, Rem.\ 6.2]{BIP} and~\cite[Cor.\ 1.3]{BIP2}.) %
  The first paragraph of the proof of~\cite[Prop.~II.2.15]{MR2642409} then shows that
  $D_{A}^{\natural}\boxtimes \Pone$ is $G$-stable, %
  and hence also~$\cO\llbracket G\rrbracket $-stable, 
  for any Artinian quotient $A$ of~$S$,
  as required.
  \end{proof}

Before stating the next main result of this section,
we note the following (presumably standard) lemma.

\begin{lem}\label{lem: finite type injects into product of CNL}
If~$A$ is a Noetherian ring,
  then for some finite set $\m_1,\dots,\m_r$ of maximal
  ideals of~$A$, the natural map $A\to\prod_{i=1}^r\widehat{A}_{\m_i}$
  is injective.  
\end{lem} %
\begin{proof}
As in the proof of Lemma~\ref{lem:Artinian replacement},
we may assume $A$ is non-zero, and then choose an element $f \in A$ and a positive 
integer $n$ so that $A \to A_f \times A/f^n$
is injective, and so that $A_f$ contains a unique associated prime~$\mathfrak{p}$.  
We then see that $A_f \to A_{\mathfrak m_1}$, and hence also $A_f \to \widehat{A}_{\m_1}$,
is injective, for any maximal ideal
$\mathfrak m_1 \supset \mathfrak p.$ 
Note that maximal ideals~$\mathfrak n$ of $A/f^n$ correspond to maximal ideals of $A$ containing~$f$,
and that for any such~$\mathfrak n$, the composite $A \to A/f^n \to (A/f^n)_{\mathfrak n}
\to \widehat{(A/f^n)}_{\mathfrak n}$
factors through the natural map $A \to \widehat{A}_{\mathfrak n}$.  
Thus we may continue by Noetherian induction and establish the lemma.
\end{proof}

\begin{thm}%
  \label{thm: D natural is G-stable} Suppose that~$p\ge 3$. Assume that~$A$ is a
  Noetherian
  $\cO/\varpi^a$-algebra for some~$a\ge 1$, and that~$D$ is a
  projective \'etale $(\varphi,\Gamma)$-module of rank~$2$ with
  $A$-coefficients and determinant $\zeta\varepsilon^{-1}$. Then $D^\natural\boxtimes\Pone$ is a
  $\cO\llbracket G\rrbracket $-stable lattice in $D\boxtimes\Pone$.
\end{thm}%
\begin{proof} Since $\cX$ is of finite presentation,
the morphism $\Spec A \to \cX$ that classifies $D$ may be
factored as $\Spec A \to \Spec A_0 \to \cX$ for some finite type $\cO/\varpi^a$-subalgebra $A_0$
of~$A$.
The morphism $\Spec A_0 \to \cX$ corresponds to a
rank~$2$  \'etale $(\varphi,\Gamma)$-module $D_{A_0}$
with $A_0$-coefficients, such that $D \iso D_{A_0} \cotimes_{A_0} A.$ 
By Lemma~\ref{lem: base change results we can use to prove G-stability}~(2),
it then suffices to show the result for $D_{A_0}$; replacing $A$ by $A_0$,
we may thus suppose that $A$ is in fact of finite type over~$\cO/\varpi^a$. 

By Lemma~\ref{lem: finite type injects
    into product of CNL} and Lemma~\ref{lem: base change results we can
    use to prove G-stability}~\eqref{item: going down for Dnatural P1
    lattice}, we can further reduce to the case that~$A$ is a finite
  product of complete local Noetherian $\cO/\varpi^a$-algebras with
  finite residue fields. This immediately reduces to the case of a
  single such complete local Noetherian $\cO/\varpi^a$-algebra, which is Theorem~\ref{thm: D natural is G stable in the CNL
    case}. %
\end{proof}

\begin{rem}
  \label{rem: short exact sequence for restriction to Qp}
Note in particular that as a consequence of Lemmas~\ref{lem: kernel of Res Qp}
and~\ref{lem: minimal stable lattice in D box Qp},
and Theorem~\ref{thm: D natural is G-stable}, we have a short exact
sequence of $B$-representations
\begin{equation}
\label{eqn:short exact sequence for restriction to Qp}
0\to \ker\ResQp = w D^{\nr} \to
D^\natural\boxtimes\Pone\to D^\natural\boxtimes\Qp\to 0. 
\end{equation}
Lemma~\ref{lem:f.g. discrete} (and Remark~\ref{rem:P box P1 lattices})
show that the induced topology on $w D^{\nr}$ is discrete.
\end{rem}

\begin{cor}\label{cor:G-stable-lattices-formal}Let $R$ be a complete local
Noetherian $\cO$-algebra with finite residue field, and let ~$D$ be a free rank 2
formal \'etale $(\varphi,\Gamma)$-module with $R$-coefficients and determinant
$\zeta\varepsilon^{-1}$.
Then
  $D^\natural \boxtimes \bP^1$ %
is~$\cO\llbracket G\rrbracket $-stable,
and we have a short exact sequence of $B$-representations \begin{equation}\label{eqn:formal-kernel-ResQp}0\to\ker\ResQp=wD^{\nr}\to D^\natural\boxtimes\Pone\to
  D^\natural\boxtimes\Qp\to 0. \end{equation}
\end{cor}
\begin{proof}
  Using the exactness of projective limits of profinite sets, this is immediate from Theorem~\ref{thm: D natural is G-stable} and Remark~\ref{rem: short exact sequence for
    restriction to Qp}.
\end{proof}

Recall from Definition~\ref{defn:Iwasawa-algebra-coefficients} and Lemma~\ref{lem:topology-Iwasawa-algebra} that if~$A$ is a finite type
  $\cO/\varpi^a$-algebra for some~$a\ge 1$, the Iwasawa algebra $A\llbracket K\rrbracket $ is a topological $A$-algebra (where $A$ has the discrete topology), whose topology agrees with the $\mathfrak a$-adic topology, where
 $\mathfrak a$ denotes
the ideal in $A\llbracket K\rrbracket $ generated by the augmentation ideal of $A\llbracket K_1\rrbracket $.
\begin{cor}%
  \label{cor: D natural box P1 is finite over K} Suppose that~$p\ge 3$. Assume that~$A$ is a finite type
  $\cO/\varpi^a$-algebra for some~$a\ge 1$, and that~$D$ is an
   \'etale $(\varphi,\Gamma)$-module of rank~$2$ with
  $A$-coefficients and determinant $\zeta\varepsilon^{-1}$. Then $D^\natural\boxtimes\Pone$ is a
 finitely generated $A\llbracket K\rrbracket $-module. Furthermore, 
the topology on $D^{\natural}\boxtimes \Pone$ coincides with
its $\mathfrak a$-adic topology.
\end{cor}
\begin{proof}
Consider the $A\llbracket K\rrbracket $-submodule $A\llbracket K\rrbracket  D^+\subseteq
D^\natural\boxtimes\Pone$. This is a lattice (since it contains $D^+$
and $wD^+$), and it is finitely generated over $A\llbracket K\rrbracket $, because $D^+$
is finitely generated over~$\bA_A^+ = A\llbracket T\rrbracket $ (being an $\bA_A^+$-lattice in~$D$) 
and thus over $A\llbracket U(\Zp)\rrbracket $. The
quotient $(D^\natural\boxtimes\Pone)/(A\llbracket K\rrbracket D^+)$ is finitely
generated over~$A$ and thus over~$A\llbracket K\rrbracket $, %
so that $D^{\natural}\boxtimes\Pone$ is indeed finitely generated over $A\llbracket K\rrbracket $.
Now Corollary~\ref{cor: O-G  action on G box P1}
shows that $D^{\natural}\boxtimes \Pone$ is in fact
a topological $A\llbracket K\rrbracket $-module, when endowed with its natural topology as a lattice in an $A$-Tate module. 
Since
this topology is complete and first countable
(and so completely metrizable), the general theory of finitely generated 
modules over Polish topological rings (cf.\ \cite[Prop.~C.6]{emertongeepicture}
as well as Remark~\ref{rem:Polish Iwasawa}) 
shows that the topology on~$D^{\natural}\boxtimes\Pone$ is its canonical
topology as a finitely generated~$A\llbracket K\rrbracket $-module, which in turn is precisely
its~$\fa$-adic topology (as was noted in Remark~\ref{rem:Polish Iwasawa}).
\end{proof}

\subsection{Non-flat base-change} 
\label{subsec:non-flat base-change}
Let $A \to B$ be a morphism of Noetherian $\cO/\varpi^a$-algebras,
and suppose given a projective rank two \'etale $(\varphi,\Gamma)$-module $D_A$ over~$\A_A$
with determinant~$\zeta \varepsilon^{-1}$.
As usual, we then write
$$D_B \coloneqq  D_A \otimes_{\A_A}\A_B = D_A\cotimes_A B.$$
Then by Lemma~\ref{lem: base change for D box P1} we have a natural base change morphism \begin{equation}
\label{eqn:non-flat base-change map for Pone}
(D_A^{\natural}\boxtimes \Pone)\cotimes_A B \to D_B^{\natural}\boxtimes \Pone.
\end{equation}

\begin{lemma}%
\label{lem:non-flat base-change} The morphism \eqref{eqn:non-flat base-change map for Pone} 
is an open mapping, whose image is a lattice in $D_B\boxtimes\Pone$.
Its kernel and cokernel 
are both finite over~$B$, and are both discrete {\em(}with respect to their
induced and quotient topologies, respectively{\em )}.
\end{lemma}
\begin{proof}

By Theorem~\ref{thm: D natural is G-stable}, the source and target of~\eqref{eqn:non-flat base-change map for Pone} are lattices in~$D_A \boxtimes \Pone$, resp.\ $D_B \boxtimes \Pone$.
By Lemma~\ref{lem: collection of facts about lattices}~(4)
(together with Remark~\ref{rem:P box P1 lattices}), the morphism~\eqref{eqn:non-flat base-change map for Pone}
is thus an open mapping, whose kernel is finite over~$B$ and discrete, and whose image is a lattice in $D_B\boxtimes\Pone$. 
Thus the cokernel of \eqref{eqn:non-flat base-change map for Pone} is a quotient of lattices,
and hence is discrete and $B$-finite. %
\end{proof}%

Our main focus in this section is on  base-change %
results for non-flat morphisms; in the flat case, we have the following result. %
\begin{lem}
  \label{lem: flat base change for D natural box P1 and D/Dnatural box P1}If $A\to B$ is a
  flat morphism of Noetherian %
$\cO/\varpi^a$-algebras, and~$D_A$ is an
   \'etale $(\varphi,\Gamma)$-module of rank~$2$ with $A$-coefficients and determinant $\zeta\varepsilon^{-1}$,
  then $(D_A^\natural\boxtimes\Pone)\widehat{\otimes}_AB=D_B^\natural\boxtimes\Pone$. %
\end{lem}%
\begin{proof}
This follows from 
  Lemma~\ref{lem: base change results we can
    use to prove G-stability}~\eqref{item: flat case} and Theorem~\ref{thm: D natural is G-stable}. %
\end{proof}

\subsubsection{Base-change in the formal context}
\label{subsubsec:formal base change}
Let $R\to S$ be a finite morphism of complete Noetherian local
$\cO$-algebras, with maximal ideals $\m$ and $\fn$ respectively,
each having finite residue field.
Let $D_{R}$ denote a projective rank $2$ formal \'etale $(\varphi,\Gamma)$-module
with $R$-coefficients
and with determinant~$\zeta \varepsilon^{-1}$, and set $D_S\coloneqq D_R\cotimes_RS=D_R\otimes_RS$.
Taking inverse limits of the morphisms~\eqref{eqn:non-flat base-change map for Pone} for $A=R/\m^n$, $B=S/\fn^{n}$, and using
Lemma~\ref{lem:tensor with finite algebra is automatically complete}
to replace a completed tensor product by a tensor product, %
we obtain a natural base-change morphism
\begin{equation}
\label{eqn:non-flat base-change map in the local case}
(D_R^{\natural}\boxtimes \Pone) \otimes_R S \to
D_S^{\natural}\boxtimes \Pone.
\end{equation}

\begin{lem}\label{lem:finite-ker-coker-completed-base-change}
  The kernel~$\cK$ and cokernel~$\cC$ of~\eqref{eqn:non-flat base-change map in the local case} are finite $S$-modules.
\end{lem}
\begin{proof}
  Consider the short exact sequence %
$$0 \to D_R^{\natural}\boxtimes \Q_p \to D_R^{\sharp}\boxtimes \Q_p
\to D_R^{\sharp}/D_R^{\natural} \to 0$$ arising from
Lemma~\ref{lem:formal-etale-phi-gamma-behave-well}.
Taking the tensor product %
of this short exact sequence with $S$ over~$R$, and taking into account %
Lemma~\ref{Dsharp box Qp flatness}~\eqref{item:66}, %
we obtain an exact sequence of $B$-representations
$$0 \to \Tor_1^R( D_R^{\sharp}/D_R^{\natural}, S) 
\to (D_R^{\natural}\boxtimes \Q_p)\otimes_R S \to (D_R^{\sharp}\boxtimes \Q_p)\otimes_R S \to (D_R^{\sharp}/D_R^{\natural})\otimes_R S \to 0.$$
Lemma~\ref{Dsharp box Qp flatness}~\eqref{item:65} shows that $(D_R^{\sharp}\boxtimes \Q_p)\otimes_RS \iso D_S^{\sharp}\boxtimes \Q_p$,
while Lemma~\ref{lem:formal D natural base change}
shows that 
    $(D_R^{\natural}\boxtimes \Qp) \otimes_R S \to D_S^{\natural}\boxtimes \Qp$ is surjective. 
Thus  we can rewrite the preceding exact sequence in the form
\begin{equation}
  \label{eqn:exact sequence for a finite morphism}
  0 \to \Tor_1^R( D_R^{\sharp}/D_R^{\natural}, S) 
  \to (D_R^{\natural}\boxtimes \Q_p)\otimes_R S
\to D_S^{\natural}\boxtimes \Q_p 
  \to 0.
\end{equation}

Next, %
the exact sequences of~\eqref{eqn:formal-kernel-ResQp}
give rise to
a diagram of $B$-representations
$$\xymatrix{
  & (w \cdot D_R^{\nr})\otimes_R S \ar[r]\ar[d] & (D_R^{\natural} \boxtimes \Pone)\otimes_R S \ar[r]
\ar^-{\text{ \eqref{eqn:non-flat base-change map in the local case}}}[d] &
  (D^{\natural}_R \boxtimes \Q_p)\otimes_R S \ar[r]
\ar[d] & 0 \\
  0 \ar[r] & w\cdot D_S^{\nr} \ar[r] & D_S^{\natural}\boxtimes \Pone \ar[r] &
  D_S^{\natural} \boxtimes \Q_p \ar[r] & 0\\
}
$$
Applying the snake lemma, we obtain an exact sequence of $B$-representations
\begin{multline}
  \label{eqn:non-flat basechange exact sequence}
  \ker\bigl( (w\cdot D_R^{\nr})\otimes_R S \to w\cdot D_S^\nr \bigr)
  \to \cK \to
  \Tor_1^R( D_R^{\sharp}/D_R^{\natural}, S) 
  \\
  \to 
  \coker\bigl( (w\cdot D_R^{\nr})\otimes_R S \to w\cdot D_S^\nr \bigr)
  \to \cC \to 0.
\end{multline}By Lemma~\ref{finiteness in formal case}, $D_R^{\nr}$ and $D_R^{\sharp}/D_R^{\natural}$ are finite $R$-modules, and~$D_S^\nr$ is a finite $S$-module. It follows from~\eqref{eqn:non-flat basechange exact sequence} that~$\cK$ and~$\cC$ are finite $S$-modules, as required.
\end{proof}

We can now exploit the correspondence between formal \'etale $(\varphi,\Gamma)$-modules and Galois representations to say something more precise about~$\cK$ and~$\cC$. Write $V\coloneqq \mathbf{V}(D_R)$,  %
so that %
$D_R \iso \bD(V).$ 
Lemma~\ref{Dnr and V^ab} %
shows that
$$D_R^{\nr} \iso \bD^{\nr}(V) \iso V^{\ab}$$
as $R$-linear $\Q_p^{\times}$-representations,
while
Lemma~\ref{Tate duality and Pontrjagin duality}
shows that
$$D_R^{\sharp}/D_R^{\natural} \iso \bD^{\sharp}(V)/\bD^{\natural}(V) 
\iso V_{\ab}\otimes \varepsilon^{-1}$$
as $R$-linear $\Qp^{\times}$-representations.
Of course, there are analogous isomorphisms for $V_S\coloneqq V\otimes_RS=\mathbf{V}(D_S)$. %

Since $\cK$ and $\cC$ are finite over~$S$, %
Lemma~\ref{trivial action} and Remark~\ref{rem: central char} show that $G$ acts on each of them
through a direct sum of $\beta\circ \det$-eigenspaces, where~
$\beta$ runs over the characters $\beta: \Q_p^{\times} \to \cO^{\times}$ with $\beta^2 = \zeta \varepsilon^{-2}$.
We can then read off $\beta$ as the restriction of the $G$-action
to the image of the embedding~\eqref{eqn:Qp-times embedding into G};
in other words, the $G$-actions on $\cK$ and~$\cC$ are entirely determined by
their restriction to this image, i.e.\ by their structure as $\Qp^{\times}$-modules. 

We can then rewrite the  exact sequences of $B$-representations ~\eqref{eqn:exact sequence for a finite morphism}
and~\eqref{eqn:non-flat basechange exact sequence} in the form
\begin{equation}
\label{eqn:exact sequence for a finite morphism V form}
0 \to \Tor_1^R( V_{\ab}\otimes \varepsilon^{-1}, S) 
\to (D_R^{\natural}\boxtimes \Q_p)\otimes_R S \to D_S^{\natural}\boxtimes \Q_p 
\to 0
\end{equation}
and
\begin{multline}
\label{eqn:non-flat basechange exact sequence V form}
\ker\bigl( (w\cdot V^{\ab})\otimes_R S \to  w\cdot (V\otimes_R S)^{\ab} \bigr)
\to \cK \to
 \Tor_1^R( V_{\ab}\otimes \varepsilon^{-1}, S) 
\\
 \to 
 \coker\bigl( (w\cdot V^{\ab})\otimes_R S \to w\cdot(V\otimes_R S)^\ab \bigr)
\to \cC \to 0.
\end{multline}

As an example of how the preceding results can be applied,
we establish the following result, which greatly constrains the 
circumstances in which $\cK$ and $\cC$ can be non-zero.

\begin{lemma}
\label{lem:basechange failure criterion}
Write $\overline{V} \coloneqq  V/\mathfrak m$, a continuous two dimensional 
representation of $G_{\Q_p}$ over the residue field $k \coloneqq  R/\m$ 
{\em (}which is a finite extension of~$\F${\em )}. Then the base-change morphism ~\eqref{eqn:non-flat base-change map in the local case} is an isomorphism unless
$\overline{V}^{\semis} \cong \chi \oplus \chi \omega$ %
for some continuous character $\chi: G_{\Q_p} \to k^{\times}$.  
\end{lemma}
\begin{proof}
As already noted,
each of $\cK$ and $\cC$ is a finite $S$-module equipped with a continuous $G$-action,
and admitting~$\zeta\varepsilon^{-1}$ as a central character,
so that, by Lemma~\ref{trivial action}, $G$ acts on each of them
through a direct sum of $\beta\circ \det$-eigenspaces, where
$\beta: \Q_p^{\times} \to \cO^{\times}$ and $\beta^2 = \zeta \varepsilon^{-2}$.

Considering~\eqref{eqn:exact sequence for a finite morphism V form} and~\eqref{eqn:non-flat basechange exact sequence V form}, we see that any such character $\beta$ must then contribute to one of~$V^{\ab}\otimes_R S,$
$(V_S)^{\ab}$, %
or $\Tor_1^R(V_{\ab}\otimes \varepsilon^{-1},S).$   
Thus the reduction $\overline{\beta}: \Q_p^{\times} \to \F^{\times}$
must be a constituent of either~$\overline{V}$ or~$\overline{V}\otimes \omega^{-1}$. 
Since $\det \overline{V} = \overline{\zeta} \omega^{-1},$
while $\overline{\beta}^2 = \overline{\zeta} \omega^{-2}$,
we find that %
$\overline{V}^{\semis} \iso \overline{\beta} \oplus \overline{\beta} \omega,$
proving the lemma.
\end{proof}

\subsubsection{Compatibility with completion}
\label{subsubsec:compatibility}
If $R$ is  a complete local Noetherian  $\cO/\varpi^a$-algebra with finite residue field (for some $a \geq 1$),
then in addition to the theory of formal \'etale $(\varphi,\Gamma)$-modules that we have 
discussed in Section~\ref{subsec:formal phi Gamma modules},
we can also consider %
the context of the theory of
projective \'etale $(\varphi,\Gamma)$-modules with $R$-coefficients,
as developed in the preceding sections.
Given a
projective \'etale $(\varphi,\Gamma)$-module
$D$ with $R$-coefficients,
we may form $\widehat{D} \coloneqq  \varprojlim_k D/\m_R^k D,$
a projective formal \'etale $(\varphi,\Gamma)$-module with $R$-coefficients.
We wish to compare some of the preceding constructions for $\widehat{D}$
with the analogous constructions for $D$ itself.

Before doing this, we need a brief discussion on topologies.
Recall that if $L$ is a lattice in an $R$-Tate module,
then in addition to its Tate-module topology,
it has a {\em weak topology}, as 
described in Section~\ref{subsubsec:weak topology}.
In particular, %
$D^{\natural}$ admits a weak topology, obtained as the inverse limit topology of the quotient topology on $D^\natural/\fm^k_R D^\natural$ under the isomorphism
$D^\natural \isom \varprojlim_k (D^\natural/\fm_R^k D^\natural)$ of Lemma~\ref{lattices in Tate modules are complete}.
On the other hand, $(\widehat D)^\natural \coloneqq  \varprojlim_k (D/\fm^k_R D)^\natural$ is endowed with the inverse limit topology with respect to the $R/\fm_R^k$-Tate module topology
on~$(D/\fm^k_R D)^\natural$.

Recalling that
$D^{\natural}\boxtimes \Q_p \coloneqq  \varprojlim_{\psi} D^{\natural}$,
we also define a weak topology on $D^{\natural}\boxtimes \Q_p$,
as the inverse limit topology arising from the weak topology
on each of the terms in the inverse limit that defines it. %

\begin{lemma}
\label{lem:compatibility}
Let~$R$ be a complete Noetherian local $\cO/\varpi^a$-algebra with finite residue field, and let~$D$ be a projective \'etale $(\varphi,\Gamma)$-module with $R$-coefficients. 
Then there are natural topological isomorphisms:
\begin{enumerate}
\item\label{item:67}  $D^{\natural} \iso (\widehat{D})^{\natural}$,
\item\label{item:69}  $D^{\natural}\boxtimes \Q_p \iso (\widehat{D})^{\natural}\boxtimes \Q_p$,
\item \label{item:70} $D^{\nr} \iso (\widehat{D})^{\nr}$,
\end{enumerate}
where~$D^{\natural}$ and~$D^{\natural}\boxtimes \Q_p$ are endowed with their weak topologies, and ~$D^{\nr}$ {\em (}a finite $R$-module,
by Lemma~{\em \ref{lem: Dnr is finite})} is endowed with its $\m_R$-adic topology.
If~$D$ furthermore has rank two and  determinant~$\zeta \varepsilon^{-1}$, then there is a natural topological isomorphism:
\begin{enumerate}[start=4]
\item\label{item:71}  $D^{\natural}\boxtimes \Pone \iso (\widehat{D})^{\natural}\boxtimes \Pone$,
\end{enumerate}
where~$D^{\natural}\boxtimes \Pone$ is endowed with its weak topology.

\end{lemma}
\begin{proof}
Throughout the proof we use the fact that completed tensor products $\text{--}\cotimes_R
R/\m_R^k$ coincide with usual tensor products $\text{--}\otimes_R R/\m_R^k 
= (\text{--})/\m_R^k$, since 
$R/\m_R^k$ is finitely presented over~$R$ (see Lemma~\ref{lem:automatic tensor completeness})

Write $D_k \coloneqq  D/\m_R^k D$ for $k \geq 1$,
so that $\widehat{D} \coloneqq  \varprojlim_k D_k$.
The natural map $D \to \widehat{D}$ is an injection (because $\A_R \to \widehat{\A}_R$
is an injection and~$D$ is projective).
Lemma~\ref{lem: D natural base change}~(2)
shows that for each $k$ we have a short exact sequence
$$0 \to M_k \to D^{\natural}/\m^k_R D^\natural \to (D_k)^{\natural} \to 0$$ 
(defining~$M_k$),
and Lemma~\ref{lem: collection of facts about lattices}~(4)
shows that each $M_k$ is a finite $R/\m_R^k$-module, 
and so is a finite set.
Note that all arrows in this exact sequence are continuous for the $R/\fm_R^k$-Tate module topology on $(D_k)^\natural$ and the quotient topology on $D^{\natural}/\m^k_R D^\natural$.
Passing to the inverse limit over~$k$, %
we obtain a continuous surjection $D^{\natural} \to (\widehat{D})^{\natural}$.
Since this map is also obtained by restricting the injection
$D \hookrightarrow \widehat{D},$    %
it is also an injection.  
Since $D^\natural$ is weakly compact, by Lemma~\ref{lattices in Tate modules are complete}, and $(\widehat D)^\natural$ is Hausdorff, this map is a topological isomorphism.
This proves~\eqref{item:67}.

Part~\ref{item:69}  follows directly from part~\eqref{item:67}.  
More precisely, we see that
$$D^{\natural}\boxtimes \Q_p \coloneqq  \varprojlim_\psi D^{\natural}
\iso \varprojlim_\psi (\widehat{D})^{\natural} \eqcolon (\widehat{D})^{\natural}\boxtimes\Qp.$$ 
We now turn to part~\eqref{item:70}. Note firstly that each surjection $D\to D_k$ restricts to a morphism
$D^\nr \to (D_k)^\nr,$ which, taken together, yield a morphism $D^\nr \to \varprojlim_k (D_k)^\nr.$
This is in fact an injection, being the restriction of the embedding $D \hookrightarrow \widehat{D}$.
Since $(D_k)^{\nr} \subseteq (D_k)^{\natural}$, we furthermore find (using part~\eqref{item:67}) that
$\varprojlim_k (D_k)^{\nr} \subseteq D^{\natural} \subset D.$ 
Now $\varphi$ acts bijectively on each $(D_k)^{\nr}$, and thus on $\varprojlim (D_k)^{\nr}$.
Thus in fact 
$(\widehat D)^{\nr} \coloneqq \varprojlim (D_k)^{\nr} \subseteq D^{\nr}$ (by the definition of $D^{\nr}$).  This gives the
opposite inclusion to the one already proved, and so establishes part~\eqref{item:70}.

Finally, we turn to part~\eqref{item:71}.
By our definition of the functors $(\text{--})^{\nr}, (\text{--})^\natural \boxtimes \Pone$ and $(\text{--})^\sharp\boxtimes \Pone$ on formal \'etale $(\varphi, \Gamma)$-modules (see 
Definition~\ref{def: various functors on formal phi Gamma modules}),
we have a commutative diagram
\[
\begin{tikzcd}
  &w \cdot D^{\nr} \ar[r] \ar[d] &
  D^{\natural} \boxtimes \mathbf{P}^1 \ar[r] \ar[d] &
  D^{\natural} \boxtimes \mathbf{Q}_p \ar[d] \arrow[r] & 0 \\
  0 \arrow[r] & w \cdot (\widehat{D})^{\nr} \ar[r] &
  (\widehat{D})^{\natural} \boxtimes \mathbf{P}^1 \ar[r] &
  (\widehat{D})^{\natural} \boxtimes \mathbf{Q}_p. &
\end{tikzcd}
\]
The first, resp.\ third vertical arrow is an isomorphism by part~\eqref{item:70}, resp.\ part~\eqref{item:69}.
The first, resp.\ second row is exact by Remark~\ref{rem: short exact sequence for restriction to Qp}, resp.\ Lemma~\ref{lem: kernel of Res Qp}.
An application of the snake lemma then shows that the middle arrow is a bijection.
Since $D^\natural \boxtimes \Pone$ is profinite in its weak topology,
by Lemma~\ref{lattices in Tate modules are complete}, 
and the middle arrow is continuous for this topology, it is therefore a topological isomorphism, as desired.
(To see this continuity, it suffices to check that each $D^\natural \boxtimes \Pone \to D_k^\natural \boxtimes \Pone$ is weakly continuous;
since this map factors through $(D^\natural \boxtimes \Pone)/\fm^k$, it is equivalent to check that it is strongly continuous; 
and this is true because it is a restriction of $D^{\oplus 2} \to D_k^{\oplus 2}$.)
\end{proof}

\subsubsection{Base-change in the finite type context}
We now %
turn to the setting of a %
morphism of finite type $\cO/\varpi^a$-algebras $A \to B$,
and a projective rank two \'etale $(\varphi,\Gamma)$-module $D_A$ over~$\A_A$
with determinant~$\zeta \varepsilon^{-1}$, corresponding to a morphism $\Spec A\to\cX$.
Recall
(Definition~\ref{def:Steinberg locus})
that we have defined a substack $\cX(\St)$ of~$\cX$.

\begin{cor}
  \label{cor:support-of-base-change-ker-coker}Let $A \to B$ be a %
morphism of finite type $\cO/\varpi^a$-algebras,
and let $D_A$ be a projective rank two \'etale $(\varphi,\Gamma)$-module with $A$-coefficients,
with determinant~$\zeta \varepsilon^{-1}$.
Then the kernel and cokernel of the base-change morphism~\eqref{eqn:non-flat base-change map for Pone} 
$$(D_A^{\natural}\boxtimes \Pone)\cotimes_A B \to D_B^{\natural}\boxtimes \Pone$$ are supported on $\Spec B\times_{\cX}\cX(\St)$. %
\end{cor}

\begin{proof}
Extend (as we may) the morphism $A \to B$ to a surjection $A[x_1,\ldots,x_n] \to B$
(for some $n\geq 1$).  Since $A \to A[x_1,\ldots,x_n]$ is flat,
the base-change morphism
$$(D_A^{\natural}\boxtimes\Pone)\cotimes_A A[x_1,\ldots,x_n] \to D_{A[x_1,\ldots,x_n]}^{\natural}
\boxtimes \Pone$$ 
is an isomorphism
(by Lemma~\ref{lem: base change results we can use to prove G-stability}~\eqref{item: flat case}),
and so we see that, in proving the lemma, we may replace $A$ by $A[x_1,\ldots,x_n]$
and $D_A$ by $D_{A[x_1,\ldots,x_n]}$, and hence assume that the morphism
$A\to B$ is furthermore surjective (and so, in particular, finite).

Now let $\cK$, resp.\ $\cC$, denote the kernel, resp.\ cokernel, of~\eqref{eqn:non-flat base-change map for Pone}. 
By Lemma~\ref{lem:non-flat base-change}, \eqref{eqn:non-flat base-change map for Pone} is open, and $\cK$ and $\cC$ are each finite over~$B$.

  We now study the completions of $\cK$ and $\cC$, in order to reduce to the results in the formal case considered above.

  Let $\fn$ be a maximal ideal of~$B$, with preimage $\m$ in~$A$.  Since $A \to \widehat{A}_\m$ is flat, Lemma~\ref{lem: base change results we can use to prove G-stability}~(3) shows that the base-change morphism
$$(D_A^{\natural}\boxtimes \Pone)\cotimes_A \widehat{A}_\m \to
D_{\widehat{A}_\m}^{\natural} \boxtimes \Pone$$ is an isomorphism.  Similarly,
$$(D_B^{\natural}\boxtimes \Pone)\cotimes_B \widehat{B}_\fn \to
D_{\widehat{B}_\fn}^{\natural} \boxtimes \Pone$$ is an isomorphism.  
Thus, tensoring~\eqref{eqn:non-flat base-change map for Pone} with $\widehat{B}_\fn$ over~$B$, we obtain an exact sequence that can be rewritten as
\[0 \to \widehat{\cK}_{\fn} \to (D_{\widehat{A}_\m}^{\natural} \boxtimes \Pone)\cotimes_{\widehat{A}_\m} \widehat{B}_{\fn} \to D_{\widehat{B}_\fn}^{\natural}\boxtimes \Pone \to \widehat{\cC}_{\fn} \to 0.\] 
(Recall that if $M$ is a finite $B$-module with the discrete topology, then $M\cotimes_B \widehat{B}_{\fn} = M\otimes_B \widehat{B}_{\fn} \iso \widehat{M}_{\fn}.$)
Taking into account the isomorphism of Lemma~\ref{lem:compatibility}~\eqref{item:71}, %
we may further replace $D_{\widehat{A}_m}$ and $D_{\widehat{B}_n}$ by their formal completions, obtaining the exact sequence
\[  0 \to \widehat{\cK}_{\fn} \to
  (\widehat{D}_{\widehat{A}_\m}^{\natural} \boxtimes \Pone)\cotimes_{\widehat{A}_\m}
  \widehat{B}_{\fn} 
  \to \widehat{D}_{\widehat{B}_\fn}^{\natural}\boxtimes \Pone \to \widehat{\cC}_{\fn} \to 0.\]
By Lemma~\ref{lem:tensor with finite algebra is automatically complete},
we can rewrite this as
\[  0 \to \widehat{\cK}_{\fn} \to
  (\widehat{D}_{\widehat{A}_\m}^{\natural} \boxtimes \Pone)\otimes_{\widehat{A}_\m}
  \widehat{B}_{\fn} 
  \to \widehat{D}_{\widehat{B}_\fn}^{\natural}\boxtimes \Pone \to \widehat{\cC}_{\fn} \to 0.\]
The result is now immediate from Lemma~\ref{lem:basechange failure criterion}.
\end{proof}

\subsection{Some explicit base-change computations}\label{Colmez functor computations}
In this section we investigate two particular instances of the base-change we studied 
in Section~\ref{subsubsec:formal base change}.
In light of Lemma~\ref{lem:basechange failure criterion},
the interesting cases are those when $D \coloneqq  \bD(V)$,  
for $V$ a free rank $2$ module over the complete Noetherian local ring~$R$
equipped with a continuous action of~$G_{\Q_p}$,
such that furthermore $\overline{V}^{\semis}$ is a twist of $\omega \oplus 1$. 
We will be primarily focused on the case when $\overline{V}$ is in fact
already semisimple, and so (up to a twist) isomorphic to~$\omega \oplus 1$; %
thus we begin with some considerations related to the corresponding versal deformation.

Throughout this section we let~$\thetabar = \omega + 1$, and we will denote the
versal deformation ring $R^{\ver}_{\thetabar}$ from Definition~\ref{defn:various versal rings}, resp.\ the versal object~$\Vver_{\thetabar}$, by~$R$, resp.\ $V$.
Then~$V$ is a versal deformation of $\omega \oplus 1$ with fixed determinant~$\varepsilon$,
and~$R$ is equal to the $\fm$-adic completion
of the ring $S$ introduced in Section~\ref{subsubsec:Steinberg CWE stack},
where $\fm$ is the maximal ideal of $S$ corresponding to the closed point
of~$\cX_{\thetabar}$ given by the representation $\omega \oplus 1$; %
namely $\fm = (a_0,a_1,b_0,b_1,c)$.  Thus
$R = \cO\llbracket a_0,a_1,b_0,b_1,c\rrbracket /(a_0b_1 + a_1 b_0)$ (and the character~$\zeta$ is implicitly fixed to be~$\varepsilon^2$, so that $\det(V) = \zeta \varepsilon^{-1}$).
The natural map $R^{\ps}_{\thetabar} \to R$ from the pseudodeformation ring is injective with image
$R^{\ps}_{\thetabar} = \cO\llbracket a_0,a_1,b_0c, b_1c\rrbracket /
(a_0b_1c+a_1b_0c).$ 
Furthermore, $V \cong R^{\oplus 2}$ endowed with a certain $G_{\Q_p}$-action,
which can be written down essentially explicitly,
see for example~\cite[Section~3]{HuTan}, \cite[Section~3.4]{JNWE} or~\cite[Appendix~B]{MR3150248},
but we won't need that here. 
Rather, we use the fact (see~\eqref{2x2 matrix order} for a discussion) that $\cO\llbracket G_{\Q_p}\rrbracket $ acts on~$V$ 
through its quotient $\tR_{\thetabar}$, and that the action map $\tld R_{\thetabar} \to M_2(R)$ is injective with image
$$
\fourmatrix{R^{\ps}_{\thetabar}}{R^{\ps}_{\thetabar}b_0+R^{\ps}_{\thetabar}b_1}
{R^{\ps}_{\thetabar}c}{R^{\ps}_{\thetabar}}.
$$

\subsubsection{Computing abelian subrepresentations and quotients}\label{abelianization}
Our key tool
for studying base-change problems involving~$V$
and its various specializations
will be the 
exact sequence~\eqref{eqn:non-flat basechange exact sequence V form}.
In order to use it, we will need to compute $V^{\ab}$ and $V_{\ab}$ (and the analogous
objects for the specializations of~$V$).   A convenient way to do this is
in terms of the ``abelianization'' of~$\tR_{\thetabar}$, i.e.\ its maximal
commutative quotient.  One easily verifies that this quotient is equal to
$$
\fourmatrix{R^{\ps}_{\thetabar}/(b_0c,b_1c)}{0}
{0}{R^{\ps}_{\thetabar}/(b_0 c,b_1 c)},
$$
which is the quotient of $\tR_{\thetabar}$ by its ideal
$$ 
J\coloneqq  \fourmatrix{(b_0c,b_1c)}{R^{\ps}_{\thetabar}b_0+R^{\ps}_{\thetabar}b_1}
{R^{\ps}_{\thetabar}c}{(b_0c,b_1 c)}.
$$
Thus, if $S$ is any quotient of~$R$, and $V_S \coloneqq  V\otimes_R S = S^{\oplus 2}$
is the corresponding $G_{\Q_p}$-representation,
then we see that
\begin{equation}
\label{eqn:V^ab formula}
(V_S)^{\ab} = V_S [J] = S[c] \oplus S[(b_0,b_1)]\subseteq S \oplus S, 
\end{equation}
while
\begin{equation}
\label{eqn:V_ab formula}
(V_S)_{\ab} = V_S / J V_S = S/(b_0,b_1)S \oplus S/c S.
\end{equation}
Since $V/\fm V = \omega \oplus 1$,
we see that $G_{\Q_p}$ acts on the first summand of~$(V_S)_{\ab}$ 
through a deformation of~$\omega$, and on the second summand through a deformation
of~$1$.
Since~$G_{\bQ_p}$ acts on both summands through its abelianization, we can regard these as representations of~$\bQ_p^\times$.

For example, if we apply the above discussion to the case of $V$ itself (i.e.\ we take $S = R$),
we obtain the following result.

\begin{lemma}
\label{lem:abelian computations}
\leavevmode
\begin{enumerate}
\item $V^{\ab}  = 0$.
\item $V_{\ab}  = R/(b_0,b_1) \oplus R/(c)$, with $\Qp^{\times}$ acting on
the first summand via
a deformation of~$\omega$,
and on the second summand via
a deformation of~$1$.
\end{enumerate}
\end{lemma} 
\begin{proof}
The first claim follows from~\eqref{eqn:V^ab formula} 
and the fact that $R$ is an integral domain,
so that $R[c] = R[(b_0,b_1)] = 0.$  (Alternatively,
it also follows from Lemma~\ref{no D^ab}.)
The second claim follows immediately from~\eqref{eqn:V_ab formula}. 
\end{proof}

If $I$ denotes an ideal in the versal ring~$R$,
and if we write
$D \coloneqq  \bD(V)$ and $D' \coloneqq  \bD(V/IV),$
then we can study the base-change morphism
\begin{equation}\label{eqn:base-change morphism on versal rings}
(D^{\natural}\boxtimes \Pone)\otimes_R R/I \to D'^{\natural}\boxtimes \Pone.
\end{equation} 
If we let $\cK$ and $\cC$ denote the kernel and cokernel respectively of 
this morphism, then (since, 
as already noted in Lemma~\ref{lem:abelian computations} above,
we have $V^{\ab} = 0$), we deduce
from~\eqref{eqn:non-flat basechange 
exact sequence V form}, together with~\eqref{eqn:V^ab formula} and~\eqref{eqn:V_ab formula},
that $\cK$ and $\cC$ admit the alternate descriptions as the 
kernel and cokernel of a morphism
\begin{equation}
\label{eqn:alternate description}
 \Tor_1^R(V/JV\otimes\varepsilon^{-1}, R/I) \to 
w\cdot (V/IV)[J]; %
\end{equation}
in particular they are finite $R/I$-modules. 

More generally, choose a descending sequence of ideals $I_n \subset I$, write
$$D_n \coloneqq  \bD(V/I_n V) = \bD(V)\otimes_R R/I_n,$$
and
let $\cK_n$ and $\cC_n$ denote the kernel and cokernel of the base-change morphism
\begin{equation}\label{eqn:base-change morphism on versal rings II}
(D_n^{\natural}\boxtimes \Pone) \otimes_{R/I_n} R/I \to  
D'^{\natural}\boxtimes \Pone.
\end{equation}
Then %
$\{\cK_n\}_{n \geq 1}$
and $\{\cC_n\}_{n \geq 1}$
are projective systems of finitely generated $R/I$-modules.
There are evident morphisms $\cK \to \varprojlim_n \cK_n$ and $\cC \to \varprojlim_n \cC_n,$
which are easily seen to be isomorphisms. In the following lemma we make the stronger observation
that, under a natural assumption on~$I_n$, these in fact are isomorphisms in the category $\Pro\Mod^{\fg}(R/I).$

\begin{lemma}%
\label{lem:coherence}%
Maintain the notation of the previous paragraph, and assume further that the ideals~$I_n$ are cofinal with the powers of~$I_1$.
Then the natural maps induce isomorphisms
$\cK \iso \quoteslim{n} \cK_n$ and $\cC \iso \quoteslim{n} \cC_n$ in~$\Pro\Mod^{\fg}(R/I)$.
In particular,
both $\quoteslim{n} \cK_n$ and $\quoteslim{n} \cC_n$, which {\em a priori} are objects
of~$\Pro \Mod^{\fg}(R/I)$, are in fact
isomorphic to objects of~$\Mod^{\fg}(R/I).$
\end{lemma}%
\begin{proof}
Since the pro-objects~$\cK_n$ and~$\cC_n$ arising from~$I_n$ and~$I_1^n$ are isomorphic, we can assume without loss of generality that $I_n = I_1^n$.
Just as in the preceding discussion,
we then apply~\eqref{eqn:non-flat basechange exact sequence V form},
together with~\eqref{eqn:V^ab formula} and~\eqref{eqn:V_ab formula},
to deduce that the $\cK_n$ and $\cC_n$ fit into the projective 
system of exact sequences
\begin{multline*}
w\cdot \ker\Bigl(\bigl((V/I_1^nV)[J]\bigr)\otimes_{R/I_1^n} R/I  \to (V/IV)[J]\Bigr)
\to \cK_n
\\
\to \Tor_1^{R/I_1^n}\Bigl(\bigl((V/I_1^nV)/J(V/I_1^n V)\bigr)\otimes \varepsilon^{-1}, R/I\Bigr)
\\
\to
w\cdot \coker\Bigl(\bigl((V/I_1^nV)[J]\bigr)\otimes_{R/I_1^n} R/I  \to (V/IV)[J]\Bigr)
\to \cC_n \to 0.
\end{multline*}
If we regard this projective system of exact sequences as giving an exact sequence
in the category $\Pro \Mod^{\fg}(R/I),$
then the Artin--Rees lemma shows %
that it is isomorphic in that category to the
exact sequence
\begin{multline}\label{eqn:sequence after Artin--Rees}
w\cdot \ker\bigl((V[J])\otimes_{R} R/I  \to (V/IV)[J]\bigr)
\to \quoteslim{n} \cK_n \to \Tor_1^{R}\bigl((V/JV)\otimes \varepsilon^{-1}, R/I\bigr)
\\
\to
w\cdot \coker\bigl((V[J])\otimes_{R} R/I  \to (V/IV)[J]\bigr)
\to \quoteslim{n} \cC_n \to 0.
\end{multline}
(More precisely, we are using the fact that the functor
\[
\Mod^{\fp}(R) \to \Pro \Mod^{\fp}(R), M \mapsto \quoteslim{n} M/I_1^n M
\]
is fully faithful and exact, which is a standard consequence of the Artin--Rees lemma.)
As already noted in Lemma~\ref{lem:abelian computations} above,
we have $V^{\ab} = V[J] = 0,$ so that~\eqref{eqn:sequence after Artin--Rees} simplifies to
\begin{multline*}
0 \to \quoteslim{n} \cK_n \to \Tor_1^{R}\bigl((V/JV)\otimes \varepsilon^{-1}, R/I\bigr)
\to
w\cdot  (V/IV)[J]
\to \quoteslim{n} \cC_n \to 0.
\end{multline*}

Since the inner two terms are in fact objects of $\Mod^{\fg}(R/I)$,
 the same is true of the outer two terms, as claimed. 
Furthermore, taking into account the description of~$\cK$ and $\cC$
as the kernel and cokernel of~\eqref{eqn:alternate description}
(and the evident naturality of the formation
of~\eqref{eqn:non-flat basechange exact sequence V form}),
we obtain the claimed isomorphisms.
\end{proof}

In the remainder of this section we consider the ``extreme case'' of base-change along
$R \to R/\m_R \iso \F$. %
Thus,
again writing $D \coloneqq  \bD(V)$, and also writing
$$\overline{D} \coloneqq  \bD(\overline{V}) = \bD(\omega\oplus 1),$$
we will study the base-change morphism
\begin{equation}
\label{eqn:versal to residual basechange}
(D^{\natural}\boxtimes \Pone)\otimes_R \F \to \overline{D}^{\natural}\boxtimes \Pone.
\end{equation}
\begin{lemma}
\label{lem:describing Dbar boxtimes Pone}
Let~$\pi_\alpha := \Ind_B^G(\omega \otimes \omega^{-1})$.
Then we have an isomorphism
$$\overline{D}^{\natural}\boxtimes \Pone 
\iso 
\pi_\alpha^\vee \oplus (\triv_G - \St)^\vee,$$
where $(\triv_G -\St)^\vee$
is the Pontrjagin dual to the unique non-split extension of~$\St$ by~$\triv_G$.
\end{lemma}
\begin{proof}
This follows
from~\cite[Section~VII.4.7]{MR2642409}, where it is shown that 
$$\bD(1)^\natural \boxtimes \bP^1 \iso (\triv_G - \St)^\vee$$ 
and 
\[\bD(\omega)^\natural \boxtimes \bP^1 \iso \pi_\alpha^\vee.\]
(Note that, since we have $\zeta = \varepsilon^2$ in this subsection, our symbol~$\boxtimes$ coincides with the~$\boxtimes$ of~\emph{loc.\ cit.})
\end{proof}
As usual, we denote the  kernel and cokernel of~\eqref{eqn:versal to residual basechange} respectively by~$\cK$ and~$\cC$.
Taking  %
Lemma~\ref{lem:abelian computations} into account,
we see that
the exact sequence~\eqref{eqn:non-flat basechange exact sequence V form} of $B$-representations
associated to the base-change morphism~\eqref{eqn:versal to residual basechange}
may be written as
\begin{equation}\label{special case of base change sequence}
0 \to \cK \to (1 \otimes 1)^{\oplus 2} \oplus (\omega^{-1}\otimes \omega)
\to (1 \otimes 1) \oplus (\omega^{-1} \otimes \omega) \to \cC \to 0.
\end{equation}
(Here we have used the fact that the $\F$-dimension of $\Tor_1^R(R/I,\F)$ is equal
to the minimal number of generators of~$I$, for any ideal $I\subseteq R$.)

Since, by Lemma~\ref{trivial action}, the $G$-action on each of $\cK$ and $\cC$ (and hence also the $B$-action)
has to be a direct sum of representations of the form $\beta\circ \det$, we see that in fact
either
$$\cK \iso \triv_G^{\vee} \quad \text{ and }\quad  \cC = 0$$
or
$$\cK \iso (\triv_G^{\vee})^{\oplus 2} \quad  \text{ and }\quad  \cC \iso \triv_G^{\vee}.$$
We write~$\triv_G^\vee$ here, rather than~$1$ as in~\eqref{special case of base change sequence}, 
to emphasize the fact that here we regard~$\cK$ and~$\cC$ as $G$-representations, as opposed to~$G_{\bQ_p}$-representations;
in other words,
since~$\cK$ and~$\cC$ are the kernel and cokernel of~\eqref{eqn:versal to residual basechange}, which is a morphism in~$\fC_{\thetabar}$, we consider them as objects of~$\fC_{\thetabar}$.

Since $\Ext^1_{\cA}(\triv_G, \pi_\alpha) = 0$ by~\cite[(165)]{MR3150248},
we deduce from Lemma~\ref{lem:describing Dbar boxtimes Pone}
together with the preceding computations of $\cK$ and $\cC$ that 
\begin{equation}
\label{eqn:preliminary description}
(D^{\natural}\boxtimes \Pone)\otimes_R \F \iso \pi_\alpha^\vee \oplus \mathbf{E}^{\vee},
\end{equation}
where $\mathbf{E}$ is either a three-step extension of the form $(\triv_G -\St-\triv_G)$
or else an extension of the form $(\St - \triv_G^{\oplus 2})$,
depending on which of the two possibilities for the structure of $\cK$ and $\cC$ holds.

We will show that in fact the second case is the one that holds, and, furthermore,
that $\mathbf{E}$ is the universal extension $E(\St)$ of $\triv_G$ by~$\St$.
We will do this by finding certain $\cO$-specializations $V'$ of $V$
such that the induced map %
$(D^{\natural}\boxtimes \Pone)\otimes_R \F \to (D(V')^{\natural}\boxtimes \Pone)\otimes_\cO \F $ 
is surjective, and such that the  target is a successive extension of $\pi_{\alpha}^{\vee}$
by $\St^{\vee}$ by~$\triv_G^{\vee}$, the latter extension having been arbitrarily
chosen in advance. 

\subsubsection{Construction of specializations}
Let $\tau \ne 0 \in \Hom_\bZ(\bQ_p^\times, \bF)$, and let~$\tau^\perp$ be its orthogonal line in $H^1(G_{\bQ_p}, \bF(\omega))$ under local Tate duality.
We write~$\rhobar_\tau$ for the non-split extension
\[
0 \to \omega \to \rhobar_\tau \to 1 \to 0
\]
corresponding to~$\tau^\perp$, and 
we choose a lift $\rho_\tau : G_{\bQ_p} \to \GL_2(\cO)$ of~$\rhobar_\tau$ with determinant~$\varepsilon$ such that~$\rho_\tau[1/p]$ is absolutely irreducible.
To construct this lift, we may need to replace~$E$ %
with a finite extension; since
we ultimately need to construct such lifts just for the two members of some chosen
$\F$-basis of~$\Hom_{\bZ}(\bQ_p^{\times}, \bF)$, we make such a replacement once and for all.
The existence of the lift is then a very special case of~\cite[Prop.\ 1.12]{BIP}, but it is easily seen to exist by a direct argument; 
indeed~$\rhobar_{\tau}$ is unobstructed, so the corresponding universal deformation ring is formally smooth, and the reducible locus is formally smooth of codimension 1.
We write~$\lbar D_\tau \coloneqq \bD(\rhobar_\tau)$ and $D_\tau \coloneqq \bD(\rho_\tau)$.

\begin{lemma}\label{Dtaunatural}
The natural map $(D_\tau^\natural \boxtimes \bP^1) \otimes_{\cO} \bF \to \lbar D_\tau^\natural \boxtimes \bP^1$ is an isomorphism.
\end{lemma}
\begin{proof}
This follows by a consideration of~\eqref{eqn:non-flat basechange exact sequence V form}
(with the map of local rings taken to be $\cO \to \F$).
Indeed, taking into account Lemma~\ref{no D^ab}, and noting that
$(\rho_{\tau})_{\ab} \otimes_{\cO}\F \iso (\rhobar_{\tau})_{\ab}\iso 1,$
we see that this sequence simplifies to
$$ 0 \to  \cK \to \omega^{-1} \otimes \omega \to \omega^{-1}\otimes\omega \to \cC \to 0,$$
from which (recalling, as always, that the action of $G$ on $\cK$ and $\cC$ has to be
via representations of the form $\beta \circ \det$) we see that $\cK = \cC = 0,$
as required.
\end{proof}

Recall from~\cite[Prop.\ VII.4.27]{MR2642409} (and the discussion preceding
the statement of that result)  that we have a short exact sequence
\[
0 \to
(\triv_G - \pi_\alpha)^\vee 
\to \lbar D_\tau^\natural \boxtimes \bP^1 \to \St^\vee \to 0,
\]
where 
$(\triv_G - \pi_\alpha)^\vee$ denotes the dual to unique non-split extension of  
$\pi_{\alpha}$ by~$\triv_G.$
Indeed, in the notation of~\emph{loc.\ cit.\ }, we have $D_1 = \bD(\omega)$ and~$D_2 = \bD(1)$,
$D_1^{\sharp}\boxtimes \Pone =
(\triv_G - \pi_\alpha)^\vee,$ 
and $(D_2^\natural \boxtimes \bP^1)_0 = \St^\vee$.
Write
\begin{equation}
  \label{eq:E-tau-defn}
  E_{\tau}^{\vee}\coloneqq (\lbar D_\tau^\natural \boxtimes \bP^1)/\pi_\alpha^\vee,
\end{equation}
so that 
 $E_{\tau}^{\vee}$ is an extension of $\St^{\vee}$ by~$\triv^{\vee}_G$, dual to an extension
$E_{\tau}$ of $\triv_G$ by~$\St$.
By ~\cite[Thm.\ VII.4.18, Prop.\ VII.4.27]{MR2642409}, the map $E_\tau\mapsto\tau$  is 
an $\bF$-linear isomorphism
\begin{equation}\label{tau isomorphism}
\Ext^1_{\cA}(\triv_G, \St) \isom \Hom_\bZ(\bQ_p^\times, \bF).
\end{equation}
Of course, the representations $\rho_{\tau}$ are not specializations of our versal representation~$V$,
since they have the wrong residual representations.
However, after enlarging $E$ if necessary (e.g.\ replacing $E$ by $E(\varpi^{1/2})$ suffices),
we may find a $G_{\Q_p}$-invariant lattice $\rho_{\semis} \subset \rho_{\tau}$
with reduction equal to $\omega \oplus 1,$
and such that 
\begin{equation}\label{adjacent lattices}
\varpi \rho_{\mss} \subset \varpi \rho_\tau \subset \rho_{\mss} \subset \rho_\tau.
\end{equation}

\begin{lemma}
\label{lem:rhoss reduction}
There is an isomorphism
$$(D(\rho_{\semis})\boxtimes \Pone)\otimes_{\cO} \F
\iso \pi_{\alpha}^{\vee} \oplus E_{\tau}^{\vee}.$$ 
\end{lemma}
\begin{proof}
We apply~\eqref{eqn:non-flat basechange exact sequence V form} to the base-change
of $D(\rho_{\semis})\boxtimes \Pone$ along the surjection $\cO\to \F$.
Taking into account Lemma~\ref{no D^ab}, and the fact
that $(\rho_{\semis})_{\ab} \otimes_{\cO} \F = \omega \oplus 1,$
we see that this sequence simplifies to 
$$0 \to \cK \to
(1 \otimes 1) \oplus (\omega^{-1}\otimes \omega) 
\to 
(1 \otimes 1) \oplus (\omega^{-1}\otimes \omega) 
\to \cC \to 0.
$$
Thus either $\cK = \cC = 0$, or $\cK = \cC = \triv_G^\vee.$
If we also take into account Lemma~\ref{lem:describing Dbar boxtimes Pone},
we conclude that, correspondingly, either 
$$(D(\rho_{\semis})\boxtimes \Pone)\otimes_{\cO} \F \iso 
\pi_{\alpha}^{\vee} \oplus (\triv_G - \St)^{\vee},$$
or else that
$$(D(\rho_{\semis})\boxtimes \Pone)\otimes_{\cO} \F \iso 
\pi_{\alpha}^{\vee} \oplus (\St - \triv_G )^{\vee},$$
for some extension $(\St - \triv_G)$ of $\triv_G$ by~$\St$.

To conclude the proof, it therefore suffices to show that there is an embedding
$$E_{\tau}^{\vee} \hookrightarrow 
\bigl(D(\rho_{\semis})\boxtimes \Pone\bigr)\otimes_{\cO} \F.$$
To this end, observe first that the functor $D \mapsto D^\natural \boxtimes \Pone$ preserves injections: 
in fact, inspecting Definition~\ref{def: various functors on formal phi Gamma modules} we see that 
there are natural injections $D^\natural \boxtimes \Pone \to D \boxtimes \Pone \to D^{\oplus 2}$,
and~$D \mapsto D^{\oplus 2}$ is exact.

We thus see that $D_\tau^\natural \boxtimes \Pone$ is $\varpi$-torsion free, and applying $(\text{--})^\natural \boxtimes \Pone$ to~\eqref{adjacent lattices}, we obtain a chain of inclusions
\[
\varpi (D^\natural_{\mss} \boxtimes \bP^1) \subset \varpi(D^\natural_{\tau} \boxtimes \bP^1) \subset D^\natural_{\mss} \boxtimes \bP^1 \subset D^\natural_{\tau} \boxtimes \bP^1,
\]
from which we deduce that $(D_{\mss}^\natural \boxtimes \bP^1) \otimes_{\cO} \bF$ has a submodule isomorphic to
\[
\coker\bigl((D^\natural_{\mss} \boxtimes \bP^1) \otimes_\cO \bF \to (D^\natural_{\tau} \boxtimes \bP^1) \otimes_\cO \bF\bigr).
\]
It therefore suffices to prove that this cokernel is isomorphic to~$E_\tau^{\vee}$.
By Lemma~\ref{Dtaunatural} and~\eqref{eq:E-tau-defn} (i.e.\ by the definition of~$E_{\tau}$), it suffices to prove that
\[
\operatorname{image}\bigl((D^\natural_{\mss} \boxtimes \bP^1) \otimes_\cO \bF \to (D^\natural_{\tau} \boxtimes \bP^1) \otimes_\cO \bF \iso \lbar D^\natural_\tau \boxtimes \bP^1 \bigr) = \pi_\alpha^\vee.
\]
We have a commutative diagram
\[
\begin{tikzcd}
(D^\natural_{\mss} \boxtimes \bP^1) \otimes_\cO \bF \arrow[r] \arrow[d] & (D^\natural_{\tau} \boxtimes \bP^1) \otimes_\cO \bF \arrow[d, "\sim"]\\
\lbar D^\natural_{\mss} \boxtimes \bP^1 \arrow[r] & \lbar D^\natural_\tau \boxtimes \bP^1
\end{tikzcd}
\]
where the lower horizontal arrow is induced by a non-zero map $\omega \oplus 1 \to \rhobar_\tau$, and so has image~$\bD(\omega)^\natural \boxtimes \bP^1 = \pi_\alpha^\vee$.
This concludes the proof, since the cokernel $\cC$ of the left vertical arrow is finite-dimensional.
\end{proof}

\begin{defn}%
Recall that
$\dim_\bF \Ext_{\cA}^1(\triv_G, \St) = 2$
(e.g.\ by~\cite[Section~10.1]{MR3150248}, or~\eqref{tau isomorphism} above).
The universal extension $E(\St)$ is defined to be the pushout along~$\St$ of any two linearly independent extensions of~$\triv_G$ by~$\St$. 
Alternatively, we may describe $E(\St)$ as the preimage in
$\operatorname{inj}_{\bF[G/Z]}(\St)$ of the $\triv_G$-isotypic component of the
first layer of the socle filtration.  
It sits in a short exact sequence
$$0 \to \St \to E(\St) \to \triv_G^{\oplus 2} \to 0.$$
Thus its Pontrjagin dual sits in a short exact sequence
$$0 \to (\triv_G^{\vee})^{\oplus 2} \to E(\St)^{\vee} \to \St^{\vee} \to 0.$$
\end{defn}
\begin{thm}\label{special fibre of versal deformation}
We have an %
isomorphism
\[
(D^\natural \boxtimes \bP^1) \otimes_{R} \bF \cong \pi_\alpha^\vee \oplus E(\St)^\vee.
\]
\end{thm}

\begin{proof}
Let $\rho_{\semis}$ be one of the lattices considered in Lemma~\ref{lem:rhoss reduction}.
Then $\rho_{\semis}$ is obtained by base-changing
$V$ along some $\cO$-algebra morphism $R \to \cO$. 
Considering the exact sequence~\eqref{eqn:non-flat basechange exact sequence V form}
for this base-change,
and again taking into account Lemma~\ref{no D^ab} (which applies to both $V$
and~$\rho_{\semis}$), we deduce that the induced map
$$D\boxtimes\Pone \to  D(\rho_{\semis})\boxtimes\Pone$$
is surjective.
Consequently
$$(D\boxtimes\Pone)\otimes_R \F  \to \bigl( D(\rho_{\semis})\boxtimes\Pone\bigr)\otimes_{\cO} \F$$
is also surjective. 
Lemma~\ref{lem:rhoss reduction}
shows that the target has a direct summand isomorphic to~$E_{\tau}^{\vee}$. %
Considering~\eqref{eqn:preliminary description},
we deduce that the direct summand~$\bE^\vee$ of the source surjects onto~$E_{\tau}^\vee$.
Letting $\tau$ range over a pair of basis vectors of $\Hom_{\Z}(\Q_p^{\times},\F),$
we find that in fact~$\mathbf{E}^{\vee}$ appearing 
in~\eqref{eqn:preliminary description} must be the dual to the universal extension~$E(\St)$,
as required.
\end{proof}

In fact, in the sequel, our application of Theorem~\ref{special fibre of versal deformation}
will be via the following corollary, which does not use its full strength,
but only its particular consequence %
\begin{equation}\label{no trivial in cosocle}
\Hom_{\fC_{\thetabar}}((D^\natural \boxtimes \bP^1)\otimes_R \bF, \triv_G^\vee) = 0.
\end{equation}

\begin{cor}
\label{cor:Steinberg D natural}Suppose that $i_S: \Spf S\to\cX_{\thetabar}$ is a
versal morphism at the closed point of~$\cX_{\thetabar}$, where~$S$ is a
complete local Noetherian $\cO$-algebra with finite residue field. Write~$D_{S}$
for the corresponding formal \'etale $(\varphi,\Gamma)$-module with
$S$-coefficients (i.e.\ $D_S \coloneqq  \bD(\ihat^*_S \cV_{\thetabar})$). 
Then $D^{\natural}_S\boxtimes \Pone$ has no
$\fC_{\thetabar}$-quotients of $\cO$-finite length.
\end{cor}
\begin{proof}
  Without loss of generality (\emph{cf.}\ \cite[\href{https://stacks.math.columbia.edu/tag/06T5}{Tag 06T5}]{stacks-project}) %
  we can and do suppose that~$S=R\llbracket x_1,\dots,x_r\rrbracket $ is a power series ring over~$R$, %
  where~$R$ is the versal ring considered above. %
If
\[
\Hom_{\fC_{\thetabar}}(D_S^\natural \boxtimes \bP^1, \triv_G^\vee) \ne 0, 
\]
then there exists~$k > 0$ such that
\[
\Hom_{\fC_{\thetabar}}\bigl((D_S^\natural \boxtimes \bP^1) \otimes_{S} S/\fm_{S}^k, \triv_G^\vee\bigr) \ne 0.
\]
(This is %
because $D_S^\natural \boxtimes \bP^1$ is a profinite $S$-module, hence has the $\fm_S$-adic topology by Lemma~\ref{properties of compact modules}~\eqref{item: compact 7}, 
and we are working with continuous homomorphisms.)
Filtering $S$ by powers of the maximal ideal, we deduce that
\begin{equation}\label{cosocle contradiction}
\Hom_{\fC_{\thetabar}}\bigl((D_S^\natural \boxtimes \bP^1) \otimes_{S} \bF, \triv_G^\vee\bigr) \ne 0.
\end{equation}
Writing $S_n \coloneqq  (R/\m_R^{n})[x_1,\dots,x_r]/(x_1 ,\dots,x_r)^{n}$, we have $S=\varprojlim_{n} S_n$, and by definition
\[
D_S^\natural \boxtimes \Pone = \varprojlim_n D_{S_n}^\natural \boxtimes \Pone.
\]
Applying part~(3) of Lemma~\ref{lem: base change results we can use to prove G-stability} 
to the finite free morphism $R/\fm_R^n \to S_n$, and using Lemma~\ref{lem:tensor with finite algebra is automatically complete} to replace~$\cotimes$ with~$\otimes$, we 
find
that
\[
D_S^\natural \boxtimes \bP^1 \cong \varprojlim_n (D_{R/\fm_R^n}^\natural \boxtimes \Pone) \otimes_{R/\fm_R^n} S_n.
\]
Since~$\bF$ is finitely presented over~$R$ and~$S$, and cofiltered limits are exact in~$\fC_{\thetabar}$, we conclude that
\[
(D_S^\natural \boxtimes \bP^1) \otimes_{S} \bF \cong \varprojlim_n ((D_{R/\fm_R^n}^\natural \boxtimes \Pone) \otimes_{R/\fm_R^n} S_n \otimes_S \bF) \cong (D^\natural \boxtimes \bP^1) \otimes_{R} \bF.
\]
Then~\eqref{cosocle contradiction} contradicts~\eqref{no trivial in cosocle}, and we
are done.
\end{proof}

\section{The functor}\label{sec: the
  functor}
We now construct our functor
$\Functor: D^b_{\fp}(\cA) \to D^b_{\coh}(\cX)$,
and prove our main theorem, i.e.\ that $\Functor$ is fully faithful.

\subsection{The definition of the functor}\label{subsec: the functor}%
In this section we construct the pro-coherent sheaf $L_{\infty}$ of $\cO\llbracket G\rrbracket _{\zeta}$-modules
over~$\cX$, and use it to define our functor~$\Functor$.

\subsubsection{Interpreting linearly topological modules as pro-coherent sheaves}
We begin by explaining a general procedure that converts certain topological modules
into pro-coherent sheaves.
We have the following (slightly {\em ad hoc}) definition.

\begin{df}
Let $A$ be a Noetherian ring endowed with the discrete topology.  We say that a linearly topological $A$-module $M$ 
is {\em pro-coherent} if $M$ is complete,
and if for any open $A$-submodule $U$ of~$M$, the quotient $M/U$ is finitely
generated over~$A$.
(Of course, it suffices to check this condition for $U$ running over a neighbourhood
basis of $0$.)
Note that the natural morphism $M \to \varprojlim_U M/U$ (the projective
limit being taken over all open submodules of~$M$, where each quotient is
endowed with its discrete topology)   
is then a topological isomorphism.  
\end{df}

\begin{remark}
Any lattice in an $A$-Tate module is pro-coherent when equipped with its induced
topology.  
\end{remark}

For any finitely generated $A$-module~$N$, we write $\widetilde{N}$ to denote
the coherent sheaf obtained by localizing $N$ over~$\Spec A$.

\begin{df} 
If $M$ is a pro-coherent linearly topological 
$A$-module, then we write
$\widetilde{M} \coloneqq  \quoteslim{U} \widetilde{M/U}$,
regarded as a pro-coherent sheaf on~$\Spec A$, i.e.\ an object of $\Pro \Coh(\Spec A)$.
\end{df}

\begin{remark}\label{rem: tilde construction on finite discrete modules}
Any finitely generated $A$-module may be regarded as a pro-coherent
$A$-module by endowing it with the discrete topology; in this case the two
possible meanings of $\widetilde{N}$ evidently coincide.
\end{remark}

\begin{remark}\label{rem: tilde construction on A[[K]]-modules}
If~$A$ is a finite type $\cO/\varpi^a$-algebra, and~$M$ is a finitely presented $A\llbracket K \rrbracket$-module, the canonical completely metrizable topology on~$M$ coincides with 
its $\fa$-adic topology, and so it is pro-coherent.
The functor $\Mod^{\fp}(A\llbracket K \rrbracket) \to \Pro \Coh(\Spec A), M \mapsto \tld M$ is isomorphic to $M \mapsto \quoteslim{n} M/\fa^n M$, and so it is exact, by Lemma~\ref{lem:topology-Iwasawa-algebra}~\eqref{item: Iwasawa 6}. 
\end{remark}

\begin{remark}
\label{rem:pro-coherent warning}
Not every pro-coherent sheaf $\quoteslim{i} \cF_i$ on $\Spec A$
need be of the form~$\widetilde{M}$ for some pro-coherent linearly topological 
$A$-module, even if the transition morphisms $\cF_{i'}\to \cF_i$ are surjective.
Indeed, if we let $M \coloneqq  \varprojlim_i \cF_i(A),$ then  
the morphisms $M \to \cF_i(A)$ need not be surjective (see e.g.~\cite[\href{https://stacks.math.columbia.edu/tag/0ANX}{Tag 0ANX}]{stacks-project}).

On the other hand,
if the indexing set $I$ is countable and the transition morphisms are surjective,
then $M$ is indeed a pro-coherent module,
and the natural morphism $\widetilde{M} \to \quoteslim{i} \cF_i$ 
is an isomorphism.
\end{remark}

The following lemma shows that the formation of~$\widetilde{M}$ gives rise to a fully
faithful functor with certain natural exactness properties.

\begin{lemma}\label{lem: properties of tilde construction}
\leavevmode
The formation of $\widetilde{M}$ induces a fully faithful $A$-linear functor from 
the $A$-linear additive category of 
pro-coherent linearly topological~$A$-modules {\em (}whose morphisms
are continuous $A$-linear morphisms{\em )} to the $A$-linear abelian category 
$\ProCoh(\Spec A)$.
Furthermore: 
\begin{enumerate}
\item If $f: M \to N$ is a topological embedding of pro-coherent linearly topological
$A$-modules, then the induced morphism
$\widetilde{M} \to \widetilde{N}$ is a monomorphism.
\item If $f:M \to N$ is a
continuous, open, and surjective morphism of pro-coherent linearly topological
$A$-modules, then the induced morphism 
$\widetilde{M} \to \widetilde{N}$ is an epimorphism.
\end{enumerate}
\end{lemma}
\begin{proof}
Let $M$ and $N$ be two pro-coherent linearly topological $A$-modules.
Then
\begin{multline*}
\Hom^{\cont}_A(M,N) \iso \varprojlim_V \Hom^{\cont}_A(M, N/V)
\\
\iso \varprojlim_V \varinjlim_U \Hom_A(M/U,N/V)
\iso \varprojlim_V \varinjlim_U
\Hom_{\Coh(\Spec A)}(\widetilde{M/U},\widetilde{N/V})
\end{multline*}
where $U$ (resp.\ $V$) runs over the open submodules of~$M$ (resp,\ $N$):
the first isomorphism follows from the isomorphism $N\iso \varprojlim_V N/V$,
 the second follows from the fact that $N/V$ is discrete, so that any continuous 
morphism
$M \to N/V$ factors through $M/U$ for some open submodule $U$ of~$M$,
and the third follows from the fact that localization induces an equivalence of categories 
$\Mod^{\fp}(A) \iso \Coh(\Spec A)$.
The right-most expression in this sequence of isomorphisms
is precisely equal to $\Hom_{\ProCoh(\Spec A)}(\widetilde{M},\widetilde{N});$
thus we do indeed obtain a fully faithful functor. 

If $M \hookrightarrow N$ is a topological embedding, then as $V$ runs over the
open submodules of $N$, the intersections $M\cap V$ run through a neighborhood basis at~$0$ consisting of open submodules
of~$M$.
Thus the induced morphism $\widetilde{M} \to \widetilde{N}$ 
arises as the formal inverse limit of the inverse system 
of injections $M/(M\cap V)
 \hookrightarrow N/V.$ This proves~(1).

If $M \to N$ is continuous, open, and surjective, then every open submodule of~$N$ is the image $\overline{U}$ %
of an open submodule $U$ of~$M$. %
In this case the induced morphism $\widetilde{M} \to \widetilde{N}$ corresponds to the formal inverse limit of the inverse system of surjections
$M/U \to N/\overline{U},$ proving~(2).
\end{proof}

In the next lemma,
if $f:\Spec B \to \Spec A$ is a morphism of affine schemes,
then
we continue to write~$f^*$ for the Pro-extension of the usual pullback $f^* : \Coh(\Spec A) \to \Coh(\Spec B)$.

\begin{lemma}
\label{lem:pro-coherent base-change}
Suppose that $M$ is a pro-coherent and first countable 
linearly topological module over the Noetherian ring~$A$.
If $B$ is a Noetherian $A$-algebra,
and $f: \Spec B \to \Spec A$ is the corresponding morphism of schemes, 
then $f^* \widetilde{M} 
\iso (M\cotimes_A B)^{\sim}.$
\end{lemma}
\begin{proof}
Let $\{U_n\}_{n \geq 0}$ be a cofinal sequence of open submodules of~$M$, which
exists by virtue of our assumption that $M$ is first countable. 
Then $M\cotimes_A B = \varprojlim_n (M/U_n)\otimes_A B$
by~\eqref{eqn: recognizing completed tensor};
to ease notation, let $N$ denote this (first countable, pro-coherent) 
linearly topological $B$-module.
Since the transition morphisms $(M/U_{n+1})\otimes_A B \to (M/U_n)\otimes_A B$
are surjective, we see that each of the canonical morphisms $N \to (M/U_n)\otimes_A B$
is surjective, and that if we let $V_n$ denote the kernel of this morphism,
then $V_n$ forms a cofinal sequence of open submodules of~$N$.

Now $\widetilde{M} \coloneqq  \quoteslim{n} \widetilde{M/U_n},$
so that
$$f^* \widetilde{M} = \quoteslim{n} f^* \widetilde{M/U_n}
= \quoteslim{n} \widetilde{\bigl( (M/U_n)\otimes_A B \bigr)}
 \iso \quoteslim{n} \widetilde{ N/V_n} \eqcolon \widetilde{N},
$$
as claimed.
\end{proof}

If $L$ is a lattice in a Tate module over a complete Noetherian local $\cO$-algebra,
then $L$ is equipped not only with its Tate-module topology,
but also with its weak topology,
as defined in Section~\ref{subsubsec:weak topology}.
We let $L^w$ denote $L$ equipped with its weak topology;
this is again a pro-coherent topology on~$L$,
and the aim of the following lemma is to give a description of the associated
pro-coherent sheaf.

\begin{lem}\label{lem: change of topology}
Let~$S$ be a complete Noetherian local $\cO$-algebra, 
and let~$L$ be a lattice in a Tate $S$-module. %
Then, under the identification $\Pro \Mod^{\fp}(S) \isoto \Pro \Coh(\Spec S)$, $\widetilde{L^w}$ is an object of $\Pro \Mod^{\fl}(S) \subset \Pro \Mod^{\fp}(S)$,
and it is isomorphic to the image of $\widetilde L$ under the map
\[
\Pro \Mod^{\fp}(S) \to \Pro \Mod^{\fl}(S)
\]
induced by~\eqref{eqn:mod-fp-to-mod-c} and~\eqref{improved compact modules and Pro}.
\end{lem}
\begin{proof}
By Lemma~\ref{lattices in Tate modules are complete}, $L^w$ is an object of~$\Mod_c(S)$, and so every open $S$-submodule of~$L^w$ has cofinite $S$-length.
By definition, this shows that $\widetilde {L^w }\in \Pro \Mod^{\fl}(S)$. %
This proves the first claim.

For the second claim, choose a neighborhood basis of the origin~$U_i \subset L$ consisting of open submodules such that $L/U_i$ is $S$-finite, for all~$i$.
Then, by definition, $\tld L \in \Pro \Mod^{\fp}(S)$ is $\quoteslim{i}L/U_i$, and so its image in $\Pro \Mod^{\fl}(S)$ is
$\quoteslim{i, n}L/(U_i, \fm^n L)$.
Since $(U_i, \fm^n L)$ forms a basis of the weak topology on~$L$, we see that~$\widetilde{L^w}$ is isomorphic to the image of~$\widetilde L$ in~$\Pro \Mod^{\fl}(S)$.
\end{proof}

\subsubsection{Interpreting \texorpdfstring{$D^{\natural}\boxtimes\Pone$}{D natural box P1} as a pro-coherent
  sheaf}%
\label{subsubsec:D natural box P one as a pro-coherent sheaf}
Suppose now that $A$
is a Noetherian %
$\cO/\varpi^a$-algebra for some~$a\ge 1$, and that~$D$ is an
   \'etale $(\varphi,\Gamma)$-module of rank~$2$ with
  $A$-coefficients and determinant $\zeta\varepsilon^{-1}$. Then we have seen (in Theorem~\ref{thm: D natural is G-stable}) 
  that $D^\natural\boxtimes\Pone$ is an
  $\cO\llbracket G\rrbracket $-stable lattice in $D\boxtimes\Pone$.
In particular, it is pro-coherent 
as a topological $A$-module, 
and so we may form the associated pro-coherent sheaf~$\widetilde{D^{\natural}\boxtimes\Pone}.$
Since $\cO\llbracket G\rrbracket $
acts on $D^{\natural}\boxtimes\Pone$
by continuous endomorphisms (by Corollary~\ref{cor: O-G  action on G box P1}),
we obtain an induced action of~$\cO\llbracket G\rrbracket $ on~$\widetilde{D^{\natural}\boxtimes\Pone},$
making it a left $\cO\llbracket G\rrbracket $-module object in $\ProCoh(\Spec A)$.

We intend %
to apply the construction of
$\widetilde{D^{\natural}\boxtimes \Pone}$
 to the universal case, i.e.\ over~$\cX$.  Since
$\cX$ is only a formal algebraic stack, rather than a scheme, %
we will do this via descent.
We first describe the construction over algebraic stacks, 
before moving to the case of the formal algebraic stack~$\cX$. 

\subsubsection{Descent for algebraic stacks.}\label{subsubsec: descent for algebraic stacks}
Suppose given a Noetherian algebraic $\cO/\varpi^a$-stack~$\cZ$ (for some $a \geq 1$) with affine diagonal,
equipped with a morphism $\cZ \to \cX$.  
Pulling back the universal rank $2$ \'etale $(\varphi,\Gamma)$-module on~$\cX$ to~$\cZ$
gives rise to a rank $2$ \'etale $(\varphi,\Gamma)$-module $D_{\cZ}$ on~$\cZ$,
and we wish to construct the associated pro-coherent sheaf of~$\cO\llbracket G\rrbracket $-modules
$\widetilde{D_{\cZ}^{\natural} \boxtimes \Pone}$ on~$\cZ$.

Let $\Spec A \to \cZ$ be a flat surjection of finite presentation
from the spectrum of a Noetherian $\cO/\varpi^a$-algebra.
Write $R \coloneqq  \Spec A \times_\cZ \Spec A$.
Then~$R$ is an affine scheme with the structure of an {\em fppf}  groupoid over~$\Spec A$.
The morphism $\cZ\to\cX$ then corresponds to  %
a rank 2 \'etale $(\varphi, \Gamma)$-module~$D$ with $A$-coefficients and determinant~$\zeta \varepsilon^{-1}$, endowed with descent data with respect to~$R$.
Concretely, if we let $s,t: R \rightrightarrows \Spec A$ denote the two projections,
and if we let 
$s^*D,$ respectively $t^*D$,
denote the base-changes of $D$, computed with respect to the corresponding (flat) morphisms $A 
\rightrightarrows B$ (where we write $R = \Spec B$),
then we have an isomorphism
$s^*D \iso t^*D$ of \'etale $(\varphi,\Gamma)$-modules over $B$, satisfying an appropriate
cocycle condition. 

Lemma~\ref{lem: base change for D box P1} then shows
that $D\boxtimes \Pone$ is again equipped with descent data with respect to~$R$,
and Lemma~\ref{lem: base change results we can use to prove G-stability}~\eqref{item: flat case}
shows that these descent data restrict to descent data
on the $\cO\llbracket G\rrbracket $-stable sublattice~$D^{\natural}\boxtimes\Pone.$
(We will write that $D^{\natural}\boxtimes\Pone$ is {\em preserved} by the descent data.)
These descent data are necessarily~$G$-equivariant, since the $G$-action
on $D^{\natural}\boxtimes \Pone$ is computed in terms of the $(\varphi,\Gamma)$-module
structure on~$D$, and the descent data on $D^{\natural}\boxtimes \Pone$ are induced
by descent data of $(\varphi,\Gamma)$-modules on~$D$.
Since $\cO[K]$ is dense in~$\cO\llbracket K\rrbracket $, and acts continuously on $D^{\natural}\boxtimes \Pone$,
we find that the descent data are furthermore~$\cO\llbracket K\rrbracket $-equivariant,
and hence also~$\cO\llbracket G\rrbracket $-equivariant. 

\begin{lemma}\label{lem:equivariant lattices for algebraic stacks}
In the context of the preceding paragraph, there exists a neighborhood basis $\fM_n \subset D^\natural \boxtimes \Pone$ of the origin
consisting of $A\llbracket K \rrbracket$-stable lattices that are preserved by the descent data. %
\end{lemma}
\begin{proof}
Assume first that there is a neighborhood basis~$\fN_n$ of the origin consisting of lattices that are preserved by the descent data.
Since $D^\natural \boxtimes \Pone$ is finitely generated over $A\llbracket K \rrbracket$, and the action of $A\llbracket K \rrbracket$ is continuous, 
for all~$n$ there exists~$i$ such that $\fN_n$ is an $A\llbracket K_i \rrbracket$-module.
Hence the (finite) intersection $\fM_n \coloneqq \bigcap_{g \in K}g(\fN_n)$ is an $A\llbracket K \rrbracket$-stable lattice, 
and it is preserved by the descent data, since these data are $\cO\llbracket K \rrbracket$-equivariant.
Thus~$\fM_n$ is the required neighborhood basis.

There remains to construct a neighborhood basis~$\fN_n$ of the origin consisting of $R$-equivariant lattices.
Since $D^\natural \boxtimes \Pone$ is an $R$-equivariant lattice in $D \boxtimes \Pone$, it suffices to do this for $D \boxtimes \Pone$.
By Lemma~\ref{lem: seeing D box P1 as a Tate module} we have an isomorphism of $R$-equivariant Tate $A$-modules
\[
D \boxtimes \Pone \isoto D \oplus (D \boxtimes p\bZ_p),
\]
and $D \boxtimes p\bZ_p$ is a direct summand of~$D$.
So it suffices to construct a neighborhood basis of the origin $\fN'_n \subset D$
such that each $\fN'_n$ is preserved by the descent data.
We can achieve this by setting $\fN'_n \coloneqq  T^n \cdot D^+$;
these are lattices in $D$, by Corollary~\ref{cor: Dplus is a lattice},
and they are preserved by the descent data,
by Lemma~\ref{lem:flat base-change for Dplus}.
\end{proof}

Lemma~\ref{lem:equivariant lattices for algebraic stacks}
shows that the descent data on $D^{\natural}\boxtimes \Pone$ preserve a neighbourhood
basis of the origin,
and hence the associated pro-coherent sheaf
$\widetilde{D^\natural \boxtimes \Pone}$ may again be equipped with descent data,
and so descends to a pro-coherent sheaf on~$\cZ$,
which we will denote by $\widetilde{D^\natural_\cZ \boxtimes \Pone}$.
Concretely, if we choose a basis~$\fM_n$ as in Lemma~\ref{lem:equivariant lattices for algebraic stacks}, and define $X_n \coloneqq (D^\natural \boxtimes \Pone)/\fM_n$, then
we have an isomorphism
\begin{equation}\label{eqn:presentation over cover}
\widetilde{D^\natural \boxtimes \Pone} \isoto \quoteslim{n} \widetilde{X_n},
\end{equation}
where each $\widetilde{X_n}$ is a coherent sheaf with $\cO\llbracket K \rrbracket$-action (i.e.\ a left $\cO\llbracket K \rrbracket$-module in $\Coh(\Spec A)$)
and descent data to~$\cZ$.
Note that~$X_n$ is a smooth representation of~$\cO\llbracket K \rrbracket$, and so the $\cO\llbracket K \rrbracket$-action on~$X_n$ factors through 
$\cO\llbracket K \rrbracket/\fa^i$ for some~$i$, since $X_n$ is finitely generated over~$A$.
In other words, using terminology from Remark~\ref{rem:Polish Iwasawa III}, $\tld X_n$ is an $\fa$-power torsion $\cO\llbracket K \rrbracket$-module in~$\Coh(\Spec A)$.

The right-hand side of~\eqref{eqn:presentation over cover} descends to a pro-system of $\fa$-power torsion $\cO\llbracket K \rrbracket$-modules~$\cF_n$ in $\Coh(\cZ)$,
and so
\begin{equation}\label{eqn:presentation over cZ}
\widetilde{D^\natural_\cZ \boxtimes \Pone} \coloneqq \quoteslim{n} \cF_n
\end{equation}
acquires the structure of an $\cO\llbracket K \rrbracket$-module in~$\Pro\Coh(\cZ)$.
Since the descent data on $D^\natural \boxtimes \Pone$ are also~$\cO\llbracket G\rrbracket $-equivariant (as we noted above),
this $\cO\llbracket K \rrbracket$-module structure on $\widetilde{D^\natural_\cZ \boxtimes \Pone}$ extends to an $\cO\llbracket G \rrbracket$-module structure.

\begin{rem}\label{rem:descended module is independent of presentation}
One easily verifies that the pro-coherent sheaf of~$\cO\llbracket G\rrbracket $-modules $\widetilde{D^{\natural}_{\cZ}\boxtimes \Pone}$
is canonically independent of the choice of presentation of $\cZ$ of the form~$[\Spec A/ R].$ 
Indeed, any two such presentations may be dominated by a third, and by  Lemma~\ref{lem: base change results we can use to prove G-stability}~\eqref{item: flat case}, if $\Spec A' \to \Spec A \to \cZ$ are {\em fppf} morphisms,
giving rise to rank $2$ \'etale $(\varphi,\Gamma)$-modules $D$ and $D'$ over $A$ and $A'$
respectively, then $D^{\natural} \boxtimes \Pone$ is obtained from $(D')^{\natural}\boxtimes \Pone$
via descent along the first arrow. 
This implies the analogous fact for the corresponding
pair of pro-coherent sheaves of~$\cO\llbracket G\rrbracket $-modules
$\widetilde{(D')^{\natural}\boxtimes \Pone}$ and $\widetilde{D^{\natural}\boxtimes \Pone}$,
and hence the pro-coherent sheaves
of $\cO\llbracket G\rrbracket $-modules on~$\cZ$ obtained by descending each of 
these are canonically isomorphic, as required.
\end{rem}

The next lemma describes the base-change properties of $\widetilde{D^{\natural}_{\cZ}\boxtimes \Pone}$. %

\begin{lemma}
\label{lem:base-change in the algebraic stack case}
Let $\cZ$ be a Noetherian algebraic stack over~$\cO/\varpi^a$, with affine diagonal, and equipped with
a morphism $\cZ \to \cX$, inducing the
pro-coherent sheaf of $\cO\llbracket G\rrbracket $-modules~$\widetilde{D^{\natural}_{\cZ}\boxtimes \Pone}$ on~$\cZ$,
as above. 
Further, let $f: \cW \to \cZ$ be a morphism of Noetherian algebraic stacks,
giving rise to the composite morphism $\cW \to \cZ \to \cX$, and hence similarly inducing
a pro-coherent sheaf of $\cO\llbracket G\rrbracket $-modules~$\widetilde{D^{\natural}_{\cW}\boxtimes \Pone}$
on~$\cW$. 
Then there is a base-change morphism of pro-coherent sheaves of $\cO\llbracket G\rrbracket $-modules
$$f^*(\widetilde{D^{\natural}_{\cZ}\boxtimes\Pone)} \to \widetilde{D^{\natural}_{\cW}\boxtimes
\Pone},$$
whose kernel and cokernel are coherent sheaves on~$\cW$,
and which is furthermore an isomorphism if~$f$ is flat.

If $\cZ$ and $\cW$ are furthermore of finite type over~$\cO/\varpi^a$ {\em (}for some
$a \geq 1${\em )}, then the kernel and cokernel of the base-change 
morphism are set-theoretically supported on~$\cW\times_{\cX} \cX(\St).$
\end{lemma}
\begin{proof}
We may find a commutative diagram 
$$\xymatrix{\Spec B \ar^g[r]\ar[d] & \Spec A \ar[d] \\
\cW \ar^f[r] & \cZ}
$$
in which the vertical arrows are {\em fppf} surjections,
and then use these vertical arrows to construct each of
$\widetilde{D^{\natural}_{\cZ}\boxtimes \Pone}$ 
and
$\widetilde{D^{\natural}_{\cW}\boxtimes \Pone}$
via descent 
from the corresponding pro-coherent sheaves
$\widetilde{D^{\natural}_{A}\boxtimes \Pone}$ 
and
$\widetilde{D^{\natural}_{B}\boxtimes \Pone}$
(using evident notation).

Lemmas~\ref{lem:pro-coherent base-change}
and~\ref{lem: base change for D box P1} then yield
a canonical morphism
$$g^*
(\widetilde{D^{\natural}_{A}\boxtimes \Pone}) 
\to
\widetilde{D^{\natural}_{B}\boxtimes \Pone},$$
which,
by Lemma~\ref{lem: base change results we can use to prove G-stability}~\eqref{item: flat case},
is even an isomorphism if~$f$ is flat (since then the morphism $g:\Spec B \to \Spec A$
is also flat).
This morphism is %
compatible with the descent data, 
and so induces the desired morphism
$$f^*(\widetilde{D^{\natural}_{\cZ}\boxtimes\Pone)} \to \widetilde{D^{\natural}_{\cW}\boxtimes
\Pone}.$$
The claim about its kernel and cokernel can be checked after pulling back to~$\Spec B$
(i.e.\ after ``undoing'' the descent),
where it follows from Lemma~\ref{lem:non-flat base-change}.
The final claim in the finite type case similarly follows from
Corollary~\ref{cor:support-of-base-change-ker-coker}.
\end{proof}

\subsubsection{Descent for~\texorpdfstring{$\cX$}{X}.}\label{subsubsec: descent for cX}
We now wish to define the pro-coherent sheaf of~$\cO\llbracket G\rrbracket $-modules
$\widetilde{D_{\cX}^{\natural}\boxtimes \Pone}$
associated to the universal \'etale $(\varphi,\Gamma)$-module $D_{\cX}$ on the 
formal algebraic stack~$\cX$.  

To this end, write $\cX \iso \colim \cX_n$
as a colimit of algebraic stacks under thickenings; 
let $i_n: \cX_n \hookrightarrow\cX_{n+1}$ denote the $n$th transition morphism,
and let $k_n: \cX_n \hookrightarrow \cX$ be the closed immersion.
By Lemma~\ref{lem:base-change in the algebraic stack case}, %
we have pro-coherent sheaves of $\cO\llbracket G\rrbracket $-modules 
$\widetilde{D^{\natural}_{\cX_n}\boxtimes \Pone},$
endowed with $\cO\llbracket G\rrbracket $-equivariant morphisms 
$$\widetilde{D^{\natural}_{\cX_{n+1}}\boxtimes \Pone}\to
i_{n,*} (\widetilde{D^{\natural}_{\cX_n}\boxtimes \Pone})$$
(obtained by applying the adjunction between $i_n^*$ and $i_{n,*}$ to 
the base-change morphisms of Lemma~\ref{lem:base-change in the algebraic stack case}),
which we may push forward to~$\cX$ so as to obtain~$\cO\llbracket G\rrbracket $-equivariant morphisms
$$k_{n+1,*}(\widetilde{D^{\natural}_{\cX_{n+1}}\boxtimes \Pone})\to
k_{n+1,*}i_{n,*}( \widetilde{D^{\natural}_{\cX_n}\boxtimes \Pone})
\cong
k_{n,*}( \widetilde{D^{\natural}_{\cX_n}\boxtimes \Pone}).$$
We thus obtain a projective system
$\bigl(k_{n,*}( \widetilde{D^{\natural}_{\cX_n}\boxtimes \Pone})\bigr)$
of pro-coherent sheaves of $\cO\llbracket G\rrbracket $-modules on~$\cX$, and we define
\begin{equation}\label{eqn:definition of tilde D natural box Pone}
\widetilde{D_{\cX}^{\natural}\boxtimes \Pone}
\coloneqq 
\limcommand_n k_{n,*}(\widetilde{D_{\cX_n}^{\natural}\boxtimes \Pone}),
\end{equation}
the inverse limit being formed in the category $\Pro \Coh(\cX)$.
It has the structure of a left $\cO\llbracket G\rrbracket $-module in~$\Pro \Coh(\cX)$.
Since any two descriptions of $\cX$ as a colimit of closed algebraic substacks are mutually
cofinal, it is well-defined independently of the choice of such a description used in its
construction. 

\begin{rem}\label{rem:Iwasawa module structure on D natural box Pone}
If we combine~\eqref{eqn:definition of tilde D natural box Pone} and~\eqref{eqn:presentation over cZ}, we find that there is a countable sequence~$\cF_n$ of $\fa$-power torsion $\cO\llbracket K\rrbracket$-modules
in~$\Coh(\cX)$, and an $\cO\llbracket K \rrbracket$-linear isomorphism
\begin{equation}\label{eqn:presentation over cX}
\widetilde{D^\natural_\cX \boxtimes \Pone} \isoto \quoteslim{n} \cF_n
\end{equation}
in $\Pro \Coh(\cX)$.
\end{rem}

Our next results describe some of the base-change properties
of~$\widetilde{D^{\natural}_{\cX}\boxtimes\Pone}$.
We first consider its behaviour with respect to versal rings.
In the statement of the next lemma, we make use of the completed pullback functor~$\ihat^*_x$ from Section~\ref{subsubsec:completed-pullback-to-completed}.

\begin{lem}\label{lem: L_infty to versal rings algebraic case}
Let $\cZ \to \cX$ be a morphism of finite type with $\cZ$ an algebraic stack with affine diagonal.
Let~$S$ be a complete local Noetherian $\cO$-algebra with finite residue field, and let $i_x: \Spf S \to \cZ$ be a 
  versal morphism at a finite type point~$x$ of $|\cZ|$, corresponding to a formal \'etale $(\varphi, \Gamma)$-module~$\widehat D_S$ with coefficients in~$S$.
Then $\ihat^*_{x}\bigl(\widetilde{D^\natural_{\cZ} \boxtimes \Pone}\bigr)$, which is an object of $\Pro \Coh(\Spf S) = \Mod_c(S)$, is naturally isomorphic to
$\widehat D^\natural_S \boxtimes \Pone$, as defined in Section~{\em \ref{subsec:formal phi Gamma modules}}.
\end{lem}
\begin{proof}
Since $\cZ$ is algebraic, $i_x$ is the $\fm_S$-adic completion of a morphism $i_x : \Spec S \to \cZ$, i.e.\ the formal \'etale $(\varphi, \Gamma)$-module $\widehat D_{S}$ 
is the $\fm_S$-adic completion of
an \'etale $(\varphi, \Gamma)$-module with coefficients in~$S$, which we denote by $D_{S}$.
By Lemma~\ref{lem: completed pullback at effective versal morphism}, $\ihat^*_x$ is the composite
\[
\Pro \Coh(\cZ) \xrightarrow{i_x^*} \Pro \Mod^{\fp}(S) \xrightarrow{\eqref{eqn:mod-fp-to-mod-c}} \Mod_c(S)
\]
where the first map is the $\Pro$-extension of the usual coherent pullback $i_x^* : \Coh(\cZ) \to \Mod^{\fp}(S)$.
Since $\cZ$ is of finite type over~$\cO/\varpi^a$ for some~$a$,
we see from \cite[\href{https://stacks.math.columbia.edu/tag/0DR2}{Tag~0DR2}]{stacks-project} that~$i_x$ is flat.
Hence Lemma~\ref{lem:base-change in the algebraic stack case} 
shows that
\[
i_x^*\bigl(\widetilde{D^\natural_{\cZ} \boxtimes \Pone}\bigr) = \bigl(\widetilde{D^\natural_{S} \boxtimes \Pone}\bigr).
\]
Applying Lemma~\ref{lem: change of topology}, we conclude that
\[
\ihat^*_x \bigl(\widetilde{D^\natural_{\cZ} \boxtimes \Pone}\bigr) = (D^\natural_{S} \boxtimes \Pone)^w \in \Mod_c(S),
\]
where the right-hand side is endowed with its weak topology.
By Lemma~\ref{lem:compatibility}~\eqref{item:71}, we conclude that
\begin{equation*} %
\ihat^*_x \bigl(\widetilde{D^\natural_{\cZ} \boxtimes \Pone}\bigr) = \widehat D^\natural_{S} \boxtimes \Pone,
\end{equation*}
as desired.
\end{proof}

\begin{lem}\label{lem: L_infty to versal rings}Let~$S$ be a complete local Noetherian $\cO$-algebra with finite residue field, and let $i_x: \Spf S \to \cX$ be a 
  versal morphism at a finite type point~$x$ of $|\cX|$, corresponding to a formal \'etale $(\varphi, \Gamma)$-module~$\widehat D_S$ with coefficients in~$S$.
Then $\ihat^*_{x}\bigl(\widetilde{D^\natural_{\cX} \boxtimes \Pone}\bigr)$, which is an object of $\Pro \Coh(\Spf S) = \Mod_c(S)$, is naturally isomorphic to
$\widehat D^\natural_S \boxtimes \Pone$, as defined in Section~{\em \ref{subsec:formal phi Gamma modules}}.
\end{lem}
\begin{proof}
Fix a presentation $\cX = \colim_n \cX_n$ as a colimit of algebraic stacks under thickenings, and let $\Spf S_n \coloneqq  \Spf S \times_{\cX} \cX_n$, 
so that $\Spf S = \colim_n \Spf S_n$. 
Let $i_n : \Spf S_n \to \cX_n$ be the morphism induced by~$i_x$, which is versal to~$\cX_n$ at~$x$.
By Lemma~\ref{lem:computing-Pro-completed-pullback} we have %
\[
\ihat^*_{x}\bigl(\widetilde{D^\natural_{\cX} \boxtimes \Pone}\bigr) = \limcommand_n\ihat^*_n\bigl(\widetilde{D^\natural_{\cX_n} \boxtimes \Pone}\bigr),
\]
the limit being taken in~$\Mod_c(S)$.
By Lemma~\ref{lem: L_infty to versal rings algebraic case}, we have  natural isomorphisms %
\begin{equation}\label{eqn:versal pullback for algebraic stacks}
\ihat^*_n \bigl(\widetilde{D^\natural_{\cX_n} \boxtimes \Pone}\bigr) = \widehat D^\natural_{S_n} \boxtimes \Pone.
\end{equation}
Applying~$\limcommand_n$ to~\eqref{eqn:versal pullback for algebraic stacks} concludes the proof.
\end{proof}

We now consider the base-change properties of
$\widetilde{D_{\cX}^\natural\boxtimes \Pone}$
with respect to %
finite type morphisms of algebraic stacks into~$\cX$.

\begin{remark}
In anticipation of the statement of  Lemma~\ref{lem:non-flat base-change, universal case},
we recall the stack~$\cX(\St)$ from Definition~\ref{def:Steinberg locus}, and
we note that if $\cY \to \cX$ is a finite type morphism
whose domain is an algebraic stack, then $\cY\times_{\cX}\cX(\St)$
is a formal algebraic substack of~$\cY$, equal to the completion
of $\cY$ along its closed algebraic substack $\cY\times_{\cX} \cX(\St)_{\red}.$
Thus $\Coh(\cY\times_{\cX} \cX(\St))$ embeds as a full subcategory 
of~$\Coh(\cY)$,
with essential image equal to the full subcategory 
of $\Coh(\cY)$ consisting of sheaves set-theoretically supported
on~$\cY\times_{\cX} \cX(\St)_{\red}.$
Consequently, 
$\Pro \Coh(\cY\times_{\cX} \cX(\St))$ embeds as a full subcategory 
of~$\Pro \Coh(\cY)$.
\end{remark}

\begin{lemma}
\label{lem:non-flat base-change, universal case}
Let $f:\cY \to \cX$ be a finite type morphism 
whose domain is an algebraic stack with affine diagonal
{\em (}so that $\cY$ is necessarily of finite type over~$\cO/\varpi^a$ for some $a \geq 1${\em )},
and let $D_{\cY}$ denote the pullback to~$\cY$ of the universal $(\varphi,\Gamma)$-module
over~$\cX$.
Then there is a base-change morphism
$${f}^*\bigl(\widetilde{D_{\cX}^\natural\boxtimes \Pone})
\to \widetilde{D_{\cY}^\natural\boxtimes \Pone} $$
{\em(}where ${f}^*$ denotes the pullback functor on pro-coherent sheaves
of Remark~{\em \ref{rem:Pro-extensions-denoted-same-symbols})}
of~$\cO\llbracket G\rrbracket $-module objects in $\Pro \Coh (\cY),$
whose kernel and cokernel, which {\em a priori} are
$\cO\llbracket G\rrbracket $-module objects of~$\Pro \Coh(\cY)$,
are in fact $\cO\llbracket G\rrbracket $-module objects of
$\Pro\Coh\bigl(\cY\times_{\cX}\cX(\St)\bigr)$.

If $f$ is furthermore a closed immersion, then this kernel and cokernel are in fact
objects of~$\Coh\bigl(\cY\times_{\cX}\cX(\St)\bigr).$
\end{lemma}
\begin{proof}
As in the discussion at the beginning of Section~\ref{subsubsec: descent for cX}, write $\cX$ as a colimit of thickenings of closed algebraic substacks
$k_n:\cX_n \hookrightarrow \cX$, so that $\cX \coloneqq  \colim \cX_n$. 
Since $\cY$ is %
an algebraic stack, and the morphism $f$ is of finite type, %
we see that $f$ factors through $\cX_n$ for $n$ sufficiently
large, and so relabelling the $\cX_n$ if necessary, we may assume
that $f$ factors through a %
morphism
$f_n: \cY \to \cX_n$ for each~$n$.

From~\eqref{eqn:definition of tilde D natural box Pone} we obtain isomorphisms
\begin{equation}
\label{eqn:base-change from formal to algebraic}
f^*
\widetilde{D^{\natural}_{\cX} \boxtimes \Pone} 
\iso \limcommand_n f^* k_{n,*}(\widetilde{D^{\natural}_{\cX_n}\boxtimes \Pone})
\iso \limcommand_n f_n^*(\widetilde{D^{\natural}_{\cX_n}\boxtimes \Pone}).
\end{equation}
(To see the first isomorphism, recall that $f^*$ is {\em defined} to be the pro-extension
of pullback on coherent sheaves, and so commutes with arbitrary cofiltered limits;
the second isomorphism holds simply because $f = k_n \circ f_n$ and~$k_{n, *}$ is fully faithful.)
Now Lemma~\ref{lem:base-change in the algebraic stack case} gives
us, for each~$n$, a base-change morphism
\begin{equation}
\label{eqn:nth base-change}
 f_n^*(\widetilde{D^{\natural}_{\cX_n}\boxtimes \Pone}) \to
\widetilde{D^{\natural}_{\cY}\boxtimes \Pone},
\end{equation}
whose kernel and cokernel, which we denote by $\cK_n$ and $\cC_n$ respectively,
are coherent sheaves on~$\cY$ which are set-theoretically supported on~$\cY(\St)\coloneqq  \cY\times_{\cX} \cX(\St)$.
By definition, this means we may regard~$\cC_n$ and~$\cK_n$ as coherent sheaves on the formal
algebraic stack~$\cY(\St).$
The morphisms~\eqref{eqn:nth base-change} form a projective system as~$n$ varies,
and thus so do the $\cC_n$ and $\cK_n$;
we then let $\cC \coloneqq  \limcommand_n \cC_n$ and $\cK \coloneqq  \limcommand_n \cK_n,$
the inverse limits being formed in~$\Pro \Coh(\cY(\St))$, or equivalently in~$\Pro \Coh(\cY)$.  

Passing to the inverse limit over the morphisms~\eqref{eqn:nth base-change},
and taking into account the isomorphism~\eqref{eqn:base-change from formal to algebraic},
we obtain the desired base-change morphism
$$
 f^*(\widetilde{D^{\natural}_{\cX}\boxtimes \Pone}) \to
\widetilde{D^{\natural}_{\cY}\boxtimes \Pone},$$
whose kernel and cokernel are isomorphic to~$\cK$ and $\cC$ respectively.

Assume now that  $f$ is a closed immersion; it remains to prove that in this case
 $\cC$ and $\cK$ are in fact objects of~$\Coh(\cY(\St))$
(and not merely pro-objects). Since $f$ is a closed immersion, the stack $\cY(\St)$ is a closed formal algebraic substack of
the formal algebraic stack~$\cX(\St)$, and so
$\Coh(\cY(\St))$ fully faithfully embeds into~$\Coh\bigl(\cX(\St)\bigr).$
Thus, we may regard each of $\cC_n$ and $\cK_n$ as objects
of~$\Coh\bigl(\cX(\St)\bigr),$
and hence regard $\cC$ and $\cK$ as objects of~$\Pro\Coh\bigl(\cX(\St)\bigr).$
It then suffices to show that $\cC$ and $\cK$ are in fact objects
of $\Coh\bigl(\cX(\St)\bigr).$

If $\cX(\St)$ is empty, then of course this entire discussion is vacuous and
there is nothing to prove.  Otherwise,
after making an unramified quadratic base-change in our coefficients and a twist, if necessary,
we may write %
\[
\cX(\St) = \coprod %
\cX_{\thetabar}
\]
where~$\thetabar$ runs through the four quadratic twists of~$1+\omega^{-1}$.
By the results of Section~\ref{subsubsec:Steinberg CWE stack}, for each~$\thetabar$
there exists a finite type $\bZ$-graded $R^{\ps}_{\thetabar}$-algebra~$S$,
and an isomorphism
\[
\colim_n [\Spec (S/\fm_{R^{\ps}_{\thetabar}}^nS)/\bG_m] \isoto \cX_{\thetabar}.
\]
Writing $i_{x}: \cX_{x} \to \cX_{\thetabar}$ for the completion of~$\cX_{\thetabar}$ at its unique closed point, we
thus deduce from Theorem~\ref{thm:summary-of-AHL}~\eqref{item:BGG-derived-formal-functions}
that the completed pullback is an equivalence
\[
\widehat i_{x} : \Coh(\cX_{\thetabar}) \isoto \Coh(\cO_{\cX_{x}}).
\]
Thus, in order to show that $\quoteslim{n} \cK_n$ and $\quoteslim{n} \cC_n$ are objects of
$\Coh(\cX(\St))$, it suffices to show (for each~$\thetabar$) that
$\quoteslim{n} \ihat_{x}\cK_n$ and $\quoteslim{n} \ihat_{x}\cC_n$ are objects of
$\Coh(\cO_{\cX_{x}})$.

Let $v: \Spf R^{\ver} \to \cX_{\thetabar}$ be the versal ring at the closed point 
induced by the completion of~$S$ at its maximal homogeneous ideal~$\fm_{S}$, so that there is an isomorphism
\begin{equation}\label{eqn:presentation of versal ring}
\colim_n [\Spec (R^{\ver}/\fm_{R^{\ver}}^nR^{\ver})/\bG_m] \isoto \cX_{x}.
\end{equation}
To ease notation, we write
$R \coloneqq  R^{\ver}$ for the remainder of the proof, and let~$I \subset R$, resp.\ $I_n \subset R$, 
be the ideal of the pullback of~$\cY \to \cX$ to $\Spf R$, resp.\ the ideal of the pullback of $\cX_n \to \cX$
to~$\Spf R$. 
Since $\cX$ is Noetherian, %
the $I_n$ are cofinal with the
powers~$I_1^n$.

Since $v : \Spf R \to \cX_{x}$ is smooth and surjective, and the coherence of a pro-coherent sheaf can be tested after pulling back to a
flat cover (by Remark~\ref{arem:checking coherence on smooth cover}),
it now suffices to 
prove that
the pullbacks $v^*\quoteslim{n} \ihat_{x}^*\cK_n$ and $v^*\quoteslim{n}\ihat_{x}^*\cC_n$ are objects of $\Coh(\Spf R)$.
Bearing in mind the isomorphism~\eqref{eqn:presentation of versal ring},
we see that the morphism $v : \Spf R \to \cX_{x}$ is representable by algebraic stacks, and so the pullback
$v^*$ coincides with the completed pullback $\widehat v^*$, by Remark~\ref{rem:pullback and completed pullback}.
For simplicity, in the following we will also write $v$ for the induced morphisms $\Spf R/I \to \cY_{x}, \Spf R/I_n \to \cX_{n, x}$.

Applying~$v^* \ihat^*_{x}$ to~\eqref{eqn:nth base-change}, we obtain an exact sequence
\[
0 \to v^*\ihat_{x}^*\cK_n \to v^* \ihat^*_{x} f_n^*(\widetilde D^\natural_{\cX_n} \boxtimes \Pone) \to v^*\ihat^*_{x} \widetilde{D^\natural_\cY \boxtimes \Pone} \to v^*\ihat_{x}^*\cC_n \to 0.
\]
Lemma~\ref{lem: L_infty to versal rings algebraic case} shows that
$v^*\ihat^*_{x} \widetilde{D^\natural_\cY \boxtimes \Pone} = \widehat D_{R/I}^\natural \boxtimes \Pone$.
Similarly, using the fact that $v^* = \widehat v^*$, the commutative diagram~\eqref{to prove commutative for compatibility of completion and pushforward} implies that
\[
v^* \ihat^*_{x} f_n^*(\widetilde D^\natural_{\cX_n} \boxtimes \Pone) = (\widehat D_{R/I_n}^\natural \boxtimes \Pone) \otimes_{R/I_n} R/I,
\]
and so we have an exact sequence
\[
0 \to v^*\ihat_{x}^*\cK_n \to (\widehat D_{R/I_n}^\natural \boxtimes \Pone) \otimes_{R/I_n} R/I \to \widehat D_{R/I}^\natural \boxtimes \Pone \to v^*\ihat_{x}^*\cC_n \to 0.
\]
By Lemma~\ref{lem:coherence}, it follows that 
$\quoteslim{n}v^* \ihat_{x}^*\cK_n$ and $\quoteslim{n}v^*\ihat_{x}^*\cC_n$ are objects of $\Coh(\Spf R)$, as desired.\qedhere

\end{proof}

\subsubsection{Defining~\texorpdfstring{$L_\infty$}{L-infty} and~\texorpdfstring{$\Functor$}{F}}
Recall from Remark~\ref{rem: central char} that the centre~$Z$ of~$G$ acts on~$D^\natural \boxtimes \Pone$ by
multiplication by~$\zeta\varepsilon^{-2}: \bQ_p^\times \to \cO^\times$, i.e.\ $D^\natural \boxtimes \Pone$ is a left $\cO\llbracket G\rrbracket _{\zeta\varepsilon^{-2}}$-module.

\begin{df}\label{definition of Linfty} 
We let $L_\infty$
denote the right $\cO\llbracket G\rrbracket _\zeta$-module in $\Pro \Coh(\cX)$ %
obtained by twisting the left $\cO\llbracket G\rrbracket _{\zeta\varepsilon^{-2}}$-action on
$\widetilde{D_{\cX}^{\natural}\boxtimes \Pone}$
by~$\zeta^{-1}\varepsilon \circ \det$, and then turning the resulting left 
$\cO\llbracket G\rrbracket _{\zeta^{-1}}$-action into a right $\cO\llbracket G\rrbracket _\zeta$-action by precomposing with $g \mapsto g^{-1}$.%
\end{df}

Lemma~\ref{lem:derived-tensor-product-EW}~\eqref{item:87}
(and the equivalence $D^b\bigl(\Coh(\cX)\bigr)\iso D^b_{\coh}(\cX)$ of~\eqref{eqn:coherent equivalence}) 
 then allows us to define a right $t$-exact functor
\begin{equation}\label{eqn:derived level functor for O-G-modules}  D^b_{\fp}(\cO\llbracket G\rrbracket _{\zeta})\to
  \Pro D^b_{\coh}(\cX), \quad \pi \mapsto L_\infty \otimes^L_{\cO\llbracket G\rrbracket _\zeta}
  \pi.
\end{equation} 
\begin{df}\label{defn:definition of F}
We let
\begin{equation}\label{eqn:defn-of-our-functor}\Functor: D^b_{\fp}(\cA)\to\Pro
  D^b_{\coh}(\cX)  %
\end{equation}
denote the right $t$-exact functor \[\pi \mapsto %
  L_\infty \otimes^L_{\cO\llbracket G\rrbracket _\zeta}
  \pi\] given by the composite of~\eqref{eqn:derived level functor for O-G-modules} and the $t$-exact fully faithful functor \[D^b_{\fp}(\cA)\into
 D^b_{\fp}(\cO\llbracket G\rrbracket _{\zeta}) \] induced by~\eqref{eq:Afp-in-Modfp-OG}.
\end{df}

\begin{rem}\label{rem: F factors through D^b_coh}
We will see below in Proposition~\ref{finiteness of
  cohomology} that $\Functor$ %
factors through a
functor~$D^b_{\fp}(\cA)\to D^b_{\coh}(\cX)$, though this does not seem obvious from its definition.
\end{rem}

The following result will allow us to get something of a handle
on the pro-coherent structure of~$L_{\infty}$. 

\begin{lemma} 
\label{lem:explicit pro-coherent structure}
If $i: \cY \hookrightarrow \cX$ is the inclusion of
a closed algebraic substack, 
then:
\begin{enumerate}
\item\label{item:53} each $(i^*L_{\infty})_{K_n}$ is a coherent sheaf on~$\cY$, and 
\item\label{item:52} $i^*L_{\infty} \iso \quoteslim{n} (i^*L_{\infty})_{K_n}$. %
\end{enumerate}
\end{lemma}
\begin{proof}

Evidently the twists and duals in Definition~\ref{definition of Linfty} %
play no role in these assertions, and so it is enough to prove the
statements of the lemma with~$L_\infty$ replaced by $\widetilde{D_{\cX}^{\natural}\boxtimes \Pone}$.
Let $p : \Spec A \to \cY$ be an \emph{fppf} and surjective morphism.
By Remark~\ref{arem:checking coherence on smooth cover}, the coherence of a pro-coherent sheaf can be established after pulling back to a flat cover, 
so it suffices to prove that
\begin{enumerate}
\item each $(p^*i^*\widetilde{D_{\cX}^\natural \boxtimes \Pone})_{K_n}$ is a coherent sheaf on~$\Spec A$, and 
\item $p^*i^*\widetilde{D_{\cX}^\natural \boxtimes \Pone} \iso \quoteslim{n} (p^*i^*\widetilde{D_{\cX}^\natural \boxtimes \Pone})_{K_n}$.
\end{enumerate}
Both of these statements follow immediately from the assertion that $p^*i^*\widetilde{D_{\cX}^\natural \boxtimes \Pone}$
is isomorphic to $\tld M$ for some finitely presented $A\llbracket K \rrbracket$-module~$M$, which we now prove.

Consider the base-change morphism
\begin{equation}\label{eqn:specialized base-change morphism}
p^*i^*(\widetilde{D^{\natural}_{\cX}\boxtimes \Pone}) \to
\widetilde{D^{\natural}_{A} \boxtimes \Pone}
\end{equation}
of Lemma~\ref{lem:non-flat base-change, universal case}.
That lemma shows that we have an exact sequence  
of right $\cO\llbracket G \rrbracket$-modules in $\Pro \Coh(\Spec A)$
of the form
\begin{equation}\label{eqn:base change exact sequence for explicit pro-coherent structure}
0 \to \cK \to 
 p^*i^*(\widetilde{D^{\natural}_{\cX}\boxtimes \Pone}) \to
\widetilde{D^{\natural}_{A} \boxtimes \Pone} \to \cC \to 0,
\end{equation}
where~$\cK$ and~$\cC$ are furthermore coherent sheaves.

Recall from Remarks~\ref{rem:Polish Iwasawa III} and~\ref{rem:Polish Iwasawa IV}
that we write~$\AKModfpA$ for the category of $\fa$-power torsion $\cO\llbracket K \rrbracket$-modules in $\Coh(\Spec A)$, 
and ~$\OKProModfpA$ for the category of $\cO\llbracket K \rrbracket$-modules in $\Pro \Coh(\Spec A)$.
The forgetful functor $\AKModfpA \to \Coh(\Spec A)$ is exact and faithful, and  by Remark~\ref{rem:Polish Iwasawa IV} it extends to an exact, fully faithful, cofiltered limit-preserving functor
\begin{equation}\label{eqn:Pro modules to modules in Pro}
\Pro\AKModfpA \to \OKProModfpA.
\end{equation}
The functor
\begin{equation}\label{eqn:tilde functor again}
\Mod^{\fp}(A \llbracket K \rrbracket) \to \OKProModfpA,\ M \mapsto \tld M \coloneqq \quoteslim{i} M/\fa^i M 
\end{equation}
factors through~\eqref{eqn:Pro modules to modules in Pro} via the functor %
\begin{equation}\label{eqn:tilde functor again II}
\Mod^{\fp}(A\llbracket K \rrbracket) \to \Pro\AKModfpA,\ M \mapsto  \quoteslim{i} M/\fa^i M 
\end{equation} obtained from ~\eqref{eqn: tilde functor on Iwasawa modules},
 and Corollary~\ref{cor: D natural box P1 is finite over K} implies that
$\widetilde{D_A^\natural \boxtimes \Pone}$, which is \emph{a priori} an object of $\OKProModfpA$,
is contained in 
the essential image of~\eqref{eqn:tilde functor again II}.
By~\eqref{eqn:presentation over cX}, the same is true for $p^*i^*(\widetilde {D^\natural_\cX \boxtimes \Pone})$.
The exact sequence~\eqref{eqn:base change exact sequence for explicit pro-coherent structure} can therefore be viewed as an exact sequence in~$\Pro \AKModfpA$.

We are now in a position to apply Lemma~\ref{lem:topology-Iwasawa-algebra}~\eqref{item: Iwasawa 7} (bearing in mind Remark~\ref{rem:Polish Iwasawa III}),
which shows that $p^*i^*(\widetilde{D^{\natural}_{\cX}\boxtimes \Pone})$ is contained in the essential image of~\eqref{eqn:tilde functor again II}
if this is true for~$\cK$ and~$\cC$ (recalling that we have already seen that this is true of  $\widetilde{D^{\natural}_A \boxtimes \Pone}$).
Finally, since~$\cK$ and~$\cC$ are objects of $\Coh(\Spec A)$ in the essential image of~\eqref{eqn:Pro modules to modules in Pro}, they are the coherent sheaves underlying
$\fa$-adically discrete $\cO\llbracket K \rrbracket$-modules in $\Mod^{\fp}(A)$, and so  (again bearing in mind Remark~\ref{rem:Polish Iwasawa III}) they are contained in the essential image of~\eqref{eqn:tilde functor again II}, as desired.\qedhere
\end{proof}

\begin{remark}
\label{rem:explicit pro-coherent structure}
As already alluded to above,
the significance of Lemma~\ref{lem:explicit pro-coherent structure} is that it gives us a concrete description
of the pro-coherent structure of~$L_\infty$. Namely, if we write $\cX = \colim_m \cX_m$
as an Ind-algebraic stack with closed transition morphisms, and write $i_m:\cX_m \hookrightarrow
\cX$ for the corresponding closed immersion, 
then (bearing in mind Lemma~\ref{lem:pro-coherent-as-limit-of-pullbacks}) we have %
\begin{equation}\label{eqn:explicit-pro-coherent-via-coinvariants-of-pullback}
  L_{\infty} \iso \limcommand_{m,n} (i_m^*L_{\infty})_{K_n}.
\end{equation}
Actually, it follows from Proposition~\ref{concentration and support} below that %
that $(L_{\infty})_{K_n}$ is supported on a $\varpi$-adic formal algebraic
closed substack of~$\cX$, and so~\eqref{eqn:explicit-pro-coherent-via-coinvariants-of-pullback} does not convey the whole truth
of the matter.  Nevertheless, it is a key intermediate step in the ultimate 
analysis of~$L_{\infty}$.
\end{remark}

\subsection{Completing \texorpdfstring{$L_{\infty}$}{L infinity}}
\label{subsec:completing L-infinity}
A key tool in our analysis of the functor~$\Functor$ %
will be completion: 
both in the sense of completion at a block $\cA_{\thetabar}$ %
of the category $\cA$, and completion along
the various substacks $\cX_{\thetabar}$ of~$\cX$.
In fact, we will show that these two notions of completion are compatible,
in an appropriate sense,
with respect to~$\Functor$ (see Proposition~\ref{prop:completing our functor on Y}). 
Relatedly, we will introduce a ``completed-at-$\thetabar$'' version of~$\Functor$,
and prove that it is fully faithful (see Theorem~\ref{thm:derived-F-thetabar-nonSt-case}).
This will be a key ingredient in the proof that $\Functor$ itself is fully
faithful. 

\subsubsection{A completed version of~\texorpdfstring{$\Functor$}{F}}
We begin by defining a completed version~$\Functor_{\thetabar}$ of~$\Functor$. It is in
fact convenient to do this for any  closed subspace  $Y \subseteq |X|$. Recall
that in~\eqref{eqn:defn-of-cts-map-f} we have defined a continuous map $\piss:|\cX|\to|X|$. %
\begin{defn} %
  \label{defn:Y-substacks-i-prime}
  \leavevmode
  \begin{enumerate}
  \item If~$Y\subseteq|X|$ is a closed subset,
  then 
  we write $\cX_{Y}$ for %
  the completion of
  $\cX$ along~$\pi_{\mathrm{ss}}^{-1}(Y)$. 
  By Proposition~\ref{aprop:Coh-set-theoretically-supported} %
  we identify~ $D^b_{\coh}(\cX_Y)$ with~$D^b_{\coh,Y}(\cX)$, i.e.\ the full subcategory of complexes whose cohomology sheaves are set-theoretically supported on~$Y$. 
  \item We write $\ihat'_{Y,*}:\Pro D^b_{\coh}(\cX_Y)\to \Pro D^b_{\coh}(\cX)$,
  $\ihat'^{*}_Y:\Pro D^b_{\coh}(\cX)\to \Pro D^b_{\coh}(\cX_Y)$ for the
  adjoint pair of $t$-exact functors defined in Appendix~\ref{subsubsec:completed-pullback-to-completed}. %
  \end{enumerate}
  \end{defn}%

\begin{rem}\label{why i'}
We remind the reader that we have already defined functors~$i_{Y,*}$ and~$\widehat i_{Y,*}$, on derived categories of $G$-representations,
in~\eqref{definition of i_Y*} and~\eqref{definition of i_Y*-hat}, and similarly for~$\widehat i_Y^*$.
This is the reason for the notation~$i'_Y$.
\end{rem}

\begin{rem}
  \label{rem:compatible-Y-thetabar-notation}
  In the case that~$Y=\{\thetabar\}$ is a single closed point (identified with a conjugacy class of $\cbF_p$-valued pseudorepresentations as in Section~\ref{chain of projective lines})
  we also have the stack $\cX_{\thetabar} \to \cX$ introduced in~\eqref{eqn:CWE to EG}.
  By Theorem~\ref{thm: thetabar substacks}, $\cX_{\thetabar}$ is canonically isomorphic to $\cX_{Y}$,
  which we will therefore also denote by~$\cX_{\thetabar}$.
  Furthermore, we will tacitly use results and notation from Section~\ref{subsec:CWE stacks} and Section~\ref{subsec:sheaves on CWE stacks}, so that, for example,
  $k_{\thetabar} : \cX_{\thetabar} \to \fX_{\thetabar}$ is the $\fm$-adic completion morphism.
\end{rem}

\begin{lem}
  \label{lem:pullback-Linfty-associativity}For any closed subset~$Y\subseteq
  |X|$, we have a natural isomorphism of right $t$-exact functors $D^b_{\fp}(\cO\llbracket G\rrbracket _{\zeta})\to \Pro D^b_{\coh}(\cX_{Y})$ %
  \[\ihat'^*_{Y}
(L_{\infty}\otimes^L_{\cO\llbracket G\rrbracket _{\zeta}}\text{--}) 
\iso \ihat'^*_{Y}L_{\infty} \otimes^{L}_{\cO\llbracket G\rrbracket _{\zeta}}\text{--}.\]
\end{lem}
\begin{proof}By Lemma~\ref{lem:derived-tensor-product-EW}~\eqref{item:84},
it suffices to note that %
each functor 
takes $\cO\llbracket G\rrbracket _{\zeta}$ to %
$\ihat'^*_{Y}L_{\infty}$. %
(Note that this argument is just a particular case 
of an evident variant, in the derived context, of
Lemma~\ref{lem:abelian-composition-right-exact-tensor}.)
\end{proof}

\begin{df}%
\label{def:F thetabar}For any closed subset~$Y\subseteq |X|$, we let  $\Functor_{Y}$ denote the functor
$$\Functor_{Y} \coloneqq  \ihat'^*_Y \Functor \, i_{Y,*} :
\Dfp^b(\cA_{Y}) \to \Pro D^b_{\coh}(\cX_{Y}).$$ 
In particular, for any $2$-dimensional $\cbF_p$-valued pseudorepresentation~$\thetabar,$
we let~  $\Functor_{\thetabar}$ denote the functor
\begin{equation}\label{eqn:F_thetabar}\Functor_{\thetabar} \coloneqq  \ihat'^*_{\thetabar} \Functor \, i_{\thetabar,*} :
 \Dfp^b(\cA_{\thetabar}) \to \Pro D^b_{\coh}(\cX_{\thetabar}).\end{equation}
\end{df}

By Lemma~\ref{lem:pullback-Linfty-associativity} we have a natural isomorphism of functors $\Dfp^b(\cA_{\thetabar}) \to \Pro D^b_{\coh}(\cX_{\thetabar})$
\begin{equation}
\label{eqn:F thetabar formula}
\Functor_{\thetabar}(\text{--}) \coloneqq  \ihat'^*_{\thetabar}\Functor(i_{\thetabar, *}(\text{--})) = \ihat'^*_{\thetabar}
(L_{\infty}\otimes^L_{\cO\llbracket G\rrbracket _{\zeta}}i_{\thetabar, *}(\text{--})) 
\iso \ihat'^*_{\thetabar}L_{\infty} \otimes^{L}_{\cO\llbracket G\rrbracket _{\zeta}}i_{\thetabar, *}(\text{--}).
\end{equation}

\subsubsection{An alternative description of
  \texorpdfstring{$\ihat'^*_{\thetabar}L_{\infty}$}{the completion of L infinity}}\label{subsubsec:alternative-description-Linfty}

Our next goal is to give an alternate description of the functors~$\Functor_{\thetabar}$,
which will allow us to relate them to the functors studied in~\cite{JNWE}, and
thus also to the results of Colmez~\cite{MR2642409} and
Pa\v{s}k\={u}nas~\cite{MR3150248}. In view of~\eqref{eqn:F thetabar formula}, the key
step in giving such a description is to describe the
sheaf~$\ihat'^*_{\thetabar}L_{\infty}$, which we do in Proposition~\ref{prop:completed L-infinity} below, by relating it to the sheaf $\cV_{\thetabar} \cotimes_{\tE_{\thetabar}} \tP_{\thetabar}$,
whose definition we now recall.

We have defined in Definition~\ref{defn: tld P} a projective
object~$\tld P_{\thetabar}$ of $\fC_{\thetabar}$, whose endomorphism algebra~$\tE_{\thetabar}$
is equipped with a canonical isomorphism $\tE_{\thetabar} \isoto \tR_{\thetabar}^{\op}$. %
If $\thetabar$ is not of type~\ref{item: Steinberg pseudorep},
then $\tP_{\thetabar}$ is even a projective generator of~$\fC_{\thetabar}.$
If $\thetabar$ is of type~\ref{item: Steinberg pseudorep}, then $\tP_{\thetabar}$ is not
a projective generator of~$\fC_{\thetabar}$;
rather, we have the projective generator $\bP_{\thetabar}$, with endomorphism algebra~$\bE_{\thetabar}$, introduced in
Definition~\ref{defn: bP in Steinberg case}. 
In Section~\ref{subsubsec:the ideal J}, we have also introduced a two-sided ideal $J \subset \bE_{\thetabar}$. 

On the other hand, recall from~\eqref{eqn:definition of V_thetabar} that~$\cV_{\thetabar}$ denotes the universal object on~$\cX_{\thetabar}$, 
and is a complete right $\tE_{\thetabar}$-module in $\Pro\Coh(\cX_{\thetabar})$.
When~$\thetabar$ has type~\ref{item: Steinberg pseudorep}, recall furthermore the sheaf~$\cW_{\thetabar}$ from Definition~\ref{defn:W-thetabar}: it is
a complete right $\bE_{\thetabar}$-module in~$\Pro \Coh(\cX_{\thetabar})$.

We will need two preliminary results.
For the first one, recall from~\eqref{eqn:abelian complete to pro} the exact and fully faithful functor
\begin{gather*}
\Coh(\cO_{\cX_{\thetabar}}) \to \Pro \Coh(\cX_{\thetabar}), \\
\cF \mapsto \quoteslim{n} \cF/\fm_{R_{\thetabar}^{\ps}}^{n+1}\cF.
\end{gather*}

\begin{lemma}\label{lem:finiteness of H coinvariants}
For every compact open subgroup~$H \subset G$, the pro-coherent sheaves 
\[
\cV_{\thetabar} \cotimes_{\tE_{\thetabar}} (\tP_{\thetabar})_H
\]
and
\[
(\ihat'^*_{\thetabar}L_\infty)_H
\]
on~$\cX_{\thetabar}$ are objects of $\Coh(\cO_{\cX_{\thetabar}})$. 
Furthermore, the natural maps
\begin{equation}\label{eqn:limit description of L_infty over cX}
\ihat'^*_{\thetabar} L_\infty \isoto \lim_H (\ihat'^*_{\thetabar}L_\infty)_H
\end{equation}
and
\begin{equation}\label{eqn:limit description of V tensor P tilde}
\cV_{\thetabar} \cotimes_{\tE_{\thetabar}} \tP_{\thetabar} \isoto \lim_H \cV_{\thetabar} \cotimes_{\tE_{\thetabar}} (\tP_{\thetabar})_H
\end{equation}
are isomorphisms in $\Pro \Coh(\cX_{\thetabar})$.
\end{lemma}
\begin{proof}
Let~$R \coloneqq R_{\thetabar}^{\ps}, \fm \coloneqq \fm_{R_{\thetabar}^{\ps}}$ and
$\cF \coloneqq \cV_{\thetabar} \cotimes_{\tE_{\thetabar}} (\tP_{\thetabar})_H$.
To show that~$\cF \in \Coh(\cO_{\cX_{\thetabar}})$ it suffices to prove that $\cF \isoto \lim_{n} \cF/\fm^{n+1} \cF$, 
and that each $\cF/\fm^{n+1} \cF$ is an object of $\Coh(\cX_{\thetabar})$.

For the first statement, we begin with the isomorphism $\tP_{\thetabar} \isoto \lim_n\tP_{\thetabar}/\fm^{n+1}\tP_{\thetabar}$ arising from the fact that
$\tP_{\thetabar}$ is a complete $\tE_{\thetabar}$-module in $\Mod_c(\cO)$ (by Lemma~\ref{compact modules are compact}), hence a complete $R$-module 
(by Lemma~\ref{lem:equivalence-of-R-and-E-complete}).
We then use the facts that
the functor $(\text{--})_H$ preserves cofiltered limits in $\Mod_c(\cO \llbracket H \rrbracket^{\op})$ (by Lemma~\ref{lem:tensor and projective limits}),
and the functor $\cV_{\thetabar} \cotimes_{\tE_{\thetabar}} \text{--}$ preserves cofiltered limits in $\Mod_c(\tE_{\thetabar})$ (by definition).
This implies that $\cF \isoto \lim_n \cF/\fm^{n+1} \cF$, and also that 
\begin{equation*}
\cF/\fm^{n+1} \cF = \cV_{\thetabar} \cotimes_{\tE_{\thetabar}} (\tP_{\thetabar}/\fm^{n+1} \tP_{\thetabar})_H.
\end{equation*}
Now $\tP_{\thetabar}/\fm^{n+1}\tP_{\thetabar}$ is coadmissible: in fact, it is a quotient of $\tP_{\thetabar}/\rad(\tE_{\thetabar})^j\tP_{\thetabar}$ for some~$j$,
since~$\tE_{\thetabar}$ is a finite $R$-module, and $\tP_{\thetabar}/\rad(\tE_{\thetabar})^j\tP_{\thetabar}$ is coadmissible, by Lemma~\ref{lem:projective generators are flat over their endos}~\eqref{item:projective 4}.
Hence $(\tP_{\thetabar}/\fm^{n+1}\tP_{\thetabar})_H$ is finitely presented over~$\tE_{\thetabar}$ (in fact, even over~$\cO$).
We conclude that
\[
\cF/\fm^{n+1} \cF = \cV_{\thetabar} \otimes_{\tE_{\thetabar}} (\tP_{\thetabar}/\fm^{n+1} \tP_{\thetabar})_H = (\cV_{\thetabar}/\fm^{n+1}\cV_{\thetabar}) \otimes_{\tE_{\thetabar}} (\tP_{\thetabar}/\fm^{n+1} \tP_{\thetabar})_H,
\]
and since $\cV_{\thetabar}/\fm^{n+1} \cV_{\thetabar}$ is a coherent sheaf on~$\cX_{\thetabar}$, this implies that $\cF/\fm^{n+1}\cF$ is a coherent sheaf on~$\cX_{\thetabar}$, as desired.
This concludes the proof that $\cV_{\thetabar} \cotimes_{\tE_{\thetabar}} (\tP_{\thetabar})_H \in \Coh(\cO_{\cX_{\thetabar}})$.

Similarly, let $\cG \coloneqq (\ihat'^*_{\thetabar}L_\infty)_H$.
By definition we have $\cX_{\thetabar} = \colim_n \cX_{\thetabar}/\fm^{n+1}$, and $\cX_{\thetabar}/\fm^{n+1}$ is an algebraic stack.
Writing $i'_{\thetabar, n+1} : \cX_{\thetabar}/\fm^{n+1} \to \cX$ for the restriction of $i'_{\thetabar}$,
Lemma~\ref{lem:pro-coherent-as-limit-of-pullbacks} shows that
\begin{equation}\label{eqn:intermediate step in limit description II}
\ihat'^*_{\thetabar} L_\infty = \lim_n (i'^*_{\thetabar, n+1}L_\infty),
\end{equation}
where $i'^*_{\thetabar, n+1}L_\infty = \ihat'^*_{\thetabar}L_\infty/\fm^{n+1}\ihat'^*_{\thetabar}L_\infty$.
Since $\Pro \Coh(\cX_{\thetabar})$ has exact cofiltered limits, the coinvariant functor $(-)_H$ preserves cofiltered limits of right $\cO\llbracket H \rrbracket$-modules in $\Pro \Coh(\cX_{\thetabar})$
(by Lemma~\ref{lem:tensor and projective limits}).
We thus deduce that
\[
\cG %
= \lim_n(i'^*_{\thetabar, n+1}L_\infty)_H, %
\]
and since $(i'^*_{\thetabar, n+1}L_\infty)_H = \cG/\fm^{n+1}\cG$, we conclude that $\cG \isoto \lim_n \cG/\fm^{n+1}\cG$, and 
(using Lemma~\ref{lem:explicit pro-coherent structure}~\eqref{item:53})
that~$\cG/\fm^{n+1}\cG$ is coherent.
This concludes the proof that $(\ihat'^*_{\thetabar}L_\infty)_H \in \Coh(\cO_{\cX_{\thetabar}})$.

The isomorphism~\eqref{eqn:limit description of L_infty over cX} now follows from~\eqref{eqn:intermediate step in limit description II} and 
Lemma~\ref{lem:explicit pro-coherent structure}~\eqref{item:52}.
Finally, the isomorphism~\eqref{eqn:limit description of V tensor P tilde} follows from the isomorphism $\tP_{\thetabar} \isoto \lim_H (\tP_{\thetabar})_H$ 
(which is valid for every object of~$\fC$) 
and the fact that $\cV_{\thetabar} \cotimes_{\tE_{\thetabar}} \text{--}$ preserves cofiltered limits.
\end{proof}

Before stating the next result, we recall from Remark~\ref{rem:contragredient again} 
that objects of~$\fC_{\thetabar}$, such as $\tP_{\thetabar}$, can be regarded as right $\cO\llbracket G\rrbracket _\zeta$-modules
in a natural way.
\begin{prop}\label{prop:versal-ring-completed-L-infinity}
  Let~$\thetabar$ be a 2-dimensional $\cbF_p$-valued pseudorepresentation of~$G_{\bQ_p}$ of determinant~$\zeta\varepsilon^{-1}$.
  \begin{enumerate}
   \item\label{item:59} 
   If~$S$ is a complete Noetherian local $\cO$-algebra with finite residue field, and $i'_S: \Spf S\to\cX_{\thetabar}$ 
   is a versal morphism at the closed point of~$\cX_{\thetabar}$, 
   then there is a natural isomorphism of right $\cO\llbracket G\rrbracket _\zeta$-modules %
    \begin{equation}\label{eqn:sought-for-isomorphism-on-versal-rings}
(\ihat'^*_S \cV_{\thetabar}\cotimes_{\tld E_{\thetabar}} \tP_{\thetabar}) /
(\ihat'^*_S \cV_{\thetabar}\cotimes_{\tld E_{\thetabar}} \tP_{\thetabar})^{\SL_2(\Q_p)}
\iso
\ihat'^*_S L_{\infty},
    \end{equation}
    where the left-hand side has the right $\cO\llbracket G\rrbracket _\zeta$-action induced from that on~$\tP_{\thetabar}$, and the right-hand side has the 
    right $\cO\llbracket G\rrbracket _\zeta$-action 
induced from that on~$L_{\infty}$.
{\em (}Here, ``natural'' means ``natural in the versal morphism~$i'_S$''.{\em )}
  \item\label{item:60}   If~$\thetabar$ is not of type~\emph{\ref{item: Steinberg pseudorep}},
  then~\eqref{eqn:sought-for-isomorphism-on-versal-rings} induces a natural isomorphism of right $\cO\llbracket G\rrbracket _\zeta$-modules
  \begin{equation}\label{eqn:sought-for-isomorphism-on-versal-rings-non-St}
\ihat'^*_S \cV_{\thetabar}\cotimes_{\tld E_{\thetabar}} \tP_{\thetabar} \iso
\ihat'^*_S L_{\infty}.
  \end{equation}
  \item\label{item:62} If~$\thetabar$ is of type~\emph{\ref{item: Steinberg pseudorep}},
  then~\eqref{eqn:sought-for-isomorphism-on-versal-rings} induces %
a natural isomorphism of right $\cO\llbracket G\rrbracket _{\zeta}$-modules 
  \begin{equation}\label{eqn:sought-for-isomorphism-on-versal-rings-St}
(\ihat'^*_S \cW_{\thetabar}/\ihat'^*_S \cW_{\thetabar}[J])\cotimes_{\bE_{\thetabar}}
\bP_{\thetabar} \iso 
\ihat'^*_S L_{\infty},
  \end{equation}
\end{enumerate}
where
$J \subset \bE_{\thetabar}$
is the two-sided ideal defined in Section~{\em \ref{subsubsec:the ideal J}},
and 
$\ihat'^*_S \cW_{\thetabar}[J]$ denotes
the $\bE_{\thetabar}$-submodule of $\ihat'^*_S\cW_{\thetabar}$ consisting of
elements annihilated by~$J$.
\end{prop}
\begin{proof}
  The existence of these
  isomorphisms is a consequence of
  Colmez's results in~\cite[IV.4]{MR2642409}, as we now explain. 
  It will be convenient to follow the exposition of~\cite{MR3267142}.

We first note that the isomorphism $i'_S$ gives rise to a formal \'etale~$(\varphi,
\Gamma)$-module $D_S$ with coefficients in~$S$, 
and Lemma~\ref{lem: L_infty to versal rings}, together with
Definition~\ref{definition of Linfty}, shows that
$$\ihat'^*_S L_{\infty} \iso 
   (D_{S}^{\natural}\boxtimes\Pone)(\zeta^{-1}\varepsilon\circ\det).
$$

Secondly, we recall that
the morphism $\cX_{\thetabar} \to \cX$ is defined so the induced family of
\'etale $(\varphi,\Gamma)$-modules over $\cX_{\thetabar}$ is the one that 
corresponds to the universal Galois representation $\cV_{\thetabar}$.
Consequently, we see that $\mathbf{V}(D_S)$ (the $2$-dimensional representation
of $G_{\Q_p}$ over $S$ determined by~$D_S$) coincides with $\ihat'^*_S \cV_{\thetabar}$. 
Our desired isomorphism~\eqref{eqn:sought-for-isomorphism-on-versal-rings}
can then be rewritten in the form
$$    \bigl(\mathbf{V}(D_{S})\cotimes_{\tld E_{\thetabar}} \tP_{\thetabar}\bigr)/\bigl(\mathbf{V}(D_{S})\cotimes_{\tld E_{\thetabar}} \tP_{\thetabar}\bigr)^{\SL_2(\Q_p)} \isoto 
    (D_{S}^{\natural}\boxtimes\Pone)(\zeta^{-1}\varepsilon\circ\det),$$
and it is in this form that we will prove it.

  Let~$A$ be an Artinian quotient
  of~$S$, and set~$D_A\coloneqq D_S\otimes_SA$. 
  Then $\mathbf{V}(D_{A})\otimes_{\tld
    E_{\thetabar}} \tP_{\thetabar}$ is a finite length object of~$\mathfrak{C}_{\thetabar} \subset \fC$,
  so it is the Pontrjagin dual~$\Pi^{\vee}$ of an object~$\Pi$ of the
  category~$\Rep_{\operatorname{tors}}(\zeta)$ considered
  in~\cite[Section~III.1]{MR3267142}.

  The~$(\varphi,\Gamma)$-module denoted~$\mathbb{D}(\Pi)$
  in~\cite{MR2642409} is, by the definition just preceding the statement of~\cite[Thm.\ IV.2.12]{MR2642409},
  isomorphic to the $(\varphi, \Gamma)$-module of $V(\Pi)^\vee \otimes \varepsilon$, where~$V$ is normalized as in~\cite{MR2642409}.
  This is the same as $V(\Pi)^\vee$ in our notation.
  Hence $\bD(\Pi)$ is equal to (in our notation)
  \[
  \bD(V(\Pi)^\vee) = \bD(\Vcheck(\Pi^\vee)\otimes \zeta^{-1}\varepsilon) = \bD(\Vcheck(\bV(D_A)\otimes_{\tE_{\thetabar}}\tP_{\thetabar})\otimes \zeta^{-1}\varepsilon).
  \]
    By Lemma~\ref{needed for universal deformation}~(1)
    (i.e.\ the very definitions of
$\tP_{\thetabar}$ and $E_{\thetabar}$), 
we have 
\[
\Vcheck(\mathbf{V}(D_{A})\otimes_{\tld E_{\thetabar}} \tP_{\thetabar})=\mathbf{V}(D_{A}),
\] 
and so we conclude that
\[
\bD(\Pi) = D_{A} \otimes \zeta^{-1}\varepsilon.
\]
Now~\cite[Prop.\ III.44]{MR3267142} produces a natural map of left $\cO\llbracket G\rrbracket _{\zeta^{-1}}$-modules
\[
\beta_{\Pone} : \Pi^\vee \to \bD(\Pi)^\sharp \boxtimes_{\zeta^{-1}} \Pone
\]
whose kernel is equal to $(\Pi^\vee)^{\SL_2(\bQ_p)}$.
Furthermore, by~\cite[Prop.\ II.1.11]{MR2642409}, we have
\[
\bD(\Pi)^\sharp \boxtimes_{\zeta^{-1}} \Pone = (D_A^\sharp \boxtimes \Pone) (\zeta^{-1}\varepsilon \circ \det),
\]
since by our conventions, $\boxtimes$ is shorthand for~$\boxtimes_{\zeta\varepsilon^{-2}}$ (compare Remark~\ref{rem: central char}). %

We thus obtain a map of left $\cO\llbracket G\rrbracket _{\zeta^{-1}}$-modules
\[\beta_{\Pone}:\mathbf{V}(D_{A})\otimes_{\tld
      E_{\thetabar}} \tP_{\thetabar}\to 
      (D_A^\sharp \boxtimes \Pone) (\zeta^{-1}\varepsilon \circ \det),\] 
  whose
  kernel is equal to $(\mathbf{V}(D_{A})\otimes_{\tld
      E_{\thetabar}} \tP_{\thetabar})^{\SL_2(\Qp)}$. We claim that~$\beta_{\Pone}$
  factors through $(D_{A}^{\natural}\boxtimes\Pone)(\zeta^{-1}\varepsilon \circ \det)$. 
  By~\cite[Prop.\ III.44(3)]{MR3267142}, this holds provided that~$\Pi^{\SL_2 (\Qp)}=0$, i.e.\
  provided that $(\mathbf{V}(D_{A})\otimes_{\tld E_{\thetabar}} \tP_{\thetabar})_{\SL_2
    (\Qp)}=0$.
 By construction,
$\mathbf{V}(D_{A})\otimes_{\tld E_{\thetabar}} \tP_{\thetabar}$ is a quotient of products
of copies of $\tP_{\thetabar}$.  By definition such a product never contains the projective envelope of a character
as a factor, and thus
$\mathbf{V}(D_{A})\otimes_{\tld E_{\thetabar}} \tP_{\thetabar}$ 
admits no non-zero map to
a character, and hence has no non-zero $\SL_2(\Q_p)$-invariant
quotients, as required. %

Thus we have obtained a natural morphism of left $\cO\llbracket G\rrbracket _{\zeta^{-1}}$-modules
\[\beta_{\Pone}:\mathbf{V}(D_{A})\otimes_{\tld E_{\thetabar}} \tP_{\thetabar}\to
  (D_{A}^{\natural}\boxtimes\Pone)(\zeta^{-1}\varepsilon\circ\det),\]whose kernel is equal to
$(\mathbf{V}(D_{A})\otimes_{\tld E_{\thetabar}}
\tP_{\thetabar})^{\SL_2(\Qp)}$. 
Furthermore, it is shown in the second paragraph
of the proof of~\cite[Thm.\ III.45]{MR3267142} that the cokernel
of~$\beta_{\Pone}$ is $\SL_2 (\Qp)$-invariant and of finite $\cO$-length.

Applying the above discussion to $A \coloneqq S/\fm_S^k$, we have for each~$k\ge 1$ a natural
exact sequence
\[0\to \ker_{k}\to \mathbf{V}(D_{S/\m_S^{k}})\otimes_{\tld E_{\thetabar}} \tP_{\thetabar}\stackrel{\beta_{\Pone}}\to
  (D_{S/\m_S^{k}}^{\natural}\boxtimes\Pone)(\zeta^{-1}\varepsilon\circ\det)\to\coker_{k}\to 0 \]where~$\ker_k$ and
~$\coker_k$ are of finite $\cO$-length. 
Thus by Mittag-Leffler we can pass
to the inverse limit and obtain an exact sequence
\[0\to \varprojlim_k\ker_k\to \mathbf{V}(D_{S})\cotimes_{\tld E_{\thetabar}}
  \tP_{\thetabar}\to (D_S^{\natural}\boxtimes\Pone)(\zeta^{-1}\varepsilon\circ\det)\to\varprojlim_k\coker_k\to
  0 \]whose kernel is equal to $(\mathbf{V}(D_{S})\otimes_{\tld E_{\thetabar}}
\tP_{\thetabar})^{\SL_2 (\Qp)}$, and whose cokernel is
$\SL_2(\Qp)$-invariant. This
cokernel vanishes by Corollary~\ref{cor:Steinberg D natural}, so after turning the left $\cO\llbracket G\rrbracket _{\zeta^{-1}}$-action into a right $\cO\llbracket G\rrbracket _\zeta$-action, we obtain the
isomorphism~\eqref{eqn:sought-for-isomorphism-on-versal-rings}.

If ~$\thetabar$ is not of type~\ref{item: Steinberg pseudorep}, then $(\mathbf{V}(D_{S})\cotimes_{\tld E_{\thetabar}}
  \tP_{\thetabar})^{\SL_2 (\Qp)}=0$, because there are no $\SL_2(\Qp)$-invariant objects in~$\mathfrak{C}_{\thetabar}$, so we obtain~\eqref{eqn:sought-for-isomorphism-on-versal-rings-non-St}.

  Finally, if ~$\thetabar$ of type~\ref{item: Steinberg pseudorep}, then 
  by~\eqref{eqn:Steinberg case rewrite bis} and the definition of~$\cW_{\thetabar}$
(i.e.\ Definition~\ref{defn:W-thetabar}),
we have a natural isomorphism of left $\cO\llbracket G\rrbracket _{\zeta^{-1}}$-modules 
  \[
\ihat'^*_S \cW_{\thetabar}
\cotimes_{\bE_{\thetabar}} \bP_{\thetabar}\iso
\ihat'^*_S \cV_{\thetabar}
\cotimes_{\bE_{\thetabar}} \bP_{\thetabar}\iso
\mathbf{V}(D_S)\cotimes_{\tE_{\thetabar}} \tP_{\thetabar}.\]
Thus, to prove~\eqref{eqn:sought-for-isomorphism-on-versal-rings-St},
it suffices to note the natural isomorphism
\[\ihat'^*_S\cW_{\thetabar}[J] \cotimes_{\bE_{\thetabar}} \bP_{\thetabar}\iso
(\ihat'^*_S\cW_{\thetabar} \cotimes_{\bE_{\thetabar}} \bP_{\thetabar})^{\SL_2 (\Qp)}\]
provided by Lemma~\ref{lem:SL_2 invariants via J torsion}.
\end{proof}

We can now give our alternative description of $\ihat'^*_{\thetabar}L_\infty$.

\begin{prop}
\label{prop:completed L-infinity}\leavevmode %
\begin{enumerate}%
\item\label{item:50} If~$\thetabar$ is {\em not} of type~{\em \ref{item: Steinberg block}},
  then there is an $\cO\llbracket G\rrbracket _{\zeta}$-equivariant isomorphism
  \begin{equation}
    \label{eqn:completed iso non Steinberg}
    \cV_{\thetabar} \cotimes_{\tld E_{\thetabar}} \tP_{\thetabar}\iso  \ihat'^*_{\thetabar}L_{\infty}.
  \end{equation}
  \item\label{item:51} If~$\thetabar$ is of type~{\em \ref{item: Steinberg block}},
    then there is an $\cO\llbracket G\rrbracket _{\zeta}$-equivariant isomorphism
\begin{equation}\label{eqn:completed iso Steinberg}
  (\cW_{\thetabar}/\cW_{\thetabar}[J])
    \cotimes_{\bE_{\thetabar}} \bP_{\thetabar}\iso \ihat'^*_{\thetabar}L_{\infty}.
  \end{equation}%
     \end{enumerate}
\end{prop}

\begin{proof}%
  We give the argument in case~\eqref{item:50}; case~\eqref{item:51} is formally identical, and we leave it to the reader. 
  We simplify our notation by writing $R\coloneqq  R^{\ps}_{\thetabar}$, with maximal ideal~$\mf{m}$. 

Let $H$ be a compact open subgroup of~$G$. 
Recall from Lemma~\ref{lem:finiteness of H coinvariants} %
that $(\cV_{\thetabar}\cotimes_{\tE_{\thetabar}}\tP_{\thetabar})_H$
and $(\ihat'^*_{\thetabar}L_{\infty})_H$
are objects of $\Coh(\cO_{\cX_{\thetabar}})$.
As in the statement of Proposition~\ref{prop:versal-ring-completed-L-infinity}, we let ~$S$ be a complete Noetherian local $\cO$-algebra with finite residue field, and $i'_S: \Spf S\to\cX_{\thetabar}$ 
   be a versal morphism at the closed point of~$\cX_{\thetabar}$.
   Since~$\ihat'^*_S$ is a right exact functor, it commutes with the formation of $H$-coinvariants, by Lemma~\ref{lem:abelian-composition-right-exact-tensor}. 
Thus the isomorphism~\eqref{eqn:sought-for-isomorphism-on-versal-rings-non-St}
induces an isomorphism
\begin{equation}
\label{eqn:H-coinvariant isos}
\ihat'^*_S\bigl( (\cV_{\thetabar}\cotimes_{\tE_{\thetabar}}\tP_{\thetabar})_H\bigr)
\iso \ihat'^*_S\bigl((\ihat'^*_{\thetabar}L_{\infty})_H\bigr),
\end{equation}
which is natural in 
the morphism~$i'_S$.
By Proposition~\ref{prop:versal formal functions CWE stacks}, %
 this collection of isomorphisms induces an isomorphism
\begin{equation}
\label{eqn:H-coinvariant isos bis}
(\cV_{\thetabar}\cotimes_{\tE_{\thetabar}}\tP_{\thetabar})_H
\iso (\ihat'^*_{\thetabar}L_{\infty})_H.
\end{equation}
Since the isomorphisms~\eqref{eqn:H-coinvariant isos}
are deduced from the $\cO\llbracket G\rrbracket_\zeta$-module 
isomorphisms~\eqref{eqn:sought-for-isomorphism-on-versal-rings-non-St},
they are %
compatible in $H$ as $H$-varies, and indeed
compatible with the $\cO\llbracket G\rrbracket_{\zeta}$-module structure on the 
projective systems
formed (by allowing $H$ to vary) by their left and right hand sides.
Passing to the limit over~$H$ in the isomorphisms~\eqref{eqn:H-coinvariant isos bis},
and taking into account the isomorphisms~\eqref{eqn:limit description of L_infty over cX} and~\eqref{eqn:limit description of V tensor P tilde},
we obtain the desired isomorphism~\eqref{eqn:completed iso non Steinberg}.
  \end{proof}%

  \subsubsection{An alternative description of \texorpdfstring{$\Functor_{\thetabar}$}{F-thetabar}}
  Using Proposition~\ref{prop:completed L-infinity}, we can now give
  an alternate description of~$\Functor_{\thetabar}$, and establish its full faithfulness.  
\begin{thm}
  \label{thm:derived-F-thetabar-nonSt-case}
  Let~$\thetabar$ be a $2$-dimensional $\cbF_p$-valued pseudorepresentation of~$G_{\bQ_p}$.
  \begin{enumerate}
    \item\label{item:79}  The functor $\Functor_{\thetabar}:D^b_{\fp}(\cA_{\thetabar})\to
\Pro D^b_{\coh}(\cX_{\thetabar})$ is fully faithful, and moreover factors through a \emph{(}necessarily fully faithful\emph{)} functor 
$D^b_{\fp}(\cA_{\thetabar})\to D^b_{\coh}(\cX_{\thetabar})$.
\item\label{item:40}If  $\thetabar$ is  {\em not} of type~\emph{\ref{item:
      Steinberg block}}, then  $\Functor_{\thetabar}$ is $t$-exact. %
      \item\label{item:41} If  $\thetabar$ is of type~\emph{\ref{item:
      Steinberg block}}, then  $\Functor_{\thetabar}$ has amplitude $[-1,0]$.
  \end{enumerate}
\end{thm}%
\begin{proof}We begin by supposing that  $\thetabar$ is not of type~\ref{item: Steinberg block}. 
We have natural isomorphisms of functors %
\begin{multline}
\label{eq:1}
\Functor_{\thetabar}(\text{--}) \xrightarrow[\eqref{eqn:F thetabar formula}]{\sim} \ihat'^*_{\thetabar}L_{\infty} \otimes^{L}_{\cO\llbracket G\rrbracket _{\zeta}}i_{\thetabar, *}(\text{--}) 
\xrightarrow[\eqref{eqn:completed iso non Steinberg}]{\sim} (\cV_{\thetabar} \cotimes_{\tld E_{\thetabar}} \tP_{\thetabar}) \otimes^{L}_{\cO\llbracket G\rrbracket _{\zeta}}i_{\thetabar, *}(\text{--}) 
\\
\xrightarrow[\textrm{Lem.\ }\ref{lem:tensor product comparison}]{\sim} \cV_{\thetabar} \cotimes^L_{\tld E_{\thetabar}} (\tP_{\thetabar} \otimes^{L}_{\cO\llbracket G\rrbracket _{\zeta}}i_{\thetabar, *}(\text{--})).
\end{multline}
By Lemma~\ref{lem:t-exactness on cA}, the functor
\[
\tP_{\thetabar}\otimes^L_{\cO\llbracket G\rrbracket _{\zeta}} i_{\thetabar, *}(\text{--}): D^b_{\fp}(\cA_{\thetabar}) \to \Pro D^b_{\fl}(\tE_{\thetabar})
\]
is $t$-exact.
By Lemma~\ref{lem:recapitulated Morita}~\eqref{item:63}, its restriction to the heart~$\cA_{\thetabar}^{\fp}$ 
is an equivalence
\[
\tP_{\thetabar}\otimes_{\cO\llbracket G\rrbracket _{\zeta}} i_{\thetabar, *}(\text{--}): \cA_{\thetabar}^{\fp} \iso \Mod^{\fl}(\tld E_{\thetabar}).
\] 
Hence, by Corollary~\ref{cor:t-exact derived functors}, the functor $\tP_{\thetabar}\otimes^L_{\cO\llbracket G\rrbracket _{\zeta}} i_{\thetabar, *}(\text{--})$ factors through a $t$-exact equivalence
\[
\tP_{\thetabar}\otimes^L_{\cO\llbracket G\rrbracket _{\zeta}} i_{\thetabar, *}(\text{--}): D^b_{\fp}(\cA_{\thetabar}) \iso D^b_{\fl}(\tE_{\thetabar}).
\]
It therefore suffices to prove that 
\[
\cV_{\thetabar} \cotimes^L_{\tE_{\thetabar}} \text{--} : \Pro D^b_{\fl}(\tE_{\thetabar}) \to \Pro D^b_{\coh}(\cX_{\thetabar})
\]
restricts to a $t$-exact and fully faithful functor $D^b_{\fl}(\tE_{\thetabar}) \to D^b_{\coh}(\cX_{\thetabar})$. 
To do so, we apply Lemma~\ref{lem:comparing-frak-and-cal-tensors --- derived case}, which produces a commutative diagram
\[
\begin{tikzcd}
	 D^b_{\fl}(\tE_{\thetabar}) &
        D^b_{\fp}(\tE_{\thetabar}) & \Pro D^b_{\fl}(\tE_{\thetabar}) \\
         D^b_{\coh}\bigl(\cX_{\thetabar}\bigr) & D^b_{\coh}\bigl(\fX_{\thetabar}\bigr)
      & \Pro D^b_{\coh}(\cX_{\thetabar}).
	\arrow["\eqref{eqn:fl to fp}", hook, from=1-1, to=1-2]
	\arrow["\cV_{\thetabar}\cotimes^{L}_{\tE_{\thetabar}}\text{--}",from=1-3, to=2-3]
	\arrow[hook, from=1-2, to=1-3, "\eqref{eqn:derived fp to pro fl}"]
	\arrow["\fV_{\thetabar}\otimes^{L}_{\tE_{\thetabar}}\text{--}",from=1-2, to=2-2]
	\arrow["\cV_{\thetabar}\otimes^{L}_{\tE_{\thetabar}}\text{--}",from=1-1, to=2-1]
        	\arrow[hook, from=2-1, to=2-2, "k_{\thetabar, *}"]
        \arrow["\cV_{\thetabar}\otimes^L_{\tE_{\thetabar}}\text{--}", from=1-2, to=2-3]
	\arrow[hook, from=2-2, to=2-3, "\widehat k_{\thetabar}^*"]
\end{tikzcd}
\]
Commutativity of the outer rectangle shows that $\cV_{\thetabar} \cotimes^L_{\tE_{\thetabar}}\text{--}$ restricts to a functor
\[
\cV_{\thetabar} \otimes^L_{\tE_{\thetabar}} \text{--}: D^b_{\fl}(\tE_{\thetabar}) \to D^b_{\coh}(\cX_{\thetabar}). 
\]
Bearing in mind the commutativity of the upper triangle, we now see that this restriction is fully faithful,
by Theorem~\ref{thm:JNWE-properties-of-V}~\eqref{item:20}, and $t$-exact, by Theorem~\ref{thm:JNWE-properties-of-V}~\eqref{item:22}.
This concludes the proof in the case that~$\thetabar$ does not have type~\ref{item: Steinberg pseudorep}.

We now turn to the case that $\thetabar$ is of type~\ref{item: Steinberg block}.
Using~\eqref{eqn:completed iso Steinberg} and arguing as above,
  we see that
  there is a natural isomorphism of functors from $D^b_{\fp}(\cA_{\thetabar})$ to
$\Pro D^b_{\coh}(\cX_{\thetabar})$
  \begin{equation}\label{eqn:Fthetabar-WJ}\Functor_{\thetabar} (\text{--}) \iso (\cW_{\thetabar}/\cW_{\thetabar}[J])
    \cotimes^L_{\bE_{\thetabar}} (\bP_{\thetabar} \otimes^L_{\cO\llbracket G\rrbracket _{\zeta}}
    i_{\thetabar, *}(\text{--})).\end{equation} 
    Again, the functor $\bP_{\thetabar} \otimes^L_{\cO\llbracket G\rrbracket _{\zeta}}
    i_{\thetabar, *}(\text{--})$ induces a $t$-exact equivalence
    $D^b_{\fp}(\cA_{\thetabar}) \iso D^b_{\fl}(\bE_{\thetabar})$. 
    We see %
    by Lemma~\ref{lem:comparing-frak-and-cal-tensors --- derived case}
    that there is a commutative diagram
\begin{equation}\label{eqn:specialization of tensor product commutative diagram to Steinberg case}
\begin{tikzcd}
	 D^b_{\fl}(\bE_{\thetabar}) &
        D^b_{\fp}(\bE_{\thetabar}) & \Pro D^b_{\fl}(\bE_{\thetabar}) \\
         D^b_{\coh}(\cX_{\thetabar}) & D^b_{\coh}(\fX_{\thetabar})
      & \Pro D^b_{\coh}(\cX_{\thetabar}).
	\arrow[hook, from=1-1, to=1-2]
	\arrow["\cW_{\thetabar}/\cW_{\thetabar}\text{$[J]$} \cotimes_{\bE_{\thetabar}}^L\text{--}",from=1-3, to=2-3]
	\arrow[hook, from=1-2, to=1-3, "\eqref{eqn:derived fp to pro fl}"]
	\arrow["\fW_{\thetabar}/\fW_{\thetabar}\text{$[J]$} \otimes_{\bE_{\thetabar}}^L\text{--}",from=1-2, to=2-2]
	\arrow["\cW_{\thetabar}/\cW_{\thetabar}\text{$[J]$} \otimes_{\bE_{\thetabar}}^L\text{--}",from=1-1, to=2-1]
        	\arrow[hook, from=2-1, to=2-2, "k_{\thetabar, *}"]
	\arrow[hook, from=2-2, to=2-3, "\widehat k_{\thetabar}^*"]
\end{tikzcd}
\end{equation}
This immediately implies that
  \begin{equation}
    \label{eqn:Fthetabar-WJ II}
    \Functor_{\thetabar} (\text{--}) \iso (\cW_{\thetabar}/\cW_{\thetabar}[J])
    \otimes^L_{\bE_{\thetabar}} (\bP_{\thetabar} \otimes^L_{\cO\llbracket G\rrbracket _{\zeta}}
    i_{\thetabar, *}(\text{--})).\end{equation} 
factors through~$D^b_{\coh}(\cX_{\thetabar})$.
Furthermore, by Lemma~\ref{lem:X* tor-dim'n}, the functor %
\[\bigl(\fW_{\thetabar}/\fW_{\thetabar}[J]\bigr)\otimes^L_{\bE_{\thetabar}} : D^b_{\fp}(\bE_{\thetabar}) \to
D^b_{\coh}(\fX_{\thetabar})\] has amplitude~$[-1,0]$,
so we see that indeed~$\Functor_{\thetabar}$ has amplitude~$[-1,0]$. %
Finally, the full faithfulness of~$\Functor_{\thetabar}$ follows as above, using~\eqref{eqn:Fthetabar-WJ} and Theorem~\ref{thm:alg-W-REnd-JWNE}. %
\end{proof}

In the context of Theorem~\ref{thm:derived-F-thetabar-nonSt-case}~\eqref{item:41} we have the following more precise description of the values of~$\Functor_{\thetabar}$ on irreducible objects of~$\cA_{\thetabar}$.
In its statement, we will identify coherent sheaves on $\fX_{\thetabar}$ with graded $S$-modules, where~$S$ is the graded ring in Section~\ref{subsubsec:Steinberg CWE stack},
and we will continue to denote by $\widehat k^*_{\thetabar} : D^b_{\coh}(\fX_{\thetabar}) \to \Pro D^b_{\coh}(\cX_{\thetabar})$ 
the $t$-exact, fully faithful completion functor.

\begin{cor}\label{cor:JNWE-calculation-of-functor-on-irreps-Steinberg-case}
  Let~$\thetabar = 1+\omega^{-1}$.
  Then: 
  \[\Functor_{\thetabar}(\pi_{\alpha})\cong \widehat k^*_{\thetabar}S(1)/(a_0, a_1, b_0, b_1,
  \varpi)[0],\] 
  \[\Functor_{\thetabar}(\St)  \cong \widehat k^*_{\thetabar}S(-1)/(a_0, a_1, c, \varpi)[0],\] 
  \[\Functor_{\thetabar}(\triv_G)\cong \widehat k^*_{\thetabar}S(-3)/(a_0, a_1, c, \varpi)[1].\]
\end{cor}
\begin{proof}
By Theorem~\ref{thm:alg-W-REnd-JWNE},
the graded $S$-module corresponding to $\fW_{\thetabar}/\fW_{\thetabar}[J]$ is the module~$X^*$ from Definition~\ref{defn:X*}. 
Bearing in mind the isomorphism~\eqref{eqn:Fthetabar-WJ II}, 
and the commutativity of the leftmost square in~\eqref{eqn:specialization of tensor product commutative diagram to Steinberg case}, %
the corollary thus follows from the computations of $\Tor_i^{\bE_{\thetabar}}(X^*, \text{--})$ in Example~\ref{example:functor on irreducible objects in the Steinberg case}. 
\end{proof}

We next note the evident extension of
Theorem~\ref{thm:derived-F-thetabar-nonSt-case}
to the case of a finite closed subset  $Y \subset X$.
\begin{cor}%
\label{cor:overview fully faithful on closed}
For any finite set of closed points $Y \subset X$,
the functor $\Functor_Y:D^b_{\fp}(\cA_{Y}) \to \Pro D^b_{\coh}(\cX_{Y})$ %
factors through a  fully faithful functor $D^b_{\fp}(\cA_{Y}) \to  D^b_{\coh}(\cX_{Y})$.
\end{cor}
\begin{proof}We can write~$\Functor_{Y}:D^b_{\fp}(\cA_{Y}) \to \Pro
  D^b_{\coh}(\cX_{Y})$ as a product 
  \[
  \prod_{\thetabar\in Y}\Functor_{\thetabar}: \prod_{\thetabar \in Y}D^b_{\fp}(\cA_{\thetabar}) \to \prod_{\thetabar \in Y}\Pro D^b_\coh(\cX_{\thetabar}),
  \]so the result 
   follows immediately from Theorem~\ref{thm:derived-F-thetabar-nonSt-case}.
\end{proof}

In view of Theorem~\ref{thm:derived-F-thetabar-nonSt-case}~\eqref{item:79}
and Corollary~\ref{cor:overview fully faithful on closed},
we will freely regard~ $\Functor_{\thetabar}$ (resp.\ $\Functor_Y$)
as a fully faithful functor $D^b_{\fp}(\cA_{\thetabar})\to D^b_{\coh}(\cX_{\thetabar})$
(resp.\ a fully faithful functor $D^b_{\fp}(\cA_Y)\to D^b_{\coh}(\cX_Y)$)
from now on.

\subsubsection{The relationship to Colmez's functor~\texorpdfstring{$V$}{V}}
Recall from Section~\ref{ColmezV}
that Colmez's functor~$V$ is exact, so in particular it induces a $t$-exact functor $D^b_{\fp}(\cA_{\thetabar})\to D^b_{\fl}(\tE_{\thetabar}) $. 
If $\thetabar$ is not of type~\ref{item: Steinberg block}, then by Lemma~\ref{needed for universal
deformation} %
we have a natural isomorphism of $t$-exact functors 
\[\tP_{\thetabar}\otimes^{L}_{\cO\llbracket G\rrbracket _{\zeta}}(\text{--})\iso V^\dagger(\text{--}),\]
 so we see from~\eqref{eq:1} that there is a natural isomorphism of $t$-exact
    functors from $D^b_{\fp}(\cA_{\thetabar})$ to
$D^b_{\coh}(\cX_{\thetabar})$
\begin{equation}
\label{eq:2}\Functor_{\thetabar} (\text{--}) \iso
\cV_{\thetabar} \cotimes^L_{\tld E_{\thetabar}} (\tP_{\thetabar}
\otimes^L_{\cO\llbracket G\rrbracket _{\zeta}} \text{--}) \iso \cV_{\thetabar} \cotimes^L_{\tld
  E_{\thetabar}} V^\dagger(\text{--}).
\end{equation}

If $\thetabar$ is of type~\ref{item: Steinberg
      block} then~\eqref{eq:2} does not hold in general, and~$\Functor_{\thetabar}$ is no longer $t$-exact.
      However, we have the following related result, %
which will be used in the proof of
Proposition~\ref{concentration and support} below. %
\begin{cor}%
  \label{cor:evaluating-H0-functor-no-trivial-quotient}
Suppose that $\thetabar$ is of type~\emph{\ref{item: Steinberg block}}, and
that~$\pi\in\cA_{\thetabar}^{\fp}$ admits no $\SL_2(\Qp)$-invariant quotient. Then
there is a natural isomorphism 
  \[H^0\Functor_{\thetabar} (\pi) \iso \cV_{\thetabar} \cotimes_{\tld E_{\thetabar}}V^\dagger(\pi).\]
\end{cor}%
\begin{proof}By 
  Lemma~\ref{right exactness in the module variable} (noting that products are exact in the category of pro-coherent sheaves on a formal algebraic stack)
  and~\eqref{eqn:Fthetabar-WJ}, %
  we have a right exact sequence %
\[\cW_{\thetabar}[J]
    \cotimes_{\bE_{\thetabar}} (\bP_{\thetabar} \otimes_{\cO\llbracket G\rrbracket _{\zeta}}
    \pi)\to \cW_{\thetabar}
    \cotimes_{\bE_{\thetabar}} (\bP_{\thetabar} \otimes_{\cO\llbracket G\rrbracket _{\zeta}}
    \pi) \to H^0\Functor_{\thetabar} (\pi)\to 0.\]
The leftmost term vanishes, because
  \begin{align*}
\cW_{\thetabar}[J]
  \cotimes_{\bE_{\thetabar}} (\bP_{\thetabar} \otimes_{\cO\llbracket G\rrbracket _{\zeta}} \pi)
  &= (\cW_{\thetabar}[J] \cotimes_{\bE_{\thetabar}} \bP_{\thetabar})
     \otimes_{\cO\llbracket G\rrbracket _{\zeta}} \pi \\
  &= (\cW_{\thetabar}[J] \cotimes_{\bE_{\thetabar}} \bP_{\thetabar}/J)
     \otimes_{\cO\llbracket G\rrbracket _{\zeta}} \pi \\
  &= \cW_{\thetabar}[J]
     \cotimes_{\bE_{\thetabar}} (\bP_{\thetabar}/J \otimes_{\cO\llbracket G\rrbracket _{\zeta}} \pi) \\
  &= \cW_{\thetabar}[J]
     \cotimes_{\bE_{\thetabar}} (\cO \otimes_{\cO\llbracket G\rrbracket _{\zeta}} \pi),
  \end{align*}
  which vanishes by our assumption that~$\pi$ admits no $\SL_2(\Qp)$-invariant quotient.
 (The second equality holds because, by construction, the map $\bE^{\op}_{\thetabar} \to \End_{\Pro \Coh(\cX_{\thetabar})}(\cW_{\thetabar}[J])$
  factors through $\bE_{\thetabar}^{\op}/J$, hence the functor $\cW_{\thetabar}[J] \cotimes_{\bE_{\thetabar}}$ factors through 
  $\Mod_c(\bE_{\thetabar}) \to \Mod_c(\bE_{\thetabar}/J), M \mapsto M/J$.)
  
  It follows from the discussion above that $H^0\Functor_{\thetabar}(\pi)$ is naturally isomorphic to $\cW_{\thetabar}
    \cotimes_{\bE_{\thetabar}} (\bP_{\thetabar} \otimes_{\cO\llbracket G\rrbracket _{\zeta}}
    \pi)$, which in turn is isomorphic to $\cV_{\thetabar}
    \cotimes_{\tE_{\thetabar}} (\tP_{\thetabar} \otimes_{\cO\llbracket G\rrbracket _{\zeta}}
    \pi)$, by~\eqref{comparing cV and cW}.
  The corollary now follows from Lemma~\ref{needed for universal deformation}~(2).
\end{proof}

\subsection{Finiteness properties of~\texorpdfstring{$\Functor$}{F}}\label{finiteness for Functor}

Our next goal is to prove Proposition~\ref{finiteness of cohomology}, which shows that
our functor~$\Functor$ takes finitely presented representations to (bounded complexes of) coherent
sheaves. 
Since one of the main steps of the argument will be to establish finiteness and concentration properties of the completions of $\Functor(\cInd_{KZ}^G \sigma)$
at closed points, we begin
by establishing a compatibility between $\Functor$ and
~$\Functor_{Y}$, where $Y$ is a finite set of closed points of~$X$. %

\subsubsection{Compatibility between~\texorpdfstring{$\Functor$}{F} and~\texorpdfstring{$\Functor_Y$}{FY}.}\label{subsubsec:compatibility of functor and completion}
Let~$Y$ be a finite set of closed points of~$X$.
Then Corollary~\ref{cor:overview fully faithful on closed} shows that the functor~$\Functor_Y$ factors through a fully faithful functor $D^b_{\fp}(\cA_Y) \to D^b_{\coh}(\cX_{Y})$, 
 and we write  $\ProFunctor_{Y}: \Pro D_{\fp}^b(\cA_{Y}) \to \Pro D^b_{\coh}(\cX_{Y})$ for its $\Pro$-extension \begin{equation}
  \label{eqn:ProF_Y}
  \ProFunctor_{Y} = \Pro(\ihat'^*_Y \Functor \, \ihat_{Y,*})= \ihat'^*_Y \Pro (\Functor) \ihat_{Y,*},
\end{equation} %
 which is also  fully faithful by %
Lemma~\ref{lem:F-fully-faithful-iff-IndProF}.
Note that since $\Functor = L_\infty \otimes_{\cO\llbracket G\rrbracket _\zeta}^L\text{--}$, the  Pro-extension of~$\Functor$ is the functor
\[
L_\infty \ccotimes_{\cO\llbracket G\rrbracket _\zeta}^L\text{--}: \Pro D^b_{\fp}(\cA) \to \Pro D^b_{\coh}(\cX)
\]
from Lemma~\ref{lem:derived-tensor-product-EW}~\eqref{item:93}.
If we precompose~\eqref{eqn:ProF_Y}
 with~$\ihat_{Y}^*$, the unit of the adjunction gives a natural transformation of functors $D^b_{\fp}(\cA)\to\Pro
  D^b_{\coh}(\cX_{Y}) $ %
\begin{equation} \label{unit of adjunction on functor}
\ihat'^*_{Y} \Functor \to (\ProFunctor_{Y}) \ihat_{Y}^*
\end{equation} 
which identifies with the natural transformation %
\[
\ihat'^*_Y(L_\infty \otimes_{\cO\llbracket G\rrbracket _\zeta}^L (\text{--})) \to \ihat'^*_Y(L_\infty \ccotimes^L_{\cO\llbracket G\rrbracket _\zeta} \ihat_{Y, *}\ihat_{Y}^*(\text{--}))
\]
induced by the restriction to $D^b_{\fp}(\cA)$ of the unit of adjunction $\id_{\Pro D^b_{\fp}(\cA)} \to \ihat_{Y, *} \ihat^*_Y$.
(Note that $L_\infty \ccotimes^L_{\cO\llbracket G\rrbracket _\zeta}$ restricts to $L_\infty \otimes^L_{\cO\llbracket G\rrbracket _\zeta}$ on~$D^b_{\fp}(\cA)$, by the discussion above, 
or equivalently by~\eqref{eqn:second-commutative-diagram-derived-tensor}.)

Applying Lemma~\ref{lem:pullback-Linfty-associativity}, we can rewrite \eqref{unit of adjunction on functor} as a natural transformation
\begin{equation}\label{to rewrite using derived completeness}
\ihat'^*_YL_\infty \otimes_{\cO\llbracket G\rrbracket _\zeta}^L (\text{--}) \to \ihat'^*_YL_\infty \ccotimes^L_{\cO\llbracket G\rrbracket _\zeta} \ihat_{Y, *}\ihat_{Y}^*(\text{--}).
\end{equation}
Assume now that~$Y = \{\thetabar\}$ is a singleton, 
and apply Proposition~\ref{prop:completed L-infinity}. 
This gives us a projective generator~$P$ of~$\fC_{\thetabar}$ with finite cosocle, with compact endomorphism ring~$E$, and a
complete right $E$-module $V$ in $\Pro \Coh(\cX_{\thetabar})$, such that %
\[
V \cotimes^L_{E} P \isoto V \cotimes_{E} P \isoto \ihat'^*_{\thetabar}L_\infty
\](where the first isomorphism follows from the fact that~$P$ is projective in $\Mod_c(E)$, by Lemma~\ref{lem:projective generators are flat over their endos}~\eqref{item:projective 3}).
By Lemma~\ref{lem:tensor product comparison}, the natural transformation~\eqref{to rewrite using derived completeness} is thus isomorphic to
\begin{equation}\label{to be used in proof of completing on Y}
V \cotimes_E^L (P \otimes_{\cO\llbracket G\rrbracket_{\zeta} }^L (\text{--})) \to V \cotimes_E^L (P \ccotimes_{\cO\llbracket G\rrbracket _\zeta}^L  \ihat_{Y, *}\ihat_{Y}^*(\text{--})).
\end{equation}

\begin{prop}
\label{prop:completing our functor on Y}%
For any finite set of closed points $Y \subset X$, 
the natural transformations~\eqref{to be used in proof of completing on Y} and $\eqref{unit of adjunction on functor}: \ihat'^*_{Y} \Functor \to (\ProFunctor_{Y})\ihat_{Y}^*$
are isomorphisms.
\end{prop}
\begin{proof} As in the proof of Corollary~\ref{cor:overview fully faithful on closed}, we may immediately reduce to the case that~$Y=\thetabar$ is a singleton.
Since~\eqref{unit of adjunction on functor} is naturally isomorphic to~\eqref{to be used in proof of completing on Y}, it suffices to prove that~\eqref{to be used in proof of completing on Y}
is an isomorphism.
To do so, it suffices to prove that for all
projective $P \in \fC_{\thetabar}$ with finite cosocle,
the natural transformation
\begin{equation}\label{to be used in proof of completing on Y II}
P \otimes^L_{\cO\llbracket G\rrbracket _\zeta} (\text{--}) \to P \ccotimes_{\cO\llbracket G\rrbracket _\zeta}^L \ihat_{\thetabar,*}\ihat_{\thetabar}^*(\text{--})
\end{equation}
is an isomorphism.
We know that $P \otimes_{\cO\llbracket G\rrbracket _\zeta}^L\text{--} : D^b_{\fp}(\cA) \to \Pro D^b_{\fl}(E)$ is $t$-exact, by 
Lemma~\ref{lem:t-exactness on cA}.
Hence its $\Pro$-extension $P \ccotimes_{\cO\llbracket G\rrbracket _\zeta}^L\text{--}$ is also $t$-exact, and so both sides of~\eqref{to be used in proof of completing on Y II} are $t$-exact functors. 
By Corollary~\ref{cor:t-exact derived functors},
it thus suffices to prove that~\eqref{to be used in proof of completing on Y II} becomes an
isomorphism after restriction to~$\cA^{\fp}$.
This is Lemma~\ref{tensor and completion}. %
\end{proof}

We can now prove that $\Functor(\cInd_{KZ}^G \sigma)$ is pure of degree zero (in the sense of Definition~\ref{defn:pure-degree-zero}) after completed pullback to~$\cX_{\thetabar}$.

\begin{lemma}\label{concentration}
Let~$\sigma$ be a Serre weight.
Let~$\thetabar$ be a %
2-dimensional $\cbF_p$-valued pseudorepresentation.
Then 
\[
\ihat'^*_{\thetabar}\Functor(\cInd_{KZ}^G\sigma) \buildrel \textrm{\emph{Prop.}\ }\ref{prop:completing our functor on Y} \over = \ProFunctor_{\thetabar}(\ihat_{\thetabar}^*\cInd_{KZ}^G\sigma)
\]
is a pro-coherent sheaf, i.e.\ it is pure of degree zero.
Furthermore, $\Functor_{\thetabar}(\cInd_{KZ}^G\sigma/f_{\thetabar}\cInd_{KZ}^G\sigma)$ is pure of degree zero unless $\thetabar$ is of type~\emph{\ref{item: Steinberg pseudorep}} and~$\sigma$ is Steinberg.
\end{lemma}
\begin{proof}
This is immediate if $\thetabar$ is not of type~\ref{item: Steinberg pseudorep}, since $\Functor_{\thetabar}$ (and thus~$\ProFunctor_{\thetabar}$) is
$t$-exact by
Theorem~\ref{thm:derived-F-thetabar-nonSt-case}~\eqref{item:40}. %
Note that in this case the only piece of information we need about
$\ihat_{\thetabar}^*\cInd_{KZ}^G\sigma$ 
and
$\ihat_{\thetabar}^*(\cInd_{KZ}^G\sigma / f_{\thetabar} \cInd_{KZ}^G\sigma)$ 
is that %
they are objects
of~$\Pro \cA_{\thetabar}^{\fp}$, a consequence of the $t$-exactness of~$\ihat_{\thetabar}^*$. 

Suppose now that~$\thetabar$ has type~ \ref{item: Steinberg pseudorep}.
After twisting we may suppose that $\thetabar = 1+\omega^{-1}$, so that by
Lemma~\ref{lem:localization-of-cInd-sigma}~\eqref{item:support of cInd sigma}, 
we only need to consider  $\sigma
\in \{\Sym^0, \det \otimes \Sym^{p-3},\Sym^{p-1}\}$, as otherwise
$\ihat_{\thetabar}^*\cInd_{KZ}^G\sigma=0$. To handle the cases $\sigma = \Sym^0$ and $\det\otimes\Sym^{p-3}$,
we first recall that by definition we have~$f_{\thetabar}=\HeckeT-1$, and
  by~\eqref{eqn:recollection-of-explicit-completion-of-cInd-sigma} we have
  \[\ihat_{\thetabar}^{*}\cInd_{KZ}^G \sigma\iso \quoteslim{n}(\cInd_{KZ}^G
    \sigma )/(\HeckeT-1)^{n}.\]  
Then, by  an obvious d\'evissage,
it suffices to show that
  $\Functor_{\thetabar}\bigl((\cInd_{KZ}^G \sigma )/(\HeckeT-1)\bigr)$
is pure of degree zero.
If~ $\sigma=\det\otimes\Sym^{p-3}$, we
have \[\Functor_{\thetabar}\bigl(\cInd_{KZ}^G \det \otimes
  \Sym^{p-3}/(\HeckeT-1)\bigr)=\Functor_{\thetabar}(\pi_\alpha),\] which is pure of
degree zero by
Corollary~\ref{cor:JNWE-calculation-of-functor-on-irreps-Steinberg-case}.

In the
case~$\sigma=\Sym^0$, we have a non-split short exact sequence
\begin{equation}\label{eqn:JH-seq-of-cInd-trivial}0 \to \St \to
\cInd_{KZ}^G\Sym^0/(\HeckeT-1) \to \triv_G \to 0,\end{equation} so that
$\cInd_{KZ}^G\Sym^0/(\HeckeT-1)$ is the cone of a non-zero map $\triv_G[-1]\to\St$.  It then follows from
Corollary~\ref{cor:JNWE-calculation-of-functor-on-irreps-Steinberg-case} and the
full faithfulness of~$\Functor_{\thetabar}$ that~$\Functor_{\thetabar}(\cInd_{KZ}^G\Sym^0/(\HeckeT-1))$ is
quasi-isomorphic to the %
image under the $t$-exact, fully faithful functor $\widehat k^*_{\thetabar} : D^b_{\coh}(\fX_{\thetabar}) \to \Pro D^b_{\coh}(\cX_{\thetabar})$
of a complex (in degrees~$-1$ and~$0$) %
\[
S(-3)/(a_0, a_1, c, \varpi) \xrightarrow{\delta} S(-1)/(a_0, a_1, c, \varpi),
\]
where~$\delta$ is non-zero, and we regard graded~$S$-modules as objects of $\Coh(\fX_{\thetabar})$. %
Since every non-zero $S$-linear endomorphism of $S/(a_0, a_1, c, \varpi) =
\bF[b_0, b_1]$ is injective (such an endomorphism is given by multiplication by a non-zero
element of the integral domain~$\bF[b_0,b_1]$), we deduce that
~$\Functor_{\thetabar}(\cInd_{KZ}^G\Sym^0/(\HeckeT-1))$ is pure of degree
zero.

Finally, suppose that~$\sigma=\Sym^{p-1}$. %
Since~$\ihat_{\thetabar}^{*}$ is exact, applying it 
to~\eqref{Steinberg change of weight}
yields a short exact sequence
\begin{equation}\label{ses-comparing-0-p-minus-1}
0\to \ihat_{\thetabar}^*\cInd_{KZ}^G\Sym^0\to
\ihat_{\thetabar}^*\cInd_{KZ}^G\Sym^{p-1}\to
\St\to 0. \end{equation}
Now 
Corollary~\ref{cor:JNWE-calculation-of-functor-on-irreps-Steinberg-case}
shows that $\Functor_{\thetabar}(\St)$ is pure of degree zero,
and since we have already seen that~$\ProFunctor_{\thetabar}(\ihat_{\thetabar}^*\cInd_{KZ}^G \Sym^0)$ is pure of degree zero,
the same is true of~$\ProFunctor_{\thetabar}(\ihat_{\thetabar}^*\cInd_{KZ}^G \Sym^{p-1})$.
\end{proof}

\subsubsection{Completed stalks of \texorpdfstring{$\Functor(\cInd_{KZ}^G \sigma)$}{F on compact inductions of Serre weights}.}
Recall from Section~\ref{subsubsec:points of X} that the closed points of~$\cX$ are in bijection with $\Gal(\cbF_p/\bF)$-conjugacy classes of
2-dimensional $\cbF_p$-valued pseudorepresentations~$\thetabar$, or equivalently of
semisimple Galois representations $\rhobar: G_{\bQ_p} \to \GL_2(\cbF_p)$.
Given~$\rhobar$ with trace~$\thetabar$, we write %
$i'_{\rhobar} : (\cX/\varpi)^\wedge_{\rhobar} \to \cX/\varpi$ for the completion at the closed subset $\{\rhobar\} \subset |\cX/\varpi|$.
In this subsection we study
the pullback $\ihat'^*_{\rhobar} \Functor(\cInd_{KZ}^G \sigma)$.

We begin by specializing the framework of Section~\ref{subsubsec:support}
to the situation at hand.
In Definition~\ref{defn:various versal rings} we have constructed a morphism
\[
\Spf R^{\ver}_{\thetabar} \to \cX_{\thetabar}
\]
mapping the unique point of the source to the unique closed point~$\rhobar$ of the target.
It thus induces a morphism
\begin{equation}\label{eqn:versal morphism to CWE stacks}
v: \Spf R^{\ver}_{\thetabar}/\varpi \to (\cX/\varpi)^\wedge_{\rhobar}.
\end{equation}

If~$\sigma$ is a Serre weight, we have defined
a reduced closed algebraic substack $z : \cZ(\sigma) \to \cX/\varpi$ in Definition~\ref{defn:special-fibre-crystalline-weight-sigma}.
We also have introduced in Definition~\ref{defn:various versal rings} a quotient~$\Rsigma_{\thetabar}$ of~$R^{\ver}_{\thetabar}$ such that 
$\Spf \Rsigma_{\thetabar} = \Spf R^{\ver}_{\thetabar}/\varpi \times_{\cX/\varpi} \cZ(\sigma)$.
Noting that 
we may factor~$i'_{\rhobar}$ as a composite %
\[
(\cX/\varpi)^\wedge_{\rhobar} \xrightarrow{i'_{\rhobar,\thetabar}} \cX_{\thetabar}/\varpi \xrightarrow{i'_{\thetabar}} \cX/\varpi,
\]
we see that we have a commutative diagram   %

\begin{equation}\label{completed stalks diagram}
  \begin{tikzcd}
\Spf R^{\ver}_{\thetabar}/\varpi \arrow[r, "v"] & 
(\cX/\varpi)^\wedge_{\rhobar} \arrow[r, "i'_{\rhobar,\thetabar}"] \arrow[rr, "i'_{\rhobar}", bend left] & 
\cX_{\thetabar}/\varpi \arrow[r, "i'_{\thetabar}"] & 
\cX/\varpi \\
\Spf \Rsigma_{\thetabar}/\varpi \arrow[u] \arrow[r, "v_{\sigma}"] & 
\cZ(\sigma)^\wedge_{\rhobar} \arrow[rr, "i'_{\rhobar,\sigma}"] \arrow[u, "z_{\rhobar}"] & 
& 
\cZ(\sigma) \arrow[u, "z"],
\end{tikzcd}
\end{equation}
with Cartesian squares, where ~$i'_{\rhobar,\sigma}$ is the completion at $\{\rhobar\} \cap |\cZ(\sigma)|$, 
and $z_{\rhobar}$ is the completion of~$z$ at~$\{\rhobar\}$.
(If~$\rhobar \not \in |\cZ(\sigma)|$, then $\cZ(\sigma)^\wedge_{\rhobar}$ is empty.)

We now note that $\ihat'^*_{\rhobar} \Functor(\cInd_{KZ}^G\sigma)$ is in fact a pro-coherent sheaf,
i.e.\ pure of degree zero. Indeed, 
Lemma~\ref{concentration} shows that
$\ihat'^*_{\thetabar} \Functor(\cInd_{KZ}^G\sigma)$ is pure of degree zero,
while the discussion of Section~\ref{subsubsec:derived passage}
shows that $\ihat'^*_{\rhobar,\thetabar}$ 
is $t$-exact.
By 
 Remark~\ref{completed pullback is functorial} we have
$\ihat'^*_{\rhobar} = \ihat'^*_{\rhobar,\thetabar}\ihat'^*_{\thetabar}$, hence  
\begin{equation}\label{eqn:compare-ihatprime-theta-rho}
  \ihat'^*_{\rhobar} \Functor(\cInd_{KZ}^G\sigma)
  \iso \ihat'^*_{\rhobar,\thetabar}\ihat'^*_{\thetabar}\Functor(\cInd_{KZ}^G\sigma)
\end{equation}
is again pure of degree zero.

As a special case of Lemma~\ref{lem: versal pullback is exact and faithful}, we have an exact and faithful pullback functor
\[
\vhat^* : \Pro\Coh\bigl( (\cX/\varpi)^\wedge_{\rhobar}\bigr)
\to \Mod_c(R^{\ver}_{\thetabar}/\varpi).
\]
We  apply $\vhat^*$ to the object 
$\ihat'^*_{\rhobar} \Functor(\cInd_{KZ}^G\sigma)$
of $\Pro\Coh\bigl( (\cX/\varpi)^\wedge_{\rhobar}\bigr),$
so as to obtain an object of
$\Mod_c(R^{\ver}_{\thetabar}/\varpi)$.
The following proposition 
provides some quite precise information about this latter object.

\begin{prop}\label{computation of completed stalk}
Let~$\sigma$ be a Serre weight, and let $\rhobar : G_{\bQ_p} \to \GL_2(\cbF_p)$ be a semisimple Galois representation.
Let~$\thetabar$ be the trace of~$\rhobar$.
\begin{enumerate}
\item If $\rhobar \not\in |\cZ(\sigma)|$, then $\ihat'^*_{\rhobar} \Functor(\cInd_{KZ}^G\sigma) = 0$.
\item If $\rhobar \in |\cZ(\sigma)|$, then $\vhat^*\ihat'^*_{\rhobar} \Functor(\cInd_{KZ}^G\sigma)$ is {\em (}the compact $R_{\thetabar}/\varpi$-module associated to{\em )}
a finitely presented $\Rsigma_{\thetabar}/\varpi$-module, and the natural map
\[
\Rsigma_{\thetabar}/\varpi \to \End_{\Rsigma_{\thetabar}/\varpi}\bigl(\vhat^*\ihat'^*_{\rhobar} \Functor(\cInd_{KZ}^G\sigma)\bigr)
\]
is injective.
\item If $\rhobar \in |\cZ(\sigma)|$ and $\sigma \ne \det^a \otimes \Sym^{p-2}$ for any~$a$, then 
\[
\vhat^*\ihat'^*_{\rhobar} \Functor(\cInd_{KZ}^G\sigma) \cong \Rsigma_{\thetabar}/\varpi.
\]
\end{enumerate}
\end{prop}
\begin{proof}

The first statement is immediate from~\eqref{eqn:compare-ihatprime-theta-rho} together with Proposition~\ref{prop:completing our functor on Y}, since if~$\rhobar \not \in |\cZ(\sigma)|$, then~$\sigma$ is not a Serre weight of~$\rhobar$, 
by Theorem~\ref{thm: fixed determinant stack for GL2 Qp}~\eqref{item:points and Serre weights}, 
and so $\ihat^*_{\thetabar}\cInd_{KZ}^G\sigma = 0$, by Lemma~\ref{lem:localization-of-cInd-sigma}~\eqref{item:support of cInd sigma}.

We now prove the remaining two statements, initially under the additional assumption that
 if~$\thetabar$ has type~\ref{item: Steinberg pseudorep}, then $\sigma$ is not a determinant twist of~$\Sym^0$.
Under this additional assumption, we claim that
\begin{equation}
\label{eqn:non Sym^0 case}
\ProFunctor_{\thetabar}(\ihat_{\thetabar}^*\cInd_{KZ}^G \sigma) 
\iso \cV_{\thetabar}\cotimes_{\tE_{\thetabar}} \Pro V^\dagger(\ihat_{\thetabar}^*
\cInd_{KZ}^G\sigma).
\end{equation}
Indeed, we first note that Lemma~\ref{concentration} shows that
$$
\ProFunctor_{\thetabar}(\ihat_{\thetabar}^*\cInd_{KZ}^G \sigma) 
\iso
H^0 \ProFunctor_{\thetabar}(\ihat_{\thetabar}^*\cInd_{KZ}^G \sigma),$$
so that we may replace the left-hand side of~\eqref{eqn:non Sym^0 case}
by its~$H^0$.

If $\thetabar$ does not have type~\ref{item: Steinberg pseudorep},
we then see that the
isomorphism~\eqref{eqn:non Sym^0 case} follows by passing to $H^0$
in the (Pro-extension of) the isomorphism~\eqref{eq:2}.
On the other hand,
if $\thetabar$ does have type~\ref{item: Steinberg pseudorep},
then we are assuming that $\sigma$ is a twist of either~$\Sym^{p-3}$ 
or~$\Sym^{p-1}$, 
in which case
$\cInd_{KZ}^G \sigma/(\HeckeT-1)^n$ does not have any $\SL_2(\bQ_p)$-invariant quotient, for any~$n > 0$.
We may thus apply 
Corollary~\ref{cor:evaluating-H0-functor-no-trivial-quotient}
to find that
\begin{multline*}
H^0\ProFunctor_{\thetabar}(\ihat^*_{\thetabar}\cInd_{KZ}^G \sigma)
\iso
\quoteslim{n}H^0\Functor_{\thetabar}\bigl(\cInd_{KZ}^G \sigma/(\HeckeT-1)^n\bigr)
\\
\iso
\quoteslim{n}\cV_{\thetabar} \cotimes_{\tE_{\thetabar}}V^\dagger\bigl(\cInd_{KZ}^G \sigma/(\HeckeT-1)^n\bigr) 
\iso
\cV_{\thetabar} \cotimes_{\tE_{\thetabar}}\Pro V^\dagger(\ihat^*_{\thetabar}\cInd_{KZ}^G \sigma),
\end{multline*}
again yielding the desired isomorphism~\eqref{eqn:non Sym^0 case}.

We then find that %
\begin{multline}\label{eqn:functor-as-tensor-over-Galois-group-ring}
\vhat^*\ihat'^*_{\rhobar} \Functor(\cInd_{KZ}^G \sigma) = \vhat^*\ihat'^*_{\rhobar, \thetabar} \ihat'^*_{\thetabar}\Functor(\cInd_{KZ}^G \sigma) %
\\
\buildrel \textrm{Prop.\ }\ref{prop:completing our functor on Y} \over = 
\vhat^*\ihat'^*_{\rhobar,\thetabar}\ProFunctor_{\thetabar}(\ihat_{\thetabar}^*\cInd_{KZ}^G\sigma)
\buildrel
\eqref{eqn:non Sym^0 case}
\over = 
\vhat^*\ihat'^*_{\rhobar,\thetabar} 
\bigl(\cV_{\thetabar} \cotimes_{\tld E_{\thetabar}}\Pro V^\dagger(\ihat_{\thetabar}^*\cInd_{KZ}^G\sigma)\bigr) \\
\buildrel \textrm{Lem.\ \ref{lem:abelian-composition-right-exact-tensor}} \over =
(\vhat^*\ihat'^*_{\rhobar,\thetabar} 
\cV_{\thetabar} )\cotimes_{\tld E_{\thetabar}}\Pro V^\dagger(\ihat_{\thetabar}^*\cInd_{KZ}^G\sigma)
=
\Vver_{\thetabar} \cotimes_{\tld E_{\thetabar}}\Pro V^\dagger(\ihat_{\thetabar}^*\cInd_{KZ}^G\sigma), %
\end{multline}
where we have followed Definition~\ref{defn:various versal rings} in writing~$\Vver_{\thetabar}$ for the versal object on $R^{\ver}_{\thetabar}$.
Under our running assumption that if $\thetabar$ has type~\ref{item: Steinberg pseudorep}, then $\sigma$ is not a determinant twist of~$\Sym^0$, the proposition follows from Proposition~\ref{prop:support of M(sigma)} (noting that~$\Vsigma^\dagger$ there is by definition equal to $\Pro V^\dagger(\ihat_{\thetabar}^*\cInd_{KZ}^G\sigma)$).

It remains to treat the case when $\thetabar$ is of type~\ref{item:
    Steinberg pseudorep} and $\sigma$ is a twist of~$\Sym^0$,
where, as usual, we make a twist so as to assume that in fact
$\thetabar=1+\omega^{-1}$ and~$\sigma = \Sym^0$.
We need to show that 
\[
\vhat^*\ihat'^*_{\rhobar}\Functor(\cInd_{KZ}^G \Sym^0) \cong R^{(1, 0), \crys}/\varpi.
\]
Recall %
that applying~$\ihat^*_{\thetabar}$ to~\eqref{Steinberg change of weight} and~\eqref{trivial change of weight} 
gives rise to a sequence of injections
\[
\ihat_{\thetabar}^*\cInd_{KZ}^G \Sym^{p-1} \to \ihat_{\thetabar}^*\cInd_{KZ}^G \Sym^0 \to \ihat_{\thetabar}^*\cInd_{KZ}^G \Sym^{p-1}
\]
such that the composition is a unit multiple of $\HeckeT-1$, the cokernel of the first map is~$\triv_G$, and the cokernel of the second map is~$\St$.
Then Corollary~\ref{cor:JNWE-calculation-of-functor-on-irreps-Steinberg-case} and Lemma~\ref{concentration}
 yield a sequence of morphisms %
\begin{multline}\label{eqn:profunctor change of weight}
\ProFunctor_{\thetabar}(\ihat_{\thetabar}^*\cInd_{KZ}^G\Sym^{p-1}) \onto \ProFunctor_{\thetabar}(\ihat_{\thetabar}^*\cInd_{KZ}^G\Sym^0)
\\ 
\into 
\ProFunctor_{\thetabar}(\ihat_{\thetabar}^*\cInd_{KZ}^G\Sym^{p-1})
\end{multline}
whose composite is a unit multiple of~$\HeckeT-1$;
note that the second arrow is (as indicated) an injection, since $\Functor_{\thetabar}(\St)$ is pure of degree~$0$, and the
first arrow is (as indicated) a surjection, since $\Functor_{\thetabar}(\triv_G)$ is pure of degree~$-1$.

Applying~$\vhat^*\ihat'^*_{\rhobar, \thetabar}$ to~\eqref{eqn:profunctor change of weight}, and using the case~$\sigma=\Sym^{p-1}$ of the proposition, 
which we have already proved, we obtain  (in the notation of Section~\ref{subsubsec:Steinberg CWE stack})
 a sequence of maps of $R^{\ver}_{\thetabar}$-modules
\[
R^{(p, 0), \crys}/\varpi \cong \bF\llbracket a_0, b_0, b_1\rrbracket /(a_0b_1) \onto \vhat^*\ihat'^*_{\rhobar}\Functor(\cInd_{KZ}^G \Sym^0) \into 
\bF\llbracket a_0, b_0, b_1\rrbracket /(a_0b_1)
\]
whose composite %
is a unit multiple of~$a_0$, by Remark~\ref{a_0 and T}.
Hence 
\[
\vhat^*\ihat'^*_{\rhobar}\Functor(\cInd_{KZ}^G \Sym^0) \cong R^{(p, 0), \crys}/(\varpi, b_1) 
\buildrel {\eqref{cor:explicit presentation of crystalline deformation ring II}} \over = R^{(1, 0), \crys}/\varpi.\qedhere
\]
\end{proof}

\subsubsection{\texorpdfstring{$\Functor$}{F} on finitely presented representations.}

If~$\sigma$ is a Serre weight, we continue to write $z : \cZ(\sigma) \to \cX/\varpi$ for the closed immersion defined in Definition~\ref{defn:special-fibre-crystalline-weight-sigma}.

\begin{prop}\label{concentration and support}
Let~$\sigma$ be a Serre weight.
Then  $\Functor(\cInd_{KZ}^G\sigma)$ is pure of degree zero, and is in fact a coherent sheaf, whose scheme-theoretic support  is~$\cZ(\sigma)$.
If furthermore $\sigma \ne \Sym^{p-2} \otimes \det^a$ for any~$a$, then $z^* \Functor(\cInd_{KZ}^G \sigma)$ is an invertible sheaf on~$\cZ(\sigma)$.
\end{prop}
\begin{proof}
In the discussion preceding the statement and proof of Proposition~\ref{computation
of completed stalk},
we noted that 
 for every closed point~$\rhobar$ of~$\cX$, the pullback $\ihat'^*_{\rhobar}\Functor(\cInd_{KZ}^G \sigma)$ is pure of degree zero. %
Since~$\Functor$ is right $t$-exact, it then follows from  Lemma~\ref{lem:pro-coh-vanishes-if-pullbacks-do}~(2)
that the same is true of $\Functor(\cInd_{KZ}^G \sigma)$, which is therefore a pro-coherent sheaf.

We will now show that  $\Functor(\cInd_{KZ}^G \sigma)$ is scheme-theoretically supported on~$\cZ(\sigma)$, in the sense of Definition~\ref{defn:scheme-theoretic-support-pro-coherent}~(1).
Recall the commutative diagram~\eqref{completed stalks diagram}.
By Lemma~\ref{lem:testing support with completions},
it suffices to prove that for every closed point $\rhobar$ of~$\cX$, 
the pro-coherent sheaf $\ihat'^*_{\rhobar}\Functor(\cInd_{KZ}^G\sigma)$ on~$(\cX/\varpi)^\wedge_{\rhobar}$
is scheme-theoretically supported on %
  the completion $\cZ(\sigma)^{\wedge}_{\rhobar} \to \cZ(\sigma)$ at $\{\rhobar\} \cap |\cZ(\sigma)|$.
We thus need to prove that the unit of the adjunction 
\[
(z_{\rhobar}^*, z_{\rhobar, *}) : \Pro\Coh\bigl((\cX/\varpi)^\wedge_{\rhobar}\bigr) \to
\Pro \Coh \bigl( \cZ(\sigma)^\wedge_{\rhobar}\bigr)
\]
of~\eqref{eqn: adjoint pair for closed immersion II} is an isomorphism at~$\ihat'^*_{\rhobar}\Functor(\cInd_{KZ}^G \sigma)$.
By Proposition~\ref{computation of completed stalk}, 
$\ihat'^*_{\rhobar} \Functor(\cInd_{KZ}^G \sigma)$ is an object of the full subcategory $\Coh(\cO_{(\cX/\varpi)^\wedge_{\rhobar}})$ 
of $\Pro \Coh\bigl((\cX/\varpi)^\wedge_{\rhobar}\bigr)$.
The completed pullback functor~$\widehat v^*$ is faithful on this subcategory, by Lemma~\ref{lem:versal pullback is exact and faithful II} (applied with $\widehat \cX = (\cX/\varpi)^\wedge_{\rhobar}$,
which is a completion of $(\fX_{\thetabar}/\varpi)$, and so falls within the scope of Lemma~\ref{lem:versal pullback is exact and faithful II}).
It is therefore enough to show
that the natural map
\[
\vhat^*\ihat'^*_{\rhobar} \Functor(\cInd_{KZ}^G \sigma) \to \bigl(\vhat^*\ihat'^*_{\rhobar} \Functor(\cInd_{KZ}^G \sigma)\bigr) \otimes_{R^{\ver}_{\thetabar}} \Rsigma_{\thetabar}
\]
is an isomorphism.
This is a consequence of Proposition~\ref{computation of completed stalk}.

At this point in the proof, we have shown that $\Functor(\cInd_{KZ}^G \sigma)$ is a pro-coherent sheaf, scheme-theoretically supported on ~$\cZ(\sigma)$.
In particular, the unit of adjunction
$\Functor(\cInd_{KZ}^G\sigma)\iso z_{*}z^{*}\Functor(\cInd_{KZ}^G\sigma)$
is an isomorphism, and 
 $\Functor(\cInd_{KZ}^G \sigma)$ is 
 isomorphic to
$L_{\infty} \otimes_{\cO\llbracket G\rrbracket _{\zeta}} \cInd_{KZ}^G \sigma$
(i.e.\ we can replace the derived tensor product in the definition of~$\Functor$  by a non-derived one).
We then find that %
\begin{multline*}
z^*\Functor(\cInd_{KZ}^G\sigma)  = z^*(L_{\infty}\otimes_{\cO\llbracket G\rrbracket _{\zeta}}\cInd_{KZ}^G \sigma) 
\buildrel \textrm{Lem.\ \ref{lem:abelian-composition-right-exact-tensor}} \over = z^{*}L_{\infty} \otimes_{\cO\llbracket G\rrbracket _{\zeta}} \cInd_{KZ}^G \sigma \\
  = z^{*}L_{\infty} \otimes_{\cO\llbracket K\rrbracket } \sigma 
  = (z^{*}L_{\infty})_{K_1} \otimes_{\cO\llbracket K/K_1\rrbracket } \sigma.
\end{multline*}
Since $(z^{*}L_{\infty})_{K_{1}}$ is coherent,
by Lemma~\ref{lem:explicit pro-coherent structure}~\eqref{item:53}, we see that $\Functor(\cInd_{KZ}^G\sigma)$ is indeed coherent.

In order to show that the scheme-theoretic support of $\Functor(\cInd_{KZ}^G\sigma)$  is~$\cZ(\sigma)$, we need to prove that the natural map
\begin{equation}\label{to prove injective for support}
\cO_{\cZ(\sigma)} \to \underline{\mathrm{End}}_{\cO_{\cZ(\sigma)}}\bigl(z^*\Functor(\cInd_{KZ}^G\sigma)\bigr)
\end{equation}
is injective.
Applying Lemma~\ref{lem:pro-coh-vanishes-if-pullbacks-do} to the kernel of~\eqref{to prove injective for support}, we see that it suffices to show that~\eqref{to prove injective for support}
becomes an injection after applying $\ihat'^*_{\rhobar,\sigma}$, for any closed point~$\rhobar$ of~$\cZ(\sigma)$.
By Lemma~\ref{lem: versal pullback and internal Hom}, applying $\vhat^*_{\sigma} \ihat'^*_{\rhobar,\sigma}$ to~\eqref{to prove injective for support} yields %
the natural map
\begin{multline}\label{eqn:action map on internal Hom}
R^\sigma_{\thetabar}/\varpi \to \underline{\mathrm{End}}_{\Mod_c(R^\sigma_{\thetabar})}(\vhat^*_{\sigma} \ihat'^*_{\rhobar,\sigma}z^*\Functor(\cInd_{KZ}^G \sigma))  \\ = 
\underline{\mathrm{End}}_{\Mod_c(R^\sigma_{\thetabar})}(\vhat^* \ihat'^*_{\rhobar}\Functor(\cInd_{KZ}^G \sigma));
\end{multline}
here the equality is obtained by an application of~\eqref{to prove commutative for compatibility of completion and pushforward} 
to the outer rectangle in the commutative diagram~\eqref{completed stalks diagram}, taking into account the fact that
$\vhat^* \ihat'^*_{\rhobar}\Functor(\cInd_{KZ}^G \sigma)$ is an $R^\sigma_{\thetabar}$-module (by Proposition~\ref{computation of completed stalk}).
The map~\eqref{eqn:action map on internal Hom} is injective, again by Proposition~\ref{computation of completed stalk}.
Since~$\ihat'^*_{\rhobar}$ pulls back the morphism~\eqref{to prove injective for support} to a morphism in 
$\Coh(\cO_{(\cX/\varpi)^\wedge_{\rhobar}})$ (by Lemma~\ref{lem:completed pullback preserves O_X modules}),
and $\widehat v^*$ is faithful on $\Coh(\cO_{(\cX/\varpi)^\wedge_{\rhobar}})$ (as we have already seen earlier in this proof),
this concludes the proof that~\eqref{to prove injective for support} is injective.

Finally, assume that $\sigma \ne \Sym^{p-2} \otimes \det^a$ for any~$a$. 
It remains to prove that $z^* \Functor(\cInd_{KZ}^G \sigma)$ is an invertible sheaf on~$\cZ(\sigma)$. 
Since~$\cZ(\sigma)$ is reduced, %
it suffices to show that its fibre dimension is equal to~$1$ at all finite type points of~$\cZ(\sigma)$.
Since~\eqref{to prove injective for support} is injective, the fibre dimension is positive at all finite type points of~$\cZ(\sigma)$, 
hence the loci of fibre dimension~$=1$ and~$\leq 1$ coincide.
By semicontinuity, this locus is an open subset of~$|\cZ(\sigma)|$, 
and by Proposition~\ref{computation of completed stalk}~(3), it contains all closed points.
It is therefore equal to $|\cZ(\sigma)|$.
\end{proof}

\begin{prop}\label{finiteness of cohomology}%
The functor  %
\[
\Functor : D^b_{\fp}(\cA)\to \Pro D^b_\coh(\cX)
\]
has essential image contained in $D^b_\coh(\cX)$. Furthermore it is of amplitude~$[-1,0]$.
\end{prop}%
\begin{proof}Since~$D^b_{\fp}(\cA)=D^b(\cA^{\fp})$, it follows from
  Lemma~\ref{lem:detecting-coh-in-pro-coh} that it suffices to show that for
  each~$\pi\in\cA^{\fp}$, 
  $\Functor(\pi)$ has  cohomology in~$\Coh(\cX)$, which furthermore vanishes
  outside of degrees $-1,0$.

We begin by showing the vanishing.
By Lemma~\ref{lem:pro-coh-vanishes-if-pullbacks-do}, it suffices to prove that $\ihat'^*_{\thetabar}H^{i}\Functor(\pi)=0$ for all closed points~$\thetabar$ of~$\cX$ and all $i \ne -1, 0$.
Since $\ihat'^*_{\thetabar}$ is exact, this is equivalent to $H^i(\ihat'^*_{\thetabar}\Functor(\pi)) = 0$.
Since~\eqref{unit of adjunction on functor} is an isomorphism (by Proposition~\ref{prop:completing our functor on Y})
we reduce to showing that $H^i(\ProFunctor_{\thetabar}(\ihat^*_{\thetabar}\pi))$ vanishes for~$i \ne -1, 0$.
This follows from Theorem~\ref{thm:derived-F-thetabar-nonSt-case}.

We now prove that $H^0\Functor(\pi)$ and $H^{-1}\Functor(\pi)$ are coherent, for all~$\pi \in \cA^{\fp}$.
Recall that~$\pi$ admits a presentation
\[\cInd_{KZ}^GW\to \cInd_{KZ}^GV\to \pi\to 0\]
with $V, W$ finite $\cO$-modules.
A d\'evissage using Proposition~\ref{concentration and support} shows  that
$H^{-1}\Functor(\cInd_{KZ}^G V)=0$, and $H^{0}\Functor(\cInd_{KZ}^GV)$ is coherent.
The right $t$-exactness of~$\Functor$ then implies that $H^0\Functor(\pi)$ is coherent for all~$\pi \in \cA^{\fp}$.
Now, since~$\cA^{\fp}$ is abelian, we have a short exact sequence \[0\to \pi'\to \cInd_{KZ}^GV\to \pi\to 0\] with~$\pi'\in\cA^{\fp}$.
We thus obtain an exact sequence \[0\to H^{-1}\Functor(\pi)\to H^0\Functor(\pi')\to H^0\Functor(\cInd_{KZ}^GV)\to H^0\Functor(\pi)\to 0, \]and we conclude that~$H^{-1}\Functor(\pi)$ is coherent, as required.
\end{proof}

\subsubsection{\texorpdfstring{$\Functor$}{F} preserves the recollement.}
We now show that for any (not necessarily finite) closed subset~$Y\subseteq X$, the functor~$\Functor_{Y}$ takes~$D^b_{\fp}(\cA_Y)$ to~$D^b_{\coh}(\cX_{Y})$.

\begin{prop}%
\label{prop:existence of FY FU overview}Let~$Y$ be a closed subset of~$X$ with
open complement~$j:U\into X$, and write %
$\cU\coloneqq \piss^{-1}(U)$ and $j'^*_\cU: D^b_{\coh}(\cX) \to D^b_{\coh}(\cU)$ for the pullback functor.
\begin{enumerate}
\item\label{item:17}The essential image of the functor ~$\Functor_{Y}:D^b_{\fp}(\cA_Y)\to
  \Pro D^b_{\coh}(\cX_{Y})$ is contained in~$ D^b_{\coh}(\cX_{Y})$, and there is a commutative
  diagram %
\[\begin{tikzcd}
	 D^b_{\fp}(\cA_Y) & D^b_{\fp}(\cA) \\
	 D^b_{\coh}(\cX_{Y}) & D^b_{\coh}(\cX)
	\arrow["i_{Y,*}", from=1-1, to=1-2]
	\arrow["\Functor_{Y}"', from=1-1, to=2-1]
		\arrow["\Functor", from=1-2, to=2-2]
		\arrow["i'_{Y,*}", from=2-1, to=2-2]
	\end{tikzcd}\]
\item\label{item:18} %
The composite $\Dfp^b(\cA) \buildrel \Functor \over \longrightarrow
D^b_{\coh}(\cX) \buildrel j'^{*}_{\cU} \over \longrightarrow
D^b_{\coh}(\cU)$ factors through a functor $$\Functor_U: \Dfp^b(\cA_U) \to D^b_{\coh}(\cU).$$ Furthermore, if~$\thetabar$ is a closed point of~$U$ then there is a natural isomorphism
\begin{equation}\label{eqn:completing our functor on Y II}
\ihat'^*_{\thetabar}\Functor_U \isoto (\Pro \Functor_{\thetabar})\ihat^*_{\thetabar} : D^b_{\fp}(\cA_U) \to \Pro D^b_{\coh}(\cX_{\thetabar}).
\end{equation}
\end{enumerate}
\end{prop}
\begin{proof}
We begin by noting that,
  for every~$\thetabar\not\in Y$, %
 it follows from
  Proposition~\ref{prop:completing our functor on Y} that
  \begin{equation}
    \label{eq:vanishing-of-FY-at-theta}
    \ihat'^*_{\thetabar}\Functor i_{Y,*}=(\ProFunctor_{\thetabar})\ihat^*_{\thetabar} i_{Y,*}=0
  \end{equation}(where we have used that $\ihat^*_{\thetabar}i_{Y,*}=0$).
It follows from~\eqref{eq:vanishing-of-FY-at-theta} and Lemma~\ref{lem:pro-coh-vanishes-if-pullbacks-do}
that $j'^*_{\cU} \Functor i_{Y, *} = 0$;  equivalently, bearing in mind Proposition~\ref{aprop:Coh-set-theoretically-supported}, we see that
the composite $\Functor i_{Y,*}: D^b_{\fp}(\cA_{Y})\to
D^b_{\coh}(\cX)$
has essential image contained in the full subcategory
$D^b_{\coh}(\cX_{Y})=D^b_{\coh, Y}(\cX)$ consisting of
complexes whose cohomology sheaves are set-theoretically supported on~$Y$.
In particular, we may write $\Functor i_{Y,*}=\ihat'_{Y,*}
\Functor'_{Y}$ for some functor $\Functor'_Y:  D^b_{\fp}(\cA_{Y})\to
D^b_{\coh}(\cX_{Y})$.

We therefore have the following commutative
diagram, which %
immediately gives~\eqref{item:17}.   %

\[\begin{tikzcd}
	D^b_{\fp}(\cA_Y) & D^b_{\fp}(\cA) & D^b_{\coh}(\cX) & \Pro D^b_{\coh}(\cX_{Y}) \\
	& D^b_{\coh}(\cX_{Y})
        \arrow["\Functor_Y",from=1-1, to=1-4,curve={height=-30pt}]
	\arrow["i_{Y,*}", from=1-1, to=1-2]
	\arrow["\Functor'_{Y}"', from=1-1, to=2-2]
	\arrow["\Functor", from=1-2, to=1-3]
	\arrow["\ihat'^{*}_{Y}", from=1-3, to=1-4]
	\arrow["i'_{Y,*}", from=2-2, to=1-3]
	\arrow[hook,from=2-2, to=1-4]
\end{tikzcd}\]
The existence of~$\Functor_U$ is a formal consequence of ~\eqref{eq:vanishing-of-FY-at-theta}, Corollary~\ref{cor:all-the-nice-properties-of-D(A)} and the  commutativity of the
  diagram. %
Finally, \eqref{eqn:completing our functor on Y II} follows from~ \eqref{unit of adjunction on functor}. 
\end{proof}%

\begin{defn}\label{defn:Ind-extended-functors}%
  By Proposition~\ref{prop:existence of FY FU overview} and
  Corollary~\ref{cor:all-the-nice-properties-of-D(A)}, the $\Ind$-extensions of
  the functors~$\Functor,\Functor_Y,\Functor_U$ give continuous, right $t$-exact
  functors, that we continue to denote by the same symbols %
  \[\Functor: D(\cA)\to \Ind D^b_\coh(\cX),\]
  \[\Functor_Y: D(\cA_Y)\to \Ind D^b_\coh(\cX_Y),\]
  \[\Functor_U: D(\cA_U)\to \Ind D^b_\coh(\cU).\]
\end{defn}

\subsection{Localization to \texorpdfstring{ $\cUgood$}{Ugood}}\label{subsection:localization-to-Ugood}%

Recall that in Definition~\ref{defn:f(t)} we defined a finite
subset~$\Ybad\subset |X|$ with dense open complement~$\Ugood$, and in Definition~\ref{defn:closed substack of bad points} we defined
$\cUgood$ as the open substack of~$\cX$ with underlying topological space $\pi_{\mathrm{ss}}^{-1}(\Ugood)$.
Our goal in this section is to establish some properties of our functor~$\Functor$ after restriction to~$\Ugood$ and~$\cUgood$; in particular, we will show that the functor~$\Functor_{\cUgood}$ is fully faithful (see Theorem~\ref{thm:full faithfulness on U}).

By Lemma~\ref{lem:reduced reducible
  connected components}, $\cUgood$ decomposes as the disjoint union of open
subsets $\cU(\sigmasigmacomp),$ and correspondingly $\Ugood$ decomposes as the
disjoint union of open subsets $U(\sigmasigmacomp)$.
Accordingly, we will be able to work with each open subset  $U(\sigmasigmacomp)$ separately, and a key tool in our analysis will be the explicit description of the stack  $\cU(\sigmasigmacomp)_{\red}$ in Proposition~\ref{prop: explicit description of underlying
      reduced of U}. We will need to study the coherent sheaf $\Functor(\cInd_{KZ}^G \sigma)$; 
if ~$\{\sigma, \sigmacomp\}$ is not of type~\ref{item: non p-distinguished block}, then this is an invertible sheaf on~$\cZ(\sigma)$ (see Proposition~\ref{concentration and support}), and 
our analysis could be simplified, but we give a uniform treatment of all~$\sigmasigmacomp$ below.

We now fix a companion pair~$\{\sigma,\sigmacomp\}$.  To
simplify notation, we make the following definitions.
\begin{defn}\label{defn:notation-for-localization-to-Ugood}\leavevmode
\begin{enumerate}
\item We set~$Y\coloneqq \Ybad\cap X(\sigmasigmacomp)$, $U\coloneqq U(\sigmasigmacomp) = X(\sigmasigmacomp)\setminus Y$,
  $\cU \coloneqq  \cU(\sigmasigmacomp)=\piss^{-1}(U)$. We also set $\cU(\sigma)\coloneqq \cU\cap\cZ(\sigma)$. 
  \item 
  We write
  $j^{*}:D(\cA)\to D(\cA_{U})$,
  $j'^{*}:\Ind D^b_{\coh}(\cX)\to \Ind D^b_{\coh}(\cU)$ for the $t$-exact localization
  functors, and $j_{*},j'_{*}$ for their $t$-exact fully faithful right
  adjoints (see Corollary~\ref{cor:all-the-nice-properties-of-D(A)} and Lemma~\ref{lem:vanishing-higher-derived-functors-j}). 
  \item 
  We write~$f=f(t)$ for the
polynomial defined in Definition~\ref{defn:f(t)}, so that $f_\sigma^{-1}\bigl(U(\sigmasigmacomp)\bigr) \subset \Spec \cH_G(\sigma)$ is the non-vanishing locus of 
$f(\HeckeT) \in \cH_G(\sigma) = \F[\HeckeT]$.
Note that this implies
\[
j_*j^* \cInd_{KZ}^G \sigma \isoto \colim_{\times f}\cInd_{KZ}^G\sigma,
\]
by Lemma~\ref{lem:localization-of-cInd-sigma}.
\end{enumerate}
\end{defn}%

We begin by specializing Corollary~\ref{cor:explicit-generators-D(A)} to the situation at hand.

\begin{lem}
  \label{lem:generators-of-quotient-category}The localizations
  $j^*\cInd_{KZ}^G\sigma$, $j^*\cInd_{KZ}^G\sigmacomp$ form a set of compact generators of
  $D(\cA_U)$.  %
\end{lem}
\begin{proof}
  By Corollary~\ref{cor:explicit-generators-D(A)}, the set of
  localizations $\{j^*\cInd_{KZ}^G\sigma'\}$, where~$\sigma'$ ranges over the $\zeta$-compatible Serre weights, forms a set of generators of
  $D(\cA_{U})$.  If $\sigma'\ne\sigma,\sigmacomp$, and ~$\sigma'$ is
  not a twist of~$\Sym^{p-1}$, then $j^{*}\cInd_{KZ}^G\sigma'=0$ by~\cite[Prop.\
  3.1.13]{DEGlocalization} (i.e.\ essentially by the very definition
  of~$j^{*}$). It remains to note that by~\eqref{Steinberg change of weight} and the
  definition of~$\Ugood$, we have (for any~$r$) \[j^{*}\cInd_{KZ}^G\det{}^r\otimes\Sym^{p-1}=j^{*}\cInd_{KZ}^G\det{}^r\otimes\Sym^0.\qedhere\]
  \end{proof}

In preparation for the next result, we recall
that $\Functor(\cInd_{KZ}^G\sigma/f\cInd_{KZ}^G\sigma)$
is pure of degree zero. 
(Indeed, 
$\cInd_{KZ}^G\sigma/f\cInd_{KZ}^G \sigma$ is an object of~$\cA_Y$,
and
Proposition~\ref{prop:existence of FY FU overview}~\eqref{item:17} identifies
$\Functor(\cInd_{KZ}^G\sigma/f\cInd_{KZ}^G\sigma)$
with~$\Functor_Y(\cInd_{KZ}^G\sigma/f\cInd_{KZ}^G\sigma),$
which Lemma~\ref{concentration} shows is pure of degree zero; note that~$\sigma$ is not Steinberg, by the definition of a companion pair.)
The following lemma describes the underlying topological space of the scheme-theoretic support of this coherent sheaf, and shows that its endomorphism~$\Functor(f)$ is nilpotent.
In its proof, we will make use of the fact that $\Functor(\cInd_{KZ}^G\sigma)$ is scheme-theoretically supported on $\cZ(\sigma)$, by Proposition~\ref{concentration and support},
hence the set-theoretic support of its quotient $\Functor(\cInd_{KZ}^G\sigma/f \cInd_{KZ}^G \sigma)$ (in the sense of Definition~\ref{defn:scheme-theoretic-support-pro-coherent})
is \emph{a priori} contained in $|\cZ(\sigma)|$.

\begin{lem}\label{support of quotients of functor on weights}\leavevmode
  Let $i_\cY' : \cY \to \cX$
be the scheme-theoretic support of $\Functor(\cInd_{KZ}^G\sigma/f\cInd_{KZ}^G\sigma)$. %
  Then
  \begin{enumerate}
  \item 
      $|\cY| = \coprod_{\thetabar \in \Ybad} |\cZ(\sigma)| \cap |\cX_{\thetabar}|= |\cZ(\sigma)|\setminus |\cU(\sigma)|.$ 
  \item The endomorphism $\Functor(f)$ of the coherent sheaf $i'^{*}_{\cY}\Functor(\cInd_{KZ}^G\sigma)$ is nilpotent.
  \end{enumerate}
\end{lem}
\begin{proof}
By Definition~\ref{defn:f(t)}, we have $f = \prod_{\thetabar \in \Ybad}f_{\thetabar}$, where    $f_{\thetabar} \in \cH(\sigma)$ is as in Definition~\ref{f_Y}, so that $\cInd_{KZ}^G \sigma/f_{\thetabar}\cInd_{KZ}^G \sigma$ is the maximal multiplicity-free quotient of $\cInd_{KZ}^G \sigma$ which is an object of~$\cA_{\thetabar}$,
and
\[
\cInd_{KZ}^G \sigma/ f \cInd_{KZ}^G \sigma = \prod_{\thetabar \in \Ybad}
\cInd_{KZ}^G \sigma/f_{\thetabar} \cInd_{KZ}^G \sigma.
\]
Thus, if we let $i'_{\cY_{\thetabar}}: \cY_{\thetabar} \hookrightarrow \cX$
denote the scheme-theoretic support of the coherent sheaf
$\Functor(\cInd_{KZ}^G \sigma/f_{\thetabar}\cInd_{KZ}^G \sigma)$,
we find that $\cY = \coprod_{\thetabar \in \Ybad} \cY_{\thetabar}.$
It therefore suffices to prove, for each~$\thetabar \in \Ybad$,
that $|\cY_{\thetabar}| = |\cZ(\sigma)| \cap |\cX_{\thetabar}|$, 
and that $\Functor(f_{\thetabar})$ acts nilpotently on $i'^{*}_{\cY_{\thetabar}}\Functor(\cInd_{KZ}^G\sigma)$.

By Proposition~\ref{prop:existence of FY FU overview} we have
\begin{multline*} \Functor(\cInd_{KZ}^G \sigma/f_{\thetabar}\cInd_{KZ}^G \sigma)
=i'_{\thetabar,*}\Functor_{\thetabar}(\cInd_{KZ}^G \sigma/f_{\thetabar}\cInd_{KZ}^G \sigma)\\
=i'_{\thetabar,*}\ihat'^{*}_{\thetabar}\Functor(\cInd_{KZ}^G \sigma/f_{\thetabar}\cInd_{KZ}^G \sigma),
\end{multline*}
so we see that
$|\cY_{\thetabar}|$
(which is \emph{a priori} contained in~$|\cZ(\sigma)|$, as recalled before the statement of the lemma)
is contained in $|\cZ(\sigma)| \cap |\cX_{\thetabar}|$.

On the other hand, by Proposition~\ref{prop:completing our functor on Y} and Lemma~\ref{lem:localization-of-cInd-sigma} we have the identifications
\begin{align*}
\ihat'^{*}_{\thetabar} \Functor(\cInd_{KZ}^G \sigma)
&= \Pro\Functor_{\thetabar}\bigl(\quoteslim{n} (\cInd_{KZ}^G \sigma)/f_{\thetabar}^n \cInd_{KZ}^G\sigma\bigr) \\
&= \quoteslim{n}
\Functor_{\thetabar}\bigl((\cInd_{KZ}^G \sigma)/f_{\thetabar}^n( \cInd_{KZ}^G\sigma)\bigr).
\end{align*}%
Applying~$\ihat'_{\thetabar,*}$, %
we find that
\begin{equation}\label{eqn:notation for source of surjection}
\ihat'_{\thetabar,*}\ihat'^{*}_{\thetabar} \Functor(\cInd_{KZ}^G \sigma)
= \quoteslim{n}
i'_{\thetabar,*}\Functor_{\thetabar}\bigl((\cInd_{KZ}^G \sigma)/f_{\thetabar}^n( \cInd_{KZ}^G\sigma)\bigr).
\end{equation}
Since $i'_{\thetabar,*}\Functor_{\thetabar}\bigl((\cInd_{KZ}^G \sigma)/f_{\thetabar}^n( \cInd_{KZ}^G\sigma)\bigr)$ is filtered by copies of the coherent sheaf
$i'_{\thetabar,*}\Functor_{\thetabar}\bigl((\cInd_{KZ}^G \sigma)/f_{\thetabar}( \cInd_{KZ}^G\sigma)\bigr)$, 
we see that the set-theoretic support of every coherent quotient of 
$\ihat'_{\thetabar,*}\ihat'^{*}_{\thetabar} \Functor(\cInd_{KZ}^G \sigma)$
is contained
in that of 
\[
i'_{\thetabar,*}\Functor_{\thetabar}(\cInd_{KZ}^G \sigma/f_{\thetabar}\cInd_{KZ}^G \sigma) = \Functor(\cInd_{KZ}^G \sigma/f_{\thetabar} \cInd_{KZ}^G \sigma),
\]
which is~$|\cY_{\thetabar}|$.
By Proposition~\ref{concentration and support}, the scheme-theoretic support of $\Functor(\cInd_{KZ}^G \sigma)$ is~$\cZ(\sigma)$, so in particular 
$\Functor(\cInd_{KZ}^G \sigma)$ has a quotient~$\cG$ whose
set-theoretic support is exactly equal to~$|\cZ(\sigma) |\cap |\cX_{\thetabar}|$.
The unit $\cG \to \ihat'_{\thetabar, *}\ihat'^*_{\thetabar} \cG$ is then necessarily an isomorphism, and so~$\cG$
is a coherent quotient of
$\ihat'_{\thetabar,*}\ihat'^{*}_{\thetabar} \Functor(\cInd_{KZ}^G \sigma)$.
It follows from the previous paragraph that~$\cG$ has set-theoretic support contained 
in that of $\Functor(\cInd_{KZ}^G \sigma/f_{\thetabar}\cInd_{KZ}^G \sigma)$, which is~$|\cY_{\thetabar}|$. 
Thus 
$|\cY_{\thetabar}|$ is
equal to $|\cZ(\sigma)| \cap |\cX_{\thetabar}|$, as claimed.
This concludes the proof of the first part.

Turning to the second part, write $i'_{\cX_{\thetabar, \red}}: \cX_{\thetabar, \red} \hookrightarrow \cX$ for the closed immersion.
Since~$\cX_{\thetabar}$ is the completion of~$\cX$ at~$\cX_{\thetabar, \red}$,
we see from 
Remark~\ref{rem:ihat surjects onto i upper star} 
that there is a natural surjection
\[
\ihat'_{\thetabar,*}\ihat'^{*}_{\thetabar} \Functor(\cInd_{KZ}^G \sigma) \onto
i'_{\cX_{\thetabar, \red}, *}i'^{*}_{\cX_{\thetabar,\red}}\Functor(\cInd_{KZ}^G \sigma)
\]
in $\Pro \Coh (\cX).$
We noted in~\eqref{eqn:notation for source of surjection} that the source of this surjection is equal to 
\[
\quoteslim{n} 
i'_{\thetabar,*}\Functor_{\thetabar}\bigl((\cInd_{KZ}^G \sigma)/f_{\thetabar}^n( \cInd_{KZ}^G\sigma)\bigr), 
\]
so we see that the coherent sheaf $i'_{\cX_{\thetabar, \red}, *}i'^{*}_{\cX_{\thetabar,\red}}\Functor(\cInd_{KZ}^G \sigma)$
is a quotient of  
\[
i'_{\thetabar,*}\Functor_{\thetabar}\bigl((\cInd_{KZ}^G \sigma)/f_{\thetabar}^n( \cInd_{KZ}^G\sigma)\bigr) = 
\Functor\bigl((\cInd_{KZ}^G \sigma)/f_{\thetabar}^n( \cInd_{KZ}^G\sigma)\bigr)
\]
for~$n$ sufficiently large, whence it is killed by $\Functor(f_{\thetabar})^{n}=\Functor(f_{\thetabar}^{n})$. 

We have already seen that $|\cY_{\thetabar}| \subseteq |\cX_{\thetabar}|,$
and so $\cY_{\thetabar,\red} \subseteq \cX_{\thetabar,\red}.$
Thus, %
letting $i'_{\cY_{\thetabar,\red}}$ 
 denote the closed immersion
$\cY_{\thetabar,\red} \hookrightarrow \cX,$ 
the result of the previous paragraph shows that
$\Functor(f_{\thetabar})$ acts nilpotently on
$i'^*_{\cY_{\thetabar,\red}}\Functor(\cInd_{KZ}^G \sigma).$
This immediately implies that
$\Functor(f_{\thetabar})$ also acts nilpotently on
$i'^*_{\cY_{\thetabar}}\Functor(\cInd_{KZ}^G \sigma),$
as required.
\end{proof}

We now study the interaction of our functor~$\Functor$ with the localization functors in Definition~\ref{defn:notation-for-localization-to-Ugood}.
By construction, we have an isomorphism $\Functor_U j^* \isoto j'^* \Functor$.
Postcomposing the inverse of this isomorphism
with $j'_*$, the unit of adjunction $1 \to j'_*j'^*$ induces a morphism
$\Functor \to j'_*\Functor_U j^*$.
Precomposing with~$j_*$, %
the counit $j^*j_* \to 1$ induces a natural transformation
\begin{equation}\label{eqn:Fj*-j*F}\Functor j_* \to j'_*\Functor_{U}.\end{equation}

\begin{prop}
\label{prop:F commutes with j lower star}
The natural transformation~\eqref{eqn:Fj*-j*F} is a natural isomorphism
of functors from $\Dfp^b(\cA_{U})$ to $\Ind D^b_{\coh}\bigl(\cX)$.
\end{prop}%
\begin{proof}By the continuity of all functors involved, it is equivalent to show that~\eqref{eqn:Fj*-j*F} is a natural isomorphism
of the corresponding Ind-extended functors  $D(\cA_{U})\to \Ind D^b_{\coh}\bigl(\cX)$.
By Lemma~\ref{lem:generators-of-quotient-category},  $j^*\cInd_{KZ}^G\sigma$, $j^*\cInd_{KZ}^G\sigmacomp$ are compact generators of  $D(\cA_{U})$, so by
the symmetry between~$\sigma$ and~$\sigmacomp$, it suffices to show that the
morphism induced by ~\eqref{eqn:Fj*-j*F} 
\begin{equation}\label{to prove isomorphism for recollement}
\Functor (j_*j^*\cInd_{KZ}^G \sigma) \to
  j'_*\Functor_{U}(j^*\cInd_{KZ}^G \sigma) %
\end{equation}
  is an isomorphism.

Recall (as noted in Definition~\ref{defn:notation-for-localization-to-Ugood}) that the functors~$j_{*},j^{*}$ and $j'_{*},j'^{*}$ are $t$-exact, and $\Functor(\cInd_{KZ}^G \sigma)$ is a sheaf (i.e.\ is in the heart of $\Ind D^b_{\coh}\bigl(\cX)$), so the same is true of $j'_*\Functor_{U}(j^*\cInd_{KZ}^G \sigma) = j'_*j'^*\Functor(\cInd_{KZ}^G \sigma)$.
Similarly (as also noted in Definition~\ref{defn:notation-for-localization-to-Ugood}) %
we have $j_*j^*\cInd_{KZ}^G\sigma = (\cInd_{KZ}^G \sigma)[1/f]=\colim_{\times f} \cInd_{KZ}^G\sigma$,
and since $\Functor$ commutes with colimits, we have
\begin{equation}\label{eqn:functor-applied-to-colim}\Functor(j_*j^*\cInd_{KZ}^G\sigma) = \Functor(\cInd_{KZ}^G\sigma)[1/\Functor(f)]\coloneqq \colim_{\times \Functor(f)}\Functor(\cInd_{KZ}^G\sigma),\end{equation}which is also a sheaf.
Thus all the objects under consideration are in the hearts of the respective $t$-structures, and we work at the abelian level in what follows.
Bearing in mind~\eqref{eqn:functor-applied-to-colim}, we see that we need to prove that the natural map %
\[\Functor(\cInd_{KZ}^G\sigma)[1/\Functor(f)]\to j'_*j'^*\Functor(\cInd_{KZ}^G \sigma) \]is an isomorphism.

Since $\Functor(\cInd_{KZ}^G\sigma)$ is scheme-theoretically supported on the closed substack~$\cZ(\sigma)$ of~$\cX$, it is easy to check that we can replace~$\cX$ by~$\cZ(\sigma)$, and $j'$ by its restriction to the open immersion $\cU(\sigma)\into\cZ(\sigma)$. We claim that the result now follows from 
Lemma~\ref{lem:localization identification} below, which we apply with  $\cZ=\cZ(\sigma)$, $\cF=\Functor(\cInd_{KZ}^G\sigma)$, and~$f$ equal to our~$\Functor(f)$.
Indeed, it is only necessary to remark that in this case the coherent sheaf $\cF/f\cF$ is (by the right $t$-exactness of~$\Functor$) equal to~$\Functor(\cInd_{KZ}^G\sigma/f\cInd_{KZ}^{G}\sigma)$, so the hypotheses of Lemma~\ref{lem:localization identification} are supplied by Lemma~\ref{support of quotients of functor on weights}.
\end{proof}

\begin{lemma} %
  \label{lem:localization identification}
Let $\cF$ be a coherent sheaf
on a Noetherian algebraic stack~$\cZ$,
and let $f$ be an endomorphism of~$\cF$.
Let~$i : \cY \to \cZ$ be the scheme-theoretic support of~$\cF/f\cF$ (a closed algebraic substack of~$\cZ$)
and let $j: \cU %
\hookrightarrow \cZ$ %
be the open immersion of the complement of~$\cY$.  
Suppose further that $f$ acts nilpotently on $i^*\cF$.
Then $f$ becomes invertible on $j^*\cF$, and the natural map $\cF[1/f] \iso j_*j^*\cF$ is an isomorphism in $\Ind \Coh \cZ$. 
\end{lemma}
\begin{proof}
There is an exact sequence
\[
0 \to \cF[f] \to \cF \xrightarrow{f} \cF \to \cF/f\cF \to 0.
\]
Note that the morphism $\cF[f] \to \cF/f^N$ is injective for $N$ sufficiently large,
so that there is an inclusion of set-theoretic supports
$$\Supp(\cF[f]) \subseteq \Supp(\cF/f^N) = \Supp(\cF/f\cF) = |\cY|.$$ 
Thus multiplication by $f$ on $\cF$ becomes an isomorphism after applying~$j^*$, %
whence the unit of adjunction $\cF \to j_*j^*\cF$ induces a morphism $\cF[1/f] \to j_*j^*\cF$ in $\Ind \Coh(\cZ)$, which we will show is an isomorphism.
In order to do this we can replace~$\cZ$ by a smooth cover by a
scheme~$Z$, and~$\cU$ and~$\cF$ by their pullbacks to~$Z$, which we denote
by~$U$ and~$\cF$.

Let $\cA$ be the coherent algebra of endomorphisms of $\cF$ generated by $\cO_Z$ together
with~$f$,
and write $W = \underline{\Spec} \, \cA$, equipped with its canonical finite affine
morphism $\pi: W \to Z.$
As usual, $\cF$ gives rise to a coherent sheaf
$\cG \coloneqq  \cO_{W}\otimes_{\pi^{-1}\cA}\pi^{-1}\cF$ 
on $W$ such that $\cF \iso \pi_*\cG.$ 
Since $\cA$ acts faithfully on~$\cF$, we see that the support of~$\cG$ is
all of $W$.  Recall that this implies that for any coherent ideal sheaf
$\cI\subseteq \cO_W$, we have that $\Supp(\cG/\cI\cG) = V(\cI).$
In particular, $\Supp(\cG/f\cG) = V(f)$.

Now $\pi_*(\cG/f\cG) = \cF/f\cF,$ and so $Y \coloneqq  \Supp(\cF/f\cF)  = \pi\bigl(V(f)\bigr).$
We claim that the inclusion $V(f) \subseteq \pi^{-1}(Y)$ is an equality.
Given this, we find that
$\pi^{-1}(U)$ is equal to the distinguished open subset $D(f)$ of~$W$, 
and so, if $j':D(f) \hookrightarrow W$ is the corresponding open immersion,
we find (recalling that~$\pi$ is affine) that indeed
$$j_*j^{*}\cF =  \pi_* j'_*j'^{*}\cG = \pi_* \cG[1/f] = \cF[1/f].$$ 

To prove the claim,
write $\cI' \coloneqq  (\pi^{-1}\cI_{Y}) \cO_W$, so that $\cI'$ cuts out a closed 
subscheme of $W$ with underlying reduced equal to~$\pi^{-1}(Y)$.
Then $i^*\cF = \pi_*(\cG/\cI'\cG)$,
and our assumption that   $f$ acts nilpotently on $i^*\cF$ implies that
$f$ acts nilpotently on $\cG/\cI'\cG$. Hence (again using 
the fact that $\cG$ has full support on~$W$) we see that
$$\pi^{-1}(Y) = V(\cI') = \Supp(\cG/\cI'\cG) \subseteq V(f).$$
We've already noted the reverse inclusion, and the desired equality follows. 
\end{proof}

We now make use of the explicit description of~$\cU_{\red}$ proved in Section~\ref{subsec: defining U}.
Write~$B\coloneqq  \F[t,f(t)^{-1},x,y]/(xy)$,  so that by Proposition~\ref{prop: explicit description of underlying
      reduced of U} we may write $\cU_{\red}=[\Spec B / \bG],$ where the reductive group~$\bG$ and its action on~$B$ is specified in the statement of Proposition~\ref{prop: explicit description of underlying
      reduced of U}.
From this explicit description we see that  %
 $B^{\bG}=\F[t,f(t)^{-1}]$ unless~$\sigmasigmacomp$ is of type~\ref{item: non p-distinguished block}, in which case $B^{\bG}=\F[s,(s^2-4)^{\pm 1}]$ where $s=t+t^{-1}$. %
In either case there is a finite morphism of $\F$-algebras
\begin{equation}
\label{eqn:invariants mapping to Hecke-type ring}
B^{\bG}\to \F[t,f(t)^{-1}],
\end{equation}
and the maximal ideals of~$B^{\bG}$ are in bijection with the closed points~$\thetabar$ of~$U$. More precisely, if  
$\sigma=\sigma_{a,b}$, and~$\thetabar$ corresponds to the  semisimple
Galois representation %
\[\rhobar_{\thetabar} \coloneqq \nr_{\alpha} \zetabar\omega^{-a} \oplus
\nr_{\alpha^{-1}} \omega^{a-1},\] then %
the corresponding maximal ideal of~$B^{\bG}$ is 
$\mathfrak n_{\thetabar}\coloneqq (t-\alpha) \subseteq B^{\bG}$ (respectively $\mathfrak n_{\thetabar}\coloneqq (s-(\alpha+\alpha^{-1})) \subseteq B^{\bG}$ if 
$\sigmasigmacomp$ is of type~\ref{item: non p-distinguished block}); and we have
\begin{equation}
\label{eqn:ideal identification}
\mathfrak{n}_{\thetabar}\F[t,f(t)^{-1}]=(f_{\thetabar}(t)).
\end{equation}

  We may now apply the results of Appendix~\ref{subsubsec:reductive} to the stack $\cU_{\red}=[\Spec B/\bG]$, taking~$G$ there to be our~$\bG$, and~$R$ to be~$\F$.

\begin{prop}\label{prop:computation-End-U}%
The functor
$\Functor_{U}$ induces an isomorphism
$$\End_{\cA^{\fp}_{U}}\bigl(j^*\cInd_{KZ}^G \sigma\bigr)
\iso
\End_{\Coh(\cU)}\bigl(\Functor_U(j^*\cInd_{KZ}^G \sigma)\bigr).$$%
\end{prop}
\begin{proof}
By Proposition~\ref{concentration and support},
$\Functor_U(j^*\cInd_{KZ}^G \sigma)=j'^*\Functor(\cInd_{KZ}^G \sigma)$ is a coherent
sheaf on~$\cU$. %
Write
\[
S: = \End_{\Coh(\cU)}\bigl( \Functor_U(j^*\cInd_{KZ}^G \sigma)\bigr),
\]
so that $S$ is a (not necessarily commutative) $B^{\bG}$-algebra which,
by Corollary~\ref{cor:coherent exts for affine mod reductive-I},
is of finite type as a $B^{\bG}$-module.
By~\eqref{eqn:recollection-of-description-of-j*-cInd-sigma}, we have
\begin{equation}\label{eqn:description-of-j*-cInd-sigma}(\cInd_{KZ}^G \sigma )[1/f] \iso j_*j^*\cInd_{KZ}^G \sigma.\end{equation}
Consequently, bearing in mind the full faithfulness of~$j_*$, we find that
\begin{equation}\label{eqn:end-of-j*sigma-explicit}\End_{\cA_{U}^{\fp}}\bigl(j^*\cInd_{KZ}^G \sigma\bigr)
= \F[\HeckeT, f(\HeckeT)^{-1}].\end{equation}
The morphism
$$\Functor_U : \End_{\cA_{U}^{\fp}}\bigl(j^*\cInd_{KZ}^G \sigma\bigr)
\to \End_{\Coh(\cU)}\bigl( \Functor_U(j^*\cInd_{KZ}^G \sigma) \bigr) %
$$
induced by $\Functor_{U}$
is thus a morphism
\begin{equation}
\label{eqn:morphism of endos}
\F[\HeckeT, f(\HeckeT)^{-1}] \to S.
\end{equation}
Our goal is to show that this morphism is an isomorphism.

As above, we choose a closed point of~$U$, corresponding to a pseudorepresentation~$\thetabar$ and a maximal ideal $\mathfrak n_{\thetabar}$ of~$B^{\bG}$.
Recall from Theorem~\ref{thm:derived-F-thetabar-nonSt-case} that
\[
\Pro\Functor_{\thetabar} : \Pro D^b_{\fp}(\cA_{\thetabar}) \to \Pro D^b_{\coh}(\cX_{\thetabar})
\]
is fully faithful, and recall from Lemma~\ref{concentration} that %
$\ProFunctor_{\thetabar}(\ihat^*_{\thetabar} \cInd_{KZ}^G \sigma) = \ihat'^*_{\thetabar}\Functor(\cInd_{KZ}^G \sigma) = \ihat'^*_{\thetabar}\Functor_U(j^{*}\cInd_{KZ}^G \sigma)$ is pure of degree zero.
We thus obtain an isomorphism
\begin{equation}
\label{eqn:completion iso}
\ProFunctor_{\thetabar} : \End_{\Pro \cA_{\thetabar}^{\fp}}( \ihat_{\thetabar}^* \cInd_{KZ}^G \sigma)  
\iso \End_{\Pro \Coh \cX_{\thetabar}}\bigl( \ihat'^*_{\thetabar}\Functor(\cInd_{KZ}^G \sigma)\bigr).
\end{equation}
By Corollary~\ref{cor:completing RHoms},
the morphism
$$\ihat'^*_{\thetabar}: S: = \End_{\Coh\cU}\bigl( \Functor_U(j^*\cInd_{KZ}^G \sigma)\bigr)
\to \End_{\Pro \Coh\cX_{\thetabar}}\bigl( \ihat'^*_{\thetabar}\Functor(\cInd_{KZ}^G \sigma)\bigr)$$
is a completion at~$\fn_{\thetabar}$, i.e.\ it induces an isomorphism
$$\widehat{S}_{\mathfrak{n}_{\thetabar}} \iso
\End_{\Pro \Coh\cX_{\thetabar}}\bigl( \ihat'^*_{\thetabar}\Functor(\cInd_{KZ}^G \sigma)\bigr).$$ 
Moreover, by Proposition~\ref{Ext and completion}, 
  the map~\eqref{eqn: Ext and completion} 
induces the first of the identifications 
\begin{equation}
\label{eqn:completion identification}
\End_{\Pro \cA_{\thetabar}^{\fp}}(\ihat^*_{\thetabar}\cInd_{KZ}^G\sigma) =\bF[\HeckeT, f(\HeckeT)^{-1}]^{\wedge}_{f_{\thetabar}(\HeckeT)}
= \bF[\HeckeT, f(\HeckeT)^{-1}]^{\wedge}_{\fn_{\thetabar}},
\end{equation}
the second being induced by~\eqref{eqn:ideal identification}.
Finally, the natural isomorphism $\ihat'^*_{\thetabar} \Functor_U \isoto (\ProFunctor_{\thetabar}) \ihat^*_{\thetabar}$ 
of Proposition~\ref{prop:existence of FY FU overview} %
produces a commutative diagram
\begin{equation}\label{compatibility of F_U and Pro F_thetabar}
\begin{tikzcd}
\bF[\HeckeT, f(\HeckeT)^{-1}] = \End_{\cA_U^{\fp}}(j^* \cInd_{KZ}^G \sigma)
\arrow[r, "\Functor_U"', "\eqref{eqn:morphism of endos}"] 
\arrow[d, "\ihat^*_{\thetabar}"] & \End_{\Coh \cU}(\Functor_U(j^*\cInd_{KZ}^G \sigma)) \eqcolon S \arrow[d, "\ihat'^*_{\thetabar}"]\\
\bF[\HeckeT, f(\HeckeT)^{-1}]^{\wedge}_{\fn_{\thetabar}} \buildrel \eqref{eqn:completion
identification} \over =
\End_{\Pro \cA_{\thetabar}^{\fp}}(\ihat^*_{\thetabar}\cInd_{KZ}^G\sigma) \arrow[r, "\sim", "\ProFunctor_{\thetabar}"'] & \End_{\Pro \Coh \cX_{\thetabar}}(\ihat'^*_{\thetabar}\Functor(\cInd_{KZ}^G\sigma))=\widehat{S}_{\mathfrak{n}_{\thetabar}}
\end{tikzcd}
\end{equation}
whose top horizontal arrow is the morphism~\eqref{eqn:morphism of endos} that is the focus
of our attention.

We now make some deductions.  %
Firstly, we see that the completion of $S$ (which we recall is of finite type as a $B^{\bG}$-module)
at each maximal ideal of $B^{\bG}$ (which is itself a finite type $\F$-algebra)
is commutative and reduced.  Thus $S$ itself is commutative and
reduced, and  is in particular a finite commutative $B^{\bG}$-algebra.     
We have already observed that~\eqref{eqn:invariants mapping to Hecke-type ring}
makes $\F[\HeckeT,f(\HeckeT)^{-1}]$ a finite $B^\bG$-algebra.
Somewhat confusingly, we have not shown %
that~\eqref{eqn:morphism of endos}
is a morphism of~$B^\bG$-algebras. 
However, we see from ~\eqref{compatibility of F_U and Pro F_thetabar}
that~\eqref{eqn:morphism of endos} induces isomorphisms of the completions of its source and target
at each maximal ideal $\fn_{\thetabar}$ of~$B^\bG$. 
Now the finiteness
of each of $\F[\HeckeT, f(\HeckeT)^{-1}]$ and $S$ over~$B^\bG$ shows that 
the maximal ideals of each of them are partitioned into finite sets
according to the maximal~$\fn_{\thetabar}$ of $B^\bG$ that they lie over, 
and then the isomorphisms given by the bottom arrow 
of~\eqref{compatibility of F_U and Pro F_thetabar}
show that~\eqref{eqn:morphism of endos}
respects these partitions of the maximal ideals in its source
and target, and furthermore induces bijections on each of the finite sets
of maximal ideals given by these partitions. 
Thus~\eqref{eqn:morphism of endos}
in fact induces a bijection between the maximal ideals in its source and target,
and an isomorphism between the corresponding completions of each of its source 
and target.
It follows from Lemma~\ref{lem:isomorphism lemma} below that \eqref{eqn:morphism of endos} is an isomorphism, as desired.
\end{proof}

\begin{lemma}\label{lem:isomorphism lemma}
Let~$k$ be a field, and let $f: A \to B$ be a morphism between finite type $k$-algebras.
Assume that $f^{-1}$ %
induces a bijection on maximal ideals, and that for every maximal ideal $\fm \subset B$, the induced map 
$\widehat f : \widehat A_{f^{-1}(\fm)} \to \widehat B_{\fm}$ is an isomorphism.
Then~$f$ is an isomorphism.
\end{lemma}
\begin{proof}
Since $f^{-1}$ %
is bijective on maximal ideals, to prove that~$f$ an isomorphism, it suffices to prove that~$\Spec f$ is an open immersion. 
By~\cite[\href{https://stacks.math.columbia.edu/tag/02LC}{Tag~02LC}]{stacks-project}, it in turn suffices  to prove that~$f$ is \'etale and universally injective.

For each maximal ideal~$\m$ of~$B$, it follows from our hypotheses together with~\cite[\href{https://stacks.math.columbia.edu/tag/039M}{Tag~039M}]{stacks-project} that the map $A_{f^{-1}(\fm)} \to B_{\fm}$ is an \'etale homomorphism of local rings, and is in particular flat and unramified.
Hence $f: A \to B$ is flat, by~\cite[\href{https://stacks.math.columbia.edu/tag/00HT}{Tag~00HT}]{stacks-project}
and furthermore it is unramified at all maximal ideals of~$B$ 
(i.e.\ $(\Omega^1_{B/A})_\fm = 0$ for all maximal $\fm \subset B$), by~\cite[\href{https://stacks.math.columbia.edu/tag/039G}{Tag~039G}]{stacks-project}.
Thus $f: A \to B$ is \'etale.

It remains to show  that~$f$ is universally injective.
By~\cite[\href{https://stacks.math.columbia.edu/tag/01S4}{Tag~01S4}]{stacks-project}, it suffices to prove that $B \otimes_A B \to B$ induces a surjection on~$\Spec$.
Since the image of $\Spec(B) \to \Spec(B \otimes_A B)$ is closed, it suffices to show that for each finite extension $k'/k$, any morphism $\phi : B \otimes_A B \to k'$ necessarily factors through~$B$. %
Given such a map~$\phi$, we have two maps $\phi_1, \phi_2: B \to k'$ (defined by $\phi_1(b) = \phi(b \otimes 1)$ and $\phi_2(b) = \phi(1 \otimes b)$) that agree on~$A$, and we must show $\phi_1 = \phi_2$.
Since~$\phi_1 ,\phi_2 $ agree on~$A$, and~$f^{-1}$ is bijective on maximal ideals by hypothesis, we deduce that $\ker(\phi_1) = \ker(\phi_2) = \fm$ for some maximal ideal $\fm \subset B$.
The map $\widehat A_{f^{-1}(\fm)} \to \widehat B_{\fm}$ is an isomorphism by hypothesis, so it induces an isomorphism on residue fields $A/f^{-1}(\fm) \xrightarrow{\sim} B/\fm$; so $\phi_1 = \phi_2$, as required.
\end{proof}

\begin{rem}\label{rem:t=T}
Although we won't need this in the sequel, the morphism ~\eqref{eqn:morphism of endos}
is in fact a morphism of~$B^\bG$-algebras. %
To see this, writing $\mSpec$ for the set of maximal ideals of a ring,
it follows from the proof of Proposition~\ref{prop:computation-End-U}
(in particular,
from the commutativity of~\eqref{compatibility of F_U and Pro F_thetabar})
that
the diagram
\begin{equation*}
\begin{tikzcd}
&\mSpec B^\bG \\
\mSpec \F[\HeckeT, f(\HeckeT)^{-1}] \arrow[ur]& \arrow[l, "\Functor_U"] \mSpec S  \arrow[u]
\end{tikzcd}
\end{equation*}
commutes.
Since the formation of~\eqref{compatibility of F_U and Pro F_thetabar} is compatible with finite base change in~$\F$, we see that in fact the diagram
\begin{equation*}
\begin{tikzcd}
&\mSpec(B^\bG \otimes_\F \F')\\
\mSpec \F'[\HeckeT, f(\HeckeT)^{-1}] \arrow[ur]& \arrow[l, "\Functor_U"] \mSpec(S \otimes_\F \F') \arrow[u]
\end{tikzcd}
\end{equation*}
commutes whenever~$\F'$ is a finite extension of~$\F$, and so it also commutes when~$\F' = \cbF_p$.
Since the functor $A \mapsto \mSpec(A \otimes_\F \cbF_p)$ is faithful on reduced $\F$-algebras of finite type (with $\F$-linear morphisms)
we conclude that the diagram
\begin{equation*}
\begin{tikzcd}
& B^\bG \arrow[dl] \arrow[d] \\
 \F[\HeckeT, f(\HeckeT)^{-1}] \arrow[r, "\Functor_U"] &   S  
\end{tikzcd}
\end{equation*}
commutes, as claimed.
\end{rem}

The following result is a mild reformulation of results from
Section~\ref{subsec:Ext-group-computations}. %

\begin{prop}\label{properties of Hecke modules sigma version statements we use}%
For any choice of $\sigma_1,\sigma_2\in\{\sigma,\sigmacomp\}$,  and for
each~$i\ge 0$:
\begin{enumerate}
 \item the $\F[\HeckeT, f(\HeckeT)^{-1}]$-module
\[
\Ext^i_{D^b_{\fp}(\cA_U)}( j^*\cInd_{KZ}^G \sigma_1,
j^*\cInd_{KZ}^G \sigma_2) 
\]  
is finitely generated, where the $\F[\HeckeT, f(\HeckeT)^{-1}]$-action is defined through the action of
$\cH_G(\sigma_2) = \bF[\HeckeT]$ on the second factor.
\item if~$\thetabar$ corresponds to a point of~$U$,
  then the map 
  \[\ihat^*_{\thetabar}:
  \Ext^i_{D^b_{\fp}(\cA_{U})}(j^*\cInd_{KZ}^G\sigma_1, j^*\cInd_{KZ}^G\sigma_2) \to \Ext^i_{\Pro D^b(\cA_{\thetabar}^{\fp})}(\ihat_{\thetabar}^*\cInd_{KZ}^G\sigma_1, 
  \ihat_{\thetabar}^*\cInd_{KZ}^G\sigma_2)\]
  exhibits the target as the completion of the source at the principal ideal of~$\F[\HeckeT, f(\HeckeT)^{-1}]$ generated by~$f_{\thetabar}(\HeckeT)$.
\end{enumerate}
\end{prop}
\begin{proof}
  Note that %
  \[
\begin{tikzcd}
\Ext^i_{D^b_{\fp}(\cA_U)}( j^*\cInd_{KZ}^G \sigma_1,
j^*\cInd_{KZ}^G \sigma_2) \arrow[r, "\sim", "j_*"'] & 
\Ext^i_{D(\cA)}( j_*j^*\cInd_{KZ}^G \sigma_1,
j_*j^*\cInd_{KZ}^G \sigma_2)
\end{tikzcd}
\]
(because $j_*$ is fully faithful on the level of derived categories),
and that
$\Ext^i_\cA = \Ext^i_{D(\cA)}$ by definition.
Part~(1) is then part of Proposition~\ref{properties of Hecke modules}, while 
part~(2) is  Proposition~\ref{Ext and completion}. %
\end{proof}

We can now deduce the following key statement.
\begin{prop}%
\label{prop:full faithfulness on localized weights bis} %
For any choice of $\sigma_1,\sigma_2\in\{\sigma,\sigmacomp\}$, and any~$i\ge 0$, the
functor~$\Functor_{U}$ %
induces an isomorphism \[\Ext^i_{D^b_{\fp}(\cA_{U})}\bigl( j^*\cInd_{KZ}^G \sigma_1,
j^*\cInd_{KZ}^G \sigma_2\bigr)\isoto\Ext^i_{ D^b_{\coh}(\cU)}\bigl(\Functor_{U}(
j^{*}\cInd_{KZ}^G \sigma_1),
\Functor_{U}(j^{*}\cInd_{KZ}^G \sigma_2)\bigr). \]%
\end{prop}
\begin{proof} 
The case $i  = 0$ and $\sigma_1 = \sigma_2$ is  immediate from
Proposition~\ref{prop:computation-End-U} (bearing in mind the symmetry
between~$\sigma$ and~$\sigmacomp$), which shows that
\begin{multline}\label{recalling isomorphism between End}
\F[\HeckeT, f(\HeckeT)^{-1}] \buildrel \eqref{eqn:end-of-j*sigma-explicit} \over = \End_{\cA_U^{\fp}}(j^*\cInd_{KZ}^G \sigma_2)
\\
\iso \End_{\Coh(\cU)}\bigl(\Functor_U(j^*(\cInd_{KZ}^G \sigma_2)\bigr).
\end{multline}
By Proposition~\ref{properties of Hecke modules sigma version statements we use}
and Corollary~\ref{cor:coherent exts for affine mod reductive}, it follows that
for each $i \geq 0$ and pair $\sigma_1,\sigma_2$,
both
\[
\begin{tikzcd}
\Ext^i_{D^b_{\fp}(\cA_U)}\bigl( j^*\cInd_{KZ}^G \sigma_1,
j^*\cInd_{KZ}^G \sigma_2\bigr) 
\end{tikzcd}
\]
and
\[\Ext^i_{D^b_{\coh}(\cU)}\bigl(\Functor_U(j^*
\cInd_{KZ}^G \sigma_1),
\Functor_U(j^*\cInd_{KZ}^G \sigma_2)\bigr)\]
are finitely generated over the finite type 
$\F$-algebra $\F[\HeckeT, f(\HeckeT)^{-1}]$ (acting through the second factor). 

Since~\eqref{recalling isomorphism between End} is induced by~$\Functor_U$,
the morphism
\[
\Functor_U: \Ext^i_{D^b_{\fp}(\cA_U)}\bigl(j^*\cInd_{KZ}^G \sigma_1, j^*\cInd_{KZ}^G \sigma_2\bigr)
\to
\Ext^i_{D^b_{\coh}(\cU)}\bigl(\Functor_U(j^*\cInd_{KZ}^G \sigma_1), \Functor_U(j^*\cInd_{KZ}^G \sigma_2)\bigr) %
\]
is $\F[\HeckeT, f(\HeckeT)^{-1}]$-linear, and so
(since we have now seen that it is an $\F[\HeckeT,f(\HeckeT)^{-1}]$-linear morphism
between finitely generated $\F[\HeckeT,f(\HeckeT)^{-1}]$-modules)
to show that it is an isomorphism,
it suffices to prove that it
becomes one %
after completing at each of the principal ideals~
$f_{\thetabar}(\HeckeT)$ of $\F[\HeckeT, f(\HeckeT)^{-1}].$ %
By~\eqref{eqn:completing our functor on Y II}
there is a commutative diagram %
\[
 \begin{adjustbox}{max width = \textwidth}
\begin{tikzcd}
\Ext^i_{D^b_{\fp}(\cA_U)}(j^*\cInd_{KZ}^G \sigma_1, j^*\cInd_{KZ}^G \sigma_2) \arrow[r, "\Functor_U"] \arrow[d, "\ihat^*_{\thetabar}"] &
\Ext^i_{D^b_{\coh}(\cU)}\bigl(\Functor_U(j^*\cInd_{KZ}^G \sigma_1), \Functor_U(j^*\cInd_{KZ}^G \sigma_2)\bigr) \arrow[d, "\ihat'^*_{\thetabar}"]\\
\Ext^i_{\Pro D^b_{\fp}(\cA_{\thetabar})}(\ihat^*_{\thetabar}\cInd_{KZ}^G \sigma_1, \ihat^*_{\thetabar}\cInd_{KZ}^G \sigma_2) \arrow[r, "\ProFunctor_{\thetabar}"] & %
\Ext^i_{\Pro D^b_{\coh}(\cX_{\thetabar})} \bigl(\ProFunctor_{\thetabar}(\ihat^*_{\thetabar}\cInd_{KZ}^G \sigma_1), \ProFunctor_{\thetabar}(\ihat^*_{\thetabar}\cInd_{KZ}^G \sigma_2)\bigr)
\end{tikzcd}
 \end{adjustbox}
\]
Proposition~\ref{properties of
  Hecke modules sigma version statements we use} shows that the left vertical arrow is a completion at the principal ideal generated by $f_{\thetabar}$. 
  Corollary~\ref{cor:completing RHoms} shows the same is true of the right vertical arrow. 
It only remains to show that the lower horizontal arrow is an isomorphism, and this is immediate from Theorem~\ref{thm:derived-F-thetabar-nonSt-case}.
\end{proof}

\begin{prop}%
\label{prop:full faithfulness on U(sigmasigmacomp)}For each companion
pair~$\{\sigma,\sigmacomp\}$, the functor 
\[
\Functor_{U(\sigmasigmacomp)}: \Dfp^b(\cA_{U(\sigmasigmacomp)})
\to D^b_{\coh}\bigl(\cU(\sigmasigmacomp)\bigr)
\]
is fully faithful.
\end{prop}
\begin{proof}
  By Lemma~\ref{lem:F-fully-faithful-iff-IndProF}, it suffices to show that 
  \[
  \Functor_{U(\sigmasigmacomp)}: D(\cA_{U(\sigmasigmacomp)})
\to \Ind D^b_{\coh}\bigl(\cU(\sigmasigmacomp)\bigr) 
  \]
is fully faithful. Using Proposition~\ref{prop:full faithfulness on
    localized weights bis}, this is immediate from an application of the implication ``$(4)\implies (1)$''  %
  of Proposition~\ref{prop:check-full-faithful-on-compact-generators}, using the compact generators $j^*\cInd_{KZ}^G\sigma$, $j^*\cInd_{KZ}^G\sigmacomp$ %
  of $D(\cA_{U(\sigmasigmacomp)})$. %
\end{proof}  

Finally, we deduce the main results of this section.
\begin{thm}%
  \label{thm:full faithfulness on U}\leavevmode
  \begin{enumerate}
  \item\label{item:25} The functor $\Functor_{\Ugood}: \Dfp^b(\cA_{\Ugood}) \to D^b_{\coh}(\cUgood)$ is fully faithful.
\item\label{item:26} The natural transformation 
\[\Functor j_{\Ugood,*} \to j'_{\cUgood,*}\Functor_{\Ugood}\]
of functors from $\Dfp^b(\cA_{\Ugood})$ to $\Ind D^b_{\coh}\bigl(\cX)$
is a natural isomorphism.
  \end{enumerate}
\end{thm}
\begin{proof}
Recalling  the
decomposition  $\Ugood = \coprod U(\sigmasigmacomp)$, this is immediate from 
Propositions~\ref{prop:full faithfulness on U(sigmasigmacomp)} and~\ref{prop:F commutes with j lower star}.
\end{proof}

\subsection{The main theorem %
}\label{subsec:proving-that-F-is-fully-faithful}
We are now in a position to deduce our main theorem, showing that if~$p\ge 5$, then there is a fully faithful functor $D(\cA) = \Ind \Dfp^b(\cA) \to \Ind D^b_{\coh}(\cX)$. %

\begin{thm}%
  \label{thm:F-is-fully-faithful} The functors
  \begin{equation}\label{eqn: to prove fully faithful for main theorem}
  \Functor:\Dfp^b(\cA) \to D^b_{\coh}(\cX)
  \end{equation}
  and
  \begin{equation}\label{eqn: to prove fully faithful for main theorem Ind version}
  \Functor: D(\cA) = \Ind \Dfp^b(\cA) \to \Ind D^b_{\coh}(\cX) 
  \end{equation}
  are fully
   faithful. 
\end{thm}
\begin{proof}We put ourselves in the situation of Section~\ref{subsec:semiorthogonal and Ind Pro}, taking the  category  $\catD^{c}$  there to be our $D^b_{\fp}(\cA)$; so (in a somewhat
unfortunate clash of notation) it follows from Corollary~\ref{cor:all-the-nice-properties-of-D(A)} that $\catD$ there is our
$D(\cA)$. %
We take  $(\catD')^{c}$ there 
  to be our $D^b_{\coh}(\cX)$. 
We let
$\catB^{c}$ be~$D^b_{\fp}(\cA_{\Ybad})$, and $(\catB')^{c}$ be $D^b_{\coh,\Ybad}(\cX) =  D^b_{\coh}(\cX_{\Ybad})$. %
  
  The theorem will be an immediate consequence of Proposition~\ref{prop:gluingfull-faithfulness-Ind-Pro-version-compact-version}, once we check that the functor~\eqref{eqn: to prove fully faithful for main theorem} satisfies the hypotheses of that result.
   Note firstly that by Corollary~\ref{cor:all-the-nice-properties-of-D(A)} and \eqref{eqn: computing quotient category in coherent recollement compact version}, we have equivalences 
   \[
   D^b_{\fp}(\cA)/D^b_{\fp}(\cA_{\Ybad})\to \Dfp^b(\cA_{\Ugood})
   \] 
   and 
   \[
   D^b_{\coh}(\cX)/D^b_{\coh}(\cX_{\Ybad}) \to D^b_{\coh}(\cUgood).
   \]
   Then Hypothesis~\eqref{item:78} is the claim that the functors
    \[
    \Functor_{\Ugood}: \Dfp^b(\cA_{\Ugood}) \to D^b_{\coh}(\cUgood)
    \]
    and
    \[
    \Functor_{\Ybad}: \Dfp^b(\cA_{\Ybad}) \to D^b_{\coh}(\cX_{\Ybad})
    \]
    are fully
    faithful, which is Theorem~\ref{thm:full faithfulness on U}~\eqref{item:25} and
    Corollary~\ref{cor:overview fully faithful on closed}.
Hypothesis~\eqref{item:76} is precisely Theorem~\ref{thm:full faithfulness on U}~\eqref{item:26}, and finally Hypothesis~\eqref{item:77} is Proposition~\ref{prop:completing our functor on Y}    (with~$Y=\Ybad$).
\end{proof}

\appendix
\section{Category theory}\label{app:category-theory}%
We begin by briefly collecting some standard results in category theory. In
addition to the material recalled here, we refer the reader to~\cite[App.\
A.3]{DEGlocalization} for some recollections on~$\Ind$ and~$\Pro$ categories.
\subsection{1-categories}\label{subsec:1-categories} As
we discuss in more detail at the beginning of Section~\ref{subsec: category
  theory background},  we typically fix a 
Grothendieck universe, %
and sets are called \emph{small} if they belong to this fixed universe.
\subsubsection{Subcategories of abelian categories}
\label{subsubsec:Serre subcategories}
Let~$\cA$ be an abelian category.
By definition, 
a \emph{Serre subcategory} of~$\cA$ is a full subcategory $\cB$ of~$\cA$ such
that, for any exact
sequence \(X\to Y\to Z\) for which $X$ and $Z$ are objects of $\cB$, we also have that~$Y$
is an object of $\cB$; equivalently, $\cB$ is closed under passing to subquotients 
and extensions in~$\cA$ (cf.~\cite[\href{https://stacks.math.columbia.edu/tag/02MP}{Tag 02MP}]{stacks-project}). 
A \emph{localizing subcategory} of~$\cA$ is a Serre subcategory $\cB$ of $\cA$
for which the quotient functor 
$\cA\to\cA/\cB$ admits a right adjoint, which is then necessarily fully faithful, by~\cite[Prop.\ III.2.3(a)]{Gabrielthesis}. 

Following
\cite[\href{https://stacks.math.columbia.edu/tag/02MN}{Tag
02MN}]{stacks-project}, we define %
a {\em weak Serre subcategory}
of~$\cA$ to be a full subcategory~$\cB$ of~$\cA$ 
with the additional property that, whenever the four outermost terms of a five term exact sequence
in~$\cA$ are objects of~$\cB$, so is the middle term. Any weak Serre
subcategory~$\cB$ (and thus any Serre subcategory) of~$\cA$ is an \emph{exact abelian
  subcategory} of~$\cA$; that is, it is a full subcategory of~$\cA$
which is an abelian category in its own right, and for which the inclusion
$\cB\hookrightarrow \cA$ is exact. (Indeed, a weak Serre subcategory of~$\cA$
is precisely an exact abelian subcategory of~$\cA$ which is also 
closed under the formation of extensions in~$\cA$; 
c.f.\ \cite[\href{https://stacks.math.columbia.edu/tag/0754}{Tag 0754}]{stacks-project}.)

\subsubsection{Grothendieck categories} Recall that a Grothendieck
category~$\cA$ is an abelian category which is cocomplete, admits a set
of generators, and for which the formation of filtered colimits in $\cA$ is
exact. It then follows that~$\cA$ is also complete~\cite[Prop.\ 8.3.27]{MR2182076}.

We say that an object~$X$ of a Grothendieck category~$\cA$ is compact
if $\Hom_\cA(X,\text{--})$ commutes with filtered colimits, and we write~$\cA^c$ for the full
subcategory of compact objects. We say that the
 category~$\cA$ is compactly generated if it admits a set of compact
generators, in which case the natural functor $\Ind(\cA^c)\to\cA$ is an equivalence. (In the literature, compact objects are often called ``finitely
presented'', and a compactly generated category is ``locally finitely
presented''.)%

Recall that an abelian category is called \emph{locally Noetherian} if it is
a Grothendieck category, and furthermore admits a set of Noetherian generators. A locally
Noetherian abelian category is compactly generated, and its compact objects are
precisely the Noetherian objects  (see e.g.\ \cite[Prop.\
A.1.1]{DEGlocalization}).

An abelian category~$\cA$ is  {\em locally coherent} if it is a
compactly generated   Grothendieck category and if furthermore~$\cA^c$ is abelian, in which case~$\cA^{c}$ is a
weak Serre subcategory, hence an
exact abelian subcategory,
of~$\cA$ %
(for example by Lemma~\ref{lem:locally-coherent-implies-weak-Serre} below). %
Any locally Noetherian category is
locally coherent.

A Serre subcategory~$\cB$ of a Grothendieck category~$\cA$ is localizing
if and only if it is closed under arbitrary direct sums in~$\cA$. Suppose
furthermore that~$\cA$ is locally Noetherian; then~$\cB$ and ~$\cA/\cB$ are also
locally Noetherian, and the right adjoint $\cA/\cB\to\cA$ preserves filtered
colimits. Furthermore the compact objects of $\cB$ are those objects of~$\cB$
which are compact in~$\cA$, and the compact objects of~$\cA/\cB$ are the images
of the compact objects in~$\cA$, so that the sequence of
functors \[\cB\to\cA\to\cA/\cB\] is the $\Ind$-extension of the induced
sequence \[\cB^c\to\cA^c\to(\cA/\cB)^c=\cA^c/\cB^{c}.\](See for example~\cite[Lem.\
A.2.7]{DEGlocalization} for these facts.)

\begin{lem}\label{lem:locally-coherent-implies-weak-Serre}
  If~$\cA$ is a locally coherent abelian category, then~$\cA^c$ is a weak Serre
  subcategory of~$\cA$. If~$\cA$ is furthermore locally Noetherian, then~$\cA^{c}$ is a
  Serre subcategory of~$\cA$.
\end{lem}
\begin{proof}
For the first statement, we have
to check that~$\cA^c$ is an exact abelian subcategory of~$\cA$ which is closed under extensions.
  This is consequence of the equivalence of~(1) and~(2) in~\cite[Thm.\ 1.6]{MR1434441}, together with~\cite[Prop.\ 1.5]{MR1434441}.
  
  For the second statement,
we have to check that~$\cA^{c}$ is closed under
  subquotients and extensions, which is clear, since in this case (i.e.\ when $\cA$ is locally
Noetherian)~$\cA^c$ coincides with the
  Noetherian objects in~$\cA$.
  \end{proof}

 \subsubsection{Ind and Pro categories}\label{subsubsec:Ind and Pro categories}
 If~$\cC$ is a small category, then, as usual,
we write $\Ind(\cC)$ for the category of small filtered diagrams of objects of~$\cC$.
We typically consider directed systems $\{X_i\}_{i\in I}$ of objects of~$\cC$,
indexed by a directed set~$I$, and write $\quotescolim{i \in I} X_i$ to denote
the associated object of~$\Ind(\cC)$.  (Of course, one can equally well consider
systems of objects indexed by a filtered category, and we will do so on occasion,
using the same notation.) 
Regarding objects of~$\cC$ as constant directed systems, we obtain a fully
faithful embedding~$\cC \hookrightarrow \Ind(\cC)$.
The object $\quotescolim{i\in I} X_i$ is
then
the colimit of the $X_i$ in $\Ind \cC$.

Recall that morphisms in $\Ind \cC$ are computed via the formula
 \[
 \Hom_{\Ind(\cC)}(\quotescolim{i \in I}X_i, \quotescolim{j \in J}Y_j) = \lim\nolimits_{i \in I}\colim\nolimits_{j \in J}\Hom_\cC(X_i, Y_j).
 \]
(This simply encodes the fact, already noted, that the natural functor $\cC \to \Ind(\cC)$ is
fully faithful, %
together with the additional fact that the objects of $\cC$ are compact in $\Ind \cC$.)

 We define~$\Pro(\cC)$ in a dual fashion, so that objects of ~$\Pro(\cC)$ can be written as cofiltered limits
 $\quoteslim{i \in I} X_i$ of  objects of~$\cC$. %
 We will usually apply these constructions when~$\cC$ is a small abelian category,
in which case~$\Ind(\cC)$ and~$\Pro(\cC)$ are also abelian (but no longer small).

We have the following standard lemmas. %

 \begin{lem}
  \label{lem:exactness-of-limits}
  Let $\cC$ and~$\cC'$ be small abelian categories, %
  and let~$f: \cC \to \cC'$ be an additive functor.
  Then:
  \begin{enumerate}
  \item 
  $f$ is left exact (resp.\ right exact, resp.\ exact) if and only if its cofiltered limit-preserving extension $\Pro(f):\Pro(\cC)\to\Pro(\cC')$ is left exact (resp.\ right exact, resp.\ exact).  
  \item
  If~$\cC'$ has cofiltered limits, then $f$ is left exact if and only if its cofiltered limit-preserving extension $\Pro(f): \Pro(\cC) \to \cC'$ is left exact;
  \item If~$\cC'$ has exact cofiltered limits, then $f$ is right exact if and only if~$\Pro(f) : \Pro(\cC) \to \cC'$ is right exact.
  \end{enumerate}
\end{lem}
 \begin{proof}
  By for example \cite[Prop.\ 8.6.6(a)]{MR2182076} any exact
  sequence $0\to A\to B\to C \to 0$
  in ~$\Pro(\cC)$ can be written as a cofiltered limit of  exact sequences
  in~$\cC$. 
  By~\cite[Theorem\ 8.6.5(ii)]{MR2182076}, the functor $\cC \to \Pro(\cC)$ is fully faithful and exact.
  Part~(1) follows, and then %
 part~(2), resp.\ part~(3), follows from part~(1) and the fact that the inverse limit functor $\Pro(\cC') \to \cC'$ is left exact, resp.\ exact under the assumptions of part~(3).
\end{proof}

  \begin{lem}
    \label{lem:left-adjoint-to-Pro-functor-abelian}Suppose that $g:\cC'\to\cC$ is an exact functor between small
  abelian categories, and continue to write $g:\Pro\cC'\to\Pro\cC$ for its
  $\Pro$-extension. Then:
  \begin{enumerate}
  \item $g:\Pro\cC'\to\Pro\cC$ admits a left adjoint $f:\Pro \cC \to \Pro \cC'$,
    which is right exact and cofiltered limit-preserving.
  \item $f$ is exact if and only if the restriction $f|_{\cC}:\cC\to\Pro\cC'$ is exact.
  \end{enumerate}
  \end{lem}
  \begin{proof}
    By assumption, $g:\Pro\cC'\to\Pro\cC$ is exact and cofiltered limit-preserving, so it is limit-preserving, and thus admits a left adjoint 
    $f:\Pro \cC \to \Pro \cC'$ by the adjoint functor theorem. 
    Being a left
    adjoint, $f$ is right exact. 
    We now prove that it preserves cofiltered limits.
    Letting $X = \lim_{i\in I}X_i$ be a cofiltered limit of objects of~$\Pro \cC$, and $Y$ be an object of~$\cC'$, %
    we have
    \begin{multline*}
      \Hom_{\Pro\cC'}(f(\lim_{i}X_i),Y)=\Hom_{\Pro\cC}(\lim_{i}X_i,g(Y))=\colim_i\Hom_{\Pro\cC}(X_i,g(Y))\\ =\colim_i\Hom_{\Pro\cC'}(f(X_i),Y)=\Hom_{\Pro\cC'}(\lim_if(X_i),Y),
    \end{multline*}
    where the second equality is because $g(Y) \in \cC$.
    Since every object of~$\Pro \cC'$ is a limit of objects of~$\cC'$, this equality remains true for arbitrary objects $Y \in \Pro \cC'$.
    Hence $f(\lim_iX_i)=\lim_if(X_i)$, as required.
The second part follows from parts~(2) and~(3) of Lemma~\ref{lem:exactness-of-limits}, since~$\Pro \cC'$ has exact cofiltered limits. %
  \end{proof}
  \begin{rem}
    \label{rem:ind-pro-extension-of-adjoint}In the setting of
    Lemma~\ref{lem:left-adjoint-to-Pro-functor-abelian}, we may form the
    $\Ind$-extensions of~$f$ and~$g$, which are a pair of adjoint functors $\Ind
    f:\Ind\Pro\cC\to\Ind\Pro\cC'$, $\Ind g:\Ind\Pro\cC'\to\Ind\Pro\cC$ (as can be checked directly from the definition of $\Hom$-spaces in $\Ind$-categories). Since
    $\Ind$-completion preserves exactness, we see that~$\Ind g$ is always exact,
    and  if
    $f|_{\cC}:\cC\to\Pro\cC'$ is exact, then so is~$\Ind f$.
  \end{rem}

We now recall a standard result about 
projective objects of~$\Pro \cC$.
(The dual statement then gives a result about injective objects in~$\Ind \cC$.) 
We begin with the following technical lemma, which is Exercise 6.11 in~\cite{MR2182076}.

\begin{lem}\label{lem:subobjects-of-Ind-category}%
  Let $\cC'$ be a full subcategory of a small category~$\cC$. %
  Then:
  \begin{enumerate}
  \item if $A= \quotescolim{i} X_i\in\Ind\cC$ is such that any $X\to A$ with $X\in \cC$ factors through some $Y \to A$ with $Y\in \cC'$, then $A\in \Ind\cC'$.
  \item  If $A= \quoteslim{i} X_i\in\Pro\cC$ is such that any $A\to X$ with $X\in \cC$  factors through some $A\to Y$ with $Y\in \cC'$, then $A\in \Pro\cC'$.
  \end{enumerate}
\end{lem}
\begin{proof}  The second statement follows from the first by passage to opposite categories, so it suffices to prove the first.
To this end,
write $\cD$ to denote the category of morphisms (in $\Ind \cC$) $Y \to A$, with $Y$ an
object of $\cC'$. It follows from the hypothesis (by a very similar argument to that in the last paragraph of this proof) that the category~$\cD$ is filtered. Accordingly, the functor $\alpha: \cD \to \cC'$ mapping the morphism $Y\to A$ 
to $Y$ gives a filtered diagram in $\cC'$, 
so it has colimit $\quotescolim{\cD} Y$ in $\Ind \cC'$.

The morphisms $Y \to A$ that define the objects of $\cD$ give rise to a tautological
morphism $\quotescolim{\cD} Y \to A,$ which we will show is an isomorphism. By Yoneda's lemma, it suffices to show that for any~$X\in \cC$, the induced morphism %
\begin{equation}
  \label{eq:quoteslim-to-A}
  \Hom_{\Ind \cC}(X, \quotescolim{\cD} Y) \to \Hom_{\Ind \cC}(X, A)
\end{equation}
is an isomorphism.

Since~$X\in \cC$, we have 
\[
\Hom_{\Ind \cC}(X, \quotescolim{\cD} Y)=\colim_{\cD}\Hom_{\cC}(X,Y),
\] 
and the surjectivity of~\eqref{eq:quoteslim-to-A} follows from the hypothesis on~$A$. For the injectivity, suppose that $X \rightrightarrows \quotescolim{\cD} Y$ are two morphisms %
that induce the same morphism to~$A$. These morphisms then factor through some common~$Y$,
and we have that the two composites $X \rightrightarrows Y \to A$ coincide.
Now the morphism $Y \to A$ factors through some~$X_i$, and if we choose $i$ large enough,
the composites $X \rightrightarrows Y \to X_i$ will coincide.  
Finally, the map $X_i \to A$ factors through some $Z \to A$ with $Z$ an object of~$\cC'$,
by hypothesis.  So our morphisms $X \rightrightarrows \quotescolim{\cD} Y$
factor through the tautological map $Z \to \quotescolim{\cD} Y$, 
via $X \rightrightarrows Y \to X_i \to Z,$ and already coincide as maps to~$Z$.  Hence
they coincide, and we've proved injectivity.
\end{proof}

\begin{lem}
  \label{lem:quasi-projectives-in-Pro-cat}
  Let~$\cC$ be a small abelian category.
  Then~$\Pro \cC$ has enough projectives.
  If~$\cC$ has enough projectives, then every projective object in~ $\Pro\cC$ is of the form $\quoteslim{i} P_i$ with the $P_i$ projective objects of~$\cC$.
\end{lem}
\begin{proof}
The first claim is dual to~\cite[Corollary~9.6.5]{MR2182076}.
We now prove the second claim.
The analogous result for injective objects of Ind-categories is immediate from~\cite[Prop.~15.2.3]{MR2182076}, and we follow the proof of that result.
Let~$P$ be a projective object of~$\Pro\cC$, and let $P\to X$ be a morphism in~$\Pro\cC$, with ~$X$  an object of $\cC$.
Since~$\cC$ has enough projectives, there exists a surjection $Q\to X$, with~$Q$ a projective object of~$\cC$.
Since~$P$ is projective, the morphism $P\to X$ factors as $P\to Q\to X$. The result follows from Lemma~\ref{lem:subobjects-of-Ind-category} (with~$\cC'$ there being the full subcategory of projective objects of~$\cC$). %
\end{proof}

We finish up this discussion of Pro-categories by recording a technical result
related to countably indexed pro-objects.

\begin{df} If $\cC$ is a small category,
we say that an object $X$ of $\Pro \cC$
is countably indexed if we can write 
$X = \quoteslim{n} X_n$, for some sequence of morphisms
$\cdots X_{n+1} \to X_n \to \cdots \to X_0$ in~$\cC$.
\end{df}

\begin{lem}
\label{lem:countably indexed}
Let $\cC$ be a small abelian category,
and $0 \to A \to B \to C \to 0$ a short exact sequence in $\Pro \cC$. 
If each of $A$ and $C$ is countably indexed,
then the same is true of~$B$.
Furthermore,
if $A = \quoteslim{n} A_n$ for some sequence of morphisms
$\cdots \to A_{n+1} \to A_n \to \cdots \to A_0$ in~$\cC$,
then we can write $0 \to A \to B \to C \to 0$
as a projective limit $0 \to A_n \to B_n \to C_n \to 0$
of short exact sequences in~$\cC$, with the same transition maps $A_{n+1} \to A_n$ as before. 
\end{lem}
\begin{proof}
If we push out the given exact sequence along the various morphisms
$A\to A_n$, we obtain exact sequences
\begin{equation}
\label{eqn:pushed out ses}
0 \to A_n \to B_n \to C \to 0
\end{equation}
which are compatible as $n$ varies.
Since the formation of cofiltered limits is exact in~$\Pro \cC$,
we find that the limit of these sequences recovers the original
short exact sequence; in particular, $B \iso \lim_n B_n$.

The class of $B_n$ gives an element in
$\Ext^1_{\Pro \cC}(C,A_n)$.  If we write $C = \quoteslim{m} C_m$,
then~\cite[Lems.~A.3.13, A.3.14(2)]{DEGlocalization} 
give an isomorphism $$\colim_m \Ext^1_{\cC}(C_m,A_n)
\iso \Ext^1_{\Pro \cC}(C,A_n).$$
In particular, for some $m_n$ there is an extension
\begin{equation}
\label{eqn:bottom level ses}
0 \to A_n \to B_{n,m_n} \to C_{m_n} \to 0
\end{equation}
from which~\eqref{eqn:pushed out ses} is obtained by pulling back along
the morphism $C\to C_{m_n}$.  If
\begin{equation}
\label{eqn:finite level ses}
0 \to A_n \to B_{n,m} \to C_m \to 0
\end{equation}
denotes the short exact sequence in~$\cC$, obtained 
by pulling back~\eqref{eqn:finite level ses} along the transition
morphism $C_m \to C_{m_n}$ (for $m\geq m_n$),
then we find that~\eqref{eqn:pushed out ses} maps isomorphically
to the limit of these short exact sequences,
and in particular, that $B_n \iso \quoteslim{m}B_{n,m}.$
Consequently, we see that $B$ itself is countably indexed.

Choosing an appropriate strictly increasing function $m = m(n)$,
and considering the corresponding diagonal system of short 
exact sequences
$$
0 \to A_n \to B_{n,m(n)} \to C_{m(n)} \to 0 %
$$
we obtain an inverse system of short exact sequences in~$\cC$
whose limit over~$n$ coincides with the given short exact sequence
$0 \to A \to B \to C \to 0.$
(More precisely, suppose we have chosen~$m(n)$ suitably.
The composite of the morphism $B_{n+1} \to B_n$ (that arises via their construction 
as pushouts) with the morphism $B_n \to B_{n,m(n)}$ factors through
some $B_{n+1,m(n+1)},$ inducing a corresponding morphism of
short exact sequences~\eqref{eqn:finite level ses}
from level $\bigl(n+1,m(n+1)\bigr)$ to $\bigl(n,m(n)\bigr).$)
This concludes the proof.\qedhere
\end{proof}

\subsubsection{Socle and radical filtrations.}\label{socle and radical}
We briefly recall some facts about socle and radical filtrations in an abelian category~$\cC$.
We say that an object~$M \in \cC$ is \emph{simple} if its only subobjects (equivalently, quotients) are~$0$ and~$M$.
We say that~$M$ is \emph{semisimple} if it is a direct sum of simple objects, and that~$M$ has \emph{finite length} if it has a finite composition series, i.e.\ a finite filtration
whose graded pieces are simple.
The length of a composition series is then independent of the choice of filtration, and defined to be the length of~$M$.
We introduce the following notation.

\begin{defn}\label{Cfl}
Let~$\cC$ be an abelian category.
We write~$\cC^{\fl}$ for the full subcategory of~$\cC$ of objects of finite length.
\end{defn}

When~$\cC = \Mod(R)$ for a ring~$R$, we will also write~$\Mod^{\fl}(R)$ for~$\cC^{\fl}$.
Assume now that~$\cC$ is complete and cocomplete (e.g.\ a Grothendieck category). 
If~$M \in \cC$, we let
\[\rad(M) \coloneqq \ker \left (\prod_{q: M \to N : N \text{ is simple }}q \right ) = \bigcap_{q: M \to N : N \text{ is simple }}\ker(q)\]
and
\[\soc(M) \coloneqq \operatorname{im}\left ( \coprod_{q: N \to M: N \text{ is simple }} q\right ).\]
Then we define %
\[
\rad^0(M) \coloneqq M, \rad^{i+1}(M) \coloneqq \rad(\rad^i(M))
\]
and
\[\soc_{-1}(M)\coloneqq 0, %
  \soc_{i+1}(M) \coloneqq \text{preimage in $M$ of } \soc(M/\soc_i(M)).\]
Finally, we set
\[\cosoc(M) \coloneqq M/\rad(M).\]
If we write~$\vee : \cC \to \cC^\op$ for the natural anti-equivalence, an induction on~$i$ then shows that
\[
(M/\rad^i M)^\vee = \soc^i(M^\vee).
\]
As the following lemmas show, these notions behave mostly as expected in Grothendieck categories: for example, since coproducts are then exact, we find that
\begin{equation}\label{eqn:socle in Grothendieck categories}
\soc(M) = \sum_{q: N \to M: N \text{ is simple }} q(N),
\end{equation}
hence~$\soc(M)$ is the sum of all simple subobjects of~$M$.
Lemma~\ref{properties of semisimple objects}~(2) then shows that
$\soc(M)$ is a semisimple subobject of~$M$ (indeed, the maximal semisimple subobject of~$M$). 
However, $\cosoc(M)$ need not be a semisimple quotient of~$M$: 
see Remark~\ref{counterexample for semisimple objects} for a counterexample, which also shows that the next lemma need not hold in an arbitrary abelian category.

\begin{lemma}\label{properties of semisimple objects}
Let~$\cC$ be an abelian category, and let~$M$ be an object of~$\cC$.
Assume furthermore either that $\cC$ is a Grothendieck category, or that $M$ is of finite length.
\begin{enumerate}
\item If $M$ is semisimple,
then every subobject, resp. quotient, of~$M$ is semisimple, and every exact sequence
$0 \to M' \to M \to M'' \to 0$
splits.
\item
If $M = \sum_{i \in I} M_i$ for a set~$I$ of simple subobjects,
then $M$ is semisimple.
\end{enumerate}
\end{lemma}
\begin{proof}
Assume that $M$ is semisimple,
and consider an exact sequence as in the statement of~(1). 
By assumption, $M = \bigoplus_{i \in I}M_i$ for a set~$I$ of simple objects of~$\cC$.
Furthermore, either $\cC$ is a Grothendieck category, or~$I$ is finite. 
In either case, there is a (possibly empty) subset $J \subseteq I$ of~$I$ maximal with respect
to the property that
\[M' \cap \bigoplus_{i \in J} M_i = 0.\]
If $I$ is finite, this is obvious.  Otherwise~$\cC$ is a Grothendieck category by assumption,
so the formation of filtered colimits is exact in~$\cC$.
Thus, for any $J \subseteq I,$
we have that
\begin{equation}\label{intersection and filtered colimits}
M' \cap \bigoplus_{i \in J}M_i = \varinjlim_{K \subset J \text{ finite}}\left ( M' \cap \bigoplus_{i \in K} M_i \right ),
\end{equation}
and so
the existence of a maximal $J$ follows from Zorn's lemma.

Having chosen such a subset~$J$, we define $N \coloneqq M' + \bigoplus_{i \in J} M_i$.
Note that this sum is direct.
Now if~$N \ne M$ then there exists~$i \in I$ such that $M_i \not \subset N$,
and so $N + M_i$ is a direct sum, which contradicts the maximality of~$J$.
We have thus proved that the exact sequence in~(1) splits, and that~$M''$ is semisimple.

As an aside,
note that this already implies~(2), since (momentarily adopting the notation of 
that statement) we have a surjection $\bigoplus M_i \to M$ with semisimple
source.  

Returning to the proof of~(1),
it remains to show that~$M'$ is semisimple.
We first note that 
if~$M' \ne 0$ then it contains a simple subobject.
Indeed, if $M' \ne 0,$
then~\eqref{intersection and filtered colimits} shows that there exists a finite $K \subset I$
such that $M' \cap \bigoplus_{i \in K} M_i \ne 0$.
Now $M' \cap \bigoplus_{i \in K} M_i$ has finite length (since~$K$ is finite) and so it necessarily contains a simple subobject, which is then a simple subobject of~$M'$.

We now prove that~$M'$ is semisimple.  By~(2) (which we've already proved) it suffices
to show that $M'$ is equal to the sum of its simple subobjects.  To this end, then,
let~$M'_0$ be the sum of the simple subobjects of~$M'$.
By what we've already shown, we know that the inclusion $M'_0 \hookrightarrow M$ is split, and
so the inclusion $M'_0 \hookrightarrow M'$ is also split
(by the restriction to $M'$ of a retraction of
$M'_0 \hookrightarrow M$).  Thus we can write $M' = M'_0 \oplus M'_1$ for some~$M'_1 \subset M'$.
By the construction of~$M'_0,$ we see 
that $M_1'$ contains no simple subobjects. 
By what we've already proved, we conclude that $M_1' = 0$, so that $M_0' = M',$
as required.
\end{proof}

\begin{cor}\label{associated graded of socle filtration is semisimple}
Let~$\cC$ be a complete and cocomplete abelian category, and let~$M \in \cC$.
\begin{enumerate}
\item If~$\cC$ is a Grothendieck category, then~$\soc_i M / \soc_{i-1} M$ is semisimple for all~$i \geq 0$.
\item If~$M$ has finite length, then $\soc_i M / \soc_{i-1} M$ and $\rad^i M / \rad^{i+1} M$ are semisimple for all~$i \geq 0$.
\end{enumerate}
\end{cor} 
\begin{proof}
The statements about~$\soc_i$ are immediate by %
Lemma~\ref{properties of semisimple objects}~(2),
which shows that~$\soc M$ is semisimple, together with the fact that $\soc_iM/\soc_{i-1}M \isoto \soc(M/\soc_{i-1}M)$ for all~$i$, by definition.
Assume therefore that~$M$ has finite length.
Arguing as in the previous case, it suffices to prove that~$M/\rad M$ is semisimple. 
By construction, $M/\rad M$ is a finite length subobject of a product $\prod_{i \in I} M_i$ of simple objects.
Since
\[\prod_{i \in I} M_i = \varprojlim_{K \subset I \text{ finite }}\prod_{i \in K} M_i\]
and~$M/\rad M$ is Artinian, there exists a finite $K \subset I$ such that 
\[
M/\rad M \to \prod_{i \in K} M_i = \bigoplus_{i \in K} M_i 
\]
is a monomorphism.
Since~$K$ is finite, %
Lemma~\ref{properties of semisimple objects}~(1) %
now implies that~$M/\rad M$ is semisimple, as desired.
\end{proof}

\begin{rem}\label{counterexample for semisimple objects}
Lemma~\ref{properties of semisimple objects} need not hold true if~$\cC$ is not a Grothendieck category.
For example, let~$\cC = \Mod(\bZ)^\op$ be the opposite of the category of abelian groups.
The simple objects in~$\cC$ are the groups~$\bZ/q$ for~$q$ a prime number.
Then $G \coloneqq \prod_{q} \bZ/q$ is a semisimple object in~$\cC$, but the subgroup generated by~$1$ in all coordinates is a quotient of~$G$ in~$\cC$ which is isomorphic to~$\bZ$, 
which is not a semisimple object of~$\cC$.

Similarly, now viewing~$G$ as an object of $\cC^\op$ (i.e.\ of $\Mod(\Z)$),
we see that $\rad_{\Mod(\bZ)}(G) = 0$.
However, $G$ is not semisimple, since it contains non-torsion elements;
hence the cosocle of an object need not be semisimple (even in a Grothendieck category).
\end{rem}

\begin{lemma}\label{general Loewy facts}
Let~$\cC$ be an abelian category, and let $M, N$ be objects of~$\cC$.
\begin{enumerate}
\item If~$\cC$ is a Grothendieck category, and $N$ is a subobject of~$M$, then $\soc_i N = N \cap \soc_i M$ for all~$i \geq 0$.
\item If~$\cC^\op$ is a Grothendieck category, and $q: M \to N$ is an epimorphism, then $\rad^i N = q(\rad^i M)$ for all~$i$.
\end{enumerate}
\end{lemma}
\begin{proof}
We begin by proving part~(1), by induction on~$i$. 
For the base case~$i = 0$, the inclusion $\soc N \subset N\cap \soc M$ is true by~\eqref{eqn:socle in Grothendieck categories}. 
On the other hand, the inclusion $N \cap \soc M \subset \soc N$ holds because $N\cap \soc M $ is semisimple, by Lemma~\ref{properties of semisimple objects}~(1),
and so it is the sum of its simple subobjects; and by~\eqref{eqn:socle in Grothendieck categories}, these are also subobjects of~$\soc(N)$.
By inductive assumption, we now have an inclusion $N/\soc_i N \to M/\soc_i M$, and the case~$i = 0$ implies that
\[\soc(N/\soc_i N) = (N/\soc_i N) \cap \soc(M/\soc_i M).\]
Hence the preimage  in~$N$ of $\soc(N/\soc_i N)$, which is $\soc_{i+1} N$, coincides with the preimage in~$N$ of $\soc(M/\soc_i M)$, which is $N \cap \soc_{i+1}M$.
This concludes the proof of part~(1).

Part~(2) now follows by duality, since it is equivalent to the statement that $q|_{\rad^i M}$ is an epimorphism onto~$\rad^i N$, and we already know that
$N^\vee \to M^\vee \to M^\vee/\soc_i M^\vee$ factors through a monomorphism $N^\vee/\soc_i N^\vee \to M^\vee /\soc_i M^\vee$.\qedhere
\end{proof}

We conclude this subsection with the following well-known fact about socle and radical filtrations of objects of finite length.

\begin{lemma}\label{Loewy length}
Let~$\cC$ be a complete and cocomplete abelian category, and let~$M \in \cC$ be an object of finite length.
Then the socle and radical filtration of~$M$ have the same length {\em (}i.e.\ number of nonzero graded pieces{\em )}, called the \emph{Loewy length} of~$M$.
\end{lemma}
\begin{proof}
Let~$m$ be the length of the radical filtration, and~$n$ the length of the socle filtration.
Then~$m$ is minimal such that $\rad^m(M) = 0$, and~$n$ is minimal such that $\soc_{n-1}(M) = M$.
Since, by Corollary~\ref{associated graded of socle filtration is semisimple}, $\soc_i M/\soc_{i-1} M$ is semisimple for all~$i$, 
we see by induction on~$i$ that $\rad^i(M) \subseteq \soc_{n-1-i}(M)$.
This implies that $\rad^{n}(M) = 0$, and so $m \leq n$. 
Similarly, Corollary~\ref{associated graded of socle filtration is semisimple} and induction on~$j$ implies that
$\rad^{m-1-j}(M) \subseteq \soc_j(M)$ for all~$j$. 
Hence $\soc_{m-1}(M) = M$, and so $n \leq m$.
\end{proof}

\subsubsection{Locally finite categories}\label{subsec:locally finite categories}
A {\em locally finite category} is a Grothendieck category admitting
a set of generators of finite length. 

\begin{lemma}\label{locally finite categories}\leavevmode
\begin{enumerate}
\item If~$\cA$ is a locally finite category, and~$M \in \cA$, then~$M$ is compact if and only if it is Noetherian, if and only if it has finite length.
\item If~$\cA$ is a locally finite category, and~$\cA^{\fl}$ is the full subcategory of objects of finite length, 
then the inclusion $\cA^{\fl} \to \cA$ induces an equivalence $\Ind \cA^{\fl}
\to \cA$. 
\item If~$\cA^\op$ is a locally finite category, and~$\cA^{\fl}$ is the full subcategory of~$\cA$ whose objects are of finite length, 
then the inclusion $\cA^{\fl} \to \cA$ induces an equivalence $\Pro \cA^{\fl} \to \cA$.
\end{enumerate}
\end{lemma}
\begin{proof}
For parts~(1) and~(2),
observe that the locally finite category $\cA$ is locally Noetherian, hence its compact objects coincide with its Noetherian objects by~\cite[Proposition~A.1.1]{DEGlocalization}, 
and (writing~$\cA^{\mathrm{Noeth}}$ for the subcategory of Noetherian objects) we have $\Ind \cA^{\mathrm{Noeth}} \iso \cA$ by~\cite[Proposition~6.3.4]{MR2182076}.
So it suffices to prove that~$\cA^{\mathrm{Noeth}} = \cA^{\fl}$.
Since it is immediate that~$\cA^{\fl} \subseteq \cA^{\mathrm{Noeth}}$, it suffices to prove that every compact object~$X \in \cA$ has finite length,
which is immediate since~$X$ is a filtered colimit of subobjects of finite length.
Part~(3) is dual to part~(2).
\end{proof}

If $\cA$ is a locally finite category, then 
a {\em block} of $\cA$ is an equivalence class of isomorphism classes of simple objects under the equivalence relation generated by
\[
S_1 \sim S_2 \text{ if } \Ext^1_{\cA}(S_1, S_2) \ne 0 \text{ or } \Ext^1_{\cA}(S_2, S_1) \ne 0.
\]
If $\fB$ is a block of~$\cA$,
then we let $\cA_{\fB}$ denote the full subcategory of~$\cA$ containing
precisely those objects all of whose irreducible subquotients 
lie in~$\fB$.

We let $\cA_{\fB}^{\fl}$ denote the full subcategory of objects in~$\cA_{\fB}$
that are of finite length in~$\cA_{\fB}$, and note the following evident result.  

\begin{lem}\label{lem:fl equals Noetherian}
For any block~$\fB$, the category $\cA_{\fB}^{\fl}$ coincides with 
the subcategory of objects of~$\cA_\fB$ that are finite length, equivalently compact, equivalently
Noetherian, in~$\cA$. %
\end{lem}
\begin{proof}
Since $\cA_{\fB}$ is a localizing subcategory of~$\cA$, an object of $\cA_{\fB}$
is of finite length in $\cA$ if and only it is of finite length when regarded as an object
of~$\cA_{\fB}$.
The present lemma then follows from Lemma~\ref{locally finite categories}.
\end{proof}

We have the following general structure theorem for locally finite categories.

\begin{prop}\label{Gabriel decomposition}
Let~$\cA$ be a locally finite category.
Then there is a canonical direct product decomposition
\begin{equation}\label{eqn:product-decomposition-blocks} \cA \cong \prod_{\thetabar} \cA_{\fB},
\end{equation}
where the product is over all the blocks~$\fB$ of~$\cA$.
For each block~$\fB$, there exists a pseudocompact topological ring~$E_{\fB}$ 
and an equivalence between $\cA_{\fB}^{\op}$ and the category
of right-pseudocompact modules over~$E_{\fB}$. 
\end{prop}
\begin{proof}
This is a consequence of the results of~\cite[Section~IV.2]{Gabrielthesis}.
\end{proof}

The rings~$E_{\fB}$ appearing in Proposition~\ref{Gabriel decomposition}
are only well-defined up to their categories of pseudocompact modules --- that is,
up to some form of Morita equivalence. 
In our applications,
we will in fact be able to choose
the rings $E_{\fB}$ to be profinite,
so that the categories of right-pseudocompact modules will simply be the categories
$\Mod_c(E_{\fB}^\op)$ of profinite right $E_{\fB}$-modules, as we now describe.

\subsubsection{Profinite modules over profinite rings.}\label{compact modules}
Recall that
if $R$ is a topological ring whose underlying space
is profinite, and which is equipped with an $\cO$-algebra structure as
an abstract ring,
then $R$ is in fact a profinite $\cO$-algebra in the strongest possible sense,
namely we may write $R \iso \varprojlim_i R_i$, where the $R_i$ are %
$\cO$-algebras of finite cardinality, the transition maps are $\cO$-algebra homomorphisms,
and $R$ is endowed with the inverse limit topology.  
This is a consequence of the equivalence of \cite[Prop.\ 5.1.2~(b), (e)]{MR2599132}.
We will refer to such $R$ as {\em profinite topological $\cO$-algebras}
from now on.
(Note that by~\emph{loc.\ cit.}, any compact (Hausdorff) ring is necessarily
profinite, %
and so we could equally well speak of {\em compact topological
$\cO$-algebras}; but {\em profinite} seems clearer.)
The examples we care about in the main body of the paper
include the completed group ring $\cO\llbracket \Gamma\rrbracket $ of
a profinite group~$\Gamma$,
the Cayley--Hamilton algebras of Section~\ref{subsec:pseudoreps-recollections},
and the endomorphisms algebras of various profinite $\cO\llbracket G\rrbracket _{\zeta}$-modules.

Any profinite topological $\cO$-algebra $R$ 
satisfies the requirements of~\cite[\S~IV.3]{Gabrielthesis}; i.e.\ $R$ is a {\em left- and right-pseudocompact}
ring. %
We can thus consider the associated category of {\em left-pseudocompact $R$-modules},
as introduced there, and recalled in the next definition.

\begin{defn}\label{category of compact modules}
We define $\Mod_c(R)$ to be the category whose objects are separated and complete topological left $R$-modules
which admit a neighborhood basis at zero consisting of open $R$-submodules of finite $R$-colength.
The morphisms are the continuous $R$-linear maps.
(These are precisely the left-pseudocompact $R$-modules in the sense of~\cite[Section~IV.3]{Gabrielthesis}.)
\end{defn}

The category $\Mod_c(R)$ has the following standard properties.
Note that, as a consequence of Lemma~\ref{properties of compact modules}~\eqref{item: compact 3}\eqref{item: compact 4}, 
every object of~$\Mod_c(R)$ is topologically isomorphic to an inverse limit of $R$-modules of finite cardinality, and so it has a profinite topology.

\begin{lemma}\label{properties of compact modules}
Let $R$ be a profinite topological $\cO$-algebra.
\begin{enumerate}
  \item\label{item: compact 1} The category $\Mod_c(R)$ is abelian, and
  the forgetful functor 
  \begin{equation}\label{forgetful functor on compact modules}
  \Mod_c(R) \to \Mod(R)
  \end{equation}
  is faithful and exact.
  \item\label{item: compact 2} $\Mod_c(R)^\op$ is a locally finite category.
  \item\label{item: compact 3} There is an equivalence \begin{equation}\label{compact modules and Pro}
    \Pro \Mod_c(R)^{\fl} \isoto \Mod_c(R)
    \end{equation}
    through the functor that sends a cofiltered diagram of objects of $\Mod_c(R)^{\fl}$ to its limit in~$\Mod_c(R)$.
  \item\label{item: compact 4} The forgetful functor~\eqref{forgetful functor on compact modules} induces a functor
  \begin{equation}\label{forgetful functor on finite length}
  \Mod_c(R)^{\fl} \to \Mod(R)^{\fl}.
  \end{equation}
  The objects of $\Mod_c(R)^{\fl}$ carry the discrete topology, and they are finite sets.
  The functor~\eqref{forgetful functor on finite length} is fully faithful.
  \item\label{item: compact 5} Any finitely presented $R$-module has a unique topology making
  it an object of~$\Mod_c(R)$,
  referred to as its {\em canonical topology}.
  Passing to the canonical topology yields a fully faithful, exact embedding  
  \begin{equation}\label{eqn:mod-fp-to-mod-c}\Mod^{\fp}(R) \hookrightarrow \Mod_c(R).\end{equation}
  \item\label{item: compact 6} The Jacobson radical~$\Rad(R)$ of~$R$ is the intersection of all maximal closed left ideals of~$R$, and so it is closed.
  \item\label{item: compact 8} If~$\Rad(R) = 0$, then every $M \in \Mod_c(R)$ is projective, and %
  can be written as a direct product
  \[
  M = \prod_{\fm \subset R} (R/\fm)^{I_\fm}
  \]
  over maximal closed left ideals~$\fm \subset R$ {\em (}with each $I_{\fm}$ being
  an appropriately chosen indexing set{\em )}.
  \item\label{item: compact 7} If $R$ is Noetherian, then $\Rad(R)^n$ is open for all~$n \geq 0$, 
  the canonical topology on any $M \in \Mod^{\fp}(R)$ is the $\Rad(R)$-adic topology, and
  the forgetful functor~\eqref{forgetful functor on finite length} is an equivalence. 
  The equivalence~\eqref{compact modules and Pro} may thus be rewritten as an equivalence
  \begin{equation}\label{improved compact modules and Pro}
  \Pro \Mod^{\fl}(R) \isoto \Mod_c(R).
  \end{equation}
  \end{enumerate}
\end{lemma}
\begin{proof}
  Part~\eqref{item: compact 1} is proved in \cite[Section~IV.3]{Gabrielthesis}.
  Part~\eqref{item: compact 2} is~\cite[Section~IV.3, Th.\ 3]{Gabrielthesis}, and it implies part~\eqref{item: compact 3} by Lemma~\ref{locally finite categories}.
  
  We now prove part~\eqref{item: compact 4}. 
  The existence of~\eqref{forgetful functor on finite length} follows from~\cite[Lemma~3.4(1)]{MR1474172}.
  Furthermore, any $M \in \Mod_c(R)^{\fl}$ is discrete:
  this is because 
  $M$ is Artinian in $\Mod_c(R)$, and so it has a minimal open submodule~$M'$, which must be unique (as can be seen by considering the intersection of two minimal open submodules);
  and since~$M$ is separated, necessarily~$M' = 0$, which shows that~$M$ is discrete.
  Hence~\eqref{forgetful functor on finite length} is fully faithful.
  Finally, to see that~$M$ is a finite set, we can assume without loss of generality that~$M$ is simple, and then it suffices to note that there exists a continuous surjection $R \to M$,
  hence~$M$ is both compact and discrete. This concludes the proof of part~\eqref{item: compact 4}.
  
  We now prove part~\eqref{item: compact 5}.
  To prove existence of the canonical topology,
  note that $\End_{\Mod_c(R)}(R) = \End_{\Mod(R)}(R) = R^{\op}$.
  By (for example) Proposition~\ref{prop:Eilenberg-Watts},
  we thus obtain a functor $\Mod^{\fp}(R) \to \Mod_c(R)$ whose composite with the forgetful functor~\eqref{forgetful functor on compact modules}
  is the embedding $\Mod^{\fp}(R) \hookrightarrow \Mod(R).$
  This functor is fully faithful, by the ``automatic continuity'' result in~\cite[Proposition~3.5]{MR1474172}.
  The same result implies the uniqueness of the canonical topology.
  This concludes the proof of part~\eqref{item: compact 5}.

  Part~\ref{item: compact 6} follows from~\cite[Section~IV, Proposition~12]{Gabrielthesis}, or can be proved directly.
  We now turn to part~\ref{item: compact 8}, and so assume that~$\Rad(R) = 0$,
  so that $R \to \prod_{\fm} R/\fm$  is a monomorphism (the product ranging over closed maximal left ideals
  of~$R$).
  Since~$\Mod_c(R)^\op$ is locally finite (hence a Grothendieck category) the dual of Lemma~\ref{properties of semisimple objects} shows that
  $R$ is a product of simple objects of~$\Mod_c(R)$.
  Now any $M \in \Mod_c(R)$ admits a surjection $\prod_{s \in S} R \to M$ for some index set~$S$: 
  to see this, it suffices to 
  use part~\eqref{item: compact 4} to
  write~$M = \lim_{i \in I} M_i$ as a limit of finite-cardinality objects of~$\Mod_c(R)$,
  and then take the limit of the system of surjections
  \[
  \prod_{m \in M_i} R = \bigoplus_{m \in M_i} R \to M_i.
  \]
  This produces a surjection $\prod_{m \in M} R \to M$.
  Applying Lemma~\ref{properties of semisimple objects} to this surjection,
  we conclude that~$M$ is a direct product of simple objects (because so is~$R$), and a direct factor of $\prod_{m \in M} R$.
  Since~$\prod_{s \in S} R$ is indeed
  projective in $\Mod_c(R)$ for any~$S$, by e.g.~\cite[Cor.\ 1.3]{Brumerpseudocompact}, this concludes the proof of~\eqref{item: compact 8}.

  The first two statements of part~\eqref{item: compact 7} are~\cite[Corollary~3.14]{MR1474172} (note that~$M$ is a Noetherian object of~$\Mod(R)$, hence of~$\Mod_c(R)$).
  Finally, to see that~\eqref{forgetful functor on finite length} is an equivalence, 
  it suffices to note that, when~$R$ is Noetherian, the forgetful functor~\eqref{forgetful functor on compact modules} actually induces an equivalence on the full subcategories of
  Noetherian objects, by~\cite[Proposition~3.21]{MR1474172}.
\end{proof}

There are the usual tensor product functors
\begin{gather*}
\text{--}\otimes_R\text{--} : \Mod(R^\op) \times \Mod(R) \to \Mod(\cO)\\
\text{--} \cotimes_R\text{--} : \Mod_c(R^\op) \times \Mod_c(R) \to \Mod_c(\cO).
\end{gather*}
From the perspective of the equivalence~\eqref{compact modules and Pro},
the completed tensor product $\text{--}\cotimes_R\text{--}$ 
can be interpreted as the Pro-extension of the usual tensor product,
when we restrict it to a functor
$\text{--}\otimes_R\text{--}: \Mod_c(R^\op)^{\fl}\times \Mod_c(R)^{\fl} \to \Mod^{\fl}(\cO).$
If $M$ and $N$ are objects of $\Mod_c(R^{\op})$ and $\Mod_c(R)$ respectively,
then there is a natural map of $\cO$-modules
\begin{equation}
\label{eqn:tensor to completed tensor}
M\otimes_R N \to M\cotimes_R N,
\end{equation}
where the source is the usual tensor product, formed without regard to
the topologies on $M$ and~$N$, and the target is the completed tensor product,
regarded as an $\cO$-module by forgetting its topology.
The image of~\eqref{eqn:tensor to completed tensor} is dense in the defining topology on $M \cotimes_R N$.
Looking ahead, from the perspective of Lemma~\ref{lem:EW-gives-completed-tensor-product}, if~$R$ is furthermore Noetherian, the functor~$M \cotimes_R -$ is associated to the 
complete right $R$-module~$M$ in~$\Mod_c(\cO)$, by Lemma~\ref{compact modules are compact}.

If~$R$ is a profinite $\cO$-algebra, the equivalence~\eqref{compact modules and Pro} shows that~$\Mod_c(R)$
is the Pro-category of an abelian category, and so it has enough projectives (by Lemma~\ref{lem:quasi-projectives-in-Pro-cat}; alternatively, this fact was established in the proof of 
Lemma \ref{properties of compact modules}~\eqref{item: compact 8}).
We write $\widehat{\Tor}_i^R$ for the left derived functors of $\text{--} \cotimes_R\text{--}$.
\begin{defn}\label{defn:topologically-flat}
  We say that a module $M \in \Mod_c(R^\op)$, respectively $M \in \Mod_c(R)$, is \emph{topologically flat} if 
$M\cotimes_R \text{--}$ (resp.\ $\text{--}\cotimes_R M$) is an exact functor.
\end{defn}

\begin{lemma}\label{topologically flat implies flat over Noetherian}
Let~$R$ be a Noetherian profinite $\cO$-algebra, and let $M \in \Mod_c(R^\op)$.
Then $M$ is topologically flat if and only if $M$ is a flat $R$-module {\em (}in the usual sense that $M \otimes_R (\text{--})$ is exact on~$\Mod(R)${\em )}.
\end{lemma}
\begin{proof}
Since Definition~\ref{defn:topologically-flat} is a specialization of Definition~\ref{topologically flat object}, 
this lemma is a consequence of Lemma~\ref{topologically flat equivalent to flat}.
\end{proof}

Since a right exact functor on an abelian category with enough projectives is exact
if and only if its left derived functors vanish, we see that an object $M$
of~$\Mod_c(R^{\op})$ (resp.\ of~$\Mod_c(R)$) is topologically flat if and only if
$\widehat{\Tor}_i^R(M, \text{--}) = 0$ for all~$i > 0$, respectively $\widehat{\Tor}_i^R(\text{--}, M) = 0$ for all~$i > 0$.
From the very construction of left derived functors,
we see that projective objects of $\Mod_c(R^\op)$ and~$\Mod_c(R)$ are topologically flat;
we will see in Lemma~\ref{lem: topologically flat implies projective}
below that the converse also holds.
The following lemma shows that the (derived) tensor product with finitely generated $R$-modules 
coincides with the (derived) completed tensor product.

\begin{lemma}
\label{lem:automatic tensor completeness}
Let~$R$ be a Noetherian profinite $\cO$-algebra.
If $M$ is an object of $\Mod_c(R^{\op})$ and $N$ is an object of~$\Mod^{\fp}(R),$
then the natural map~{\em \eqref{eqn:tensor to completed tensor}} induces isomorphisms $\Tor_i^R(M, N) \to \widehat{\Tor}_i^R(M, N)$. %
\end{lemma}
\begin{proof}
The statement for~$i = 0$ is immediate from the fact that the restriction of $M \cotimes_R \text{--}$ to $\Mod^{\fp}(R)$ is naturally isomorphic to $M \otimes_R -$
(see Corollary~\ref{cor:restricting-completed-tensor} for a generalization).
In turn, this implies the statement for all~$i$, since we can compute the left derived functors via a resolution of~$N$ by finite free $R$-modules.
\end{proof}

\begin{remark}
\label{rem:automatic tensor completeness}
In the context of Lemma~\ref{lem:automatic tensor completeness}, we will often regard $M\otimes_R N$ as
a compact $R$-module, by identifying it with~$M\cotimes_R N$.
Note that the case~$i = 0$ of Lemma~\ref{lem:automatic tensor completeness} is true without the Noetherian assumption on~$R$, 
and in fact applies in the more general context of pseudocompact rings and modules (compare~\cite[Lemma~2.1]{Brumerpseudocompact}).
\end{remark}

\begin{remark}
\label{rem:computing the cosocle}
Let~$R$ be a profinite $\cO$-algebra.
Since $\Rad(R)$ is closed in~$R$ (by Lemma~\ref{properties of compact modules}~\eqref{item: compact 6}),
it inherits a profinite topology making it a subobject of $R$ in~$\Mod_c(R)$.
Thus we may consider the quotient $R/\Rad(R)$, 
and the associated completed tensor product functor
$(R/\Rad(R)) \cotimes_R -$ on~$\Mod_c(R)$.  
Evidently, $R/\Rad(R)$ is again a profinite $\cO$-algebra,
and this functor takes values in $\Mod_c\bigl( R/\Rad(R)\bigr).$
Lemma~\ref{properties of compact modules}~\eqref{item: compact 8}
shows that if $M$ is an object of $\Mod_c(R)$, then $(R/\Rad(R)) \cotimes_R M$
is a product of simple discrete $R$-modules.
In fact, $(R/\Rad(R)) \cotimes_R M$
may be characterized as the maximal quotient of $M$ in $\Mod_c(R)$
which is a product of simple discrete
$R$-modules.
To see this, note that
if~$N$ is a product of simple discrete $R$-modules, and $M \to N$ is a surjection in~$\Mod_c(R)$, then 
$M \to N$ factors through $M/\Rad(R) M$. 
Since $N \in \Mod_c(R)$, it actually factors through the quotient of~$M$
by the closure of $\Rad(R) M$ in~$M$. 
Using the fact that~\eqref{eqn:tensor to completed tensor} has dense image (applied to $\rad(R) \otimes_R M \to \rad(R) \cotimes_R M$) 
we see that this quotient is $(R/ \Rad(R)) \cotimes_R M$.
\end{remark}

We next %
show that topological flatness is equivalent to being projective in~$\Mod_c(R)$.
(This is standard, see for example~\cite[Lemma~2.4]{MR4350140}, but we provide a proof for convenience.)

\begin{lemma}\label{lem: topologically flat implies projective}
Let~$R$ be a profinite $\cO$-algebra, and let~$M \in \Mod_c(R)$.
Then~$M$ is topologically flat over~$R$ if and only if it is projective in~$\Mod_c(R)$.
\end{lemma}
\begin{proof}
If~$M$ is projective, then it is topologically flat, since $\widehat \Tor_i^R(\text{--}, M) = 0$ for all~$i > 0$.
We now prove the converse.
By Remark~\ref{rem:computing the cosocle},
we know that~$(R/\Rad(R)) \cotimes_R M$ can be written as a direct product
\[
(R/\Rad(R)) \cotimes_R M = \prod_\fm (R/\fm)^{I_\fm}
\]
over maximal closed left ideals~$\fm \subset R$, for some appropriately chosen indexing sets~$I_\fm$.
Now choose a projective envelope $R_\fm \to R/\fm$ of~$R/\fm$ in~$\Mod_c(R)$ (whose existence is guaranteed by the fact that $\Mod_c(R)^{\op}$ is a Grothendieck category,
which has injective envelopes, by~\cite[Thm.\ II.6.2]{Gabrielthesis}).
By definition, $R_\fm$ has a unique discrete simple quotient (namely $R/\fm$), %
and so Remark~\ref{rem:computing the cosocle} shows that $(R/\Rad(R)) \cotimes_R R_\fm
\iso R/\fm$.
So we obtain a morphism
\begin{equation}\label{isomorphism by Nakayama's lemma}
\prod_\fm R_\fm^{I_\fm} \to M
\end{equation}
which induces an isomorphism after applying $R/\Rad R \cotimes_R (\text{--})$.
We now apply Nakayama's lemma for~$\Mod_c(R)$,
which asserts that if~$N$ is an object of $\Mod_c(R)$
for which $\bigl(R/\Rad(R)\bigr) \cotimes_R N = 0$, then~$N = 0$.
(See for example~\cite[Lemma~3.23]{MR1474172}.)
Together with the assumption that~$M$ is topologically flat, this implies that~\eqref{isomorphism by Nakayama's lemma} is an isomorphism, and so~$M$ is projective, as desired.
\end{proof}

\subsubsection{Tensoring objects of abelian categories by modules}
\label{subsubsec:tensor products in abelian categories}
In preparation for the statements of Morita theory, we now discuss  a general formalism for tensor products in abelian categories.
\begin{defn}\label{defn:R-module-objects}
  Let~$R$ be a ring, and let~$\cA$ be an abelian category. Then a \emph{right $R$-module  in~$\cA$} is a pair consisting of an object $M$ of~$\cA$, and a ring homomorphism $R^{\op} \to \End_{\cA}(M)$. The category of right $R$-modules in $\cA$ is the category whose objects are right $R$-modules  in~$\cA$, and whose morphisms are the morphisms of underlying of objects of~$\cA$ which are compatible with the homomorphisms from ~$R^{\op}$.
\end{defn}
Given a functor $F:\Mod^{\fp}(R)\to\cA$, we see that $F(R)$ is naturally a right $R$-module in~$\cA$. We have the following standard result.
\begin{prop}
  \label{prop:Eilenberg-Watts}
Let~$R$ be a ring, and let~$\cA$ be an abelian category.
  \begin{enumerate}
  \item\label{item:56} The functor $F\mapsto F(R)$ is an equivalence of categories between the category of right exact functors
$\Mod^{\fp}(R)\to\cA$ and the category of right $R$-modules ~$M$ in~$\cA$;
we denote the functor $F$ associated to $M$ by~$M\otimes_R\text{--}$. 
  \item\label{item:54} If~$\cA$ is cocomplete, then the functor $F\mapsto F(R)$ is an equivalence of categories between the category of colimit-preserving functors %
$\Mod(R)\to\cA$ and the category of right $R$-modules~$M$ in~$\cA$; we denote
the functor $F$ associated to $M$ by~$M\otimes_R\text{--}$. Furthermore, the functor~$F$ is then
left adjoint to the functor $\Hom_{\cA}(M,\text{--}\, ): \cA \to \Mod(R)$.
    \item\label{item:55}  If~$\cA$ is complete, then the functor $F\mapsto F(R)$ is an equivalence of categories between the category of right exact, cofiltered limit-preserving functors
$\Pro\Mod^{\fp}(R)\to\cA$ and the category of right $R$-modules in~$\cA$;
we denote the functor $F$ associated to $M$ by $M\cotimes_R\text{--}$.
    \item The functor $F\mapsto F(R)$ is an equivalence of categories between the category of right exact, cofiltered limit-preserving functors
$\Pro\Mod^{\fp}(R)\to \Pro\cA$ 
    sending~$R$ to an object of~$\cA$, and the category of right $R$-modules in~$\cA$;
we denote the functor $F$ associated to $M$ by $M\cotimes_R\text{--}$.
  \end{enumerate}
\end{prop}
\begin{proof}If~$R$ is commutative, then~\eqref{item:56} is \cite[\href{https://stacks.math.columbia.edu/tag/0GNN}{Tag 0GNN}]{stacks-project}, and the same proof works for general~$R$.
  Part~\eqref{item:54} then follows by passage to filtered colimits, see~\cite[\href{https://stacks.math.columbia.edu/tag/0GNQ}{Tag 0GNQ}]{stacks-project} and~\cite[\href{https://stacks.math.columbia.edu/tag/0GNS}{Tag 0GNS}]{stacks-project}.

  Alternatively, \eqref{item:54} is the special case $k=\Z$ of~\cite[Thm.\ 3.1]{MR3507176}. Part~\eqref{item:56} can then be deduced by applying~ \eqref{item:54} with~$\cA$ replaced by~$\Ind\cA$, noting that because $M$ is compact in $\Ind \cA$ (being an object of~$\cA$),
the functor $\Hom_{\Ind \cA}(M,\text{--}):\Ind\cA\to\Mod(R)$ is filtered colimit-preserving,
from which one deduces that its left adjoint~$F$
preserves compact objects.  (Cf.\ Lemma~\ref{lem: adjoint continuous} below for the analogous
statement in the context of functors between stable $\infty$-categories.)

Finally~\eqref{item:55} follows from~\eqref{item:56} by passing to cofiltered limits, and part~(4) follows from parts~\eqref{item:55} and~\eqref{item:56}. 
\end{proof}%

\begin{rem}
  \label{rem:algebra-version-EW}The equivalences of categories in Proposition~\ref{prop:Eilenberg-Watts} give the usual bifunctoriality of the tensor product, and even the trifunctoriality of the tensor-$\Hom$ adjunction (see for example~\cite[Chapter VI, Thm.\ 3.1]{MR202787} for this last point). In addition, if~$S$ is a commutative ring, and $R$ is an $S$-algebra, then there is an obvious $S$-linear variant of Proposition~\ref{prop:Eilenberg-Watts}; see~\cite[Cor.\ 2.3]{MR3507176}.
\end{rem}
\begin{rem}
  \label{rem:Hom-version-of-EW}
  Proposition~\ref{prop:Eilenberg-Watts} also applies to left $R$-modules~$M$ in~$\cA$, which (by definition) are right $R^{\op}$-modules.
  As a consequence, if~$M$ is a left $R$-module in~$\cA$, we will sometimes write $\text{--}\otimes_R M: \Mod^{\fp}(R^{\op}) \to \cA$ for the right exact functor previously denoted $M \otimes_{R^{\op}} -$.
  If~$\cA$ is cocomplete, it extends to a functor $\Mod(R^{\op})\to\cA$ which is left adjoint to $\Hom_{\cA}(M, \text{--}) : \cA \to \Mod(R^{\op})$.

  On the other hand, a left $R$-module~$M$ in~$\cA$ can also be regarded as a right $R$-module in~$\cA^{\op}$, and so it defines a (covariant) right exact functor $\Mod^{\fp}(R) \to \cA^{\op}$.
  We will write $\Hom_{R}(\text{--}, M) : \Mod^{\fp}(R) \to \cA$ for the corresponding \emph{contravariant} right exact functor.
  If~$\cA$ is complete, so that Proposition~\ref{prop:Eilenberg-Watts}~\eqref{item:54} applies to~$\cA^{\op}$, then the functor $\Hom_R(\text{--}, M)$ extends to a contravariant functor $\Mod(R)\to\cA$ (sending colimits to limits) 
  and there is a bifunctorial isomorphism
  \[
  \Hom_{\cA}(X, \Hom_R(N, M)) \cong \Hom_R(N, \Hom_\cA(X, M))
  \]
  for all~$X \in \cA$ and all~$N\in\Mod(R)$.
  
  In applications, we will have a finitely generated two-sided ideal~$J$ of~$R$, and a left $R$-module~$M$ in~$\cA$, and we will write $M[J] \coloneqq  \Hom_R(R/J, M)$.
  Then $M[J]$ is naturally a subobject of~$M$ in the category of left $R$-modules in~$\cA$, and the $R$-action on $M[J]$ factors through~$R/J$.
  The right exact functor $\otimes_R M[J]$ is therefore naturally isomorphic to the composition 
  \[
  \Mod^{\fp}(R^{\op}) \xrightarrow{\otimes_R R/J} \Mod^{\fp}((R/J)^{\op}) \xrightarrow{\otimes_{R/J} M[J]}  \cA. 
  \]
\end{rem}

\begin{lemma}
\label{lem:tensor and projective limits}
Let~$R$ be a ring, and let $\cA$ be a complete abelian category with exact cofiltered limits.
Let~$M_i$ be a cofiltered system of right $R$-modules in~$\cA$, and let $M \coloneqq  \varprojlim_i M_i$.
Then the natural map $(\varprojlim_i M_i \otimes_R \text{--}) \to \varprojlim_i(M_i \otimes_R \text{--})$ is an isomorphism of functors $\Mod^{\fp}(R)\to \cA$.
\end{lemma}
\begin{proof}
Since both sides send~$R$ to~$\varprojlim_i M_i$, it suffices (by Proposition~\ref{prop:Eilenberg-Watts}) to prove that $\varprojlim_i(M_i \otimes_R \text{--})$ is a right-exact functor on~$\Mod^{\fp}(R)$, 
which is a consequence of the exactness of cofiltered limits in~$\cA$.
\end{proof}

We will also need the following variant of the above constructions, which produces a formalism of ``completed tensor products'' in complete abelian categories.

\begin{defn}%
  \label{defn:complete-object-abelian-category}
Let~$R$ be a Noetherian profinite  $\cO$-algebra, and let~$\cA$ be an abelian category. 
  Then we say that a right $R$-module~$M$ in ~$\cA$ is \emph{complete} if the natural morphisms $M\to M\otimes_R\bigl(R/\Rad(R)^n\bigr)$ give rise to an isomorphism $M\to\varprojlim_nM\otimes_{R} \bigl(R/\Rad(R)^n\bigr)$ 
  (i.e.\ this limit exists in~$\cA$, and is given by~$M$).
\end{defn}

\begin{lem}
  \label{lem:EW-gives-completed-tensor-product}Let $R$ be a Noetherian profinite topological $\cO$-algebra, and let~$\cA$ be a complete abelian category.
Then the assignment $F\mapsto F(R)$ is an equivalence between the category of right exact and cofiltered limit-preserving functors 
$\Mod_c(R) \to \cA$, and the category of complete right $R$-modules~$M$ in~$\cA$; we denote the functor~$F$ associated to~$M$ by~$M \cotimes_R \text{--}$.
\end{lem}
\begin{proof}
Suppose that~$F:\Mod_c(R) \to \cA$ is right exact and cofiltered limit-preserving, and write $M\coloneqq F(R)$. 
Then the restriction of~$F$ to~$\Mod^{\fp}(R)$ is given by $M\otimes_R\text{--}$, and since 
by Lemma~\ref{properties of compact modules}~\eqref{item: compact 7} we have $R\isoto \varprojlim_{n} \bigl(R/\Rad(R)^n\bigr)$ in $\Mod_c(R)$, we see that~$M$ is 
   necessarily a complete right $R$-module.
This shows that $F \mapsto F(R)$ defines a functor between the categories in the statement of the lemma.

We now construct a quasi-inverse to $F \mapsto F(M)$.
If $M\in\cA$ is a complete right $R$-module, then we have the right exact and cofiltered limit-preserving functor $M\cotimes_R\text{--}:\Pro\Mod^{\fp}(R)\to \cA$ 
of Proposition~\ref{prop:Eilenberg-Watts}~\eqref{item:55}. %
Restricting this functor to $\Pro\Mod^{\fl}(R)$, and recalling ~\eqref{improved compact modules and Pro}, we obtain a right exact and cofiltered limit-preserving functor $\Mod_c(R)\isoto\Pro\Mod^{\fl}(R)\to\cA$.

Given~$M$, the functor $M \cotimes_R -$ takes~$R$ to $\varprojlim_nM\otimes_R\bigl(R/\Rad(R)^n\bigr)$, which coincides with~$M$, since~$M$ is assumed to be complete.
On the other hand, given~$F$, the functors~$F$ and $F(M) \cotimes_R - : \Mod_c(R) \to \cA$ preserve cofiltered limits, and are isomorphic to~$F(M) \otimes_R -$ 
after restriction to $\Mod^{\fl}(R)$;
they are therefore naturally isomorphic.
This concludes the proof that $F \mapsto F(M)$ is an equivalence, with quasi-inverse $M \mapsto M \cotimes_R -$.
\end{proof}

\begin{cor}
  \label{cor:restricting-completed-tensor}Let $R$ be a Noetherian profinite $\cO$-algebra. 
  Let~$\cA$ be a complete abelian category, and let $M\in\cA$ be a complete right $R$-module.
Then the restriction of $M\cotimes_R\text{--}:\Mod_c(R)\to\cA$ to~$\Mod^{\fp}(R)$ is $M\otimes_R\text{--}$.
\end{cor}
\begin{proof}
  This is immediate by construction (or by the defining properties of the tensor products).
\end{proof}

Next, we record a compatibility between our tensor product functors, and right exact functors between abelian categories.

\begin{lem}%
\label{lem:abelian-composition-right-exact-tensor}\leavevmode
\begin{enumerate}  
 \item\label{item:49} If~$\cA, \cB$ are abelian categories, $R$ is a ring, $F:\cA\to\cB$ is right exact, and~$M$ is an $R$-module in~$\cA$, then there is a natural isomorphism of functors %
$\Mod^{\fp}(R)\to\cB$ \begin{equation}\label{eqn:tensor-product-composed-with-F} (F(M)\otimes_R\text{--})\to F(M\otimes_R\text{--}).\end{equation} %
\item\label{item:99}  If~$\cA, \cB$ are complete abelian categories, $R$ is a Noetherian profinite topological $\cO$-algebra, $F: \cA \to \cB$ is right exact and cofiltered limit-preserving,
and~$M$ is a complete $R$-module in~$\cA$,
then there is a natural isomorphism of functors 
$\Mod_c(R)\to\cB$ \begin{equation}\label{eqn:completed-tensor-product-composed-with-F} (F(M)\cotimes_R\text{--})\to F(M\cotimes_R\text{--}).\end{equation}
\end{enumerate}
\end{lem}
\begin{proof}
  The existence of the natural morphism of functors~\eqref{eqn:tensor-product-composed-with-F} follows from adjunction (i.e.\ from the deduction of part ~\eqref{item:56} of Proposition~\ref{prop:Eilenberg-Watts} from part~\eqref{item:54})  Since the functor $F(M\otimes_R\text{--})$ is right exact and takes~$R$ to~$F(M)$, this morphism  is an isomorphism by Proposition~\ref{prop:Eilenberg-Watts}~\eqref{item:56}.
  
  For part~\eqref{item:99}, since $\Mod_c(R) \isom \Pro \Mod^{\fl}(R)$ (by Lemma~\ref{properties of compact modules}~\eqref{item: compact 7}), and both sides of~\eqref{eqn:completed-tensor-product-composed-with-F} preserve cofiltered limits, it suffices
  to construct the isomorphism~\eqref{eqn:completed-tensor-product-composed-with-F} after restricting to~$\Mod^{\fl}(R)$, which is a consequence of part~\eqref{item:49} and Corollary~\ref{cor:restricting-completed-tensor}. 
\end{proof}

Now let~$\cA$ be a complete abelian category, and let~$R$ be a Noetherian profinite $\cO$-algebra.
If~$X \in \Mod_c(R)$, then Lemma~\ref{lem:EW-gives-completed-tensor-product} shows that the completed tensor product with~$X$ gives rise to a functor
\[
\text{--}\cotimes_R X : (\text{complete right $R$-modules in~$\cA$}) \to \cA.
\]
More precisely, given a morphism $M_1 \to M_2$ of complete right $R$-modules in~$\cA$, the morphism $M_1 \cotimes_R X \to M_2 \cotimes_R X$ is the evaluation at~$X$
of the natural transformation $M_1 \cotimes_R - \to M_2 \cotimes_R -$ arising from Lemma~\ref{lem:EW-gives-completed-tensor-product}.

\begin{lem}\label{right exactness in the module variable}
Let~$\cA$ be a complete abelian category in which products are exact, and let~$R$ be a Noetherian profinite $\cO$-algebra.
  Then for any object $X$ of~$\Mod_c(R)$,
and for any sequence of morphisms~$M_1 \to M_2 \to M_3 \to 0$ of complete right $R$-modules in~$\cA$ which is exact when viewed as a sequence in~$\cA$,
the resulting sequence $M_1 \cotimes_R X \to M_2 \cotimes_R X \to M_3 \cotimes_R X \to 0$ is exact; 
or, more succinctly,
for any~$X \in \Mod_c(R)$,
the functor
$\text{--}\cotimes_R X$ {\em(}regarded as a functor from the category of
complete right $R$-modules in~$\cA$ to~$\cA${\em )} is right exact.
\end{lem}
\begin{proof}%
  Note firstly that the claim is immediate if~$X = R$.  For general~$X$, we may choose a presentation $\prod_{i \in I} R \to \prod_{j \in J} R \to X \to 0$, so that for any right exact sequence $M_1 \to M_2 \to M_3 \to 0$ of $R$-modules, we have a commutative diagram
  \[
    \begin{tikzcd}
      \prod_{i \in I}M_1 \arrow[r] \arrow[d] & \prod_{i \in I} M_2 \arrow[r] \arrow[d] & \prod_{i \in I} M_3 \arrow[r] \arrow[d]& 0\\
      \prod_{j \in J}M_1 \arrow[r] \arrow[d]& \prod_{j \in J} M_2 \arrow[r] \arrow[d]& \prod_{j \in J} M_3 \arrow[r] \arrow[d]& 0\\
      M_1 \arrow[r] \arrow[d]\cotimes_R X & M_2 \arrow[r] \arrow[d]\cotimes_R X & M_3 \cotimes_R X \arrow[r] \arrow[d]& 0\\
      0 & 0 & 0 &
    \end{tikzcd}
  \]
  with all rows and columns exact, except possibly the last row, which is therefore also exact.
\end{proof}

\begin{lem}\label{compact modules are compact}
  Let~$R$ be a Noetherian profinite $\cO$-algebra.
  Let~$M$ be an object of $\Mod_c(R^{\op})$.
  Then~$M$ %
  is a complete right $R$-module in $\Mod_c(\cO)$.
  The functor $M \cotimes_R - : \Mod_c(R) \to \Mod_c(\cO)$ obtained from
  Lemma~{\em \ref{lem:EW-gives-completed-tensor-product}}
  coincides with the usual completed tensor product with~$M$.
\end{lem}
\begin{proof}  
  As a consequence of Lemma~\ref{properties of compact modules}~\eqref{item: compact 3}\eqref{item: compact 4}, the canonical topology on~$M$ is profinite, and so~$M$ is 
  a right $R$-module in $\Mod_c(\cO)$.
  We now show that~$M$ is complete, i.e.\ that the natural map $M \to \varprojlim_n M \otimes_R R/\rad(R)^n$ is an isomorphism. 
  This map is part of a natural transformation of functors $\Mod_c(R^{\op}) \to \Mod_c(\cO)$, which induces the identity after evaluating at~$M = R$, by Lemma~\ref{properties of compact modules}~\eqref{item: compact 7}. 
  So it suffices to prove that both functors preserve cofiltered limits, and to do so, it suffices to prove that 
  $M \mapsto M \otimes_R R/\rad(R)^n$ preserves cofiltered limits.
  Since~$R/\rad(R)^n$ is finitely presented, this is a consequence of Lemma~\ref{lem:tensor and projective limits}. %
  This concludes the proof that~$M$ is complete.
  The last statement of the corollary now follows, because 
  the functor $M \cotimes_R -$ and the usual completed tensor product are right exact, cofiltered limit-preserving functors mapping~$R$ to~$M$.
\end{proof}

\begin{rem}
  \label{rem:completed-tensor-product}
  Let~$R$ be a Noetherian profinite $\cO$-algebra, and let~$\cA$ be a complete abelian category.
If~$M$ is a right $R$-module in~$\cA$, we have associated to~$M$ a pro-extended functor $M \cotimes_R -: \Pro \Mod^{\fp}(R) \to \cA$ 
(Proposition~\ref{prop:Eilenberg-Watts}~\eqref{item:55}).
It is right exact, and cofiltered limit-preserving.
Since $\Mod_c(R) = \Pro \Mod^{\fl}(R)$, there is an exact and cofiltered limit-preserving inclusion functor \[g: \Mod_c(R) \to \Pro \Mod^{\fp}(R),\]
which sends~$R$ to $\quoteslim{n} R/\Rad(R)^n$.
Restricting~$M \cotimes_R -$ through~$g$ gives rise to a right exact, cofiltered-limit preserving functor $\Mod_c(R) \to \cA$.
It sends~$R$ to the right $R$-module $\varprojlim_{n} (M \otimes_R R/\Rad(R)^n)$.
If~$M$ is complete, we thus obtain the functor associated to~$M$ in Lemma~\ref{lem:EW-gives-completed-tensor-product}.
This explains why we have denoted both functors by $M \cotimes_R -$.

The above discussion can also be understood in terms of the left adjoint~$f$ to~$g$,
which is pro-extended from the inclusion $\Mod^{\fp}(R) \to \Mod_c(R)$ described in Lemma~\ref{properties of compact modules}~\eqref{item: compact 5}:
in fact, if $\quoteslim{i \in I}M_i$ is a cofiltered system of finitely presented $R$-modules, and $N = \varprojlim_{j \in J} N_j$ is a compact $R$-module written as a cofiltered limit of finite length $R$-modules, then
\begin{multline*}
\Hom_{\Mod_c(R)}(\varprojlim_{i} M_i, N) = \varprojlim_j\Hom_{\Mod_c(R)}(\varprojlim_i M_i, N_j) = \varprojlim_j \varinjlim_i \Hom_{\Mod_c(R)}(M_i, N_j)\\
= \varprojlim_j \varinjlim_i \Hom_{\Mod^{\fp}(R)}(M_i, N_j) = \Hom_{\Pro \Mod^{\fp}(R)}(\quoteslim{i}M_i, \quoteslim{j}N_j)\\ 
= \Hom_{\Pro \Mod^{\fp}(R)}(\quoteslim{i}M_i, g(N)).
\end{multline*}
In more detail, the second equality is because objects of~$\Mod^{\fl}(R)$ are compact in $\Mod_c(R)^{\op} = \Ind (\Mod^{\fl}(R)^{\op})$, and the third equality is because finite length objects are finitely presented, and
$\Mod^{\fp}(R) \to \Mod_c(R)$ is fully faithful.
It follows immediately that the counit $fg \to \id_{\Mod_c(R)}$ is an isomorphism.
Furthermore, by Lemma~\ref{lem:left-adjoint-to-Pro-functor-abelian}, the functor~$f$ is exact. 
We are therefore in the situation of~\cite[Proposition~5, \S III.2]{Gabrielthesis}, and we conclude that~$f$ is a Serre quotient functor.

The condition that~$M \in \cA$ be ``complete'' can therefore be understood as an explicit characterization
 of the condition that the functor $M \cotimes_R - : \Pro \Mod^{\fp}(R) \to \cA$ 
factors through~$f$
(and, indeed, this characterization then
shows that $f$ is the quotient of $\Pro \Mod^{\fp}(R)$
by the Serre subcategory generated by
the kernel of the canonical morphism $R \to \quoteslim{n} R/\Rad(R)^n$).
If this condition holds, then we can compute the resulting functor $\Mod_c(R) \to \cA$ by restricting through~$g$, since $fg \isom \id_{\Mod_c(R)}$.
\end{rem}

\subsubsection{Morita theory for abelian categories.}
\label{subsubsec:theta Morita}
Now let $\cC$ be a cocomplete abelian category,
and let $P$ be an object of~$\cC$.
Let~$E \coloneqq \End_\cC(P)$.
Then $E$ is a right $E^{\op}$-module in~$\cA$,
and so the functor
\begin{equation}\label{equivalence without topologies}
\Hom_\cC(P, \text{--}): \cC \to \Mod(E^\op)
\end{equation}
admits a left adjoint $\text{--}\otimes_E P$, by Proposition~\ref{prop:Eilenberg-Watts}~\eqref{item:54}. 

Assume now that~$P$ is projective.
Recall that the right orthogonal to~$P$, which we denote by~$\cT$, is defined to be
the full subcategory
of $\cC$ consisting of objects $X$ for which $\Hom_{\cC}(P,X) = 0.$
The object~$P$ is then a generator of~$\cC$ if and only if $\cT = 0$. 
We have the following result.

\begin{lemma}\label{properties of projective objects}
Let~$\cC$ be a cocomplete abelian category, let $P \in \cC$ be projective, and let~$\cT$ be the right orthogonal to~$P$. Then:
\begin{enumerate}
\item
 $\cT$ is a Serre subcategory of~$\cC$.
\item If, for any object $X$ of~$\cC$, we let $\overline{X}$ denote the
image of $X$ in~$\cQ \coloneqq  \cC/\cT$, then $\Hom_{\cC}(P,X) \iso \Hom_{\cQ}(\overline{P},\overline{X}).$
\item
The image~$\overline{P}$ of $P$ in $\cC/\cT$ is a projective generator 
of~$\cQ$.
\end{enumerate}
\end{lemma}
\begin{proof}
Part~(1) is because~$\cT$ is the kernel of the exact functor $\Hom_\cC(P, \text{--})$.
Part~(2) can be proved in the same way as~\cite[Lem.\ III.2.1(c)]{Gabrielthesis}.
We now prove part~(3), starting with the claim that~$\lbar P$ is projective in~$\cQ$.
If~$\lbar X \to \lbar Y$ is a surjection in~$\cQ$, then by~\cite[Cor.\ III.1.1]{Gabrielthesis} we can replace it by an isomorphic surjection, and 
assume that it is represented by a surjection $X \to Y$ in~$\cC$.
Since~$P$ is projective in~$\cC$, this shows that 
$\Hom_{\cC}(P,X) \to \Hom_{\cC}(P,Y)$
is surjective.  
Now (2) shows that $\Hom_{\cQ}(\overline{P}, \overline{X}) \to
\Hom_{\cQ}(\overline{P},\overline{Y})$ is surjective.
Hence~$\lbar P$ is projective in~$\cQ$.
Now it remains to prove that the orthogonal to~$\lbar P$ in~$\cQ$ is zero, which is a direct consequence of part~(2).
\end{proof}

In the context of Lemma~\ref{properties of projective objects}, 
assume furthermore that~$\cC$ is a Grothendieck category, 
and that~$P$ is a compact projective object of~$\cC$ (or equivalently, $\Hom_\cC(P, \text{--})$ is exact and commutes with direct sums).
If $E \coloneqq  \End_{\cC}(P),$ then from Lemma~\ref{properties of projective objects}~(2) 
we see that also $E \coloneqq  \End_{\cQ}(\overline{P})$,
and then Lemma~\ref{properties of projective objects}~(3) together with~\cite[Cor.\ V.1.1]{Gabrielthesis} implies that 
\[
\Hom_{\cQ}(\lbar P,\text{--}) : \cQ \to \Mod(E^\op)
\]
is an equivalence.
Its left adjoint $\text{--}\otimes_E \overline{P}$, as constructed in Proposition~\ref{prop:Eilenberg-Watts}~\eqref{item:54}, %
is then a quasi-inverse to this equivalence.
In particular, if $P$ is a compact projective generator 
of~$\cC$,
then $\Hom_{\cC}(P,\text{--})$ induces an equivalence $\cC \iso \Mod(E^{\op})$ with quasi-inverse $\text{--}\otimes_E P$.

Some of the Grothendieck categories of interest to us do not admit projective generators
(compact or not),
and so the preceding
discussion 
will not apply.  However,   
Lemma~\ref{locally finite categories} shows that, if $\cA$ is locally finite,
then its opposite category is naturally the Pro-category of an abelian category, and so in particular does contain
projective objects.  
We now describe a variant of the preceding results which applies in this context. 

Let~$\cC$ be the opposite category to a locally finite category, and let~$P$ be
a projective generator of~$\cC$.
Then, following~\cite[Section~IV.4]{Gabrielthesis}, the ring $E \coloneqq  \End_\cC(P)$ can be endowed with the structure of a right-pseudocompact topological ring, 
and the functor $\Hom_\cC(P, \text{--})$ factors through $\Mod_c(E^{\op})$. 
We make the additional assumption that~$P$ has finite cosocle, and that~$E$ is a Noetherian profinite $\cO$-algebra; this will be true in all cases of interest in this paper.
By~\cite[Thm.\ IV.4.4]{Gabrielthesis}, 
the functor 
\begin{equation}\label{equivalence from Gabriel's thesis}
\Hom_{\cC}(P, \text{--}) : \cC \to \Mod_c(E^\op)
\end{equation}
is then an equivalence: this is in fact how Proposition~\ref{Gabriel decomposition} is proved.

\begin{prop}\label{compact Morita theory}
Under the assumptions in the previous paragraph,
$P$ is a complete left $E$-module in~$\cC$, and the functor $\text{--}\cotimes_E P : \Mod_c(E^{\op}) \to \cC$ {\em (}as defined in Lemma~{\em \ref{lem:EW-gives-completed-tensor-product})}
is a quasi-inverse to~\eqref{equivalence from Gabriel's thesis}.
\end{prop}
\begin{proof}
We begin by proving that~$P$ is complete, i.e.\ $P \isoto \varprojlim_n P/\Rad(E)^n P$.
Write $\vee: \cC \to \cC^\op$ for the natural anti-equivalence.
Let~$I \coloneqq  P^\vee \in \cC^{\op}$, and let $I_n \coloneqq  (P/\Rad(E)^n P)^\vee$.
Choose generators~$x_1, \ldots, x_t$ of~$\rad(E)^n$ as a left ideal.
Then
\[
I_n = \bigcap_{i=1}^t \ker(x_i : I \to I), 
\]
and we need to prove that $I = \varinjlim_n I_n$.
By assumption, $\cC^{\op}$ is a locally finite category, so it suffices to prove that for each finite length subobjects $M \subset I$,
there exists~$n$ such that 
\[
\Rad(E)^n \subset \Ann_{E^\op}(M) \coloneqq  \{x \in E^\op: x|_M = 0\}.
\]
The right-pseudocompact topology on~$E = \End_\cC(P)$, which coincides with the left-pseudocompact topology on~$E^\op = \End_{\cC^\op}(I)$, 
is the $\Rad(E)$-adic topology, 
by Lemma~\ref{properties of compact modules}~\eqref{item: compact 7}. 
So we need to prove that $\Ann_{E^\op}(M)$ is open in the
left-pseudocompact topology on~$E^\op$.
This is true by definition of this topology, see \cite[Section~IV.4, end of page~396]{Gabrielthesis}, where $\Ann_{E^\op}(M)$ is denoted $\mathfrak{l}(M)$.
This concludes the proof that~$P$ is a complete left $E$-module in~$\cC$.

It now suffices to prove that $\text{--} \cotimes_E P$ is left adjoint to the equivalence $\Hom_{\cC}(P, \text{--})$.
Recall that $\cC = \Pro \cC^{\fl}$ and $\Mod_c(E^\op) = \Pro \Mod^{\fl}(E^\op)$.
We claim that~$\Hom_{\cC}(P, \text{--})$ and $\text{--} \cotimes_E P$ are Pro-extended from functors between $\cC^{\fl}$ and $\Mod^{\fl}(E^\op)$.

Assuming the claim, it suffices to prove that these restricted functors are adjoint.
By Corollary~\ref{cor:restricting-completed-tensor}, $\text{--} \cotimes_E P$
and $\text{--} \otimes_E P: \Mod(E^{\op}) \to \Ind(\cC)$ have the same restriction to~$\Mod^{\fl}(E^{\op})$.
By Proposition~\ref{prop:Eilenberg-Watts}, $\text{--} \otimes_E P$ is left adjoint to $\Hom_{\Ind(\cC)}(P, \text{--})$, whose restriction to~$\cC^{\fl}$ is $\Hom_{\cC}(P, \text{--})$.
Putting these together, we see that indeed $\text{--} \cotimes_E P$ and~$\Hom_{\cC}(P, \text{--})$ restrict to adjoint functors between $\Mod^{\fl}(E^{\op})$ and $\cC^{\fl}$. 
This concludes the proof of the proposition. %

There remains to prove the claim.
Since $\Hom_{\cC}(P, \text{--})$ and $\text{--} \cotimes_E P$ commute with cofiltered limits, it suffices to prove that they preserve objects of finite length.
This is true for $\Hom_\cC(P, \text{--})$ because it is an equivalence.
By d\'evissage, it now suffices to prove that $M \cotimes_E P$ has finite length if~$M$ is a simple $E^\op$-module.
Now $M \cotimes_E P$ is a quotient of $(E/\Rad(E)) \cotimes_E P$, which equals~$P/\Rad(E)P$, by Corollary~\ref{cor:restricting-completed-tensor}.
Since~$P$ has cosocle of finite length, it suffices to prove that $\Rad(E)P = \rad(P)$.
By~\cite[Prop.\ IV.4.12]{Gabrielthesis}, we have
\[\Rad(E) = \{\varphi \in E: \varphi(P) \subseteq \rad(P)\},\]
and so $\Rad(E)P \subseteq \rad(P)$.
If~$\rad(P)/\Rad(E)P \ne 0$, then there is a non-zero map $\lbar \varphi : P \to \rad(P)/\Rad(E)P$, because~$P$ is a generator of~$\cC$; 
and~$\lbar \varphi$ lifts to a map $\varphi : P \to \rad(P)$, because~$P$ is a projective object of~$\cC$.
But then~$\varphi \in \Rad(E)$, which contradicts the fact that~$\lbar \varphi \ne 0$.
\end{proof}

\subsection{\texorpdfstring{$\infty$}{infinity}-categories}\label{subsec: category theory background}
We now recall some basic results about
$\infty$-categories, and especially
stable $\infty$-categories, that
we will use below. %
Our main references are
ultimately~\cite{MR2522659, LurieHA, kerodon}, but we find it convenient to
refer to~\cite{MR3070515, MR3904731, MR3931682} for %
some of the results we use. In order to deal with set-theoretic issues, 
we follow the approach of~\cite{MR2522659}. %
In particular, we typically fix a 
Grothendieck universe, %
and sets are called \emph{small} if they belong to this fixed universe. Furthermore,
all limits and
colimits are assumed to be small. %
At one point in the sequel we employ the technical device of enlarging the universe 
(allowing us to regard an {\em a priori} large $\infty$-category  as
in fact being small), and so in Subsection~\ref{subsec:universes} below we
briefly recall how the process of enlarging the universe interacts with the
$\infty$-categorical constructions in which we're interested.
However, throughout the present subsection, the universe remains fixed.

It will sometimes be useful to refer to the more traditional literature of
triangulated categories. To this end, we recall that the homotopy category of a
stable $\infty$-category is a triangulated category~\cite[Rem.\
1.1.2.15]{LurieHA}, and that a functor $\cA\to\cB$ between stable
$\infty$-categories is an equivalence (resp.\ is fully faithful) if and only if
the induced functor on homotopy categories is  an equivalence (resp.\ is fully faithful).

The $\infty$-categories that we consider typically come in two sizes: {\em small},
and {\em presentable}.
More precisely, we say that an $\infty$-category~$\cC$ is \emph{small} if it is equivalent to
an $\infty$-category whose underlying simplicial set is small, i.e.\ if it is
essentially small in the sense of \cite[\href{https://kerodon.net/tag/03SM}{Tag~03SM}]{kerodon}.
Following
\cite[\href{https://kerodon.net/tag/06K6}{Tag~06K6}, \href{https://kerodon.net/tag/06NF}{Tag~06NF}]{kerodon}, 
for every small regular cardinal~$\kappa$, we say that~$\cC$ is \emph{$\kappa$-accessible} if
there exists a small $\infty$-category~$\cC_0$ such that
$\cC$ is equivalent to the $\kappa$-completion $\Ind_\kappa(\cC_0)$.
We then say that~$\cC$ is \emph{$\kappa$-presentable} if it is $\kappa$-accessible and cocomplete, and we say that~$\cC$ is \emph{presentable}
if it is $\kappa$-presentable for some small regular cardinal~$\kappa$.
By~\cite[\href{https://kerodon.net/tag/06PU}{Tag~06PU}]{kerodon}, presentable $\infty$-categories are complete.
In applications, we will usually be able to take $\kappa = \omega$, and
we refer to Section~\ref{subsubsec:Ind-cats} for a review of the $\Ind_\omega$-construction (which we will simply denote by~$\Ind$).

If~$\cC$ is a stable $\infty$-category, then it is $\omega$-presentable if and only if it is $\omega$-accessible: 
in fact, if it is $\omega$-accessible, then its cocompleteness follows from the fact that
it has filtered colimits and finite colimits.
Hence, if~$\cC_0$ is a small stable $\infty$-category, then $\Ind(\cC_0)$ is $\omega$-presentable, because it is stable, as explained in Remark~\ref{rem:stable indization}, 
and $\omega$-accessible, by definition.

We will refer to %
colimit-preserving functors between cocomplete stable $\infty$-categories as~\emph{continuous}.
In many texts, such as~\cite[Def.\ 5.3.4.5]{MR2522659}, {\em continuous} signifies the preservation of filtered colimits,
and so the functors between cocomplete stable $\infty$-categories
that we call continuous would there be called exact and continuous.
Since we have no occasion to consider non-exact functors, we have opted 
to incorporate the exactness condition into our definition of
continuous functors.

Finally, we will say that an arrow $f : x \to y$ in an $\infty$-category~$\cC$ is an \emph{isomorphism} if it represents an isomorphism
in the homotopy category~$h\cC$.
In the literature, this is often referred to as~$f$ being an \emph{equivalence}, and we will sometimes use this terminology interchangeably.

\subsubsection{Anima, and morphisms in \texorpdfstring{$\infty$}{infinity}-categories}
In the theory of $\infty$-categories, %
the $\infty$-category of {\em anima}\footnote{This is the terminology 
of~\cite{MR4681144}, and is short for {\em animated set}.} (referred to as the $\infty$-category of
{\em spaces} in \cite{MR2522659, LurieHA}, and as the $\infty$-category of
\emph{$\infty$-groupoids} in other references), which we denote $\operatorname{Ani}$,
plays the role that the category of sets plays in ordinary category theory.
 In particular,
if $x$ and $y$ are objects of an $\infty$-category~$\cC$, 
then we may form $\Maps_{\cC}(x,y)$, the anima %
of morphisms from $x$ to~$y$. 
The basic constructions related to this notion are described e.g.\ in~\cite[\S 1.2.2]{MR2522659}
and~\cite[\S 3.7]{MR3931682}.

The formation of $\Maps_{\cC}(\text{--},\text{--})$ is functorial in each of its arguments.  One way to formulate this statement
is via the existence of the Yoneda embedding (see \cite[Prop.\ 5.1.3.1]{MR2522659} or \cite[Thm.\ 5.8.13]{MR3931682}), which is a fully faithful functor
\begin{equation}
\label{eqn:Yoneda embedding}
\cC \hookrightarrow \Fun(\cC^{\op}, \operatorname{Ani});
\end{equation}
here the target denotes the $\infty$-category of functors from~$\cC^{\op}$ to~$\operatorname{Ani}$.

\subsubsection{Spectra}
\label{subsubsec:spectra}
Since the theory of stable $\infty$-categories inevitably bumps up against the theory
of {\em spectra} (in the algebraic topological sense), we give some brief reminders
regarding this latter theory for the non-expert. 
In fact, in the present paper, all the stable $\infty$-categories
that we consider will be $R$-linear (in a suitable sense)
for some commutative ring~$R$, and hence, as we explain below,
all the spectra that we encounter will
actually be objects of $D(R)$ (the derived category of $R$-modules).
Nevertheless, the references we give for $\infty$-categorical results typically work
in a framework that involves the consideration of general spectra;
thus we hope the following recollections will be helpful.

Intuitively, the $\infty$-category
$\Sp$ of spectra is obtained from the $\infty$-category of pointed anima by
making the suspension functor $\Sigma$ become invertible.  In particular, to each pointed anima~$x$,
there is associated a ``suspension spectrum'' $\Sigma^{\infty} x$, which is the image of~$x$
in the $\infty$-category of spectra with respect to the canonical functor
from pointed anima to spectra. 

The category $\Sp$ carries a canonical right complete $t$-structure\footnote{We recall
the notion of $t$-structure in Section~\ref{subsec:t-structures} below,
and the notion of right completeness in Section~\ref{subsubsec:completeness}.},
whose connective part $\Sp^{\leq 0}$ contains the full subcategory of suspension spectra,
and is obtained from this by taking the closure under colimits and extensions~\cite[Rem.~1.4.3.5]{LurieHA}.
Right completeness means that any spectrum $s$ can be written as a colimit
$$\colim s_n[-n] \iso s,$$
where the $s_n$ are connective spectra.  In fact, we may even arrange things 
so that the $s_n$ are suspension spectra,
say $s_n = \Sigma^\infty x_n$.
Eliding the difference between a pointed anima and its associated suspension spectrum,
we may then write this as
$$\colim \Sigma^{-n} x_n \text{``$\iso$''} s,$$ 
which in turn prompts the most naive notion of spectrum: namely, a spectrum
consists of a sequence $(x_n)$ of pointed anima equipped with transition morphisms
$\Sigma x_n \to x_{n+1}$.  (From this perspective,
the suspension spectrum of a pointed anima $x$ is given
by the particular sequence $(\Sigma^n x)$.)
These transition morphisms may be expressed in adjoint form 
as maps $x_n \to \Omega x_{n+1}$, and we say that the sequence~$(x_n)$ is an
{\em $\Omega$-spectrum} if these maps are homotopy equivalences.
It turns out that any ``naive'' spectrum is equivalent to an $\Omega$-spectrum,
and so one may restrict attention to $\Omega$-spectra when developing the theory.

Another way of thinking of an $\Omega$-spectrum is that we start with a pointed anima~$x_0$,
which we then enhance to a sequence $(x_n)$ equipped with equivalences $x_n \iso
\Omega x_{n+1}.$ In other words, $(x_n)$ is a sequence of successive deloopings
of~$x_0$, endowing~$x_0$ with the structure of an infinite loop space.  Motivated 
by this, we write $x_0 \coloneqq  \Omega^{\infty} s;$ this construction then yields
a functor $\Omega^{\infty}: \Sp \to \operatorname{Ani}_*$ (the codomain
denoting the $\infty$-category of pointed anima) which is right adjoint
to~$\Sigma^{\infty}$.

The preceding discussion is of course rather informal, but one rigorous version
of it can be found
in \cite[\S 1.4]{LurieHA}.  The discussion of \cite[\S 1.4.2]{LurieHA}  presents
a general explanation of how to
``stabilize'' an $\infty$-category~$\cC$ admitting finite limits,
yielding a stable $\infty$-category~$\Sp(\cC)$.  
If $\cC$ is furthermore pointed, then
\cite[Prop.~1.4.2.24]{LurieHA}
shows that
$\Sp(\cC)$ may be described as the $\infty$-category of ``$\Omega$-spectrum objects''
in~$\cC$.
(If $\cC$ admits  finite limits, and so in particular a final object,
then we may consider the $\infty$-category $\cC_*$ of pointed objects
of~$\cC$.  The forgetful functor $\cC_* \to \cC$ then induces
an equivalence $\Sp(\cC_*) \to \Sp(\cC)$~\cite[Rem.~1.4.2.18]{LurieHA},
and so it is in fact no loss of generality to assume that~$\cC$ is pointed;
cf.\ ~\cite[Rem.~1.4.2.25]{LurieHA}.)
Another discussion of the description of $\Sp(\cC)$ in terms of spectrum objects,
for a pointed $\infty$-category $\cC$ admitting finite limits, can be found
in~\cite[\S 2.3]{MR3070515}. 

In particular, taking $\cC$ to  be the $\infty$-category $\operatorname{Ani}_*$
of pointed anima,
we obtain the stable $\infty$-category $\Sp$ as its stabilization in the sense just discussed.
(The functor $\Omega^{\infty}$ is described in~\cite[1.4.2.20]{LurieHA}.
The construction of the $t$-structure on~$\Sp$ is given by~\cite[Prop.~1.4.3.4]{LurieHA},
and the left adjoint~$\Sigma^{\infty}$ to $\Omega^{\infty}$ is introduced
in the course of proving that result (and again in~\cite[Prop.~1.4.4.4]{LurieHA}).
The description of $\Sp^{\leq 0}$
is given in~\cite[Rem.~1.4.3.5]{LurieHA}.)

The heart of the $t$-structure on~$\Sp$ is the category of abelian groups.
If $A$ is an abelian group, then one typically denotes the corresponding spectrum
by~$HA$; it is the so-called {\em Eilenberg--MacLane spectrum} of~$A$.
We may furthermore describe the homotopy groups of spectra in terms of the $t$-structure
on~$\Sp$.
Indeed, if $s$ is a spectrum, then we write
$\pi_0(s) \coloneqq  \tau^{\geq 0}\tau^{\leq 0}s,$
and more generally
$\pi_{-n}(s) \coloneqq  \pi_0(s[n]) =(\tau^{\geq n}\tau^{\leq n} s)[n].$ 
Of course, the homotopy groups of spectra, defined in this manner,
also admit an interpretation as stable homotopy groups.\footnote{For example,
if we describe $s$ as an $\Omega$-spectrum $(x_m)$,
then $\pi_{-n}(s) = \pi_{m-n}(x_m)$ (for any $m \geq n$),
where $\pi_{\bullet}$ on the right hand side denotes the usual homotopy groups,
because the  
equivalence $x_m \iso \Omega x_{m+1}$ induces an isomorphism $\pi_{m-n}(x_m)
\iso \pi_{m -n +1}(x_{m+1}).$
More generally, if $s$ is the spectrum associated to a sequence
$(x_m)$ of pointed anima equipped with morphisms $\Sigma x_m \to x_{m+1},$
then $\pi_{-n}(s) = \colim_m \pi_{m-n}(x_m)$,
the transition morphisms being induced by the morphisms
$\pi_{m-n}(x_m) \to \pi_{m-n+1}(\Sigma x_m) \to \pi_{m-n+1}(x_{m+1})$.
So, if $x$ is a pointed anima,
then $\pi_{-n}\Sigma^{\infty} x$ is the $(-n)$th stable homotopy group of~$x$.}

If $R$ is a commutative ring, then its associated Eilenberg--MacLane spectrum~$HR$
has a natural $E_{\infty}$-structure.
By~\cite[Thm.\ 7.1.2.13]{LurieHA}, which is an unbounded version of the Dold--Kan correspondence,
we may identify the stable $\infty$-category of $HR$-module spectra
with the $\infty$-categorical version of the derived category~$D(R)$.
As indicated at the beginning of this section, these sorts of spectra
are the only ones that we will actually have occasion to consider in the present paper.

\subsubsection{\texorpdfstring{$\RHom$}{RHom} in stable \texorpdfstring{$\infty$}{infinity}-categories}
If $\cC$ is %
a {\em stable} $\infty$-category
and $x$ and $y$ are two objects of~$\cC$,
then in addition to forming the anima of morphisms from~$x$ to~$y$,
we may also form the
{\em mapping spectrum} from $x$ to $y$,
denoted by $\RHom_{\cC}(x,y)$.  
One way to express this is by saying that $\cC$ is canonically enriched over spectra.
(The reference~\cite[App.~C]{heyer20246functorformalismssmoothrepresentations} gives
a rather general discussion of enriched categories, which applies in particular
in our present context.) In the following discussion we state, and sketch the proofs,
of the basic results that we need related to this enriched structure.   

For a terse explanation of this enriched structure, see e.g.\ \cite[Rem.~7.1.2.2]{LurieHA}.
More details can be found in the discussion of \cite[\S 2.3]{MR3070515};
see also the discussion following~\cite[Def.~2.15]{MR3070515} (noting that
those authors write $\operatorname{maps}$ where we write~$\Maps$,
and write $\Maps$ where we write~$\RHom$).
Concretely,
the discussion following  ~\cite[Rem.~1.1.2.8]{LurieHA}  
shows that if $x$ and $y$ are objects of the stable $\infty$-category~$\cC$,
then we have 
$\Maps_{\cC}(\Sigma^{-n} x, y) \iso \Omega \Maps_{\cC}(\Sigma^{-n-1}x, y),$
and so the sequence $\bigl(\Maps_{\cC}(\Sigma^{-n}x ,y)\bigr)_{n \in \Z}$
forms a spectrum, or more specifically, an $\Omega$-spectrum.
This $\Omega$-spectrum is the mapping spectrum from $x$ to~$y$,
which (as already indicated) we denote by~$\RHom_{\cC}(x,y)$.
(Note that when $\cC$ is a stable $\infty$-category,
the mapping anima $\Maps_{\cC}(\text{--},\text{--})$ are canonically pointed, by the zero
morphism.) %

We may recover $\Maps$ from $\RHom$ by applying the functor~$\Omega^{\infty}$:
$$\Maps_{\cC}(x,y) \iso \Omega^{\infty}\RHom_{\cC}(x,y).$$
Also, we find that 
\begin{equation}
\label{eqn:homotopy of RHom}
\pi_{-n}\RHom_{\cC}(x,y) = \pi_0\Maps_{\cC}(x[-n],y) = \Hom_{h\cC}(x[-n],y),
\end{equation}
where $h\cC$ denotes the underlying homotopy category of~$\cC$.

In the examples we care about, our stable $\infty$-categories will (in a suitable sense)
be $R$-linear for some commutative ring~$R$;
indeed, they will be derived categories of $R$-linear
abelian categories, or closely related to such. 
The corresponding
mapping spectra $\RHom_{\cC}(x,y)$ will thus be module spectra over the Eilenberg--MacLane spectrum~$HR$ of~$R$.
As already noted above, \cite[Thm.\ 7.1.2.13]{LurieHA}
identifies the stable $\infty$-category of $HR$-module spectra
with the $\infty$-categorical version of the derived category~$D(R)$.
Accordingly, we regard $\RHom_{\cC}(x,y)$ as taking values in~$D(R)$. 

In traditional homological algebra, one often defines $\RHom$ 
on the derived category of an $R$-linear abelian category
as a derived functor (again with values in $D(R)$),
and so {\em a priori} $\RHom$ may be ambiguously defined.
However, the exactness properties that characterize derived functors in
the $\infty$-categorical context (see e.g.\ the various results discussed in
Section~\ref{subsubsec:derived functors} below) will show that $\RHom$ as
defined above will coincide with any definition of $\RHom$ as a derived 
functor that comes up in practice.

Relatedly,
if $\cA$ is an abelian category and $\cC \coloneqq  D(\cA)$ (in the sense of Definition~\ref{def:derived category} below),
then for any two objects $x, y$ 
of~$\cA$ (thought of as a full subcategory of~$\cC$ in the canonical way),  
we deduce from~\eqref{eqn:homotopy of RHom}
and standard homological algebra
(see e.g.~\cite[\href{https://stacks.math.columbia.edu/tag/06XU}{Tag~06XU}]{stacks-project})
that $\pi_{-n}\RHom_{\cC}(x,y)$
agrees with the Yoneda Ext-group $\Ext^n_{\cA}(x,y).$
In keeping with
this observation,
and with~\cite[Notation~1.1.2.17]{LurieHA}, we employ the usual notation
\begin{equation}\label{definition of Ext}
\Ext^n_{\cC}(x, y) \coloneqq  \pi_{-n}\RHom_{\cC}(x, y) 
\end{equation}
for any objects $x,y$ of an arbitrary stable $\infty$-category~$\cC$.\footnote{As indicated by our notation $\Sp^{\leq 0}$, we have employed cohomological conventions 
for the $t$-structure on~$\Sp$. This contrasts with~\cite{LurieHA}, which uses homological conventions, and defines $\pi_n(x) := \tau_{\leq 0}\tau_{\geq 0}x[-n]$ \cite[Def.\ 1.2.1.11]{LurieHA}.
However, the underlying stable $\infty$-category of spectra is the same, hence so is
the suspension functor on it; and since we have defined $\pi_{-n}(x)$ to be $\tau^{\geq 0}\tau^{\leq 0}x[n]$, which equals $\tau_{\leq 0} \tau_{\geq 0}x[n]$, our notion of $\pi_{-n}(x)$ 
coincides with the one in~\cite{LurieHA}.
Note also that, because of this, our $\Sp^{\leq 0}$ coincides with the category of spectra whose stable homotopy is concentrated in \emph{nonnegative} degrees.}

\subsubsection{Spectral variant of the Yoneda embedding}
\label{subsubsec:spectral Yoneda}
If $\cC$ is a stable $\infty$-category, then we may use $\RHom_{\cC}$ to construct a variant
of the usual Yoneda embedding~\eqref{eqn:Yoneda embedding}.
Namely, $y \mapsto \RHom_{\cC}(\text{--},y)$ gives a fully faithful exact functor
\begin{equation}
\label{eqn:spectral Yoneda}
\cC \hookrightarrow \Fun^{\mathrm{ex}}(\cC^{\op},\Sp)
\end{equation}
where $\Sp$ denotes the stable $\infty$-category of spectra, and
$\Fun^{\mathrm{ex}}$ denotes the stable $\infty$-category of exact functors between
stable $\infty$-categories.  
One construction of this functor is given in~\cite[Def.~2.15]{MR3070515},
where it is obtained as the stabilization (in the sense of~\cite[\S 1.4.2]{LurieHA})
of the usual Yoneda embedding~\eqref{eqn:Yoneda embedding},
bearing in mind that $\Sp(\cC) \iso \cC$ when $\cC$ is itself a stable $\infty$-category. 

\subsubsection{\texorpdfstring{$\RHom$}{RHom} as a right adjoint}

If $\cC$ is a presentable stable $\infty$-category, %
then we have an action of spectra
on~$\cC$~\cite[Rem.~7.1.2.2, Prop.~4.8.2.18, Rem.~4.8.2.20]{LurieHA}.
More precisely, the $\infty$-category ${\mathcal Pr}^L$ whose objects
are  presentable $\infty$-categories
and whose morphisms are left adjoint functors admits a symmetric monoidal structure
$(\cC,\cD) \mapsto \cC\otimes \cD$, where $\cC\otimes \cD$ is the presentable
$\infty$-category characterized by the property that it receives a functor
$\cC\times  \cD \to \cC\otimes\cD$ which is colimit preserving in each
variable separately, and which is universal for this property.  Furthermore,
the stable $\infty$-category $\Sp$ of spectra is idempotent with respect
to this structure, and the $\Sp$-local objects of ${\mathcal Pr}^L$, i.e.\
the presentable $\infty$-categories for which $\cC \iso \cC \otimes \Sp$, 
are precisely the presentable stable $\infty$-categories \cite[Prop.~4.8.2.18]{LurieHA}.
The inverse equivalence $\cC \otimes \Sp \to \cC$ then gives the
action of $\Sp$ on~$\cC$; in particular, this action is colimit preserving in either
variable.

The $\RHom_{\cC}$ construction may also be interpreted as right adjoint 
to the action of $\Sp$ on~$\cC$. In a little more detail,
for any object $x$ of~$\cC$, the functor $\text{--}\otimes x : \Sp \to \cC$
is colimit preserving, and so admits a right adjoint, which we (temporarily)
denote by~$G$.  If $\mathbf{S}$ denotes
the sphere spectrum (i.e.\ $\Sigma^\infty S^0$, where $S^0$ denotes the pointed
$0$-sphere), then
$\mathbf{S}$ is a unit object for $\text{--}\otimes x,$  
and so, for any object $y$ of~$\cC$, we find that
$$\Maps_{\cC}(x, y) 
\iso \Maps_{\cC}(\mathbf{S}\otimes x, y) \iso \Maps_{\Sp}\bigl(\mathbf{S}, G(y)\bigr).$$
Now for any object $t$ of~$\Sp,$ we have that
$$
\Maps_{\Sp}(\mathbf{S}, t) \iso \Maps_{\Sp}(\Sigma^{\infty} S^0 , t)
\iso \Maps_{\operatorname{Ani}_*}(S^0, \Omega^{\infty} t).
$$
Thus
$$\Maps_{\cC}(x,y) \iso \Omega^{\infty} G(y),$$
and replacing $y$ by $\Sigma^{n}y,$ we find that
$$\Maps_{\cC}(\Sigma^{-n} x,y) \iso \Omega^{\infty} \Sigma^{n}G(y).$$
Now any spectrum~$t$ is identified with
the $\Omega$-spectrum given by the sequence
$(\Omega^{\infty}\Sigma^n t)_{n \geq 0}.$ 
Thus we see that $G(y)$ is identified with the $\Omega$-spectrum 
$\Maps_{\cC}(\Sigma^{-n} x, y)$, which is to say with $\RHom_{\cC}(x,y)$. 
To summarize:
if $s$ is a spectrum and~$x, y \in \cC$, we have a natural isomorphism
\begin{equation}
\label{eqn:spectral action adjunction}
\Maps_{\cC}( s \otimes x, y) \iso \Maps_{\Sp}\bigl(s, \RHom_{\cC}(x,y)\bigr).
\end{equation}

\subsubsection{Formal properties of~\texorpdfstring{$\RHom$}{RHom}}\label{RHom and colimits}
We recall (and verify) some properties of~$\RHom_{\cC}$,
which are analogues of the usual formal properties for~$\Maps_{\cC}$. 
(The notion of a stable $\infty$-category being compactly generated
is recalled in some detail below; in particular, it implies presentability.)

\begin{lemma}\label{lem:formal properties of RHom}
Let~$\cC$ be a compactly generated stable $\infty$-category, and let~$x \in \cC$ be an object of~$\cC$. 
Then:
\begin{enumerate}
\item The functor $\RHom_\cC(x, \text{--})$ commutes with small limits.
\item If~$x = \colim_{i \in I}x_i$, then $\RHom_\cC(x, \text{--}) = \lim_{i \in I}\RHom_\cC(x_i, \text{--})$.
\item If~$x$ is compact in~$\cC$, then $\RHom_\cC(x, \text{--})$ commutes with small colimits.
\end{enumerate}
\end{lemma}
\begin{proof}
The description~\eqref{eqn:spectral action adjunction} of $\RHom_{\cC}(x,\text{--})$ as a right adjoint 
immediately implies part~(1).

We now prove part~(2).
Since $s \otimes \text{--}$ preserves colimits for any~$s \in \Sp$, we obtain a natural isomorphism
\[
\Maps_{\Sp}(s, \RHom_\cC(\colim x_i, \text{--})) \isoto \lim_i\Maps_{\Sp}(s, \RHom_\cC(x_i, \text{--}))
\]
and so part~(2) follows from the Yoneda lemma in~$\Sp$.

Finally, we prove part~(3).
By Lemma~\ref{lem: adjoint continuous} below,
applied to~$F \coloneqq \text{--}\otimes x$ and $G \coloneqq \RHom_\cC(x, \text{--})$,
it suffices to prove that~$\text{--}\otimes x$ preserves compact objects.
The compact objects in~$\Sp$ are precisely the colimits of finite diagrams, all of whose objects are shifts of the sphere spectrum~$\mathbf{S}$, or zero.
Since~$\mathbf{S} \otimes x = x$, and~$\text{--}\otimes x$ preserves colimits, we see that~$\text{--}\otimes x$ sends compact objects to 
colimits of finite diagrams, all of whose objects are shifts of~$x$, or zero. 
Since such colimits are compact if~$x$ is compact, this concludes the proof.
\end{proof}

\subsubsection{\texorpdfstring{$\RHom$}{RHom} and adjoints}\label{RHom and adjoints}
If $F: \cC \to \cD$ and $G: \cD \to \cC$ are (left, resp.\ right) adjoint functors between 
$\infty$-categories, then they are necessarily exact (being adjoints) and so are compatible
with suspension.
One then immediately verifies that $\RHom_{\cD}\bigl( F(x), y\bigr) \iso \RHom_{\cC}\bigl( x, G(y) \bigr)$
for any object $x$ of~$\cC$ and $y$ of~$\cD$.

If $F: \cC \to \cD$ is a continuous %
functor between presentable stable~$\infty$-categories,
then from the canonical nature of the $\Sp$-action on each of $\cC$ and~$\cD$,
one deduces that $F$ is compatible with these actions. Since $F$ admits a right adjoint~$G$,
one can also deduce this from the compatibility of adjunctions with $\RHom$ noted in the
preceding paragraph, together with the adjunction~\eqref{eqn:spectral action adjunction}.

\subsubsection{Ind-categories}\label{subsubsec:Ind-cats}
We recall %
the notion of Ind-completion of $\infty$-categories, 
following the exposition in \cite{MR2522659} and~\cite[\S 2.4]{MR3070515}. 
If~$\cC$ is a
small $\infty$-category, the Ind-category $\Ind(\cC)$ is defined to be
the formal closure of~$\cC$ under filtered colimits. %
Dually, we have
the Pro-category $\Pro(\cC)$, which may be defined as
$\Ind(\cC^\op)^\op$, and enjoys analogous %
properties. 
The $\infty$-category~$\Ind(\cC)$ is
characterised by the property that it admits filtered colimits, and
admits a fully faithful functor $\cC\into\Ind(\cC)$ which induces 
(via restriction) an
equivalence of $\infty$-categories of functors
\begin{equation}
\label{eqn:universal pro-property}
\Fun'\bigl(\Ind(\cC),\cD) \iso \Fun(\cC,\cD),
\end{equation}
for any $\cD$ admitting filtered colimits, where the domain denotes the full sub-$\infty$-category of
$\Fun\bigl(\Ind(\cC),\cD\bigr)$ consisting  
of those functors that preserve filtered colimits. 
See~\cite[Prop.\ 5.3.5.10]{MR2522659}.
The embedding
$\cC\into\Ind(\cC)$ preserves all limits, and if~$\cC$ admits
finite limits then so does $\Ind(\cC)$.

If we take the target $\infty$-category $\cD$ in~\eqref{eqn:universal pro-property}
to be~$\operatorname{Ani}^{\op}$, and compose the resulting equivalence with
the Yoneda embedding~\eqref{eqn:Yoneda embedding} for~$\Ind(\cC)$, we obtain a functor
$$\Ind(\cC) \hookrightarrow \Fun'\bigl(\Ind(\cC),\operatorname{Ani}^{\op}\bigr)^{\op}
\iso \Fun(\cC, \operatorname{Ani}^{\op})^{\op} = \Fun(\cC^{\op}, \operatorname{Ani}),$$
which is again fully faithful, and preserves filtered colimits (since the objects of~$\cC$ are compact in~$\Ind(\cC)$).
This functor identifies $\Ind(\cC)$ with the full subcategory
of~$\Fun(\cC^{\op},\operatorname{Ani})$ generated by~$\cC$ under the formation of filtered
colimits~\cite[Cor.~5.3.5.4]{MR2522659}.

\begin{remark}
\label{rem:stable indization}
If $\cC$ is furthermore a stable $\infty$-category,
then so is~$\Ind(\cC)$~\cite[Prop.~1.1.3.6]{LurieHA},  
and the canonical fully faithful embedding $\cC \hookrightarrow \Ind \cC$ is exact (see~\cite[Prop.~5.1.3.2, 5.3.5.14]{MR2522659}).
Furthermore,
if~$\cC$ and~$\cD$ are stable $\infty$-categories, then~\eqref{eqn:universal pro-property}
restricts to an equivalence
\begin{equation}
\label{eqn:exact Ind property}
\Fun^{\mathrm{ex}\prime}(\Ind \cC,\cD) \iso \Fun^{\mathrm{ex}}(\cC,\cD)
\end{equation}
(the source denoting exact and filtered colimit preserving --- or
equivalently, continuous --- functors, and the target denoting exact functors).
Indeed, since the inclusion $\cC \hookrightarrow \Ind \cC$ is exact,
the restriction of an exact functor $\Ind \cC \to \cD$ is again exact;
while conversely,
if $F: \cC \to \cD$ is exact, then $\Ind(F) : \Ind(\cC) \to \cD$ preserves cofibre sequences,
and so is also exact.
(To see this claim regarding preservation of cofibre sequences,
note that by the proof of~\cite[Prop.\ 1.1.3.6]{LurieHA},
every cofibre sequence in $\Ind(\cC)$ is a filtered colimit of cofibre sequences in $\cC$.
Since~$F$ preserves filtered colimits, and a colimit of cofibre sequences in~$\cD$ is a cofibre sequence, this implies that~$\Ind(F)$ preserves cofibre sequences, as claimed.)

In particular, 
taking~$\cD$ to be $\Sp^{\op}$ in~\eqref{eqn:exact Ind property}
and composing with the spectral Yoneda embedding~\eqref{eqn:spectral Yoneda}
for~$\Ind(\cC)$, we obtain an exact and fully faithful functor   
$\Ind \cC \hookrightarrow \Fun^{\mathrm{ex}}(\cC,\Sp^{\op})^{\op} = \Fun^{\mathrm{ex}}(\cC^{\op},\Sp)$.
In this way we may alternatively regard $\Ind \cC$ as the full subcategory of $\Fun(\cC^{\op},\Sp)$
generated by $\cC$ under filtered colimits.
\end{remark}

\begin{remark}
\label{rem:Indization as an adjunction}
In fact, the equivalence~\eqref{eqn:universal pro-property}
is part of an adjunction, as we now briefly explain.
Namely, we may consider the functor
\begin{equation}\label{eqn: indization restriction}
\Fun\bigl(\Ind(\cC),\cD) \to \Fun(\cC,\cD)
\end{equation}
given by restriction. Then ``Indization'', i.e.\ the inverse to the
equivalence~\eqref{eqn:universal pro-property},
is a left adjoint to this functor. 
This is a consequence of~\cite[Prop.\ 4.3.3.7]{MR2522659}, because 
the indization of a functor $F: \cC \to \cD$ is its left Kan extension through the canonical embedding $\cC \to \Ind(\cC)$, by~\cite[Lem.\ 5.3.5.8]{MR2522659}.
The left adjoint to~\eqref{eqn: indization restriction} is furthermore fully faithful, %
since it induces an equivalence onto its essential image $\Fun'\bigl(\Ind(\cC),\cD\bigr)$.
\end{remark}

\begin{rem}\label{rem:comparison of Ind-completions}
If~$\cC$ is a 1-category, 
then~\cite[\href{https://kerodon.net/tag/065F}{Tag~065F}]{kerodon},
shows that the formation of~$\Ind \cC$ is independent (up to canonical equivalence) 
of whether we regard $\cC$ as a $1$-category and form its Ind-completion in
the sense  of Section~\ref{subsubsec:Ind and Pro categories},
or whether we regard it as an $\infty$-category and form its Ind-completion
in the sense of the present discussion.
\end{rem}

If~$\cC$ admits filtered colimits then we let~$\cC^c$ denote the full sub-$\infty$-category of~$\cC$ consisting of
compact objects.%

\begin{lemma}
\label{lem:Ind compact}
Suppose that $\cC$ is an $\infty$-category admitting filtered colimits for which
$\cC^c$ is small.  If $\cC'$ is any sub-$\infty$-category of~$\cC^c$,
then the inclusion $\cC'\hookrightarrow \cC$ extends
{\em (}essentially{\em )} uniquely
to a filtered colimit preserving functor $\Ind \cC' \hookrightarrow \cC$,
which is again fully faithful.
\end{lemma}
\begin{proof}
The defining property~\eqref{eqn:universal pro-property}
(taking $\cC$ there to be~$\cC'$ and $\cD$ to be~$\cC$)
shows the existence and essential uniqueness of the extension.
Its full faithfulness 
is then an immediate application of~\cite[Prop.\ 5.3.5.11(1)]{MR2522659}
(taking the $\kappa$ there to be~$\omega$),
and is in any case easily verified directly. (It is here
that we use the compactness in~$\cC$ of the objects of~$\cC'$.)
\end{proof}

We say that~$\cC$ is \emph{compactly generated} if it is $\omega$-accessible, or equivalently
if it admits filtered colimits, if $\cC^c$ is small, and if the natural functor
$\Ind(\cC^c)\to\cC$
of Lemma~\ref{lem:Ind compact}
(sending a filtered diagram in~$\cC^c$ to its colimit in~$\cC$) is an
equivalence.\footnote{Note that 
in~\cite[Def.\ 5.5.7.1]{MR2522659} and \cite[\S 2.4]{MR3070515},
compactly generated is taken to mean $\omega$-presentable (rather than merely $\omega$-accessible),
while~\cite[\href{https://kerodon.net/tag/0673}{Tag~0673}]{kerodon} gives
essentially the same definition as the one we give here, but omits the requirement
that $\cC^c$ be small. In the stable case, $\omega$-accessible and $\omega$-presentable
coincide, so in this case (which is our primary focus) our definition coincides with
those of~\cite{MR2522659} and~\cite{MR3070515}.}
If $\cA$ is a small $\infty$-category,
then $\Ind (\cA)$ is idempotent-complete~\cite[Cor.~4.4.5.16]{MR2522659},
and the compact objects of~$\Ind(\cA)$ are the idempotent completion
of~$\cA$ \cite[Lem.\ 5.4.2.4]{MR2522659}.
We thus obtain a correspondence \cite[Prop.\ 5.4.2.17]{MR2522659} between small
idempotent-complete $\infty$-categories and compactly generated
$\infty$-categories. %
\footnote{By \cite[Prop.\
5.4.2.17]{MR2522659}, this correspondence can be upgraded to an equivalence of
$\infty$-categories, where the functors between compactly generated
$\infty$-categories are those which preserve filtered colimits and
compact objects, but we will not need this result; and indeed we will
consider functors between compactly generated categories which at
least \emph{a priori} do not preserve compact objects.}
If $\cC$ is compactly generated and furthermore stable,
then it is in fact cocomplete (since by the definitions of compact generation and
stability respectively, it admits both filtered and finite colimits).
Hence, if~$\cC$ is a stable $\infty$-category, then it is compactly generated if and only if it is $\omega$-presentable.

We say that a sub-$\infty$-category $\cA$ of a stable $\infty$-category $\cB$  
is a {\em stable sub-$\infty$-category} if it is full as a sub-$\infty$-category,
stable as an $\infty$-category, and if the inclusion $\cA \subseteq \cB$
is exact.  (This is what in~\cite{MR3904731} is called a stable subcategory.) If
$\cC$ is stable and admits filtered colimits, then $\cC^c$ is a stable
sub-$\infty$-category of~$\cC$, for example by~\cite[Lem.~1.1.3.3]{LurieHA}.
On the other hand, as we already noted, the Ind-completion of a small stable $\infty$-category is
stable. %
In particular, we see that
a compactly generated $\infty$-category~$\cC$ is stable if and only
if~$\cC^c$ is stable.

Any functor $F:\cC\to\cD$ between small $\infty$-categories induces a filtered colimit-preserving
functor $F_{\Ind(\cC)}: \Ind(\cC)\to\Ind(\cD)$, corresponding
under~\eqref{eqn:universal pro-property} to the composite functor
$\cC\stackrel{F}{\to}\cD\to\Ind(\cD)$;
we will sometimes denote $F_{\Ind(\cC)}$ simply
by~$F$ if this will not cause confusion. If~$F$ is fully faithful,
then so is ~$F_{\Ind(\cC)}$ \cite[Prop.\ 5.3.5.11(1)]{MR2522659}. If~$F:\cC\to\cD$ is an exact functor between small stable $\infty$-categories
then~$F_{\Ind(\cC)}:\Ind(\cC)\to\Ind(\cD)$ preserves all colimits, because it
preserves filtered colimits and %
finite colimits, 
and thus it is continuous.

\begin{remark}
\label{rem:compact generation}
Recall that there is a general notion of what it means
for a set of objects $X$ belonging to a cocomplete
stable $\infty$-category $\cC$ to {\em generate} $\cC$;
see e.g.\ \cite[Def.~A.8.5]{emerton2023introduction}.
Briefly: we define $X^{\perp}$ to be the full sub-$\infty$-category
of~$\cC$ whose objects are those $y$ for which $\RHom_{\cC}(x,y) = 0$
for all $x \in X$, and we say that $X$ {\em generates} ~$\cC$
if $X^{\perp} = 0.$ Equivalently, by e.g.\ \cite[Cor.\
A.8.7]{emerton2023introduction}, $X$ generates~$\cC$ if and only if the only
cocomplete stable subcategory of~$\cC$ containing all objects of~$X$ is~$\cC$
itself. We say that~$X$ is a \emph{set of compact generators} of~$\cC$ if every
object of~$X$ is compact, and if $X$ generates~$\cC$.

Lemma~\ref{lem:compact generation} below shows that the notion of compact generation
introduced above is compatible with this
more general notion of generation. In particular, it follows immediately from
Lemma~\ref{lem:compact generation}~\eqref{item:44} that if~$\cC$ has a set of compact generators, then~$\cC$
is compactly generated (because the existence of a set of compact generators implies that $(\cC^{c})^{\perp}=0$).
\end{remark}

\begin{lemma}
\label{lem:compact generation}
Let $\cC$ be a cocomplete 
stable $\infty$-category.

\begin{enumerate}
\item If $\cC$ is compactly generated,
then $(\cC^c)^{\perp} = 0.$
\item\label{item:44} 
Suppose that~$\cC'$ is a small
stable sub-$\infty$-category of~$\cC^c$,
with the property that
$(\cC')^{\perp} = 0$.  
Then $\cC^c$ is equal to the idempotent completion of~$\cC'$,
and $\cC$ is compactly generated.
In particular, if 
$(\cC^c)^{\perp} = 0$
then $\cC$ is compactly generated.
\end{enumerate}
\end{lemma}
\begin{proof}
Since $\Ind \cC^c$ is defined to be a subcategory of the presheaf category
of~$\cC^c$, we see that the right perpendicular of $\cC^c$ in $\Ind \cC^c$ 
is trivial; this proves~(1).

We now prove~(2).
Since $\cC$ is cocomplete,
Lemma~\ref{lem:Ind compact} shows that
the inclusion $\cC' \hookrightarrow \cC$
Ind-extends to a filtered colimit preserving
fully faithful functor
\begin{equation}
\label{eqn:Ind ff}
\Ind \cC' \hookrightarrow \cC.
\end{equation}
Part~(2) will follow if we further show
that~\eqref{eqn:Ind ff} is an equivalence, since we have
already observed that the analogue  
of~(2) holds for~$\Ind \cC'$.

Since the inclusion of~$\cC'$ into $\cC$ is an exact functor
between stable $\infty$-categories, the same is true of its Ind-extension~\eqref{eqn:Ind ff}, by Remark~\ref{rem:stable indization}.
Thus this Ind-extension is continuous. %
The source is furthermore cocomplete (being compactly generated and stable),
and hence its essential image is a cocomplete and stable sub-$\infty$-category
of~$\cC$, containing~$\cC'$.  Our assumption that
$(\cC'){}^\perp = 0$, interpreted in light of Remark~\ref{rem:compact generation},
shows that this essential image coincides with~$\cC$, as required.
\end{proof}

The following standard result gives a convenient criterion to check the full
faithfulness of a functor by checking it on compact generators.

\begin{prop}
  \label{prop:check-full-faithful-on-compact-generators}Let $F:\cC\to\cD$ be a
  continuous exact functor between compactly generated stable $\infty$-categories $\cC$ and~$\cD$, and
  let~$X$ be a set of compact generators of~$\cC$. Suppose that~$F$ preserves
  compact objects. %
  Then the following conditions are equivalent:
  \begin{enumerate}
  \item\label{item:45} $F$ is fully faithful,
i.e.\ 
for all objects $x,y\in \cC$, the induced
    map $$\Maps_{\cC}(x,y)\to\Maps_{\cD}(F(x),F(y))$$
    is an isomorphism.  
      \item\label{item:48}
For all objects $x,y\in \cC$, the induced
    map \begin{equation}\label{eqn:induced-map-on-RHom}\RHom_{\cC}(x,y)\to\RHom_{\cD}(F(x),F(y))\end{equation}
    is an isomorphism.  %
  \item\label{item:46} \eqref{eqn:induced-map-on-RHom} is an isomorphism for all $x,y\in X$.
  \item\label{item:47} For all objects $x,y\in X$ and each~$n\in\Z$, the induced
    map \[\Ext^{n}_{\cC}(x,y)\to\Ext^{n}_{\cD}(F(x),F(y))\]
    is an isomorphism. 
  \end{enumerate}
\end{prop}
\begin{proof}
The first and second conditions are equivalent by the discussion in Section~\ref{RHom and colimits}.
More precisely, $\RHom$ is constructed as an $\Omega$-spectrum from~$\Maps$,
while $\Maps$ is recovered from $\RHom$ by applying~$\Omega^{\infty}$.

The third and fourth conditions are equivalent by~\cite[Remark~1.4.3.8]{LurieHA}.
The second condition trivially implies the third, and it remains to prove the converse.

    Assume then that hypothesis~\eqref{item:46} holds, and
  consider the (full, stable) subcategory~$Y$ of~ $\cC$ consisting of the $y\in \cC$ with the property
  that~\eqref{eqn:induced-map-on-RHom} is an isomorphism for all~$x\in X$. By
  assumption, $Y$ contains every $y\in X$. In order to show that~$Y=\cC$, it
  therefore suffices to show that~$Y$ is closed under colimits in~$\cC$, which
  follows by noting that for any $x\in X$ and any colimit $\colim_{i\in
    I}y_{i}$, we have (by the compactness of~$x$ and of~$F(x)$, the
  continuity of~$F$, and Lemma~\ref{lem:formal properties of RHom}~(3))
  \begin{align*}
    \RHom_{\cC}(x,\colim_{i \in I} y_i) 
    &= \colim_{i \in I} \RHom_{\cC}(x, y_i) \\
    &= \colim_{i \in I} \RHom_{\cD}(F(x), F(y_i)) \\
    &= \RHom_{\cD}(F(x), \colim_{i \in I} F(y_i)) \\
    &= \RHom_{\cD}(F(x), F(\colim_{i \in I} y_i)).
\end{align*}

  We now consider the (full, stable) subcategory~$Z$ of~ $\cC$ consisting of the $x\in \cC$ with the property
  that~\eqref{eqn:induced-map-on-RHom} is an isomorphism for all~$y\in \cC$. We
  have just seen that~$Z$ contains~$X$, so to show that~$X=\cC$, it again
  suffices to show that~$Z$ is closed under colimits in~$\cC$.

  To see this, we have (using Lemma~\ref{lem:formal properties of RHom}~(2) and the continuity
  of~$F$) %
  \begin{align*}
    \RHom_{\cC}(\colim_{i \in I} x_i, y) 
    &= \lim_{i \in I} \RHom_{\cC}(x_i, y) \\
    &= \lim_{i \in I} \RHom_{\cD}(F(x_i), F(y)) \\
    &= \RHom_{\cD}(\colim_{i \in I} F(x_i), F(y)) \\
    &= \RHom_{\cD}(F(\colim_{i \in I} x_i), F(y)).
\end{align*}
  Thus~\eqref{eqn:induced-map-on-RHom} is an isomorphism for all~$x,y\in\cC$, as required.
  \end{proof}

\subsubsection{Quotients}\label{subsubsec: quotients}
An important role in our arguments will be played by {\em Verdier quotients} 
of stable $\infty$-categories. Given a stable sub-$\infty$-category~$\cA$ of a stable
$\infty$-category~$\cB$ (or, what is essentially the same data,
an exact fully faithful functor $\cA \hookrightarrow \cB$ between stable $\infty$-categories),
the Verdier quotient $\cB/\cA$ is a stable $\infty$-category receiving a functor
$\cB \to \cB/\cA$ whose composite with $\cA \hookrightarrow \cB$ is the zero functor,
and which is universal (i.e.\ initial) for these properties.    
In what follows, we recall the construction of the Verdier quotient,
and then discuss its interaction with the formation of~$\Ind$-categories.

The key point in the construction of $\cB/\cA$ comes by noting that,
given an exact functor of stable $\infty$-categories $\cB \to \cC$,
the objects in (the essential image of) $\cA$ all map to zero
precisely if all the arrows in $\cB$ with fibres lying in (the essential image of)
$\cA$ become invertible in~$\cC$.  Thus the construction of
Verdier quotients becomes a special case of the construction of
$\infty$-categorical localizations, a construction that we now briefly recall.

If~$\cC$ is an $\infty$-category, and $W$ is a collection of arrows of~$\cC$, 
then the {\em localization} $\cC[W^{-1}]$ is an $\infty$-category
receiving a functor $\cC \to \cC[W^{-1}]$ under which all the elements
of~$W$ become invertible, and which is universal (i.e.\ initial) for 
these properties.
This construction (under the slightly more circumlocutious name of
``the $\infty$-category obtained from $\cC$ by inverting the 
set of morphisms~$W$'') is the subject of~\cite[Def.\ 1.3.4.1, Rem.~1.3.4.2]{LurieHA},
and in the context of
small $\infty$-categories
it is also discussed in~\cite[\S I.3]{MR3904731},
where it is  called ``Dwyer--Kan localization''.
It is also discussed in~\cite[\href{https://kerodon.net/tag/01M4}{Tag [01M4}]{kerodon},
where it is simply called ``localization'' (and it this last terminological
choice that we follow).

In~\cite[Thm.~I.3.3]{MR3904731} it is shown that if $\cA$ is a stable sub-$\infty$-category
of the stable $\infty$-category~$\cB$, and if $W$ denotes the set of arrows in~$\cB$
whose fibre lies in~$\cA$, then the localization $\cB[W^{-1}]$ is again stable,
and hence (for the reason already explained above)
satisfies the universal property of the Verdier quotient.

\begin{remark}
\label{rem:localization ambiguity}
As already mentioned, in~\cite[\S~I.3]{MR3904731}, a smallness hypothesis is imposed.
In our setting, in which we use universes, this is not essential,
since we can always enlarge the universe to make any given $\infty$-category
small. One could worry that the universal property that $\cC[W^{-1}]$ satisfies
is universe dependent (i.e.\ that it ceases to hold if one enlarges the
universe), but in fact this does not happen (i.e.\ the localization of
$\cC$ at $W$ remains such after any enlargement of the universe).
To see this, one  notes~\cite[\href{https://kerodon.net/tag/05ZN}{Tag [05ZN}]{kerodon}
that the size of any localization of~$\cC$ is bounded in terms of the size of~$\cC$
itself. 
\end{remark} 

Among all localizations of $\infty$-categories, a particularly important
class is given by the so-called {\em Bousfield localizations}.
These are localizations for which the canonical functor $\cC \to \cC[W^{-1}]$
admits a fully faithful right adjoint. This right adjoint
then realizes $\cC[W^{-1}]$ as a full subcategory of~$\cC$ itself. 
In~\cite{MR2522659} and~\cite{LurieHA}, these
are referred to simply as ``localizations'', and are characterized in terms
of the resulting full subcategories of~$\cC$; in~\cite{kerodon} 
they are called ``reflective localizations''.
See~\cite[Ex.~I.3.4.3]{LurieHA} for one explanation of why these ``localizations'' (defined
as they are in~\cite{MR2522659}, in terms of certain full subcategories of~$\cC$)
are localizations in the sense that we are using here.

In practice, the right adjoint required for a Bousfield localization will be
constructed via an application of the adjoint functor theorem (see the discussion
below for more on this), and so Bousfield
localizations typically arise in the context of large, rather than small, $\infty$-categories.
We note, though, that~\cite[Lem.~I.3.4]{MR3904731} shows that
one can use Yoneda embeddings to
map (the opposite of) any localization functor $\cC \to \cC[W^{-1}]$ 
into an associated  Bousfield localization of $\Fun(\cC,\operatorname{Ani})
\to \Fun(\cC[W^{-1}],\operatorname{Ani})$; thus Bousfield localizations
can be used as a tool to study arbitrary localizations (such as localizations
of small categories, which is the actual context under consideration in~\cite{MR3904731}). 

Since $\Ind$-categories especially lend themselves to the construction of right adjoints
(again, see the discussion below),
for our purposes, the role of Bousfield
localizations will be to aid in the analysis of Verdier quotients of Ind-categories
of stable $\infty$-categories. 
For example, if $\cA \hookrightarrow \cB$ is a fully faithful 
continuous functor of presentable stable~$\infty$-categories,
then it follows from~\cite[Prop.~5.6]{MR3070515}
that the Verdier quotient $\cB/\cA$ is a Bousfield localization of~$\cB$.

The following proposition describes how the formation of Ind-categories
and the formation of Verdier quotients interact.

\begin{prop}
\label{prop:Ind is exact}
Let $\cA$ be an idempotent complete stable sub-$\infty$-category of the small idempotent
complete stable $\infty$-category~$\cB$.
Then the $\Ind$-extension $\Ind(\cB) \to \Ind(\cB/\cA)$ 
of the quotient functor $\cB \to \cB/\cA$ (of small stable $\infty$-categories)
is a %
Bousfield localization at the collection of arrows of~$\Ind(\cB)$ whose cofibre lies in~$\Ind(\cA)$.
Furthermore,
$\cA = (\Ind \cA )\cap \cB$ {\em (}the intersection
taking place in $\Ind \cB${\em )}. 
\end{prop}
\begin{proof}
The claim that $\Ind(\cB) \to \Ind(\cB/\cA)$ is a Bousfield localization at the collection of arrows with cofibre in~$\Ind(\cA)$ is part of~\cite[Prop.\ I.3.5]{MR3904731}.
It thus follows that the objects of $\Ind(\cA) \cap \cB$ are precisely the objects of~$\cB$ that become equivalent to~$0$ in~$\cB/\cA$. 
Since~$h(\cB) \to h(\cB/\cA)$ is by construction a Verdier quotient of triangulated categories, the equality $\Ind(\cA) \cap \cB = \cA$ thus follows from
\cite[\href{https://stacks.math.columbia.edu/tag/05RK}{Tag~05RK}]{stacks-project}, since~$\cA$ is idempotent complete.
\end{proof}%

\begin{rem}
  \label{rem:quotients-compactly-generated}
Suppose now that
$\cC \hookrightarrow\cD$ is a fully faithful continuous functor between
compactly generated stable $\infty$-categories, which furthermore preserves compact objects,
so that it can be identified as
the Ind-ification of its restriction $\cC^c \hookrightarrow \cD^c$
to the sub-$\infty$-categories of compact objects.
Then Proposition~\ref{prop:Ind is exact} can be rephrased as the statement %
that the induced functor
$$\cD^c/\cC^c \to \cD/\cC = \Ind \cD^c / \Ind \cC^c $$
induces an equivalence
$$\Ind (\cD^c/\cC^c) \iso \cD/\cC.$$
It thus follows that, in this case, the Verdier quotient $\cD/\cC$ is compactly generated.
\end{rem}

\subsubsection{Adjoint functors}\label{app-subsub-adjoint-functors}
By the
adjoint functor theorem \cite[Cor.\ 5.5.2.9(1)]{MR2522659},
a functor $F: \cC \to \cD$ between cocomplete and compactly generated
stable $\infty$-categories is continuous if and only if it admits a
right adjoint~$G$. 
In particular, if~$F:\cC\to\cD$ is an exact functor
between small stable $\infty$-categories then~$F:\Ind(\cC)\to\Ind(\cD)$
always admits a right adjoint; and passing to opposite categories, we
see that the induced functor $\Pro(\cC)\to\Pro(\cD)$ always admits a
left adjoint. In addition, if $F:\cC\to\cD$ is left adjoint to
$G:\cD\to\cC$, then the induced functors $F:\Ind(\cC)\to\Ind(\cD)$ and
$G:\Ind(\cD)\to\Ind(\cC)$ are an adjoint pair. This can be seen either by the defining properties
of the Ind extended functors, or by~\cite[Prop.\
5.3.5.13]{MR2522659}.

The
following is %
a special case of~\cite[Prop.\ 5.5.7.2]{MR2522659}. %
\begin{lem}
  \label{lem: adjoint continuous}If $F: \cC \to \cD$ is a
  continuous functor between compactly generated
  stable $\infty$-categories, then $F$ preserves compact objects if
  and only if its right adjoint $G$ is continuous.
\end{lem}

\subsection{Recollements}
We will use the notion of recollements, which in the $\infty$-category setting
were defined by Lurie in~\cite[\S A.8]{LurieHA}.  We recall the definition here,
in the context of recollements of stable $\infty$-categories.

\begin{defn}%
  \label{def:recollement}
Let $\catA$, $\catB$, and $\catD$ be stable $\infty$-categories,
and let $i_*:\catB \to \catD$ and $j_*:\catA \rightarrow \catD$
be fully faithful exact functors.
We will often regard~$\catB$ and~$\catA$ as full subcategories of~$\cA$ by identifying them with their essential images. %
We say that $\catD$ is the {\em recollement} of $\catB$ and $\catA$
if the following hypotheses hold:
\begin{enumerate}
\item $i_*$ and $j_*$ admit (necessarily exact) left adjoints $i^*$ and $j^*.$
\item $j^*i_* = 0$.
\item $i^*$ and $j^*$ are jointly conservative, in the sense that if 
$i^*X = 0$ and $j^*X = 0$ for some object $X$ of~$\catD$, then $X = 0.$
\end{enumerate}
\end{defn}

We now describe a method for recognizing recollements.
We will usually assume the following hypothesis.

\begin{hypothesis}
  \label{hyp:usual hyp for semiorthogonal}
Let $i_*:\catB\into\catD$ be a  fully faithful functor between compactly
generated stable
$\infty$-categories, which is continuous and preserves compact
objects. (Equivalently, $i_{*}$ is the Ind-extension of a fully
faithful exact functor $i_{*}:\catB^c\into\catD^c$.) We will freely regard~$\catB$ as a full subcategory of~$\catD$, by identifying it with the essential image of~$i_{*}$.
 Set
$\catA\coloneqq \catD/\catB$, and write $j^*:\catD\to\catA$ for the quotient
functor. 
\end{hypothesis}

The following lemma contains the specializations of some results from Section~\ref{subsubsec: quotients} to the present context.

\begin{lem}\label{lem: abstract semiorthogonal decomposition}%
Assume Hypothesis~{\em \ref{hyp:usual hyp for semiorthogonal}}.
  Then:
 \begin{enumerate}
 \item $\catB \cap \catD^c = \catB^c$.
\item  $\catA$ is compactly generated, and the functor $j^*$ is continuous and preserves compact objects. The kernel of~$j^{*}$ is~$\catB$.
\item $j^*$ has a fully faithful continuous right adjoint $j_*:\catA\to\catD$.

\end{enumerate}
In particular, the functors~$i_{*},j^{*},j_{*}$ are the Ind-extensions of their restrictions
 $i_{*}:\catB^c\into\catD^c$, $j^{*}:\catD^c\to\catA^c$, $j_*:\catA^{c}\into\catD$.
\end{lem}
\begin{proof} 
Applying
Proposition~\ref{prop:Ind is exact}
to the
inclusion $\catB^c\subseteq\catD^c$, 
we see that $\catB \cap \catD^c = \catB^c$, that~$j^*$ is the $\Ind$-extension 
of the quotient functor $\catD^c \to \catD^c/\catB^c$ (hence it preserves compact objects), and
that the kernel of~$j^*$ is $\catB$.
Furthermore, $j^*$ has a fully faithful right adjoint~$j_*$ (and so
in particular is continuous),
because it is a Bousfield localization, and $\catA$ is compactly generated, because
it is equivalent to $\Ind(\catD^c/\catB^c)$. 
  \end{proof}

We can now state our criterion for recognizing recollements.

\begin{lemma}
\label{lem:semiorthogonal to recollement}
Assume Hypothesis~{\em \ref{hyp:usual hyp for semiorthogonal}}. Then the following conditions are equivalent:
\begin{enumerate}
\item\label{item:72}  $\catD$ is the recollement of $\catB$ and~$\catA$.
  \item\label{item:68}  $i_*$ admits a left adjoint~$i^*$.  
  \item\label{item:73} $i_{*}|_{\catB^c}:\catB^c\to\catD^c$ admits a left adjoint.
\end{enumerate}
\end{lemma}
\begin{proof}
By definition, \eqref{item:72} implies \eqref{item:68}.
For the converse implication, by Lemma~\ref{lem: abstract semiorthogonal decomposition}, the only point to verify is that if~ $X$ is an object of~$\catD$ with $i^*X = 0$
  and $j^*X = 0$ then~$X= 0.$ 
  But if $j^* X = 0$, then $X$ is isomorphic to~$i_* Y$ for some $Y \in \cA_Z$. %
  Now $i_*$ is fully faithful, so the counit $i^*i_* Y \to Y$ is an isomorphism.
  If furthermore $i^* X = 0$, this implies $Y = 0$, hence $X = 0$, as desired.  
  
  We now turn to the equivalence of~  \eqref{item:68} and~ \eqref{item:73}.
  If~\eqref{item:73} holds, then passing to $\Ind$-categories gives  ~\eqref{item:68},
  by the discussion of %
  Section~\ref{app-subsub-adjoint-functors}.
Conversely, by Lemma~\ref{lem: adjoint continuous} since~$i_{*}$ is continuous, any left adjoint~$i^{*}$ preserves compact objects, and so restricts to a left adjoint of~$i_{*}|_{\catB^c}$.
\end{proof}

Consider now a continuous functor
$F:\catD \to \catD'$ of compactly generated stable $\infty$-categories,
each of which satisfies
Hypothesis~\ref{hyp:usual hyp for semiorthogonal};
that is, we have fully faithful, continuous, and compact
object preserving functors of compactly generated stable $\infty$-categories
 $i_*:\catB \hookrightarrow \catD$ and
 $i'_*:\catB' \hookrightarrow \catD'$.
 We say that $F(\catB) \subseteq \catB'$ if the essential image of~$Fi_*$ is a %
 subcategory of the essential image of~$i'_*$.

  \begin{lem}
    \label{lem:defn-of-F(B)-B'}The following conditions are equivalent:
    \begin{enumerate}
    \item $F(\catB) \subseteq \catB'$.
    \item There exists a continuous functor $F_{\catB}:\catB \to \catB'$,
characterized by the existence of
a natural isomorphism
\begin{equation}
\label{eqn:F and i lower star}
Fi_* \iso i'_* F_{\catB}. 
\end{equation}
\item There exists a continuous functor $F_{\catA}: \catA \coloneqq  \catD/\catB 
\to \catA' \coloneqq  \catD'/\catB',$
characterized by the existence of
a natural isomorphism
\begin{equation}
\label{eqn:F and j upper star}
j'^{*}F \iso F_{\catA} j^*.
\end{equation}
\item $j'^{*}Fi_{*}=0$.
\end{enumerate}
  \end{lem}
  \begin{proof} The first, second and fourth conditions are equivalent by definition.
    The third implies the fourth, because
    $j^{*}i_{*}=0$. Conversely, if the fourth condition holds, then
   $j'^{*}F:\catD\to\catA'$ is the zero functor %
   on~$\catB$, so induces a
   functor $F_{\catA}: \catA\to\catA'$ and a natural isomorphism as in~(3),
by the universal property of the quotient. \end{proof}%

Suppose now given $\catD$, $\catD'$, and $F:\catD \to \catD'$ as above,
and suppose that $F(\catB) \subseteq \catB'$,
so that the equivalent conditions of 
Lemma~\ref{lem:defn-of-F(B)-B'} hold. 
We will write $\eta, \eta'$ (resp.\ $\epsilon, \epsilon'$) for the units (resp.\ the counits) of the adjunctions $(j^*, j_*), (j'^*, j'_*)$.
Precomposing  the unit $\eta': \id_{\catD'}\to j'_{*}j'^{*}$ with $Fj_{*}$,
and taking into account~\eqref{eqn:F and j upper star}, 
we obtain a natural transformation
\begin{equation}
\label{eqn:F and j lower star}
F j_* \xrightarrow{\eta' F j_*} j'_*j'^*Fj_* \xrightarrow{j_*'\eqref{eqn:F and j upper star}j_*} j_*'F_{\catA}j^*j_* \xrightarrow{j_*'F_{\catA}\epsilon} j'_* F_{\catA},
\end{equation}
which need not be an isomorphism in general. %

Now
suppose %
further that
that each of $i_*$ and $i'_*$ admits a left adjoint $i^*$ and~$i'^*$, so that
by Lemma~\ref{lem:semiorthogonal to recollement}, each of~$\catD$ and~$\catD'$ is a recollement.
Postcomposing~ \eqref{eqn:F and i lower star} with
$i'^{*}$, and recalling that $i'_{*}$ is fully faithful (and so the counit is an isomorphism), we obtain an isomorphism
\[
i'^{*}Fi_* \iso i'^*i'_* F_{\catB} \iso F_{\catB}. 
\]
It follows that if we postcompose the unit of
adjunction $1_{\catD}\to i_{*}i^{*}$ with $i'^{*}F$, then we obtain
 a natural transformation %
\begin{equation}
\label{eqn:F and i upper star}
i'^*F \to i'^*Fi_*i^* \isom F_{\catB} i^*,
\end{equation}
which again need not be an isomorphism in general.

\begin{rem}\label{rem:morphism of recollements}
Using terminology from~\cite[Def.\ 2.3, Def.\ 2.6]{shah2022recollements}, if each of~$\catD$, $\catD'$ is a recollement, then a functor~$F$ that satisfies the conditions of Lemma~\ref{lem:defn-of-F(B)-B'} 
defines a \emph{morphism of recollements} if and only if~\eqref{eqn:F and i upper star} is an
isomorphism,
and a \emph{strict morphism of recollements} if and only if both~\eqref{eqn:F and j lower star} and~\eqref{eqn:F and i upper star} are isomorphisms. 
\end{rem}

In the framework just described, we have the following criterion for~$F$ to be fully faithful. 
\begin{prop}
\label{prop:gluing full faithfulness, recollement version}
Suppose that $F:\catD \to \catD'$ is a continuous functor between
stable $\infty$-categories each satisfying Hypothesis~{\em \ref{hyp:usual hyp for semiorthogonal}},
and suppose that  $F(\catB) \subseteq \catB'$. %
Suppose furthermore that 
\begin{enumerate}
\item each of $i_*:\catB \to \catD$ and $i'_*:\catB' \to \catD'$
admits a left adjoint $i^*$ and~$i'^{*}$ respectively;
\item the natural transformations~{\em \eqref{eqn:F and j lower star}} and {\em \eqref{eqn:F and i upper star}}
are natural isomorphisms.
\end{enumerate}
Then if each of $F_{\catA}$ and $F_{\catB}$ are fully faithful, so is~$F$.
\end{prop}
\begin{proof}
Write $F(\catD), F_{\catB}(\catB), F_{\catA}(\catA)$ for the essential images of~$F, F_{\catB}, F_{\cA_U}$.
The natural isomorphisms~\eqref{eqn:F and i lower star} and~\eqref{eqn:F and j upper star} imply that we have a commutative diagram
\begin{equation}\label{eqn:restricting the recollement}
\begin{tikzcd}
\catB \arrow[r, "i_*"] \arrow[d, "F_{\catB}"] & \catD \arrow[r, "j^*"] \arrow[d, "F"] & \catA \arrow[d, "F_{\catA}"] \\
F_{\catB}(\catB) \arrow[r, "i'_*"] & F(\catD) \arrow[r, "j'^*"] & F_{\catA}(\catA).
\end{tikzcd}
\end{equation}
Assumption~(1) and Lemma~\ref{lem:semiorthogonal to recollement} show that~$\catD$ is the recollement of~$\catB$ and~$\catA$ via~$i_*$ and~$j_*$.
Assumption~(2) implies in particular that~$j'_*$ and~$i'^*$ restrict to functors
$j'_*: {F_{\catA}}(\catA) \to F(\catD)$ and $F(\catD) \to F_{\catB}(\catB)$, which are then adjoint to (the appropriate restrictions of)~$i'_*$ and~$j'^*$.
We therefore see that all the properties in Definition~\ref{def:recollement} are satisfied, and so $F(\catD)$ is the recollement of $F_{\catB}(\catB)$ and~$F_{\catA}(\catA)$ via~$i'_*$ and~$j'_*$.
Furthermore, assumption~(2) also shows that~\eqref{eqn:restricting the recollement} is a ``strict morphism of recollements'' in the sense of~\cite[Def.\ 2.6]{shah2022recollements} (compare Remark~\ref{rem:morphism of recollements}).
We can therefore replace~$\catD'$, $\catB'$, $\catA'$ by~$F(\catD)$, $F_{\catB}(\catB)$, $F_{\catA}(\catA)$ respectively,  thus reducing ourselves to the case that~$F_{\catB}$ and~$F_{\catA}$ are equivalences of categories.
We may then apply~\cite[Rem.\ 2.7]{shah2022recollements}, or equivalently~\cite[Prop.\ A.8.14]{LurieHA}, to conclude that $F$ is also an equivalence, as required. %
\end{proof}

\subsection{Enlarging the universe, and the Ind Pro construction}
\label{subsec:universes}
In the sequel we will use Ind Pro categories.
In order to deal with the size issues that this entails, we make a digression
on change of universe.

\subsubsection{Change of universe}
Writing~$\uf{U}$ for the universe (implicitly)
considered in Section~\ref{subsec: category
  theory background}, we now choose another universe~$\uf{V}$ with
$\uf{U}\in\uf{V}$.
We wish to compare size-related notions, such as smallness, small colimits, and compact objects,
with respect to the two universes $\uf{U}$ and~$\uf{V}$.

To this end,
we now say that an $\infty$-category $\cC$ is $\uf{U}$-small
if it is small with respect to~$\uf{U}$.
Given a $\uf{U}$-small $\infty$-category~$\cC$, we let $\Ind^{\uf{U}} \cC$ denote
the Ind category of $\cC$ built by formally adjoining 
$\uf{U}$-small filtered colimits.
Given an $\infty$-category~$\cC$,
we say that an object $X$ in $\cC$ is $\uf{U}$-compact if
$\Maps_{\cC}(X, \text{--})$ commutes with the formation of $\uf{U}$-small filtered colimits,
and write $\cC^{\uf{U}-c}$
to denote the full sub-$\infty$-category of $\cC$ consisting of $\uf{U}$-compact objects.
We say that $\cC$ is $\uf{U}$-compactly generated if
it admits $\uf{U}$-small filtered colimits, if $\cC^{\uf{U}-c}$ is $\uf{U}$-small,
and if $\Ind^{\uf{U}} \cC^{\uf{U}-c} \to \cC$ is an equivalence.

All these notions have evident $\uf{V}$-analogues,
given by replacing the universe $\uf{U}$ by the universe $\uf{V}$ in the definitions.
It is evident that, on a literal level, the notions of $\uf{U}$-compactly
generated and $\uf{V}$-compactly generated objects are distinct.  However,
there is in fact a sense in which compactly generated categories are insensitive to the choice 
of universe.\footnote{We could phrase this in a precise mathematical manner by describing 
the construction $\cC \mapsto \cC^{\uf{V}}$ that we define below as inducing an equivalence
of certain $\infty$-categories of $\infty$-categories, although
we don't do that here.}  
Namely, if $\cC$ is $\uf{U}$-compactly generated, 
so that 
$$\Ind^{\uf{U}} \cC^{\uf{U}-c} \iso \cC,$$
then we define
$$\cC^{\uf{V}} \coloneqq  \Ind^{\uf{V}} \cC^{\uf{U}-c}.$$
This definition makes sense, because the $\uf{U}$-small category $\cC^{\uf{U}-c}$ is in particular
$\uf{V}$-small.
By construction, $\cC^{\uf{V}}$ is $\uf{V}$-compactly generated,
and $(\cC^{\uf{V}})^{\uf{V}-c} = \cC^{\uf{U}-c}$.
Thus $\cC^{\uf{V}}$ ``promotes'' $\cC$ to a $\uf{V}$-compactly generated $\infty$-category
whose $\uf{V}$-compact objects coincide with the $\uf{U}$-compact objects in the original 
$\infty$-category~$\cC$.
One can think of $\cC^{\uf{V}}$ as a ``version'' of $\cC$ which is adapted to the larger 
universe~$\uf{V}$.

Similarly, any $\uf{U}$-continuous functor (i.e.\ a functor compatible with
$\uf{U}$-small filtered colimits)  $F:\cC \to \cD$ between $\uf{U}$-compactly generated
$\infty$-categories extends canonically to a $\uf{V}$-continuous
functor $F^{\uf{V}}: \cC^{\uf{V}} \to \cD^{\uf{V}}$.  If $F$ preserves $\uf{U}$-compact
objects, then $F^{\uf{V}}$ preserves $\uf{V}$-compact objects.

\subsubsection{Ind Pro categories}
\label{subsubsec:Ind Pro categories}
As we already explained,
our motivation for introducing change-of-universe considerations is so that we
may consider $\Ind \Pro$ categories.  More precisely, if $\cC$ is a $\uf{U}$-small
$\infty$-category, then we may define $\Pro^{\uf{U}}\cC$ in the obvious manner.
This is no longer a $\uf{U}$-small category, but it {\em is}
$\uf{V}$-small (e.g.\ by~\cite[Thm.\ 1(i)]{MR3731483}), %
and so we may then define the Ind completion
$$\Ind^{\uf{V}} \Pro^{\uf{U}} \cC,$$
a $\uf{V}$-compactly generated $\infty$-category.

We note the following lemma related to computing morphisms in $\Ind^{\uf{V}}\Pro^{\uf{U}} \cC$.

\begin{lemma}
\label{lem:Ind Pro maps}
Let $\cC$ be a $\uf{U}$-small $\infty$-category,
let $X = \lim_{i\in I} X_i$ be an object
of $\Pro^{\uf{U}} \cC$, written as a $\uf{U}$-small cofiltered limit
of objects $X_i$ of~$\cC$,
and let $Y$ be an object of $\Ind^{\uf{V}} \cC$.
Then the fully faithful embeddings
$\Pro^{\uf{U}} \cC \hookrightarrow \Ind^{\uf{V}} \Pro^{\uf{U}} \cC$
and
$\Ind^{\uf{V}} \cC \hookrightarrow \Ind^{\uf{V}} \Pro^{\uf{U}} \cC$
allow us to regard both these objects as belonging to~$\Ind^{\uf{V}} \Pro^{\uf{U}} \cC,$ 
and we have
$$\Maps_{\Ind^{\uf{V}} \Pro^{\uf{U}} \cC}( X, Y)
\iso \colim_i \Maps_{\Ind^{\uf{V}} \cC}(X_i, Y).$$
\end{lemma}
\begin{proof}
We may write 
$Y = \colim_{j\in J} Y_j$, a $\uf{V}$-small filtered colimit of objects $Y_j$ of~$\cC$. 
Furthermore, the embedding %
$\Ind^V \cC \hookrightarrow \Ind^V \Pro^U \cC$ 
preserves %
$V$-small filtered colimits.
Then
\begin{multline*}
\Maps_{\Ind^{\uf{V}} \Pro^{\uf{U}} \cC}(X, \colim_j Y_j)
\iso \colim_j \Maps_{\Pro^{\uf{U}} \cC}(X,  Y_j)
\\\iso \colim_j \colim_i \Maps_{\cC}(X_i, Y_j)
\iso \colim_i \colim_j \Maps_{\cC}(X_i,Y_j)
\\
\iso \colim_i \Maps_{\Ind^{\uf{V}} \cC}(X_i, \colim_j Y_j),
\end{multline*}
as claimed.
\end{proof}

Suppose now that $F:\cC \to \cD$ is an exact functor between $\uf{U}$-small
(and hence $\uf{V}$-small)
stable $\infty$-categories, 
so that the induced functor
$F_{\Ind^{\uf{V}} \cC} : \Ind^{\uf{V}} \cC \to \Ind^{\uf{V}} \cD$ 
admits a right adjoint~$G: \Ind^{\uf{V}} \cD \to \Ind^{\uf{V}} \cC$.
The functor $F$ also induces an exact functor between $\uf{V}$-small stable $\infty$-categories
$F_{\Pro^{\uf{U}} \cC}: \Pro^{\uf{U}} \cC \to \Pro^{\uf{U}} \cD,$
whose Ind-ification
$$F_{\Ind^{\uf{V}}\Pro^{\uf{U}} \cC}: \Ind^{\uf{V}}\Pro^{\uf{U}}\cC 
\to \Ind^{\uf{V}}\Pro^{\uf{U}}\cD$$
similarly admits a right adjoint
$$G_{\Ind^{\uf{V}}\Pro^{\uf{U}} \cD}: 
\Ind^{\uf{V}}\Pro^{\uf{U}}\cD\to
\Ind^{\uf{V}}\Pro^{\uf{U}}\cC.$$

The following lemma shows that 
$G_{\Ind^{\uf{V}}\Pro^{\uf{U}} \cD}$ and~ $G$ are compatible in the
evident way. 

\begin{lemma}
\label{lem:Ind Pro adjoints}
The functor $G_{\Ind^{\uf{V}} \Pro^{\uf{U}} \cD}$ is continuous, and so is determined by
its restriction to $\Pro^{\uf{U}}\cD$.  Furthermore, 
if $Y\coloneqq\lim_{j\in J} Y_j$ is an object of $\Pro^{\uf{U}}\cD$, written as a cofiltered limit 
of objects $Y_j$ of $\cD$, then 
$$G_{\Ind^{\uf{V}} \Pro^{\uf{U}} \cD}(Y) 
\iso \lim_{j\in J} G(Y_j).$$
{\em (}The limit exists in $\Ind^{\uf{V}} \Pro^{\uf{U}} \cC$.{\em
  )}
In particular, the restrictions of
$G_{\Ind^{\uf{V}} \Pro^{\uf{U}} \cD}$ and of~$G$ to $\cD$ coincide.
\end{lemma}
\begin{proof}
Since $F_{\Ind^{\uf{V}} \Pro^{\uf{U}} \cC}$
is defined as the Ind-ification of~$F_{\Pro^{\uf{U}} \cC}$,
it preserves $\uf{V}$-compact objects, and so its right adjoint
$G_{\Ind^{\uf{V}} \Pro^{\uf{U}} \cD}$
is continuous by Lemma~\ref{lem: adjoint continuous}.  Furthermore, right adjoints necessarily commute with limits (which also implies the claimed existence of 
$\lim_{j\in J} G(Y_j)$ in $\Ind^V \Pro^U \cC$). 
Thus it remains to prove the final claim of the lemma, namely that the
restrictions to $\cD$ of
$G_{\Ind^{\uf{V}} \Pro^{\uf{U}} \cD}$ and of~$G$ coincide. %

Let $Y$ be an object of~$\cD$.  Then $G(Y)$ is an object of $\Ind^{\uf{V}} \cC$,
and we must exhibit a canonical equivalence
$$\Maps_{\Ind^{\uf{V}} \Pro^{\uf{U}} \cC} \bigl(X,G(Y)\bigr)
\iso 
\Maps_{\Ind^{\uf{V}} \Pro^{\uf{U}} \cD} \bigl(F_{\Ind^{\uf{V}} \Pro^{\uf{U}}\cC} (X),Y\bigr)
$$
for any object $X$ of~$\Ind^{\uf{V}}\Pro^{\uf{U}} \cC$.
Since $F_{\Ind^{\uf{V}} \Pro^{\uf{U}} \cC}$
is continuous, the formation of both mapping spaces is compatible with the formation 
of colimits in~$X$, and so it suffices to consider the case when
$X = \lim_{i\in I} X_i$ is an object of~$\Pro^{\uf{U}} \cC$, written as a cofiltered limit
of objects $X_i$ of~$\cC$.
We then find that
\begin{multline*}
\Maps_{\Ind^{\uf{V}} \Pro^{\uf{U}} \cC} \bigl(X,G(Y)\bigr)
\iso\colim_i \Maps_{\Ind^{\uf{V}} \cC} \bigl( X_i, G(Y) \bigr) 
\\
\iso\colim_i \Maps_{\Ind^{\uf{V}} \cD} \bigl( F(X_i), Y \bigr) 
\iso \Maps_{\Ind^{\uf{V}} \Pro^{\uf{U}} \cD} \bigl(\lim_{i\in I} F(X_i), Y\bigr)
\\
\iso \Maps_{\Ind^{\uf{V}} \Pro^{\uf{U}} \cD} \bigl(F_{\Ind^{\uf{V}}\Pro^{\uf{U}}\cC} (\lim_{i\in I} X_i),
Y\bigr),
\end{multline*}
the first and third equivalences being provided by Lemma~\ref{lem:Ind Pro maps},
the second equivalence being provided by the adjunction between
$F_{\Ind^{\uf{V}} \cC}$ and $G$ (and recalling that $F_{\Ind^{\uf{V}} \cC}$ 
is the Ind-ification of~$F$, and so coincides with $F$ on objects of~$\cC$),
and the final equivalence being provided by
the definition of $F_{\Ind^{\uf{V}}\Pro^{\uf{U}} \cC}$ as the Ind Pro-ification of~$F$.
\end{proof}
We also have the following useful lemmas.
\begin{lem}
  \label{lem:F-IndPro-commutes-limits-colimits}$F_{\Ind^{\uf{V}}\Pro^{\uf{U}}
    \cC}$ commutes with all limits and colimits in~$\Ind^{\uf{V}}\Pro^{\uf{U}}
    \cC$.
  \end{lem}
  \begin{proof}
    We already observed that~$F_{\Ind^{\uf{V}}\Pro^{\uf{U}} \cC}$ admits a right
    adjoint~$G_{\Ind^{\uf{V}} \Pro^{\uf{U}} \cD}$, so it commutes with
    colimits. Similarly, the functor
    $F_{\Pro^{\uf{U}} \cC}: \Pro^{\uf{U}} \cC \to \Pro^{\uf{U}} \cD$ admits a
    left adjoint, and thus so does its Ind-extension
    $F_{\Ind^{\uf{V}}\Pro^{\uf{U}} \cC}$; so
    $F_{\Ind^{\uf{V}}\Pro^{\uf{U}} \cC}$ also commutes with limits.
  \end{proof}

\begin{lem}
  \label{lem:F-fully-faithful-iff-IndProF}Let $F:\cC \to \cD$ be an exact functor between $\uf{U}$-small stable $\infty$-categories. Then the following conditions are
  equivalent.
  \begin{enumerate}
  \item\label{item:6} $F$ is fully faithful.
  \item\label{item:7} $F_{\Ind^{\uf{V}} \cC}$ is fully faithful.
  \item\label{item:8} $F_{\Pro^{\uf{U}} \cC}$ is fully faithful.
  \item\label{item:5} $F_{\Ind^{\uf{V}}\Pro^{\uf{U}} \cC}$ is fully faithful.
  \end{enumerate}%
\end{lem}
\begin{proof}%
  Since the restriction of~$F_{\Ind^{\uf{V}} \cC}$ to~$\cC$ agrees with~$F$ by
  definition, the second condition implies the first. Conversely, suppose that
  $F$ is fully faithful, and let~$X,Y$ be arbitrary objects
  of~$\Ind^{\uf{V}} \cC$, which we write as filtered colimits
  $X=\colim_{i\in I}X_i$, $Y = \colim_{j\in J} Y_j$ of objects $X_i,Y_j$
  of~$\cC$. Then $F_{\Ind^{\uf{V}} \cC}(X)=\colim_{i\in I}F(X_i)$, $F_{\Ind^{\uf{V}} \cC}(Y)=\colim_{j\in
    J}F(Y_j)$, so we have 
    \begin{align*}
      \Maps_{\Ind^{\uf{V}} \cD}(F_{\Ind^{\uf{V}} \cC}(X),F_{\Ind^{\uf{V}} \cC}(Y))=
      \lim_{i\in I}\colim_{j\in J}\Maps_{\cD}(F(X_i),F(Y_j))\\
      =\lim_{i\in I}\colim_{j\in J}\Maps_{\cC}(X_i,Y_j)=\Maps_{\Ind^{\uf{V}} \cC}(X,Y).
    \end{align*} 
    This shows that the first two
  conditions are equivalent; the remaining equivalences follow by analogous
  arguments which we leave to the reader.
\end{proof}

\subsubsection{Notational conventions}
When we apply the preceding material,
we will again suppress the choice of $\uf{V}$ (as well as of the original
universe~$\uf{U}$), and will simply write $\Ind \Pro \cC$.  %
When we implement this construction,
it will be furthermore understood that all 
$\uf{U}$-compactly generated categories $\cD$ under consideration, and any $\uf{U}$-continuous
functors $F$ between them,  are
replaced by their ``$\uf{V}$-versions'' $\cD^{\uf{V}}$ and~$F^{\uf{V}}$.   
This stipulation is required because, for example, we will want to view a $\uf{U}$-compactly generated category~$\cA$ as
a full subcategory of $\Ind \Pro (\cA^c)$, and we will do this by means of the canonical embedding
\[
\cA^V = \Ind^V(\cA^c) \to \Ind^V \Pro^U(\cA^c).
\]

\subsection{An application of Ind Pro categories to full
  faithfulness}%
\label{subsec:semiorthogonal and Ind Pro}
We prove a criterion for full faithfulness
(Proposition~\ref{prop:gluingfull-faithfulness-Ind-Pro-version-compact-version} below)
which will in turn be used in the proof of our main theorem. 

\subsubsection{Constructing a recollement}
\label{subsubsec:constructing a recollement}
We return to the setting of Hypothesis~\ref{hyp:usual hyp for semiorthogonal},
so that we assume given a fully faithful, continuous, compact object
preserving functor $i_*:\catB \hookrightarrow \catD$
of compactly generated stable $\infty$-categories,
with a corresponding Verdier quotient functor $j^*: \catD \to \catB$
which (by Lemma~\ref{lem: abstract semiorthogonal decomposition})
is also continuous and compact object preserving.
In fact, Lemma~\ref{lem: abstract semiorthogonal decomposition} shows that
this data is obtained by Ind-extending the corresponding functors
$i_*:\catB^c \hookrightarrow \catD^c$ and $j^*:\catD^c \to \catA^c$
obtained by restriction to the various full sub-stable-$\infty$-categories 
of compact objects.

Starting with these functors on the categories of compact objects,
we now take Pro-extensions, to obtain
functors that we denote
$\widehat{i}_*: \Pro \catB^c \hookrightarrow \Pro \catD^c$
and $\widehat{j}^* : \Pro \catD^c \to \Pro \catA^c.$
The evident Pro-analogue of Proposition~\ref{prop:Ind is exact}
shows that $\widehat{j}^*$ realizes $\Pro \catA^c$ as 
the Verdier quotient $\Pro \catD^c/\Pro \catB^c$.

Next we Ind-extend these Pro-extended functors (using a change of universe,
as discussed in the preceding section, in order to regard these Pro-categories
as being small) to obtain functors (notated in the same manner)
$$
\widehat{i}_*: \Ind \Pro \catB^c \hookrightarrow \Ind \Pro \catD^c
$$
and 
$$\widehat{j}^* : \Ind \Pro \catD^c \to \Ind \Pro \catA^c.$$
Again, Proposition~\ref{prop:Ind is exact} shows that $\widehat{j}^*$ 
realizes $\Ind \Pro \catA^c$ as the Verdier quotient $\Ind \Pro \catD^c/\Ind \Pro \catB^c$.
Furthermore, we are now again in the context of Hypothesis~\ref{hyp:usual hyp for semiorthogonal},
with the the role of $\catB^c$, $\catD^c$, and $\catA^c$ now being played
by $\Pro \catB^c$, $\Pro \catD^c$, and $\Pro \catA^c$.

We can now apply
Lemma~\ref{lem: abstract semiorthogonal decomposition} 
twice.  Firstly, returing to the context of $i_*:\catB \hookrightarrow \catD$ and 
its quotient~$\catA$, we obtain the fully faithful right adjoint to~$j^*$,
which (following the notation of that Lemma) we denote by
$j_*:\catA \hookrightarrow \catD.$ 
But secondly, considering the context of $\widehat{i}_*: \Ind \Pro \catB^c
\hookrightarrow \Ind \Pro \catD^c$ and its quotient $\Ind \Pro \catA^c,$
we obtain the fully faithful right adjoint to~$\widehat{j}^*$, 
which we denote by
$$\widehat{j}_*: \Ind \Pro \catA^c \hookrightarrow \Ind \Pro \catD^c.$$
We can compare these two constructions.  Indeed,
the canonical fully faithful 
embeddings $\catA^c \hookrightarrow \Pro \catA^c$ 
and $\catD^c \hookrightarrow \Pro \catD^c$
induce fully faithful embeddings
$$\catA \iso \Ind \catA^c \hookrightarrow \Ind \Pro \catA^c$$
and 
$$\catD \iso \Ind \catD^c \hookrightarrow \Ind \Pro \catD^c,$$
and we see that we are in a particular  
instance of the general setting of
Lemma~\ref{lem:Ind Pro adjoints}.
That lemma then shows that $\widehat{j}_*$ restricts to $j_*$.

In the present situation, where we have made a preliminary Pro-extension before performing 
an Ind-extension, there is one more %
functor that arises (that
is not part of the general package coming from Hypothesis~\ref{hyp:usual hyp for semiorthogonal}).
Namely, the Pro-extended functor $\widehat{i}_*$ has a left adjoint
$$\widehat{i}^*: \Pro \catD^c \to \Pro \catB^c.$$ Then, since Ind-extension preserves
adjunctions, we see that its Ind-extension
$$\widehat{i}^*: \Ind \Pro \catD^c \to \Ind \Pro \catB^c$$
(which we denote by the same symbol) is left adjoint to (the Ind-extended version of)~$\widehat{i}_*$.
Lemma~\ref{lem:semiorthogonal to recollement} then shows
 that the category $\Ind \Pro \cA^c$ is the recollement
of $\Ind \Pro \cA_Z^c$ and~$\Ind \Pro \cA_U^c$, via~$\ihat_*$ and~$\jhat_*$.

\subsubsection{A criterion for full faithfulness}
Suppose now that we have a functor $F_{\catD^c}: \catD^c \to (\catD')^c$ with $\Ind$-extension   $F:\catD \to \catD'$, %
$\Pro$ extension %
$F_{\Pro \catD^c}: \Pro\catD^c \to \Pro (\catD')^c$, and $\Ind\Pro$ extension %
$F_{\Ind \Pro \catD^c}: \Ind\Pro\catD^c \to \Ind\Pro (\catD')^c$. 
Suppose also that we have a stable sub-$\infty$-category~$i'_{*}:(\catB')^c\into(\catD')^c$, %
and write~$\catA'$, $\ihat'_{*}$ etc.\ for the corresponding
constructions; and suppose further that $F_{\catD^c}(\catB^c) \subseteq (\catB')^c$, so that $F_{\catD^c}$ restricts to a functor $F_{\catB^c}:\catB^c\to(\catB')^c$. 
Then %
$F_{\catB^c}$ induces functors $F_{\catB}:\catB\to\catB'$, $F_{\Ind\Pro\catB^c}:\Ind\Pro\catB^c\to \Ind\Pro(\catB')^c$,  and functors $F_{\catA^c}:\catA^c\to(\catA')^c$, $F_{\catA}:\catA\to\catA'$
and $F_{\Ind\Pro\catA^c}:\Ind\Pro\catA^c\to \Ind\Pro(\catA')^c$.

The following criterion for~$F$ to be fully faithful is used to prove the main theorem of the paper (see Theorem~\ref{thm:F-is-fully-faithful}).
\begin{prop}
  \label{prop:gluingfull-faithfulness-Ind-Pro-version-compact-version}
  Suppose, with notation as above, that the functor $F_{\catD^c}: \catD^c \to (\catD')^c$ satisfies:
\begin{enumerate}
  \item\label{item:78}  $F_{\catA^c}$ and $F_{\catB^c}$ are fully faithful.
\item\label{item:76} The natural transformation of functors $\catA^{c}\to \cA'$
\[F j_* \to j'_* F_{\catA}\]  \emph{(}arising from precomposing the unit of adjunction $ \id_{\catD'}\to j'_{*}j'^{*}$  with  $Fj_{*}$\emph{)} is
  a natural isomorphism. 
\item\label{item:77}  The natural transformation of functors   $\catD^c\to \Pro(\catB')^c$ 
\[\ihat'^*F_{\catD^c} \to
F_{\Pro\catB^c}\ihat^{*}\]
  \emph{(}arising from postcomposing $\ihat'^{*}F_{\Pro\catD^c}$ with the unit of
adjunction $1_{\catD}\to \ihat_{*}\ihat^{*}$\emph{)} %
is
  a natural isomorphism.  
    \end{enumerate}
Then~$F$ and~$F_{\catD^c}$ are fully faithful.
\end{prop}
\begin{proof} By Lemma~\ref{lem:F-fully-faithful-iff-IndProF}, assumption~\eqref{item:78} implies that~$F_{\Ind \Pro \catA^c}$ and~$F_{\Ind \Pro \catB^c}$ are
  fully faithful, and that the proposition is equivalent to proving that ~$F_{\Ind \Pro \catD^c}$ is fully faithful.
  We will deduce this from Proposition~\ref{prop:gluing full
    faithfulness, recollement version}, applied with the functor $F:\cA\to\cA'$ in that proposition being our functor~$F_{\Ind\Pro\catD^c}:\Ind\Pro\catD^c\to\Ind\Pro(\catD')^c$.
In order to apply this result, we need to verify that Hypothesis~\ref{hyp:usual hyp for semiorthogonal} holds (which we have already done, in the paragraph preceding the statement of this proposition) 
and to show that the natural transformations 
\begin{equation}\label{eqn:F-IndProAc-jstar}
  F_{\Ind\Pro\catD^c} \jhat_* \to \jhat'_* F_{\Ind\Pro\catA^c}
\end{equation}
and 
\begin{equation}
\label{eqn:F on completions}
\ihat'^*F_{\Ind\Pro\catD^c} \to
F_{\Ind\Pro\catB^c}\ihat^{*}
\end{equation} are natural isomorphisms.
(These are the natural transformations~\eqref{eqn:F and j lower star} and~\eqref{eqn:F and i upper star}, with~$\catD$ there being our
$\Ind\Pro\catD^c$.)

These statements follow easily from our assumptions ~\eqref{item:76} and~\eqref{item:77}, using  Lemma~\ref{lem:F-IndPro-commutes-limits-colimits} (for $F_{\Ind\Pro\catA^c}$ and $F_{\Ind\Pro\catB^c}$), Lemma~\ref{lem:Ind Pro adjoints} (for $\jhat_*$ and $\jhat'_*$), and the compatibility of $\ihat^{*},\ihat'^{*}$ with cofiltered limits. %
\end{proof}

\subsection{\texorpdfstring{$t$}{t}-structures}\label{subsec:t-structures}
We recall some of the key facts
about $t$-structures that we will use.
We refer to \cite[App.~A]{emerton2023introduction} for a more complete recollection, 
as well as further references.

The data of a $t$-structure on the stable $\infty$-category $\cD$ is determined by the
specification
of a full sub-$\infty$-category $\cD^{\leq 0} \hookrightarrow \cD$.
The full sub-$\infty$-category $\cD^{\geq 1} \hookrightarrow \cD$ is then defined
to be the right orthogonal to~$\cD^{\leq 0}$.  
Of course $\cD^{\leq 0}$ should satisfy some conditions (e.g.\
$\cD^{\leq 0}[1] \subseteq \cD^{\leq 0}$),
which amount
to the requirement that the image of $\cD^{\leq 0}$ and $\cD^{\geq 1}$
in the homotopy category of $\cD$ define a $t$-structure on this homotopy category
in the usual sense.

For any integer $n$,
we define $\cD^{\leq n} \coloneqq  \cD^{\leq 0}[-n],$ and $\cD^{\geq n} \coloneqq  \cD^{\geq 1}[1-n].$
For any pair of integers $a \leq b,$ we write $\cD^{[a,b]} \coloneqq  \cD^{\geq a} \cap
\cD^{\leq b}.$ 
We say that an object $X$ of $\cD$ is {\em bounded above} (resp.\ {\em below})
if it lies in  $\cD^{\leq n}$ (resp.\ $\cD^{\geq n}$) for some~$n$.

We say that the $t$-structure on $\cD$ is {\em bounded} (resp.\ {\em bounded above}, resp.\
{\em bounded below})
if every object of~$\cD$ is bounded (resp.\ bounded above, resp.\ bounded below).

The heart of the $t$-structure, often denoted $\cD^{\heartsuit}$,
is defined to be $\cD^{\leq 0} \cap \cD^{\geq 0}  \eqcolon \cD^{[0,0]}.$
It is an abelian category.

Each inclusion $\cD^{\leq n} \hookrightarrow \cD$ admits a right adjoint~$\tau^{\leq n}:
\cD\to \cD^{\leq n}$,
while each inclusion $\cD^{\geq n} \hookrightarrow \cD$ admits a left adjoint
$\tau^{\geq n}: \cD \to \cD^{\geq n}.$
These satisfy various relations; for example, if $a \leq b$ then
there is a natural isomorphism $\tau^{\geq a} \tau^{\leq b} \iso \tau^{\leq b} \tau^{\geq a}$
of functors $\cD \to \cD^{[a,b]}$. For any integer $n$ and any object $X$ of~$\cD$,
$\tau^{\geq n}\tau^{\leq n} X$ 
is an object of $\cD^{[n,n]} = \cD^{\heartsuit}[-n].$
We then define the functor\footnote{We follow traditional homological algebra notation
(with {\em cohomological} conventions) by writing~$H^n$.  This suits our purposes,
since the stable $\infty$-categories that we consider in the body of the paper
will typically be derived categories of various abelian categories of modules,
or variants thereof. In the case of the stable $\infty$-category of spectra,
and other stable $\infty$-categories of a more homotopical nature,
it is more usual to write $\pi_{-n}$. We have applied this convention
in the discussion of Section~\ref{subsubsec:spectra} above, but in all
other cases that we consider, we use the cohomological notation established here.} 
 $H^n: \cD \to \cD^{\heartsuit}$
via $H^n(X) \coloneqq  (\tau^{\geq n} \tau^{\leq n} X)[n].$ Recall also that  for each
$X\in\cD$ and~$n\in\Z$ there is a
fibre sequence \begin{equation}\label{eqn:fibre-sequence-truncations}\tau^{\leq n-1} X \to X \to \tau^{\geq n}X. \end{equation}

\subsubsection{Completeness}
\label{subsubsec:completeness}
If $\cD$ is a stable $\infty$-category endowed with a $t$-structure,
then we have a functor
$$\cD \to \lim \cD^{\leq n}$$
defined by $X \mapsto (\tau^{\leq n} X)$ (the transition morphisms being 
given by the obvious truncations).   
We say that $\cD$ is {\em right complete} if this functor is an equivalence.
Unwinding this definition, the full faithfulness of this morphism
is equivalent to the requirement that for any object $X$ of~$\cD$,
the canonical morphisms
$\tau^{\leq n} X \to X$ 
induce an isomorphism $\colim \tau^{\leq n} X \iso X.$
The essential surjectivity is equivalent to the requirement 
that if $\{X_n\}$ is any sequence of objects endowed with identifications
$X_n \iso \tau^{\leq n} X_{n+1}$, then $\colim_n X_n$ exists in~$\cD$.

There is a dual notion of {\em left complete}; namely, this is the condition
that the canonical functor
$$\cD \to \lim \cD^{\geq n},$$
defined by $X \mapsto (\tau^{\geq n} X)$, should be an equivalence. (Note that
this is a limit with $n\to-\infty$.)
The full faithfulness now amounts to the requirement that for any object $X$ of~$\cD$,
the canonical morphisms
$X \to \tau^{\geq n} X$ induce an isomorphism
$$X \iso \lim_n \tau^{\geq n} X.$$
The essential surjectivity is equivalent to the requirement 
that if $\{X_n\}$ is any sequence of objects endowed with identifications
$\tau^{\geq n}X_{n-1} \iso  X_{n}$, then $\lim_n X_n$ exists in~$\cD$.

\subsubsection{\texorpdfstring{$t$}{t}-exactness of functors}
If $\cD$ and $\cD'$ are two stable $\infty$-categories endowed with $t$-structures,
then a functor $F:\cD \to \cD'$ is {\em right $t$-exact} (resp.\ {\em left $t$-exact})
if $F(\cD^{\leq 0}) \subseteq (\cD')^{\leq 0}$
(resp.\ $F(\cD^{\geq 0}) \subseteq (\cD')^{\geq 0}$),
and {\em $t$-exact}
if it is both right and left $t$-exact.

\subsubsection{Some cautions}
If $\cD$ is a stable $\infty$-category endowed with a $t$-structure,
it is natural to ask whether
if $X$ is an object of $\cD$ for which
$H^i(X) = 0$ for all~$i$,
is it necessarily the case that $X = 0$?
While we recall some positive results in this direction,
Example~\ref{ex:Ind/Pro objects with vanishing cohomology}
below shows that it need not hold in general.

\begin{lemma}
\label{lem:cohomology is conservative on bounded objects}
If $X$ is an object of $\cD$ which is bounded, and $H^i(X) = 0$ for all~$i$,
then $X = 0$.
\end{lemma}
\begin{proof}
By assumption, $X \in \cD^{[a,b]}$ for some $a \leq b$.  
We have the fibre sequence
$$ \tau^{\leq b-1} X \to X \to H^b(X)[-b]  = 0,$$
so that 
$\tau^{\leq b-1} X \iso X$. Now $\tau^{\leq b-1}X$ is an object of $\cD^{[a,b-1]},$
and the lemma follows by descending induction on the amplitude $b -a \geq 0$.
\end{proof}

\begin{lemma}
\label{lem:cohomology is conservative on bounded above objects for left complete}
If the $t$-structure on $\cD$ is left complete, and $X$ is an object of $\cD$
which is bounded above and for which $H^i(X) = 0$ for all~$i$, then $X = 0.$
\end{lemma} 
\begin{proof}
Since $H^i(X) = 0$ for all $i$, we have $H^i(\tau^{\geq n}X) = 0$ for all~$i, n$, so Lemma~\ref{lem:cohomology is conservative on bounded objects} shows that~$\tau^{\geq n}X=0$ for all~$n$.
Since~$\cD$ is left complete, we deduce that~$X=0$, as required.
\end{proof}

\begin{remark}
There is an obvious analogue of 
Lemma~\ref{lem:cohomology is conservative on bounded above objects for left complete}
for right complete $t$-structures and objects that are bounded below.
\end{remark}

\subsubsection{\texorpdfstring{$t$}{t}-structures and Ind/Pro constructions}
For later reference, we recall
(from \cite[Prop.~2.13]{MR3935042} and~\cite[Lem.~C.2.4.3]{LurieSAG})
the following general facts
about $t$-structures and Ind-categories,
as well as the dual version for Pro-categories.

\begin{prop}
\label{prop:Ind t-structures}
Let $\cD$ be a stable $\infty$-category endowed with a $t$-structure.

\begin{enumerate}
\item
$\Ind \cD$ inherits a $t$-structure, characterized by the requirements
that $\cD  \hookrightarrow  \Ind \cD$ is $t$-exact,
and that the inclusions
$\cD^{\leq 0}  \hookrightarrow \cD$, 
$\cD^{\geq 0}  \into \cD$, 
and $\cD^{\heartsuit} \hookrightarrow \cD$
induce equivalences
$\Ind(\cD^{\leq 0})  \iso  \Ind(\cD)^{\leq  0},$
$\Ind(\cD^{\geq 0})  \iso  \Ind(\cD)^{\geq  0},$
$\Ind(\cD^{\heartsuit})  \iso  \Ind(\cD)^{\heartsuit},$
respectively.
If the $t$-structure on $\cD$ is furthermore bounded above,
then the $t$-structure  on $\Ind \cD$ is  right  complete.
\item
$\Pro \cD$ inherits a $t$-structure, characterized by the requirements
that $\cD  \hookrightarrow  \Pro \cD$ is $t$-exact,
and that the inclusions $\cD^{\leq 0}  \hookrightarrow \cD$, $\cD^{\geq 0}  \into \cD$, 
and $\cD^{\heartsuit} \hookrightarrow \cD$
induce equivalences
$\Pro(\cD^{\leq 0})  \iso  \Pro(\cD)^{\leq  0},$
$\Pro(\cD^{\geq 0})  \iso  \Pro(\cD)^{\geq  0},$
$\Pro(\cD^{\heartsuit})  \iso  \Pro(\cD)^{\heartsuit},$
respectively.
If the $t$-structure on $\cD$ is furthermore bounded below,
then the $t$-structure  on $\Pro \cD$ is  left complete.
\end{enumerate}
\end{prop}%
\begin{remark}
Since neither of the references cited above explicitly state the result about
the heart, we recall a proof here in the case of Ind categories; the Pro category
case is formally dual.  

On the one hand, the inclusion
$\cD^{\heartsuit} \hookrightarrow \cD$
induces a corresponding inclusion
$\Ind \cD^{\heartsuit} \hookrightarrow \Ind \cD$,
which lies in $(\Ind \cD)^{\heartsuit}$,
since it lies in both $(\Ind \cD)^{\leq 0} \iso \Ind (\cD^{\leq 0}) \hookrightarrow \Ind \cD$ and $ (\Ind \cD)^{\geq 0} \iso \Ind (\cD^{\geq 0}) \hookrightarrow \Ind \cD$.

On the other hand,
the functors
$\tau^{\leq 0}: \Ind \cD \to
(\Ind \cD)^{\leq 0} \iso \Ind (\cD^{\leq 0})$
and
$\tau^{\geq 0}: \Ind \cD \to (\Ind \cD)^{\geq 0} \iso \Ind (\cD^{\geq 0})$
are the Ind-extensions of the corresponding functors for~$\cD$ (as can be seen from the second displayed equation in the proof of~\cite[Prop.\ 2.13]{MR3935042}) and 
in particular commute with filtered colimits.
Thus, if $X$ is any object of~$(\Ind \cD)^{\heartsuit}$,
so that $\tau^{\geq 0}\tau^{\leq 0} X \iso X,$
then, if we write $X = \colim_i X_i$ for some objects $X_i$ of~$\cD$,
we find that
$X \iso \tau^{\geq 0}\tau^{\leq 0} X \iso \colim_i \tau^{\geq 0}\tau^{\leq 0} X_i,$
which describes $X$ as an object of~$\Ind \cD^{\heartsuit},$
so that $(\Ind \cD)^{\heartsuit} \subseteq \Ind \cD^{\heartsuit}.$
\end{remark}

\begin{example}
\label{ex:Ind/Pro objects with vanishing cohomology}
Let $\cC$ be the category of finitely generated modules over $k[\epsilon]$
(the dual numbers over a field~$k$), and let $D^b(\cC)$ be the bounded derived category of the abelian category~$\cC$, in the sense of Definition~\ref{def:derived category}
and Remark~\ref{rem:subcategories of the derived category}.
As usual, regard $k$ as a $k[\epsilon]$-module by having $\epsilon$ act by zero.
Recall then that $\Ext^{1}_{k[\epsilon]} (k,k)$
is one-dimensional over~$k$, and that if we let $x$ denote a basis vector,
then Yoneda product of $\Ext$s induces an isomorphism
$k[x] \iso \Ext^{\bullet}_{k[\epsilon]} (k,k)$.
Now consider the sequence
$$k \buildrel x \over \longrightarrow k[1]
\buildrel x \over \longrightarrow k[2]
\buildrel x \over \longrightarrow \dots
\buildrel x \over \longrightarrow k[n]
\buildrel x \over \longrightarrow \dots
$$
of morphisms in $D^b(\cC)$.
By the facts just recalled, the composite of any $n$ successive morphisms
in this 
sequence corresponds to the (non-zero!) element $x^n \in \Ext^n_{k[\epsilon]}(k,k).$
Thus none of the transition morphisms in this sequence are zero, and it induces
a non-zero object of
$\Ind D^b(\cC).$
This is then an example (probably the most standard example)
of 
an object of $\Ind D^b(\cC)$
which is non-zero (as one can see by computing its endomorphisms), but has all cohomologies equal to zero.
Note this object is bounded above, but not bounded below.

Similarly, if we instead consider the sequence 
$$
\dots \buildrel x \over \longrightarrow 
k[-n] \buildrel x \over \longrightarrow \dots
\buildrel x \over \longrightarrow k[-2]
\buildrel x \over \longrightarrow k[-1]
\buildrel x \over \longrightarrow k
$$
of morphisms in $D^b(\cC)$,
it gives rise to an object of $\Pro D^b(\cC)$ which is non-zero, but has all cohomologies
equal to zero; it is bounded below, but not bounded above.
\end{example}

\subsection{Derived categories}
We recall a range of notation and results related 
to derived categories of abelian categories, working in the framework of stable
$\infty$-categories.

\begin{df}
\label{def:derived category}
If $\cA$ is an abelian category, then we define the {\em derived category}
$D(\cA)$ of~$\cA$ to be the %
localization of the category
$\CoCh(\cA)$ 
of cochain complexes valued in~$\cA$ with respect to the class of
quasi-isomorphisms. %
In symbols, 
if we let $W$ denote the collection of quasi-isomorphisms in~$\CoCh(\cA)$,
then $D(\cA) \coloneqq  \CoCh(\cA)[W^{-1}].$
\end{df}%

\begin{remark}\label{rem:homotopy-category}
If we continue to let $W$ denote the collection of quasi-isomorphisms in~$\CoCh(\cA)$,
and let $W_0$ denote the collection of cochain homotopy equivalences,
then \cite[Prop.~1.3.4.5]{LurieHA} shows that $\CoCh(\cA)[W_0^{-1}] \iso
N_{\mathrm{dg}}\bigl(\CoCh(\cA)\bigr)$ (the dg-nerve of the dg-category $\CoCh(\cA)$).
Then 
\begin{equation}\label{localizing in stages}
\CoCh(\cA)[W^{-1}] \iso \CoCh(\cA)[W_0^{-1}][W^{-1}] \iso N_{\mathrm{dg}}
\bigl(\CoCh(\cA)\bigr)[W^{-1}].
\end{equation}
\end{remark}

\begin{remark}
Definition~\ref{def:derived category}
is the most general applicable definition of the $\infty$-categorical version
of the derived category of an abelian category of which we are aware. 
It is a stable $\infty$-category, by~\eqref{localizing in stages} and~\cite[Theorem~I.3.3]{MR3904731}, and 
the %
homotopy category of~$D(\cA)$ 
is the unbounded triangulated derived category of~$\cA$, as defined e.g.\ 
in~\cite[\href{https://stacks.math.columbia.edu/tag/05RU}{Tag 05RU}]{stacks-project}.  

If $\cA$ is a Grothendieck category, %
then an alternative definition of $D(\cA)$ is given in~\cite[Def.~1.3.5.8]{LurieHA}.
The equivalence of that definition with Definition~\ref{def:derived category}
is provided by~\cite[Prop.\ 1.3.5.15]{LurieHA}.
By~\cite[Prop.\ 1.3.5.21]{LurieHA}, we thus see that~$D(\cA)$ is presentable (hence cocomplete) whenever~$\cA$ is a Grothendieck category.
\end{remark}

\begin{remark}\label{rem:subcategories of the derived category}
If~$\cA$ is an abelian category,
then the usual $t$-structure on the unbounded triangulated derived category
of~$\cA$, i.e.\ the one given by the vanishing of cohomology groups (see e.g.\
\cite[\S IV.4]{MR1950475}),
induces a $t$-structure on~$D(\cA)$. 
We can then define $D^+(\cA)$, $D^-(\cA)$, and $D^b(\cA)$ in the usual manner.
Their homotopy categories are the corresponding triangulated subcategories of the unbounded triangulated derived category of~$\cA$.
Since $\infty$-categorical enhancements of these triangulated categories are unique up to equivalence~\cite[Corollary~6.5, Corollary~6.10]{Antieauuniqueness},
our definitions will coincide with most variants appearing in the literature, such as \cite[Def.~1.3.2.7]{LurieHA}, which gives an alternative definition of~$D^-(\cA)$
under the assumption that~$\cA$ admits enough projectives.
\end{remark}

The following definition uses the notion of a weak Serre subcategory, which was
recalled in Section~\ref{subsubsec:Serre subcategories}.
\begin{defn}\label{defn:D-sub-A-B}
  Let~$\cA$ be an abelian category, and let $\cB$ be a weak Serre subcategory
  of~$\cA$. Then we write~$D_{\cB}(\cA)$ for the full subcategory of $D(\cA)$
  consisting of those objects~$x$ all of whose cohomologies~$H^n(x)$ lie
  in~$\cB$. We define~$D^*_{\cB}(\cA)$ analogously, when $*$ is any of the boundedness
conditions $b, +,$ or~$-$.
\end{defn}

\begin{remark}
Recall that, by definition \cite[\S~1.2.11]{MR2522659},
(full) subcategories of an $\infty$-category correspond to (full) subcategories
of its %
homotopy category.  Under this correspondence,
the subcategory $D_{\cB}(\cA)$ of Definition~\ref{defn:D-sub-A-B}
corresponds to the subcategory of the triangulated derived category of~$\cA$
given 
by~\cite[\href{https://stacks.math.columbia.edu/tag/06UP}{Tag 06UP}]{stacks-project}.
\end{remark}

\begin{lem}\label{lem:quotient-of-bounded}Let~$\cA$ be an abelian category,
  and let~$\cB$ be a Serre subcategory. Then the natural functor
  $D^b(\cA)/D_{\cB}^{b}(\cA)\to D^{b}(\cA/\cB)$ is a $t$-exact equivalence.  
\end{lem}
\begin{proof}This is~\cite[Thm.\ 3.2]{MR1125707}.  
\end{proof}

Recall that a Serre subcategory~$\cB$ of an abelian category~$\cA$ is \emph{localizing} if the quotient functor
  $\cA\to\cA/\cB$ admits a right adjoint. This right adjoint functor is
  necessarily fully faithful~\cite[\S III.2 Prop.\ 3]{Gabrielthesis}.
\begin{lem}
  \label{lem:quotient-of-unbounded-with-adjoint}Let~$\cA$ be an abelian
  category, and let~$\cB$ be a localizing subcategory. Then the natural functor
  $ D(\cA)/D_{\cB}(\cA)\to D(\cA/\cB)$ is a $t$-exact equivalence.
\end{lem}
\begin{proof}This is~\cite[Lem.\ 5.9]{MR3398723}.  
\end{proof}

\subsubsection{A criterion for compact
  generation}\label{subsubsec:compact-derived}

By Lemma~\ref{lem:locally-coherent-implies-weak-Serre},  if~$\cA$ is a locally
coherent abelian category, then~$\cA^c$ is a weak Serre subcategory of~$\cA$, so
that~$D^b_{\cA^c}(\cA)$ is defined. Then we have the following lemma.
\begin{lem}%
  \label{lem:compact-objects-bounded-Serre-category}If~$\cA$ is a locally
  coherent abelian category, then the natural functors $D^b(\cA^c)\to
  D^b_{\cA^c}(\cA)$ and $D^-(\cA^c) \to D^-_{\cA^c}(\cA)$ are $t$-exact equivalences.
\end{lem}
\begin{proof}
  Bearing in mind that $\cA$ is (by definition) compactly generated, and so $\Ind(\cA^c) \isoto \cA$,
  the lemma is an immediate consequence of~\cite[Thm.\ 15.3.1]{MR2182076}. 
\end{proof}

If~$R$ is a coherent ring, the category~$\Mod(R)$ is locally coherent, and its compact objects coincide with the finitely presented modules.
We introduce the notation 
\[
D^b_{\fp}(R)\coloneqq D^b_{\Mod^{\fp}(R)}(\Mod(R))
\]
and
\[
D^-_{\fp}(R)\coloneqq D^-_{\Mod^{\fp}(R)}(\Mod(R)),
\]
and note the following immediate corollary of Lemma~\ref{lem:compact-objects-bounded-Serre-category}.

\begin{cor}\label{finitely presented derived category for coherent ring}
Let~$R$ be a coherent ring.
Then the natural functors 
$D^b(\Mod^{\fp}(R)) \to D^b_{\fp}(R)$
and $D^-(\Mod^{\fp}(R)) \to D^-_{\fp}(R)$
are $t$-exact equivalences.
\end{cor}

The following proposition is closely related to~\cite[Lem.\
4.5]{MR2157133}. By ``a set of weak generators'' of an abelian category~$\cA$,
we mean a set~$X$ of objects of~$\cA$ having the property that for each~$a\in \cA$, there
exists an~$x\in X$ and a non-zero morphism $x\to a$.  %

\begin{prop}\label{prop:abelian-compact-implies-derived-compactly-generated}%
Let $\cA$ be a locally coherent Grothendieck category, and suppose furthermore
that $\cA^c \subseteq D(\cA)^c$.  Then:

\begin{enumerate}
\item $D(\cA)$ is compactly generated, and $D^b_{\cA^c}(\cA)$ coincides with $D(\cA)^c$.
\item Any subset $X\subseteq \cA^{c}$ which is a set of weak generators of~$\cA$ is a set of
  compact generators of~$D(\cA)$.
  \item The standard $t$-structure on~$D(\cA)$ is the one induced by the standard
  $t$-structure on $D^b_{\cA^c}(\cA) $ via Proposition~\ref{prop:Ind t-structures} and the equivalence 
  $\Ind D^b_{\cA^c}(\cA) \iso D(\cA)$ induced by part~(1).  
\end{enumerate}
\end{prop}
\begin{proof}
By Lemma~\ref{lem:compact-objects-bounded-Serre-category}, 
the full sub-$\infty$-category $D^b_{\cA^c}(\cA)$
of $D(\cA)$ is well-defined, and equivalent to $D^b(\cA^c)$.
Furthermore, the assumption that $\cA^c \subset D(\cA)^c$, together with an induction on amplitude, shows that 
$D^b_{\cA^c}(\cA) \subset D(\cA)^c$.
Since $D^b_{\cA^c}(\cA)$ is evidently idempotent complete,
we see that part~(1) of the proposition will follow
from Lemma~\ref{lem:compact generation}~(2), applied to $\cC' := D^b_{\cA^c}(\cA)$ and~$\cC := D(\cA)$,
provided that we prove that $\bigl(D^b_{\cA^c}(\cA)\bigr)^{\perp} = 0.$
To prove this, it suffices to show that if $a \in D(\cA)$ is not zero, and $X \subseteq \cA^c$ is a set of weak generators of~$\cA$, then there exists $x \in X$ such that
$\RHom_{D(\cA)}(x, a) \ne 0$; note that this will also establish part~(2).

Suppose then that $a$ is a non-zero object of $D(\cA)$.   %
Since the homotopy category of~$D(\cA)$ is the unbounded derived category of~$\cA$ with the usual $t$-structure,  
there exists
some~$n$ such that $H^n(a) \neq 0$. %
We may represent the truncation $\tau^{\leq n}a$ by a complex $ \cdots \to a_{n-1} \to a_n \to 0 \to 0 \to \cdots$ in~$\CoCh(\cA)$.
Since $\coker(a_{n-1} \to a_n) \eqcolon H^n(\tau^{\leq n} a) =H^n(a)\neq 0,$
and since 
$X$ is a set of weak generators of~$\cA$ by assumption,
we may find an object $x$ of $X$ and a morphism $x \to a_n$ such that the composite
$x\to H^n(\tau^{\leq n} a)$ is non-zero.
Then the composite $x[-n] \to \tau^{\leq n} a \to a$ (the first arrow being induced
by the chosen morphism $x \to a_n$) is non-zero (since it induces a non-zero morphism
on~$H^n$), and so in particular $\RHom_{D(\cA)}(x, a) \neq 0,$ as required. 

Finally, the claim in part~(3) about $t$-structures is equivalent to the $t$-structure
on~$D(\cA)$ being compatible with filtered colimits, which is~\cite[Prop.\ 1.3.5.21]{LurieHA}.
\end{proof}

\subsubsection{Derived functors}
\label{subsubsec:derived functors}
A very general mechanism for constructing left derived functors is given by the following 
result~\cite[Thm.~1.3.3.2]{LurieHA},
which characterizes $D^-(\cC)$ (for abelian categories with enough projectives)
by a mapping property. %

\begin{thm}
\label{thm:left derived functors}
If $\cC$ is an abelian category with enough projectives,
and if $\cD$ is a stable  $\infty$-category equipped with
a left complete $t$-structure,
then $F \mapsto \tau^{\geq 0} F_{|\cC}$ {\em (}the restriction
being taken by identifying $\cC$ with the heart of~$D^-(\cC)${\em )}
is an equivalence between the $\infty$-category of right $t$-exact
functors $D^-(\cC)\to\cD$ which carry projective objects of~$\cC$ into $\cD^{\heartsuit},$
and the ordinary category of right exact functors $\cC\to \cD^{\heartsuit}$.

Furthermore, by restriction from~$D^-(\cC)$ to~$D^b(\cC)$, these categories are equivalent to the $\infty$-category of right $t$-exact
functors $D^b(\cC)\to\cD$ which carry projective objects of~$\cC$ into $\cD^{\heartsuit}.$
\end{thm}%
\begin{proof}
As already noted, the 
first claim is simply a restatement of~\cite[Thm.~1.3.3.2]{LurieHA}.
The second claim is a consequence of the proof of that result, and is a manifestation
of the role that left completeness plays in that result, as we now explain.

Namely, suppose that $F:D^-(\cC) \to \cD$ is a right $t$-exact functor.
Write $D^{-,\geq n}(\cC)$ to denote the full subcategory of $D^-(\cC)$
consisting of objects whose cohomology in degrees $<n$ vanishes (i.e.\ the ``$\geq n$''
part of the $t$-structure on~$D^-(\cC)$),
and write
$$F_n\coloneqq  \tau^{\geq n} F_{| D^{-,\geq n}(\cC)}: 
D^{-,\geq n}(\cC) \to \cD^{\geq n}.$$
The fact that $F$ is right $t$-exact shows that the canonical functor
\begin{equation}
\label{eqn:truncation iso}
\tau^{\geq n} F \to \tau^{\geq n} F \tau^{\geq n}
\end{equation}
is an isomorphism, from which one easily verifies that the diagram 
of functors
$$\xymatrix{
D^{-,\geq n-1}(\cC) \ar^-{F_{n-1}}[r] \ar^-{\tau^{\geq n}}[d] &
\cD^{\geq n-1} \ar^-{\tau^{\geq n}}[d] \\
D^{-,\geq n}(\cC) \ar^-{F_n}[r] & \cD^{\geq n}
}
$$
commutes. 
Thus, taking into account the left completeness of each of $D^-(\cC)$
(by~\cite[Prop.~1.3.3.16]{LurieHA}) and $\cD$ (by hypothesis),
we may form the functor 
$$\lim_n F_n: D^-(\cC) \to \cD.$$
Again, the fact that~\eqref{eqn:truncation iso} is an isomorphism
shows that we have natural transformations (for $X$ an object of $D^-(\cC)$)
$$F(X) \to \tau^{\geq n} F(X) \iso \tau^{\geq n}F (\tau^{\geq n} X) = F_n(\tau^{\geq n} X),$$
evidently compatible as $n$ varies, and hence a natural transformation
\begin{equation}
\label{eqn:limit transformation}
F(X) \to \lim_n \bigl ( F_n(\tau^{\geq n}X)\bigr) = (\lim_n F_n) (X). 
\end{equation}

If $F$ furthermore takes projective objects of $\cC$ to objects of~$\cD^{\heartsuit}$,
then~\eqref{eqn:limit transformation} is a natural isomorphism.
(Indeed, both source and target are right $t$-exact functors
from $D^-(\cC)$ to $\cD$ which coincide on~$\cC$, and take projective objects
of~$\cC$ to~$\cD^{\heartsuit}$.)
Thus $F$ can be recovered (and constructed) from the various~$F_n$, and thus from its
restriction to~$D^b(\cC)$.
\end{proof}

\begin{cor}
\label{cor:left derived functors exact case}
If $\cC$ is an abelian category with enough projectives,
and if $\cD$ is a stable  $\infty$-category equipped with
a left complete $t$-structure,
then $F \mapsto F_{|\cC}$ 
induces an equivalence between the $\infty$-category of  $t$-exact
functors $D^-(\cC)\to\cD$ %
and the ordinary category of  exact functors $\cC\to \cD^{\heartsuit}$.
\end{cor}
\begin{proof}
This follows from  Theorem~\ref{thm:left derived functors}, noting that  the $t$-exactness assumption implies that {\em all} objects
of~$\cC$, projective or not, are carried into
objects of~$\cD^{\heartsuit}$. See also ~\cite[Rem.~1.3.3.6]{LurieHA}.
\end{proof}%

In the case of $t$-exact functors, we can in fact state
results  %
valid for any abelian category (related, morally,
to the fact that we don't
need projective resolutions to compute the derived functors
of exact functors).
This is because of the following very general mapping property satisfied by
bounded derived categories~\cite[Cor.~7.4.12]{bunke2024controlledobjectsleftexactinftycategories}.

\begin{thm}
\label{thm:exact derived functors}
If $\cC$ is an abelian category, and 
$\cD$ is a stable $\infty$-category,
then $F \mapsto F_{|\cC}$ induces an equivalence
between the $\infty$-category of exact functors 
$D^b(\cC) \to \cD$
and the $\infty$-category of finite coproduct-preserving %
functors $\cC \to \cD$
which take exact sequences in~$\cC$ to fibre sequences in~$\cD$.
\end{thm}

This has the following immediate corollary:

\begin{cor}
\label{cor:t-exact derived functors}
If $\cC$ is an abelian category, 
and if $\cD$ is a stable  $\infty$-category equipped with
a $t$-structure,
then there
is an equivalence between the $\infty$-category of $t$-exact
functors $D^b(\cC)\to\cD$ 
and the ordinary category of exact functors $\cC\to\cD^{\heartsuit}$, given by $F \mapsto F_{|\cC}$ {\em (}the restriction
being taken by identifying $\cC$ with the heart of~$D^b(\cC)${\em )}.
\end{cor}

We also note the following variant of the preceding corollary.

\begin{cor}
\label{cor:derived functors of exact functors}
If $\cC$ is an abelian category, 
if $\cD$ is a stable $\infty$-category equipped with
a $t$-structure,
and if $F: D^b(\cC) \to \cD$ is an exact functor for which $F_{|\cC}$ 
takes values in $\cD^{\heartsuit}$,
then $F$ is $t$-exact. 
\end{cor}
\begin{proof}
The hypotheses imply that $F_{|\cC}$ is an exact functor $\cC \to \cD^{\heartsuit}$,
which by Corollary~\ref{cor:t-exact derived functors} arises by restriction
from a $t$-exact functor $G: D^b(\cC) \to \cD$.
Taking into account the equivalence of categories described in the 
statement of Theorem~\ref{thm:exact derived functors},
the equality $F_{|\cC} = G_{|\cC}$ is induced by an isomorphism
$F \iso G$.  Thus $F$ is $t$-exact, as claimed.
\end{proof}

We now consider right derived functors. We begin with the following simple lemma.
\begin{lem}\label{lem:exact-functor-induces-t-exact}
  If~$F:\cA\to\cB$ is an exact functor between  abelian categories, then~$F$
  induces a $t$-exact functor $D(\cA)\to D(\cB)$, whose restriction to the heart
  of~$D(\cA)$ is~$F$.
\end{lem}
\begin{proof}We obtain a functor $\CoCh(\cA)\to\CoCh(\cB)$ by applying~$F$
  degreewise. Since~$F$ is exact, it preserves quasi-isomorphisms, so it induces
  the required functor $D(\cA)\to D(\cB)$ by the universal property of
  localization.
\end{proof}%

Suppose now that~$F:\cA\to\cB$ is an additive functor between Grothendieck
categories. Then~$F$ induces a functor $\CoCh(\cA)\to\CoCh(\cB)$, and
thus a functor $\overline{F}:\CoCh(\cA)\to D(\cB)$.
Write~$Q:\CoCh(\cA)\to D(\cA)$ for the localization
functor. Following~\cite[7.5.23]{MR3931682}
we make the following definition,
which is a ``lifting'' to the $\infty$-categorical setting of a definition
on the level of homotopy categories
that goes back
to Quillen %
(see~\cite[Def.~1.4.1]{MR223432},
or~\cite[Def.~2.3.1]{MR3931682} for a more recent treatment).

\begin{defn}\label{def:right derived functor}
Let $F: \cA \to \cB$ be an additive functor between Grothendieck categories.
A \emph{right derived functor~$RF$}
of~$F$ is a functor $RF:D(\cA)\to D(\cB)$ which is equipped with a natural
transformation $\overline{F}\to RF\circ Q$, %
and is such that for any other functor
$G:D(\cA)\to D(\cB)$ together with a natural transformation
$\overline{F}\to G\circ Q$, there is a unique (up to a contractible space of
choices) natural transformation $RF\to G$ giving rise to the given
$\overline{F}\to G\circ Q$.
\end{defn}

\begin{lem}%
  \label{lem:right-derived-of-exact-is-exact}If~$F:\cA\to\cB$ is an exact
  functor between Grothendieck
categories, then the right derived functor~$RF$ exists.
Furthermore $RF$ is $t$-exact, and agrees with
the extension of~$F$ in Lemma~{\em \ref{lem:exact-functor-induces-t-exact}}.
\end{lem}
\begin{proof}
  This is immediate from the definitions.
\end{proof}

In fact, 
the right derived functor~$RF$ exists for any additive functor~$F$, 
although it is only well-behaved if~$F$ is left exact, in a sense made precise by Lemma~\ref{lem:right-derived-left-exact} below.                                                          
Right derived functors can be computed using resolutions by $K$-injective complexes of injective objects:
see~\cite[\href{https://stacks.math.columbia.edu/tag/079P}{Tag
079P}]{stacks-project} and
\cite[\href{https://stacks.math.columbia.edu/tag/070K}{Tag
  070K}]{stacks-project} in the context of triangulated categories.
Alternatively, in the context of $\infty$-categories and Definition~\ref{def:right derived functor},
we can give $\CoCh(\cA)$ the model category structure of~\cite[Proposition~1.3.5.3]{LurieHA}.
The fibrant objects are then precisely the $K$-injective complexes of injective objects: compare~\cite[Remark~2.3.18]{Hoveybook} for this fact, 
bearing in mind that the ``dg-injective complexes" in \emph{loc.\ cit.} are precisely the $K$-injective complexes of injective objects, 
by~\cite[1.2.I]{AvramovFoxby}
(where $K$-injective complexes are called $\pi$-injective).
Then, as explained in~\cite[Section~7.5.25]{MR3931682}, the existence of~$RF$, and the fact that it can be computed by $K$-injective complexes of injective objects,
both follow from the fact that $F$ sends weak equivalences between fibrant objects of $\CoCh(\cA)$ to weak equivalences in $\CoCh(\cB)$, which in turn follows from~\cite[Proposition~1.3.5.14]{LurieHA},
asserting that every weak equivalence between fibrant objects of $\CoCh(\cA)$ is a chain homotopy equivalence. %

\begin{lem}
  \label{lem:right-derived-left-exact}Let ~$F:\cA\to\cB$ be an additive functor between Grothendieck
categories. Then the right derived functor~$RF$ is left $t$-exact. Furthermore, the natural
map $F\to H^0RF$ \emph{(}induced by the natural transformation
$\overline{F}\to RF\circ Q$\emph{)}  is an isomorphism if and only if~$F$ is left exact.  %
\end{lem}%
\begin{proof}
  Both claims follow easily from the existence of injective
  resolutions, and the fact that bounded below complexes of injective objects are $K$-injective, see e.g.\
  \cite[\href{https://stacks.math.columbia.edu/tag/070J}{Tag 070J}]{stacks-project} and
  \cite[\href{https://stacks.math.columbia.edu/tag/05TD}{Tag 05TD}]{stacks-project}. %
  \end{proof}

We next quote two results from~\cite{emerton2023introduction} concerning the interplay between derived functors and adjunctions.

\begin{lem}\label{lem:derived-adjunction-unbounded-Grothendieck}%
  Let $F:\cA\to\cB$ be
  an exact, colimit-preserving functor between Grothendieck categories,
  and write $G:\cB\to\cA$ for its right adjoint. Then the right-derived functor
  $RG: D(\cB)\to D(\cA)$ is right adjoint to $F:D(\cA)\to D(\cB)$. In
  particular, $F:D(\cA)\to D(\cB)$ is continuous.
\end{lem}
\begin{proof}
This is~\cite[Prop.\
A.7.1]{emerton2023introduction}.
\end{proof}

\begin{thm}
  \label{thm:unbounded-adjoints-localizing}
Let $F:\cA \hookrightarrow \cB$
be the inclusion of a localizing  subcategory into a Grothendieck
category, with right adjoint $G:\cB \to \cA$. %
 Suppose further that:
\begin{enumerate}
\item\label{item:weaker-version-of-preserving-injectives} For any objects $X$ and $Y$ of $\cA$ with $Y$ injective,
there is an epimorphism $Z \to X$ %
such that
$\Ext^i_{\cB}(F(Z),F(Y)) = 0$ for $i > 0$. %
\item The formation of products in~$\cB$ is exact.
\item The derived right adjoint $RG$ has finite cohomological dimension.
\end{enumerate}
 Then $F:D(\cA) \to D(\cB)$ is fully faithful, 
with essential image equal to $D_{\cA}(\cB)$.
\end{thm}
\begin{rem}
 Note in
particular that condition~\eqref{item:weaker-version-of-preserving-injectives}
is satisfied if~$F$ preserves injectives.
\end{rem}
\begin{proof}[Proof of Theorem~{\em \ref{thm:unbounded-adjoints-localizing}}]
  This is a combination of parts~(1) and~(2) of~\cite[Prop.\
  A.7.3]{emerton2023introduction}. %
\end{proof}

We end this section with a result
on the compatibility
of localization of locally Noetherian categories with passage to derived
categories. Rather than consider the general case, we will assume various
additional hypotheses that will hold in our application of this material in Section~\ref{subsec:derived-categories-smooth-representations}.

Suppose  that~$\cB$ is
a locally Noetherian abelian category and that ~$i_{*}:\cA\to\cB$ is an
inclusion of a localizing subcategory;  so~$i_{*}$ is exact and is compatible
with colimits, and by Lemma~\ref{lem:derived-adjunction-unbounded-Grothendieck},
we have a continuous $t$-exact extension $i_{*}:D(\cA)\to D(\cB)$. 
Write
$j^{*}:\cB \to\cB/\cA$ for the quotient functor, which is also exact and compatible
with colimits, and $j_{*}:\cB/\cA\to\cB$ for its fully faithful right adjoint,
which (as recalled in Appendix~\ref{subsec:1-categories}) preserves filtered colimits.  If we suppose further that~$j_{*}$ is exact, it is
furthermore compatible with all colimits, and so we have continuous $t$-exact
extensions $j^{*}:D(\cB)\to D(\cB/\cA)$, $j_{*}: D(\cB/\cA)\to D(\cB)$.

\begin{prop}%
  \label{prop:localizing-derived-compact-objects}Suppose as above that~$\cB$ is
  a locally Noetherian abelian category and that ~$i_{*}:\cA\to\cB$ is an
  inclusion of a localizing subcategory. Write $j^{*}:\cB \to\cB/\cA$ for the
  quotient functor, and $j_{*}:\cB/\cA\to\cB$ for its right adjoint.

  Suppose that: 
  \begin{enumerate}[label=(\alph*)]
  \item $\cB^c \subseteq D(\cB)^c$.
  \item $j_{*}:\cB/\cA\to\cB$ is exact.
  \item  $i_{*}:D(\cA)\to D(\cB)$ is fully faithful, with essential image
    $D_{\cA}(\cB)$. %
  \end{enumerate}
  Then 
  \begin{enumerate}
   \item \label{item:13}$j^{*}$ induces equivalences $D(\cB)/D(\cA)\isoto D(\cB/\cA)$ and $D^b(\cB^c)/D^b(\cA^c)\isoto D^b(\cB^c/\cA^c)$.
    \item \label{item:10} $j_{*}:D(\cB/\cA)\to D(\cB)$ is fully faithful, and is right adjoint
    to~$j^{*}:D(\cB)\to D(\cB/\cA)$.
  \item\label{item:11} $D(\cB)$, $D(\cA)$ and $D(\cB/\cA)$ are compactly
    generated.
  \item\label{item:12} $D(\cB)^c$ coincides with $D^b_{\cB^c}(\cB)$, which in turn is
    equivalent to $D^b(\cB^{c})$. The analogous statements hold for
    $D(\cA)^{c}$ and $D(\cB/\cA)^{c}$.
  \end{enumerate}
\end{prop}
\begin{proof}%
  By assumption we have $\cB^c \subseteq D(\cB)^c$.  Since~$i_{*}:\cA\to\cB$
  preserves colimits and compact objects, and (by assumption)
  $i_{*}:D(\cA)\to D(\cB)$ is continuous and fully faithful, we see also that
  $\cA^c \subseteq D(\cA)^c$.  Then the statements of \eqref{item:11} and
  \eqref{item:12} for~$D(\cB)$ and~$D(\cA)$ follow from
  Proposition~\ref{prop:abelian-compact-implies-derived-compactly-generated} and
  Lemma~\ref{lem:compact-objects-bounded-Serre-category}; furthermore
  $i_{*}:D(\cA)\to D(\cB)$ preserves compact objects.

  We now prove part~\eqref{item:13}.
  By assumption, $i_{*}$ induces an equivalence between $D(\cA)$ and
  $D_{\cA}(\cB)$, hence the first statement of part~\eqref{item:13} follows from
  Lemma~\ref{lem:quotient-of-unbounded-with-adjoint};
  note that this implies that Hypothesis~\ref{hyp:usual hyp for semiorthogonal} is satisfied by~$i_*$ and~$j^*$.
  The second statement of part~\eqref{item:13} follows similarly from Lemma~\ref{lem:quotient-of-bounded}, provided we show
  that~$i_{*}$ induces an equivalence between $D^b(\cA^c)$ and
  $D^b_{\cA^c}(\cB^c)$. 
  Certainly~$i_{*}$ induces a fully faithful functor
  $i_{*}:D^b(\cA^c)\to D^b_{\cA^c}(\cB^c)$.
  It thus suffices to note that any object of
  $D^b_{\cA^c}(\cB^c)$ is compact in~$D(\cB)$ (by the already proved statement of part~\eqref{item:12} for~$D(\cB)$) and contained in~$D_{\cA}(\cB)$, and thus in the essential
  image of~$D(\cA)$; and $D(\cB)^c\cap D(\cA)=D(\cA)^c=D^b(\cA^c)$ by
  Lemma~\ref{lem: abstract semiorthogonal decomposition} (and the already proved statement of part~\eqref{item:12} for~$D(\cA)$).
  This concludes the proof of part~\eqref{item:13}.

  Finally, by Lemmas~\ref{lem:derived-adjunction-unbounded-Grothendieck}
  and~\ref{lem:right-derived-of-exact-is-exact}, $j_{*}:D(\cB/\cA)\to D(\cB)$ is
  right adjoint to~$j^{*}:D(\cB)\to D(\cB/\cA)$.
  Hence the statements of part~\eqref{item:10}, as well as part~\eqref{item:11} for $D(\cB/\cA)$, are immediate consequences of 
  Lemma~\ref{lem: abstract semiorthogonal decomposition}, which also shows that~$j^*$ preserves compact objects.
  Thus the objects of $(\cB/\cA)^c = j^*(\cB^c)$ are compact in $D(\cB/\cA)^c$,
  and so the statement of part~\eqref{item:12} for~$D(\cB/\cA)$ follows from 
  Proposition~\ref{prop:abelian-compact-implies-derived-compactly-generated} and Lemma~\ref{lem:compact-objects-bounded-Serre-category}.\qedhere
\end{proof}

\subsubsection{Colimits of abelian and derived categories}
\label{subsec:colimits}
If $\{\cA_i\}_{i \in \cI}$ is a filtered system of abelian categories indexed
by~$\cI$, with exact transition functors,
then we can form $\cA \coloneqq  \colim_i \cA_i,$
which is again an abelian category.
By Lemma~\ref{lem:exact-functor-induces-t-exact}
we get an induced system $\{ D^b(\cA_i) \}_{i \in \cI}$
with $t$-exact transition functors,
and the universal property of Corollary~\ref{cor:t-exact derived functors}
shows that 
$$ \colim_i D^b(\cA_i) \iso D^b(\cA),$$
the colimit being formed in the $\infty$-category of small %
$\infty$-categories. %

Similarly, a consideration of the universal property of the Ind construction shows that
we obtain an equivalence
$$\colim_i \Ind D^b(\cA_i) \iso \Ind D^b (\cA),$$
the colimit now being formed in the $\infty$-category of compactly generated stable
$\infty$-categories (whose morphisms are the continuous functors). 

\subsection{Left derived functors via Pro-categories}
\label{subsec:left derived}
The goal of this section is to study how certain constructions of derived
functors interact with the formation of Pro-categories.
The reason for these considerations is as follows: if $\cC$ is any
abelian category, then $\Pro \cC$ has enough projectives, and so passing
to $\Pro \cC$ can be a convenient first step in the construction of left derived functors.
On the other hand, these derived functors are then {\em a priori} defined
on $D^-(\Pro \cC)$, which is an inconvenient source category for the applications 
we have in mind; we would prefer to work with $\Pro D^b(\cC)$.  
More precisely, then, the following discussion will explain how we can 
form derived functors whose source is this latter category.

Another motivation for passing to $\Pro \cC$ is to construct left adjoints which
might otherwise not exist.    
With this in mind,
we also show that our construction of derived functors is compatible with the formation
of such adjoints.

Finally,
we note that our discussion is closely related to the discussion
of~\cite[Ch.~15]{MR2182076},
which considers the analogous constructions for Ind (rather than Pro) categories,
in the setting of triangulated (rather than stable $\infty$-) categories.  
For a more precise explanation of the connection,
see Remark~\ref{rem:KS} below.

\subsubsection{The Pro-categorical context}
Throughout this discussion, we let $\cC$ and $\cC'$ denote
two small abelian categories. We begin with the following consequence of Theorem~\ref{thm:exact derived functors}.
\begin{cor}
\label{cor:pro t-exact derived functors}
If $\cC$ and $\cC'$ are two small abelian categories,
then an exact functor $D^b(\cC) \to \Pro D^b(\cC')$ is $t$-exact if and only if its limit-preserving extension $\Pro D^b(\cC) \to \Pro D^b(\cC')$ is $t$-exact.
Furthermore, the functors $F \mapsto F_{|\Pro \cC}$ and $F\mapsto F_{|\cC}$ induce equivalences between the following categories:
\begin{enumerate}
\item The $\infty$-category of $t$-exact
  limit-preserving functors $\Pro D^b(\cC)\to \Pro D^b(\cC')$.
  \item The $\infty$-category of $t$-exact
 functors $D^b(\cC)\to \Pro D^b(\cC')$.
\item The ordinary category of exact limit-preserving
functors from $\Pro \cC$ to~$\Pro \cC'$.
\item The ordinary category of exact functors
$\cC \to \Pro \cC'$.
\end{enumerate}
\end{cor}
\begin{proof}
The first statement of the corollary is immediate from Proposition~\ref{prop:Ind t-structures}, which also implies that the $\infty$-category of $t$-exact
limit-preserving functors $\Pro D^b(\cC)\to \Pro D^b(\cC')$ is equivalent to the
category of $t$-exact functors $D^b(\cC)\to \Pro D^b(\cC')$. The result is then immediate from Corollary~\ref{cor:t-exact derived functors}.
\end{proof}

The category $\Pro \cC$ has enough projectives, by Lemma~\ref{lem:quasi-projectives-in-Pro-cat},
while $\Pro D^b(\cC')$ is 
left complete with respect to its natural $t$-structure, by
Proposition~\ref{prop:Ind t-structures}. 
Thus we have the following particular case of Theorem~\ref{thm:left derived functors}
and Corollary~\ref{cor:left derived functors exact case}.

\begin{prop} 
\label{prop:left derived functors}
\leavevmode
\begin{enumerate}
\item\label{item:80} 
$F \mapsto \tau^{\geq 0} F_{|\Pro \cC}$
is an equivalence between the $\infty$-category of right $t$-exact
functors $D^-(\Pro \cC)\to \Pro D^b(\cC')$
which carry projective objects of~$\Pro \cC$ into $\Pro \cC'$,
and the ordinary category of right exact functors $\Pro \cC\to \Pro \cC'$.
\item\label{item:61} 
$F \mapsto F_{|\Pro \cC}$
is an equivalence between the $\infty$-category of $t$-exact
functors $D^-(\Pro \cC)\to \Pro D^b(\cC')$
and the ordinary category of exact functors $\Pro \cC\to \Pro \cC'$.
\end{enumerate}
\end{prop}

By Proposition~\ref{prop:left derived functors}~\eqref{item:61}, the natural identification $\Pro \cC \iso \bigl(\Pro D^b(\cC)\bigr)^{\heartsuit}$
induces  a $t$-exact functor
$p:D^-(\Pro\cC) \to \Pro D^b(\cC).$ 
Of course, the canonical inclusion $\cC \hookrightarrow \Pro \cC$
(which is exact)
induces  $t$-exact functors $i:D^b(\cC) \to D^-(\Pro \cC)$ %
and $j:D^-(\cC)\to D^-(\Pro\cC)$.
There are also the canonical inclusions
$D^b(\cC)\hookrightarrow D^-(\cC)$ and
$D^b(\cC) \hookrightarrow \Pro D^b(\cC)$.

\begin{lemma}
\label{lem:checking Pro comparison commutatitivy}We have a commutative diagram of $t$-exact functors:
\begin{equation}
\label{eqn:Pro comparison}
\begin{tikzcd}[column sep=large,row sep=large]
  & D^-(\mathcal{C}) \ar[d,"j"]  
\\
  D^b(\mathcal{C}) \ar[r,"i"] \ar[dr] \ar[ur]
    & D^-(\mathrm{Pro}\,\mathcal{C}) \ar[d,"p"] 
\\
  & \mathrm{Pro}\,D^b(\mathcal{C})
\end{tikzcd}
\end{equation}
The horizontal and diagonal arrows in~\eqref{eqn:Pro comparison} are fully faithful.
\end{lemma}%
\begin{proof}%
The commutativity follows directly from Corollary~\ref{cor:t-exact derived functors},
since all the functors are $t$-exact, and their restrictions to the hearts are either
the identity functors or the %
natural embedding $\cC \hookrightarrow \Pro \cC$.
The full faithfulness of the diagonal arrows holds by definition.
The full faithfulness of~$i$ is a consequence of~\cite[Thm.\ 15.3.1]{MR2182076}.
\end{proof}

As was already remarked on above,
in applications, it is much more convenient to 
have the domain of our left-derived functors be equal to~$\Pro D^b(\cC)$,
rather than the $\infty$-category~$D^-(\Pro \cC)$, which is rather hard to work with.
In other words, under appropriate hypotheses, we would like
to canonically factor the left derived functors constructed
by Proposition~\ref{prop:left derived functors} through the
functor~$p$ of~\eqref{eqn:Pro comparison}. %
It turns out that this is possible whenever the functor we are deriving is cofiltered limit-preserving, i.e.\ it
  is  Pro-extended from~ $\cC$. %
  
To this end, suppose that $f:\cC\to\Pro(\cC')$ is right exact, so that $\Pro(f):\Pro(\cC)\to\Pro(\cC')$ is right exact and cofiltered limit-preserving.
Write  
$$F:D^-(\Pro \cC) \to \Pro D^b(\cC')$$
for the right $t$-exact functor corresponding to $\Pro(f)$ via Proposition~\ref{prop:left derived functors}~\eqref{item:80}.
We will also consider the composite
\[
(F\circ i)_{|\cC} : \cC \subset D^b(\cC) \xrightarrow{i} D^-(\Pro \cC) \xrightarrow{F} \Pro D^b(\cC');
\]
recalling from Remark~\ref{rem:comparison of Ind-completions} that~$\Pro\cC$ is naturally identified with the $\infty$-categorical $\Pro$-completion of~$\cC$, 
this functor can be $\Pro$-extended to a functor
\[
\Pro((F\circ i)_{|\cC}) : \Pro\cC \to \Pro D^b(\cC').
\]

\begin{lemma}
  \label{lem:truncated cofiltered limit-preserving}
  With the notation of the preceding paragraph, the functor $F_{| \Pro \cC} : \Pro \cC \to \Pro D^b(\cC')$ 
preserves cofiltered limits; equivalently, it is naturally equivalent to $\Pro((F\circ i)_{|\cC})$.
\end{lemma}%
\begin{proof}
Set
\[
G \coloneqq \Pro\bigl( (F\circ i)_{|\cC} \bigr) : \Pro\cC \longrightarrow \Pro D^b(\cC').
\]
By the Pro-version of the adjunction remarked upon in
Remark~\ref{rem:Indization as an adjunction},
there is a canonical natural transformation
\[
\eta : F_{|\Pro\cC} \longrightarrow G,
\]which we need to prove is an isomorphism.
Concretely, if $X = \quoteslim{i} X_i$ with the $X_i$ objects
of~$\cC$, then $G(X) \coloneqq  \quoteslim{i} F(X_i)$; and so the morphisms
$F(X) \to F(X_i)$ induce a morphism $\eta_X: F(X) \to G(X)$.
From this description,
we see that~$\eta_X$ is an isomorphism whenever~$X \in \cC$.

Since the $t$-structure on $\Pro D^b(\cC')$ is left complete,
it suffices to prove for every $n\le 0$, the truncation
\begin{equation}\label{eqn:nth truncation}
\tau^{\ge n}\eta : \tau^{\ge n}F_{|\Pro\cC} \longrightarrow \tau^{\ge n}G
\end{equation}
is an equivalence.
When~$n=0$, we have $\Pro(f)=\tau^{\geq 0}F_{| \Pro \cC} \iso \tau^{\geq 0}G,$ so by induction it suffices to show that if~\eqref{eqn:nth truncation} is an isomorphism for some~$n\le 0$, then it is also an isomorphism for~$n-1$.
Considering the fibre sequences \[
H^{n-1}(X)[1-n] \longrightarrow  \tau^{\ge n-1}X \longrightarrow \tau^{\ge n}X,
\]we see that it suffices in turn to show that the induced morphism \begin{equation}\label{eq:cuipx67a7m}H^{n-1}(F(X))\to H^{n-1}(G(X))\end{equation}
is an isomorphism for all objects~$X$ of~$\Pro\cC$.
We already know that it is an isomorphism when~$X \in \cC$.

Suppose that we have a short exact sequence
\[
  0\to Y\to P\to X\to 0
\]in~$\Pro\cC$, where~$P$ is projective.
Bearing in mind that~$F$ sends projective objects of~$\Pro \cC$ into $\Pro \cC'$, so that in particular $H^{n-1}(F(P)) = 0 $,
we obtain a diagram of %
exact sequences in~$\Pro(\cC')$
$$ \xymatrix{
  0 \ar[r] \ar[d] & H^{n-1}\bigl( F(X) \bigr) \ar[r]\ar[d] &
  H^n\bigl( F(Y)\bigr) \ar^-{\sim}[d] \\
  H^{n-1}\bigl( G(P)\bigr) \ar[r] & H^{n-1}\bigl( G(X) \bigr) \ar[r] &
  H^n\bigl( G(Y)\bigr)
}
$$
in which the last vertical map is an isomorphism by inductive assumption.
Since~$\Pro\cC$ has enough projectives, it follows that in order to show that~\eqref{eq:cuipx67a7m} is an isomorphism for all~$X$, it is enough to show that $H^{n-1}\bigl( G(P)\bigr)=0$ for all projective objects~$P$ of~$\Pro\cC$.

Suppose furthermore that~$X$ is an object of~$\cC$; then the morphism $H^{n-1}\bigl( G(P)\bigr) \to H^{n-1}\bigl( G(X) \bigr)$ vanishes (because the natural morphism $H^{n-1}\bigl( F(X) \bigr)\to H^{n-1}\bigl( G(X) \bigr)$ is an isomorphism).
If now $Q\to X$ is another (not necessarily epi) morphism with~$Q$ a projective object of~$\Pro\cC$, then by lifting to a morphism~$Q\to P$, we see that the  morphism  $H^{n-1}\bigl( G(Q) \bigr)\to H^{n-1}\bigl( G(X) \bigr)$ also vanishes.

In particular, if we write~$P=\quoteslim{i} P_i$ for objects $P_i$ of~$\cC$
(which is possible by Lemma~\ref{lem:quasi-projectives-in-Pro-cat}), 
and recall that the truncation functors on $\Pro D^b(\cC')$ are compatible with cofiltered limits, we have \[H^{n-1}\bigl( G(P)\bigr)\iso \quoteslim{i} H^{n-1}\bigl(G(P_i)\bigr). \] Since (as shown in the previous paragraph, taking~$Q$ to be~$P$ and~$X$ to be~$P_i$) the morphisms $H^{n-1}\bigl( G(P)\bigr)\to H^{n-1}\bigl( G(P_i)\bigr) $ are all zero, we deduce that
$H^{n-1}\bigl( G(P)\bigr)=0$, as required.
\end{proof}

\begin{lem}
  \label{lem:Pro commutative diagram}
  Let $f : \cC \to \Pro(\cC')$ be right exact, and let $F: D^-(\Pro \cC) \to \Pro D^b(\cC')$ be the right $t$-exact functor
  corresponding to $\Pro(f)$ via Proposition~{\em \ref{prop:left derived functors}}~\eqref{item:80}.
  Then we can extend the commutative diagram~\eqref{eqn:Pro comparison}
to a commutative diagram
\begin{equation}
\label{eqn:Pro diagram}
\begin{tikzcd}[column sep=large,row sep=large]
  & D^-(\mathcal{C}) \ar[d,"j"]  \ar[dr,"F\circ j"]
\\
  D^b(\mathcal{C}) \ar[r,"i"] \ar[dr,hookrightarrow] \ar[ur]
  & D^-(\mathrm{Pro}\,\mathcal{C}) \ar[d,"p"] \ar[r,"F"]
  & \mathrm{Pro}\,D^b(\mathcal{C}')
\\
  & \mathrm{Pro}\,D^b(\mathcal{C}) \ar[ur,"\Pro(F\circ i)"']
\end{tikzcd}
\end{equation}%
\end{lem}
\begin{proof}
We only have to prove that the bottom right triangle commutes.
Since $p$ restricts to the identity on~$\Pro \cC$,
it follows from Lemma~\ref{lem:truncated cofiltered limit-preserving}  that
$$\bigl(\Pro(F\circ i) \circ p \bigr)_{| \Pro \cC} 
= \Pro(F\circ i)_{|\Pro \cC} = \Pro\bigl( (F\circ i)_{|\cC})
= F_{| \Pro \cC}.$$
Proposition~\ref{prop:left derived functors}
now shows that 
$\Pro(F\circ i) \circ p$ and $F$ are naturally isomorphic, as required.
\end{proof}

\begin{thm}
\label{thm:pro left derived functors}
If ~$\cC$ and~ $\cC'$ are two small abelian categories, then the above constructions induce equivalences between the following categories.
\begin{enumerate}
  \item  The ordinary category of right exact functors
    $f:\cC \to \Pro \cC'$.
    \item  The ordinary category of right exact cofiltered limit-preserving
functors $\Pro(f):\Pro \cC \to \Pro \cC'$.
\item\label{item:81} The $\infty$-category of right $t$-exact functors
  $F:D^-(\Pro \cC)\to \Pro D^b(\cC')$ which carry projective objects of~$\Pro \cC$
  into $\Pro \cC'$ and for which $\tau^{\geq 0} F_{| \Pro \cC}$ is cofiltered
  limit-preserving.
\item\label{item:82} The $\infty$-category of right $t$-exact
limit-preserving functors $F':\Pro D^b(\cC)\to \Pro D^b(\cC')$ which carry
projective objects of~$\Pro\cC$ into $\Pro \cC'$.
\end{enumerate}
Under these equivalences, we have $F' = \Pro(F \circ i)$, and $F = F' \circ p$.
\end{thm}

\begin{proof}
The equivalences between the first three categories are immediate from Proposition~\ref{prop:left derived functors}~\eqref{item:80}, so it remains to see that the functors $F \mapsto \Pro(F\circ i)$ and $ F' \mapsto F'\circ p$ are inverse equivalences of categories between ~\eqref{item:81} and~\eqref{item:82}.
This follows from %
Lemma~\ref{lem:Pro commutative diagram},
which shows that the composites in
each direction of the purported equivalences are naturally equivalent to the
identity functors. Indeed, to show that $F' \mapsto F'\circ p \mapsto \Pro
(F'\circ p \circ i)$
is isomorphic to the identity, we need to check that
$F' = \Pro(F'\circ p \circ i) = \Pro F'_{| D^b(\cC)}$.  This is precisely the
statement that $F'$ is cofiltered 
limit-preserving (or equivalently, limit-preserving, since~$F'$ is exact). %
On the other hand, the statement that $F \mapsto \Pro(F \circ i) \mapsto \Pro(F \circ i) \circ p$ is isomorphic to the identity 
follows from the commutativity of the bottom right triangle in~\eqref{eqn:Pro diagram}.\qedhere

\end{proof}

\begin{rem}
\label{rem:KS}
The functor $p:D^-(\Pro \cC) \to \Pro D^b(\cC)$ (or, more precisely, its restriction
to $D^b(\cC)$) constructed above is a Pro-analogue, in the $\infty$-categorical
context, of the morphism
$J: D^b(\Ind \cC) \to \Ind D^b(\cC)$ constructed in~\cite[Thm.~15.4.3]{MR2182076}:
see Section~\ref{subsubsec:Pro co-Yoneda} for an elaboration of this point.
Our Theorem~\ref{thm:pro left derived functors}
is then a Pro- (and $\infty$-categorical) analogue of~\cite[Prop.~15.4.7]{MR2182076}.
\end{rem}

\subsubsection{A compatibility}

Suppose now that~$\cC$  itself has enough projectives, so that we can apply Theorem~\ref{thm:left derived functors} to construct derived functors whose source is~$D^-(\cC)$.
We will show that this is compatible with the equivalences of Theorem~\ref{thm:pro left derived functors}.
Note firstly that by Corollary~\ref{cor:left derived functors exact case}, the $t$-exact functor 
\[
j:D^-(\cC)\to D^-(\Pro\cC)
\]
appearing in~\eqref{eqn:Pro comparison}, 
which is the functor given by the functoriality of the formation of the derived category~$D^{-}$, 
can also be thought of as being the derived functor arising from the exact functor $\cC\to \Pro\cC$. %

As in Lemma~\ref{lem:Pro commutative diagram}, we suppose that $f:\cC\to\Pro(\cC')$ is a right exact functor, and
we write $F:D^-(\Pro \cC) \to \Pro D^b(\cC')$
for the right $t$-exact functor corresponding to $\Pro(f)$ via Proposition~\ref{prop:left derived functors}~\eqref{item:80}.
We also write $F_{\cC}:D^-(\cC) \to \Pro D^b(\cC')$, $F^{-}_{\cC}:D^-(\cC) \to D^-(\Pro \cC')$ for the right $t$-exact functors corresponding to $f$ via Theorem~\ref{thm:left derived functors}.

We can then consider the following variant on
 the commutative diagram~\eqref{eqn:Pro diagram}:
\begin{equation}
\label{eqn:Pro diagram enough projectives}
\begin{tikzcd}[column sep=large,row sep=large]
  & D^-(\mathcal{C}) \ar[d,"j"]  \ar[dr, "F_{\cC}"] \ar[r,"F^{-}_{\cC}"] &  D^-(\Pro\mathcal{C'})\ar[d,"p'"]
\\
  D^b(\mathcal{C}) \ar[r,"i"] \ar[dr,"b"' hookrightarrow] \ar[ur,"a"]
  & D^-(\mathrm{Pro}\,\mathcal{C}) \ar[d,"p"] \ar[r,"F"]
  & \mathrm{Pro}\,D^b(\mathcal{C}')
\\
  & \mathrm{Pro}\,D^b(\mathcal{C}) \ar[ur,"\Pro(F\circ i)"']
\end{tikzcd}
\end{equation}

\begin{lem}
  \label{lem:Pro commutative diagram enough projective}Suppose as above that~ $\cC$ has enough projectives, and that~$f:\cC\to\Pro(\cC')$ is right exact.
Then the diagram ~{\em \eqref{eqn:Pro diagram enough projectives}} commutes; in particular, $\Pro(F\circ i)=\Pro(F_{\cC}\circ a)$. %
\end{lem}
\begin{proof}Note firstly that~$F_{\cC}$ is naturally isomorphic to $p'\circ F^{-}_{\cC}$, by Theorem~\ref{thm:left derived functors}
  and the $t$-exactness of~$p'$.  
  We next show that $\Pro(F\circ i)=\Pro(F_{\cC}\circ a)$.
  Since $F_{\cC}\circ a$ is right $t$-exact, its Pro-extension
is a right $t$-exact functor 
\[
\Pro(F_{\cC}\circ a): \Pro D^b(\cC) \to \Pro D^b(\cC')
\]
satisfying
$$(\tau^{\geq 0} \Pro(F_{\cC}\circ a))_{| \Pro \cC} 
\iso 
\Pro (\tau^{\geq 0} F_{\cC}\circ a)_{| \Pro \cC}
\iso \Pro f.$$
Thus, to show that $\Pro(F_{\cC}\circ a)$ coincides with $\Pro(F\circ i)$, 
we have (according to Theorem~\ref{thm:pro left derived functors})
to show that $\Pro(F_{\cC}\circ a)$ takes projective objects of $\Pro \cC$ to
the heart $\Pro \cC'$ of $\Pro D^b(\cC')$.

By Lemma~\ref{lem:quasi-projectives-in-Pro-cat}, any projective object in~ $\Pro\cC$ is of the form $\quoteslim{i} X_i$, with the $X_i$ being projective objects of~$\cC$. 
Thus, for such an object, 
$$\Pro(F_{\cC}\circ a) (\quoteslim{i} X_i) \coloneqq  \quoteslim{i} F_{\cC}(X_i)$$
indeed lies in~$\Pro \cC'$, since each $F_{\cC}(X_i)$ lies in~$\Pro \cC'$ (as each
$X_i$ is projective in~$\cC$, and $F_{\cC}$ is associated to~$f$ via the construction
of Theorem~\ref{thm:left derived functors}).
This concludes the proof that $\Pro(F_\cC \circ a) = \Pro(F \circ i)$.

At this point, bearing in mind Lemma~\ref{lem:Pro commutative diagram}, in order to conclude the proof that~\eqref{eqn:Pro diagram enough projectives} commutes, 
there remains to check that in fact $F_{\cC}=F\circ j$, which we will do by applying Theorem~\ref{thm:left derived functors}.
Note first that the $t$-exactness of~$j$ implies that
\[
\tau^{\geq 0}(F\circ j)|_{\cC} = \tau^{\geq 0}F|_{\Pro \cC} \circ j|_{\cC} = f;
\]
so by Theorem~\ref{thm:left derived functors}, it suffices to show that $F\circ j$ takes projective objects of~$\cC$ to $\Pro \cC'$.
To this end, note that
\[ (F\circ j)\circ a=   F\circ i=\Pro(F\circ i)\circ b=\Pro(F_{\cC}\circ a)\circ b=F_{\cC}\circ a.\]
Since the restriction of~$a$ to~$\cC$ is the identity functor, this implies that the restriction of $F\circ j$ to~$\cC$ agrees 
with the restriction of~$F_{\cC}$ to~$\cC$, and since $F_{\cC}$ takes projective objects of~$\cC$ to $\Pro \cC'$ by construction, we are done.
\end{proof}

\subsubsection{Deriving adjoint functors}

Let $g: \cC' \to \cC$ be an exact functor between small abelian categories.
We continue to write~$g$ for its $\Pro$-extension $g:\Pro(\cC')\to\Pro(\cC)$, which by Lemma~\ref{lem:left-adjoint-to-Pro-functor-abelian} has a left adjoint $f:\Pro \cC \to \Pro \cC'$,
    which is right exact and cofiltered limit-preserving.

\begin{prop}%
  \label{prop:pro-adjoints}As above, let $g:\cC'\to\cC$ be an exact functor between small
  abelian categories, let $g:\Pro\cC'\to\Pro\cC$ be its $\Pro$-extension, and
  let $f:\Pro\cC\to\Pro\cC'$ be the left adjoint of~$g$. Let
  $F:\Pro D^b(\cC) \to \Pro D^b(\cC')$ and $G: \Pro D^b(\cC') \to \Pro D^b(\cC)$
  be the limit-preserving functors determined by~$f,g$ respectively by
  Theorem~{\em \ref{thm:pro left derived functors}} and  Corollary~{\em \ref{cor:pro t-exact derived functors}}, which are respectively right
  $t$-exact and $t$-exact.

  Then: 
  \begin{enumerate}
  \item\label{item:14} $F$ is left adjoint to~$G$. 
\item\label{item:15}   If furthermore~$f$ is exact, then~$F$ is
  $t$-exact, and is the Pro-extension of the $t$-exact functor $D^b(\cC) \to \Pro D^b(\cC')$
  determined by~$f$. %
  \item\label{item:16} Continuing to assume that~$f$ is exact, the
    $\Ind$-extensions of~$F$ and~$G$ are adjoint $t$-exact functors. 
  \end{enumerate}

\end{prop}
\begin{proof}
The composite $G\circ F: \Pro D^b(\cC) \to \Pro D^b(\cC)$ is again right $t$-exact
and limit-preserving, and takes 
projective objects of~$\Pro \cC$ into $ \Pro \cC$.
Theorem~\ref{thm:pro left derived functors} (in particular,
the equivalence of $\infty$-categories that it provides)
thus implies that the unit
of adjunction $\id_{\Pro \cC} \to g\circ f$ induces a morphism
$\id_{\Pro D^b(\cC)}\to G\circ F.$ 
For any objects $X$ of $\Pro D^b(\cC)$ and $Y$ of $\Pro D^b(\cC')$,
this morphism induces functorial maps
\begin{multline*}
\Maps_{\Pro D^b(\cC')}\bigl( F(X), Y\bigr)
\\
\xrightarrow{\text{ apply } G} 
\Maps_{\Pro D^b(\cC)}\bigl( GF(X), G(Y)\bigr)
\\
\longrightarrow \Maps_{\Pro D^b(\cC)}\bigl(X,G(Y)\bigr)
.
\end{multline*}
We claim that this composite is an isomorphism, thus showing that $F$ and $G$
are indeed adjoints.

Writing each of $X$ and $Y$ as the limit of objects in $D^b(\cC)$ and $D^b(\cC')$
respectively, and using the fact that $F$ and $G$ are limit-preserving,
we reduce to the case when $X$ and $Y$ are in fact objects of $D^b(\cC)$
and~$D^b(\cC')$. An obvious d\'evissage using the fibre
sequences~\eqref{eqn:fibre-sequence-truncations} for~$X$ and~$Y$ (and the
exactness of~$F$ and~$G$) reduces us the 
case that~$X$ is an object of~$\cC$ and~$Y$ is an object of~$\cC'$. %

We are thus reduced to a consideration of the morphism
$$
\Maps_{\Pro D^b(\cC')}\bigl( F(X), Y\bigr) 
\to
\Maps_{\Pro D^b(\cC)}\bigl( X, G(Y)\bigr) 
$$
with $X$ and $Y$ objects of~$\cC$ and $\cC'$ respectively.
The adjunction property of~$\tau^{\geq 0}$ lets us rewrite this as
$$
\Maps_{\Pro D^b(\cC')}\bigl( \tau^{\geq 0}F(X), Y\bigr) 
\to
\Maps_{\Pro D^b(\cC)}\bigl( X, G(Y)\bigr), 
$$
and then (using the $t$-exactness of~$G$) as
$$
\Maps_{\Pro \cC'}\bigl( f(X), Y\bigr) 
\to
\Maps_{\Pro \cC}\bigl( X, g(Y)\bigr). 
$$
This is precisely the map induced by the unit of adjunction for $f$ and~$g$,
and so is indeed an isomorphism, as required. This concludes the proof of the first point.
The second point then follows from
Corollary~\ref{cor:pro t-exact derived functors}, and the third point is clear,
since $\Ind$-extension preserves $t$-exactness and adjoints.
\end{proof}

\subsection{Pro-categories and the co-Yoneda embedding}\label{subsubsec:Pro co-Yoneda}%
Let $\cC$ be a small abelian category.
We may restrict the $t$-exact functor $p$ of~\eqref{eqn:Pro comparison} to a
$t$-exact functor
\begin{equation}
\label{eqn:bounded pro functor}
p: D^b(\Pro \cC) \to \Pro D^b(\cC).
\end{equation}
If $\cB$ is an exact full abelian subcategory of $\Pro \cC$, containing the image of the canonical embedding $\cC \to \Pro \cC$,
then the composite of~\eqref{eqn:bounded pro functor} with
the canonical $t$-exact functor $D^b(\cB) \to D^b(\Pro \cC)$ induces
a $t$-exact functor
\begin{equation}\label{eqn:defn-of-q}
  q: D^b(\cB) \to \Pro D^b(\cC).
\end{equation}
The functor~\eqref{eqn:bounded pro functor} need not be fully faithful in general,
by e.g.\ \cite[Rem.\ 15.4.5, Ex.\ 15.2]{MR2182076},
and so the functor $q$ certainly need not be fully faithful in general.
Nevertheless, it is sometimes the case that~$q$ is fully faithful: for example, by Lemma~\ref{lem:Pro commutative diagram}, 
this is true when $\cB = \cC$. 
Our goal in this section is to give a useful 
sufficient condition for full faithfulness of~$q$
(see Lemma~\ref{lem:full faithfulness criterion}). %

Before doing so,
it will be useful to give another description of~\eqref{eqn:bounded pro functor}.
In order to do this, we %
first explain the ``dual'' set-up, where we consider Ind-categories
rather than Pro-categories.  
This will let us speak of presheaves, Yoneda embeddings,
and so on, rather than their ``co'' variants. 
 Passing to opposite categories
at the end of the discussion will then give the description we need. 

If~$\cC$ is a small abelian category,
the functor $p : D^-(\Pro \cC^\op) \to \Pro D^b(\cC^\op)$
induces, by passage to opposite categories, 
a $t$-exact functor
$\check p: D^+(\Ind \cC) \to \Ind D^b(\cC)$, 
which we may restrict to a functor   
$$\check p: D^b(\Ind \cC) \to \Ind D^b(\cC).$$
By Corollary~\ref{cor:t-exact derived functors}, $\check p$ is the unique (up to equivalence) $t$-exact functor
$D^b(\Ind \cC) \to \Ind D^b(\cC)$ inducing the identity on hearts.
Using the Yoneda embedding, we can also give the following alternative description of~$\check p$:
the canonical functor $D^b(\cC) \to D^b(\Ind \cC)$
gives rise to an exact 
functor
\begin{equation}
\label{eqn:Yoneda} 
D^b(\Ind \cC) \to \Fun^{\mathrm{ex}}\bigl(D^b(\cC)^{\op},\Sp\bigr),
\quad X \mapsto \RHom_{D^b(\Ind \cC)}(\text{--}, X).
\end{equation}
If we regard $\Ind D^b(\cC)$ as a stable sub-$\infty$-category
of~$\Fun^{\mathrm{ex}}\bigl(D^b(\cC)^{\op},\Sp\bigr)$
via the discussion of Remark~\ref{rem:stable indization}, %
then we see that $\check p$ and~\eqref{eqn:Yoneda}
are exact functors  
\[
D^b(\Ind \cC) \to \Fun^{\mathrm{ex}}\bigl(D^b(\cC)^{\op},\Sp\bigr),
\]
both of which induce the canonical equivalence
$\Ind \cC \iso \bigl(\Ind D^b(\cC)\bigr)^{\heartsuit}$.
By Theorem~\ref{thm:exact derived functors},
they must coincide.

Concretely, this means that if $X$ is an object of $D^b(\Ind \cC)$ and $Y$ is an
object of~$D^b(\cC)$,
then
\begin{equation}\label{eqn:Yoneda equivalence}
\RHom_{D^b(\Ind \cC)}(Y,X) \iso \RHom_{\Ind D^b(\cC)}\bigl(Y, \check p(X) \bigr).
\end{equation}
Returning now to the ``un-dualized'' setting of $\Pro$-categories,
we find that the restriction~\eqref{eqn:bounded pro functor} of $p$
admits a co-Yoneda description. 
Concretely, if $X$ %
is an object of $D^b(\Pro \cC)$
and $Y$ is an object of~$D^b(\cC)$,
then the canonical morphism
\begin{multline}
\label{eqn:co-Yoneda iso}
\RHom_{D^b(\Pro \cC)}(X,Y)
=\RHom_{D^-(\Pro \cC)}(X,Y)
\\
\to \RHom_{\Pro D^b(\cC)}\bigl (p(X),Y\bigr) 
\end{multline}
induced by~$p$ 
is an isomorphism.
(Here we have used the commutativity of~{\em \eqref{eqn:Pro comparison}}
to identify the restriction of $p$ to $D^b(\cC)$ with the canonical embedding of
$D^b(\cC)$ into $\Pro D^b(\cC)$, and so to identify~$Y$ with~$p(Y)$.)

We then have the following lemma, which specializes, in the case when $\cB$ equals $\Pro \cC$ and~$p = q$, 
to the fact that~\eqref{eqn:co-Yoneda iso} is an isomorphism.

\begin{lem}
\label{lem:pro co-yoneda}
Let~$\cC$ be a small abelian category, and let $\cC \subseteq \cB \subseteq \Pro \cC$ be an exact full abelian subcategory of~$\Pro \cC$.
If $X$ and $Y$ are objects of $D^b(\cB)$ and~$D^b(\cC)$ respectively,
then $q$ induces a natural isomorphism
$$\RHom_{D^b(\cB)}(X,Y) \iso \RHom_{\Pro D^b(\cC)} \bigl( q(X), Y),$$
where we used the commutativity of~{\em \eqref{eqn:Pro comparison}}
to identify the restriction of $p$ to $D^b(\cC)$ with the canonical embedding of
$D^b(\cC)$ into $\Pro D^b(\cC)$, and so to identify~$Y$ with~$q(Y)$. 
\end{lem}
\begin{proof}
A standard truncation argument 
reduces to the case when $X$ is an object of $\cB$ and $Y$ is an object
of~$\cC$.
We may and do write $X = \quoteslim{i} X_i$ for some objects $X_i$ of~$\cC$.
We then consider the sequence of morphisms
\begin{multline}\label{eqn:composition for Pro co-Yoneda}
\colim_i \RHom_{D^b(\cC)}(X_i,Y) \to 
\colim_i \RHom_{D^b(\cB)}(X_i,Y) \to 
\RHom_{D^b(\cB)}(X,Y) 
\\
\to
\RHom_{D^b(\Pro \cC)}(X,Y) \iso
\RHom_{\Pro D^b(\cC)}\bigl (q(X),Y\bigr) 
\iso \colim_i \RHom_{D^b(\cC)}(X_i,Y)
\end{multline}
(with the first three morphisms being the evident ones,
the first isomorphism being
that of~\eqref{eqn:co-Yoneda iso},
and the second
isomorphism arising from the definition
of morphisms in the Pro category). 
Evidently the composite of this sequence of morphisms is the identity. %

Passing to cohomology, we see that the natural morphism
$$\Ext^i_{D^b(\cB)}(X,Y) \to \Ext^i_{D^b(\Pro \cC)}(X,Y)$$
is surjective for all $i\geq 0$, for any objects $X$ and $Y$ as above. 
Since this morphism is
an isomorphism when $i = 0$, a standard argument 
then shows that it is in fact an isomorphism for all~$i \geq 0.$
(For the sake of completeness,
we have recalled this argument in Lemma~\ref{lem:standard homological algebra} below; note
that Yoneda $\Ext$ functors computed in any abelian category are element-wise effaceable
in the sense of that lemma.)
Thus in fact~$\RHom_{D^b(\cB)}(X,Y) \iso \RHom_{D^b(\Pro \cC)}(X,Y),$
and the lemma follows by another consideration of~\eqref{eqn:composition for Pro co-Yoneda}.
\end{proof}

\begin{lem}
\label{lem:standard homological algebra} 
Let $\eta^{\bullet}: F^{\bullet} \to G^{\bullet}$ 
be a natural transformation of abelian group-valued $\delta$-functors on an abelian category~$\cA$.
Suppose that $\eta^0$ is an isomorphism, and that
$\eta^n$ is an epimorphism for each~$n$. 
Suppose also that $F^n$ is {\em element-wise effaceable}
for each $n > 0$ {\em (}in the sense that if $\alpha \in F^n(X)$ for
some object $X$ of~$\cA$, then there exists a monomorphism
$X \hookrightarrow Y$ in~$\cA$ with respect to which the image of~$\alpha$ 
in $F^n(Y)$ vanishes{\em )}.
Then $\eta^{\bullet}$
is in fact a natural isomorphism.
\end{lem}
\begin{proof}
We prove that $\eta^n$ is an isomorphism for each~$n \geq 0$
via induction on~$n$, the case $n = 0$ holding by hypothesis.
Thus, we assume that $\eta^i$ is an isomorphism for all $0 \leq i < n$,
and prove that $\eta^n$ is an isomorphism.  In fact, we need only
show that $\eta^n$ is a monomorphism, since it is an epimorphism by assumption.
Choose an object $X$ in $\cA$, and an element $\alpha \in \ker\bigl(\eta^n(X):
F^n(X) \to G^n(X)\bigr).$

By hypothesis we may find a monomorphism $X \hookrightarrow Y$ such 
that the image of $\alpha$ in $F^n(Y)$ vanishes. Let $Z$ denote the cokernel of 
this embedding; we then have a commutative diagram with exact rows
$$\xymatrix{
F^{n-1}(Y) \ar[r]\ar[d] & F^{n-1}(Z) \ar[r]\ar[d]
&  F^n(X) \ar[r]\ar[d] &  F^n(Y) \ar[d] \\
G^{n-1}(Y) \ar[r] & G^{n-1}(Z) \ar[r]
&  G^n(X) \ar[r] &  G^n(Y)  }
$$
A straightforward diagram chase (using that the left two vertical arrows
are isomorphisms) shows that $\alpha$ itself vanishes.  Thus we see that
$\eta^n$ is indeed a monomorphism, as required.
\end{proof}

The next lemma provides the promised criterion for the functor $q$ to be fully faithful.

\begin{lem}
\label{lem:full faithfulness criterion}Suppose that $\cC \subseteq \cB \subseteq \Pro \cC$ is an exact full abelian subcategory of $\Pro \cC$. %
Suppose further that for any objects $X, Y$ of~$\cB$, 
writing $Y = \quoteslim{i}Y_i$ as a cofiltered limit of objects of~$\cC$,
the canonical morphism %
\begin{equation}\label{eqn:assumption in full faithfulness criterion}
\RHom_{D^b(\cB)}(X,Y) \to \lim_i \RHom_{D^b(\cB)}(X,Y_i)
\end{equation}
is an isomorphism. 
Then $q: D^b(\cB) \to \Pro D^b(\cC)$ is fully faithful.
\end{lem}
\begin{proof}
An easy argument with truncations (taking into account the
exactness of~$q$) %
reduces the proposition to verifying that the natural morphism
\begin{equation}\label{eqn:to prove isomorphism for full faithfulness criterion}
\RHom_{D^b(\cB)}(X,Y)
\to
\RHom_{\Pro D^b(\cC)}\bigl(q(X),q(Y)\bigr)
\end{equation}
is an isomorphism
whenever $X$ is an object of $D^b(\cB)$ and $Y$ is an object of~$\cB$.
Similarly, our assumption that~\eqref{eqn:assumption in full faithfulness criterion} is an isomorphism when~$X \in \cB$ implies that~\eqref{eqn:assumption in full faithfulness criterion}
is also an isomorphism when~$X \in D^b(\cB)$.

Write $Y =\quoteslim{i} Y_i$. 
Since~\eqref{eqn:assumption in full faithfulness criterion} is an isomorphism, and $q(Y)$ coincides with~$Y$ (viewed as an object of~$\Pro \cC$), 
we see that it furthermore suffices to prove that~\eqref{eqn:to prove isomorphism for full faithfulness criterion} 
is an isomorphism when~$Y$ is replaced by each of the~$Y_i$;
that is, we reduce to the case when $Y$ is an object of~$\cC$.
The result in this case follows from Lemma~\ref{lem:pro co-yoneda}.
\end{proof}

\subsection{Derived tensor products}%
\label{subsec:derived tensor products in abelian categories}
We develop derived analogues of the various constructions considered in 
Section~\ref{subsubsec:tensor products in abelian categories}.

\subsubsection{The derived category of finite length $R$-modules}
Let~$R$ be a Noetherian profinite $\cO$-algebra.
By ~\eqref{improved compact modules and Pro}, we have
   $\Mod_c(R) \isoto \Pro\Mod^{\fl}(R)$.
\begin{df}
\label{def:D^b_fl(R)}
Write
$D^b_{\fl}(R)$ to denote the full subcategory of $D^b\bigl(\Mod_c(R)\bigr)$ consisting of objects whose cohomologies are of finite length, i.e.\ lie in~$\Mod^{\fl}(R)$.
\end{df}
\begin{remark}
Since $\Mod^{\fl}(R)$ is a Serre subcategory of~$\Mod_c(R)$, Definition \ref{def:D^b_fl(R)}
makes sense, as a particular instance of Definition~\ref{defn:D-sub-A-B}.
Despite what the notation may suggest, $D^b_{\fl}(R)$ {\em a priori} depends on
the structure of $R$ as a topological ring, since it is defined as a full subcategory
of $D^b\bigl(\Mod_c(R)\bigr)$.
Part~(1) of Lemma~\ref{lem:finite-length-derived-cat} shows in particular that it admits another description
depending only on the structure of $R$ as an abstract ring.
\end{remark}

Recall that by Lemma~\ref{properties of compact modules} %
we have natural fully faithful exact embeddings
$\Mod^{\fl}(R) \hookrightarrow \Mod_c(R)$, 
$\Mod^{\fp}(R) \hookrightarrow \Mod_c(R)$, and $\Mod^{\fl}(R)\hookrightarrow\Mod^{\fp}(R)$.
The following lemma shows that the derived extensions of these functors remain fully faithful, and, 
in the case of the first embedding, identifies the essential image with the
stable $\infty$-category $D^b_{\fl}(R)$ introduced in Definition~\ref{def:D^b_fl(R)}.

   \begin{lem}
     \label{lem:finite-length-derived-cat}
\leavevmode
\begin{enumerate}
\item
The natural $t$-exact functor $D^b\bigl(\Mod^{\fl}(R)\bigr)\to D^b\bigl(\Mod_c(R)\bigr)$ is
fully faithful,
and induces an equivalence
$D^b(\Mod^{\fl}(R))\iso D^b_{\fl}(R).$
\item
The natural $t$-exact functor $D^b\bigl(\Mod^{\fp}(R)\bigr)\to D^b\bigl(\Mod_c(R)\bigr)$ is
fully faithful.
\item The natural $t$-exact functor 
$D^b\bigl(\Mod^{\fl}(R)\bigr)\to D^b\bigl(\Mod^{\fp}(R)\bigr)$ is
fully faithful. %
\end{enumerate}
   \end{lem}
   \begin{proof}By Lemma~\ref{properties of compact modules}~\eqref{item: compact 2},
$\Mod_c(R)^\op$ is a locally finite category, in the sense recalled in~\ref{subsec:locally finite categories},
i.e.\ it is a Grothendieck category with compact objects ~$\Mod^{\fl}(R)^{\op}$
(see~Lemma~\ref{locally finite categories}~(1) for this assertion regarding compact objects).
Since $\Mod^{\fl}(R)^{\op}$ is abelian, $\Mod_c(R)^\op$ is a locally coherent abelian category,
and part~(1) of the lemma is therefore immediate from Lemma~\ref{lem:compact-objects-bounded-Serre-category} by passing to opposite categories.

To prove~(2), it follows from Proposition~\ref{prop:check-full-faithful-on-compact-generators}~(2)
that it suffices to show that the embedding
$\Mod^{\fp}(R) \to \Mod_c(R)$ induces an isomorphism on all~$\Ext^i$.  
This follows from~\cite[Prop.~3.16]{MR1474172}. (Note that Corollary 3.12 of the
same reference shows that the objects of~$\Mod^{\fp}(R)$ are Noetherian
when regarded as objects of~$\Mod_c(R)$, so that the cited Proposition does
indeed apply.)

Part~(3) follows immediately from parts~(1) and~(2).
\end{proof}

\begin{remark}
\label{rem:finite-length-derived-cat}
Lemma~\ref{lem:finite-length-derived-cat}
gives rise to the following convenient $t$-exact fully faithful embedding:
\begin{equation}
\label{eqn:fl to fp}
D^b_{\fl}(R) \iso D^b\bigl(\Mod^{\fl}(R)\bigr) \hookrightarrow D^b\bigl(\Mod^{\fp}(R)\bigr)
\iso D^b_{\fp}(R),
\end{equation}
with the first equivalence being the inverse of the equivalence
of part~(1) of the lemma, the fully faithful embedding being that
of part~(3) of the lemma, and the final equivalence being that
of Corollary~\ref{finitely presented derived category for coherent ring}. 
The utility of~\eqref{eqn:fl to fp} is that while the source is defined with reference to 
the topology of~$R$, the target depends only on $R$ as an abstract ring.
\end{remark}

\begin{remark}\label{rem:derived fp to pro fl}
By the discussion in Section~\ref{app-subsub-adjoint-functors}, %
the Pro-extension of~\eqref{eqn:fl to fp} has a left adjoint
\begin{equation}\label{eqn:derived fp to pro fl}
\Pro D^b_{\fp}(R) \to \Pro D^b_{\fl}(R).
\end{equation}
This left adjoint is $t$-exact, and maps the constant object $R \in \Pro \Mod^{\fp}(R)$ to $\quoteslim{n}R/\rad(R)^n \in \Pro \Mod^{\fl}(R)$. 
To see this, we apply Proposition~\ref{prop:pro-adjoints} with~$g : \Mod^{\fl}(R) \to \Mod^{\fp}(R)$ taken to be the inclusion, noting that the left adjoint 
$\Pro \Mod^{\fp}(R) \to \Pro \Mod^{\fl}(R)$ to~$\Pro(g)$ is exact and sends~$R$ to $\quoteslim{n}R/\rad(R)^n$, by Remark~\ref{rem:completed-tensor-product}.

Note furthermore that the composite
\[
D^b_{\fl}(R) \xrightarrow{\eqref{eqn:fl to fp}} D^b_{\fp}(R) \subset \Pro D^b_{\fp}(R) \xrightarrow{\eqref{eqn:derived fp to pro fl}} \Pro D^b_{\fl}(R)
\]
is the natural embedding of $D^b_{\fl}(R)$ in~$\Pro D^b_{\fl}(R)$. 
(These functors being $t$-exact, this statement can be checked on restriction to the hearts.)
\end{remark}

\subsubsection{Derived tensor products}
We now return to the setting of Definition~\ref{defn:R-module-objects}.
We assume in addition that the ring~$R$ is coherent, so that~$\Mod(R)$ is locally coherent, 
with compact objects given by the abelian category~$\Mod^{\fp}(R)$.
By Corollary~\ref{finitely presented derived category for coherent ring}, 
the natural functor 
\[
D^b(\Mod^{\fp}(R))\to D^b_{\fp}(R)\coloneqq D^b_{\Mod^{\fp}(R)}(\Mod(R)) 
\]
is an equivalence.
Similarly, bearing in mind the fact that $\Mod^{\fp}(R) \to \Mod(R)$ preserves projectives, 
it follows from~\cite[Prop.\ 1.3.3.7]{LurieHA} (see also~\cite[Prop.\ A.5.7]{emerton2023introduction}) that 
the natural functor 
\[
D^-(\Mod^{\fp}(R))\to D^-_{\fp}(R)\coloneqq D^-_{\Mod^{\fp}(R)}(\Mod(R)) 
\]
is an equivalence.

We now apply the machinery of Appendix~\ref{subsec:left derived} to construct various derived versions of the  tensor products provided by Proposition~\ref{prop:Eilenberg-Watts} and Lemma~\ref{lem:EW-gives-completed-tensor-product}.
\begin{lem}%
  \label{lem:derived-tensor-product-EW}Let~$R$ be a coherent ring, %
  and let~$\cA$ be an abelian category.
  \begin{enumerate}
  \item\label{item:83}  The functor $F\mapsto F(R)$ gives equivalences of categories between the category of right $R$-modules~$M$ in~$\cA$ and the following categories of functors.
    \begin{enumerate}[series=tensorfunctors]
    \item\label{item:95} The category of right exact functors \[M\otimes_R\text{--}:\Mod^{\fp}(R)\to\cA.\]
    \item\label{item:96} The category of right exact cofiltered limit-preserving functors \[M\cotimes_R\text{--}:\Pro\Mod^{\fp}(R)\to\Pro(\cA)\]taking~$R$ to an object of~$\cA$.
    \item\label{item:89}  The category of right $t$-exact functors \[M\otimes_R^{L}\text{--}:D^b_{\fp}(R)\to D^-(\cA)\]taking~$R$ to an object of~$\cA$.  %
    \item\label{item:90}  The category of right $t$-exact functors \[M\otimes_R^{L}\text{--}:D^-_{\fp}(R)\to D^-(\cA)\] taking~$R$ to an object of~$\cA$.
    \item\label{item:91}  The category of right $t$-exact limit-preserving functors \[M\cotimes_R^{L}\text{--}:\Pro D^b_{\fp}(R)\to\Pro D^b(\cA)\] taking~$R$ to an object of~$\cA$. %
    \end{enumerate}
    Furthermore, we have a commutative diagram
    \begin{equation}\label{eqn:first-commutative-diagram-derived-tensor}\begin{tikzcd}
	D^b_{\fp}(R) & D^-_{\fp}(R) & D^-(\cA) & \Pro D^b(\cA) \\
	& \Pro D^b_{\fp}(R)
	\arrow[from=1-1, to=1-2]
	\arrow["\eqref{item:89}", curve={height=-25pt}, from=1-1, to=1-3]
	\arrow[from=1-1, to=2-2]
	\arrow["\eqref{item:90}", from=1-2, to=1-3]
	\arrow[from=1-2, to=2-2]
	\arrow[from=1-3, to=1-4]
	\arrow["\eqref{item:91}"', from=2-2, to=1-4]
\end{tikzcd}\end{equation}
    \item\label{item:84} The functor $F\mapsto F(R)$ gives equivalences of categories between the category of right $R$-modules~$M$ in~$\Pro(\cA)$ and the following categories of functors.
    \begin{enumerate}[resume=tensorfunctors]
    \item\label{item:86a} The category of right exact functors \[M\otimes_R\text{--}:\Mod^{\fp}(R)\to\Pro(\cA).\]
    \item\label{item:86} The category of right exact cofiltered limit-preserving functors \[M\ccotimes_R\text{--}:\Pro\Mod^{\fp}(R)\to\Pro(\cA).\]
    \item\label{item:87} The category of right $t$-exact functors \[M\otimes_R^{L}\text{--}:D^b_{\fp}(R)\to \Pro D^b(\cA)\] taking~$R$ to an object of~$\Pro(\cA)$.
    \item\label{item:92}  The category of right $t$-exact functors \[M\otimes_R^{L}\text{--}:D^-_{\fp}(R)\to \Pro D^b(\cA)\] taking~$R$ to an object of~$\Pro(\cA)$.
    \item\label{item:93}  The category of right $t$-exact limit-preserving functors \[M\ccotimes_R^{L}\text{--}:\Pro D^b_{\fp}(R)\to\Pro D^b(\cA)\] taking~$R$ to an object of~$\Pro(\cA)$.
    \end{enumerate}Furthermore, we have a commutative diagram
    \begin{equation}\label{eqn:second-commutative-diagram-derived-tensor}\begin{tikzcd}
	D^b_{\fp}(R) \\
	D^-_{\fp}(R) & \Pro D^b(\cA) \\
	\Pro D^b_{\fp}(R)
	\arrow[from=1-1, to=2-1]
	\arrow["\eqref{item:87}", from=1-1, to=2-2]
	\arrow["\eqref{item:92}", from=2-1, to=2-2]
	\arrow[from=2-1, to=3-1]
	\arrow["\eqref{item:93}"', from=3-1, to=2-2]
\end{tikzcd}\end{equation}
    \item\label{item:85} Suppose furthermore that $R$ is in fact a Noetherian profinite topological $\cO$-algebra.  
    Then the functors $F\mapsto F(R)$ and $F\mapsto \varprojlim_nF(R/\Rad(R)^{n})$ give equivalences of categories between 
    the category of complete right $R$-modules in~$\Pro(\cA)$ {\em (}in the sense of Definition~{\em \ref{defn:complete-object-abelian-category})} and the following categories of functors. %
    \begin{enumerate}[resume=tensorfunctors]
    \item\label{item:95a} The category of right exact functors \[M\otimes_R\text{--}:\Mod^{\fl}(R)\to\Pro(\cA).\]
    \item\label{item:96a} The category of right exact cofiltered limit-preserving functors
      \[M\cotimes_R\text{--}:\Mod_c(R)\to\Pro(\cA).\]
    \item\label{item:88}  The category of right $t$-exact functors \[M\otimes_R^{L}\text{--}:D^b_{\fl}(R)\to \Pro D^b(\cA)\] for which $\limcommand_n \Bigl(M\otimes_R^L \bigl(R/\Rad(R)^n\bigr)\Bigr)$ {\em(}computed in $\Pro D^b(\cA)${\em)} lies in~$\Pro(\cA)$ {\em(}the latter being regarded as the heart of~$\Pro D^b(\cA)${\em )}.
    \item\label{item:94} The category of right $t$-exact limit-preserving functors \[M\cotimes_R^{L}\text{--}:\Pro D^b_{\fl}(R)\to\Pro D^b(\cA)\] taking~$R$ to an object of~$\Pro(\cA)$.
    \end{enumerate}
  \end{enumerate}
\end{lem}
\begin{proof}

We begin with part~\eqref{item:83}.
The equivalences between the category of right $R$-modules in~$\cA$ and the categories of functors in~\eqref{item:95} and ~\eqref{item:96}  are established in  Proposition~\ref{prop:Eilenberg-Watts}, and the existence of the derived functors in \eqref{item:89}--\eqref{item:91} follows from Theorem~\ref{thm:left derived functors} and Proposition~\ref{prop:left derived functors}.
By Theorem~\ref{thm:pro left derived functors} and Lemma~\ref{lem:Pro commutative diagram enough projective}, in order to see that these constructions induce the claimed equivalences of categories (and to see the commutativity of~\eqref{eqn:first-commutative-diagram-derived-tensor}), we must show that in each of \eqref{item:89}--\eqref{item:91}, the assumption that~$R$ is taken to an object of~$\cA$ implies that all projective objects in the heart of the source are taken to objects in the heart
of the target.
In~\eqref{item:89} and~\eqref{item:90} this is clear, since the projective objects of $\Mod^{\fp}(R)$ are direct summands of finite free $R$-modules.
The case of~\eqref{item:91} then follows from Lemma~\ref{lem:quasi-projectives-in-Pro-cat}.

The proof of part ~\eqref{item:84} is almost identical.
Turning to part ~\eqref{item:85}, the equivalence between the category of complete right $R$-modules in~$\Pro(\cA)$ and the category of functors in ~\eqref{item:96a} is established in Lemma~\ref{lem:EW-gives-completed-tensor-product}.
The equivalence between~\eqref{item:95a} and  ~\eqref{item:96a} is immediate from the definition of the $\Pro$-category (recalling also the equivalence $\Mod_c(R) \iso \Pro \Mod^{\fl}(R)$
of~\eqref{improved compact modules and Pro}),
as is the equivalence between \eqref{item:88} and  ~\eqref{item:94}.
By Theorem~\ref{thm:pro left derived functors}, we are reduced to showing that any functor as in~\eqref{item:94} necessarily takes projective objects of $\Mod_c(R)$ to objects of~$\Pro(\cA)$.
This follows from the fact that every projective object of~$\Mod_c(R)$ is a product of finite projective~$R$-modules (see~\cite[\S IV.3, Cor.\ 1]{Gabrielthesis}), and the assumption that the functor in~\eqref{item:94} is limit-preserving.
\end{proof}

We now describe the cohomology of the derived tensor product with a right module in~$\Pro(\cA)$.

\begin{lemma}\label{homology of derived tensor product}
Let~$R$ be a coherent ring, let $\cA$ be an abelian category, and let~$\{M_i: i \in I\}$ be a cofiltered system of right $R$-modules in~$\cA$.
Let $M \coloneqq  \quoteslim{i}M_i$, a right $R$-module in~$\Pro(\cA)$.
Then for all $X \in D^-_{\fp}(R)$ and~$p \in \bZ$, we have
\[
H^p(M \otimes_R^L X) = \quoteslim{i} H^p(M_i \otimes^L_R X)
\] 
as objects of $\Pro D^b(\cA)^\heartsuit = \Pro(\cA)$, where~$M \otimes_R^L -$ and $M_i \otimes_R^L -$ are defined as in
  Lemma~{\em \ref{lem:derived-tensor-product-EW}~\eqref{item:92}}.
\end{lemma}
\begin{proof}
Since truncation functors commute with cofiltered limits in $\Pro D^b(\cA)$, 
the functor 
\[
\quoteslim{i} (M_i \otimes_R^L\text{--}) : D^-_{\fp}(R) \to \Pro D^b(\cA) 
\]
is right $t$-exact. 
It takes~$R$ to $\quoteslim{i} M_i$, computed in $\Pro D^b(\cA)$, which coincides with the limit computed in~$\Pro(\cA)$,
which is~$M$.
(This is because $\Pro(\cA) \to \Pro D^b(\cA)$ is Pro-extended, and so preserves cofiltered limits.)
By the uniqueness assertion in~Lemma~\ref{lem:derived-tensor-product-EW}~\eqref{item:92},
$\quoteslim{i} (M_i \otimes_R^L\text{--})$ is isomorphic to $M \otimes_R^L -$.
Hence $H^p(M \otimes_R^L X) = \quoteslim{i} H^p(M_i \otimes_R^L X)$, as desired.
\end{proof}
\begin{rem}
In the context of Lemma~\ref{homology of derived tensor product}, the cohomology %
$H^p(M_i \otimes_R^L X)$ does not change depending on whether we regard $M_i \otimes_R^L$
as a functor $F: D^-_{\fp}(R) \to D^-(\cA)$
as in~Lemma~\ref{lem:derived-tensor-product-EW}~\eqref{item:90}
or as a functor
$F': D^-_{\fp}(R) \to \Pro D^b(\cA)$
as in~Lemma~\ref{lem:derived-tensor-product-EW}~\eqref{item:92}.
In fact, in~\eqref{eqn:Pro diagram} we have described a $t$-exact functor
$pj: D^-(\cA) \to \Pro D^b(\cA)$ inducing the natural inclusion on hearts.
Now $pjF$ is right $t$-exact and sends~$R$ to~$M_i$, and so it is isomorphic to~$F'$.
\end{rem}

There are several compatibilities between the functors in Lemma~\ref{lem:derived-tensor-product-EW}.
Most of these follow straightforwardly from our constructions, in the manner of the following lemma.

\begin{lem}\label{upper triangle in diagram}
  Let~$R$ be a Noetherian profinite $\cO$-algebra, and let~$\cA$ be an abelian category.
Let~$M$ be a complete right $R$-module in $\Pro(\cA)$.
Then the composite
\[
D^b_{\fp}(R) \xrightarrow{\eqref{eqn:derived fp to pro fl}} \Pro D^b_{\fl}(R) \xrightarrow{\eqref{item:94}} \Pro D^b(\cA)
\]
is naturally isomorphic to~\eqref{item:87}.
\end{lem}
\begin{proof}
As explained in Remark~\ref{rem:derived fp to pro fl}, the first arrow is $t$-exact and maps~$R$ to $\varprojlim_n R/\rad(R)^n$. 
By construction, the second arrow is right $t$-exact and sends $\varprojlim_n R/\rad(R)^n$ to~$M$.
So the composite is right $t$-exact and sends~$R$ to~$M$, as desired.
\end{proof}

There is however one subtlety, relating to the compatibility of the functors of parts~\eqref{item:87} and ~\eqref{item:88} of Lemma~\ref{lem:derived-tensor-product-EW}.
Considering this compatibility leads to the following definition and lemmas.
\begin{defn}
  \label{defn:derived-complete-object}Let~$R$ be a Noetherian profinite  $\cO$-algebra, and let~$\cA$ be an abelian category. Let~$M$ be a right $R$-module in ~$\Pro(\cA)$.
  Then we say that $M$ is \emph{derived complete} if the natural morphisms $M\to M\otimes^L_R\bigl(R/\Rad(R)^n\bigr)$  (where the derived tensor product is computed by the functor in part~\eqref{item:87} of Lemma~\ref{lem:derived-tensor-product-EW}) %
  give rise to an isomorphism $M\isoto\limcommand_nM\otimes^L_{R} \bigl(R/\Rad(R)^n\bigr)$  in $\Pro D^b(\cA)$.
\end{defn}%

\begin{lem}
  \label{lem:derived-complete-implies-complete}Let~$R$ be a Noetherian profinite  $\cO$-algebra, and let~$\cA$ be an abelian category. Let~$M$ be a derived complete right $R$-module in ~$\Pro(\cA)$.
  Then~$M$ is complete.  
\end{lem}
\begin{proof}
  By the definition of the $t$-structure on $\Pro(\cA)$, the functor~$\tau^{\geq 0}$ preserves cofiltered limits, and so we have
  \begin{align*}
    M=\tau^{\ge 0}M & \isoto \tau^{\ge 0}\bigl(\limcommand_nM\otimes^L_{R} \bigl(R/\Rad(R)^n\bigr)\bigr)\\
    & \isoto \limcommand_n\tau^{\ge 0}\bigl(M\otimes^L_{R} \bigl(R/\Rad(R)^n\bigr)\\
    & \isoto \limcommand_nM\otimes_{R} \bigl(R/\Rad(R)^n\bigr),
  \end{align*}as required.
\end{proof}

We now consider the converse to Lemma~\ref{lem:derived-complete-implies-complete}.
We begin with the following lemma.
\begin{lem}
  \label{lem:equivalence-of-R-and-E-complete}Let~$R$ be a commutative Noetherian profinite $\cO$-algebra, and let $E$ be a finite $R$-algebra \emph{(}i.e.\ $E$ is finite as an $R$-module\emph{)}.
Let~$\cA$ be an abelian category, and let $M$ be a right $E$-module in $\Pro(\cA)$.
\begin{enumerate}
\item\label{item:97} $M$ is complete if and only if it is complete as an $R$-module.
\item\label{item:98} Suppose furthermore that $E$ is flat over~$R$. Then $M$ is derived complete if and only if it is derived complete as an $R$-module.
\end{enumerate}
\end{lem}
\begin{proof}
  Writing \[M\otimes_RR/\Rad(R)^n=M\otimes_E(E\otimes_RR/\Rad(R)^{n})=M\otimes_EE/\Rad(R)^{n}E,\] 
  and bearing in mind Lemma~\ref{properties of compact modules}~\eqref{item: compact 7}, 
  the first part follows from the cofinality of~$(\Rad(R)^nE)_{n \geq 0}$ and $(\Rad(E)^{n})_{n \geq 0}$ (which in turn follows from the finiteness of~$E$ over~$R$).

  The second part follows in the same way, noting that if~$E$ is flat over~$R$ then we have
  \[M\otimes^L_RR/\Rad(R)^n=M\otimes^L_E(E\otimes^L_RR/\Rad(R)^{n})=M\otimes^{L}_EE/\Rad(R)^{n}E. \qedhere\]
\end{proof}

The following lemma presumably holds in greater generality, but for simplicity, we will only prove it in the case of profinite $\cO$-algebras that are finite flat over a commutative ring.

\begin{lem}\label{complete if and only if derived complete}
Let~$R$ be a commutative Noetherian profinite $\cO$-algebra, and let $E$ be a finite flat $R$-algebra.
Let~$\cA$ be an abelian category, and let~$M$ be a complete right $E$-module in~$\Pro(\cA)$.
Then~$M$ is derived complete.
\end{lem}
\begin{proof}
By Lemma~\ref{lem:equivalence-of-R-and-E-complete}, we can assume without loss of generality that~$E = R$ is commutative.
Let~$\fm \coloneqq  \rad(R)$.
We need to prove that 
$\lim_n M \otimes_R^L R/\fm^n$ is concentrated in degree zero.
Note that this coincides with $M \ccotimes_R^L \lim_n R/\fm^n$ (where~$M \ccotimes_R^L$ is the functor in Lemma~\ref{lem:derived-tensor-product-EW}~\eqref{item:93}).
Choose generators of~$\fm$, say $\fm = (f_1, \ldots, f_r)$, and write~$\fm_n \coloneqq  (f_1^n, \ldots, f_r^n)$.
Then the pro-objects $\lim_n R/\fm^n$ and $\lim_n R/\fm_n$ are isomorphic.
So it suffices to prove that
$M \ccotimes_R^L \lim_n R/\fm_n = \lim_n M \otimes_R^L R/\fm_n$
is concentrated in degree zero, i.e.\ the pro-object 
$\lim_n H^p(M \otimes_R^L R/\fm_n)$ is zero for all~$p\ne 0$.
Since~$M$ is complete, Lemma~\ref{homology of derived tensor product} allows us to replace~$M$ by~$M/\fm^N M$, and so to assume that
$M$ is an $R/\fm^N R$-module, for some~$N>0$.

Let~$K_n$ be the Koszul complex on $(f_1^n, \ldots, f_r^n)$, normalized as a cochain complex 
in degrees $-r, \ldots, -1, 0$
with $\wedge^t R^{\oplus r}$ in degree~$-t$.
There is a map of complexes $K_{r+i} \to K_r$, given by $\wedge^t(\diag(f_1^i, \ldots, f_r^i))$ in degree~$-t$, and inducing the natural map
$R/\fm_{r+i} \to R/\fm_r$ on~$H^0$.
Consider now the inverse system of complexes in $\Pro(\cA)$ given by $M \otimes_R K_r$, 
Since $M$ is an $R/\fm^N R$-module, 
the transition map $M \otimes_R K_{r+N} \to M \otimes_R K_r$ is zero for any~$r$ in all degrees~$t \ne 0$,
because it is given by a matrix with entries in $\fm^N$.
Hence the pro-object $\lim_rH^t(M \otimes_R K_r) \in \Pro(\cA)$ is zero for all~$t \ne 0$.
Since~$K_r$ is a complex of finite projective $R$-modules, the complex~$M\otimes_R K_r$ represents
the image of~$K_r$ under the derived functor
$D^-(\Mod^{\fp}(R)) \to D^-(\Pro(\cA))$ of $M \otimes_R\text{--}$. %
Since $M \otimes_R^L$ is the composite of this derived functor with $D^-(\Pro(\cA)) \to \Pro D^b(\cA)$,
which is $t$-exact and induces an equivalence on hearts,
we conclude that $\lim_rH^t(M \otimes^L_R K_r) = 0$ for all~$t \ne 0$.

Of course, $K_r$ need not be a projective resolution of~$R/\fm_r$.
However, we can define an inverse system of fibre sequences $K_r' \to K_r \to R/\fm_r[0]$ in $D^b_{\fp}(R)$,
and it is proved in \cite[\href{https://stacks.math.columbia.edu/tag/0921}{Tag~0921}]{stacks-project}
that for all~$r$ there exists~$i$ such that $K_{r+i}' \to K_r'$ is zero.
This shows that $\lim_r K_r' = 0$, and so $\lim_r K_r \isoto \lim_r R/\fm_r[0]$ in $\Pro D^b_{\fp}(R)$.
Hence 
\[
\lim_r (M \otimes_R^L R/\fm_r[0]) = M \ccotimes_R^L \lim_r R/\fm_r[0] = M \ccotimes_R^L \lim_r K_r = \lim_r M \otimes_R^L K_r.
\]
The previous paragraph shows that $\lim_r M \otimes_R^L K_r$ is concentrated in degree zero, hence so is $\lim_r (M \otimes_R^L R/\fm_r[0])$,
as desired.
\end{proof}

The following result is the reason for introducing Definition~\ref{defn:derived-complete-object}.

\begin{lem}
  \label{lem:derived-complete-commutative-diagram}Let~$R$ be a Noetherian profinite  $\cO$-algebra, and let~$\cA$ be an abelian category. Let~$M$ be a derived complete right $R$-module in ~$\Pro(\cA)$.
Then we have a commutative diagram
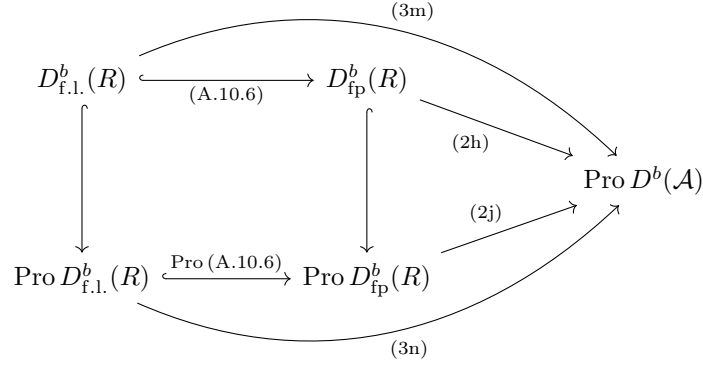
\[\begin{tikzcd}
D^b_{\fl}(R)	 && D^b_{\fp}(R) \\
	&&&& \Pro D^b(\cA) \\
	\Pro D^b_{\fl}(R) && \Pro D^b_{\fp}(R)
	\arrow["\eqref{eqn:fl to fp}"',hook, from=1-1, to=1-3]
	\arrow["\eqref{item:88}", curve={height=-50pt}, from=1-1, to=2-5]
	\arrow[from=1-1,hook, to=3-1]
	\arrow["\eqref{item:87}"', from=1-3, to=2-5]
	\arrow[from=1-3,hook, to=3-3]
	\arrow["\eqref{item:94}"', curve={height=50pt}, from=3-1, to=2-5]
	\arrow["\Pro\eqref{eqn:fl to fp}",hook, from=3-1, to=3-3]
	\arrow["\eqref{item:93}", from=3-3, to=2-5]
\end{tikzcd}\]
\end{lem}
\begin{proof}%
  Note firstly that by Lemma~\ref{lem:derived-complete-implies-complete}, all of the functors in the diagram are defined. 
  It evidently suffices to check that the functor \[M\cotimes_R^{L}\text{--}:\Pro D^b_{\fl}(R)\to\Pro D^b(\cA)\] of~\eqref{item:94} is given by the composite of the embedding
\[\Pro D^b_{\fl}(R) \into \Pro D^b_{\fp}(R)
\] deduced from~\eqref{eqn:fl to fp} and the functor \[M\cotimes_R^{L}\text{--}:\Pro D^b_{\fp}(R)\to\Pro D^b(\cA)\] of~\eqref{item:93}.
Since each of these functors $\Pro D^b_{\fl}(R)\to\Pro D^b(\cA)$ is right $t$-exact and limit-preserving, it suffices by Lemma~\ref{lem:derived-tensor-product-EW}~\eqref{item:85} to show that they agree on the object~$R$ of~$\Mod_c(R)=(\Pro D^b_{\fl}(R))^{\heartsuit}$.
By their definitions, they respectively take~$R$ to~$M$ and to $\limcommand_nM\otimes^L_{R} \bigl(R/\Rad(R)^n\bigr)$, and since~$M$ is derived complete by assumption, we are done.
\end{proof}

We have the following  criterion for the full faithfulness of a derived tensor product.

 \begin{lem}\label{lem:fp-tensor-full-faithfulness-criterion}%
Let ~$R$ be a coherent ring, and let $\cA$ be an abelian category.
Suppose that~$M$ is a right $R$-module in~$\Pro(\cA)$, and that the natural morphism 
\[
R\to\REnd_{\Pro D^b(\cA)}(M) 
\]
is an isomorphism.
Then the functor $M\otimes_{R}^{L}\text{--}:D^b_{\fp}(R)\to \Pro D^{b}(\cA)$ defined in {\em Lemma~\ref{lem:derived-tensor-product-EW}}~\eqref{item:87} is fully faithful.
\end{lem}
\begin{proof}By Lemma %
  ~\ref{lem:F-fully-faithful-iff-IndProF}, it is equivalent to prove that the %
  induced continuous functor \[D(R)=\Ind D^b_{\fp}(R)\to \Ind\Pro D^{b}(\cA)\] is fully faithful. 
Since ~$D(R)$ is compactly generated by~$R$, the result is immediate from
Proposition~\ref{prop:check-full-faithful-on-compact-generators}, together with the full faithfulness of $\Pro D^b(\cA) \to \Ind \Pro D^b(\cA)$.  
\end{proof}

\subsubsection{Topological flatness}
We now introduce various notions of flatness, in the context of Lemma~\ref{lem:derived-tensor-product-EW}.

\begin{defn}\label{topologically flat object}
Let~$R$ be a Noetherian profinite $\cO$-algebra, and let~$\cA$ be an abelian category.
\begin{enumerate}
\item Let~$M$ be a right $R$-module in~$\cA$. The \emph{$\Tor$-dimension of~$M$} is the amplitude of the right $t$-exact functor 
$M \otimes_R^L\text{--} : D^-_{\fp}(R) \to D^-(\cA)$, i.e.\ the smallest~$n\ge 0$
such that the restriction of $M \otimes_R^L\text{--}$ to~$D^-_{\fp}(R)^{\heartsuit}$
(i.e.\ the abelian category $\Mod^{\fp}(R)$) factors through $D^-(\cA)^{\geq -n}$.   
\item Let~$M$ be a right $R$-module in~$\cA$. We say that~$M$ is \emph{flat} if $M \otimes_R\text{--} : \Mod^{\fp}(R) \to \cA$ is exact, or equivalently, if $M$ has $\Tor$-dimension zero.
\item Let~$M$ be a complete right $R$-module in  $\Pro(\cA)$.
We say that~$M$ is \emph{topologically flat} if $M \cotimes_R\text{--} : \Mod_c(R) \to \Pro(\cA)$ in Lemma~\ref{lem:EW-gives-completed-tensor-product} is exact.
\end{enumerate}
\end{defn}

 \begin{lem}\label{topologically flat equivalent to flat}
Let~$R$ be a Noetherian profinite $\cO$-algebra, and let~$\cA$ be an abelian category. 
Let~$M$ be a complete right $R$-module in~$\Pro (\cA)$.
Then~$M$ is flat if and only if~$M$ is topologically flat.
\end{lem}
\begin{proof}
  By Corollary~\ref{cor:restricting-completed-tensor}, the restriction of $M \cotimes_R\text{--}$ to $\Mod^{\fp}(R)$ is $M \otimes_R\text{--}$.
  Hence, if $M$ is topologically flat, then~$M$ is flat.
  On the other hand, the inclusion $\Mod^{\fl}(R) \to \Mod_c(R)$ factors through $\Mod^{\fp}(R)$.
  Hence, if $M$ is flat, then $M \cotimes_R\text{--}$ restricts to an exact functor on $\Mod^{\fl}(R)$, and since $\Mod_c(R) = \Pro \Mod^{\fl}(R)$, we see that $M \cotimes_R\text{--}$ is %
  exact, by Lemma~\ref{lem:exactness-of-limits}~(2) and~(3).
  Thus~ $M$ is topologically flat, as desired.
\end{proof}

\begin{lem}\label{topological flatness and t-exactness}
Let~$R$ be a Noetherian profinite $\cO$-algebra, and let~$\cA$ be an abelian category.
Let~$M$ be a complete right $R$-module in~$\Pro(\cA)$.
If~$M$ is topologically flat, then the functor 
\[
M \cotimes_R^L\text{--}: \Pro D^b_{\fl}(R) \to \Pro D^b(\cA),
\] 
defined in {\em Lemma~\ref{lem:derived-tensor-product-EW}}~\eqref{item:94}, is $t$-exact.
\end{lem}
\begin{proof}
By definition, $M \cotimes_R^L \text{--}$ is 
the Pro-extension of the composite
\[
D^b_{\fl}(R) = D^b(\Mod^{\fl}(R)) \xrightarrow{i} D^-(\Mod_c(R)) \to \Pro D^b(\cA), 
\]
where the first arrow is induced by the inclusion $\Mod^{\fl}(R) \subset \Mod_c(R)$ (and is $t$-exact) and the second arrow is the derived functor of
$M \cotimes_R\text{--}$.
Since~$M$ is topologically flat, $M \cotimes_R\text{--}$ is exact, and so both arrows are $t$-exact.
The lemma thus follows from Corollary~\ref{cor:pro t-exact derived functors}.
\end{proof}

\begin{lem}\label{derived tensor and finite global dimension}
Let~$R$ be a Noetherian profinite $\cO$-algebra, and assume that~$R$ has finite global dimension (i.e.\ there exists~$n$ such that every left $R$-module has a projective resolution
of length~$\leq n$).
Let~$\cA$ be an abelian category, and let $M$ be a right $R$-module in~$\cA$.
Then the functor $M \otimes_R^L\text{--} : D^b_{\fp}(R) \to D^-(\cA)$ defined in~\eqref{item:89} factors through $D^b(\cA)$.
\end{lem}
\begin{proof}
An induction on the amplitude of~$X \in D^b_{\fp}(R)$ shows
that it suffices to check that $M \otimes_R^L X \in D^b(\cA)$ whenever~$X \in \Mod^{\fp}(R)$.
However, by assumption, $M \otimes_R^L X$ is represented by a complex concentrated in degrees~$[-n, 0]$.
\end{proof}

\section{Coherent sheaves on formal algebraic stacks}
\label{sheaves on formal algebraic stacks}
In this appendix, after recalling some of the basic framework and results related
to the theory of coherent sheaves on algebraic stacks,
we develop various constructions and results on
coherent and pro-coherent
sheaves on formal algebraic stacks. These include assorted pullback and pushforward
functors, as well as results related to the theorem on formal functions.

\subsection{Coherent sheaf theory on algebraic stacks}
\label{subsec:coherent sheaves}
There are various competing approaches to defining both the abelian and
the derived category of quasi-coherent sheaves
on an algebraic stack; see \cite[App.~B.2]{emerton2023introduction} for a discussion of
some of them.  In the context of a Noetherian algebraic stack~$\cX$
with affine diagonal, it follows from ~\cite[Prop.\ 1.3, Rem.\
1.5]{Hall_2017} that
these various approaches lead to an unambiguous definition of~$D^b_{\coh}(\cX)$,
the stable $\infty$-category of bounded complexes of $\cO_{\cX}$-modules having
coherent cohomology. A similar discussion applies to quasi-coherent cohomology, although
the coherent case will be the most important one for us in the present paper.

For definiteness, we will use sheaves on the lisse-\'etale
site of~$\cX$ as our model for $D^b_{\coh}(\cX)$,
and the various related categories of sheaves that we will have to consider.
We make this choice primarily to make contact with various pieces 
of literature (e.g.~\cite{BCGAGA} and~\cite{alper2023coherently})
whose results we are going to cite.

As in the preceding discussion, we assume that~$\cX$ is a Noetherian algebraic
stack with affine diagonal, and we let $\cO_{\cX}$ denote the usual structure sheaf
on the lisse-\'etale site of~$\cX$.
We let $\Mod(\cO_{\cX})$ denote the abelian category 
of sheaves of $\cO_{\cX}$-modules on the lisse-\'etale site of~$\cX$, 
and we let $D(\cO_{\cX})$ denote the derived $\infty$-category of~$\Mod(\cO_{\cX})$.

We let $\Coh(\cX)$, resp.\ ~$\QCoh(\cX)$, denote the abelian category of coherent, 
resp.\ ~quasi-coherent, sheaves of $\cO_{\cX}$-modules (see e.g.\ \cite[\S 6]{MR2312554}).
These are both weak Serre subcategories of~$\Mod(\cO_{\cX})$, %
and $\Coh(\cX)$ is a Serre subcategory of $\QCoh(\cX)$ (see e.g.\ \cite[\href{https://stacks.math.columbia.edu/tag/07B4}{Tag 07B4}]{stacks-project}, \cite[\href{https://stacks.math.columbia.edu/tag/0GRB}{Tag 0GRB}]{stacks-project}).
We have the following standard lemma about the structure of~$\QCoh(\cX)$.

\begin{lemma}\label{QCoh of a Noetherian algebraic stack is locally coherent}
Let~$\cX$ be a Noetherian algebraic stack.
Then~$\QCoh(\cX)$ is a locally coherent Grothendieck category, and $\QCoh(\cX)^c = \Coh(\cX)$.
\end{lemma}
\begin{proof}
$\QCoh(\cX)$ is a Grothendieck category by~\cite[\href{https://stacks.math.columbia.edu/tag/0781}{Tag~0781}]{stacks-project}. 
To see that coherent sheaves are compact in~$\QCoh(\cX)$, choose a smooth surjective morphism $p: X \to \cX$  with~$X$ a Noetherian scheme. 
Then~$p^*$ preserves colimits, and sends coherent sheaves on~$\cX$ to coherent sheaves on~$X$.
It then follows from~\cite[\href{https://stacks.math.columbia.edu/tag/06WT}{Tag~06WT}]{stacks-project} that every $\cF \in \Coh(\cX)$ is compact in $\QCoh(\cX)$, since~$p^*\cF$ is compact
in~$\QCoh(X)$.
Conversely, compact objects of~$\QCoh(\cX)$ are coherent, because every object of $\QCoh(\cX)$ is the filtered colimit of its coherent subsheaves (\cite[Proposition~15.4]{MR1771927} or \cite[\href{https://stacks.math.columbia.edu/tag/0GRF}{Tag 0GRF}]{stacks-project}).
\end{proof}

In particular, the embeddings $\QCoh(\cX) \hookrightarrow \Mod(\cO_{\cX})$
and $\Coh(\cX) \hookrightarrow \Mod(\cO_{\cX})$
are both exact,
and so induce $t$-exact functors of associated derived $\infty$-categories
$D\bigl(\QCoh(\cX)\bigr) \to D_{\qc}(\cX)$ 
and
$D\bigl(\Coh(\cX)\bigr) \to D_{\coh}(\cX),$ 
where the targets denote the full subcategories of $D(\cO_{\cX})$
consisting of objects whose cohomology sheaves lie in~$\QCoh(\cX)$, resp.\ ~$\Coh(\cX)$.
The first of these functors restricts to an equivalence
\begin{equation}
\label{eqn:quasi-coherent equivalence}
D^+\bigl(\QCoh(\cX)\bigr) \iso D^+_{\qc}(\cX),
\end{equation}
by~\cite[Prop.~1.6]{MR3720855} or~\cite[Thm.~C.1]{MR3995721}; 
note however that the first reference uses a slightly different framework
to the one we are using here, which is explained and compared with our framework of lisse-\'etale sites in~\cite{MR3431476}. 
By Lemma~\ref{lem:compact-objects-bounded-Serre-category} (whose hypotheses hold because of Lemma~\ref{QCoh of a Noetherian algebraic stack is locally coherent})
the equivalence~\eqref{eqn:quasi-coherent equivalence}
restricts to a $t$-exact equivalence
\begin{equation}
\label{eqn:coherent equivalence}
D^b\bigl(\Coh(\cX)\bigr) \iso D^b_{\coh}(\cX).
\end{equation}
In the body of the text, we will not work directly with quasi-coherent 
sheaves that are not coherent. Rather, we will work with Ind-coherent sheaves;
that is, we will consider the stable $\infty$-category~$\Ind D^b_{\coh}(\cX).$
Proposition~\ref{prop:Ind t-structures}
shows that the $t$-structure on $D^b_{\coh}(\cX)$
extends canonically to a $t$-structure on~$\Ind D^b_{\coh}(\cX)$. 
Since $D_{\qc}(\cX)$ %
is cocomplete, there is a
$t$-exact functor $\Ind D^b_{\coh}(\cX) \to D_{\qc}(\cX),$
which restricts to an equivalence
\begin{equation}
\label{eqn:ind-coh to qcoh}
\bigl(\Ind D^b_{\coh}(\cX)\bigr)^+ \iso D^+_{\qc}(\cX)
\end{equation}(see Remark~\ref{rem:more-general-IndCoh} below).
Then, when it comes to a consideration of
functors such as $Rj_*$ for an open immersion $j: \cY \to \cX$,
which are usually defined on the categories $D^+_{\qc}$,
we will extend them
to the categories $\Ind D^b_{\coh}(\cX)$; furthermore, we will often construct such functors
using this Ind-coherent perspective, as in the discussion of Section~\ref{subsubsec:closed
  pushforward} below.

\begin{rem}
  \label{rem:more-general-IndCoh}In the literature one finds a definition of a stable infinity category~$\Ind\Coh(\cX)$ which is not necessarily compactly generated, but is rather defined by a universal mapping property, which incorporates~\eqref{eqn:ind-coh to qcoh}; see for example~\cite[Defn.\ 5.10]{cautis2024indgeometricstacks}. 
  However, in our setting this more sophisticated definition agrees with our $\Ind D^b_{\coh}(\cX)$; for example, all of the stacks that we consider in this paper are “coherent ind-geometric stacks” in 
  the sense of~ \cite{cautis2024indgeometricstacks}, so this agreement follows from \cite[Prop.\ 5.30]{cautis2024indgeometricstacks}.
\end{rem}

\subsubsection{Operations on coherent and quasi-coherent sheaves}
\label{subsubsec:operations}
Let $\cX$ and $\cY$ be algebraic stacks,
both assumed to be Noetherian and having affine diagonal,
and let $f: \cY \to \cX$ be a morphism.
In this context one can define pushforward and pullback functors on both the abelian
and derived categories of quasi-coherent sheaves, as well on various related categories.
We now recall some details of these constructions.

\begin{remark}
The hypotheses we have placed on $\cX$ and $\cY$ are certainly
unnecessarily restrictive; our goal here is simply to give some details
regarding, and references for, the construction
of these functors in the specific framework that we have adopted. 
One motivation for taking care on this point is that we wish to combine
citations to works on sheaf theory on algebraic stacks, such as~\cite{alper2023etale, alper2023coherently} 
with the $\infty$-categorical framework, and so
we want to ensure that the functors we compare in these various different contexts
are indeed the ones we imagine they are.  Another is that we will extend some
of these constructions to the context of formal algebraic stacks, where there 
are fewer existing references to rely on, and we want to be sure that
our constructions in this context are well-founded (for {\em some} choice
of foundations).
\end{remark}

The construction of pullbacks and pushforwards is simplest in the case
when $f$ is representable by schemes. (One simple but important case when this holds
is when $f$ is an immersion.)
In this case $U \mapsto U\times_{\cX} \cY$ gives a continuous functor~$u$
(in the sense discussed
in~\cite[\href{https://stacks.math.columbia.edu/tag/00WU}{Tag 00WU}]{stacks-project})
from the lisse-\'etale site of $\cX$ to the lisse-\'etale site of~$\cY$.
Adopting notation from \emph{loc.\ cit.}, we thus obtain a functor $u^s$
from the category of sheaves on the lisse-\'etale site
of~$\cY$ to the category of sheaves on the lisse-\'etale site of~$\cX$.

This functor $u$ is the restriction to lisse-\'etale sites of
the analogously defined functor on big flat sites. This latter functor
induces a morphism of (ringed) 
topoi~\cite[\href{https://stacks.math.columbia.edu/tag/06W8}{Tag 06W8}]{stacks-project}.
Unfortunately, the functor $u$ on lisse-\'etale sites
typically {\em does not} define a morphism of topoi:
the left adjoint $u_s$ to $u^s$  
is typically not exact,
see e.g.~\cite[\href{https://stacks.math.columbia.edu/tag/07BF}{Tag 07BF}]{stacks-project}.

Nevertheless, the morphism $f$ induces functors
$f_*:\QCoh(\cY) \to \QCoh(\cX)$ and $f^*:\QCoh(\cX) \to \QCoh(\cY)$. 
In the Stacks Project, these are defined using big sites: the pullback
$f^*$ is induced by the left-adjoint of the ``big site version''
of the functor $u^s$ considered above, and the push-forward $f_*$ 
being defined via slightly elaborate manipulations, for which 
see~\cite[\href{https://stacks.math.columbia.edu/tag/070A}{Tag 070A}]{stacks-project}.
In the context of small lisse-\'etale sites,
however, we can describe $f_*$ rather directly.
Namely, there is a canonical morphism $\cO_{\cX} \to u^s(\cO_{\cY})$, 
so that $u^s$ induces a functor $f_*:\Mod(\cO_{\cY}) \to \Mod(\cO_{\cX})$,
which by~\cite[Lem.~6.5]{MR2312554}
restricts to a functor $\QCoh(\cY)\to \QCoh(\cX)$.
(Note that our assumption that $\cY$ and $\cX$ are Noetherian
implies that $f$ is necessarily quasi-compact.) 
If $f$ is furthermore proper, then $f_*$ restricts to a functor $\Coh(\cY) \to \Coh(\cX)$:
see~\cite[Thm.~10.13]{MR2312554} for a proof of the (stronger) derived statement,
although note that (as the proof 
of this result indicates) under the assumption that $f$ is representable by schemes,
this result follows directly from the corresponding result for schemes.  

We may then construct the left adjoint $f^*:\Mod(\cO_{\cX}) \to 
\Mod(\cO_{\cY})$ of~$f_*$, via composing the left adjoint $u_s$ to $u^s$ with the extension of
scalars from $u_s\cO_{\cX}$ to $\cO_{\cY}$. %
The functor $f^*$ then restricts
to functors $\QCoh(\cX) \to \QCoh(\cY)$ and $\Coh(\cX) \to \Coh(\cY).$ 
(For this material, see again~\cite[Lem.~6.5]{MR2312554}.)

We may right derive the left exact functor $f_*$ to obtain a functor
$Rf_*: D^+(\cY) \to D^+(\cX),$
which restricts to a functor
$Rf_*: D^+_{\qc}(\cY) \to D^+_{\qc}(\cX)$,
and even to a functor
$Rf_*: D^+_{\coh}(\cY) \to D^+_{\coh}(\cX)$
if~$f$ is proper. (See~\cite[Lem.~6.20, Thm.~10.13]{MR2312554}.) Recall 
that in general $Rf_*$ need not have bounded amplitude, and consequently
need not take $D_{\qc}(\cY)$ to $D_{\qc}(\cX)$;
see e.g.\ ~\cite[\href{https://stacks.math.columbia.edu/tag/07DC}{Tag 07DC}]{stacks-project}.

One may also consider the left-derived functor $Lf^*$. 
The treatment of this functor in the lisse-\'etale formalism of~\cite{MR2312554}
is complicated by the failure of $u$ to induce a morphism
of topoi. 
However,
in the very particular case when $f$ is a smooth morphism (e.g.\ an open immersion, which is the case used in Section~\ref{subsection:localization-to-Ugood}),
$u$ {\em does} induce
a morphism of topoi; indeed, 
the functor $u^s$ is then simply given by restriction along~$f$, and so {\em is} exact
(see the remark at the bottom of~\cite[p.~60]{MR2312554}).
In this case $f^*$ is exact, as it coincides with the exact functor $u_s$,
and by Lemma~\ref{lem:exact-functor-induces-t-exact} it induces
a $t$-exact functor
$D(\cO_{\cX}) \to D(\cO_{\cY})$,
which restricts to $t$-exact functors
$D_{\qc}(\cX) \to D_{\qc}(\cY)$
and
$D_{\coh}(\cX) \to D_{\coh}(\cY).$

If $f:\cX \to \cY$ is not representable by schemes,  
then the construction of pullbacks and pushforwards is more involved;
cf.\ \cite[Rem., p.~488]{MR3431476}. %
The reference~\cite{MR3431476} gives a construction in this case, using a similar framework
to the lisse-\'etale site (namely, a version of the flat site), and building on the methods of~\cite{MR2312554}. 

\subsubsection{Pushforward along an immersion}
\label{subsubsec:closed pushforward}
In the main body of the paper, we will frequently need to consider 
the derived functors of~$f^*$ or~$f_*$ in the case when $f$ is a closed or open immersion, which we will
typically denote by $i$ or $j$ respectively --- possibly with additional decorations
--- rather than~$f$.
If $i$ is a closed immersion, then
$i_*$ is exact: this is easily checked directly, by working locally on the
target and so reducing to the case of schemes (see also~\cite[Cor.~6.6]{MR2312554}).
Thus, by Lemma~\ref{lem:right-derived-of-exact-is-exact},
$Ri_*$ is simply the $t$-exact functor induced by $i_*$
following the prescription of Lemma~\ref{lem:exact-functor-induces-t-exact}.
We will thus write simply $i_*$ rather than $Ri_*$.

If $j:\cU \to \cX$ is an open immersion, we can 
give an alternative construction of $Rj_*$ directly as a functor on Ind-coherent complexes. 
To begin with, note that  $j_* : \QCoh(\cU) \to \QCoh(\cX)$ is fully faithful, hence the counit $j^*j_* \to \id_{\QCoh(\cU)}$ is an isomorphism. Since  $j^*:\QCoh(\cX) \to \QCoh(\cU)$ is exact, it follows from ~\cite[Proposition~5, \S III.2]{Gabrielthesis} that~$j^*$ is a Serre quotient 
by its kernel, which is a localizing subcategory of~$\QCoh(\cX)$.
Since (by Lemma~\ref{QCoh of a Noetherian algebraic stack is locally coherent}) the categories $\QCoh(\cX)$ and~$\QCoh(\cU)$ are locally coherent with compact
objects $\Coh(\cX)$ and~$\Coh(\cU)$ respectively, \cite[Theorem~2.6]{MR1426488} shows that~$j^*$
realizes $\Coh(\cU)$ as the quotient of ~$\Coh(\cX)$ 
by its kernel, which is the Serre subcategory~$\Coh_Z(\cX)$, whose objects are those coherent sheaves
on~$\cX$ with set-theoretic support contained in~$Z \coloneqq |\cX| \setminus |\cU|$. %
Lemma~\ref{lem:quotient-of-bounded} then provides a $t$-exact equivalence
\[
D^b\bigl( \Coh(\cX) \bigr) / D^b_{\Coh_Z(\cX)}\bigl(\Coh(\cX) \bigr)
\iso D^b\bigl( \Coh(\cU)\bigr),
\]
which~\eqref{eqn:coherent equivalence}
allows us to interpret as a $t$-exact equivalence
\begin{equation}\label{quotient of D^b_coh by D^b_coh, Z}
D^b_{\coh}(\cX) / D^b_{\coh,Z}(\cX) \iso D^b_{\coh}(\cU),
\end{equation}
where~$D^b_{\coh, Z}(\cX)$ is the full subcategory of complexes whose cohomology is set-theoretically supported on~$Z$.

The functor $D^b_{\coh}(\cX) \to D^b_{\coh}(\cU)$ implicit in~\eqref{quotient of D^b_coh by D^b_coh, Z} is of course simply
the $t$-exact extension of
$j^*$.
The Ind-extension of this functor  
is then a continuous %
$t$-exact functor $j^*: \Ind D^b_{\coh}(\cX) \to \Ind D^b_{\coh}(\cU)$
which restricts
(via the equivalence~\eqref{eqn:ind-coh to qcoh} for each of~$\cX$ and~$\cU$)
to the functor $j^*: D^+_{\qc}(\cX) \to D^+_{\qc}(\cU)$
considered above.   
The right adjoint of $j^*$ is then a functor
\begin{equation}
\label{eqn:Ind derived pushforward}
\Ind D^b_{\coh}(\cU) \to \Ind D^b_{\coh}(\cX)
\end{equation}
which restricts 
(again via the equivalence~\eqref{eqn:ind-coh to qcoh} for each of~$\cX$ and~$\cU$)
to the functor $Rj_*$ considered above.
Since $j^*$ preserves compact objects (equivalently, is constructed via Ind extension
from the corresponding functor on compact objects), its right adjoint~\eqref{eqn:Ind derived pushforward} 
is also continuous, %
by Lemma~\ref{lem: adjoint continuous}.
Thus~\eqref{eqn:Ind derived pushforward} is the Ind-extension of its restriction
$D^b_{\coh}(\cU) \to \Ind D^b_{\coh}(\cX),$ which in turn may be regarded as the composite
$$D^b_{\coh}(\cU) \buildrel Rj_* \over \longrightarrow D^+_{\qc}(\cX) \iso
\bigl(\Ind D^b_{\coh}(\cX)\bigr)^+ \hookrightarrow \Ind D^b_{\coh}(\cX).$$ 
Consequently, the functors $Rj_*$ and~\eqref{eqn:Ind derived pushforward} 
canonically determine one another.

\begin{remark}
\label{rem:cohomologically affine case}
In fact, in the cases of interest to us in the main body of the paper,
the relevant open immersions $j$ will be furthermore cohomologically affine. 
Hence $Rj_*$ will be $t$-exact,
and so, just as in the case of closed
immersions, will be obtained from $j_*$ by an appropriate application
of Lemma~\ref{lem:exact-functor-induces-t-exact}.
From the preceding explanation of how to recover~\eqref{eqn:Ind derived
pushforward} from~$Rj_*$,  
we conclude that~\eqref{eqn:Ind derived pushforward} is also $t$-exact in this case.
Consequently, in this case we will denote both $Rj_*$ and its continuous extension~\eqref{eqn:Ind
derived pushforward} simply by~$j_*$.
\end{remark}

\subsection{Coherent sheaf theory on formal algebraic stacks}
\label{subsec:formal coherent sheaves}
Let $\cX$ be a Noetherian formal algebraic stack with affine diagonal, 
written as a colimit
\begin{equation}\label{eqn:present-formal-alg-as-colimit}\cX \iso \colim_n \cX_n\end{equation}
with the $\cX_n$ being algebraic stacks
and the transition maps being thickenings. 
Note that the $\cX_n$ are then closed substacks
of~$\cX$, and are thus also Noetherian with affine diagonal. 
We then write
\begin{equation}\label{eqn:defn-of-coherent-sheaves}D^b_{\coh}(\cX) \coloneqq  \colim_n D^b_{\coh}(\cX_n),\end{equation}
the transition maps being given by pushforward of sheaves as described in
Section~\ref{subsubsec:closed pushforward} above,
and the colimit being formed in the $\infty$-category of stable
$\infty$-categories.  Since any two presentations of~$\cX$ of the
form~\eqref{eqn:present-formal-alg-as-colimit} are mutually cofinal,
we see that the resulting colimit is (up to canonical equivalence) independent
of the choice of presentation~\eqref{eqn:present-formal-alg-as-colimit}.
Since the transition maps are $t$-exact,
the colimit $D^b_{\coh}(\cX)$ is endowed with a canonical $t$-structure.
Furthermore, each $D^b_{\coh}(\cX_n)$ is idempotent complete, hence the same is true of~$D^b_{\coh}(\cX)$, by~\cite[\href{https://kerodon.net/tag/042P}{Tag~042P}]{kerodon}:
to see the idempotent-completeness of $D^b_{\coh}(\cX_n)$, by~\cite[Lemma~1.2.4.6]{LurieHA} it suffices to prove that the homotopy category of each $D^b_{\coh}(\cX_n)$ is idempotent complete,
which is a consequence of the equivalence~\eqref{eqn:coherent equivalence} and the fact that the bounded derived category of an abelian category is
idempotent complete~\cite[Corollary~2.10]{BalmerSchlichting}.

We let $\Coh(\cX)$ denote the heart of the $t$-structure on $D^b_{\coh}(\cX)$.
The discussion of Section~\ref{subsec:colimits} shows that there are canonical equivalences
\begin{equation}
  \label{eq:cohX-colim-cohXn}
  \colim_n \Coh(\cX_n) \iso \Coh(\cX) \quad \text{ and } \quad
  D^b\bigl(\Coh(\cX)\bigr) \iso D^b_{\coh}(\cX).
\end{equation}
The presentation~\eqref{eqn:present-formal-alg-as-colimit} allows us to define the topological space $|\cX| \coloneqq  \colim_n |\cX_n|$, and then $|\cX_n| \isoto |\cX|$ for any choice of~$n$.
If~$\cF \in \Coh(\cX)$, then it can be represented by an object $\cF_n \in \Coh(\cX_n)$ for some~$n$, and the set-theoretic support $\supp \cF_n \subset |\cX_n|$ is independent of the choice of
representative.
Using the $t$-structure on~$D^b_{\coh}(\cX)$ we can thus associate to each closed subset $Z \subset |\cX|$ of the underlying topological space of~$\cX$ 
(or equivalently, of~$|\cX_n|$, for any~$n$)
a full sub-$\infty$-category $D^b_{\coh, Z}(\cX)$ consisting of objects whose cohomology sheaves are set-theoretically supported on~$Z$.
By construction, \eqref{eqn:defn-of-coherent-sheaves} induces an equivalence
\begin{equation}\label{colimit with fixed support}
\colim_n D^b_{\coh, Z}(\cX_n) \isoto D^b_{\coh, Z}(\cX).
\end{equation}

As well as considering the bounded derived category $D^b_{\coh}(\cX)$
of a Noetherian formal algebraic stack~$\cX$ with affine diagonal,
we will also have occasion to consider the associated Ind and Pro categories
(and even the associated Ind Pro category).
If we present $\cX$ as in~\eqref{eqn:present-formal-alg-as-colimit},
then (as noted in Section~\ref{subsec:colimits}) we have an equivalence 
$$
\colim_n \Ind D^b_{\coh}(\cX_n)
\iso
\Ind D^b_{\coh}(\cX), 
$$
the colimit now being formed in the $\infty$-category of compactly generated stable
$\infty$-categories (whose morphisms are the continuous functors). 
Since $D^b_{\coh}(\cX)$ is idempotent complete,
it is the subcategory of compact
objects in $\Ind D^b_{\coh}(\cX)$.
We also note the following result. 

\begin{lem}
  \label{lem:detecting-coh-in-pro-coh}Let~$\cX$ be a Noetherian formal algebraic
  stack with affine diagonal. Then $D^b_{\coh}(\cX)$ is the full subcategory of
  $\Ind D^b_{\coh}(\cX)$ {\em (}resp.\ ~$\Pro D^b_{\coh}(\cX)${\em )} given by those objects
  whose cohomology is coherent and concentrated in finitely many degrees.
\end{lem}
\begin{proof}
  We give the argument in the Ind case, the argument in the Pro case being dual to it. %
  It suffices to prove that if $x \in \Ind D^b_{\coh}(\cX)^{[a, b]}$ is a bounded object with coherent cohomology, then~$x$ is compact in $\Ind D^b_{\coh}(\cX)$.
  When~$a = b$ we have an isomorphism $x \isoto H^b(x)[-b]$, and $H^b(x) \in \Ind \Coh(\cX)$ is assumed to be contained in $\Coh(\cX) \subset D^b_{\coh}(\cX)$.
  Hence~$H^b(x)$ is compact, and so is~$x$.
  In general, by induction on~$b-a$, we see that~$x$ is part of a cofibre sequence $\tau^{\leq (b-1)}x \to x \to H^b(x)[-b]$
  where the first and third term are compact, so~$x$ is also compact, as required.
\end{proof}

\subsubsection{Completion along a closed substack}
\label{subssubsec:completion-along-substack-formal-algebraic}
Let $\cX$ be a Noetherian formal algebraic stack with affine diagonal,
and let $\cZ \hookrightarrow \cX$ be a closed substack of~$\cX$.
We define the completion of $\cX$ along~$\cZ$ to be
the substack $\widehat{\cX}=\widehat{\cX}_{\cZ}$ of $\cX$ that classifies
morphisms $T \to \cX$ for which the induced map on underlying topological spaces
$|T| \to |\cX|$ has image landing in $|\cZ|$. By definition, then,
$\widehat{\cX}$ depends only on the closed subset $Z\coloneqq  |\cZ|$ of~$|\cX|$,
not on its particular realization as a closed substack $\cZ$ of~$\cX$, and we will sometimes write $\widehat{\cX}_{Z}$ rather than~$\widehat{\cX}_{\cZ}$.
Note that $\widehat{\cX}$ is again a Noetherian formal algebraic stack.
To see this, note for example that~\cite[Lem~8.13]{Emertonformalstacks} implies that $\widehat{\cX}$ is locally
Noetherian; since $\cX$ is furthermore assumed Noetherian,
i.e.\ to be quasi-compact and quasi-separated in addition to being locally Noetherian,
it is straightforward to deduce that $\widehat{\cX}$ is also Noetherian.

If we choose a presentation~$\widehat{\cX} \iso \colim_n \cZ_n$
for $\widehat{\cX}$ as a colimit with respect to thickenings of Noetherian algebraic
stacks~$\cZ_n$,
then we may choose our presentation~\eqref{eqn:present-formal-alg-as-colimit}
for~$\cX$ so that each composite %
$$\cZ_n \hookrightarrow \widehat{\cX} \rightarrow \cX$$
factors through~$\cX_n$, inducing a closed immersion
$$i_n:\cZ_n \hookrightarrow \cX_n.$$ 
The canonical monomorphism 
$i: \widehat{\cX} \hookrightarrow \cX$ may then be written as the colimit
of the closed immersions~$i_n$.
Passing to the colimit,
the ($t$-exact) pushforward functors $i_{n,*}: D^b_{\coh}(\cZ_n) \to D^b_{\coh}(\cX_n)$
then induce a $t$-exact functor
\begin{equation}
\label{eqn:formal pushforward}
i_*: D^b_{\coh}(\widehat{\cX}) \to D^b_{\coh}(\cX).
\end{equation}

\begin{aprop}%
\label{aprop:Coh-set-theoretically-supported}
Let $\cX$ be a Noetherian formal algebraic stack with affine diagonal,
let $\cZ$ be a closed substack of~$\cX$,
and let $\widehat{\cX}$ denote the completion of~$\cX$ along~$\cZ$.
Then~{\em \eqref{eqn:formal pushforward}} %
is a $t$-exact fully faithful functor %
$D^b_{\coh}(\widehat{\cX}) \hookrightarrow D^b_{\coh}(\cX),$
whose essential image coincides with
$D^b_{\coh,Z}(\cX),$ i.e.\
the full sub-$\infty$-category of $D^b_{\coh}(\cX)$ consisting of
objects whose cohomology sheaves are set-theoretically supported on~$Z \coloneqq  |\cZ|$.
\end{aprop}
\begin{proof}
We may replace $\cZ$ by $\cZ_{\red}$, which is a closed algebraic substack of~$\cX$.
Then, if we choose a presentation $\cX \iso \colim_n \cX_n$
as in the discussion above,
we see that $\cZ$ is in fact a closed substack of each~$\cX_n$,
that 
\[
\widehat{\cX} \iso \colim_n \widehat{\cX_n},
\]
where all the indicated completions are taken along~$\cZ$,
and that 
\[
D^b_{\coh}(\widehat{\cX}) \iso \colim_n D^b_{\coh}(\widehat{\cX_n}),
\]
where the colimit is formed in the $\infty$-category of stable $\infty$-categories.

Since a filtered colimit of fully faithful $t$-exact morphisms is fully faithful
and $t$-exact, and
\[
D^b_{\coh, Z}(\cX) = \colim_n D^b_{\coh, Z}(\cX_n)
\]
by~\eqref{colimit with fixed support}, to prove that~\eqref{eqn:formal pushforward} is fully faithful with essential image $D^b_{\coh, Z}(\cX)$, it then suffices to show that pushforward induces a $t$-exact equivalence
\begin{equation}\label{equivalence at finite stages}
D^b_{\coh}(\widehat{\cX_n}) \isoto D^b_{\coh, Z}(\cX_n)
\end{equation}
for each~$n$.
This allows us to assume that $\cX = \cX_n$ for some~$n$, i.e.\ that $\cX$ 
is algebraic.

In this case, let $\cI$ denote the coherent ideal sheaf cutting out~$\cZ$ in~$\cX$, %
and let $\cZ^{{(n)}}$ denote the closed substack of $\cX$ cut out by~$\cI^n$.
Then $\widehat{\cX} \iso \colim_n \cZ^{{(n)}}$, and so we have to show that pushforward
induces a fully faithful $t$-exact functor
$$\colim_n D^b_{\coh}(\cZ^{{(n)}}) \hookrightarrow D^b_{\coh}(\cX),$$
with essential image equal to~$D^b_{\coh,Z}(\cX).$
If~$\cX$ is an affine scheme, this is proved in~\cite[Proposition~B.1.9, Proposition~B.1.17]{emerton2023introduction}:
it amounts to the statement that the first two arrows of~\cite[(B.1.8)]{emerton2023introduction}
are equivalences.
The general case can be proved in the same way.
Alternatively, it can be deduced as a special case of~\cite[Theorem~2.2.3]{H-LP}.
\end{proof}

We use the same notation (namely~$i_*$) to denote the extension
of~\eqref{eqn:formal pushforward}
to the corresponding Ind categories;
this is a fully faithful embedding
\begin{equation}
\label{eqn:ind-coherent embedding} 
i_*: \Ind D^b_{\coh}(\widehat{\cX}) \iso
\Ind D^b_{\coh,Z}(\cX)\hookrightarrow \Ind D^b_{\coh}(\cX). 
\end{equation}

\begin{lemma}\label{lem:standard hypothesis for sheaves}
The functor~\eqref{eqn:ind-coherent embedding} satisfies the conditions of
Hypothesis~\ref{hyp:usual hyp for semiorthogonal},
taking $\cA_Z^{c} \coloneqq  D^b_{\coh}(\widehat{\cX})$
and $\cA^{c} \coloneqq   D^b_{\coh}(\cX)$.
\end{lemma}
\begin{proof}
This is an immediate consequence of the construction of~\eqref{eqn:ind-coherent embedding} as an Ind-extension of a fully faithful functor.
\end{proof}
Write $\cU \coloneqq  \cX \setminus \cZ$ (an open substack of~$\cX$),
and let $j: \cU \hookrightarrow \cX$ denote the canonical open immersion.
If we again choose a presentation~\eqref{eqn:present-formal-alg-as-colimit} for~$\cX$,
and write $\cU_n \coloneqq  \cU\times_{\cX} \cX_n,$
then we obtain a compatible presentation $\colim_n \cU_n \iso \cU$ of~$\cU$.
As we noted in~\eqref{quotient of D^b_coh by D^b_coh, Z}, %
restriction induces equivalences
$$D^b_{\coh}(\cX_n)/D^b_{\coh,Z}(\cX_n) \iso D^b_{\coh}(\cU_n).$$
Passing to the colimit over~$n$, and recalling the statement of Proposition~\ref{aprop:Coh-set-theoretically-supported},
we obtain an equivalence
\begin{equation}\label{eqn: computing quotient category in coherent recollement compact version}
D^b_{\coh}(\cX)/D^b_{\coh}(\widehat{\cX}) \iso D^b_{\coh}(\cU).
\end{equation} 

Passing to Ind-categories, and taking into account Proposition~\ref{prop:Ind is exact},
we obtain an equivalence
\begin{equation}\label{eqn: computing quotient category in coherent recollement}
\Ind D^b_{\coh}(\cX)/\Ind D^b_{\coh}(\widehat{\cX}) \iso \Ind D^b_{\coh}(\cU),
\end{equation} 
which we use to identify its left and right hand sides from now on.
Taking into account this identification,
as well as the prescription of 
Hypothesis~\ref{hyp:usual hyp for semiorthogonal},
we then write $j^*$ to denote the composite
\begin{equation}\label{eqn:j-upper-star}
j^*: D^b_{\coh}(\cX) \to D^b_{\coh}(\cX)/D^b_{\coh}(\widehat{\cX}) \iso D^b_{\coh}(\cU),
\end{equation}
as well as its Ind-extension 
\begin{equation}
\label{eqn:Ind-extended j lower star}
j^{*}:\Ind D^b_{\coh}(\cX) \to \Ind D^b_{\coh}(\cX)/\Ind D^b_{\coh}(\widehat{\cX})
\iso \Ind D^b_{\coh}(\cU).
\end{equation}
By Lemma~\ref{lem: abstract semiorthogonal decomposition},
the functor~\eqref{eqn:Ind-extended j lower star}  
admits a right adjoint
$j_*:\Ind D^b_{\coh}(\cU) \to \Ind D^b_{\coh}(\cX).$

\begin{remark}
Although we don't need it in the present paper,
we remark that the Pro-extension of $j^{*}$ also admits %
a left adjoint $j_!$,  which %
is the Pro version of ``extension by zero from~$\cU$''.  As far as
we understand, defining $j_!$ in this way was the original motivation
for introducing pro-coherent sheaves (by Deligne in~\cite[Appendix]{MR222093}). %
\end{remark}

\begin{rem}\label{notational remark about j star}
The functors $j^*$ in~\eqref{eqn:j-upper-star} and~\eqref{eqn:Ind-extended j lower star} %
are induced by the functors $j_n^*$ arising from
the open immersions $j_n:\cU_n \hookrightarrow \cX_n$, 
in the sense that if we write each of $D^b_{\coh}(\cX)$ and $D^b_{\coh}(\cU)$
as colimits over~$n$, following the prescription of~\eqref{eqn:defn-of-coherent-sheaves},
then $j^*$ is obtained as colimit of the functors~$j_n^*$. 
Thus the notation~$j^*$ is compatible with our notation for operations on derived categories
of coherent sheaves on algebraic stacks.

Each of the functors $j_n^*$ has a corresponding right adjoint 
as in~\eqref{eqn:Ind derived pushforward},
which we denote by
$j_{n,*}:D^b_{\coh}(\cU_n) \to \Ind D^b_{\coh}(\cX_n).$
Since the functor $j_*$ is compatible with the formation of (filtered, and hence all)
colimits (since it is right adjoint to the compact object-preserving functor~$j^*$),
one finds that $j_*$ is similarly the colimit of the~$j_{n,*}$.

Although it might be better to write $Rj_*$ rather than~$j_*$, 
to indicate the derived nature of this functor,
in the discussion of Subsection~\ref{subsubsec:closed pushforward}
we have reserved the former notation to denote the derived functor of~$j_*$
in the context of derived categories of complexes of quasi-coherent sheaves,
rather than in the context of Ind coherent complexes.  
Thus we content ourselves with writing simply~$j_*$ in the Ind coherent
case, trusting that context will resolve any possible ambiguity.
This also allows us to maintain
notational consistency with the general discussion of Appendix~\ref{app:category-theory}
related to~Hypothesis~\ref{hyp:usual hyp for semiorthogonal}.
In any event,
in our eventual applications,
the relevant open immersions $j$ will be cohomologically affine,
and so the functors $j_*$ that we have to consider will in fact be $t$-exact;
thus no confusion should arise from our choice of notation.
\end{rem}

\begin{aremark}
\label{rem:our definition of coherent}
Our definition of the categories
$\Coh(\cX)$ and
$D^b_{\coh}(\cX)$ is  less expansive than some others in the literature,
since, by definition, any object of either category %
is supported on a
closed algebraic substack of~$\cX$.  Thus, for example, the structure sheaf
$\cO_{\cX}$ will not be an object of $\Coh(\cX)$ or $D^b_{\coh}(\cX)$, unless $\cX$ happens to
be an algebraic stack, rather than merely formally algebraic.  Rather,
in our setting, the structure sheaf will be realized as an object
of $\Pro \Coh(\cX)$, or $\Pro D^b_{\coh}(\cX)$, namely the pro-object 
$\quoteslim{n} \cO_{\cX_n}.$.
In other papers, the $\infty$-category that we denote by $D^b_{\coh}(\cX)$
might instead be denoted $D^b_{\coh,\cX_{\red}}(\cX)$: the category of coherent
complexes whose cohomologies are set-theoretically supported on the closed
substack $\cX_{\red}$ of $\cX$. %

We use the definition that we do because it fits well with the general definition
of Ind-coherent sheaves on Ind-schemes, and also because it ensures that there are no hidden
issues related to completion or adic topologies lurking in the background;
objects of $D^b_{\coh}(\cX)$ and $\Ind D^b_{\coh}(\cX)$ are, by definition,
adically discrete.
Issues of adic topologies and completions are then handled by the use of Pro categories.

We also note that we are, in fact,  unaware of a general
treatment of stable $\infty$-categories of coherent sheaves on    
formal algebraic stacks in the more expansive sense that includes the structure 
sheaf~$\cO_{\cX}$ as one its objects, although presumably it can be realized
as a special
case of Clausen--Scholze's theory of stable $\infty$-categories of solid modules on analytic stacks.
The one exception is the case of formal algebraic stacks obtained as
completions of algebraic stacks; in this case, \cite{BCGAGA} introduces the
corresponding triangulated category.
We refer to
Proposition~\ref{prop:completion comparison}
below for a discussion of how to relate this more standard approach 
to our approach via Pro categories.
\end{aremark}

\subsection{Pullback of pro-coherent sheaves} 
\label{subsubsec:restriction and support}

In this section  we consider various pullback and pushforward functors defined on abelian categories of coherent and pro-coherent sheaves. %

\subsubsection{Representable morphisms}\label{subsubsec: representable morphisms}%
Assume that~$\cZ, \cX$ are Noetherian formal algebraic stacks with affine diagonal, and assume that $f: \cZ \to \cX$ is representable by algebraic stacks.
Choose a presentation $\cX = \colim_{n} \cX_n$ as in~\eqref{eqn:present-formal-alg-as-colimit}, and let $\cZ_n \coloneqq  \cZ \times_{\cX} \cX_n$, which is a closed algebraic substack of~$\cZ$.
By~\cite[Lemma~4.8]{Emertonformalstacks}, the natural map $\colim_{n} \cZ_n \to \cZ$ is an isomorphism, hence it is a presentation of~$\cZ$ as in~\eqref{eqn:present-formal-alg-as-colimit}.
Write~$x_{mn} : \cX_m \to \cX_n$ for the transition map (which is a closed immersion), and similarly for $z_{mn} : \cZ_m \to \cZ_n$.
The map~$f$ is then the colimit over~$n$ of morphisms $f_n : \cZ_n \to \cX_n$ of algebraic stacks. 
Since the transition maps~$x_{mn}, z_{mn}$ are closed immersions, the diagrams
\[
\begin{tikzcd}
\Coh(\cX_n) \arrow[r, "f_n^*"] & \Coh(\cZ_n)\\
\Coh(\cX_m) \arrow[r, "f_m^*"] \arrow[u, "x_{mn, *}"] & \Coh(\cZ_m) \arrow[u, "z_{mn, *}"]
\end{tikzcd}
\]
commute.
Since $\Coh(\cX) = \colim_n \Coh(\cX_n)$ and
$\Coh(\cZ) = \colim_n \Coh(\cZ_n)$,
we can thus define
\[
f^* = \colim_n f_n^*: \Coh(\cX) \to \Coh(\cZ).
\]
If~$f$ is furthermore proper (in the sense of~\cite[Definition~3.11, Remark~3.13]{Emertonformalstacks}, i.e.\ its pullback to any algebraic stack is proper)
the pushforwards~$f_{n*}$ preserve coherence, as recalled in Section~\ref{subsubsec:operations}. %
We can thus define
\[
f_* = \colim_n f_{n*}: \Coh(\cZ) \to \Coh(\cX).
\]
Up to natural equivalence, the functors~$f^*, f_*$ are independent of the choice of presentation of~$\cX$.
They form an adjoint pair.

\begin{rem}\label{rem:Pro-extensions-denoted-same-symbols}
We will denote the Pro-extensions of~$f^*$ and~$f_*$ by the same symbols, to avoid confusion with the completed pullbacks to be defined in Section~\ref{subsubsec:completed pullbacks}.
\end{rem}

\begin{rem}\label{rem: closed immersions}
  If $i:\cZ \to \cX$ is a closed immersion, then it is representable by algebraic stacks and proper, and by the discussion above we have an adjunction
  \[
  (i^*, i_*): \Pro\Coh(\cX) \to \Pro\Coh(\cZ).
  \]
  The right adjoint $i_*$ is exact and fully faithful, and the unit $1 \to i_*i^*$ is surjective.

  If~$\cZ$ is furthermore algebraic, and we let~$\widehat \cX$ denote the completion of~$\cX$ along~$\cZ$,
  then $\Coh(\cZ)$ is a full subcategory of~$\Coh(\widehat{\cX})$,
  and $i_*$ is the restriction to~$\Coh(\cZ)$
  of the functor $\Coh(\widehat{\cX}) \to \Coh(\cX)$ 
  obtained by restricting 
  the $t$-exact functor~\eqref{eqn:formal pushforward}
  to the heart of its domain.
  This can be checked by computing~$i_*$ with respect to a presentation $\cX = \colim_{n} \cX_n$ such that $\cZ$ is a closed substack of~$\cX_n$ for all~$n$.
\end{rem}

\subsubsection{Presentation of pro-coherent sheaves}
The functors~$i^*$ and~$i_*$ %
can be used to give presentations for pro-coherent sheaves, as in the following lemma.

\begin{lem}
  \label{lem:pro-coherent-as-limit-of-pullbacks}
Let~$\cX$ be a Noetherian formal algebraic stack with affine diagonal,
and choose a presentation
$\cX \iso \colim_n \cX_n$ as in~\eqref{eqn:present-formal-alg-as-colimit},
Write $i_n: \cX_n \hookrightarrow \cX$ to denote the canonical closed
immersion.  
Then, for any pro-coherent sheaf~$\cF$ on~$\cX$, we have $\cF\isoto\limcommand_n i_{n,*} i_n^{*}\cF$ in $\Pro \Coh(\cX)$. 
\end{lem}
\begin{proof}
We first consider the case when $\cF$ is actually coherent.  Then by definition
$\cF = (i_m)_*\cF_m$ for some $m$ and some object $\cF_m$ of~$\Coh(\cX_m)$.
Thus
$$\limcommand_n\, i_{n,*} i_n^*\cF = \limcommand_{n \geq m} i_{n,*}i_n^* (i_m)_* \cF_m
\iso \limcommand_{n\geq m} (i_m)_* \cF_m \iso (i_m)_* \cF_m = \cF.
$$
In the general case that $\cF$ is pro-coherent, 
write $\cF \iso \limcommand_{i\in I}\cF_i$ as a cofiltered limit of coherent sheaves~$\cF_i$.
  Then \[\limcommand_n i_{n,*} i_n^{*}\cF
\iso
\limcommand_n\limcommand_{i\in I}i_{n,*} i_n^{*}\cF_i
\iso
\limcommand_{i\in I}\limcommand_n i_{n,*} i_n^{*}\cF_i=\limcommand_{i\in I}\cF_i, \]
as required.
In more detail, the first isomorphism holds because $i_{n,*}$ and $i_n^*$ are pro-extended,
and so compatible with the formation of cofiltered limits; 
the second isomorphism is an interchange of cofiltered limits;
and the third isomorphism follows from the case of coherent sheaves
that we've already proved.
\end{proof}

\subsubsection{Completed pullbacks}\label{subsubsec:completed pullbacks}
We continue to assume that~$f:\cZ\to \cX$ is a morphism of Noetherian formal algebraic stacks with affine diagonal, but we no longer assume that~$f$ is representable by algebraic stacks. %

Suppose firstly that~$\cX$ is algebraic, and choose a presentation $\cZ = \colim_{n} \cZ_n$ as in~\eqref{eqn:present-formal-alg-as-colimit}.
Let~$f_n : \cZ_n \to \cX$ be the restriction of~$f$ to~$\cZ_n$.
We can then define a functor
\begin{align*}
  \widehat f^* : \Coh(\cX) &\to \Pro \Coh(\cZ)\\ \cF &\mapsto \quoteslim{n}f^*_n \cF.
\end{align*} Up to natural equivalence, this functor is independent of the choice of presentation of~$\cZ$.

\begin{lem}
  \label{lem:cartesian-square-completed-pullback}
Let~$f: \cZ\to \cX$ be a morphism of Noetherian formal algebraic stacks with affine diagonal. 
Suppose that~$\cX$ is in fact an algebraic stack, let $i:\cY\to\cX$ be a closed immersion, and let $f': \cZ\times_{\cX} \cY \to \cY$
and $i':\cZ\times_{\cX} \cY \to \cZ$
denote the base-changes of $f$ and $i$ respectively.
Then the natural diagram
{\em (}where the right hand vertical arrow 
is defined following Remark~{\em \ref{rem:Pro-extensions-denoted-same-symbols})}
\begin{equation*}
\begin{tikzcd}
\Coh(\cX) \arrow[r, "\widehat f^*"] & \Pro \Coh(\cZ) \\
\Coh(\cY) \arrow[r, "\widehat f'^*"] \arrow[u, "i_{*}"] & \Pro \Coh(\cZ\times_{\cX}\cY) \arrow[u, "i'_{*}"]
\end{tikzcd}
\end{equation*}
commutes.
\end{lem}
\begin{proof}By the definition of~$\widehat f^*$, we are reduced to noting that  the pushforward functor for a closed immersion of algebraic stacks commutes with arbitrary base change.  
\end{proof}
In the more general case that~$\cX$ is not necessarily algebraic, we as usual choose a presentation~$\cX = \colim_n \cX_n$.
Then we can form the fibre product $\cZ_n \coloneqq  \cZ \times_\cX \cX_n$, a closed formal algebraic substack of~$\cZ$ with a map $f_n : \cZ_n \to \cX_n$.
By Lemma~\ref{lem:cartesian-square-completed-pullback}, the diagrams
\begin{equation*}%
\begin{tikzcd}
\Coh(\cX_n) \arrow[r, "\widehat f_n^*"] & \Pro \Coh(\cZ_n) \\
\Coh(\cX_m) \arrow[r, "\widehat f_m^*"] \arrow[u, "x_{mn, *}"] & \Pro \Coh(\cZ_m) \arrow[u, "z_{mn, *}"]
\end{tikzcd}
\end{equation*}
commute, where the vertical arrows are induced by the transition maps $x_{mn} : \cX_m \to \cX_n$.
We can thus pass to the colimit and obtain a functor %
\begin{equation}\label{eqn:definition of completed pullback functor}
\widehat f^* : \Coh(\cX) \to \Pro \Coh(\cZ),
\end{equation}
which is independent (up to natural isomorphism) of the choice of presentation of~$\cX$.
We continue to denote by $\widehat f^{*}$ the $\Pro$-extension \[
\widehat f^* : \Pro\Coh(\cX) \to \Pro \Coh(\cZ).
\]
We can also consider the composite functors 
\[
\Pro \Coh(\cX) \xrightarrow{x_n^*}
\Pro \Coh(\cX_n) \xrightarrow{\widehat f_n^*} \Pro \Coh(\cZ_n) \xrightarrow{z_{n*}} \Pro\Coh(\cZ).
\]
where $x_n: \cX_n \to \cX$ and $z_n: \cZ_n \to \cZ$ denote the closed immersions. 
(Noting that the target of~$f_n: \cZ_n \to \cX_n$ is algebraic,
the functor $\widehat f_n^*$ is defined by the discussion above,
while $x_n^*$ and ~$z_{n*}$ are defined in Remark~\ref{rem:Pro-extensions-denoted-same-symbols}.)  %
\begin{lem}%
  \label{lem:computing-Pro-completed-pullback}With notation as above, the functor $\widehat f^* : \Pro\Coh(\cX) \to \Pro \Coh(\cZ)$ is naturally isomorphic %
to~$\lim_n z_{n,*}\widehat f_n^{*}x_n^*$.
\end{lem}
\begin{proof}
Lemma~\ref{lem:pro-coherent-as-limit-of-pullbacks}
gives a natural isomorphism of functors 
$$\id \iso \lim_{n} x_{n,*} x_n^*$$
on~$\Pro \Coh(\cX)$. 
Applying $\widehat f^{*}$, we find that
$$\widehat{f}^{*} \iso \lim_{n} \widehat f^* x_{n,*} x_n^*
\iso \lim_n z_{n,*} \widehat f^*_n x_n^*$$
(the first isomorphism following from the fact that $\widehat f$ commutes with cofiltered limits
by its construction as a Pro-extension, and the second claim
following from the definition of~\eqref{eqn:definition of completed pullback functor}), %
as claimed.
\end{proof}

\begin{rem}\label{rem:map of adjunctions}
Let $f : \cZ \to \cX$ be a morphism of Noetherian formal algebraic stacks with affine diagonal, 
and let $i : \cY \to \cX$ be a closed immersion.
As in Lemma~\ref{lem:cartesian-square-completed-pullback},
we can form the Cartesian diagram
\begin{equation}\label{eqn: Cartesian diagram closed immersions}
\begin{tikzcd}
\cZ \arrow[r, "f"] &
\cX\\
\cZ \times_\cX \cY \arrow[r, "f'"] \arrow[u, "i'"]&
\cY \arrow[u, "i"].
\end{tikzcd}
\end{equation}
Analyzing the definitions, we see that we have commutative diagrams
\begin{equation}\label{to prove commutative for compatibility of completion and pushforward}
\begin{array}{cc}
\begin{tikzcd}
  \Pro\Coh(\cX) & \Pro\Coh(\cZ)\\
  \Pro \Coh(\cY) & \Pro\Coh(\cZ \times_\cX \cY)
  \arrow["\widehat f^*", from=1-1, to=1-2]
  \arrow["\widehat f'^*", from=2-1, to=2-2]
  \arrow["i_*", from=2-1, to=1-1]
  \arrow["i'_*", from=2-2, to = 1-2]
\end{tikzcd}
&
\begin{tikzcd}
\Pro\Coh(\cX) &  \Pro\Coh(\cZ) \\
\Pro\Coh(\cY) & \Pro \Coh(\cZ \times_\cX \cY)
  \arrow["\widehat f^*", from=1-1, to=1-2]
  \arrow["\widehat f'^*", from=2-1, to=2-2]
  \arrow["i'^*", from=1-2, to=2-2]
  \arrow["i^*", from=1-1, to = 2-1]
\end{tikzcd}
\end{array}
\end{equation}
and furthermore, 
for all~$\cF \in \Pro\Coh(\cX)$, the image under~$\widehat f^*$ of the unit $\cF \to i_{*}i^*\cF$ 
is naturally isomorphic to the unit of $(i'^*, i'_*)$ evaluated at~$\widehat f^* \cF$.
Note that commutativity of the leftmost square in~\eqref{to prove commutative for compatibility of completion and pushforward} is a
generalization of Lemma~\ref{lem:cartesian-square-completed-pullback} to the case where~$\cX$ is not necessarily algebraic.

We will often apply this construction in a context where~$\cY$ is furthermore algebraic,
and $f: \widehat \cX_Z \to \cX$ is the completion of~$\cX$ at a closed subset $Z \subset |\cX|$.
Then~\eqref{eqn: Cartesian diagram closed immersions} becomes a Cartesian diagram
\[
\begin{tikzcd}
\widehat \cX_Z \arrow[r, "f"] &
\cX\\
\widehat \cY_{|\cY| \cap Z} \arrow[r, "f'"] \arrow[u, "i'"]&
\cY \arrow[u, "i"]
\end{tikzcd}
\]
since the pullback $\widehat \cX_Z \times_\cX \cY$ is the completion of~$\cY$ at $|\cY| \cap Z$.
We will sometimes write $i_Z$ for the closed immersion~$i'$, so that we have an adjunction
\begin{equation}\label{eqn: adjoint pair for closed immersion II}
(i_Z^*, i_{Z, *}) : \Pro \Coh(\widehat \cX_Z) \to \Pro\Coh(\widehat \cY_{|\cY| \cap Z}),
\end{equation}
as in Remark~\ref{rem: closed immersions}, and the image under $\widehat f^*$ of the unit of $(i^*, i_*)$ is the unit of $(i_Z^*, i_{Z, *})$.
\end{rem}

\begin{rem}\label{completed pullback is functorial}
If~$f : \cX \to \cY$ and~$g : \cY \to \cZ$ are morphisms of Noetherian formal algebraic stacks with affine diagonal, 
then $\widehat{gf}^*$ is naturally isomorphic to~$\widehat f^* \widehat g^*$.
To see this, we can assume without loss of generality that $\cZ$ is algebraic.
Choose a presentation $\cY = \colim_n \cY_n$, and let $\cX_n \coloneqq  \cX \times_\cY \cY_n$.
Then~$\cX_n$ is a closed formal algebraic substack of~$\cX$, and $\cX = \colim_n \cX_n$ (note that any map $\Spec A\to \cX$ factors through some~$\cX_n$, because the composite $\Spec A\to\cX\to\cY$ factors through some~$\cY_n$).%

Choose also a presentation $\cX = \colim_m \cX'_m$, and define $\cX_{n,m} = \cX_n \times_\cX \cX_m'$.
Then $\cX_n = \colim_m \cX_{n, m}$ for all~$n \geq 0$, by~\cite[Lemma~4.8]{Emertonformalstacks},
so that $\cX = \colim_{\bZ_{\geq 0}^2} \cX_{n, m}$, where the colimit is taken with respect to the product order on~$\bZ_{\geq 0}^2$.

Now, if~$\cF \in \Coh(\cZ)$, then by construction we have
\[
\widehat f^* \widehat g^*(\cF) = \widehat f^*(\lim_n g_n^*\cF) = \lim_{n, m} f^*_{n, m}g^*_n \cF
\]
where~$g_n : \cY_n \to \cZ, f_{n, m}: \cX_{n, m} \to \cY_n$ are the natural maps.
This also coincides with $\widehat{gf}^* \cF$ computed
with respect to the presentation $\cX = \colim_{\bZ_{\geq 0}^2} \cX_{n, m}$, as desired.
\end{rem}

\begin{rem}\label{rem:pullback and completed pullback}
Let~$\cZ, \cX$ be Noetherian formal algebraic stacks with affine diagonal, and assume that~$f : \cZ \to \cX$ is representable by algebraic stacks (which is always the case if~$\cZ$ is
algebraic, by our assumption on the diagonal of~$\cX$).
It then follows from the definitions that~$\widehat f^* : \Coh(\cX) \to \Pro\Coh(\cZ)$ takes values in~$\Coh(\cZ)$, and is naturally isomorphic to~$f^*$.
It follows that their Pro-extensions are also naturally isomorphic, as functors $f^*, \widehat f^* : \Pro \Coh(\cX) \to \Pro \Coh(\cZ)$.
\end{rem}

\subsubsection{Completed pullbacks to completed stacks}\label{subsubsec:completed-pullback-to-completed}

As in Section~\ref{subssubsec:completion-along-substack-formal-algebraic},
we let $\cX$ be a Noetherian formal algebraic stack with affine diagonal,
and let $\cZ \hookrightarrow \cX$ be a closed substack of~$\cX$.
We write $i:\widehat{\cX}\into\cX$ for the completion of $\cX$ along~$\cZ$.
Following as above the notational conventions of Appendix~\ref{subsec:semiorthogonal and Ind Pro}, we write 
\begin{equation}
\label{eq:ihatlower*}
\ihat_{*}:\Pro\Coh(\cXhat)\to\ProCoh(\cX),
\end{equation}
\begin{equation}
\label{eq:ihatupper*}
\widehat i^{*}:\Pro\Coh(\cX)\to\ProCoh(\cXhat)
\end{equation}for the $\Pro$-extensions of $i_{*}:\Coh(\cXhat)\to\Coh(\cX)$ and $\widehat i^{*}:\Coh(\cX)\to\ProCoh(\cXhat)$.

\begin{lem}
  \label{lem:completed-Pro-pullback-completion-formula}With notation as above, write $\cXhat=\colim_n\cZ_n$ as a colimit of thickenings of closed algebraic substacks of~$\cX$, 
  and let $k_n:\cZ_n\into\cXhat$
and $i_n: \cZ_n \into \cX$ denote the corresponding closed immersions.
Then the functor $\widehat i^{*}:\Pro\Coh(\cX)\to\ProCoh(\cXhat)$ agrees with the functor~$\lim_n 
k_{n,*} i_n^{*}$.
\end{lem}
\begin{proof}
Since~$\cZ_n$ is algebraic, the functor $\widehat k_n^*$, resp.\ $\ihat^*_n$, coincides with~$k_n^*$, resp.\ $i_n^*$, by Remark~\ref{rem:pullback and completed pullback}.
  Lemma~\ref{lem:pro-coherent-as-limit-of-pullbacks}
(applied to~$\cXhat$) shows that
$$\widehat i^* \iso  \lim_n  k_{n,*} k_n^* \widehat i^* 
\iso \lim_n  k_{n,*} i_n^*,$$
where we have applied the natural isomorphisms $k_n^* \widehat i^* = \widehat k_n^* \widehat i^* \iso \ihat^*_n = i_n^*$
of Remark~\ref{completed pullback is functorial}, taking into account that
$i_n = i k_n$ for each $n$.
\end{proof}

\begin{rem}
If~$\cX$ is furthermore algebraic, the preceding lemma recovers the definition of the completed pullback~$\ihat^*$.
On the other hand, if we take~$\cX = \widehat \cX$, $i_n = k_n$, and $\ihat^*$ to be the identity of $\Pro \Coh(\widehat \cX)$, then we recover
Lemma~\ref{lem:pro-coherent-as-limit-of-pullbacks}.
\end{rem}

\begin{prop}
  \label{prop:exact-adjoint-ihat-upper-star}Write $i:\widehat{\cX}\into\cX$ for the inclusion of the completion~$\widehat{\cX}$ of a Noetherian formal algebraic stack~$\cX$, 
  with affine diagonal, along a closed algebraic substack~$\cZ$. 
  Then the functor $\widehat i^{*}:\Pro\Coh(\cX)\to\ProCoh(\cXhat)$ is exact, and is left adjoint to $\ihat_{*}:\Pro\Coh(\cXhat)\to\ProCoh(\cX)$.
\end{prop}
\begin{proof}
  Since the completion $\widehat \cX$ only depends on the underlying topological space of~$\cZ$, 
  if~$\cX = \colim_n \cX_n$ is a presentation as in~\eqref{eqn:present-formal-alg-as-colimit}, then we have $\widehat \cX = \colim_n \cXhat_n$. %
  It is therefore easy (and formal) to reduce to the case that~$\cX$ is algebraic, which we assume from now on.

  We begin by showing exactness, for which it suffices to show that the restriction $\widehat i^{*}:\Coh(\cX)\to\ProCoh(\cXhat)$ is exact.
Write $\cI$ for the coherent ideal sheaf of~$\cZ$ in~$\cX$, so that for any coherent sheaf~$\cF$ we have
\begin{equation}
  \label{eq:completion-AR-version}
  \widehat i^{*}\cF=\quoteslim{n} \cF/\mathcal I^n.
\end{equation}
 The exactness of~$\ihat^{*}$ then follows from a standard Artin--Rees argument, \emph{cf.}\ \cite[\href{https://stacks.math.columbia.edu/tag/0G9M}{Tag 0G9M}]{stacks-project}.

We now turn to the adjunction.
Since both~$\ihat^{*}$ and~ $\ihat_{*}$ are $\Pro$-extended, it suffices to show that if~ $\cF$ is an object of $\Coh(\cX)$ and $\cG$ is an object of $\Coh(\cXhat)$ then \[\Hom_{\Pro\Coh(\cXhat)}(\ihat^{*}\cF,\cG)=\Hom_{\Pro\Coh(\cX)}(\cF,\ihat_{*}\cG),\] or equivalently that \[\Hom_{\Pro\Coh(\cXhat)}(\ihat^{*}\cF,\cG)=\Hom_{\Coh(\cX)}(\cF,i_{*}\cG).\]   %
As usual, we write $\widehat{\cX}=\colim_n\cX_n$, so that $i_n:\cX_n\into\cX$ is a closed immersion of an algebraic substack.
Relabeling if necessary, we can assume that~$\cG$ is supported on~$\cX_1 $, so that we can (slightly abusively) write~$i_{*}\cG=i_{n,*}\cG$ for each~$n$.
Then we compute that
\begin{multline*} \Hom_{\Pro\Coh(\cXhat)}(\ihat^{*}\cF,\cG)= \\
\Hom_{\Pro\Coh(\widehat{\cX})}(\quoteslim{n} i_n^{*}\cF, \cG) = 
\colim_{n}
\Hom_{\Coh(\widehat{\cX})}(i_n^{*}\cF,\cG) = \\
\colim_{n}
\Hom_{\Coh(\cX_n)}(i_n^{*}\cF,\cG)  = 
\colim_{n}
\Hom_{\Coh(\cX)}(\cF,i_{n,*}\cG)  =\\ \Hom_{\Coh(\cX)}(\cF,i_{*}\cG),
\end{multline*}
as required. (In more detail, the first isomorphism is the definition of~$\ihat^{*}$, the second 
holds by virtue of how morphisms are computed in the Pro-category, the third is by the definition of~$\Coh(\cXhat)$ and our assumption that~$\cG$ is supported on~$\cX_1 $, the fourth is the adjunction between $i^{*}_n$ and $i_{n,*}$, and the fifth is again by our assumption that~$\cG$ is supported on~$\cX_1 $.)
\end{proof}

\begin{rem}\label{rem:ihat surjects onto i upper star}
We put ourselves in the context of Lemma~\ref{lem:completed-Pro-pullback-completion-formula}, 
so that $i_1: \cZ \hookrightarrow \cX$ is a closed immersion of an algebraic stack~$\cZ$ into
a Noetherian formal algebraic stack~$\cX$ with affine diagonal,
$i: \widehat \cX \to \cX$ is the completion of~$\cX$ along~$\cZ$, and
$k_1: \cZ \to \widehat \cX$ is the defining closed immersion.
Then for any object $\cF$ of $\Pro \Coh(\cX)$,
there is a canonical surjection
\begin{equation}\label{eqn:surjection of ihat onto i^*}
\widehat{i}^* \cF \to k_{1, *}i_1^*\cF
\end{equation}
in~$\Pro \Coh(\widehat \cX)$.  
In fact, as noted in Remark~\ref{rem: closed immersions}, the unit of adjunction 
$\cF \to i_{1, *} i_1^* \cF$ is surjective, and so it remains surjective after
applying the functor~$ \widehat{i}^*$, 
which is exact by Proposition~\ref{prop:exact-adjoint-ihat-upper-star}.
It now suffices to note that $\widehat{i}^* i_{1, *}i_1^{*}\cF=k_{1, *}i_1^{*}\cF$, %
because
the counit of adjunction $\widehat{i}^* \widehat{i}_*
\to \id$ is an isomorphism, and $i_{1, *}  = \ihat_* k_{1, *}$. %

Applying the exact functor~$\ihat_*$ to~\eqref{eqn:surjection of ihat onto i^*} we obtain
a surjection %
$$\ihat_* \widehat{i}^* \cF \to i_{1, *}i_1^*\cF$$
of objects in~$\Pro\Coh(\cX).$
\end{rem}

\begin{aremark}\label{arem:checking coherence on smooth cover}
If~$\cX$ is a Noetherian algebraic stack with affine diagonal,
and $p: \Spec A \to \cX$ is a flat and surjective morphism,
then the pullback $p^*: \Pro \Coh(\cX) \to \Pro \Coh(\Spec A)$ reflects coherence, 
in the sense that if $\cF = \quoteslim{i \in I} \cF_i \in \Pro \Coh(\cX)$, and $p^* \cF$ is an object of $\Coh(\Spec A)$, then $\cF$ is an object of~$\Coh(\cX)$.
To see this, we first remark that $p^*$ is exact and faithful (and thus conservative) on $\Coh(\cX)$, 
and so the same is true of its Pro-extension, by e.g.\ \cite[Prop.\ 6.1.10, Cor.\ 8.6.8]{MR2182076}.
Now, if $\cG \coloneqq p^* \cF$ is coherent, then $\cG$ has descent data to~$\cX$, since it's pulled back from~$\cX$ when thought of as a pro-coherent sheaf.
It thus descends to a coherent sheaf $\cG'$ on~$\cX$.
The isomorphism $p^* \cG' \isoto p^* \cF$ induces morphisms $p^* \cG' \to p^* \cF_i$ for all~$i$, which are compatible with the transition maps in~$\cF_i$,
and with the descent data on source and target.
These morphisms therefore descend to a morphism $\cG' \to \cF$ in $\Pro \Coh(\cX)$, which becomes an isomorphism after applying~$p^*$.
Since, as we already recalled, $p^*$ is conservative, we conclude that $\cG' \isoto \cF$ and so~$\cF$ is coherent, as desired.

More generally, if $\cX = \colim_n \cX_n$ is a Noetherian formal algebraic stack with affine diagonal, presented as in~\eqref{eqn:present-formal-alg-as-colimit},
and $\Spf A \to \cX$ is a smooth surjective morphism which is representable by algebraic stacks,
then the pullback $p^*: \Pro \Coh(\cX) \to \Pro \Coh(\Spf A)$ reflects coherence.
To see this, let $\cF \in \Pro \Coh(\cX)$, and assume that $p^*\cF$ is coherent.
Then, if we define $\Spec A_n \coloneqq  \Spf A \times_\cX \cX_n$, and write $k_n : \cX_n \to \cX$ for the closed immersion and for its pullback to~$\Spf A$,
there exists~$n$ such that $p^* \cF \to k_{n,*}k_n^*p^* \cF$ is an isomorphism. 
Hence $\cF \to k_{n,*}k_n^* \cF$ is also an isomorphism.
Since~$\cX_n$ is an algebraic stack, we have thus reduced to the result of the previous paragraph.
\end{aremark}

\subsubsection{Passage to derived categories}
\label{subsubsec:derived passage}

By Proposition~\ref{prop:exact-adjoint-ihat-upper-star} and Corollary~\ref{cor:pro t-exact derived
functors}, the functors~\eqref{eq:ihatlower*} and \eqref{eq:ihatupper*}
induce $t$-exact adjoint functors\begin{equation}\label{ihat lower star}
\widehat{i}_*: \Pro D^b_{\coh}(\widehat{\cX}) \hookrightarrow \Pro D^b_{\coh}(\cX)
\end{equation}
and
\begin{equation}\label{ihat upper star}
\widehat{i}^*: \Pro D^b_{\coh}(\cX) \to \Pro D^b_{\coh}(\widehat{\cX}).
\end{equation}

\begin{remark} %
\label{rem:i-hat-upper-star as quasi-inverse}
Since $\widehat{i}_*$ is fully faithful (as a consequence of Proposition~\ref{aprop:Coh-set-theoretically-supported}), the counit of adjunction
$\widehat{i}^*\widehat{i}_* \to \id$ is a natural isomorphism.
Thus the restriction of $\widehat{i}^*$ to $D^b_{\coh,Z}(\cX)$
induces an inverse to the equivalence $D^b_{\coh}(\widehat{\cX}) \iso D^b_{\coh,Z}(\cX)$
of Proposition~\ref{aprop:Coh-set-theoretically-supported}.
\end{remark}

\subsubsection{Scheme-theoretic support.}\label{subsubsec:support}
It is unclear to us in what generality it is sensible to define the notion of ``scheme-theoretic'' support for pro-coherent sheaves, but the following definitions will suffice for us. %

\begin{defn}
  \label{defn:scheme-theoretic-support-pro-coherent}
  Let~$z:\cZ\into\cX$ be a closed immersion of Noetherian formal algebraic stacks with affine diagonal. %
  \begin{enumerate}
  \item Let~$\cF$ be a pro-coherent sheaf on~$\cX$.
Then we say that~$\cF$ is \emph{scheme-theoretically supported on~$\cZ$} if~$\cF$ is in the essential image of $z_{*}:\Pro\Coh(\cZ)\to \Pro\Coh(\cX)$, or 
equivalently, if the unit of adjunction $\cF \to z_*z^*\cF$ is an isomorphism.
  \item Let~$\cF$ be a coherent sheaf on~$\cX$, and assume that $\cZ$ is algebraic.
  Then we say that~$\cZ$ is the \emph{scheme-theoretic support} of~$\cF$ if 
  $\cF$ is scheme-theoretically supported on~$\cZ$ and 
  the natural map $\cO_{\cZ} \to \underline{\End}_{\Coh(\cZ)}(z^*\cF)$ is injective, 
  where~$\underline{\End}$ denotes the sheaf-End. This being the case, we say that $|\cZ|$ is the  \emph{set-theoretic support} of~$\cF$. %
  \end{enumerate}
\end{defn}

\begin{rem}\label{rem: scheme-theoretic supports on schemes}
If~$\cX$ is a scheme, %
and $\cF$ is a coherent sheaf on~$\cX$,
then the scheme-theoretic support of~$\cF$ is the closed subscheme associated to the kernel of $\cO_\cX \to \underline{\End}_{\Coh(\cX)}(\cF)$,
which is a quasicoherent ideal sheaf.
Hence Definition~\ref{defn:scheme-theoretic-support-pro-coherent}~(2) agrees with the usual definition of the scheme-theoretic support of $\cF$ in this case. %
\end{rem}

\begin{defn}
  \label{defn:pure-degree-zero}We say that an object~$\cF$ of  $\Pro D^b_{\coh}(\cX)$   is \emph{pure of degree zero} if it is a pro-coherent sheaf: that is, it is contained in the heart $\ProCoh(\cX)$ of~$\Pro D^b_{\coh}(\cX)$.
  (Note that this is \emph{a priori} a stronger condition than asking that the cohomology sheaves of~$\cF$ vanish outside of degree zero.)
\end{defn}

\begin{lem}
  \label{lem:pro-coh-vanishes-if-pullbacks-do}
Let~$\cX$ be a Noetherian formal algebraic stack with affine diagonal.
  \begin{enumerate}
  \item Suppose that $\cF$ is an object of $\ProCoh(\cX)$. Then~$\cF=0$ if and
    only if, for every closed point $x \in |\cX|$ with completion-at-$\{x\}$ morphism $i_x : \widehat \cX_x \to \cX$, we have $\ihat_x^{*}\cF=0$.
  \item Suppose that $\cG$ is an object of $\Pro D^b_{\coh}(\cX)$ which is bounded above.
Then $\cG$ is
    pure of degree zero, resp.\ equal to zero, if and only if for every closed point $x \in |\cX|$,
    $\ihat_x^{*}\cG$ is pure of degree zero, resp.\ equal to zero.
  \end{enumerate}
\end{lem}%
\begin{proof}%
Lemma~\ref{lem:cohomology is conservative on bounded above objects for left complete}
shows that the second
  statement follows from the first by considering cohomology sheaves.

For the first statement, assume that
$\ihat_x^{*}\cF=0$
for every closed point~$x$ of~$\cX$.
Write $\cF = \quoteslim{j} \cF_j$ as a cofiltered limit of~$\cF_j\in\Coh(\cX)$.
Now fix an index~$j$.
For each point $x$, there is (by assumption)
an index $j'_x$ such that $\cF_{j'_x} \to \cF_j$ becomes zero after applying~$\ihat_x^*$. %
We claim that~$x$ has an open neighbourhood $\cU_x$ over which the restriction of the transition morphism $\cF_{j'_x} \to \cF_j$
vanishes.
Granting this,
we find that the union of these $\cU_x$ is equal to all of~$\cX$ (since it is an open
substack of~$\cX$ that contains every closed point).
Since $\cX$ is quasi-compact, 
we may therefore cover $\cX$ by finitely many of these $\cU_x$,
and thus find a single index $j'$ (some index lying above
each of the finitely many corresponding~$j'_x$)
such that $\cF_{j'} \to \cF_j$ vanishes.
This implies that~$\cF$ itself is the zero object in~$\ProCoh(\cX)$.

It remains to prove the claim.
Since~$\ihat^*_x$ is exact, by a consideration of the image sheaf, it suffices to prove that,
if~$\cG$ is a coherent sheaf on~$\cX$, and~$\ihat^*_x \cG = 0$,
then there exists an open neighborhood~$\cU_x$ of~$x$ such that $\cG|_{\cU_x} = 0$.
By definition, the coherent sheaf $\cG$ is supported on a closed algebraic substack of~$\cX$,
and replacing $\cX$ by this closed algebraic substack, we may furthermore assume
that~$\cX$ is itself algebraic.
Now $\ihat^*_x$ can be thought of as the functor of $\cI$-adic completion (with $\cI$ denoting
the coherent ideal sheaf
on~$\cX$ cutting out the underlying reduced substack structure on the closed subset~$\{x\}$),
and so this amounts to the standard fact
that if the $\cI$-adic completion of a coherent sheaf vanishes, then the coherent sheaf
itself vanishes in a neighbourhood of~$x$. 
For the sake of completeness (and since our completions are computed in the Pro-sense,
rather than literally, and since we are working on an algebraic stack rather
than a scheme), we recall the proof. To this end, we choose  a smooth surjection
$p:X = \Spec A \to \cX$ from a Noetherian affine scheme onto some open neighbourhood of~$x$;
since $p$ is open~\cite[\href{https://stacks.math.columbia.edu/tag/04XL}{Tag~04XL}]{stacks-project},
it suffices to prove that $p^*\cG$ vanishes after restriction to some open neighbourhood
of the closed subset~$p^{-1}(x)$ of~$X$.
Let $I$ denote the ideal in $A$ cutting out the underlying reduced 
subscheme structure on~$p^{-1}(x)$,
and let $M$ denote the finitely generated $A$-module corresponding to the coherent
sheaf~$p^*\cG$. 
Our assumption that $\ihat^*_x \cG  = 0$ implies that the morphism $M/I^n M \to M/IM$
coincides with the zero morphism
for sufficiently large~$n$; since this morphism is surjective, we see that its 
target vanishes, i.e.\ that $M = IM$. 
A standard application of Nakayama's lemma now implies that indeed the coherent sheaf
attached to~$M$ vanishes in some neighbourhood of~$p^{-1}(x)$, as required.
\end{proof}

\begin{lemma}\label{lem:testing support with completions}%
Let $z:\cZ \to \cX$ be a closed immersion of Noetherian formal algebraic stacks with affine diagonal.
Then a pro-coherent sheaf $\cF$ on $\cX$ is scheme-theoretically supported on~$\cZ$ %
if and only if, for every closed point~$x \in |\cX|$ with completion-at-$\{x\}$ morphism $i_x: \widehat \cX_x \to \cX$,
the pro-coherent sheaf $\ihat_x^* \cF$ is scheme-theoretically supported on~$\widehat \cZ_{x \cap |\cZ|}$. %
\end{lemma}
\begin{proof}
  If~$x \in |\cX|$ is a closed point, then the image under~$\ihat^*_x$ of the unit of adjunction $\cF \to z_* z^* \cF$
  is isomorphic to the unit $\ihat^*_{x}\cF \to z_{x, *} z^*_{x}(\ihat^*_{x} \cF)$, by Remark~\ref{rem:map of adjunctions}.
   Since~$\ihat_x^*$ is exact (by Proposition~\ref{prop:exact-adjoint-ihat-upper-star}), the kernel and cokernel of $\cF \to z_* z^*\cF$ vanish under~$\ihat_x^*$ if and only if $\ihat_x^*\cF$ 
   is scheme-theoretically supported on~$\widehat \cZ_{x \cap |\cZ|}$.
Thus the lemma follows from Lemma~\ref{lem:pro-coh-vanishes-if-pullbacks-do}~(1).
\end{proof}

\subsubsection{Pullbacks to versal rings}\label{subsubsec:versal pullbacks}

\begin{lemma}\label{lem: completed pullback at effective versal morphism}
Let~$\cX$ be a Noetherian algebraic stack with affine diagonal, let $x \in |\cX|$ be a finite type point, and let $i_x: \Spec R \to \cX$ be a morphism from the spectrum
of a complete Noetherian local ring~$R$ with finite residue field, sending the maximal ideal~$\fm$ to~$x$.
Let~$\Spf R$ be the formal $\fm$-adic spectrum of~$R$. %
Then the composite
\begin{equation}\label{eqn:coherent pullback composite}
\Coh(\cX) \xrightarrow{i_x^*} \Mod^{\fp}(R) \xrightarrow{\eqref{eqn:mod-fp-to-mod-c}} \Mod_c(R) = \Pro\Mod^{\fl}(R) = \Pro\Coh(\Spf R)
\end{equation}
is naturally isomorphic to the completed pullback for the composite $\Spf R \to \Spec R \xrightarrow{i_x} \cX$, where the first arrow is the completion at~$\{\fm\}$.
\end{lemma}
\begin{proof}
Let~$\cF \in \Coh(\cX)$, and let~$M \coloneqq  i_x^*(\cF)$, a finitely presented $R$-module.
Then the composite of~\eqref{eqn:coherent pullback composite} sends~$\cF$ to the object $\quoteslim{n} M/\fm^n M$ of $\Pro \Coh(\Spf R)$.
On the other hand, we can compute the completed pullback to~$\Spf R$ using the presentation
\[
\Spf R = \colim_n \Spec R/\fm^n.
\]
Then, by definition, we have $\widehat i_x^* \cF = \quoteslim{n} i_x^*(\cF)/\fm^n i_x^*(\cF)$, as desired.
\end{proof}

\begin{lemma}\label{lem: versal pullback is exact and faithful}
Let~$\widehat \cX$ be a Noetherian formal algebraic stack with affine diagonal. %
Assume that $|\widehat \cX| = \{x\}$ is a singleton, and that ~$\widehat \cX_{\red}$ is of finite type over~$\Spec \F$.
Let $v_x: \Spf R \to \widehat \cX$ be a morphism from the formal $\fm$-adic spectrum of a complete Noetherian local ring~$R$
with residue field~$k$,
and assume that~$v_x$ is versal to the induced morphism $\Spec k \to \cX$. 
Then
\[
\vhat_x^* : \Coh(\widehat \cX) \to \Pro\Coh(\Spf R) = \Mod_c(R)
\]
is exact and faithful.
\end{lemma}
\begin{proof}
Since the formation of versal rings and of the completed pullbacks $\vhat_x^*$ are both compatible with base change to closed immersions 
(using the commutativity of the leftmost square in~\eqref{to prove commutative for compatibility of completion and pushforward} for the latter), 
we reduce to the case that~$\widehat\cX$ is algebraic and locally of finite type over~$\Spec \cO$.
In this case $v_x$ is effective, by~\cite[\href{https://stacks.math.columbia.edu/tag/0DR1}{Tag~0DR1}]{stacks-project}, so by Lemma~\ref{lem: completed pullback at effective versal morphism}  that it suffices to prove that the usual coherent pullback for the map
\[
v_x: \Spec R \to \widehat\cX
\]
is exact and faithful.
This follows from ~\cite[\href{https://stacks.math.columbia.edu/tag/0DR2}{Tag~0DR2}]{stacks-project}, which shows that~$v_x$ is flat and surjective.
\end{proof}

\subsubsection{Completed pullback and internal~\texorpdfstring{$\Hom$}{Hom}}\label{subsubsec:completed-internal-hom}If~$\cX$ is a Noetherian algebraic stack, then by \cite[\href{https://stacks.math.columbia.edu/tag/0GQN}{Tag 0GQN}]{stacks-project}, $\Coh(\cX)$ has an internal~$\Hom$ functor $\uHom_{\Coh(\cX)}(\text{--}, \text{--})$.

\begin{lemma}\label{lem: versal pullback and internal Hom}%
Let~$\cX$ be an algebraic stack with affine diagonal, which is locally of finite type over~$\Spec \cO$.
Let $i_x: \Spf R \to \cX$ be a versal morphism at a closed point $x \in |\cX|$.
Then we have a natural $R$-linear isomorphism
\[
\widehat i_x^*\uHom_{\Coh(\cX)}(\text{--}, \text{--}) \isoto \Hom_{\Mod_c(R)}(\widehat i_x^*\text{--}, \widehat i_x^*\text{--}).
\] 
\end{lemma}
\begin{proof}
By \cite[\href{https://stacks.math.columbia.edu/tag/0DR1}{Tag~0DR1}]{stacks-project}, the map~$i_x$ is effective.
By Lemma~\ref{lem: completed pullback at effective versal morphism}, it suffices to show that we have a natural $R$-linear isomorphism %
\[
i_x^*\uHom_{\Coh(\cX)}(\text{--}, \text{--}) \isoto \Hom_{\Mod^{\fp}(R)}(i_x^*\text{--}, i_x^*\text{--}),
\] 
where~$i_x^*$ is the coherent pullback for the map $v_x: \Spec R \to \cX$.
Since~$i_x$ is flat, by \cite[\href{https://stacks.math.columbia.edu/tag/0DR2}{Tag~0DR2}]{stacks-project}, this is a consequence of \cite[\href{https://stacks.math.columbia.edu/tag/0GQP}{Tag 0GQP}]{stacks-project}. %
\end{proof}

\subsection{Coherent sheaves and formal functions
in the case of ``affine mod reductive'' quotient stacks} 
\label{subsec:affine mod reductive}
In this section we prove some results about stacks of the form $[\Spec B/G]$, where $B$ is a finite type algebra
over a Noetherian ring~$R$ and $G$ is a reductive group scheme over~$R$.
Some of our results  are deduced from those of~\cite{alper2023coherently}; since our setup is slightly
different to theirs (see Remark~\ref{rem:our definition of
  coherent}), we begin by relating our approach to completions of algebraic stacks to that taken in~\cite{BCGAGA} and~\cite{alper2023coherently}.

\subsubsection{\texorpdfstring{$\cO_{\widehat{\cX}}$}{O-X-hat}-modules}%
\label{subsec:completing an algebraic stack}
Suppose that $\cX$ is a Noetherian algebraic stack with affine diagonal, and that $\cZ$ is a
closed substack with coherent ideal sheaf~$\cI$. We write~$\cXhat$ for the 
 formal completion
$\widehat{\cX}$ of $\cX$ along~$\cZ$.
In this context,
\cite{BCGAGA} defines a sheaf of rings
\begin{equation}
\label{eqn:completed structure sheaf}
\cO_{\widehat{\cX}} \coloneqq  \limcommand_n \cO_{\cX}/\mathcal I^{n+1}
\end{equation}
on the lisse-\'etale
site of~$\cX$,
which is morally the structure sheaf of the completion~$\widehat{\cX}$,
as well as an associated abelian category $\Coh(\cO_{\widehat{\cX}})$ of coherent sheaves
of~$\cO_{\widehat{\cX}}$-modules. 
This allows one to also introduce the
stable $\infty$-category~$D^b_{\coh}(\cO_{\widehat{\cX}})$.   
Furthermore, writing $\cZ^{[n]}$ for the closed substack of~$\cX$ with ideal sheaf~$\cI^n$,
the functors $\cO_\cX/\cI^{n+1} \otimes_{\cO_{\widehat \cX}}\text{--}$ induce an equivalence of categories (see~\cite[Theorem~2.3]{BCGAGA}) %
\begin{equation}\label{eqn:BC fully faithful embedding}
\Coh (\cO_{\widehat \cX}) \iso \lim_n \Coh(\cZ^{[n]}).
\end{equation}

\begin{remark}
We note that the abelian category $\Coh(\cO_{\widehat{\cX}})$ is
denoted simply by~$\Coh(\widehat{\cX})$ in~\cite{BCGAGA}.  Similarly, the
stable $\infty$-category~$D^b_{\coh}(\cO_{\widehat{\cX}})$ ---  
or, rather, its underlying triangulated category ---
is denoted\footnote{Strictly speaking, \cite{BCGAGA} only considers
the ``+'' version of the derived category, which it denotes by~$D^+_{\coh}(\widehat{\cX})$.} 
$D^b_{\coh}(\widehat{\cX})$ in~\cite{BCGAGA}.
However, as indicated in Remark~\ref{rem:our definition of coherent},
the categories defined in~\cite{BCGAGA}
are larger than the category $D^b_{\coh}(\widehat{\cX})$
that we consider, and so to distinguish the two we always denote 
the categories of~\cite{BCGAGA} by $\Coh(\cO_{\widehat{\cX}})$ and~$D^b_{\coh}(\cO_{\widehat{\cX}})$.   
\end{remark}

Rather than forming the actual projective limit in
\eqref{eqn:completed structure sheaf} to obtain a sheaf of rings,
one can consider the formal projective limit, so as to obtain an object of
$\Pro\Coh(\widehat{\cX}).$
More precisely, 
since every $\cO_{\cX}$-module that is annihilated by a power of~$\cI$ is also an $\cO_{\widehat{\cX}}$-module, 
we have an exact and fully faithful functor
\begin{equation}\label{abelian coherent to OXhat}
\Coh(\widehat \cX)\to \Coh(\cO_{\widehat \cX}),
\end{equation}
and we can also define an exact and fully faithful functor 
\begin{equation}
\label{eqn:abelian complete to pro} 
\Coh(\cO_{\widehat{\cX}}) \hookrightarrow \Pro\Coh(\widehat{\cX}), \;\; \cF \mapsto \quoteslim{n} \cF/\mathcal I^{n+1}.
\end{equation}
(The claimed properties of this functor follow again
from~\cite[Thm.~2.3]{BCGAGA}, which implies full faithfulness, together with 
a standard Artin--Rees argument to deduce exactness, as in the proof of Proposition~\ref{prop:exact-adjoint-ihat-upper-star}.) %
The composition of~\eqref{eqn:abelian complete to pro} and~\eqref{abelian coherent to OXhat} is the natural inclusion of~$\Coh(\widehat \cX)$ in its Pro-completion~$\Pro\Coh(\widehat\cX)$.

As remarked on \cite[p.~1]{BCGAGA}, the sheaf $\cO_{\widehat{\cX}}$ is
a flat sheaf of $\cO_{\cX}$-algebras.
Thus $\cO_{\widehat{\cX}}\otimes_{\cO_{\cX}}\text{--}$ defines an
exact functor 
\begin{equation}
\label{eqn:completing abelian sheaves}
\Coh(\cX) \to \Coh(\cO_{\widehat{\cX}}).
\end{equation}
By definition, the composite of ~\eqref{eqn:completing abelian sheaves} and~\eqref{eqn:abelian complete to pro} is the functor~$\ihat^{*}:\Coh(\cX)\to\Pro\Coh(\cXhat)$ defined in Section~\ref{subsubsec:completed pullbacks}.

It follows from Corollary~\ref{cor:pro t-exact derived functors}~(2), applied to $\cC = \Coh(\cO_{\widehat \cX})$ and $\cC' = \Coh(\widehat \cX)$,
that~\eqref{eqn:abelian complete to pro}
induces a $t$-exact functor
\[
D^b(\Coh \cO_{\widehat \cX}) \to \Pro D^b(\Coh \widehat \cX).
\]
Remembering that
$D^b\bigl(\Coh \cO_{\widehat{\cX}}\bigr) \iso D^b_{\coh}(\cO_{\widehat{\cX}})$
and
$D^b\bigl(\Coh \widehat{\cX}\bigr) \iso D^b_{\coh}(\widehat{\cX})$, we obtain a $t$-exact functor
\begin{equation}
\label{eqn:derived complete to pro} 
D^b_{\coh}(\cO_{\widehat{\cX}}) \to \Pro D^b_{\coh}(\widehat{\cX}),
\end{equation}
whose composition with %
the $t$-exact functor 
\begin{equation}\label{derived coherent to OXhat}
D^b_{\coh}(\widehat \cX) \to D^b_{\coh}(\cO_{\widehat \cX})
\end{equation}
associated to~\eqref{abelian coherent to OXhat} is
again the natural inclusion of~$D^b_{\coh}(\widehat \cX)$ to its $\Pro$-completion.

\begin{prop}
\label{prop:completion comparison}
The functor~{\em \eqref{eqn:derived complete to pro}} is fully faithful.
\end{prop}
\begin{proof}%
Either by inspection, or by an application of Corollary~\ref{cor:t-exact derived functors}, one sees
that~\eqref{eqn:derived complete to pro} must coincide with the functor~$q$ of~\eqref{eqn:defn-of-q}
(taking $\cC$ to be $\Coh(\widehat\cX)$
and $\cB$ to be $\Coh(\cO_{\widehat{\cX}})$, viewed as an exact abelian subcategory of~$\Pro \cC$ via~\eqref{eqn:abelian complete to pro}).
We will prove its full faithfulness by applying the criterion of Lemma~\ref{lem:full faithfulness criterion}. 
Let $\cG$ be an object of $\Coh(\cO_{\widehat{\cX}}).$
Then $\cG$ maps to $\quoteslim{n} \cG/\cI^{n+1}$ under~\eqref{eqn:abelian complete to pro}. 
It follows from
~\cite[\href{https://stacks.math.columbia.edu/tag/0BKY}{Tag 0BKY}]{stacks-project}
and the Artin--Rees lemma that in fact
$$\cG \iso R\limcommand_n \cG/\cI^{n+1}$$
in the unbounded derived category~$D(\cO_{\widehat{\cX}})$
of~$\cO_{\widehat{\cX}}$-modules,
which is to say, the corresponding $R^1\limcommand_n$ vanishes.
Thus the natural morphism
$$
\RHom_{D(\cO_{\widehat{\cX}})}(\cF,\cG)
\to
\Rlim_n \RHom_{D(\cO_{\widehat{\cX}})}(\cF,\cG/\cI^{n+1})
$$
is an isomorphism,
for any object $\cF$ of~$D(\cO_{\widehat{\cX}}).$
Restricting to objects $\cF$ of the full subcategory $D^b_{\coh}(\cO_{\widehat{\cX}})$
of~$D(\cO_{\widehat{\cX}})$, we see that the hypothesis of Lemma~\ref{lem:full faithfulness
criterion} is satisfied, as required.
\end{proof}

It will be useful to give an alternative definition of the abelian
category~$\Coh(\cO_{\widehat{\cX}})$ (as a full subcategory of $\Pro \Coh(\widehat{\cX})$);
e.g.\ it allows us to give a simple
proof of Lemma~\ref{lem:completed pullback preserves O_X modules}
below.
This definition also has the merit of applying to a more general class of formal algebraic
stacks (although we don't have cause to use it in that greater level of generality).

\begin{defn}\label{defn:O_X modules in general} %
Let~$\cY$ be a Noetherian formal algebraic stack with affine diagonal.
We define $\Coh(\cO_\cY)$ as the full subcategory of $\Pro \Coh(\cY)$ whose objects are the pro-coherent sheaves~$\cF$ such that, for all morphisms
$i : \cZ \to \cY$, where~$\cZ$ is a quasicompact algebraic stack with affine diagonal, the pullback $i^*\cF$ is coherent.
\end{defn}

\begin{rem}
As the preceding discussion already implied, 
in the case when~$\widehat \cX$ is the completion of a Noetherian algebraic stack~$\cX$ at a closed substack~$\cZ^{[0]}$ with ideal sheaf~$\cI$, 
the category $\Coh(\cO_{\widehat \cX})$, as defined in Definition~\ref{defn:O_X modules in general}, coincides with the essential image of~\eqref{eqn:abelian complete to pro}.

In fact, we have $\cX = \colim_n \cZ^{[n]}$, where $i_{n}: \cZ^{[n]} \to \cX$ is the closed substack with ideal sheaf~$\cI^{n+1}$. 
Every morphism $\cZ \to \cX$ from a quasicompact algebraic stack factors through some $\cZ^{[n]}$ (this is a consequence of~\cite[Lem.\ 4.4]{Emertonformalstacks}, after passing to a cover by a quasicompact scheme),
and so the essential image of~\eqref{eqn:abelian complete to pro} is contained in $\Coh(\cO_{\widehat \cX})$.
Conversely,
if~$\cF \in \Coh(\cO_{\widehat \cX})$ 
(as defined in Definition~\ref{defn:O_X modules in general}),
then the sequence~$(i_n^*\cF)_n$ defines an object of the right-hand side of~\eqref{eqn:BC fully faithful embedding II}, and so
of its left-hand side.
The
isomorphism $\cF \isoto \quoteslim{n} i_{n, *}i_n^* \cF$ from Lemma~\ref{lem:pro-coherent-as-limit-of-pullbacks} implies that $\cF$ is in the essential image of~\eqref{eqn:abelian complete to pro}.
\end{rem}

\begin{lemma}\label{lem:completed pullback preserves O_X modules}
Let $f: \cY \to \cX$ be a morphism of Noetherian formal algebraic stacks with affine diagonal.
Then the completed pullback functor~$\widehat f$ restricts to a functor
\[
\widehat f : \Coh(\cO_{\cX}) \to \Coh(\cO_{\cY}).
\]
\end{lemma}
\begin{proof}
Let $i : \cZ \to \cY$ be a morphism from a quasicompact algebraic stack with affine diagonal~$\cZ$, and let $\cF \in \Coh(\cO_{\cX})$.
We need to prove that $i^* \widehat f^* \cF$ is coherent.
By Remark~\ref{rem:pullback and completed pullback}, it is isomorphic to $\ihat^* \widehat f^* \cF$, which is the same as $\widehat{(f \circ i)}^*\cF$, by Remark~\ref{completed pullback is functorial}.
Again by Remark~\ref{rem:pullback and completed pullback}, we conclude that $i^* \widehat f^* \cF \cong (f \circ i)^*\cF$, which is coherent, because $\cF \in \Coh(\cO_{\cX})$.
\end{proof}

\subsubsection{Reductive groups}\label{subsubsec:reductive}%
We will require some coherent sheaf theory on stacks of the form $[\Spec B/G]$, where $B$ is a finite type algebra over a Noetherian ring~$R$ and $G$ is a reductive group scheme over~$R$.
The definition of ``reductive group scheme'' in the literature is not always consistent, due to different conventions about disconnected groups.
However, the notion of a connected reductive group is well-defined; indeed, over an an arbitrary base scheme~$S$, a connected reductive group~$G/S$ is a smooth $S$-affine group scheme all of whose geometric fibres are connected reductive group.

Turning to the not-necessarily-connected case,  if $G/S$ is a finitely presented smooth group scheme over an arbitrary base scheme~$S$, then there is a unique open subgroup scheme~$G^{\circ}\subseteq G$ with the property that for all points~  $s\in S$, the fibre $(G^{\circ})_s$ is equal to the identity component $(G_{s})^{\circ}$  of~$G_{s}$.
(See for example ~\cite[\S 3.1]{MR3362641}.) The formation of~$G^{\circ}$ commutes with arbitrary base change.
By~\cite[Prop.\ 3.1.3]{MR3362641}, if  $G/S$ is separated, smooth, and of finite presentation,  and each fibre~$G^{\circ}_s$ is reductive, then $G^{\circ}$ is open and closed in~$G$, and the quotient~$G/G^{\circ}$ is a separated \'etale group scheme of finite presentation.

For us, a ``reductive group scheme'' will always mean the following.
\begin{defn}
  \label{defn:our-defn-reductive}A group scheme ~$G/S$ is reductive if~$G$ is smooth and affine,  $G^{\circ}$ is connected reductive, and~$G/G^{\circ}$ is finite \'etale.
\end{defn}

The following result of Alper gives an alternative characterisation of reductive group schemes, which we use in the proof of Theorem~\ref{thm:summary-of-AHL} below.
Recall that by definition~\cite[Defn.\ 9.1.1]{MR3272912}, a group scheme~$G/S$ is \emph{geometrically reductive} if and only if it is finitely presented, flat and separated over~$S$ and the morphism $BG\to S$ is adequately affine in the sense of~\cite[Defn.\ 4.1.1]{MR3272912}.
\begin{prop}\label{prop:geometrically-reductive-equivalent}Suppose that~$G/S$ is a smooth affine group scheme, and that~$G^{\circ}$ is \emph{(}connected\emph{)} reductive. %
 Then~$G$ is reductive \emph{(}in the sense of Definition~{\em \ref{defn:our-defn-reductive})} if and only if it  is geometrically reductive. %
  \end{prop}
\begin{proof}
As noted above, since~$G^{\circ}$ is connected reductive, the component group scheme ~$G/G^{\circ}$ is separated. The result is then  immediate from ~\cite[Thm.\ 9.7.6]{MR3272912}.%
\end{proof}

Now let $B$ be a finite type $R$-algebra for a Noetherian ring~$R$,
let $G$ be a  reductive group scheme over~$R$, and suppose we are given an action 
of $G$ on $\Spec B$ (in the category of $R$-schemes); that is, a morphism
$$G \times_{\Spec R} \Spec B \to \Spec B$$
of $R$-schemes satisfying the usual axioms for an action.
We then write $[\Spec B/ G]$ for the usual stack quotient.  
In more detail, the pair of morphisms 
$G \times_{\Spec R} \Spec B \rightrightarrows \Spec B$ --- the first being the obvious projection,
and the second being the action morphism --- are the source and target map for a smooth groupoid
over $\Spec B$, and $[\Spec B/ G]$  is the quotient stack of $\Spec B$ with respect to this groupoid.

The following theorem summarizes results of Alper and of Alper--Hall--Lim, but since the paper~\cite{alper2023coherently} only states them in the case that~$G$ is connected reductive, we explain the straightforward deduction of the non-necessarily-connected reductive case from their results. We remind the reader  that a finite type morphism
$f:\cY\to\cX$ of Noetherian algebraic stacks is \emph{cohomologically proper} \cite[Def.~2.1(2)]{alper2023coherently}
if for every coherent sheaf~$\cF$ on~$\cY$, the derived pushforward
$f_{*}\cF$ lies in $D^{+}_{\coh}(\cX)$. 

\begin{thm}
  \label{thm:summary-of-AHL}%
Let $G$ be a reductive group scheme over  a Noetherian ring~$R$, and let  $B$ be a finite type $R$-algebra with an action of~$G$.
Then:
\begin{enumerate}

  \item \label{item:finite-type-invariants} $B^{G}$ is of finite type over~$R$ \emph{(}and is in particular Noetherian\emph{)}.
  \item \label{item:universall-coh-proper} $[\Spec B/ G]\to \Spec B^G$ is %
    cohomologically proper. %
\item \label{item:BGG-derived-formal-functions}Suppose that~$I\subseteq B$ is a $G$-equivariant ideal, and that $B^G$ is $I^G$-adically complete. 
Write~$\cY=[\Spec B/G]$, $\cY^{[n]}=[\Spec (B/I^{n+1})/G]$, and let $\widehat{\cY} = \colim_n \cY^{[n]}$ 
be the completion of~$\cY$ along~$i:\cY^{[0]}\into \cY$. Then the functor 
\begin{equation}\label{eqn:completed pullback in AHL context}
\widehat{i}^*: D^b_{\coh}(\cY) \to \Pro D^b_{\coh}(\widehat{\cY})
\end{equation}
is $t$-exact %
and fully faithful, and the functor
\begin{equation}\label{eqn:BC fully faithful embedding II}
\cO_{\widehat \cY} \otimes_{\cO_\cY} \text{--} : \Coh(\cY) \to \Coh(\cO_{\widehat \cY})
\end{equation}
is fully faithful.
If furthermore~$G$ is a linearly reductive closed subgroup of $\GL_{n, R}$, then the essential image of~\eqref{eqn:completed pullback in AHL context}
is $D^b_{\coh}(\cO_{\widehat \cY})$.  
\end{enumerate}
\end{thm}
\begin{proof}%
   Part~\ref{item:finite-type-invariants} is immediate from~\cite[Rem.\ 9.1.5, Thm.\ 6.3.3]{MR3272912}. %

We now turn to part~\ref{item:universall-coh-proper}. By\cite[Ex.~2.7]{alper2023coherently}, the morphism $[\Spec B/ G^{\circ}]\to \Spec B^{G^{\circ}}$ is cohomologically proper. Furthermore, the morphism $\Spec B^{G^{\circ}}\to\Spec B^G$ is
finite and therefore cohomologically proper, so that the composite  $[\Spec B/ G^{\circ}]\to \Spec B^G$ is cohomologically proper. Since
the map \[[ \Spec B/ G^{\circ}]\to [\Spec B/ G]\] is finite \'etale and surjective, and $[\Spec B/G]\to \Spec B^{G}$ is of finite type, the claim follows from ~\cite[Thm.\ 6.1(1)]{alper2023coherently}. %

Finally, we consider part~\ref{item:BGG-derived-formal-functions}.
The $t$-exactness of ~$\ihat^{*}$ is Proposition~\ref{prop:exact-adjoint-ihat-upper-star}.
As noted after \eqref{eqn:completing abelian sheaves}, it follows immediately from the definitions that  we may factor~$\widehat{i}^*$
as the composite %
\[
D^b_{\coh}(\cY) \xrightarrow{\cO_{\widehat{\cY}}\otimes^L_{\cO_{\cY}}\text{--}} D^b_{\coh}( \cO_{\widehat \cY}) \xrightarrow{\eqref{eqn:derived complete to pro}} \Pro D^b_{\coh}(\widehat \cY).
\]
Both arrows are $t$-exact, and
the second is fully faithful, by Proposition~\ref{prop:completion comparison}. 
So, to prove that~\eqref{eqn:completed pullback in AHL context} is fully faithful (resp.\ an equivalence if~$G$ is linearly reductive and closed in $\GL_{n, R}$),
it suffices to show that $\cO_{\widehat \cY} \otimes_{\cO_\cY}^L\text{--}$
is fully faithful (resp.\ an equivalence if~$G$ is linearly reductive and closed in $\GL_{n, R}$).
This will also imply that~\eqref{eqn:BC fully faithful embedding II}
is fully faithful, and complete the proof of the proposition.

It suffices therefore to show that the pair~$(\cY,\cY^{[0]})$ satisfies \emph{derived formal functions} (resp.\ is \emph{derived coherently complete}) 
in the sense of~\cite[Defn.\ 3.2]{alper2023coherently}.
By~\cite[Prop.\ 3.5]{alper2023coherently}, it is in turn enough to show the pair~$(\cY,\cY^{[0]})$ satisfies \emph{formal functions}
(resp.\ satisfies formal functions and is \emph{coherently complete}) 
in the sense of~\cite[Defn.\ 3.2]{alper2023coherently}. %

We begin by showing that $(\cY, \cY^{[0]})$ satisfies formal functions.
Since the preimage of $[\Spec (B/I)/G]$ in $[\Spec B/G^{\circ}]$ is $[\Spec (B/I)/G^\circ]$, %
 \cite[Thm.\ 6.1(2)]{alper2023coherently}
reduces us to showing that the pair $([\Spec B/G^{\circ}],[\Spec (B/I)/G^\circ])$ satisfies formal functions.
Noting that~$B^{G^{\circ}}$ is $I^{G^{\circ}}$-adically complete by Lemma~\ref{lem:invariants-finite-group} below,  %
this is part of~\cite[Cor.\ 4.9]{alper2023coherently}. 

Finally, under the additional assumption that~$G$ is linearly reductive and closed in $\GL_{n, R}$, 
the fact that~$(\cY, \cY^{[0]})$ is coherently complete
follows from~\cite[Thm.\ 1.6]{alper2023etale}, %
bearing in mind that
$[\Spec B/G] \to \Spec B^G$ is a good moduli space, by~\cite[Thm.\ 13.2]{MR3237451}, and that
$\cY^{[0]} = [\Spec(B/I)/G]$ has the resolution property, by~\cite[Rem.\ 2.5]{alper2023etale}, since~$G \to \Spec(R)$ is assumed to be embeddable.
\end{proof}

The following lemma and its proof are presumably standard.
\begin{lem}
  \label{lem:invariants-finite-group}Let~$G$ be a finite group acting on a ring~$A$, and let $I\subseteq A$ be a $G$-equivariant ideal.
If~$A^{G}$ is $I^G$-adically complete, then ~$A$ is $I$-adically complete.
\end{lem}
\begin{proof}
  Since~$A$ is finite over~$A^G$ (being both finite type and integral), $A$ is~$I^GA$-adically complete.
It therefore suffices to show that $\rad(I^GA)=\rad(I)$, i.e.\ that $I\subseteq \rad(I^GA)$.
To this end, let~$x\in I$ be arbitrary.
Then~$x$ is a root of the monic polynomial $\prod_{g\in G}(X-g(x))$, all of whose non-leading coefficients are contained in~$I^{G}$.
Thus $x^{|G|}\in I^GA$, as required.
\end{proof}

We now record some consequences of Theorem~\ref{thm:summary-of-AHL}.
\begin{cor}
\label{cor:coherent exts for affine mod reductive-I} Let $G$ be a reductive group scheme over  a Noetherian ring~$R$,  let  $B$ be a finite type $R$-algebra with an action of~$G$, and let 
$\cF$ and $\cG$ be  coherent sheaves on~$\cY= [\Spec B/ G]$.  Then for each $n \geq 0,$
$\Ext^n_{\cO_{\cY}}(\cF,\cG)$ is of finite type as a $B^G$-module.
\end{cor}
\begin{proof} Since the morphism $f:\cY \to \Spec B^G$ is  cohomologically proper (by Theorem~\ref{thm:summary-of-AHL}~\ref{item:universall-coh-proper}), we see that each cohomology group
$H^{i}(\cY,\cH)$ is of finite 
type over $B^G$, for any coherent sheaf~$\cH$ on~$\cY$. %
In particular, we can apply this observation to the  Ext sheaves $\underline{\Ext}_{\cO_{\cY}}^{\bullet}(\cF,\cG)$ (which are coherent, since $\cF$ and~$\cG$ are).
The result then follows from  the sheaf-Ext-to-Ext spectral sequence
\[E_2^{p,q} \coloneqq  H^p\bigl(\cY, \underline{\Ext}^{q}_{\cO_{\cY}}(\cF,\cG)\bigr)
\implies \Ext^{p+q}_{\cO_{\cY}}(\cF,\cG). \qedhere \]
\end{proof}

\begin{acor}
\label{cor:coherent exts for affine mod reductive} Let $G$ be a reductive group scheme over  a Noetherian ring~$R$,  let  $B$ be a finite type $R$-algebra with an action of~$G$, and let 
$\cF$ and $\cG$ be  coherent sheaves on~$\cY= [\Spec B/ G]$.  Then for each $n \geq 0,$
$\Ext^n_{\cO_{\cY}}(\cF,\cG)$ is of finite type as an $\End_{\cO_{\cY}}(\cF)$-module and as an $\End_{\cO_\cY}(\cG)$-module. 
\end{acor}
\begin{proof}
The action of $B^G$ on $\cG$ makes $\End_{\cO_{\cY}}(\cG)$ a $B^G$-algebra,
and the action of $\End_{\cO_{\cY}}(\cG)$ on~$\Ext^n_{\cO_{\cY}}(\cF,\cG)$ is compatible
with this algebra structure and the {\em a priori} action of $B^G$ on the~$\Ext^n$.
The corollary thus follows from Corollary~\ref{cor:coherent exts for affine mod reductive}.
The proof for~$\End_{\cO_\cY}(\cF)$ is the same. 
\end{proof}

\begin{acor}\label{cor:completing RHoms} As above, let~$R$ be a Noetherian
  ring, and let~$B$ be a finitely generated $R$-algebra endowed with an action
  of a reductive group scheme~$G$ over~$\Spec R$. Furthermore, let $I\subseteq B$
  be a $G$-equivariant ideal, let $\widehat{B^G}$ denote the $I^G$-adic
  completion of~$B^G$, and
  let $\widehat{\cY}$ denote the completion of $\cY\coloneqq  [\Spec B/G]$ along
  the zero locus $i:[\Spec (B/I)/G]\into\cY$ of~$I$.
Let $\cF$ and $\cG$ be  coherent sheaves on~$\cY$. Then for each $n \geq 0,$  the functor
$\widehat{i}^{*}$ induces an isomorphism
\[\Ext^n_{D_{\coh}^b(\cY)}(\cF,\cG)^\wedge \iso
\Ext^n_{\Pro D_{\coh}^b(\widehat{\cY})}(\widehat{i}^* \cF, \widehat{i}^*\cG),\]
where the $(\text{--})^\wedge$ in the source indicates the $I^G$-adic completion.
\end{acor}%

\begin{proof}%

The morphism $R \to B$ factors through ~$B^G$, and so
we may rewrite the $G$-action on $\Spec B$ in the form
$$(G\times_{\Spec R} \Spec B^G) \times_{\Spec B^G} \Spec B \to \Spec B.$$  Since~$B^G$ is Noetherian (by Theorem~\ref{thm:summary-of-AHL}~\eqref{item:finite-type-invariants}), we may replace $R$ by $B^G$, and $G$ by $G\times_{\Spec R} \Spec B^G$,
and so we assume from now on that $R = B^G$.
The $I^G$-adic completion $\widehat{B^G}$ of $B^G$ is flat over $B^G$ 
(since $B^G$ is Noetherian). 
We then set 
$B'  =\widehat{B^G}\otimes_{B^G} B,$
$G' =  \Spec \widehat{B^G}\times_{\Spec B^G}G,$
and
$\cY' \coloneqq  \Spec \widehat{B^G} \otimes_{\Spec B^G} \cY = [\Spec B'/G'].$ 
From the flatness of $\widehat{B^G}$ over~$B^G$, we see that
$(B')^{G'} = \widehat{B^G} \otimes_{B^G} B^G = \widehat{B^G}.$ We write $i':[\Spec (B'/IB')/G']\into\cY'$ for the base-change of~$i$.

Let $f: \cY' \to \cY$ denote the canonical morphism,
and write $f^*$ for the induced pullback morphism on coherent sheaves.
Since $f$ is flat, %
the derived pullback ~$f^{*}$ is a $t$-exact functor
$$f^*:D^b_{\coh}(\cY) \to D^b_{\coh}(\cY').$$
For any $\cF,\cG\in\Coh(\cY)$,  flat base-change shows that  %
then $f^*$ induces  isomorphisms %
\begin{equation}
\label{eqn:f base-change}
\widehat{B^G} \otimes_{B^G} 
\Ext^n_{D^b_{\coh}(\cY)}(\cF,\cG) \iso \Ext^n_{D^b_{\coh}(\cY')}(f^*\cF, f^*\cG).
\end{equation}

If we let $\widehat{\cY'}$ denote the completion of $\cY'$ along the zero locus
of $IB',$ then $f$ induces an isomorphism $\widehat{\cY'} \iso \widehat{\cY}$.
Considering the commutative diagram %
\[\begin{tikzcd}
\Coh(\cY) &  \Coh(\cY') \\
\Pro\Coh(\widehat{\cY}) & \Pro \Coh( \widehat{\cY'})
  \arrow["f^*", from=1-1, to=1-2]
  \arrow["\widehat f^*", "\sim"', from=2-1, to=2-2]
  \arrow["\ihat'^*", from=1-2, to=2-2]
  \arrow["\ihat^*", from=1-1, to = 2-1]
\end{tikzcd}\]
and noting that~$\ihat'^*$ is fully faithful by Theorem~\ref{thm:summary-of-AHL}~\eqref{item:BGG-derived-formal-functions},
we deduce from ~\eqref{eqn:f base-change} that~$\ihat^{*}$ induces isomorphisms %
\begin{multline*}
\widehat{B^G} \otimes_{B^G} 
\Ext^n_{D^b_{\coh}(\cY)}(\cF,\cG) \iso
\Ext^n_{\Pro D^b_{\coh}(\widehat{\cY'})}(\ihat'^*f^*\cF, \ihat'^*f^*\cG)
\\
\iso 
\Ext^n_{\Pro D^b_{\coh}(\widehat{\cY'})}(\widehat f^* \ihat^* \cF, \widehat f^* \ihat^*\cG)
\iso
\Ext^n_{\Pro D^b_{\coh}(\widehat{\cY})}(\ihat^*\cF, \ihat^*\cG).
\end{multline*}
Corollary~\ref{cor:coherent exts for affine mod reductive}
shows that each $\Ext^n_{D_{\coh}^b(\cY)}(\cF,\cG)$ is of finite type as a $B^G$-module.
Thus tensoring it with~$\widehat{B^G}$ over $B^G$ coincides with $I^G$-adically completing it.
This proves the corollary.
\end{proof}

\subsection{Another version of formal functions}\label{subsec:another version of formal functions}

Let $R$ be a complete Noetherian local $\cO$-algebra with maximal ideal~$\m$ and finite residue field, and let $S$ be a finite type $R$-algebra
equipped with an action of a reductive group scheme $G$ over~$R$ (where ``reductive'' is as in Definition~\ref{defn:our-defn-reductive}).  
We may then form the algebraic stack $\fX \coloneqq  [\Spec S/G]$ over~$R$.
We furthermore assume that $S^G$ is finite over~$R$ and local (so that
$S^G$ is again a complete Noetherian local ring). 
This assumption implies
that $\fX$ has a unique closed point~$z_0$, corresponding to a $G$-equivariant maximal ideal~$\fm_S$.
This is because $|\cZ| \to |\Spec S^G|$ induces a bijection on closed points, by~\cite[Thm.\ 5.3.1~(5)]{MR3272912};
recall from~\cite[Rem.\ 9.1.5]{MR3272912} that the map $\cZ \to \Spec S^G$ is an adequate moduli space (using that $G \to \Spec R$ is geometrically reductive 
in the sense of \emph{loc.\ cit.},
by Proposition~\ref{prop:geometrically-reductive-equivalent}).

Then:
\begin{enumerate}
  \item We write $\cX$ to denote the $\fm_R$-adic completion of~$\fX$, and $\widehat \cX$ to denote the $\fm_S$-adic completion of~$\fX$. %
Note that $|\widehat \cX|$ consists of a single point (corresponding to the unique
closed point $z_0$ of $\fX$).
  \item We say that a morphism $\Spf A \to \widehat \cX$, for~$A$ a complete Noetherian local $\cO$-algebra with finite residue field, equipped with its $\fm_A$-adic topology, 
  is \emph{versal} if it is versal to the induced morphism $\Spec A/\fm_A \to \widehat \cX$, 
in the sense of e.g.~\cite[Defn.\ 2.2.9]{EGstacktheoreticimages}.
  Precisely, this means that
  given a commutative square consisting of the solid arrows in the
  following diagram
  \begin{equation}\label{eqn: versality diagram}
  \xymatrix{ \Spec B_0 \ar[r]\ar[d] & \Spec B_1 \ar[d]\ar@{-->}[dl] \\ \Spf A \ar[r] & \widehat \cX}
  \end{equation}
  where $B_1\to B_0$ is a surjection of finite type Artinian local $R$-algebras,
and where the corresponding homomorphisms $A \to B_0$ and $B_1 \to B_0$ 
  induce isomorphisms on residue fields, then the dotted arrow can be filled in. 
  
  \item We write $(\Ver/\widehat \cX)$ for the category of versal morphisms with target~$\widehat \cX$, where a morphism
  \[
  (g: \Spf B \to \widehat \cX) \to (f: \Spf A \to \widehat \cX)
  \]
  is given by a commutative diagram of morphisms of formal algebraic stacks
$$\xymatrix{\Spf B \ar[r] \ar^-g[dr] & \Spf A \ar^-f[d] \\
& \widehat \cX},$$
or, equivalently,  
a local morphism $\alpha: A \to B$ of complete Noetherian local $\cO$-algebras together with an isomorphism $f \circ (\Spf \alpha) \isoto g$.
\end{enumerate}

\begin{rem}\label{rem:what it means to be versal at the closed point}
It will sometimes be useful to consider a versal morphism $\Spf A \to \widehat \cX$ as a versal morphism to the unique closed point of~$\cX$,
by implicitly composing it with the completion morphism $\widehat \cX \to \cX$,
Accordingly, we sometimes refer to the objects of~$(\Ver/\widehat \cX)$ as \emph{versal morphisms to the closed point of~$\cX$}.
\end{rem}

In this subsection we give a description of $\Coh(\cX)$ in terms of a category of coherent sheaves on $(\Ver/\widehat \cX)$, whose definition we turn to next. %
We begin by noting that 
if $\Spf A \to \widehat \cX$ is an object of $(\Ver/ \widehat \cX)$, then we have equivalences
\[
\Coh(\cO_{\Spf A}) \isoto \lim_n \Coh(\Spec A/\fm_A^n) \xleftarrow{\sim} \Coh(\Spec A).
\]
Furthermore,
if $\varphi: \Spf B \to \Spf A$ is a morphism in $(\Ver/ \widehat \cX)$, then the completed pullback 
\[
\widehat \varphi^* : \Pro \Coh (\Spf A) \to \Pro \Coh (\Spf B)
\]
restricts to a functor
\[
\widehat \varphi^*: \Coh(\cO_{\Spf A}) \to \Coh(\cO_{\Spf B})
\]
by Lemma~\ref{lem:completed pullback preserves O_X modules}.

\begin{defn}\label{defn: versal categories}
A coherent sheaf~$\cF$ on $(\Ver/\widehat \cX)$ consists of the following data:
\begin{itemize}
\item 
for every object $f: \Spf A \to \widehat \cX$ of~$(\Ver/\widehat \cX)$, 
we have a finitely presented $A$-module $\cF(x)$ (the sections of $\cF$ over $\Spf A$);
\item given a morphism $\varphi: g \to f$ of $(\Ver/\widehat \cX)$, lying over
a morphism $\Spf B \to \Spf A$,
or in more usual terms, a $2$-commutative diagram
$$\xymatrix{\Spf B \ar^-\varphi[r]\ar^-g[dr] & \Spf A \ar^-f[d] \\
& \widehat \cX},$$
we have 
a pullback isomorphism
$c_{\varphi}: \widehat\varphi^* \cF(g) = B\otimes_A \cF(g) \iso \cF(g);$ these pullback isomorphisms satisfy the following cocycle condition:
\item 
given morphisms $\psi:h \to g$ and $\varphi: g \to x$,
lying over $\Spf C \to \Spf B \to \Spf A$,
the isomorphism $c_{\varphi\circ \psi}$ coincides 
with the composite $c_{\psi} \circ (\id_C \otimes_B c_{\varphi})$,
after we identify $C\otimes_A \cF(f)$ and $C\otimes_B (B\otimes_A \cF(f))$
in the natural manner.
\end{itemize}

If $\cF$ and $\cG$ are two coherent sheaves on $(\Ver/\widehat \cX)$,
a morphism $\alpha: \cF \to \cG$ is defined to be a collection of
morphisms $\alpha_f: \cF(f) \to \cG(f)$, compatible in an evident manner
with the pullback isomorphisms~$c_\varphi$.
We write~$\Coh(\Ver/\widehat \cX)$ for the resulting category of coherent sheaves on~$(\Ver/\widehat \cX)$.
\end{defn}

\begin{lemma}\label{lem:versal pullback is exact and faithful II}
Let $f : \Spf A \to \widehat \cX$ be an object of $(\Ver/\widehat \cX)$, and let $\cF \in \Coh(\cO_{\widehat \cX})$.
Then $\widehat f^* \cF$ is an object of $\Coh(\cO_{\Spf A}) = \Coh(\Spec A)$, and the resulting functor
\begin{equation}\label{eqn:versal pullback is exact and faithful II}
\widehat f^* : \Coh(\cO_{\widehat \cX}) \to \Coh(\Spec A)
\end{equation}
is exact and faithful.
\end{lemma}
\begin{proof}
The coherence of $\widehat f^* \cF$ is a consequence of Lemma~\ref{lem:completed pullback preserves O_X modules}.
The exactness of $\widehat f^*$ is immediate from Lemma~\ref{lem: versal pullback is exact and faithful}, which implies that~$\widehat f^*$ is exact on the larger category~$\Pro \Coh(\widehat \cX)$.
There remains to prove that~$\widehat f^*$ is faithful, or equivalently (because of the exactness) conservative. 
The isomorphism $\cF \isoto \quoteslim{i} \cF/\fm_S^{i+1} \cF$ arising from Lemma~\ref{lem:pro-coherent-as-limit-of-pullbacks} exhibits any $\cF \in \Coh(\cO_{\widehat \cX})$ as a cofiltered limit
of objects of $\Coh(\widehat \cX)$, with surjective structure maps.
If~$\cF \ne 0$, then $\cF/\fm_S^{i+1} \cF \ne 0$ for some~$i$, and then $\widehat f^*(\cF/\fm_S^{i+1}\cF) \ne 0$, by Lemma~\ref{lem: versal pullback is exact and faithful}.
Hence $\widehat f^* \cF \ne 0$, since it surjects onto $\widehat f^*(\cF/\fm_S^{i+1} \cF)$.
This concludes the proof.
\end{proof}

By Lemma~\ref{lem:versal pullback is exact and faithful II},
there is a faithful restriction functor
\begin{equation}\label{eqn: versal restriction functor II}
\Coh(\cO_{\widehat \cX}) \to \Coh(\Ver/\widehat \cX)
\end{equation}
sending
a coherent $\cO_{\widehat \cX}$-module $\cF$ to the coherent sheaf on $(\Ver/\widehat \cX)$ whose sections on 
$f: \Spf A \to \widehat \cX$ are given by $\widehat f^*\cF(\Spf A)$. 
Writing~$i : \widehat \cX \to \cX$ for the completion map, we also have the composite
\begin{equation}\label{eqn: versal restriction functor I}
\Coh(\cO_{\cX}) \xrightarrow{\ihat^*} \Coh(\cO_{\widehat \cX}) \xrightarrow{\eqref{eqn: versal restriction functor II}} \Coh(\Ver/\widehat \cX).
\end{equation}

\begin{prop}\label{prop: versal formal functions}
With the notation and assumptions introduced above,
the functor~\eqref{eqn: versal restriction functor I} is fully faithful.
\end{prop}
\begin{proof}
We will show that both arrows of~\eqref{eqn: versal restriction functor I} are fully faithful.
For the first arrow, 
taking into account the isomorphism $\cF \isoto \quoteslim{n}\cF/\fm_R^{n+1}\cF$ arising from Lemma~\ref{lem:pro-coherent-as-limit-of-pullbacks},
it suffices to prove that $\ihat^*$ is fully faithful on the full subcategory of $\fm_R$-power torsion objects, i.e.\ that the restriction
\[
\Coh(\cX) \xrightarrow{\eqref{abelian coherent to OXhat}} \Coh(\cO_{\cX}) \xrightarrow{\ihat^*} \Coh(\cO_{\widehat \cX})
\]
is fully faithful.
This restriction can be written as the composite
\[
\Coh(\cX) \xrightarrow{i_{\fm_R, *}} \Coh(\fX) \xrightarrow{\cO_{\widehat \cX} \otimes_{\cO_{\fX}} \text{--}} \Coh(\cO_{\widehat \cX}), %
\]
where $i_{\fm_R}: \cX \to \fX$ is the $\fm_R$-adic completion morphism,
and so it is fully faithful, by Proposition~\ref{aprop:Coh-set-theoretically-supported} and 
Theorem~\ref{thm:summary-of-AHL}~\eqref{item:BGG-derived-formal-functions}.
This concludes the proof that the first arrow of~\eqref{eqn: versal restriction functor I} is fully faithful.

We now prove that the second arrow, i.e.\ the functor \eqref{eqn: versal restriction functor II}, is fully faithful.
We have already seen that it is faithful, as a consequence of Lemma~\ref{lem:versal pullback is exact and faithful II}.
We now have to show that it is also full.

To this end, let $v: \Spf \widehat{S} \to \widehat \cX$ be the  versal morphism given by the colimit of the flat covers $\Spec S/\fm_S^i \to [(\Spec S/\fm_S^i)/G_i]$.

Suppose that~$\cF$ and~$\cG$ are objects of $\Coh(\cO_{\widehat \cX})$, and that we are given a morphism~$\alpha$ between their restrictions to~$\Coh(\Ver/\widehat \cX)$; that is,   for each
$f: \Spf A \to \widehat \cX$ in~$(\Ver/\widehat \cX)$ we are given a morphism $\alpha_f : \cF(f) \to \cG(f)$, compatible with all base changes in $(\Ver/\widehat \cX)$.
If we 
write $\widehat{G}\coloneq G \times_{\Spec R} \Spf \widehat{S} = \colim_n G_n$ to
denote the $\fm_S$-adic completion of $G_S\coloneq G \times_{\Spec R} \Spec S$ (compare~\cite[Lem.\ 4.8]{Emertonformalstacks}),
then since $\cF(v)$ and $\cG(v)$ are pulled back from~$\cF$ and~$\cG$,
they are naturally~$\widehat{G}$-equivariant modules over~$\widehat{S}$,
and $\cF$ and $\cG$ are recovered by descending $\cF(v)$ and $\cG(v)$ with 
respect to this equivariant structure.
We will prove that the morphism
$\alpha_v$ is also equivariant for the action of~$\widehat{G}$; 
it then descends to the desired morphism from $\cF$ to~$\cG$.

The equivariance of each of $\cF(v)$ and $\cG(v)$ 
is encoded
by isomorphisms $\imath_{\cF(v)}: \widehat p^*\cF(v) \iso \widehat a^*\cF(v)$ and $\imath_{\cG(v)}: \widehat p^*\cG(v) \iso \widehat a^* \cG(v).$ 
Our goal then is to show that
\begin{equation}\label{eqn: to prove for equivariance}
\imath_{\cG(v)} \widehat p^*(\alpha_v) = \widehat a^*(\alpha_v) \imath_{\cF(v)}
\end{equation}
(as morphisms $\widehat p^*\cF(v) \to \widehat a^*\cG(v)$ of $\cO_{\widehat{G}}$-modules).

Since~$\widehat{G}$ is the $\fm_S$-adic completion of $G_S$, which has finite type over~$\Spec S$, 
it suffices to check that~\eqref{eqn: to prove for equivariance} holds %
after pulling back along every
morphism $f: \Spf A \to \widehat{G}$
arising as the completion of $G_S$ at a closed point lying above~$\fm_S$. %
Fixing such an~$f$,
the compositions
\[
\Spf A \xrightarrow{f} \widehat{G} \rightrightarrows \Spf \widehat{S} \xrightarrow{v} \widehat \cX
\]
coincide, and the resulting morphism $u: \Spf A \to \widehat \cX$ is an object of $(\Ver/\widehat \cX)$: this is a consequence of the fact that the projection map
$p: \widehat{G} \to \Spf \widehat{S}$ is a completion of the smooth morphism $G_S \to \Spec S$. %

The definition of morphisms in $\Coh(\Ver/\widehat \cX)$ gives the commutativity of the two squares in the following diagram, 
where the horizontal morphisms are induced by the sheaf property of~$\cF$ and  the fact that $v p = v a$
as morphisms $\widehat{G} \to \widehat \cX$,
\[ \begin{tikzcd}[row sep=huge, column sep=huge]
\widehat f^*\widehat p^*\mathcal{F}(v) \arrow[d, "\widehat f^*\widehat p^*(\alpha_v)"] \arrow[r, "\sim"'] \arrow[rr, bend left=20, "\widehat f^*\imath_{\cF(v)}"]
& \mathcal{F}(u) \arrow[d, "\alpha_u"] \arrow[r, "\sim"'] & \widehat f^*\widehat a^*\mathcal{F}(v) \arrow[d, "\widehat f^*\widehat a^*(\alpha_v)"']  \\
\widehat f^*\widehat p^*\mathcal{G}(v) \arrow[r, "\sim"'] \arrow[rr, bend right=20, "\widehat f^*\imath_{\cG(v)}"'] & \mathcal{G}(u) \arrow[r, "\sim"']
& \widehat f^*\widehat a^*\mathcal{G}(v)  
\end{tikzcd} \]
That \eqref{eqn: to prove for equivariance} holds after pulling back along~$f$ is now immediate from the commutativity of the diagram, and we are done.
\end{proof}

\section{Stacks of \texorpdfstring{$(\varphi,\Gamma)$}{(phi,Gamma)}-modules}\label{sec:stacks-appendix}In this
appendix we prove some minor variations on results
of~\cite{emertongeepicture}. We freely use the notation and conventions
of~\cite{emertongeepicture}. In particular, we fix a finite extension
$K/\Qp$, %
 and for each integer~$d\ge 1$ we write $\cX_d$ for the
corresponding moduli stack of rank~$d$ projective \'etale
$(\varphi,\Gamma)$-modules.
\subsection{Fixed determinant stacks}
\label{subsec:fixed-determ-stacks}We begin by extending some of the
results of~\cite{emertongeepicture} to the case of moduli stacks of
\'etale $(\varphi,\Gamma)$-modules with fixed determinant. %
 There is a natural morphism $\wedge^d:\cX_d\to\cX_1$, sending a rank~$d$
projective \'etale $(\varphi,\Gamma)$-module $D$ to~$\wedge^dD$. 
Let~$\chi:G_K\to\cO^\times$ be a character. %
Write~$\chibar:G_K\to\F^\times$ for the reduction
of~$\chi$ modulo$~\varpi$.
\begin{adefn}\label{defn: Xdchi}
  We define $\cX_d^{\chi}$ as the pullback %
  \[ \begin{tikzcd}
      \cX_d^{\chi}\ar{d}\ar{r}\arrow[dr, phantom, "\ulcorner", very near start]&\cX_d\ar{d}{\wedge^d}\\
      \Spf\cO\ar{r}{\chi}&\cX_1
    \end{tikzcd}\]
  More explicitly, $\cX_d^\chi(A)$ is the groupoid of pairs
  $(D,\theta)$ where~$D$ is a rank~$d$ projective \'etale
  $(\varphi,\Gamma)$-module with $A$-coefficients, and $\theta$ is an
  identification of~$\wedge^dD$ with~$\chi$.
\end{adefn}

If~$\chi$ is a de Rham character, $\underline{\lambda}$ is a
Hodge type, and~$\tau$ is an inertial type,
we define \begin{equation}\label{eqn:defn-of-crys-lambda-tau-chi}\cX_d^{\crys,\lambdau,\tau,\chi}\coloneqq \cX_d^{\crys,\lambdau, \tau}\times_{\cX_d}\cX_d^{\chi}. \end{equation}
When~$\tau$ is trivial, we will often omit it from the notation.
By
construction, $\cX_d^{\crys,\lambdau,\tau,\chi}$  is a closed substack
of~$\cX_d^\chi$, and by~\cite[Thm.\ 4.8.12]{emertongeepicture},
it is
a $p$-adic formal algebraic stack. It is non-zero if and only
if~$(\lambdau, \tau)$ is compatible with ~$\chi$ in the following sense. %

\begin{adefn}
  \label{adefn: Hodge type compatible with chi }We say that ~$(\lambdau, \tau)$
  is \emph{compatible with} ~$\chi$ if %
  $\chi\det(\tau)^{-1} $ is crystalline, and for
  each $\sigma:K\into\Qpbar$ we have $\HT_\sigma(\chi\det(\tau)^{-1})=\sum_{i=1}^d\lambda_{\sigma,i}$.
\end{adefn}

\begin{acor}\label{acor: fixed determinant crystalline
    stacks}If~$(\lambdau, \tau)$ is compatible
  with~$\chi$, and~$\lambdau$ is regular, then $\cX_d^{\crys,\lambdau,\chi}\times_{\Spf
    \cO}\Spec\F$ is equidimensional of dimension~$[K:\Qp]d(d-1)/2$. 
\end{acor}
\begin{proof}The proof of ~\cite[Thm.\
  4.8.14]{emertongeepicture} goes over immediately, replacing the
crystalline   deformation rings with their fixed determinant variants
(whose dimensions are known by the proof of~\cite[Thm.\
3.3.8]{MR2373358}; for a precise statement see for example~\cite[Thm.\
A]{zbMATH07042062}), and the reductive group~$\GL_d$ by~$\SL_d$.  
\end{proof}

\begin{adefn} %
  \label{adefn: compatible Serre weight}We say that a Serre
  weight~$\underline{k}$  \emph{is compatible with~$\chi$} (or  \emph{is compatible with~$\chibar$})
  if
  \[\chibar|_{I_K}=\varepsilonbar^{-d(d-1)/2}\prod_{\sigmabar:k\into\Fpbar}\omega_{\sigmabar}^{-\sum_{i=1}^dk_{\sigmabar,i}}.\]
\end{adefn}

\begin{rem}
  \label{rem: two notions of compatible are compatible}The conventions
  of this paper for Serre weights are opposite to those
  of~\cite{emertongeepicture}. This has the slightly unfortunate
  consequence that in the
  particular case that~$K=\Qp$ and~$d=2$, Definition~\ref{adefn:
    compatible Serre weight} compares with Definition~\ref{defn:
    compatible Serre weight} as follows: a Serre weight $\sigma_{a,b}$ is
  compatible with $\zetabar$ in the sense of Definition~\ref{defn:
    compatible Serre weight} if and only if it is compatible with
  $\zetabar^{-1}\varepsilonbar^{-1}$ in the sense of Definition~\ref{adefn:
    compatible Serre weight}. %
(To
see this, take %
take
  $k_1=a+b$, $k_2=a$ (where we have suppressed the unique choice
  of~$\sigmabar$ from the notation), and note that
  $\omega=\varepsilonbar$.) 
\end{rem}

\begin{adefn}
  \label{adefn: Steinberg Serre weight}We say that a Serre
  weight~$\underline{k}$ is \emph{Steinberg} if for all
  $\sigmabar:k\into\Fpbar$ and all $1\le i\le d-1$ we have
  $k_{\sigmabar,i}-k_{\sigmabar,i+1}=p-1$. 
\end{adefn}

The following is the analogue for~$\cX_d^\chi$ of some of the results of~\cite{emertongeepicture} for~$\cX_d$.

\begin{athm}
  \label{athm: fixed determinant stack in general}$\cX_d^\chi$ is a Noetherian
  formal algebraic stack. Its underlying reduced substack $\cX_{d,\red}^\chi$
  is of finite type over~$\Fp$, and is equidimensional of
  dimension~$[K:\Qp]d(d-1)/2$. The irreducible components of~$\cX^{\chi}_{d,\red}$ admit a
  natural bijection with the Serre weights which are compatible
  with~$\chi$, except that each Steinberg weight which is
  compatible with~$\chi$ corresponds to multiple irreducible
  components, indexed by the~$d$th roots of unity in~$\Fpbar$. In each
  case, a component corresponding to a Serre weight $\underline{k}$
  has a dense open substack which is maximally non-split of niveau one
  and weight~$\underline{k}$.  %
\end{athm}
\begin{proof}The argument of
  the proof of~\cite[Thm.\ 6.5.1]{emertongeepicture} goes through
  almost unchanged to show that every irreducible component
  of~$\cX_{d,\red}^\chi$ has dimension at least~$[K:\Qp]d(d-1)/2$, %
  replacing the appeal to~\cite[Thm.\
  4.8.14]{emertongeepicture} with one to Corollary~\ref{acor: fixed
    determinant crystalline stacks}. (The one possibly subtle point is
  that~ $\cX_{d,\red}^\chi$ only depends on~$\chibar$, so in making
  the comparison to $\cX_d^{\crys,\lambdau,\chi}$ we are free to
  replace~$\chi$ by a crystalline character which is compatible with~$\lambdau$.)

  We now
  consider the morphism \begin{equation}\label{eqn: twisting from fixed
    determinant to unfixed}f:\cX_{d,\red}^\chi\times\Gm\to
  \cX_{d,\red}\end{equation} given by forgetting~$\theta$ and taking
  unramified twists, as in~\cite[\S 5.3]{emertongeepicture}. Note that
  any family of $(\varphi,\Gamma)$-modules with determinant~$\chi$ is
  essentially twistable in the sense of~\cite[Defn.\
  5.3.1]{emertongeepicture} (because any twist preserving the
  determinant would have to be by a $d$th root of unity). 

  Let~$\cZ$ be an irreducible component of~$\cX_{d,\red}^\chi$, and let~$f(\cZ)$ be
  the scheme-theoretic image of~$\cZ\times\Gm$ in
  $\cX_{d,\red}$.  Since we have already shown
  that~$\cZ$ has dimension at least~$[K:\Qp]d(d-1)/2$, it follows from~\cite[Lem.\ 5.3.2]{emertongeepicture}, that~$f(\cZ)$ has dimension at least
  $[K:\Qp]d(d-1)/2$. (Use
  ~\cite[\href{https://stacks.math.columbia.edu/tag/0DS4}{Tag
    0DS4}]{stacks-project}, together with the fact that the forgetful
  morphism $\cX_{d,\red}^\chi\to\cX_{d,\red}$ has fibres of dimension
  (at most)~$1$.) %
Since~$\cX_{d,\red}$ is
  equidimensional of this dimension, we see that~$\cZ$ must have
  dimension exactly $[K:\Qp]d(d-1)/2$, and that~$f(\cZ)$ is
  an irreducible component of~$\cX_{d,\red}$, so
  that~$f(\cZ)=\cX_{d,\red}^{\underline{k}}$ for some Serre
  weight~$\underline{k}$. Since~$\cZ$ was arbitrary, we have in
  particular shown that~$\cX_d^\chi$ is equidimensional of
  dimension~$[K:\Qp]d(d-1)/2$, as claimed.

It follows from Definition~\ref{adefn: compatible Serre weight}
and~\cite[Defn.\ 5.5.1]{emertongeepicture} (together with the defining
property of $\cX_{d,\red}^{\underline{k}}$) that~$\underline{k}$ is
compatible with~$\chi$. It therefore remains to show that the
association of~$\underline{k}$ to~$\cZ$ is
surjective, and that it is injective if~$\underline{k}$ is not
Steinberg, and to determine the Steinberg components. %

Let $\cX_d^{\underline{k},\chi}$ be  the closed substack of
~$\cX_{d,\red,\Fpbar}^{\chi}$ given by the fibre product
\[ \begin{tikzcd}
      \cX_d^{\underline{k},\chi}\ar{d}\ar{r}\arrow[dr, phantom, "\ulcorner", very near start]&\cX_{d,\red,\Fpbar}^{\underline{k}}\ar{d}{\wedge^d}\\
  \Spec\Fpbar\ar{r}{\chi}&\cX_1
\end{tikzcd}\] Then  $\cX_d^{\underline{k},\chi}$ is nonempty and   has dimension~$[K:\Qp]d(d-1)/2$ (for example by
another appeal to ~\cite[\href{https://stacks.math.columbia.edu/tag/0DS4}{Tag
    0DS4}]{stacks-project}), so %
  there is at least
  one irreducible component of~$\cX_{d,\red}^\chi$ corresponding
  to~$\underline{k}$. Note also that if
  $\cX_{d,\red,\Fpbar}^{\underline{k}}$ is generically maximally
  non-split of weight~$\underline{k}$, so is $\cX_d^{\underline{k},\chi}$. 

  If~$\underline{k}$ is Steinberg then the classification of the
  irreducible components follows easily from \cite[Thm.\
  5.5.12]{emertongeepicture} and the definition of being maximally
  non-split of weight $\underline{k}$, so we assume from now on that
  $\underline{k}$ is not Steinberg. We need to show that $\cX_d^{\underline{k},\chi}$ is
  irreducible. We prove this by induction on~$d$. The inductive argument is similar to the proof of~\cite[Thm.\
  5.5.12]{emertongeepicture} given in~\cite{EGaddenda}, but much simpler,
  as we are already using the conclusions of~\cite[Thm.\
  5.5.12]{emertongeepicture}. We content ourselves with explaining the
  key differences in the argument. The case~$d=1$ is trivial, and the
  case $d=2$ is easy and is left to the reader. If~$d\ge 3$ then after possibly replacing all representations with their duals (which
has the effect of reversing the order of the~$k_i$) we can and do
assume that $\underline{k}_{d-1}$ is also not Steinberg.

  Let
$\cU^{\underline{k},\chi}$ be a dense open substack of
$\cX_d^{\underline{k},\chi}$ which is maximally
non-split of weight~$\underline{k}$, and let
$\nu_1:\cU^{\underline{k},\chi}\to\Gm$ be the character given
by~\cite[(5.5.10)]{emertongeepicture} (which by definition determines
the unique quotient character of each $\Fpbar$-point of
$\cU^{\underline{k},\chi}$; this character is  the unramified twist
  of $\alphabar\coloneqq \varepsilonbar^{1-d}\omega_{\underline{k},1}$
  by~$\nu_1$). %
Let~$\cZ$ be
an irreducible component of $\cX_d^{\underline{k},\chi}$, and let~$\cU$  be its 
intersection with $\cU^{\underline{k},\chi}$.  Since~$\cU$ has dimension~$[K:\Qp]d(d-1)/2$, the locus in~$\Gm$ where
the fibres of~$\nu_1$ have dimension at least $[K:\Qp]d(d-1)/2-1$ is
nonempty, and since it is Zariski closed it is either finite or all
of~$\Gm$.

If this locus were finite, then since the dimension
of~$\cU^{\underline{k},\chi}$ is~$[K:\Qp]d(d-1)/2$, there would be
some fibre of dimension~$[K:\Qp]d(d-1)/2$.  By the usual computations
of the dimensions of extension groups, this would necessitate some~$\cX_{d-1}^{\chi'}$ having dimension at least
$[K:\Qp](d-1)(d-2)/2+1$, a contradiction. It follows that for each
$\Fpbar$-point of~$\Gm$, the fibre of~$\nu_1$ has dimension
$[K:\Qp]d(d-1)/2-1$. Furthermore by the inductive hypothesis, we see
that such a fibre has precisely one irreducible component
of dimension $[K:\Qp]d(d-1)/2-1$, and this irreducible component
contains a dense open substack of the universal extension of the fixed
twist of~$\alphabar$ by a dense open substack of the (irreducible, by
the inductive hypothesis) $\cX_{d-1}^{\underline{k}_{d-1},\chi'}$ (where
$\chi'$ is determined by~$\chi$ and the twist of~$\alphabar$
determined by the point of~$\Gm$). Since this analysis was independent
of the choice of irreducible component~$\cZ$, we see that~$\cZ$ is
unique, as required.
\end{proof}

\subsection{Stacks of Galois representations}\label{subsec:stacks-of-Galois-reps-appendix}
As in~\cite[\S 6.7]{emertongeepicture}, we write~$\cX_d^{\Gal}$ for Wang-Erickson's formal algebraic
stack (see~\cite[Thm.\ 3.8, Rem.\ 3.9]{MR3831282}), which is
characterised by the following property: if~$A$ is an
$\cO$-algebra in which $p$ is nilpotent, then 
 $\cX_{d}^{\Gal}(A)$ is the groupoid
of rank~$d$ projective $A$-modules~$T_A$ with a continuous action of~$G_K$
(where ~$T_A$ has the discrete topology). 
The universal object~$\cV_d$ on~$\cX_d^{\Gal}$ is an object of $\ProCoh (\cX_d^{\Gal})$.
By~\cite[Thm.\
6.7.2]{emertongeepicture}, there is a monomorphism 
(i.e.\ a fully faithful functor, compatible with the structure of fibred categories)
\begin{equation}
\label{eqn:Galois stacks to ours}
\cX_d^{\Gal}\to\cX_d. 
\end{equation}
In order to describe this monomorphism, it is convenient to recall that
$\cX_d$ can be regarded as classifying rank $d$ projective \'etale $(\varphi,G_K)$-modules
with $A$-coefficients, in the sense of \cite[Defn.\ 2.7.7]{emertongeepicture};
these objects are equivalent to rank~$d$
projective \'etale $(\varphi,\Gamma)$-modules with $A$-coefficients,
by~\cite[Prop.\ 2.7.8]{emertongeepicture}.
We adopt this perspective for the duration of this discussion.

The monomorphism~\eqref{eqn:Galois stacks to ours}
is then given by the following construction: if
$T_A\in\cX_d^{\Gal}(A)$, where~$A$ is a finite type $\cO/\varpi^a$-algebra  
then
\begin{equation}\label{eqn:phi-GK-from-Galois}
T_A \mapsto
\mathbb{D}_A(T_A)\coloneqq T_A\otimes_A W(\C^\flat)_A = T_A\otimes_A W_a(\C^\flat)_A\end{equation}
(so that $\mathbb{D}_A(T_A)$ becomes a rank~$d$ projective \'etale $(\varphi,G_K)$-module with
$A$-coefficients when equipped with the $\varphi$-action induced by the $\varphi$-action
on~$W(\C^\flat)_A$, together with the diagonal~$G_K$-action).
Recall that, by definition, 
\begin{gather*}
W(\Cflat)_A \coloneqq  \varprojlim\nolimits_i W_i(\Cflat)_A,\\
W_i(\Cflat)_A \coloneqq  W_i(\cO_\C^\flat)_A[1/v],\\
W_i(\cO_\C^\flat)_A = W_i(\cO_\C^\flat) \cotimes_{\bZ_p} A \coloneqq  (W_i(\cO_{\C}^\flat) \otimes_{\bZ_p} A)^\wedge_{\text{$v$-adic}},
\end{gather*}
where~$v \in W_i(\cO_{\C}^\flat)$ is any element of the maximal ideal of $W_i(\cO_\C^\flat)$ whose image in $\cO_\C^\flat$ is not zero.
Furthermore, the Frobenius $\varphi : W_i(\cO_\C^\flat) \to W_i(\cO_\C^\flat)$ is $v$-adically continuous,
and so extends uniquely to an $A$-linear continuous endomorphism of $W_i(\C^\flat)_A$, still denoted~$\varphi$.

\begin{alem}\label{lem: properties of phi}
  If~$A$ is a finite type $\cO/\varpi^a$-algebra, then the sequence
  \[
  0 \to A \to W_a(\Cflat)_A \xrightarrow{\varphi-1} W_a(\Cflat)_A \to 0
  \]
  is exact.
\end{alem}
\begin{proof}
The claim that $W_a(\Cflat)_A^{\varphi = 1} = A$ is~\cite[Lem\ 2.2.19]{emertongeepicture}.
There remains to prove that $\varphi - 1$ is surjective on $W_a(\Cflat)_A$. 
Since $p$ is nilpotent on $W_a(\C^\flat)_A$, it suffices to prove that $\varphi-1$ is surjective on $\Cflat_A=W_a(\Cflat)_A/p$.

The action of $\varphi-1$ on $\Cflat$ is surjective, because~$\Cflat$ is algebraically closed.
Thus it also acts surjectively on $\Cflat\otimes_{\Zp}A$.
Let~$z\in\Cflat_A$ be arbitrary.
Since $\Cflat_A=(\cO_{\Cflat}\cotimes_{\Zp}A)[1/v]$, there exists~$z'\in \Cflat\otimes_{\Zp}A$ with $y\coloneqq z-z'\in v(\cO_{\Cflat}\cotimes_{\Zp}A)$.
Then~$y=(\varphi-1)(-\sum_{n\ge 0}\varphi^n(y))$ is in the image of~$\varphi-1$, as is~$z'$ (since it is an element of $\Cflat\otimes_{\Zp}A)$, 
and hence so is $z=y+z'$, as required.\qedhere

\end{proof}

\begin{alem}
  \label{lem:Galois-stack-characterise-by-Frobenius-invariants}If~$A$
  is a finite type $\cO/\varpi^a$-algebra, and $M_A\in\cX_d(A)$ is an
  \'etale $(\varphi,G_K)$-module, then $M_A$ is in the essential image 
  of ~\eqref{eqn:Galois stacks to ours} if and only if the natural map \begin{equation}\label{eqn:phi-invariants-condition}M_A^{\varphi=1}\otimes_{A}W(\C^\flat)_A\to M_A\end{equation}  is an isomorphism. 
\end{alem}
\begin{proof}%

    If~$M_A$ is in the essential image
  of ~\eqref{eqn:Galois stacks to ours}, then~\eqref{eqn:phi-invariants-condition} is an isomorphism, by a consideration of~\eqref{eqn:phi-GK-from-Galois} and the fact that
  $W(\C^\flat)_A^{\varphi=1}=A$ (see Lemma~\ref{lem: properties of phi}).
  Conversely, if
  ~\eqref{eqn:phi-invariants-condition} is an isomorphism, it suffices
  to show that the $G_K$-module $M_A^{\varphi=1 }$ is a projective
  $A$-module of
  rank~$d$. By~\cite[\href{https://stacks.math.columbia.edu/tag/058S}{Tag
  058S}]{stacks-project}, it suffices to show that $W(\C^\flat)_A$ is a
faithfully flat $A$-algebra. %
By~\cite[Prop.\
  2.2.12]{emertongeepicture}, it suffices in turn to show that the Laurent
  series ring $A((T))$ is a faithfully flat $A$-algebra. Since $A$ is
  Noetherian, and $A((T))$ is a localization of a completion of the
  flat $A$-algebra $A[T]$, we see that $A((T))$ is a flat
  $A$-algebra. To see that it is faithfully flat, we can use the
  morphism of $A$-modules $A((T))\to A$ given by $\sum a_iT^i\mapsto a_0$. %
\end{proof}

\begin{lem}
  \label{lem:Galois-substack-nilpotent-lift}
  Let $\Spec A\to\cX_d$ be a morphism whose source is a finite type $\cO/\varpi^a$-algebra for some~$a\ge 1$, and let~$I$ be a nilpotent ideal in~$A$.
Suppose that the composite $\Spec A/I\to\Spec A\to \cX_d$ factors through~$\cX_d^{\Gal}$.
Then the morphism $\Spec A\to \cX_d$ itself factors through~$\cX_d^{\Gal}$.
\end{lem}
\begin{proof}
  Since~$A$ is Noetherian, we immediately reduce to the case where $I$ is generated by a single element~$\varepsilon$ such that $\varepsilon^n = 0$ for some~$n \ge 1$.
  We prove the lemma by induction on~$n$ (the case~$n=1$ being trivial).

  Let~$M_A$ be the finite projective \'etale  $(\varphi,G_K)$-module over~$W(\C^{\flat})_{A}$ corresponding as above to the given morphism $\Spec A\to \cX_d$.  By Lemma~\ref{lem:Galois-stack-characterise-by-Frobenius-invariants}, we need to show that the natural map~\eqref{eqn:phi-invariants-condition} is an isomorphism.  %

We now consider the exact sequence
\begin{equation}\label{eqn:exact sequence phi GK}
0 \to \varepsilon M_A \to M_A \to M_A/\varepsilon M_A \to 0.
\end{equation}
Let~$J \coloneqq  \Ann_A(\varepsilon)$, so that $\varepsilon^{n-1} \in J$.
Since $W(\C^\flat)_A$ is flat over~$A$ (as we noted in the proof of Lemma~\ref{lem:Galois-stack-characterise-by-Frobenius-invariants}), 
and $M_A$ is finite projective over $W(\C^\flat)_A$, we see that the natural map $ M_A/J\to \varepsilon M_A$ is an isomorphism.
By our inductive assumption, the morphism $\Spec A/\varepsilon^{n-1} A \to \cX_d$ classifying $M_A/\varepsilon^{n-1}M_A$ factors through $\cX_d^{\Gal}$,
hence the same is true for the morphism $\Spec A/JA \to \Spec A/\varepsilon^{n-1} A \to \cX_d$ classifying $M_A/J\isoto \varepsilon M_A$.

It follows from Lemma~\ref{lem:Galois-stack-characterise-by-Frobenius-invariants} that the natural map 
\[(\varepsilon M_A)^{\varphi=1}\otimes_{A}W(\C^\flat)_A\to \varepsilon M_{A}\]is an isomorphism. 
In particular, we see that $\varphi-1$ is surjective on $\varepsilon M_{A}$, since it is surjective on $W(\C^\flat)_A$, by Lemma~\ref{lem: properties of phi}. 
Applying the snake lemma to the endomorphism $\varphi-1$ of~\eqref{eqn:exact sequence phi GK},
we deduce that
\[
0 \to (\varepsilon M_A)^{\varphi=1} \to (M_A)^{\varphi=1} \to (M_A/\varepsilon M_A)^{\varphi=1} \to 0
\]
is exact.
Since $W(\C^\flat)_A$ is flat over~$A$,  we obtain a commutative diagram with exact rows
\[\begin{adjustbox}{max width=\textwidth} \begin{tikzcd}[column sep=small]0 \arrow[r] & (\varepsilon M_A)^{\varphi=1} \otimes_A W(\C^\flat)_A \arrow[r] \arrow[d%
] & 
(M_A)^{\varphi=1} \otimes_A W(\C^\flat)_A \arrow[r] \arrow[d%
] & 
(M_A/\varepsilon M_A)^{\varphi=1} \otimes_A W(\C^\flat)_A \arrow[r] \arrow[d%
] & 0\\
0 \arrow[r] & \varepsilon M_A \arrow[r] & M_A \arrow[r] & M_A/\varepsilon M_A \arrow[r] & 0 \end{tikzcd} \end{adjustbox}
\]
Since the left and right hand vertical arrows are isomorphisms (by the inductive hypothesis), the five lemma implies that so is the middle vertical arrow, as required.
\end{proof}

If~$L/K$ is a finite extension, then we have maps (with
obvious notation)
$\cX_{K,d}^{\Gal}\to \cX_{L,d}^{\Gal}$ and $\cX_{K,d}\to \cX_{L,d}$,
in each case given by restricting the action of~$G_K$ to its
subgroup~$G_L$.
The second morphism is analyzed in detail in~\cite[Lem.~3.7.5]{emertongeepicture}, 
where it is shown to be representable by algebraic spaces, affine, and of finite presentation.

\begin{alem}
  \label{lem:Galois-stack-Cartesian-diagram-extension}If~$L/K$ is a
  finite extension, then the following
  diagram is Cartesian. 
\[\begin{tikzcd}
	\cX_{K,d}^{\Gal} &  \cX_{K,d} \\
\cX_{L,d}^{\Gal}	 &\cX_{L,d}
	\arrow[from=1-1, to=1-2]
	\arrow[from=1-1, to=2-1]
	\arrow[from=1-2, to=2-2]
	\arrow[from=2-1, to=2-2]
\end{tikzcd}\]
\end{alem}
\begin{proof}
  This is immediate from
  Lemma~\ref{lem:Galois-stack-characterise-by-Frobenius-invariants}
  (since the criterion of that lemma does not involve the $G_K$-action or $G_L$-action at all!).
\end{proof}

\subsubsection{Stacks with fixed pseudorepresentation}\label{subsubsec:stacks with fixed pseudorepresentation}
Let~$A$ be an $\cO$-algebra in which $p$ is nilpotent, and assume $\Spec A$ is connected.
If $T_A \in \cX_d^{\Gal}(A)$, then the $\cO$-subalgebra $B \subseteq A$ generated by the values of the $A$-valued pseudorepresentation~$\theta_A$
associated to~$T_A$ is Artinian and local with finite residue field, 
by~\cite[Lem.\ 3.10, Defn.\ 3.12]{MR3444227}.
If~$\thetabar$ is a $d$-dimensional $\cbF_p$-pseudorepresentation of~$G_K$,
we define $\cX_{\thetabar}(A) \subset \cX_d^{\Gal}(A)$ to be the
full subgroupoid of objects~$T_A$ such that
$\theta_A \otimes_B \cbF_p$ %
is $\Gal(\cbF_p/\F)$-conjugate to~$\thetabar$. 
It follows from \cite[Cor.\ 3.14]{MR3444227} (see also \cite[Thm.\ 3.5]{MR3831282} for a restatement in terms of the stacks~$\cX_d^{\Gal}$)
that
\begin{equation}\label{eqn: Galois breaks up
    over
    pseudorepns}\cX_d^{\Gal}=\coprod_{\thetabar}\cX_{\thetabar},\end{equation}
where~ $\thetabar$ runs over representatives of the $\Gal(\cbF_p/\F)$-conjugacy classes of $d$-dimensional
$\Fpbar$-pseudorepresentations of~$G_K$.

We write
$\F_{\thetabar}$ for the field extension of~$\F$ generated by the values of~$\thetabar$. 
Then~$\thetabar$ can be regarded as a $\F_{\thetabar}$-valued pseudorepresentation.
The universal (pseudo)deformation ring of~$\thetabar$, to complete Noetherian local $W(\F_{\thetabar}) \otimes_{W(\F)} \cO$-algebras, will be denoted
$R_{\thetabar}^{\ps}$.
Note that there is a natural morphism 
\begin{equation}\label{eqn:CWE-stack-to-Spf-Rpseudo}
  \cX_{\thetabar}\to\Spf R^{\ps}_{\thetabar},
\end{equation}
since for any test morphism $\Spf A \to \cX_{\thetabar}$ (with $A$, as above,
an $\cO$-algebra in which $p$ is nilpotent), the inclusion
$B\subseteq A$ (where, as above, $B$ denotes the sub-$\cO$-algebra
of traces) induces the composite morphism  
$R_{\thetabar}^{\ps} \to B \subset A$,
giving rise to a morphism $\Spf A \to \Spf R_{\thetabar}^{\ps}$.
Up to natural isomorphism, the ring~$R_{\thetabar}^{\ps}$ and the map~\eqref{eqn:CWE-stack-to-Spf-Rpseudo} only depend on the $\Gal(\cbF_p/\F)$-conjugacy class of~$\thetabar$.

The following theorem confirms the expectation of~\cite[Rem.\ 6.7.4]{emertongeepicture}.
\begin{athm}\label{thm: thetabar substacks app version}
The morphism \[\cX_{\thetabar}\to\cX_d\] induced by~\eqref{eqn: Galois breaks up
    over
    pseudorepns} and~\eqref{eqn:Galois stacks to ours}
is a monomorphism which induces a closed immersion on underlying reduced substacks, and exhibits~$\cX_{\thetabar}$  as the completion of~$\cX_{d}$ along its closed substack~$\cX_{\thetabar,\red}$. %
\end{athm}
\begin{proof}
That the morphism is a  monomorphism  is immediate from
\cite[Thm.~6.7.2]{emertongeepicture}. To see that it induces a closed
immersion on underlying reduced substacks, we begin by noting that there is a finite extension~$L/K$ such that every
representation $\rhobar:G_K\to\GL_d(\Fpbar)$ corresponding to an
$\Fpbar$-point of $\cX_{\thetabar}$ factors
through~$\Gal(L/K)$. Indeed, this is a basic fact that underlies
the entire framework of~\cite{MR3831282}.  
To see it directly, note that we may firstly replace~$K$ by the
fixed field of~$\ker\thetabar$, and thus assume that $\thetabar$ is
trivial, so that~$\rhobar(G_K)$ is unipotent. We can then repeatedly
(a total of~$(d-1)$ times) replace~$K$ by its maximal elementary
abelian $p$-extension.

By~\cite[Thm.\ 6.6.3(2)]{emertongeepicture}, the trivial $G_L$-representation 
gives rise to a closed residual gerbe $Z \hookrightarrow \cX_{L,d}$, and so also of $\cX_{L,d}^{\Gal}$.
As recalled above, the right hand vertical 
arrow in the Cartesian diagram of
Lemma~\ref{lem:Galois-stack-Cartesian-diagram-extension} is representable by algebraic
spaces, and so the fibre product $Z \times_{\cX_{L,d}} \cX_{K,d}$
is a closed algebraic substack of $\cX_{K,d}$ contained in $\cX_{K,d}^{\Gal}$ which
(by virtue of our choice of~$L$)
contains~$\cX_{\thetabar,\red}$ as a closed substack (indeed, as a connected component).
Thus the induced morphism on underlying reduced substacks is indeed a closed immersion.

Finally, write $\widehat \cX_d \to \cX_d$ for the 
completion of~$\cX_d$ along $\cX_{\thetabar, \red}$.
Then the morphism
$\cX_{\thetabar}\to\cX_d$ 
factors through a monomorphism $\cX_{\thetabar}\to\widehat \cX_d$, which we need to prove is an equivalence.
Since it is a monomorphism, it is fully faithful, and the essential surjectivity is a consequence of
Lemma~\ref{lem:Galois-substack-nilpotent-lift}. 
\end{proof}

We will write~$\fX_{\thetabar}$ for the algebraic stack
denoted~$\Rep_{\lbar D}$ in~\cite{MR3831282} (with $\lbar D \coloneqq 
\thetabar$); this is of finite type over $\Spec
R_{\thetabar}^{\ps}$, and  is in fact (by definition) a quotient stack \[\Rep_{\lbar
    D}=[\Rep_{\lbar D}^{\square}/\GL_d],\]where $\Rep_{\lbar
  D}^{\square}$ is a finite type scheme over $\Spec
R_{\thetabar}^{\ps}$.

\begin{athm}\label{athm:CWE-main-result}
  The morphism~\eqref{eqn:CWE-stack-to-Spf-Rpseudo} is representable
  by algebraic stacks, and it arises as the~$\m_{R_{\thetabar}^{\ps}}$-adic
  completion of $\fX_{\thetabar}\to \Spec R_{\thetabar}^{\ps}$.    
\end{athm}%
\begin{proof}This is immediate from ~\cite[Thm.\
  3.8]{MR3831282}. 
\end{proof}

\begin{cor}\label{cor:CWE stacks satisfy assumption for formal functions}
Let 
$S \coloneqq  \Gamma(\Rep_{\lbar D}^{\square},\cO_{\Rep_{\lbar D}^{\square}})$
and $S_i \coloneqq  S\otimes_{ R^{\ps}_{\thetabar}} R^{\ps}_{\thetabar}/(\m_{R_{\thetabar}^{\ps}})^i.$
Then $S^{\GL_d}$
is a finite local
$R_{\thetabar}^{\ps}$-algebra \emph{(}so in particular a complete local
Noetherian ring\emph{)}, and each
$S_i^{\GL_d}$ is a finite local $R^{\ps}_{\thetabar}/(\m_{R_{\thetabar}^{\ps}})^i$-algebra.
\end{cor}
\begin{proof}
By~\cite[Thm.\ 3.8]{MR3831282}, the map $R^{\ps}_{\thetabar} \to S^{\GL_d}$ is an adequate homeomorphism (i.e.\ it is an integral universal homeomorphism, which is a local isomorphism at all points of residue characteristic zero; see~\cite[Defn.\ 3.3.1]{MR3272912}). 
By~\cite[Rem.\ 5.2.2]{MR3272912}, the map 
\[
S^{\GL_d} \otimes_{R^{\ps}_{\thetabar}} R^{\ps}_{\thetabar}/(\m_{R_{\thetabar}^{\ps}})^i \to (S \otimes_{R^{\ps}_{\thetabar}} R^{\ps}_{\thetabar}/(\m_{R_{\thetabar}^{\ps}})^i)^{\GL_d} = S_i^{\GL_d}
\]
is also an adequate homeomorphism.
The property of being an adequate homeomorphism is stable under base change and composition, so that $R^{\ps}_{\thetabar}/(\m_{R_{\thetabar}^{\ps}})^i \to S_i^{\GL_d}$ is an adequate homeomorphism,
and by Theorem~\ref{thm:summary-of-AHL}~\eqref{item:finite-type-invariants}, it is furthermore of finite type.

It thus suffices to prove that if $f: R \to R'$ is an adequate homeomorphism of finite type between $\cO$-algebras, and~$R$ is local, then~$R'$ is local 
and $f$ is finite (and therefore local).
That~$R'$ is local is immediate from $R\to R'$ being a homeomorphism, and that~$f$ is finite is immediate from it being integral and of finite type, as desired.
\end{proof}

We now specialize some of the framework of Appendix~\ref{sheaves on formal algebraic stacks}
to the context at hand.
In all of the following discussion, 
we freely use the equivalence
$$D^b\bigl(\Coh(\cZ)\bigr) \iso D^b_{\coh}(\cZ)$$
for any Noetherian algebraic, or formal algebraic, stack~$\cZ$ having an affine diagonal (see~\eqref{eqn:coherent equivalence}).

We write $\imm_{\thetabar} : \cX_{\thetabar} \hookrightarrow \fX_{\thetabar}$ for the completion map, and we write $\fX_{0}$ for the underlying reduced substack of~$\cX_{\thetabar}$. 
Recall from
Appendix~\ref{subsubsec:completed-pullback-to-completed} (in particular Proposition~\ref{prop:exact-adjoint-ihat-upper-star})
that we have a $t$-exact fully faithful functor
$$\imm_{\thetabar,*}: D^b_{\coh}(\cX_{\thetabar}) \to D^b_{\coh}(\fX_{\thetabar}),$$
whose Pro-extension 
$$\widehat{\imm}_{\thetabar,*}: \Pro D^b_{\coh}(\cX_{\thetabar}) \to
\Pro D^b_{\coh}(\fX_{\thetabar}),$$
which is again $t$-exact, admits a $t$-exact left adjoint
of ``$\fm$-adic completion''
$$\widehat{\imm}_{\thetabar}^*: \Pro D^b_{\coh}(\fX_{\thetabar}) \to
\Pro D^b_{\coh}(\cX_{\thetabar}).$$
We use the same notation to denote the various functors on hearts that are induced
by these $t$-exact functors; %
Corollary~\ref{cor:t-exact derived functors} shows that in turn these various $t$-exact functors are determined by 
the induced functors on hearts.

The canonical morphism $\fX_{\thetabar} \to \Spec R^{\ps}_{\thetabar}$ endows
any object $\fF$ of $\Coh(\fX_{\thetabar})$ with an $R^{\ps}_{\thetabar}$-module structure.
This structure is natural in $\fF$, i.e.\ $\Coh(\fX_{\thetabar})$ is enriched
over $R^{\ps}_{\thetabar}$-modules, and so any object of $\Pro \Coh(\fX_{\thetabar})$ 
is again endowed with an $R^{\ps}_{\thetabar}$-module structure.
Similarly, the adic morphism $\cX_{\thetabar} \to \Spf R^{\ps}_{\thetabar}$ induces
a natural $R^{\ps}_{\thetabar}$-module structure on any object~$\cF$ of $\Coh(\cX_{\thetabar})$,
or, more generally, on any object of $\Pro \Coh(\cX_{\thetabar})$. 
The functors 
$\imm_{\thetabar,*}$, $\widehat{\imm}_{\thetabar,*}$,
and $\widehat{\imm}_{\thetabar}^*$ are all compatible with these $R^{\ps}_{\thetabar}$-module
structures, i.e.\ they are functors of categories enriched over $R^{\ps}_{\thetabar}$-modules.

We summarize the situation in the following result.

\begin{athm}\label{thm:coherent-completeness-CWE-stacks}%
The functor
\[\imm_{\thetabar,*}: D^b_{\coh}(\cX_{\thetabar}) \hookrightarrow D^b_{\coh}(\mf{X}_{\thetabar}),\]
is $t$-exact and fully faithful, with essential image equal to
$D^b_{\coh,\mf{X}_0}(\mf{X}_{\thetabar}),$
i.e.\ the full sub-$\infty$-category of $D^b_{\coh}(\mf{X}_{\thetabar})$ consisting of
objects whose cohomology sheaves are set-theoretically supported on~$\mf{X}_0$. %

The Pro-extended functor
\[\widehat{\imm}_{\thetabar,*}: \Pro D^b_{\coh}(\cX_{\thetabar}) \hookrightarrow
\Pro D^b_{\coh}(\mf{X}_{\thetabar})\]
is again $t$-exact and fully faithful, as is the restriction
of its left adjoint
\[\widehat{\imm}_{\thetabar}^*: D^b_{\coh}(\mf{X}_{\thetabar})\to \Pro
    D^b_{\coh}(\cX_{\thetabar}). \] %
The restriction of $\widehat{\imm}_{\thetabar}^*$ to $D^b_{\coh,\fX_0}(\fX_{\thetabar})$
is an inverse equivalence to the equivalence $D^b_{\coh}(\cX_{\thetabar})
\iso D^b_{\coh,\fX_{0}}(\fX_{\thetabar})$ induced by~$\imm_{\thetabar,*}$.

\end{athm}%
\begin{proof}%
  The claim of the
first paragraph is a special case of Proposition~\ref{aprop:Coh-set-theoretically-supported}.
That $\widehat{\imm}_{\thetabar,*}$ is $t$-exact and fully faithful follows immediately, as the Pro-extension of a $t$-exact and fully faithful functor is again $t$-exact and fully
faithful.
Remark~\ref{rem:i-hat-upper-star as quasi-inverse} shows that  $\widehat{\imm}_{\thetabar}^*$ restricts to an equivalence 
$ D^b_{\coh,\fX_{0}}(\fX_{\thetabar})\iso D^b_{\coh}(\cX_{\thetabar})$.

It remains to show that $\widehat{\imm}_{\thetabar}^*:D^b_{\coh}(\mf{X}_{\thetabar})\to \Pro
    D^b_{\coh}(\cX_{\thetabar})$ is fully faithful.
To this end, as in Corollary~\ref{cor:CWE stacks satisfy assumption for formal functions}  we
write $\Rep_{\lbar D}^{\square}=\Spec S$, and then
apply Theorem~\ref{thm:summary-of-AHL}~\eqref{item:BGG-derived-formal-functions},
taking~$R=R_{\thetabar}^{\ps}$, 
$B=S$, $I=\m_{R_{\thetabar}^{\ps}}B$,
and $G=\GL_{d}$; note that by
Corollary~\ref{cor:CWE stacks satisfy assumption for formal functions}, $B^G$ is a complete local ring, and
is in particular $I^G$-adically
complete.    
 \end{proof}

Finally, since the assumptions of Definition~\ref{defn: versal categories} are satisfied by the algebraic stack~$\fX_{\thetabar}$, 
we have the following specialization of %
Proposition~\ref{prop: versal formal functions}. %

\begin{prop}\label{prop:versal formal functions CWE stacks}
Let~$\widehat \cX_{\thetabar}$
be the completion of $\fX_{\thetabar}$ (or equivalently, $\cX_{\thetabar}$) at its unique closed point.
Then the functor
\[
\Coh(\cO_{\cX_{\thetabar}}) \to \Coh(\Ver / \widehat \cX_{\thetabar}),
\]
sending a coherent sheaf~$\cF$ to the coherent sheaf $(f: \Spf A \to \widehat \cX_{\thetabar}) \mapsto \widehat f^*\cF(\Spf A)$, is fully faithful. 
\end{prop}

 \begin{arem}\label{arem:fixed-det-versions-of-Galois-stacks}
  There are obvious analogues of these results %
  in the fixed determinant case, which are proved in exactly the same way
  (replacing the pseudodeformation ring with the fixed determinant
  version, and $\GL_d$ with~$\SL_d$; 
see also~\cite[\S 2.1]{JNWE}.)  
We will use these analogues in the main body of the paper. 
\end{arem}

\bigskip

\emergencystretch=3em

\printbibliography

\end{document}
